\documentclass[a5paper,10pt,twoside,openany,final]{extbook}

\usepackage[utf8]{inputenc}
\usepackage[T2A, TS1,T1]{fontenc}
\usepackage[english, russian]{babel}

\usepackage[pdftex]{geometry}
\usepackage{floatrow}
\usepackage{graphicx}
\usepackage{indentfirst}
\usepackage[noadjust]{cite}
\usepackage{amssymb}
\usepackage{amsmath}
\usepackage{amsthm}
\usepackage[hyphens]{xurl}
\usepackage[colorlinks=true,linkcolor=blue,citecolor=blue]{hyperref}
\usepackage[usenames, dvipsnames]{color}
\usepackage[numbered, framed]{mcode}
\definecolor{light-gray}{gray}{0.95}
\usepackage{fancyhdr}
\renewcommand{\chaptermark}[1]{\markboth{#1}{}}
\renewcommand{\sectionmark}[1]{\markright{\thesection\ #1}}
\fancypagestyle{plain}{%
	\fancyhead{}
}

\usepackage{float}
\usepackage{algorithm}
\usepackage{algorithmic}

\numberwithin{algorithm}{chapter}
\def\nk#1{{\color{black}#1}}

\def\ag#1{{\color{black}#1}}

\def \NN {\mathbb N}

\def \R {\mathbb R}

\usepackage{enumitem}
\makeatletter
\AddEnumerateCounter{\asbuk}{\russian@alph}{}
\makeatother

\graphicspath{{images/}{images/Steiner/}} 

\title{Оптимизация}   
\author{Author Name} 
\date{\today}

\everymath{\displaystyle}
\renewcommand{\le}{\leqslant} \renewcommand{\ge}{\geqslant}

\newtheorem{lemma}{Лемма}[chapter]
\DeclareMathOperator*{\argmin}{arg\,min}

\DeclareMathOperator\rk{rk}
\DeclareMathOperator\tr{tr}
\DeclareMathOperator\ext{ext}
\DeclareMathOperator\spa{span}
\DeclareMathOperator\diag{diag}
\DeclareMathOperator\sign{sgn}
\DeclareMathOperator\conv{conv}
\DeclareMathOperator\Pos{Pos}
\DeclareMathOperator\End{End}
\DeclareMathOperator\interior{int}
\DeclareMathOperator\Supp{supp}

\usepackage{epstopdf}
\usepackage[font=small, labelsep=period]{caption}

\usepackage{epigraph}
\usepackage{graphicx}
\usepackage{epstopdf}
\usepackage{amsfonts, amsmath, amsthm}
\usepackage{amssymb}
\usepackage{url}
\usepackage{cite}
\usepackage[framemethod=1]{mdframed}
\usepackage{mathtools}
\usepackage{makecell}
\usepackage{tikz}
\usepackage{dsfont}
\usepackage{verbatim}

\usepackage{thmtools}

\declaretheoremstyle[
    headfont=\bfseries, 
    bodyfont=\normalfont\itshape,
    spacebelow=\parsep,
    spaceabove=\parsep,
    mdframed={
            hidealllines=true,
            leftline=true,
            skipabove=\parsep, 
            skipbelow=\parsep } 
]{myteostyle}

\declaretheoremstyle[
    headfont=\bfseries, 
    bodyfont=\normalfont,
    spacebelow=\parsep,
    spaceabove=\parsep,
    mdframed={
            hidealllines=true,
            leftline=true,
            skipabove=\parsep, 
            skipbelow=\parsep } 
]{mydefinstyle}

\declaretheoremstyle[
    headfont=\bfseries, 
    bodyfont=\normalfont,
    spacebelow=\parsep,
    spaceabove=\parsep,
    mdframed={
            hidealllines=true,
            skipabove=\parsep, 
            skipbelow=\parsep } 
]{myremarkstyle}

\declaretheoremstyle[
    headfont=\bfseries, 
    bodyfont=\normalfont,
    spacebelow=8pt,
    spaceabove=8pt,
    mdframed={
            leftline=false,
            rightline=false,
            splittopskip=20pt, 
            skipabove = 16pt 
            } 
]{myexamplestyle}

\declaretheoremstyle[
    headfont=\bfseries, 
    bodyfont=\normalfont,
    spacebelow=\parsep,
    spaceabove=\parsep,
    mdframed={
            hidealllines=true,
            skipabove=\parsep, 
            skipbelow=\parsep } 
]{myexercisestyle}

\declaretheorem[style=myteostyle, numberwithin=chapter, name=Теорема]{teo}

\declaretheorem[style=mydefinstyle, numberwithin=chapter, name=Определение]{defin}

\declaretheorem[style=myremarkstyle, numberwithin=chapter, name=Замечание, qed=$\blacksquare$]{remark}

\declaretheorem[style=myexamplestyle, numberwithin=chapter, name=Пример 
]{example}

\declaretheorem[style=myexercisestyle, numberwithin=chapter, name=Упражнение, qed=\ensuremath{\triangle}
]{exercise}

\declaretheorem[style=myexercisestyle, numbered=no, name=Указание, qed=$\blacksquare$
]{instruction}

\declaretheorem[style=myexercisestyle, numbered=no, name=Утверждение, qed=$\blacksquare$
]{claim}

\declaretheorem[style=myteostyle, numberwithin=chapter, name=Следствие]{corr}

\DeclareRobustCommand{\qed}{%
  \ifmmode
    \mathqed
  \else
    \leavevmode\unskip
    \penalty 9999
    \hbox{}\nobreak\quad\hbox{\qedsymbol}
  \fi
}

\usepackage{tikz}

\usepackage[final]{pdfpages}

\newcommand{\leqarg}[1]{\ensuremath{\stackrel{\text{#1}}{\leq}}}

\newcommand*\circled[1]{\tikz[baseline=(char.base)]{
            \node[shape=circle,draw,inner sep=2pt] (char) {#1};}}

\makeatletter 
\def\@biblabel#1{#1. } 
\makeatother

\usepackage{tempora} 
\usepackage{newtxmath} 

  \DeclareCaptionLabelFormat{rightline}{\rightline{\bothIfFirst{#1}{ }#2}}
\begin{document}


%
%
%

\thispagestyle{empty}


\begin{center}
	{\small
		Министерство науки и высшего образования  Российской Федерации\\[3pt]
		
		Федеральное государственное автономное  
		образовательное учреждение \\ высшего образования		
		<<Московский физико-технический институт\\ 
		(национальный исследовательский университет)>>
		
	}
	
	\vspace{2cm}
	
	
	
	{\huge {\cal Е. А. Воронцова, Р. Ф. Хильдебранд, \\[5pt] А. В. Гасников,  Ф. С. Стонякин}}

	\vspace{1.1cm}
	
	

%

{\fontsize{37}{46}\selectfont 
	{\bfseries Выпуклая\\[9pt] оптимизация}}
	
	\vspace{20mm}
	
	{\LARGE Учебное пособие
		
	}
	
	%
	%
	


	
	
	\vfill 
	
	МОСКВА \\ МФТИ \\ 2021
\end{center}


\clearpage

\thispagestyle{empty}
\begin{flushleft}
	{ УДК\hspace{3pt}519.853.3(075.8)\\
		ББК\hspace{4pt}22.21я73} \\
	\hspace{21pt}{В92}
\end{flushleft}

\vspace*{-4mm}




\begin{center}
	
	{Авторы: Е. А. Воронцова, Р. Ф. Хильдебранд,\\  А. В. Гасников,  Ф. С.~Стонякин}\\[7pt]

	{  Рецензенты:\\[3pt]
		
		Доктор физико-математических наук,\\  профессор Технологического института Джорджии (Атланта) \\ \emph{А.\,С. Немировский}\\[3pt]
		
		Доктор физико-математических наук, чл-корр. РАН,\\ профессор Высшей школы экономики (Москва) \\ \emph{В.\,Ю. Протасов}\\[3pt]
		
		
	}
	
\end{center}

\vspace{1.0mm}

{\normalsize
	
	\hspace{20.5pt}{\bfseries Выпуклая~оптимизация\,:}~учебное~пособие~/~Е.\;А.\;Воронцова,\linebreak 
	В92\hspace{11pt}Р. Ф. Хильдебранд, А. В. Гасников,  Ф. С. Стонякин. --\,Москва\,: МФТИ,\linebreak
	\hspace*{27,0pt}2021. --\,364~с. --\,ил. --\,21 см. --\,Библиогр.\,: 326~назв. 
	
	\hspace{20.5pt}ISBN 978-5-7417-0776-0 
	
	\vspace{1.5mm}
	
}

{


	
	{В основном \nk{рассмотрены} выпуклые задачи оптимизации. Связанно это с тем, что, во-первых, на данный момент только для выпуклых задач имеется достаточно богатая теория и, во-вторых (и это, наверное, главное), только для выпуклых задач  \nk{данная} теория позволяет разрабатывать эффективные методы на практике. Конечно, многие практические задачи оптимизации не являются выпуклыми. Однако даже для таких задач \nk{учебное пособие} является тем фундаментом, без которого едва ли возможно какое-то продвижение в их решении.}
	
	{\nk{Предназначено} для студентов 3 курса физтех-школы ПМИ МФТИ, изучающих \nk{курс <<Оптимизация>>}. С момента основания ФУПМ МФТИ (более 50 лет назад) оптимизационный курс занимал особое положение в большом цикле дисциплин, которые должны были освоить студенты в бакалавриате.} 
	

}


\hspace{210pt} {\small\bf УДК 519.853.3(075.8)}

\hspace{210pt} {\small\bf ББК 22.21я73}


\vfill

\begin{center}
	
	{\fontsize{8}{9}\selectfont \it	Печатается по решению Редакционно-издательского совета Московского 
		физико-технического института {\rm(}национального исследовательского университета{\rm)}
		
	}
	
\end{center}

\vfill




{\footnotesize \parindent=0mm  
	
	
	{\bf ISBN 978-5-7417-0776-0}\hspace{4.995em}\copyright\,Воронцова Е. А., Хильдебранд Р. Ф., 
	
	\hspace{16.4em}Гасников А. В.,  Стонякин Ф. С., 2021
	
	
	\hspace{15.43em}\copyright\,Федеральное государственное автономное
	
	\hspace{16.4em}образовательное учреждение высшего образования 
	
	
	\hspace{16.4em}<<Московский физико-технический институт
	
	\hspace{16.4em}(национальный исследовательский университет)>>, 2021

}

\newpage
\thispagestyle {empty}

\vspace*{1cm}

\begin{center}
	\large{\parbox{11.5cm}{
		\begin{raggedright}
{\large 
			{\parindent=6.3mm Optimization has an unusual status in~mathematics. It really has~to stand on three legs. One leg is some kind of basic theory like convex analysis. Another leg is the understanding of the various ways of formulating problems, what are the important things, not important things~--- in~other words, artful mathematical modelling. The~third leg on which optimization stands is computation. All three interact deeply. Computation is~often based on optimality conditions, which come from analysis and~especially duality.}
		}
	
		\vspace{.5cm}\hfill{{\large\it\ag{\it  R. Tyrrell} Rockafellar}\hspace{2mm}\cite{Roc_interview}}\vspace{1cm}

	
	Оптимизация занимает необычное место в 
	\ag{математике.}
\ag{В~действительности, она}
	должна стоять на трёх китах. Один кит~--- это \ag{некоторая}
	базовая теория вроде выпуклого анализа. Второй кит~--- это понимание различных 
	способов 
	\ag{постановки}
	задач и того, какие объекты важные, какие второстепенные, иначе говоря, \nk{---} искусство математического моделирования. Третий кит, на котором стоит оптимизация, --- это численные методы. 
	\nk{Все три существенно взаимосвязаны.}
	В вычислениях часто используют условия оптимальности, доказательство последних, в свою очередь,
	основано на анализе и, в частности, теории двойственности.
		
		\vspace{.5cm}\hfill{{\large\it\nk{\it Р. Тиррелл Рокафеллар} }\hspace{2mm}\cite{Roc_interview}}
	\end{raggedright}
	}
}
\end{center}

\newpage 	

\tableofcontents	


\chapter*{Введение\markboth{Введение}{Введение}}\addcontentsline{toc}{chapter}{\hspace*{5.4mm}Введение\dotfill}

Данное пособие написано на основе лекций, которые авторы читали в разное время в Физтех-школе ПМИ МФТИ (Москва), ФКН ВШЭ (Москва), ДВФУ (Владивосток), КФУ имени В.И. Вернадского (Симферополь), АГУ (республика Адыгея), университете Гренобль-Альпы (Гренобль, Франция).

В первую очередь при написании пособия авторы ориентировались на программу двухсеместрового курса лекций по оптимизации, который читается студентам Физтех-школы ПМИ МФТИ. Первая глава пособия содержит материалы первого семестра (<<Основы выпуклого анализа и оптимизации>>), вторая и третья главы~--- материалы второго семестра (<<Численные методы выпуклой оптимизации>>). В пособие не вошли некоторые материалы, которые включают сейчас в лекции второго семестра (например, элементы невыпуклой оптимизации и стохастической оптимизации). Эти материалы подробно 
изложены в книге~\cite{Gas_Pos18}, более точно заточенной на программу второго семестра университетского курса по методам оптимизации.

Пособие имеет ряд особенностей, о которых стоит заранее предупредить читателя. Во-первых, в отличие от классических пособий, в которых принято доказывать все основные факты, в данном пособии приводятся доказательства далеко не всех упоминаемых результатов. Это позволило, с~одной стороны, отметить больше связей и описать больше конструкций, с~другой стороны, сделало изложение менее самодостаточным. Вторая важная особенность~--- часть материала пособия является продвинутой и публикуется в учебной литературе, по-видимому, впервые. В-третьих, акценты, которые делаются в пособии, далеко не всегда совпадают с общепринятыми акцентами в популярных сейчас учебниках. Речь идет, прежде всего, о~достаточно продвинутом изложении конической оптимизации, в том числе робастной оптимизации, как яркой демонстрации возможностей современного выпуклого анализа, а также теории ускоренных градиентных методов. В~пособие вошли последние достижения в этой области, в том числе касающиеся разработки ускоренных тензорных методов. Особо отметим примеры, которые были собраны в отдельной главе. На наш взгляд, такие примеры должны наглядно продемонстрировать: 1) роль выпуклого анализа в постановке (формализации) задачи; 2) роль современных численных методов в~решении задач выпуклой оптимизации.  

Авторы выражают благодарность за помощь в оформлении данного пособия студентам кафедры математических методов в экономике ДВФУ Владиславу Хазову, Марии Рудь, Алексею Луценко, Валерии Волосевич, Екатерине Шавыровой, Алине Влучко, Марии Романовой, Екатерине Маникаевой, Валентине Шальтис, студентам кафедры алгебры и функционального анализа КФУ им. В.И. Вернадского Сейдамету Аблаеву и Константину Алексееву, а также старшему преподавателю кафедры алгебры и функционального анализа КФУ им. В.И. Вернадского Инне Баран. 

Ряд ценных замечаний был получен от Дарины Двинских (МФТИ, Институт Вейерштрасса, ИППИ РАН), Анастасии Ивановой (ВШЭ, университет Гренобль-Альпы), Адриена Тейлора (INRIA), Павла Двуреченского (Институт Вейерштрасса, МФТИ,  ИППИ РАН), Александра Тюрина (Научно-технологический университет имени короля Абдаллы, Саудовская Аравия), Aлександра Титова (МФТИ), Дмитрия Пасечнюка (МФТИ) и Артема Комарова (МФТИ). Также хотелось бы выразить благодарность профессору Е.\,А.~Нурминскому (ДВФУ, Владивосток), профессору А.\,С.~Немировскому 
(\ag{\nk{Технологический} институт Джорджии, Атланта, США), профессору Б.\,Т.~Поляку (ИПУ РАН, Москва), профессору Б.\,Ш.~Мордуховичу (университет Уэйна, США) и профессору А.\,В.~Назину (МФТИ)}.

Работа Е.\,А.~Воронцовой и Ф.\,С.~Стонякина над пособием была поддержана грантом РФФИ 18-29-03071 мк. Работа А.\,В.~Гасникова над пособием была выполнена при поддержке Министерства науки и высшего образования Российской Федерации (госзадание) №\,075-00337-20-03, номер проекта 0714-2020-0005.

 Если вы обнаружили
какие-либо опечатки, ошибки или неточности  в данном пособии, большая просьба сообщить
о них по адресу \href{mailto:vorontsovaea@gmail.com}{vorontsovaea@gmail.com}.



\chapter*{Обозначения и сокращения\markboth{Обозначения и сокращения}{Обозначения и сокращения}}\addcontentsline{toc}{chapter}{\hspace*{5.4mm}Обозначения и сокращения\dotfill}


\noindent В пособии используются следующие сокращения:
\begin{itemize}
\item[] АД --- автоматическое дифференцирование; 

\item[] БГМ --- быстрый градиентный метод; 

\item[] БШС --- модель вычислений Блюм--Шуба--Смейла;

\item[] ЛП --- линейное программирование;

\item[] МСГ --- методы сопряжённых градиентов;

\item[] МСП --- многослойный перцептрон;

\item[] МЦТ --- метод центров тяжести;

\item[] НС --- нейронная сеть;

\item[] ОГМ --- ограниченная машина Больцмана;

\item[] ОТ --- оптимальный транспорт;

\item[] РО --- робастная оптимизация

\end{itemize}

\noindent и обозначения:

\begin{itemize}
\item[] $\boldsymbol{1}$ --- вектор из единиц;

\item[] $\| a \|_\infty = \max\limits_i |a_i|$ --- бесконечная норма вектора~$a$;

\item[] $A \succ 0$ --- матрица $A$ -- положительно определена;

\item[] $A \succeq 0$ --- матрица $A$ -- неотрицательно определена;

\item[] $A \prec 0$ --- матрица $A$ -- отрицательно определена;


\item[] $\mathbb{A}^m$ --- множество вещественных
кососимметричных матриц размера~$m \times m$;

\item[] $\textnormal{Arg}\min_{x \in Q}\psi(x)$ --- множество точек минимума $\psi(x)$ на множестве $Q$;

\item[] $\arg\min_{x \in Q} \psi(x)$ --- некоторый элемент множества $\textnormal{Arg}\min_{x \in Q} \psi(x)$;

\item[] $B_{\|\cdot\|}(x_0, \, r)$ ---
шар с центром в точке $x_0 \, \in \, \mathbb{R}^n$ радиуса~$r$ (множество вида
$\{ x \, \in \, \mathbb{R}^n \, | \, \|x - x_0 \| \le r \}$, $r \ge 0$);

\item[] $\mathbb{C}$ --- поле комплексных чисел;

\item[] $\mbox{conv } M$ --- выпуклая оболочка множества
$M$;

\item[] $\text{diag}(x)$ --- диагональная матрица с компонентами вектора $x$ на диагонали;

\item[] $\dim \, x$ --- размерность объекта $x$;

\item[] $\text{dom } f$ --- область определения функции $f$;

\item[] $\mbox{epi } f = \{ (x, \, \alpha) \, \in \, \mbox{dom} \, f \, \times \, \mathbb{R} \, | \, \alpha \, \ge \, f(x) \}$ --- надграфик функции $f$;

\item[] $F[x]$ --- кольцо многочленов от конечного числа переменных с коэффициентами из поля $F$;

\item[] $\mathbb{H}^m$ --- множество комплексных эрмитовых матриц размера~$m \times m$;

\item[] $\mathbb{H}_+^m$ --- конус неотрицательно определённых комплексных эрмитовых матриц размера~$m \times m$;

\item[] $I$ --- единичная матрица;

\item[] $\mbox{int} \, X$, $X^o$ --- внутренность множества $X$;

\item[] $nnz \left ( A \right)$~--- число  ненулевых  элементов  в матрице $A$; 

\item[] $O(n)$ --- неотрицательная функция от $n$, которую можно ограничить сверху $Cn$, где $C$ -- числовая константа;

\item[] $\tilde{O}(n)$ --- неотрицательная функция от $n$, которую можно ограничить сверху $C\log^{r}(n)n$, где $C$ и $r$ -- числовые константы (аналогичный смысл мы будем вкладывать в $O$ и $\tilde{O}$ в случае большего числа параметров);

\item[] $\mathbb{R}^n$ --- $n$-мерное вещественное евклидово пространство;

\item[] $\mathbb{R}_{+}$ --- множество неотрицательных вещественных чисел;

\item[] $\mathbb{R}_{++}$ --- множество положительных вещественных чисел;

\item[] $\bar{\mathbb{R}} = \mathbb{R} \,  \cup \, \{-\infty\} \, \cup \, \{+\infty\}$; 

\item[] $\text{rank } A$ --- ранг матрицы~$A$;

\item[] $\mathbb{S}^m$ --- множество вещественных
симметричных матриц размера~$m \times m$;

\item[] $\mathbb{S}_+^m$ --- конус неотрицательно определённых вещественных
симметричных матриц размера~$m \times m$;

\item[] $\text{Tr}(A)$ --- след квадратной матрицы~$A$, т.е. сумма диагональных элементов~$A$;

\item[] $\langle x, y \rangle$ --- скалярное произведение
векторов ${x = (x_1, \, \ldots, \, x_n)  \in  \mathbb{R}^n}$\linebreak и~$y = (y_1, \, \ldots, \, y_n) \, \in \, \mathbb{R}^n$. Также иногда для обозначения скалярного произведения используется
и вариант $x^T y$;

\item[] $x \, \ge \, y$ --- все компоненты вектора $x \, \in \mathbb{R}^n$ не меньше соответствующих компонент вектора $y \, \in \, \mathbb{R}^n$, т.е. $x_i \, \ge \, y_i$, $i = 1, \, \ldots, \, n$ (неравенства такого рода являются покомпонентными). Часто вместо вектора $y$ в таких неравенствах стоит $0$, что обозначает нулевой вектор соответствующей размерности; 

\item[] $\mathbb{Z}$ --- множество целых чисел;

\item[] $\lambda _{\max } \left( A \right)$ --- максимальное собственное значение матрицы~$A$.

\end{itemize}




\chapter[Основы выпуклого анализа\dotfill]{Основы выпуклого анализа\markboth{Основы выпуклого анализа}{Основы выпуклого анализа}}

В данной главе рассматриваются в основном задачи \textit{выпуклой оптимизации} (минимизации выпуклой функции на выпуклом множестве). Описан \textit{принцип множителей Лагранжа} как следствие \textit{теоремы об <<отделимости>>} граничной точки выпуклого множества от самого множества (опорной) гиперплоскостью. При этом конечномерность пространств не предполагается, см., например, обоснование с помощью принципа множителей Лагранжа леммы Неймана--Пирсона в разделе~\ref{neuman_pirson} (основной леммы в проверке статистических гипотез). Обсуждаются различные моменты, которые активно используются при разработке эффективных численных методов решения задач выпуклой оптимизации (и их удобных переформулировок), о чём пойдёт речь в том числе и в последующих главах. Особенно отметим: \textit{теорему фон\;Неймана--Сиона--Какутани} в разделе~\ref{minmax} ($\min \max = \max \min$), \textit{теорему (формулу) Демьянова--Данскина--Рубинова} о дифференцировании под знаком максимума (минимума) в разделе~\ref{ch1_sect_dem-dan}, \textit{теорию двойственности} в разделе~\ref{ch1_sect_extrem} и \textit{теоремы об альтернативах} в разделе~\ref{ch1_sect_fark}. Отличительной особенностью материала данной главы от содержания многих других русскоязычных учебных пособий по выпуклой оптимизации является
наличие заключительных разделов, которые содержат весьма продвинутый материал. Если \textit{коническое программирование} уже достаточно прочно вошло в различные западные курсы по выпуклой оптимизации, то \textit{робастную оптимизацию}, насколько нам известно, в основном излагают только в виде отдельных спецкурсов. В этом пособии мы, следуя А.\,С.~Немировскому, стараемся продемонстрировать, что именно такие приложения позволяют по-настоящему прочувствовать все возможности арсенала современной выпуклой оптимизации.  

Для более глубокого погружения в затронутые в данной главе вопросы можно рекомендовать замечательные книги по выпуклой оптимизации В.\,М.~Тихомирова, А.\,С.~Немировского\footnote{У Владимира Михайловича Тихомирова и у Аркадия Семёновича Немировского довольно много разных книг на тему выпуклой оптимизации, см., например, \cite{MagTih_convan, Ben-Tal}. У Владимира\linebreak Михайловича акцент делается на общих вопросах выпуклого анализа, у Аркадия Семёновича~--- на вычислительных. В совокупности это в хорошей степени покрывает материал, отвечающий глубокому изучению данного курса.} и, конечно, книгу Стивена Бойда и~Ливена Ванденберга~\cite{Boyd_book_2004}, которую можно рекомендовать всем желающим разобраться \nk{в том}, что такое <<выпуклая оптимизация>>. Вот уже почти 20 лет эта книга служит непревзойдённым по доступности изложения источником базовых фактов о выпуклой оптимизации.
\ag{Также можно порекомендовать книги
Ю.\,Е.~Нестерова~\cite{Nest}, Б.\,Н.~Пшеничного~\cite{Pshen_80}, Б.\,Т.~Поляка~\cite{Polyak}.}

\section{Элементы топологии}
В этом параграфе собраны некоторые определения и факты~\cite{AvakovMagaril, Verb_Topol}, которые необходимы
для доказательства теорем, представленных далее.
\begin{defin}
Векторное пространство~$X$ называется вещественным {\it нормированным пространством}, если каждому элементу $x \, \in \, X$ сопоставлено вещественное число $\| x \|_X$,
называемое \emph{нормой} $x$, которое удовлетворяет
условиям:
\begin{enumerate}
\renewcommand\labelenumi{\textup{\arabic{enumi})}}
    \item {\normalfont{невырожденность}}: для любого $x  \in  X$ имеем $\|x\|_X  \ge  0$ и $\|x \|_X = 0$ тогда и только тогда, когда $x = 0$;
    \item  {\normalfont{мультипликативность}}: для любого $\alpha \in  \mathbb{R}$ и
    $x  \in  X$ имеем
    $\| \alpha x \|_X =\linebreak= |\alpha| \|x\|_X$;
    \item  {\normalfont{неравенство треугольника}}: для любых $x, \, y  \in  X$ имеем $\| x + y \|_X  \le \linebreak \le \|x \|_X + \|y \|_X$.
\end{enumerate}
\end{defin}

{\it Евклидова норма} вектора $x = (x_1, \, x_2, \, ..., \, x_n)$ (её также называют {\it квадратичной})
определяется как
$$
\|x \|_2 = (x^T x)^{1/2} = \left ( \sum_{i = 1}^n |x_i|^2 \right )^{1/2}.
$$

Неравенство треугольника для евклидовой нормы называется {\it неравенством Коши--Буняковского},
оно же в нерусскоязычной литературе называется
{\it неравенством Коши--Шварца} (Cauchy--Schwarz inequality).
Действительно, возьмём два ненулевых неколлинеарных вектора $x$, $y$ в векторном пространстве~$X$ с положительно определённым скалярным произведением~$x^T y$. Неравенство треугольника
для пространства~$X$ записывается так:
\begin{equation}
\label{ch1_eq_trian_ineq_eucl1}
\sqrt{(x + y)^T (x + y)} \, \le \, \sqrt{x^T x} + \sqrt{y^T y}.
\end{equation}
Возводим обе части~\eqref{ch1_eq_trian_ineq_eucl1}
в квадрат и раскрываем скобки, тогда получим эквивалентное утверждение
\begin{equation}
\label{ch1_eq_trian_ineq_eucl2}
x^T y \, \le \, \sqrt{(x^T x) (y^T y)}.    
\end{equation}
Рассмотрим скалярное произведение $(x - \lambda y)^T (x - \lambda y)$, $\lambda \, \in \, \mathbb{R}$, $x - \lambda y \neq 0$.
Поскольку скалярное произведение положительно, то
квадратичный полином
$$
P(\lambda) = y^T y \lambda^2  - 2  x^T y \lambda
+ x^T x
$$
не имеет корней, следовательно, его дискриминант
$4 (x^T y)^2 - 4 (y^T y) (x^T x)$ отрицателен,
что доказывает~\eqref{ch1_eq_trian_ineq_eucl2}.

Евклидова норма является частным случаем векторной {\it нормы Гёльдера}, задаваемой
выражением
$$
\|x \|_p = \left ( \sum_{i = 1}^n |x_i|^p \right )^{1/p},
$$
где $x \, \in \, \mathbb{R}^n$, $p \, \ge \, 1$.

Ещё один важный пример нормы~--- $\max$-норма (также её называют {\it бесконечная норма}, $\sup$-норма или $\ell_\infty$):
$$
\|x \|_\infty = \lim_{p \rightarrow \infty} \|x \|_p = \max_{i = 1, \, \ldots, \, n} |x_i|.
$$

\begin{exercise}
Докажите, что $\max$-норма действительно является нормой.
\end{exercise}

Наконец, $\ell_1$-норма, или
просто {\it $1$-норма} (частный случай $p$-нормы Гёльдера при $p = 1$), определяется так:
$$
\| x \|_1 = \sum_{i = 1}^n |x_i|.
$$
\begin{exercise}
Докажите, что $\ell_1$-норма действительно является нормой.
\end{exercise}

Упомянутые\footnote{Вообще говоря, не только упомянутые нормы, но и абсолютно все на конечномерном пространстве.} векторные нормы эквивалентны\footnote{Две нормы $\| \cdot \|_{k_1}$ и $\| \cdot \|_{k_2}$ в 
пространстве~$X$ называются \textit{эквивалентными}, если существуют такие положительные числа $C_1$ и $C_2$, что для любого $x \, \in \, X$ выполнены неравенства
$C_1 \|x\|_{k_1} \leq\linebreak\leq \| x \|_{k_2} \leq C_2 \| x \|_{k_1}$.}
на конечномерном пространстве.

Норма на пространстве~$X$ задаёт метрику
(\textit{расстояние между элементами $x$ и $y$}) по формуле $\| x - y \|_X$ при $x, \, y \, \in \, X$.

Расстоянием от точки $a \, \in \, X$ до
множества $B \subset X$ называется
величина $\inf\limits_{x \, \in \, B}  \| x - a\|_X$.

Пусть $Q \subset X$. Точка $x \, \in \, Q$
называется \textit{внутренней} точкой $Q$,
если она входит в $Q$ вместе с некоторой окрестностью
с центром в $x$. Множество внутренних точек множества $Q$ называется \textit{внутренностью}~$Q$ и обозначается~$\mbox{int} \, Q$ или $Q^o$.

Множество $Q$ называется \textit{открытым}, если $\mbox{int} \, Q = Q$.

Точка $x^* \, \in \, X$ называется \textit{граничной} точкой множества $Q$, если в любой окрестности точки $x^*$ содержатся точки как из множества $Q$, так и из его дополнения.
Множество называется \textit{замкнутым}, если
ему принадлежат все граничные точки.
Известно, что дополнение к \textit{замкнутому}
множеству является открытым. И наоборот,
дополнение к открытому множеству является
замкнутым.

\section[Выпуклые множества и выпуклые функции]{Выпуклые множества и выпуклые функции}\label{ch1_sect_conv}

Выпуклый анализ~--- раздел математики, в котором изучают
выпуклые множества, выпуклые функции и выпуклые
экстремальные задачи.  
Исторические корни понятия выпуклости уходят в античные времена.
Определение выпуклой кривой
присутствует, например, в сочинениях Архимеда (III~век до н.э.). В дальнейшем активное изучение выпуклых фигур началось в
XIX~веке, в трудах
О.~Коши (один из первых результатов~ --- теорема Коши о
жёсткости выпуклых многогранников, 1813~год), Я.~Штейнера и
Г.~Минковского 
возникла отдельная ветвь геометрии~--- выпуклая геометрия. Окончательно
новой областью математики 
выпуклая геометрия стала после выхода
книги Минковского\footnote{{\it Minkowski H.} Geometrie der Zahlen. Leipzig: Teubner Verlag, 
1896, 1910. Последующие переиздания: New York: Chelsea, 1953; New York: Johnson, 1968.}. 
Так были заложены геометрические основания выпуклого анализа.

В двадцатые и тридцатые годы ХХ~века бурное развитие теории выпуклых множеств продолжилось. Итоги этого периода были подведены в 
вышедшей в 1934~году монографии немецких математиков Боннесена и Фенхеля <<Теория выпуклых тел>>. 
Многие значимые идеи выпуклой геометрии 
стали основой для важных положений функционального анализа (теоремы Хана--Банаха, Крейна--Мильмана). 

Понятие выпуклости стало привлекать особое внимание исследователей в 40-е годы XX~века, когда начал происходить настоящий <<выпуклый бум>>,
по выражению В.\,М.~Тихомирова~\cite[стр.~3]{Tikh_kvant03}, поскольку выяснилось, какую огромную роль играет выпуклость в экономических задачах,
следствием решения которых стало рождение нового направления в теории экстремума~ ---
{\it линейного 
программирования}\footnote{Подробнее
о развитии линейного программирования
см. п.~\ref{subsubsect_lp} на с.~\pageref{subsubsect_lp}.}. 

Результаты первооткрывателя нетривиальных задач и теории линейного программирования
Л.\,В.~Канторовича
были переоткрыты на Западе, там было осознано
значение выпуклых экстремальных задач при решении
актуальных проблем экономики и военно-промышленного
комплекса, и многие исследователи приняли участие
в развитии новой дисциплины
на стыке математического анализа и геометрии.
В 1951~г. выходит издание
курса лекций Фенхеля по теории
выпуклых функций и множеств
({\it Fenchel W.} Convex Cones, Sets and Functions.
Princeton University Press, 1951).
Название окончательно
оформившейся новой дисциплины~--- выпуклый анализ~--- предложил
профессор Принстонского университета А.\,У.~Таккер, что описано в предисловии
к подводившей итоги 20-летнего периода
развития теории
и вышедшей в США в 1970~г.
 монографии 
американского математика Р.~Рокафеллара <<Выпуклый анализ>>~\cite{Rocaf}.
Рокафеллар был пионером в применении выпуклого анализа к оптимизации.

\begin{defin}
Множество $ Q$ называется {\it выпуклым}, если для любых двух точек, принадлежащих этому множеству, отрезок, соединяющий эти две точки, целиком принадлежит этому множеству, т.е.
$$
\forall \, x, \, y \, \in \, Q, \, \forall \, \alpha \, \in \, [0, \, 1] \, \rightarrow \, \alpha x + (1-\alpha) y \, \in \, Q. 
$$
\end{defin}
По индукции легко доказывается, что выпуклое множество~$Q$ вместе с любыми своими точками $q_1, \, \ldots, \, q_n$ содержит
и их выпуклую комбинацию $\sum\limits_{i = 1}^n \alpha_i q_i$,
где все $\alpha_i \, \in \, \mathbb{R}$, $\alpha_i \geq 0$, $i = 1, \, \ldots, \, n$, $\sum\limits_{i = 1}^n \alpha_i = 1$.

Пустое множество --- выпукло.



\begin{figure}[h!]
	\begin{floatrow}
		\ffigbox
		{\includegraphics[width=0.4\textwidth]{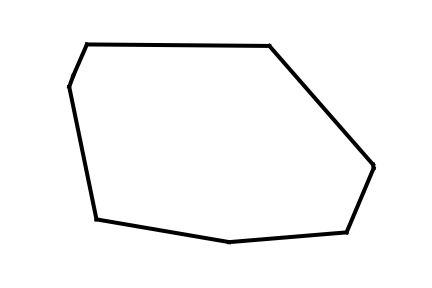}}{\caption{Пример выпуклого множества}}
		\ffigbox
		{\includegraphics[width=0.4\textwidth]{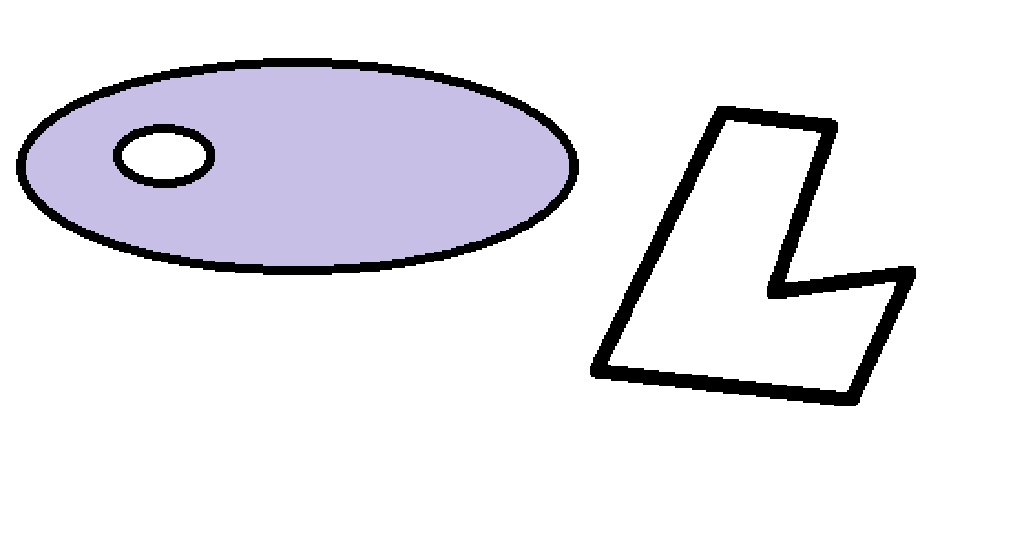}}{\caption{Пример невыпуклых множеств}}       
	\end{floatrow}
\end{figure}


\begin{example}
Шар с центром в точке $x_0 \, \in \, \mathbb{R}^n$ радиуса~$r$, с определённой
нормой $\|\cdot \|$ на пространстве~$\mathbb{R}^n$ (т.е. множество вида
$\{ x  \in \mathbb{R}^n \, | \, \|\,x -\linebreak -\, x_0 \| \le r \}$, $r \ge 0$) является
выпуклым множеством. Здесь и далее будем обозначать
такой шар $B_{\|\cdot\|}(x_0, \, r)$. 
\end{example}

\begin{exercise}
\label{excs_1_ball}
Докажите, пользуясь определением выпуклого множества, что единичный шар с центром в нуле~$B_{\|\cdot\|}(0, \, 1)$ является выпуклым.
\end{exercise}

\begin{instruction}[\nk{\cite[лекция~2]{Verb_Topol}}]
Представим отрезок $[x, \, y]$ как множество точек вида~$\alpha x + (1 - \alpha) y$, $\alpha \, \in \, [0, \, 1]$. По определению
выпуклого множества при $x, \, y \, \in \, B_{\|\cdot\|}(0, \, 1)$ все точки
 вида~$\alpha x + (1 - \alpha) y$, $\alpha \, \in \, [0, \, 1]$
 тоже должны принадлежать этому шару, т.е. должно выполняться
 неравенство~$\| \alpha x + (1 - \alpha) y\| \leq 1$.
 В силу неравенства треугольника для нормы и мультипликативности нормы
 $$
 \| \alpha x + (1 - \alpha) y\| \leq
 \| \alpha x \| + \| (1 - \alpha) y\| \leq
 \alpha \|  x \| + (1 - \alpha)  \|  y\|.
 $$
 Поскольку $x$, $y$ принадлежат единичному шару,
 $$
 \alpha \|  x \| + (1 - \alpha)  \|  y\| \leq \alpha + (1 - \alpha)
 = 1.
 $$
 Следовательно, $\| \alpha x + (1 - \alpha) y\| \leq 1$.
\end{instruction}

Образно говоря, множество является выпуклым, если до каждой точки этого множества можно добраться из любой другой точки по прямому непрерывающемуся пути, целиком лежащему в этом множестве.

Пересечение любой совокупности выпуклых множеств выпукло.
\begin{defin}
Пересечение всех выпуклых множеств, содержащих множество~$Q \subset \mathbb{R}^n$, называется \emph{выпуклой оболочкой} множества~$Q$
и обозначается $\conv(Q)$.
\end{defin}

\begin{defin}
Функция $f \, : \, \mathbb{R}^n \, \rightarrow \, \mathbb{R}$ называется \emph{выпуклой}, если
\emph{$\mbox{dom}$} $f$ --- выпуклое множество и
если для всех 
$x, \, y \, \in$ \emph{$\mbox{dom}$} $f$
выполняется неравенство Йенсена:
\begin{equation}
\label{conv_func_def}
\forall \, \alpha \, \in \, [0, \, 1] \, \rightarrow \, f \left (\alpha x + (1 - \alpha) y \right ) \le \alpha f(x) + (1 - \alpha) f(y). 
\end{equation}
\end{defin} 

Геометрически неравенство Йенсена означает, что отрезок, соединяющий точки $(x, \, f(x))$ и $(y, \, f(y))$, лежит не ниже графика функции $f(x)$. 

Сопоставим каждой функции $f : \mbox{dom} \, f \, \rightarrow \, \mathbb{R}$ множество
$$
\mbox{epi } f = \{ (x, \, \alpha) \, \in \, \mbox{dom} \, f \, \times \, \mathbb{R} \, | \, \alpha \, \ge \, f(x) \},
$$
которое называется {\it надграфиком} $f$.

Иногда говорят, что функция {\it выпукла}, если её надграфик является выпуклым множеством. Оба эти определения эквивалентны.

\begin{defin}
Функция $f$ называется \emph{строго выпуклой}, если
$$
\forall \, \alpha \, \in \, (0, \, 1) \ \rightarrow \ f \left (\alpha x + (1 - \alpha) y \right ) < \alpha f(x) + (1 - \alpha) f(y)
$$
для всех $x \neq y$. 
\end{defin}

\begin{defin}
Функция $f$ называется \emph{вогнутой}, если $-f$ --- выпуклая функция, и функция $f$ называется \emph{строго вогнутой} функцией, если $-f$ --- строго выпуклая функция.
\end{defin}

Анализ выпуклых функций является хорошо проработанной областью исследований. 
Одним из важных результатов является, например, следующий: функция $f$ является выпуклой
тогда и только тогда, когда для всех 
$x \, \in \, \mbox{dom} \, f$ и для всех $v$
функция $g(t) = f(x+tv)$ является выпуклой на своей области определения $\{t \, | \, x + tv \, \in \, \mbox{dom} \, f\}$. Это полезное свойство позволяет сводить проверку на выпуклость функции многих переменных к проверке на выпуклость функции одной переменной.

Будем обозначать здесь и далее положительную определённость матрицы $A$ таким образом: $A \succ 0$. Соответственно, если матрица неотрицательно определена, обозначение следующее: $A \succeq 0$. И наоборот, отрицательно определённая матрица: $A \prec 0$.

Пусть функция $f$ дважды дифференцируема, т.е. её матрица вторых производных (которую называют {\it гессиан}, или {\it матрица Гессе}) $\nabla^2 \, f$ существует в каждой точке области определения $\mbox{dom} \, f$. Тогда функция~$f$ выпукла тогда и только тогда, когда  $\mbox{dom} \, f$ --- выпуклое множество и гессиан функции $f$ --- неотрицательно определённая матрица:
$$
\forall \, x \, \in \,  \mbox{dom} \, f \quad \nabla^2 \, f \, \succeq \, 0.
$$
Для функции одной переменной это условие означает, что её первая производная является неубывающей функцией на интервале $\mbox{dom} \, f$: 
$$
f''(x) \, \ge \, 0.
$$

Аналогично, функция $f$ является вогнутой тогда и только тогда, когда $\mbox{dom} \, f$ --- выпуклое множество и $\nabla^2 \, f \, \preceq \, 0$
для всех $x \, \in \,  \mbox{dom} \, f$.

Для строгой выпуклости функции $f$ положительная определённость матрицы вторых производных является
только достаточным условием. Если $\nabla^2 f(x) \, \succ \, 0$  для всех $x \, \in \,  \mbox{dom} \, f$, то функция $f$ строго выпукла. Обратное утверждение неверно (контрпример --- строго выпуклая функция $f(x) = x^4, \, x \, \in \, \mathbb{R}$, имеющая в точке $x = 0$ вторую производную, равную $0$).

Как следует из определения, все аффинные\footnote{Функция $f : \, \mathbb{R}^m \, \rightarrow \, \mathbb{R}^n$
называется \textit{аффинной}, 
если существуют линейная функция 
$l: \, \mathbb{R}^m \, \rightarrow \, \mathbb{R}^n$
и вектор $b \, \in \, \mathbb{R}^n$
такие, что
$f(x) = l(x) + b$
для всех $x \, \in \mathbb{R}^m$.}, а следовательно, и линейные, функции являются одновременно и выпуклыми, и вогнутыми. Приведём ещё несколько примеров выпуклых и вогнутых функций.  
\begin{example}
Степенная функция $f(x) = x^a$ выпукла на $\mathbb{R}_{++}$, если $a \, \ge \, 1$ или $a \, \le \, 0$, и вогнута
при $0 \, \le \, a \, \le 1$.
\end{example}

\begin{example}
Экспонента $f(x) = e^{ax}$ выпукла
на $\mathbb{R}$ для любого $a \, \in \, \mathbb{R}$.
\end{example} 

\begin{example}
 Степень модуля $f(x) = |x|^p$,
при $p \, \ge \, 1$ выпукла на $\mathbb{R}$.
\end{example} 

\begin{example}
Логарифм $f(x) = \log x$ --- вогнутая функция на $\mathbb{R}_{++}$.
\end{example}

\begin{example}
Сумма $k$ максимальных компонент $f(x) = x_{[1]} + \ldots
+ x_{[k]}$ выпукла на~$\mathbb{R}^n$.
Здесь имеется в виду, что $x_{[i]}$ --- $i$-й по значению элемент вектора~$x$, $i = 1, \, 2, \, \ldots, \, n$:
$$
x_{[1]} \, \ge \, x_{[2]} \, \ge \, \ldots \,
\ge \, x_{[n]}.
$$
Данную функцию можно представить как максимум
из всех возможных сумм $k$~элементов вектора~$x$:
$$
f(x) = \sum_{i=1}^k x_{[i]} = \max \{ x_{i_1} + \ldots + x_{i_k} \, | \, 1 \, \le \, i_1 \, < \, i_2 \, < \, \ldots \, < \, i_k \, \le \, n \}.
$$
\end{example}

\begin{example}
Функция $f\left( x \right)=\log \left( {\sum\limits_{k=1}^n {\exp \left( x_k  \right)} } \right)$
выпукла на $\mathbb{R}^n$.
Это неочевидное утверждение можно доказать, проверив, что гессиан функции неотрицательно определён.
Действительно \cite{Boyd_book_2004}, вычислим гессиан~$f(x)$:
\begin{equation}
\label{ch1_hes_log-sum-exp}
\nabla^2 f(x) = 
\frac{1}{(\boldsymbol{1}^T z)^2} \left ( \left ( \boldsymbol{1}^T z \right ) \mbox{diag} (z) - z z^T \right )
, 
\end{equation}
где $z = (e^{x_1}, \, \ldots, \, e^{x_n})$, и
здесь, и далее $\text{diag}(x)$~--- диагональная матрица с компонентами вектора~$x$ на главной диагонали.

Для того чтобы доказать неотрицательную определённость гессиана $f$, нужно воспользоваться определением неотрицательно определённой матрицы и доказать следующее свойство: для всех $v$ выполняется
неравенство $v^T \nabla^2 f(x) v \, \ge \, 0$.
После подстановки $\nabla^2 f(x)$ из (\ref{ch1_hes_log-sum-exp}) получаем
$$
\frac{1}{(\boldsymbol{1}^T z)^2} \left ( \left ( \sum_{i=1}^n z_i \right ) \left ( \sum_{i=1}^n v_i^2 z_i \right ) - \left ( \sum_{i=1}^n v_i z_i \right )^2 \right ) \, \ge \, 0.
$$ 
Последнее неравенство следует из неравенства Коши--Буняковского:
$$
(a^T a)(b^T b) \, \ge \, (a^T b)^2,
$$
в котором нужно взять $a_i = \sqrt{z_i}$, $b_i = v_i \sqrt{z_i}$.
\end{example} 

\begin{example}
Максимальное собственное значение симметричной матрицы как функция, определённая в пространстве симметричных матриц, является выпуклой функцией:
$$
f\left( X \right)=\lambda _{\max } \left( X \right),
\quad 
X - \mbox{матрица}, \, X=X^T.
$$
\end{example} 

\begin{example}
Ещё один неочевидный пример --- функция $f\left( X \right)=\log \det X^{-1}$, определённая на множестве
симметричных положительно определённых матриц, является 
выпуклой. Здесь и далее $\log$ обозначает логарифм числа
по любому основанию. При этом $f\left( X \right)=\log \det X$ при
$X\succ 0$~--- вогнутая функция.
\end{example} 

\begin{example}
{\it Матрично-дробная функция} $f\left( {x,Y} \right)= {x^T Y^{-1}x} $,
где $x \, \in \, \mathbb{R}^n$, а $Y$~--- положительно
определённая симметричная матрица размера $n \, \times \, n$,
является выпуклой по совокупности переменных.
\end{example} 

Традиционно
в теории оптимизации 
используются
\textit{условия Липшица} для производных определённого порядка.

\begin{defin}
Пусть множество~$Q$ является подмножеством $\mathbb{R}^n$.
Говорят, что $p$-я производная функции \textit{удовлетворяет условию Липшица}\protect\footnotemark   
{ на~$Q$} с константой~$L_{p} > 0$, если для всех~$x, \, y \, \in \, Q$
выполняется неравенство
$$
\| \nabla^p f(x) - \nabla^p f(y) \|_2
\leq L_{p} \| x - y \|_2.
$$
\end{defin}
\footnotetext{В частном случае, при $p = 0$, функцию
называют непрерывной по Липшицу: должна существовать
константа~$L_0 > 0$, такая, что для всех~$x, \, y \, \in \, Q$ выполняется неравенство
$$
| f(x) - f(y) |
\leq L_0 \| x - y \|_2.
$$
В главе~2 константа Липшица для функции
часто обозначается через $M$.
}

\begin{exercise}
\label{excs_1_lip1}
Пусть функция $f \, : \, \mathbb{R} \rightarrow \mathbb{R}$
дифференцируема и её производная ограничена. Докажите, что
$f$ непрерывна по Липшицу.
\end{exercise}

\begin{instruction}
Если функция дифференцируема на~$\mathbb{R}$, то она непрерывна на~$\mathbb{R}$. Для непрерывной функции выполняется теорема Лагранжа о среднем значении: для любых $a, \, b \, \in \, \mathbb{R}$,
таких, что $a < b$, найдётся такая точка $c \, \in \, (a, \, b)$,
для которой
$$
\frac{f(b) - f(a)}{b - a} = f'(c).
$$
По условию $f'$ ограничена. Следовательно, существует константа
$L_0 > 0$, такая, что $|f'(c)| \leq L_0$. Следовательно,
$$
\frac{ | f(b) - f(a) | }{ | b - a | }  \leq L_0,
$$
что и требовалось доказать.
\end{instruction}

По сути, для $L_0$-липшицевой функции $L_0$ является показателем
скорости изменения функции.
Кроме того, для непрерывной по Липшицу функции всегда существует
двойной конус с вершиной, которую можно перемещать
по графику функции в любом направлении, но при этом
для такой функции её график всегда
лежит вне конуса.

\begin{example}
Не все дифференцируемые функции являются непрерывными по Липшицу.
Рассмотрим, например, дифференцируемую функцию $f(x) = x^2$:
$$
\sup_{x \neq y \, \in \, \mathbb{R}} \frac{|f(x) - f(y) |}{| x - y |}
= \sup_{x \neq y \, \in \, \mathbb{R}} \frac{|x^2 - y^2 |}{| x - y |}
=  \sup_{x \neq y \, \in \, \mathbb{R}} {| x + y |} = \infty,
$$
поэтому $f$ не является $L_0$-липшицевой.
\end{example}

\begin{example}
Дифференцируемость не является и необходимым условием липшицевости
функции. Например, функция $f(x) = |x|$ недифференцируема,
но непрерывна по Липшицу.
\end{example}

Тем не менее верна следующая теорема. 

\begin{teo}[Радемахера, \cite{MakPodkor, Ambrosio}] Пусть $Q$~--- открытое
подмножество $\mathbb{R}^n$. Функция $f \, : \, Q \, \rightarrow \, \mathbb{R}$, удовлетворяющая условию Липшица на $Q$, дифференцируема почти везде на $Q$. 
\end{teo}

\begin{exercise}
\label{excs_1_lip2}
Докажите, что функция $f(x) = |x|$ непрерывна по Липшицу с $L_0 = 1$.
\end{exercise}

\begin{instruction}
Достаточно рассмотреть обратное неравенство треугольника
$$
\left | |b| - |a| \right | \leq | b - a |.\qedhere
$$
\end{instruction}

\section[Конусы и отношения порядка]{Конусы и отношения порядка в евклидовых\\ пространствах\sectionmark{Конусы и отношения порядка}}
\sectionmark{Конусы и отношения порядка}
Кроме уже определённых основных объектов выпуклого анализа~--- выпуклых множеств и выпуклых функций, важных для решения оптимизационных задач, ещё необходимо 
упомянуть выпуклые конусы. 

\begin{defin}
Непустое множество $K$ называется \emph{конусом}, 
если для любого $x \, \in \, K$ и  
$\lambda \, \ge \, 0$
выполняется включение~$\lambda x \, \in \, K$.
Иными словами, множество в линейном пространстве является 
конусом (с вершиной в точке~$0$), если с каждой своей точкой оно  
содержит и весь луч, исходящий из начала координат и проходящий 
через эту точку.
\end{defin}

Например, всё пространство~$\mathbb{R}^n$ является
конусом. Из определения следует, что точка~$0$ всегда принадлежит любому конусу.

\begin{defin}
Множество $K$ называется \emph{выпуклым конусом}, 
если $K$ --- конус и $K$ --- выпуклое множество.
\end{defin}

Луч~$\{x_0 + \theta v \, | \, \theta  \ge  0\}$,
где $v \neq 0$, является выпуклым конусом, если
$x_0 = 0$.

Назовём конус~$K \, \subseteq \, \mathbb{R}^n$ {\it правильным} или {\it регулярным}, если он 
удовлетворяет следующим свойствам:
\begin{enumerate}
\renewcommand\labelenumi{\textup{\arabic{enumi})}}
\item $K$~--- выпуклый;
\item $K$~--- замкнутый;
\item у $K$ непустая внутренность;
\item $K$~--- остроконечный, т.е. не содержит прямых:
$$
x \, \in \, K, \, -x \, \in \, K \Rightarrow x = 0.
$$
\end{enumerate}

Правильный конус можно использовать для
{\it обобщённых неравенств}, которые являются
отношениями частичного порядка в~$\mathbb{R}^n$:
$$
\forall  x \, \in \,  \mathbb{R}^n \, \,
\forall  y \, \in \,  \mathbb{R}^n \,: \,
x \succeq_K y \Leftrightarrow x - y \, \in \, K. 
$$
Вместо $y \succeq_K x$ можно написать: $x \preceq_K y$.
Аналогично определяется строгое обобщённое
неравенство:
$$
\forall  x \, \in \, \mathbb{R}^n \, \,
\forall  y \, \in \,  \mathbb{R}^n \,: \,
x \succ_K y \Leftrightarrow x - y \, \in \, \text{int } K, 
$$
и вместо $y \succ_K x$ можно написать: $x \prec_K y$.

Отношение $\succeq_K$ является отношением частичного порядка со следующими свойствами:
\begin{itemize}
    \item {\it рефлексивность}: $\forall \, x \quad x \succeq_K x$;
    \item {\it антисимметричность}: $\forall \, x \succeq_K y, \, 
    y \succeq_K x \, \Rightarrow \, x = y$;
    \item {\it транзитивность}: $\forall \, x \succeq_K y, \, 
    y \succeq_K z \, \Rightarrow \, x \succeq_K z$;
    \item {\it совместимость с линейными операциями}:
   
    \begin{enumerate}[label=\asbuk*),ref=\asbuk*]
       \item {\it однородность}: 
       если $x \succeq_K y$ и $\lambda$~--- любое неотрицательное вещественное число, то
       $\lambda x \succeq_K \lambda y$;
        \item {\it аддитивность}:
       если $x \succeq_K y$ и $w \succeq_K z$, то $x + w \succeq_K y + z$. 
    \end{enumerate}
\end{itemize}

\begin{example}
При $K = \mathbb{R}_+$ частичный порядок $\succeq_K$ превращается
в обычное отношение порядка $\ge$ на $\mathbb{R}$ (аналогично
и со строгими неравенствами).
\end{example}

\begin{example}
Неотрицательный ортант $K = \mathbb{R}^n_+$ является правильным
конусом, а соответствующее ему обобщённое неравенство
$x \succeq_K y$ обозначает покомпонентные неравенства для
векторов~$x$ и $y$:
$$
x_i \, \ge \, y_i, \quad i = 1, \, \ldots, \, n,
$$
то же самое можно записать и так: $x \ge y$.
\end{example}

\begin{example}\label{ch1_ex_loren}
Часто встречается в подобных неравенствах
{\it конус Лоренца}\protect\footnotemark: 
    $$
    L^n = \left \{ x = (x_1, \, \ldots, \, x_{n-1}, \, x_{n})^T \, \in \, \mathbb{R}^{n} \, \left| \, x_{n} \, \ge
    \, \sqrt{x_1^2 + x_2^2 + \ldots + x_{n-1}^2} \right. \, \right \}, 
    $$
    где $n \geq 2$.
\end{example}
\footnotetext{В англоязычной литературе его иногда называют ''ice-cream cone''~--- <<конус-мороженое>>.}

\begin{example}\label{ch1_ex_matrix}
В пространстве $\mathbb{S}^n$ вещественных симметричных матриц размера $n \times n$ определён 
{\it матричный конус}: 
    $$
    \mathbb{S}_+^n = \left \{ A \, \in \, \mathbb{S}^n \, | \, \lambda_{\min}(A) \geq 0 \, \right \}.
    $$
Это конус симметричных матриц, все собственные значения которых неотрицательны.
\end{example}

\begin{defin} Пусть $K \subset V$ --- произвольный конус в вещественном векторном пространстве $V$. \emph{Двойственным} (или \emph{сопряжённым}) к $K$ конусом называется множество
\begin{equation} \label{ch1_eq_dual_cone} 
K^* = \{ p \in V^* \,|\, \langle x,p \rangle \geq 0\ \forall\ x \in K \}.
\end{equation}
\end{defin}

{
\renewcommand{\baselinestretch}{0.96}
\selectfont

Отметим, что двойственный к $K$ конус определён в двойственном к $V$ векторном пространстве. Нетрудно видеть, что имеет место включение $K \subset (K^*)^*$. 

\begin{teo}[о сопряжённом конусе \cite{Ben-Tal}] Пусть~$K \subset  \mathbb{R}^n$ является непустым множеством. Тогда множество~$K^*$, определённое в~\eqref{ch1_eq_dual_cone}, является замкнутым выпуклым конусом.
\begin{itemize}
    \item Если $\mbox{int} \, K \neq \varnothing$, то
    $K^*$ является остроконечным конусом.
    \item Если $K$~--- замкнутый выпуклый остроконечный конус, то $\textup{\mbox{int}} \, K^* \neq \varnothing$.
    \item Если $K$~--- замкнутый конус, то сопряжённый к конусу~$K^*$ совпадает
    с исходным:
    $$
    (K^*)^* = K.
    $$
\end{itemize}
\end{teo}

Таким образом, если конус $K$~--- регулярный, то двойственный конус $K^*$ также будет регулярным. Более того, в этом случае мы имеем соотношение $(K^*)^* = K$. 

Для приведённых выше семейств конусов двойственный конус совпадает с прямым, если отождествить пространства $V$ и $V^*$ путём отождествления канонического базиса в $V$ с его двойственным базисом в $V^*$. Такие конусы называются \emph{самодвойственными}.

\section{Как получать выпуклые функции?}\label{ch1_sect_preserve_conv_func}
Способность распознавать тот факт, что функция является выпуклой, либо
способность сконструировать выпуклую функцию из
функций, про которые заведомо известно, что они выпуклы, могут стать
определяющими факторами успешного решения задачи.
Поэтому в данном разделе рассмотрены некоторые операции
над выпуклыми функциями, которые сохраняют выпуклость (или вогнутость) 
полученных функций.

\subsection{Неотрицательные взвешенные суммы}
Очевидно, что сумма двух выпуклых функций выпукла. Также очевидно, что если умножить выпуклую функцию $f$ на некоторое неотрицательное число $\alpha$, то полученная функция $\alpha f$ также выпукла.
Комбинируя эти свойства, получим следующую лемму.
\begin{lemma}
Пусть функции $f_1 \, : \, \mathbb{R}^n \, \rightarrow \, \mathbb{R}, \, \ldots, \, f_m(x) \, : \, \mathbb{R}^n \, \rightarrow \, \mathbb{R}$ выпуклы, $w_1 \, \in \, \mathbb{R}_+, \, \ldots, \, w_m \, \in \, \mathbb{R}_+$. Тогда функция
$$
f = w_1 f_1 + \ldots + w_m f_m
$$
является выпуклой.
\end{lemma}

}

Аналогично, неотрицательная взвешенная сумма вогнутых функций является вогнутой функцией.

\subsection{Аффинные преобразования аргумента}
Сразу сформулируем данное свойство выпуклых функций в виде леммы.
\begin{lemma}
Пусть $f \, : \, \mathbb{R}^n \, \rightarrow \, \mathbb{R}$ --- выпуклая функция, $A \, \in \, \mathbb{R}^{n \, \times \, m}$, $b \, \in \, \mathbb{R}^n$. Тогда функция $g(x) = f(Ax + b)$ с областью определения \emph{$\mbox{dom }$} $g = \{x \, | \, Ax +\linebreak +\, b \, \in \,$ \emph{$\mbox{dom }$} $f \}$ является выпуклой. Аналогично, если $f$ --- вогнутая функция, то $g$ --- тоже вогнута. 
\end{lemma}

\subsection[Композиция, максимум и точная верхняя грань\\ функций]{Композиция, максимум и точная верхняя грань\\ функций}\label{ch1_subs_max}
Данные очень важные свойства часто используются при доказательстве выпуклости функций.
\begin{lemma}
Пусть $f \, : \, \mathbb{R}^n \, \rightarrow \, \mathbb{R}$ --- выпуклая функция, $h \, : \, \mathbb{R} \, \rightarrow \, \mathbb{R}$ --- выпуклая неубывающая функция. Тогда композиция этих двух функций $g (x) = h \left ( f(x) \right )$ с областью определения \emph{$\mbox{dom }$} $ g = 
\{ x \, \in \,$ \emph{$\mbox{dom }$} $f \, | \, f(x) \, \in \,$\linebreak  $\in\, $ \emph{$\mbox{dom }$} $ h \}$ является
выпуклой функцией.
\end{lemma}

Если $f_1$ и $f_2$ --- выпуклые функции, то функция $f(x) = \linebreak = \max \, 
\{ f_1(x), \, f_2(x) \}$ с областью определения $\mbox{dom} \, f = 
\mbox{dom} \, f_1 \, \cap \, \mbox{dom} \, f_2$ тоже выпукла. Данное свойство легко доказать, воспользовавшись неравенством Йенсена.
Понятно, что аналогичное свойство 
можно доказать для функции $f(x) = \max \{ f_1(x), \, \ldots, \, f_m(x) \}$,
которая является выпуклой, если функции $f_1$, $f_2$,~$\ldots$, $f_m$ выпуклы. 

Воспользовавшись выпуклостью функции максимума, можно доказать,
например, выпуклость функции суммы $k$ максимальных компонент (см.~предыдущий параграф настоящей книги, а также доказательство в~\cite{Boyd_book_2004}). 

Операция взятия точной верхней грани от множества, содержащего бесконечное количество выпуклых функций, также сохраняет
свойство выпуклости.
\begin{lemma}\label{ch1_lem_sup}
Пусть $A$ --- некоторое множество и $f(x, \, y)$ --- выпуклая по $x$ функция для любого $y \, \in \, A$,  тогда функция 
$$
g(x) = \sup_{y \, \in \, A} f(x, \, y)
$$
является выпуклой по $x$ с областью определения
\[ \emph{\mbox{dom }} 
 g = 
 \left\{ x \, \left| \, (x, \, y) \, \in \,
 \emph{\mbox{dom }}
 f \quad \forall 
\, y \, \in \, A, \, \sup\limits_{y \, \in \, A} f(x, \, y) \, < \, \infty  \right. \right\}.
\]
Аналогично, поточечная нижняя грань множества вогнутых функций является вогнутой функцией.
\end{lemma}
Лемму \ref{ch1_lem_sup} легко доказать, используя тот факт, что надграфик поточечной операции нахождения верхней грани выпуклых функций есть пересечение надграфиков этих функций, а пересечение выпуклых множеств выпукло.

\subsection{Минимизация}
Операции взятия поточечного максимума и нахождения точной верхней грани сохраняют выпуклость (см. параграф~\ref{ch1_subs_max}). Оказывается, в определённых случаях операция нахождения точной нижней грани также сохраняет выпуклость.
\begin{lemma}
\label{ch1_lem_conv_min}
Пусть $G\left( {x, \, y} \right)$ --- выпуклая по $(x, \, y) $ функция, $Q$ --- выпуклое непустое множество.
Если задача $\inf\limits_{y \, \in \, Q} G(x, \, y)$ разрешима для всех $x$, 
то функция\protect\footnotemark 
$$
f(x) = \inf_{y \, \in \, Q} G(x, \, y)
$$
выпукла по $x$.
При этом область определения функции $f$ есть проекция
$\textup{dom } G$ на её $x$-координаты:
$$
\textup{dom } f = \left \{ x \, | \, \exists\ y \in Q:\ (x, \, y) \, \in \, \textup{dom } G \right \}.
$$
\end{lemma}
\footnotetext{О формализованной записи экстремальных задач см. п.~\ref{ch1_sect_extrem}.}
Доказать 
эту лемму можно, воспользовавшись неравенством Йенсена для 
$x_1$, $x_2 \, \in \, \mbox{dom } f$.

\begin{example}[{\emph{расстояние до множества}}]
Расстояние от точки $x$ до множества $Q \, \subseteq \, \mathbb{R}^n$ определяется как
$$
\inf_{y \, \in \, Q} \|x - y\|.
$$
Функция $\| x - y \|$ выпукла по $(x, \, y)$ как любая функция вычисления нормы вектора, следовательно, если множество $Q$ выпукло, то функция $\inf\limits_{y \, \in \, Q} \|x - y\|$ выпукла по $x$. 
\end{example}

\begin{exercise}
Если в лемме~\ref{ch1_lem_conv_min} предположить $\mu$-сильную выпуклость $G(x,y)$ в 2-норме по совокупности переменных $(x,y)$, то верно ли, что функция $f$ будет $\mu$-сильно выпуклая в 2-норме?\protect\footnotemark
\end{exercise}
\footnotetext{См. утверждение 1 из \cite{GasDvu_TranspOptima1516} в качестве отправной точки. Отметим также (см. утверждение~3~\cite{GasDvu_TranspOptima1516}), что если $G(x,y)$ имеет константу Липшица градиента в 2-норме $L$, то константа Липшица градиента функции~$f$ в 2-норме будет не больше $L$.}

\section{Негладкие задачи}
\label{ch2_subs_nonsmooth}
При анализе задач оптимизации нельзя ограничиваться только случаем дифференцируемых функций.
Причина в том, что существует множество интересных задач оптимизации,
для которых функции, которые нужно минимизировать,
оказываются недифференцируемыми или даже разрывными. 
Понятно, что в общем случае невозможно найти минимум разрывной
функции. Но если функция состоит из нескольких гладких кусочков, с разрывами между этими кусками, можно попытаться что-то сделать, минимизируя каждый гладкий
кусок отдельно. 

В общем виде задача минимизации негладких непрерывных функций чрезвычайно сложна. Невозможно гарантировать,
что по информации, полученной о 
поведении такой функции в одной точке,
можно будет сделать правильный вывод о поведении функции
в соседних точках. Поэтому мы ограничимся важным частным случаем недифференцируемых функций~--- выпуклыми функциями. Для выпуклых функций можно ввести понятие, во многих отношениях аналогичное градиенту.
В частности, если функция непрерывна повсюду,
но не дифференцируема в определённых точках\footnote{Стандартный пример --- не дифференцируемая только в точке $x = 0$ функция $f(x) = | x |$, $x \, \in \, \mathbb{R}$.},
то можно её минимизировать, используя значения {\it субградиента}
или {\it обобщённого градиента} в таких точках.
Далее приводятся краткие сведения из субдифференциального исчисления~\cite{MagTih_convan, MordukhNumBook, Nesterov_Eff_meth89}, которое для выпуклого анализа представляет собой некоторый аналог дифференциального исчисления, где роль производной
играет {\it субдифференциал}. 

\begin{defin}\label{ch2_def_subdif}
Пусть $X \, \subset \, \mathbb{R}^n$ и $f: X \rightarrow \mathbb{R}$,
точка $x_0 \, \in \, X$.
Множество (возможно, пустое):
$$
\partial f(x_0) = \left\{g \, \in \, \mathbb{R}^n \, | \, f(x) - f(x_0) \ge g^T (x - x_0), \,
\forall \, x \, \in \, X \right\}
$$
называется субдифференциалом функции~$f$ в точке~$x_0$.
\end{defin}

Элементы $\partial f(x_0)$ называются {\it субградиентами функции~$f$ в точке~$x_0$}. Субградиент
определяется, вообще говоря, неоднозначно.
Из определения~\ref{ch2_def_subdif} ясно, что вектор $g \, \in \, \mathbb{R}^n$~--- субградиент~$f$ в точке~$x_0$
тогда и только тогда, когда 
линейная функция $l(x) = f(x_0) + g^T (x - x_0)$
всюду не превосходит $f$. Также имеется ещё одна
геометрическая интерпретация субградиента: 
 вектор $g \,  \in  \, \mathbb{R}^n$ является субградиентом в точке $x_0$ тогда и только тогда, когда $(g, \, -1)$ определяет опорную гиперплоскость\footnote{См. определение опорной гиперплоскости на с.~\pageref{ch1_def_sup_hyp}.}
 к надграфику~$f(x)$ в~точке $(x_0, \, f(x_0))$.

Если функция $f$ выпукла и дифференцируема, то её градиент
    в точке $x_0$ совпадает с субградиентом и он соответственно в этом случае единственен (доказательство см., например, в~\cite{MagTih_convan}).

Можно показать, что всякая выпуклая на
открытом множестве функция дифференцируема
почти всюду (т.е.
за исключением множества меры нуль).

\begin{example}
Рассмотрим функцию $f(x) = | x |$, $x \, \in \, \mathbb{R}$. Функция является выпуклой.
По определению субдифференциала имеем
$$
g \, \in \, \partial f(0) \Leftrightarrow
|x| - |0| \ge g \cdot (x - 0) \quad \forall x \, \in \, \mathbb{R},
$$
тогда $|x| \ge gx$, следовательно, субградиентом функции $f(x) = |x|$ в~точке~$0$ может быть любое 
число из отрезка $[-1, \, 1]$, и
$\partial f(0) = [-1, \, 1]$.
\end{example}

\begin{exercise}
Подумайте, чему будет равен субдифференциал\protect\footnotemark 
$\,$ функции
$\inf_{y \, \in \, Q} \|x - y\|$
(расстояние от точки~$x$ до множества~$Q \, \subseteq \, \mathbb{R}^n$). 
\end{exercise}
\footnotetext{Ответ на этот вопрос можно найти в \cite[с.~242]{DemVas_book_81}.}

Следующая теорема является чрезвычайно важным для выпуклых
экстремальных задач фактом.
\begin{teo}[аналог теоремы Ферма]\label{ch1_teo_ferma} Необходимым и достаточным условием глобального минимума 
выпуклой функции в точке~$x^*$ является
$$
0 \, \in \, \partial f(x^*).
$$
\end{teo}

\begin{proof}
Решается задача $\min\limits_{x \, \in \, \mathbb{R}^n} f(x)$.
Если $x^*$~--- точка минимума $f(x)$, то
для любого $x \, \in \, \mathbb{R}^n$
верно $f(x) \geq f(x^*)$, т.е.
$f(x) - f(x^*) \geq 0 = 0 \cdot (x -\linebreak - x^*)$,
и, по определению субдифференциала, $0 \, \in \, \partial f(x^*)$.
Для доказательства достаточности нужно повторить все рассуждения в обратном порядке.
\end{proof}

Приведём без доказательства набор
стандартных свойств субдифференциала,
знание которых понадобится для полного понимания
дальнейшего изложения.

\begin{lemma}
Множество субградиентов в любой точке непусто, выпукло,
замкнуто и ограничено.
\end{lemma}

Знание субградиента помогает вычислять производную по направлению.
\begin{lemma}
Производная выпуклой функции~$f$ в точке~$x_0$ 
по любому направлению $s$ всегда
существует и определяется равенством
$$
\frac{\partial f(x_0)}{\partial s} =
\lim_{t\to 0}\frac{f(x_0 +ts) - f(x_0)}{t}=
\max_{a \, \in \, \partial f(x_0)} a^T s.  
$$
\end{lemma}

Следующие две теоремы и лемма дают возможность вычислять субградиенты сложных функций так же просто, как и градиенты дифференцируемых функций.
\begin{teo}[Моро--Рокафеллара]
Если функции $f_1(x)$, $f_2(x)$
выпуклы, то для всех 
$x \, \in \, \textup{dom } f_1 \cap
\textup{dom } f_2$ справедлива
формула
$$
\partial (f_1 + f_2) (x) = \partial f_1(x) + \partial f_2(x).
$$
\end{teo}

\begin{teo}[Дубовицкого--Милютина]
\label{ch2_teo_dub_mil}
Пусть функции
$f_1$, $f_2$ выпуклые и $f_1(x_0) = f_2(x_0)$ {\rm(}причём в точке~$x_0$
обе функции непрерывны{\rm)}. Тогда
$$
\partial \max (f_1, \, f_2) (x_0) = \textup{conv } \left ( \partial f_1(x_0) \cup
\partial f_2(x_0) \right ).
$$
\end{teo}

Очевидно, что две последние теоремы распространяются
по индукции на любое конечное число функций.

\begin{lemma}
Пусть $A$~--- матрица размера~$m \times n$, $\varphi(y)$~---
выпуклая функция на~$\mathbb{R}^m$, $f(x) = \varphi (Ax)$,
$x \,\in \, \mathbb{R}^n$. Тогда
$$
\partial f(x) = A^T \partial \varphi (A x).
$$
Обозначение $\partial \varphi (A x)$ означает субдифференциал
$\varphi$ по отношению к аргументу в точке~$Ax$.
\end{lemma}

%
%
\section{Формула Демьянова--Данскина}\label{ch1_sect_dem-dan}
В приложениях, приводящих к оптимизационным задачам на минимакс,
часто возникает необходимость дифференцирования
функций под знаком максимума. Для
этих целей применяют теорему о дифференцировании (по направлению) под знаком максимума,
которую часто называют теоремой или формулой Демьянова--Данскина\footnote{В англоязычной
литературе эту теорему часто называют просто теорема Данскина (Danskin's theorem). В России также можно встретить название этой формулы (теоремы), как формулы Демьянова--Данскина--Рубинова.}, поскольку оба автора независимо пришли
к аналогичному утверждению в конце 60-х годов XX века~\cite{Dem_book_74, Danskin}. 
Сформулируем одну из форм этого утверждения.

\begin{teo}
Пусть $G(x, \, y) \, : \, \mathbb{R}^n \, \times \, Y \, \rightarrow \, \mathbb{R}$ --- дифференцируемая по $x$ для
всех $y \, \in \, Y$ функция, а $Y \, \subset \, \mathbb{R}^m$ -- компактное множество. Пусть также
$G(x, \, y)$ выпукла по $x$ для любого $y \, \in \, Y$.
Тогда функция
$$
f(x) = \max_{y\in Y} \, G(x, \, y) 
$$
\nk{выпукла, и в случае,} если существует единственный
$y\left( x \right)$, \nk{такой что} 
$G\left( {x, \, y\left( x \right)} \right) = \max\limits_y \,G\left( {x, \, y} \right)$,
 \nk{градиент} функции $f(x)$ можно вычислить
следующим образом:
$$
\nabla f\left( x \right)=\nabla_x G\left( {x,y\left( x \right)} 
\right)=\left. {{\frac{\partial G\left( {x,y} \right)}{\partial 
x }}  } \right|_{y=y\left( x \right)}.
$$
Если $\partial G / \partial x$ непрерывна по $y$ для всех $x$, то
субдифференциал\protect\footnotemark
$$
\partial f\left( x \right) = \conv\bigcup\limits_{\tilde 
{x}:\;y\left( {\tilde {x}} \right)=y\left( x \right)} {\nabla _x G\left( 
{\tilde {x},y\left( {\tilde {x}} \right)} \right)},
$$
где $\conv$ обозначает операцию взятия выпуклой
оболочки от множества.
\end{teo}
\footnotetext{Определение и свойства субдифференциала см. в п.~\ref{ch2_subs_nonsmooth}.}

Не вдаваясь в детали, поясним основную идею доказательства этой теоремы в случае, когда $y(x)$ определяется (единственным образом) из системы ($m$ уравнений):
$$G'_y(x, y) = \nabla_y G(x, y)\equiv 0.$$
Эта идея делает интуитивно понятным
основное утверждение теоремы~Демьянова--Данскина.
Поскольку 
$$
\nabla f(x)=\nabla_x G(x, y(x))+y'(x)^T\underbrace{\nabla_y G(x, y(x))}_{=0},
$$
то $\nabla f(x) = G'_x(x, y(x)) = \nabla_x G(x, y(x)) $. Здесь  $y'(x)$ --- матрица Якоби ($(i,j)$ элемент матрицы $y'(x)$ равен $\partial y_i(x)/\partial x_j$).   Отсюда также становится понятно, что аналогичное утверждение (естественное, при некоторых предположениях) можно сформулировать для
$f(x) = \min\limits_y G(x, \, y)$ и даже для $f(x) = \min\limits_y\max\limits_z G(x,y,z)$, см. раздел~\ref{subsec:inverse}.

В 1971 году Дмитрий Бертсекас в своей диссертации доказал
более общий результат, не требующий дифференцируемости $G(\cdot, \, y)$, см. \cite{Bertsekas_phd}.

Утверждения, аналогичные теореме~Демьянова--Данскина, были
известны и ещё раньше. Так, в середине 60-х годов XX века А.\,Я.~Дубовицкий и~А.\,А.~Милютин при
получении принципа максимума для задач оптимального управления с фазовыми ограничениями доказали и фактически
уже использовали теорему о дифференцировании по направлению под знаком $\sup$ \cite{DubMil},
см. также теорему~\ref{ch2_teo_dub_mil} данного пособия.

В более общем виде теорему Демьянова--Данскина
можно сформулировать следующим образом 
(см. также \cite{Danskin_Bern}).
\begin{teo}[\cite{Danskin}]
Пусть $Y$~--- компактное топологическое пространство, $G: \mathbb{R}^n \times Y \rightarrow \mathbb{R}$~--- непрерывное отображение, непрерывно дифференцируемое по первой переменной.
Пусть
$$
f(x) = \max_{y\in Y} G(x, \, y) 
$$
и
$$
\hat Y(x) = \left \{ y \, \in \, Y \, | \, G(x, \, y) = f(x) \right \}.
$$
Тогда функция $f(x)$ для каждого $x$ и $h$ из $\mathbb{R}^n$ имеет
производную в точке~$x$ по направлению $h$,
вычисляемую по формуле
$$
\frac{\partial f(x)}{\partial h} 
=\lim_{t\to 0} \frac{f(x+th) - f(x)}{t}
= \max_{y \, \in \, \hat Y} \, \sum_{i=1}^n 
h_i  \frac{\partial G\left( {x,y} \right)}{\partial 
x_i }. 
$$
\end{teo}
Существует версия теоремы Демьянова--Данскина и без условия компактности~$Y$
(см.~\cite{Danskin_Bern}).

Существует несколько так называемых теорем
об очистке, доказанных В.\,Л.~Левиным
(в том числе об очистке субдифференциалов).
Приведём основную из этих теорем.
\begin{teo}[В.\,Л. Левина об очистке, \cite{MagTih_convan}]
\label{ch1_teo_Levin_clear}
Пусть $T$~--- компакт, $X$~--- конечномерное пространство, 
$\emph{card}(T)$~--- мощность~$T$, $\emph{dim}(X)$~--- размерность пространства~$X$,
выполняется условие $\emph{\mbox{card}}(T) \geq \emph{\mbox{dim}}(X) + 1$
и функция $G: T \times X \rightarrow \mathbb{R}$ такова, что
\begin{enumerate}
\renewcommand\labelenumi{\textup{\arabic{enumi})}}
    \item отображение $G(t, \, \cdot): \, X \rightarrow \mathbb{R}$ выпукло для каждого $t \, \in \, T$;
\item отображение $G(\cdot, \, x): \, T \rightarrow \mathbb{R}$ полунепрерывно сверху для каждого $x \, \in \, X$;
\item $M = \inf_{x \, \in \, X} \max_{t \, \in \, T} \,G(t, \, x) \, > \, - \infty$.
\end{enumerate}
Тогда найдутся натуральное число $r$, такое, что $r \leq \emph{dim}(X) + 1$, и точки $\tau_1, \, \ldots, \, \tau_r \, \in \, T$, такие, что  $M = \inf_{x \, \in \, X} \max_{1 \leq i \leq r} \,G(\tau_i, \, x)$.
\end{teo}

Теорема~\ref{ch1_teo_Levin_clear} означает, что исходную задачу минимизации максимума отображения
можно заменить на задачу с меньшим
числом ограничений. Этим объясняется
такое название теоремы~--- <<очистка>>,
отбрасывание некоторых ограничений.

Связь между теоремой Демьянова--Данскина и теоремами об очистке более очевидна, если
привести субдифференциальную форму теоремы об 
очистке В.\,Л.~Левина,
которая, по сути, является обобщением теоремы
Дубовицкого--Милютина (теорема \ref{ch2_teo_dub_mil}) на произвольное число функций.

\begin{teo}[субдифференциальная форма теоремы об очистке, \cite{MagTih_convan}]
\label{ch1_teo_Levin_clear_subdif}
Пусть $T$~--- компактное топологическое пространство, $X$~--- конечномерное 
локально выпуклое топологическое векторное пространство {\rm(}т.е. оно хаусдорфово и существует база окрестностей нуля, состоящая из выпуклых 
множеств{\rm)} и функция $G: T \times X \rightarrow \mathbb{R}$ такова, что отображение $G(t, \, \cdot): X \rightarrow \mathbb{R}$ выпукло для каждого $t \, \in \, T$, а отображение $G(\cdot, \, x): T \rightarrow \mathbb{R}$ полунепрерывно сверху для каждого $x \, \in \, X$. Тогда если $y \, \in \, \partial \max_{t \, \in \, T} G(t, \, \hat x)$ для некоторого $\hat x \, \in \, X$, то найдутся натуральное число $r$, такое, что $r \leq \emph{dim}(X) + 1$, и точки $t_i \, \in \, T$, $i = 1, \, \dots, \, r$, такие, что
$$
y \, \in \, \conv \,\bigcup_{i=1}^r\, \partial_x G(t_i, \, \hat x), \quad \max_{1 \leq i \leq r} G(t_i, \, \hat x) = G(t_1, \, \hat x) =
\dots = G(t_r, \, \hat x),
$$
где $\partial_x G(t_i, \, \hat x)$ обозначает
субдифференциал функции $G(\cdot, \, x)$ в точке~$\hat x$.
\end{teo}
Для доказательства данной теоремы достаточно применить теорему Ферма (теорема~\ref{ch1_teo_ferma}), обычную теорему об очистке (теорема~\ref{ch1_teo_Levin_clear})
и теорему Дубовицкого--Милютина (теорема~\ref{ch2_teo_dub_mil}) о дифференцировании максимума из двух выпуклых функций (см.~\cite{MagTih_convan}).

С теоремами об очистке связано одно
из доказательств теоремы Хелли,
см. п.~\ref{ch1_subsect_Hilb}.

\section{Дополнение по Шуру} \label{ch1_section_schur}
Рассмотрим симметричную матрицу размера $n \times n$
следующего вида:
\begin{equation}
\label{ch1_schur_x}
X = \left[ {\begin{array}{l}
 A\quad B \\ 
 B^T\;\;C \\ 
 \end{array}} \right],
\end{equation}
где $A$ -- симметричная матрица размера $k \times k$.
Если $\mbox{det } A \neq 0$, то
матрица 
$$
S = C - B^T A^{-1} B
$$
называется дополнением по Шуру матрицы $A$ в $X$.
Матрицы такого вида возникают в различных важных теоремах и формулах, в частности, при решении систем линейных алгебраических уравнений методом исключения блоков переменных.

В этом разделе ограничимся рассмотрением применения дополнения по Шуру при решении задачи минимизации квадратичной формы:
\begin{equation}
\label{ch1_shur_quad}
G\left( {x, \, y} \right)= {x^T Cx}  +  {y^T Ay} 
 +2 { x^T B y}  ,
\end{equation}
где матрица $A$ положительно определена, т.е. ${\forall  x\neq 0}$ выполняется 
 ${{x^T Ax}   >  0}$.

Решением задачи минимизации $\min\limits_y G(x, \, y)$ является
$y = - A^{-1} B x$.
Тогда 
\begin{equation}
\label{ch1_fCBA}
f\left( x \right)= \min_y G\left( {x, \, y} \right)=G\left( 
{x, \, y\left( x \right)} \right)= {x^T  \left( {C-B^T A^{-1} B} 
\right)x} .
\end{equation}

Для того чтобы~(\ref{ch1_shur_quad}) являлась выпуклой функцией, необходимо и достаточно, чтобы
матрица~(\ref{ch1_schur_x}) была неотрицательно определённой, так как квадратичную форму~(\ref{ch1_shur_quad})
можно представить как
$$
G(x, \, y) = \left[ {\begin{array}{l}
 y \\ 
 x  
 \end{array}} \right]^T 
 \left[ {\begin{array}{l}
 A\quad B \\ 
 B^T\;\;C \\ 
 \end{array}} \right]
 \left[ {\begin{array}{l}
 y \\ 
 x  
 \end{array}} \right] .
$$

Но функция (\ref{ch1_fCBA}) является выпуклой по лемме \ref{ch1_lem_conv_min}. Так как (\ref{ch1_fCBA}) тоже является квадратичной формой, для её выпуклости необходимо и достаточно, чтобы матрица $S = {C-B^T A^{-1} B}$ была неотрицательно определённой.
Следовательно, верны следующие леммы.

\begin{lemma}
Если $A \succ 0$, то $X \succeq 0$ тогда и только тогда, когда $S \succeq 0$.
\end{lemma}

\begin{lemma}
$X \succ 0$ тогда и только тогда, когда $S \succ 0$ и $A \succ 0$.
\end{lemma}

\section{Теорема о неявной функции}\label{ch1_implicit_sect}
Теорема о неявной функции, известная из математического анализа, часто
применяется при исследовании итеративных оптимизационных методов,
не разрешённых явно относительно 
текущего приближения к решению~$x_{k + 1}$
на $k$-й итерации. 
В таких методах новое приближение $x_{k+1}$ может быть решением некоторой вспомогательной задачи минимизации, как, например, в методах регуляризации \ag{(подробнее см.~\cite{Polyak})}.
Также 
теорема о неявной функции применяется при доказательстве принципа множителей Лагранжа.

Пусть определено отображение $g(x, \, y) : \mathbb{R}^m \times \mathbb{R}^n \, \rightarrow \, \mathbb{R}^n$.
Будем обозначать $g'_x(x, \, y)$, $g'_y(x, \, y)$ матрицы Якоби ($(i, \, j)$ элемент матрицы $g'_x(x, \, y)$ равен $\partial g_i(x, \, y)/\partial x_j$), которые для краткости будем называть (используя одномерную ($n = 1$) аналогию) <<производные>> $g$ по соответствующим переменным.

\begin{teo}[о неявной функции]
Пусть $g(x^*, \, y^*) = 0$, $g(x, \, y)$
непрерывна по $\{x, \, y\}$ в окрестности $x^*$, $y^*$ и матрица $g'_y \left (x^*, \, y^* \right )$ невырождена.
Тогда существует единственная непрерывная в окрестности $x^*$ функция $y(x)$ такая, что $y^* = y (x^*)$, $g(x, \, y(x)) = 0$.
Если дополнительно существует $g'_x(x^*, \, y^*)$, то $y(x)$ дифференцируема в $x^*$ и
$$
y'(x^*) = - \left [ g'_y(x^*, \, y^*) \right ]^{-1} g'_x(x^*, \, y^*).
$$
\end{teo}

Данная теорема без больших изменений переносится с числовых функций на отображения произвольных банаховых пространств (в такой общности данная теорема используется в разделе~\ref{subsec:inverse}).
\begin{teo}[о неявной функции, бесконечномерный случай, \cite{Kolmog_FunAn}]
\label{ch1_teo_impl_f1}
Пусть $X, \, Y, \, Z$~--- банаховы пространства, $U$~--- окрестность
точки $(x^*, \, y^*) \, \in \, X \times Y$ и $g$~--- отображение $U$
в $Z$, обладающее следующими свойствами:
\begin{enumerate}
\renewcommand\labelenumi{\textup{\arabic{enumi})}}
\item $g$ непрерывно в точке $(x^*, \, y^*)$;
\item $g(x^*, \, y^*) = 0$;
\item Частная производная $g'_y(x, \, y)$ существует в $U$ и непрерывна в точке $(x^*, \, y^*)$, а оператор $g_y'(x^*, \, y^*)$
имеет ограниченный обратный.
\end{enumerate}
Тогда существуют такие $\varepsilon > 0$, $\delta > 0$ и такое
отображение
\begin{equation}
\label{ch1_eq_teo_impl_f1}
y = f(x),
\end{equation}
определённое при $\| x - x^* \| < \delta$ и непрерывное в точке~$x^*$, что каждая пара $(x, \, y)$, для которой $\| x - x^* \| < \delta$ и $y = f(x)$, удовлетворяет уравнению
\begin{equation}
\label{ch1_eq_teo_impl_f2}
g(x, \, y) = 0,
\end{equation}
и обратно, каждая пара $(x, \, y)$, удовлетворяющая уравнению~\eqref{ch1_eq_teo_impl_f2} и~условиям $\| x - x^* \| < \delta$, $\| y - y^* \| \leq \varepsilon$, удовлетворяет
и \eqref{ch1_eq_teo_impl_f1}.
\end{teo}
Доказательство теоремы основано на принципе сжимающих отображений, его
можно найти в \cite[с. 508]{Kolmog_FunAn}. Там же приведено
доказательство и нижеследующей теоремы об условиях, при которых функция, определяемая уравнением вида $g(x, \, y) = 0$,
дифференцируема.
\begin{teo}
Пусть выполнены условия теоремы~{\rm\ref{ch1_teo_impl_f1}}, и пусть, кроме того, в $U$ существует частная производная $g'_x$, непрерывная
в точке $(x^*, y^*)$. Тогда отображение $f$ дифференцируемо в точке $x^*$ и
$$
f'(x^*) = - [g_y'(x^*, \, y^*)]^{-1} g_x'(x^*, \, y^*).
$$
\end{teo}

Одно из важных применений теоремы о неявной функции~---
вопрос об эквивалентности двух определений касательных
многообразий в банаховых пространствах (теорема Люстерника \cite{Kolmog_FunAn, Ioffe}). В частном случае, на плоскости, эти определения
выглядят так. Пусть $x = (x_1, \, x_2) \, \in \, \mathbb{R}^2$ и $f(x)$~--- дифференцируемая функция.
Уравнение $f(x) = 0$ определяет на плоскости кривую~$C$.
Если $x_0 \, \in \, C$, то касательную к $C$ в точке~$x_0$ можно определить либо как совокупность векторов вида $x_0 + th$, где вектор~$h$ перпендикулярен
к градиенту $\nabla f(x_0)$, либо как совокупность точек
$x_0 + th$, расстояние от которых до кривой~$C$ есть бесконечно малая выше первого порядка относительно~$t$.

Теорема Люстерника, в свою очередь, играет важную роль в задачах оптимального управления. С помощью теоремы Люстерника правило множителей Лагранжа нахождения условного экстремума (см. п.~\ref{ch1_sect_extrem}) может быть доказано для экстремальных задач в банаховых пространствах.

\section{Экстремальные задачи и двойственность}
\label{ch1_sect_extrem}
Важность экстремальных задач, часто возникающих
в различных приложениях, объясняет необходимость
существования общей теории максимумов
и минимумов функций на заданных множествах.
Эта теория значительно упрощается, если
возможно использование свойств, связанных с выпуклостью. Так были получены следующие
результаты, относящиеся к теории экстремальных задач:
разнообразные теоремы двойственности и теоремы
о характеристике точек, в которых достигается
экстремум (необходимые условия экстремума,
достаточные условия экстремума),
теория существования решений экстремальных задач.

Точно поставленная ({\it формализованная}) экстремальная задача обязательно
включает в себя следующие элементы~\cite{Magaril_lect, MagTih_convan, Pshen_80}:
\begin{itemize}
    \item {\it целевую функцию}~--- функционал, который нужно максимизировать
    или минимизировать~--- 
${f  :  X  \rightarrow  \bar{\mathbb{R}}}$
с областью определения ${\text{dom } f = X}$, где 
$X$~--- векторное пространство (например, $\mathbb{R}^n$),
$\bar{\mathbb{R}}$~--- вещественная прямая~$\mathbb{R}$, пополненная символами $\pm\infty$,
т.е. $\bar{\mathbb{R}} = \mathbb{R} \, \cup \,  \{-\infty\} \, \cup \, \{+\infty\}$;
\item {\it ограничения}, которые в общем виде можно
представить как подмножество $D \, \subset \, X$.
\end{itemize}
Формализованную экстремальную задачу на минимум
записывают в следующей форме:
\begin{equation}
\label{ch1_eq_min_first}
\min_{x \, \in \, D} f(x),
\end{equation}
соответственно формализованную экстремальную задачу на максимум записывают так\footnote{В российской традиции формулировку
задачи~\eqref{ch1_eq_min_first} принято записывать
и так:
$$
f(x) \, \rightarrow \, \min_{x \, \in \, D}.
$$
}:
\begin{equation}
\label{ch1_eq_max_first}
\max_{x \, \in \, D} f(x).
\end{equation}
Точки $x \, \in \, D$ называют {\it допустимыми} для
задач~\eqref{ch1_eq_min_first}, \eqref{ch1_eq_max_first}, само множество~$D$ называют {\it допустимым множеством}. Если $D = X$, то задачи~\eqref{ch1_eq_min_first}\,--\,\eqref{ch1_eq_max_first} называют {\it задачами без ограничений}. 

Допустимая точка~$x^*$ называется
{\it решением} задачи~\eqref{ch1_eq_min_first} 
или {\it глобальным {\rm(}абсолютным{\rm)} минимумом} в задаче~\eqref{ch1_eq_min_first}, если
для любого $x \, \in \, D$ имеет место неравенство $f(x) \, \ge \, f(x^*)$.
Допустимая точка~$x^*$ называется
{\it решением} задачи~\eqref{ch1_eq_max_first} 
или {\it глобальным {\rm(}абсолютным{\rm)} максимумом} в задаче~\eqref{ch1_eq_max_first}, если
для любого $x \, \in \, D$ имеет место неравенство $f(x) \, \le \, f(x^*)$.
Если в задаче требуется найти и точки максимума,
и точки минимума, то в записи~\eqref{ch1_eq_min_first}
(\eqref{ch1_eq_max_first})
вместо $\min \, (\max)$ пишется $\text{extr}$ (такая задача
будет называтьcя задачей на экстремум функционала~$f$).

Ясно, что если $x^*$~--- решение задачи~\eqref{ch1_eq_min_first} на минимум (\eqref{ch1_eq_max_first} на максимум),
то $x^*$ также является решением аналогичной задачи на максимум
(минимум) с функционалом $-f(x)$ вместо $f(x)$.

 Введём также понятие локального экстремума (т.е. локального минимума и максимума).
\begin{defin}
Если в $X$ определено понятие <<окрестности
точки>>, то точка $\tilde{x} \, \in \, D$ называется \emph{локальным минимумом} (\emph{максимумом}) в задаче~\eqref{ch1_eq_min_first} (\eqref{ch1_eq_max_first}), если существует
такая её окрестность~$U$, что $f(x) \,\ge \, f(\tilde{x})$ ($f(x) \, \le \, f(\tilde{x})$)
для всех допустимых $x \, \in \, D \cap U$.
\end{defin} 
Для задач без ограничений типа~\eqref{ch1_eq_min_first}, \eqref{ch1_eq_max_first}, где $D = X$, говорят,
что точка локального экстремума $\tilde{x}$ доставляет
также локальный экстремум функции~$f(x)$.

Далее для большей наглядности мы в основном будем рассматривать конечномерные задачи. Однако все основные результаты естественным образом переносятся и на бесконечномерные пространства. В частности, как уже отмечалось в конце п.~\ref{ch1_implicit_sect},
рассмотренный далее принцип множителей Лагранжа
верен не только в конечномерном евклидовом пространстве, его также
можно распространить и на бесконечномерные банаховы пространства
с помощью теоремы Люстерника, и даже на абстрактные векторные пространства (см.~\cite{MagTih_convan}). Несколько примеров приложения описанных далее результатов в бесконечномерных пространствах мы будем разбирать в разделах~\ref{neuman_pirson} и \ref{subsec:inverse}. 

\subsection[Правило множителей Лагранжа для задачи\\ с ограничениями типа равенств]{Правило множителей Лагранжа для задачи\\ с ограничениями типа равенств}
\epigraph{Можно высказать следующий общий принцип. Если ищется максимум или минимум некоторой функции многих переменных при условии, что между этими переменными имеется связь, задаваемая одним или несколькими уравнениями, нужно прибавить к функции, экстремум которой мы ищем, функции, задающие уравнения связи, умноженные на неопределённые множители, и искать затем максимум или минимум построенной суммы, как если бы переменные были независимы. Полученные уравнения, присоединённые к уравнениям связи, послужат для определения всех неизвестных.}{\textit{Ж.\,Л.~Лагранж}
\protect\footnotemark}

\footnotetext{Цит. по~\cite[с.~126]{MinMax}.}

В этом параграфе изложение во многом следует~\cite{AvakovMagaril, MagTih_convan} и
\cite{Polyak}. Здесь рассматривается
правило множителей Лагранжа для гладких задач
с ограничениями, задаваемыми равенствами. {\it Принцип Лагранжа} (метод решения задач,
основанный на идее Лагранжа)
описан самим автором в эпиграфе, поставленном
нами к данному параграфу.
В дальнейшем оказалось, что
идея Лагранжа поразительно плодотворна~--- необходимые условия экстремума в самых разных экстремальных задачах (с ограничениями типа
равенств, неравенств и включений)
имеют форму правила множителей Лагранжа.
По~словам Р.~Рокафеллара~\cite{Rocaf},
<<теория множителей Лагранжа позволяет преобразовывать экстремальные задачи с ограничениями в задачи, имеющие меньше 
ограничений, но больше переменных>>. 
До представленного ниже формализованного описания правила
множителей Лагранжа кратко поясним, что
принцип Лагранжа заключается в идее составления функции
(теперь она называется функцией Лагранжа), которая
является суммой минимизируемого (максимизируемого) функционала
и функций-ограничений, задающих равенства, с неопределёнными множителями. Тогда
необходимые условия экстремума для исходной
задачи с ограничениями эквивалентны
таким условиям для безусловной задачи
минимизации (максимизации) функции Лагранжа.
Кроме того, для выпуклых задач решение исходной задачи доставляет глобальный минимум функции Лагранжа, т.е. необходимые условия экстремума
минимума функции Лагранжа становятся и достаточными условиями.

Рассмотрим задачу с ограничениями типа равенств вида
\begin{equation}
\label{ch1_lagr_eqs_problem}
\min_{g_i(x) = 0} f(x), \, x \, \in \, \mathbb{R}^n, \, i = 1, \, \ldots, \, m,
\end{equation}
где $f$, $g_i,  \, i = 1, \, \ldots, \, m$, --- достаточно гладкие функции.

\begin{teo}[правило множителей Лагранжа] 
Если $\nabla g_i(x^*)$, $i\!\!\!=\linebreak = 1, \, \ldots, \, m$ линейно независимы и если $x^*$~--- точка минимума в задаче~{\rm(\ref{ch1_lagr_eqs_problem})}, т.е.
$x^*  \in  Q = \{x  |   g_i(x) = 0, \, i = 1, \, \ldots, \, m \}$ и $\forall  x  \in  Q \linebreak f(x^*) \, \le \, f(x)$, то найдутся такие $\lambda_1^*, \, \ldots \, \lambda_m^*$,
что
\begin{equation}
\label{ch1_grad_Lagr_func}
\nabla f(x^*) + \sum_{i = 1}^m \lambda_i^* \nabla g_i (x^*) = 0.
\end{equation}
\end{teo}

Числа $\lambda_1^*, \, \ldots \, \lambda_m^*$ в (\ref{ch1_grad_Lagr_func}) называются множителями Лагранжа.

Составим {\it функцию Лагранжа}: 
\begin{equation}
\label{ch1_Lagr_func}
L(x, \, \lambda_1, \, \ldots, \, \lambda_m) = f(x) + \sum_{i = 1}^m \lambda_i g_i (x),
\end{equation}
определённую на $\mathbb{R}^n \, \times \, \mathbb{R}^m$.
Если обозначить $g(x) = \left (g_1(x), \, \ldots, \, g_m(x) \right )$\linebreak и~$\lambda = (\lambda_1, \, \ldots, \, \lambda_m)$, то  
(\ref{ch1_Lagr_func}) можно переписать в виде
$$
L(x, \, \lambda) = f(x) + \langle \lambda, \, g(x) \rangle.
$$
Тогда правило множителей Лагранжа (\ref{ch1_grad_Lagr_func})
переформулируется как
\begin{equation}
\label{ch1_lx_llam}
L_x'(x^*, \, \lambda^*) = 0, \quad  L_\lambda'(x^*, \, \lambda^*) = 0,
\end{equation}
где $L_x'$, $L_\lambda'$ --- производные (градиенты) по соответствующим переменным.
В~(\ref{ch1_lx_llam}) переменные $x$ называются {\it прямыми}, а $\lambda$~--- {\it двойственными}.

\subsubsection[Идея доказательства с помощью теоремы о неявной\\ функции]{Идея доказательства с помощью теоремы о неявной функции}
Данное доказательство формулы (\ref{ch1_grad_Lagr_func}) основано на идее сведения задачи с ограничениями к задаче безусловной минимизации с помощью исключения переменных. В данном разделе градиент функции понимается как вектор-строка.

Разбиваем прямые переменные $x \, \in \, \mathbb{R}^n$
на две группы $v \, \in \, \mathbb{R}^{n - m}$ и $y \, \in \, \mathbb{R}^{m}$. Далее из равенств-ограничений $g(x) = 0$
выражаем одну группу через другую: $y = y(v)$. Рассмотрим задачу безусловной минимизации для $f(v, \, y(v))$.
Применение аналога теоремы Ферма (см. \nk{теорему}~\ref{ch1_teo_ferma}) для этой функции даёт  решение исходной задачи (\ref{ch1_lagr_eqs_problem}):
\begin{equation}
\label{ch1_lagr_impl_pr}
\nabla f(v^*, \, y(v^*)) = f'_v (v^*, \, y(v^*)) + 
f'_y(v^*, \, y(v^*))y'(v^*) = 0.
\end{equation}

Поскольку $g(v^*, \, y^*) = 0$, то в условиях выполнения теоремы о неявной функции (см.~п.~\ref{ch1_implicit_sect}):
$$
y'(v^*) = - \left [ g'_y(v^*, \, y^*) \right ]^{-1} g'_v(v^*, \, y^*).
$$
Тогда \eqref{ch1_lagr_impl_pr} можно переписать так:
\begin{equation}
\label{ch1_lagr_impl_pr2}
f'_v (v^*, \, y(v^*)) 
- f'_y(v^*, \, y(v^*)) \left [ g'_y(v^*, \, y^*) \right ]^{-1} g'_v(v^*, \, y^*)
= 0.
\end{equation}

Обозначим
\begin{equation}
\label{ch1_lam}
\left [ g'_y(v^*, \, y^*)^T \right ]^{-1} 
f'_y(v^*, \, y(v^*))^T = -\lambda^*.
\end{equation}
Если объединить \eqref{ch1_lagr_impl_pr2} и \eqref{ch1_lam}, то получатся следующие равенства:
$$
f'_v (v^*, \, y(v^*))^T + {g'_v(v^*, \, y^*)}^{T} \lambda^* = 0,
\quad
f'_y(v^*, \, y(v^*))^T + g'_y(v^*, \, y^*)^T \lambda^* = 0,
$$
которые эквивалентны \eqref{ch1_lx_llam}.
Таким образом, обосновано правило множителей Лагранжа,
и \eqref{ch1_lam} даёт  явную формулу для множителей Лагранжа.

Аналогичные рассуждения можно провести, считая, что переменные принадлежат общему гильбертову пространству, см. раздел~\ref{subsec:inverse} (или даже банахову пространству, см. раздел~\ref{ch1_implicit_sect}).


\subsection[Гладкие задачи с ограничениями типа равенств\\ и неравенств]{Гладкие задачи с ограничениями типа равенств\\ и неравенств} 
В этом параграфе изложение частично следует~\cite{Boyd_book_2004}.

Рассмотрим задачу с ограничениями типа равенств и неравенств вида
\begin{eqnarray}
\label{ch1_lagr_all_problem}
& \min & f(x) \\  
& g_i(x) = 0, \, i = 1, \, \ldots, \, m, &\nonumber \\
& h_i(x) \, \le \, 0, \, i = 1, \, \ldots, \, p, & \nonumber
\end{eqnarray}
где $ x \, \in \, \mathbb{R}^n$; $f$, $g_i,  \, i = 1, \, \ldots, \, m$, $h_i,  \, i = 1, \, \ldots, \, p$~--- достаточно гладкие функции. 
Ясно, что задача~\eqref{ch1_lagr_eqs_problem}
является частным случаем задачи~\eqref{ch1_lagr_all_problem}.

Пусть область определения
задачи~\eqref{ch1_lagr_all_problem}
\begin{equation}
\label{ch1_eq1_dom_all_prob}
\mathbb{D} =  \textup{dom } f \, \cap \, \bigcap_{i = 1}^m \textup{dom } g_i(x) \, \cap \, \bigcap_{i = 1}^p \textup{dom } h_i(x)
\end{equation}
не является пустой.
Здесь и далее будем обозначать $x^*$~--- оптимальное
решение задачи~\eqref{ch1_lagr_all_problem},
а оптимальное значение $f(x^*) = f^*$. Такие
же обозначения общеприняты и для задач математического
программирования общего вида.
Как уже упоминалось, функцию~$f$ называют
{\it целевой}. 

Точку~$x' \, \in \, \mathbb{R}^n$ называют {\it допустимой}, если~$x' \, \in \, \mathbb{D}$,
$g_i(x') = 0, \, i = 1, \, \ldots, \, m$, 
$h_i(x') \, \le \, 0, \, i = 1, \, \ldots, \, p$,
т.е. если~$x'$ удовлетворяет всем ограничениям
задачи. {\it Допустимым множеством} задачи~\eqref{ch1_eq1_dom_all_prob} называют
множество всех допустимых точек.

Ограничение-неравенство $h_i(\cdot) \, \le \, 0$ называется
{\it активным} в точке~$x'$, если в этой
точке оно обращается в равенство: $h_i(x') = 0$.
Ограничения-равенства по определению всегда
активны в допустимой точке.

Как уже было показано в предыдущем параграфе, для построения
{\it лагранжиана} нужно сложить целевую функцию со взвешенной суммой функций-ограничений. Для задачи~\eqref{ch1_lagr_all_problem} {\it функция Лагранжа}
имеет вид
\begin{equation}
\label{ch1_Lagr_all_func}
L(x, \, \lambda, \, \nu) = f(x) + \sum_{i = 1}^m \lambda_i g_i (x) + \sum_{i = 1}^p \nu_i h_i (x)
\end{equation}
с областью определения $\text{dom } L = \mathbb{D} \, \times \, \mathbb{R}^m \, \times \, \mathbb{R}_+^p$.
{\it Множители Лагранжа}~$\lambda_i$ связаны соответственно
с функциями~$g_i(x)$ в ограничениях-равен\-ствах задачи~\eqref{ch1_lagr_all_problem}, а 
{\it множители Лагранжа}~$\nu_i$ связаны соответственно
с функциями~$h_i(x)$ в ограничениях-неравенствах задачи~\eqref{ch1_lagr_all_problem}. Напомним, что векторы
$\lambda = (\lambda_1, \, \ldots, \, \lambda_m)$
и $\nu = (\nu_1, \, \ldots, \, \nu_m)$
называют {\it двойственными переменными} задачи~\eqref{ch1_lagr_all_problem}.

Впервые задача с неравенствами~\eqref{ch1_lagr_all_problem}
была рассмотрена в
диссертации В.~Каруша (1939~г., Чикаго), но результаты остались
неопубликованными. В 1948~г. задача~\eqref{ch1_lagr_all_problem}
и метод её решения с помощью функции Лагранжа были
переоткрыты Ф.~Джоном. Подробно история этих открытий описана в научно-популярной статье~\cite{About_Karush}.

Существует много разных способов обоснования принципа множителей Лагранжа. Отметим здесь, в частности, работы сотрудников МФТИ~\cite{Birykov}, \cite{Evtushenko2018}.

\subsubsection{Двойственная функция и нижняя оценка}
{\it Двойственная функция Лагранжа} для задачи~\eqref{ch1_lagr_all_problem} определяется
как
$$
\varphi (\lambda, \, \nu) = \inf_{x \, \in \, \mathbb{D}}
L(x, \, \lambda, \, \nu) = \inf_{x \, \in \, \mathbb{D}}
\left \{ f(x) + \sum_{i = 1}^m \lambda_i g_i (x) + \sum_{i = 1}^p \nu_i h_i (x) \right \},
$$
область определения $\text{dom } \varphi = \mathbb{R}^m \,
\times \, \mathbb{R}_+^p$.

Если лагранжиан~\eqref{ch1_Lagr_all_func} не ограничен
снизу, двойственная функция принимает значение~$- \infty$.
Поскольку двойственная функция представляет собой поточечные нижние грани аффинных функций от $\lambda$ и $\nu$, она является вогнутой, даже
если задача~\eqref{ch1_lagr_all_problem} не является
выпуклой задачей.

\begin{lemma}[{\emph{о нижней оценке оптимального значения}}]\label{ch1_lem_low_bound}
Любое значение двойственной функции ${\varphi (\lambda,  \nu)}$ всегда не больше оптимального значения~$f^*${\rm :}
$$
\forall \, \nu \, \ge \, 0, \, \forall \, \lambda
\quad \varphi(\lambda, \, \nu) \, \le \, f^*.
$$
\end{lemma}
Здесь и далее векторные неравенства 
со знаками $\ge$ и $\le$ являются покомпонентными.

Действительно, пусть $x'$ --- {\it допустимая точка} для задачи~\eqref{ch1_lagr_all_problem}, т.е. 
$g_i(x') = 0$, $h_i(x') \, \le \, 0$. Пусть $\nu \, \ge \, 0$. Тогда
$$
\sum_{i = 1}^m \lambda_i g_i (x') + \sum_{i = 1}^p \nu_i h_i (x') \, \le \, 0,
$$
поскольку каждое слагаемое в первой сумме равно $0$,
а каждое слагаемое во второй сумме неположительно.
Отсюда
$$
L(x', \, \lambda, \, \nu) = f(x') + \sum_{i = 1}^m \lambda_i g_i (x') + \sum_{i = 1}^p \nu_i h_i (x') \, \le \, f(x').
$$
Следовательно,
$$
\varphi(\lambda, \, \nu) = \inf_{x \, \in \, \mathbb{D}}
L(x, \, \lambda, \, \nu) \, \le \, L(x', \, \lambda, \, \nu) \, \le \, f(x').
$$
Поскольку
$$
\varphi(\lambda, \, \nu) \, \le \, f(x')
$$
для любого допустимого~$x'$, то это верно и при $x' = x^*$:
$$
\varphi(\lambda, \, \nu) \, \le \, f(x^*) = f^*,
$$
что и требовалось доказать.

Итак, для каждой пары $(\lambda, \, \nu)$ c $\nu \, \ge \, 0$ двойственная функция Лагранжа даёт нижнюю оценку оптимального значения $f^*$ задачи~\eqref{ch1_lagr_all_problem}. Возникает вопрос:
а можно ли найти {\it наилучшую} нижнюю оценку
среди всех нижних оценок, которые можно получить
с помощью двойственной функции Лагранжа?
Ответом на этот вопрос служит
решение задачи
\begin{equation}
\label{ch1_eq_lagr_dual_prob}
\max_{\nu \, \ge \, 0} \varphi(\lambda, \, \nu),    
\end{equation}
которая является задачей выпуклого программирования\footnote{См. далее параграф~\ref{ch1_subsect_conv_prog}.}, 
поскольку целевая функция $\varphi(\lambda, \, \nu)$
вогнута, а ограничения задачи выпуклы.

Задача~\eqref{ch1_eq_lagr_dual_prob} называется {\it двойственной} к задаче~\eqref{ch1_lagr_all_problem}.
Последнюю задачу поэтому иногда называют {\it прямой}.

Пару $(\lambda^*, \, \nu^*)$, являющуюся решением задачи~\eqref{ch1_eq_lagr_dual_prob}, называют
{\it оптимальными} множителями Лагранжа.

\subsection{Чувствительность в оптимизации}
Множители Лагранжа часто имеют вполне конкретную интерпретацию в различных интересных практических
задачах.
Например, в экономических приложениях множители
Лагранжа часто интерпретируются как цены.

В практических приложениях нередко возникает вопрос:
насколько изменится оптимальное значение
целевой функции задачи
нелинейного программирования
при незначительных изменениях функций-ограничений?
Нахождение такой зависимости и называют 
исследованием чувствительности задачи оптимизации.
Оказывается, что для оценки чувствительности
можно использовать множители Лагранжа.

Для примера рассмотрим простую задачу с одним аффинным ограничением~\cite{Bertsekas_conv_book}:
\begin{equation}
\label{ch1_eq_sens_prob}
\min_{a^T x = b} f(x), \quad x \, \in \, \mathbb{R}^n, \quad a \neq 0.
\end{equation}
Пусть $x^*$ --- как обычно, решение задачи~\eqref{ch1_eq_sens_prob} и соответственно $\lambda^*$~--- оптимальный множитель Лагранжа.

\begin{figure}[htb]
\begin{center}
\includegraphics[width=0.7\linewidth]{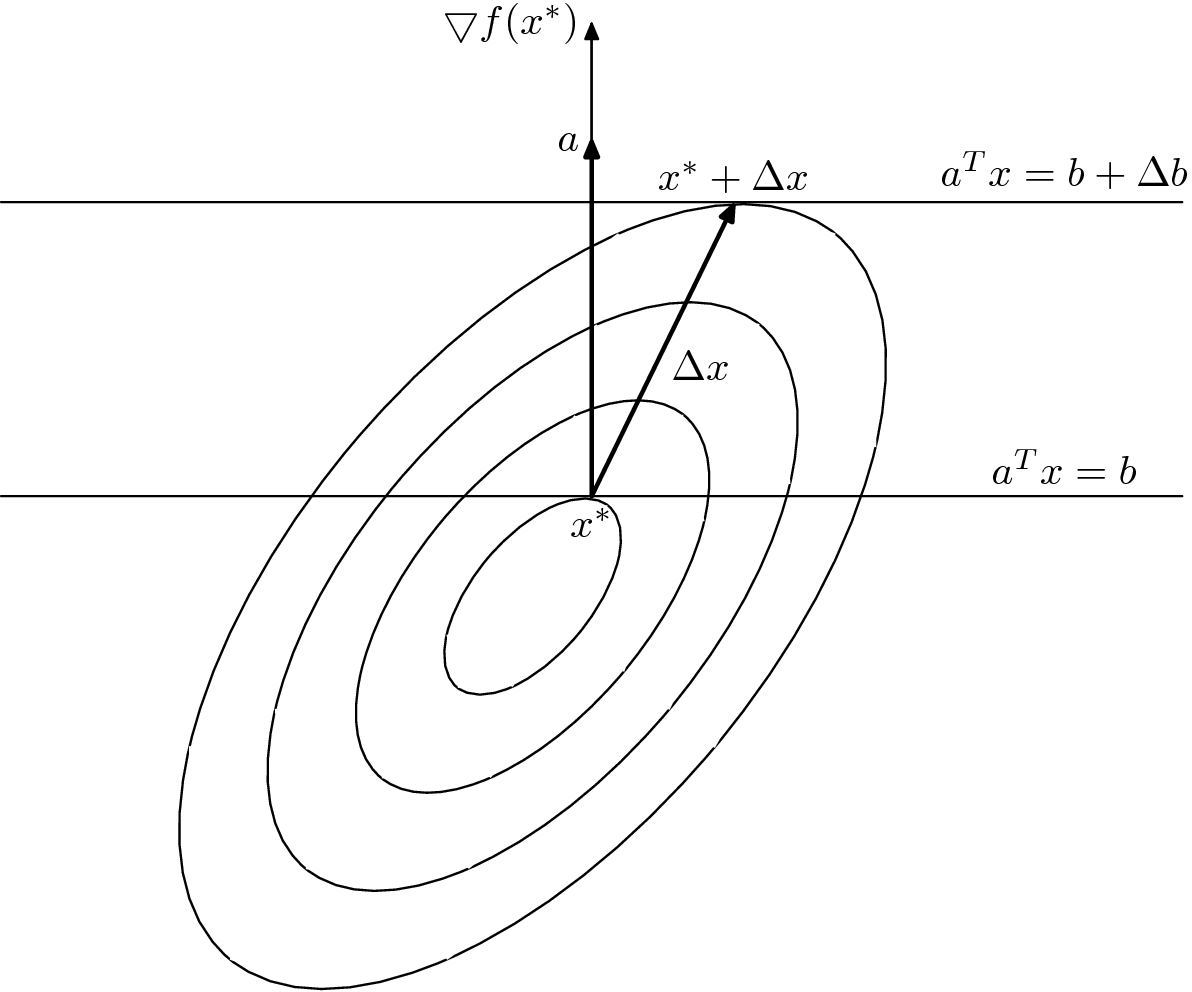}
\end{center}
\caption{Интерпретация множителей Лагранжа (\cite{Bertsekas_conv_book}, с изменениями обозначений)} \label{ch1_dr_sens}
\end{figure}

Если правая часть ограничения $a^Tx = b$ изменится на $\Delta b$ (рис.~\ref{ch1_dr_sens}):
$$
b \mapsto b + \Delta b,
$$
то оптимальное значение $x^*$ изменится на какое-то $\Delta x$.
При этом $b + \Delta b =\linebreak =  a^T (x^* + \Delta x) =a^T x^* + a^T \Delta x = b + a^T \Delta x$, т.е.
\begin{equation}
\label{ch1_eq_deltab}
\Delta b = a^T \Delta x.
\end{equation}
Используя~\eqref{ch1_grad_Lagr_func}, получаем,
что $\nabla f(x^*) = -\lambda^* a$.
Оценим соответствующее изменение оптимального значения целевой функции:
$$
\Delta f = f(x^* + \Delta x) - f(x^*) = \left ( \nabla
f(x^*) \right )^T \Delta x + o(\| \Delta x \|) = -\lambda^* a^T \Delta x + o \left (\| \Delta x \| \right ).
$$
Учитывая~\eqref{ch1_eq_deltab}, получим
$\Delta f = -\lambda^* \Delta b + o \left (\| \Delta x \| \right )$,
поэтому
с точностью до первого порядка по $\| \Delta x \|$
\begin{equation}
\label{ch1_eq_sens_lam_delta}
\lambda^* = - \frac{\Delta f}{\Delta b},
\end{equation}
т.е. величина множителя Лагранжа соответствует скорости
уменьшения~$\Delta f$ при увеличении $b$ на $\Delta b$.
Иными словами, $\lambda^*$ можно проинтерпретировать
как {\it цену ограничения}.

К соотношению, аналогичному~\eqref{ch1_eq_sens_lam_delta},
можно прийти и таким путём~\cite{boss_opt}.
Лагранжиан для задачи~\eqref{ch1_eq_sens_prob}
$$
L(x, \, \lambda) = f(x) + \lambda (a^Tx - b).
$$
Рассматривая лагранжиан как функцию от параметра~$b$
(соответственно, искомое решение $x^*(b)$, а
множитель Лагранжа~--- $\lambda^*(b)$),
получаем
$$
\frac{d L(b)}{db} = \frac{d L(x^*(b), \, \lambda^*(b))}{db} =
$$
\begin{equation}
\label{ch1_sens_dldb}
= \sum_{i=1}^n \frac{\partial f}{\partial x_i} \frac{d x_i^*}{db} + \frac{d \lambda^*}{d b} \left ( a^T x^*(b) - b \right ) +
\lambda^*(b) \left ( \sum_{i=1}^n a_i \frac{d x_i^*}{d b} - 1 \right ).
\end{equation}
Поскольку $a^T x^*(b) \equiv b$, из~\eqref{ch1_sens_dldb}
получаем
$$
\frac{dL}{d b} = \sum_{i=1}^n \frac{\partial f}{\partial x_i} \frac{d x_i^*}{db}  \overset{\text{def}}{=} 
\frac{df^*}{db}. 
$$
С другой стороны, если оставить в
правой части равенства~\eqref{ch1_sens_dldb}
только первое и третье слагаемые, можно прийти к
такому равенству:
$$
\frac{dL}{db} = \sum_{i=1}^n \left ( \frac{\partial f}{\partial x_i} + \lambda^*(b) a_i \right )
\frac{dx_i^*}{db}  - \lambda^*(b).
$$
В силу~\eqref{ch1_grad_Lagr_func} из последнего равенства получаем, что
$\frac{dL}{db} = -\lambda^*(b)$, т.е.
$$
\frac{df^*}{db} = -\lambda^*(b),
$$
что означает, что множитель Лагранжа можно
интерпретировать как показатель
чувствительности оптимума по отношению
к ограничению.

\begin{remark}
Если в задаче~\eqref{ch1_eq_sens_prob} имеются
несколько ограничений вида
$(a^i)^T x = b_i$, $i = 1, \, \ldots, \, m$,
то с помощью аналогичных преобразований можно
прийти к соотношению
$$
\Delta f = - \sum_{i = 1}^m \lambda_i^* \Delta b_i +
o \left (\| \Delta x \| \right ). 
$$
В случае существования всех нужных зависимостей
$x^*(b)$, $\lambda^*(b)$ каждый из множителей
Лагранжа интерпретируется следующим образом:
\begin{equation}
\label{ch1_eq_sens_ld_i}
\lambda_i^* = -\frac{\partial f^*}{\partial b_i},
\quad i = 1, \, \ldots, \, m.
\end{equation}
\end{remark}

\begin{remark}
Нетрудно заметить, что к абсолютно тем же результатам~\eqref{ch1_eq_sens_lam_delta}, \eqref{ch1_eq_sens_ld_i} можно
прийти и в более общем случае ограничений
вида $g_i(x) = b$ в задаче~\eqref{ch1_eq_sens_prob}.
\end{remark}

\begin{remark}
Более полное изложение современной теории чувствительности для конечномерных задач оптимизации см. в учебнике~\cite{Bertsekas_nonlinear} и монографии~\cite{Ism_Sens}. 
\end{remark}

\subsection{Задачи выпуклого программирования}\label{ch1_subsect_conv_prog}

\subsubsection{Задачи без ограничений}\label{ch1_subs_fermat}
Пусть $X$~--- векторное пространство,
$f \, : \, X \, \rightarrow \, \bar{\mathbb{R}}$ --- выпуклая функция с областью определения $\text{dom } f = X$.

Выпуклой задачей без ограничений называется задача
\begin{equation}
\label{ch1_eq_min_conv}
\min_{x \, \in \, X} f(x).
\end{equation}



\subsubsection{Постановка задачи выпуклого программирования}

{\it Задачей выпуклого программирования} называется
задача минимизации выпуклой функции на выпуклом множестве.
Пусть в нашей формулировке имеются выпуклые ограничения-неравенства:
\begin{eqnarray}
\label{ch1_eq_min_conv_in}
& \min & f(x) \\  
& h_i(x) \, \le \,  0, \, i = 1, \, \ldots, \, p, &\nonumber \\
& x \, \in \, Q & \nonumber
\end{eqnarray}
где $ x \, \in \, \mathbb{R}^n$; $f$, $h_i,  \, i = 1, \, \ldots, \, p$~--- достаточно гладкие выпуклые функции, $Q$~---~выпуклое подмножество в $\text{dom } f \, \cap \, \bigcap_{i = 1}^p\, \text{dom } h_i$.


{\it Двойственная функция Лагранжа} для задачи~\eqref{ch1_eq_min_conv_in} определяется
как
$$
\varphi (\nu) = \inf_{x \, \in \, Q}
L(x, \, \nu) = \inf_{x \, \in \, Q}
\left \{ f(x) +  \sum_{i = 1}^p \nu_i h_i (x) \right \},
$$
область определения $\text{dom } \varphi = \mathbb{R}_+^p$.

В выпуклой задаче локальный минимум всегда будет и глобальным. Этот факт можно доказать с помощью неравенства Йенсена (см., например, \cite{galeev89}).

Сделать обоснование метода множителей Лагранжа
(и вывода необходимых условий экстремума)
можно и с помощью теоремы об отделимости
непересекающихся выпуклых множеств.
В следующем параграфе 
приведены формулировки двух теорем об отделимости,
касающиеся конечномерного случая,
а на с.~\pageref{ch1_subs_teo_kkt}
приведено доказательство принципа Лагранжа для задачи выпуклого
программирования~\eqref{ch1_eq_min_conv_in}.

Следует отметить, что формулировка
задачи только с ог\-ра\-ни\-че\-ниями-не\-ра\-вен\-ства\-ми,
без ограничений-равенств, не умаляет общности
постановки задачи. Действительно, каждое ограничение-равенство
вида $Ax - b = 0$ можно заменить двумя ограничениями-неравенствами $Ax - b \, \le \, 0$ и
$Ax - b \, \ge \, 0$.

\begin{figure}[htb]
\begin{center}
\includegraphics[width=0.8\linewidth]{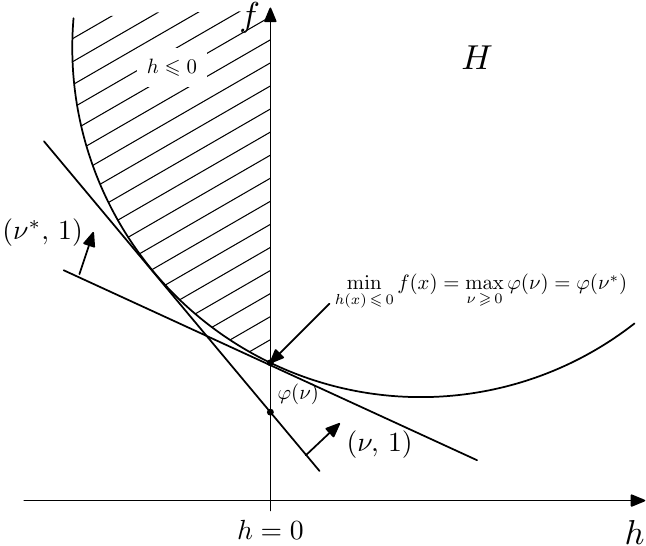}
\end{center}
\caption{Геометрическая интерпретация существования
множителей Лагранжа для задачи выпуклого программирования~\eqref{ch1_eq_min_conv_in}. В данной задаче одно ограничение ${h(x)\!\le\!0}$, множество ${Q\!=\!\mathbb{R}^n}$,
$\varphi (\nu)$~-- двойственная функция Лагранжа для задачи~\eqref{ch1_eq_min_conv_in}
} 
\label{ch1_dr_geom_int1}
\end{figure}

\begin{figure}[htb]
\begin{center}
\includegraphics[width=0.7\linewidth]{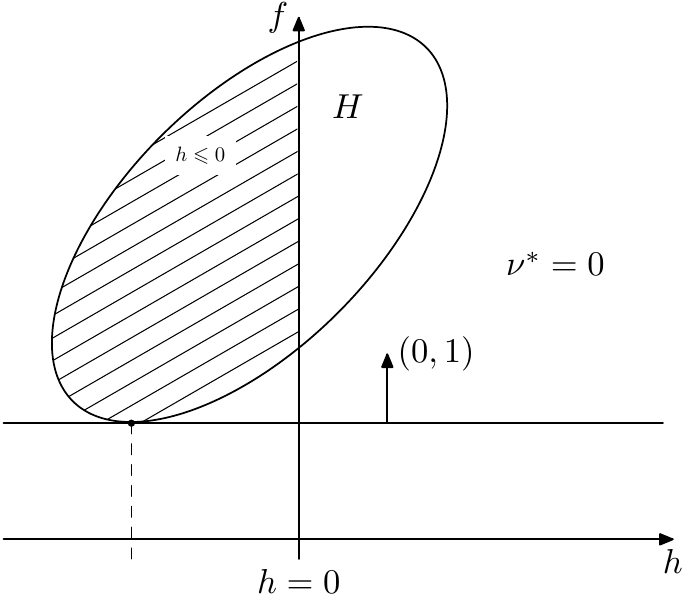}
\end{center}
\caption{Геометрическая интерпретация существования
множителей Лагранжа для задачи выпуклого программирования~\eqref{ch1_eq_min_conv_in}. В данном примере множитель Лагранжа $\nu^* = 0$, т.е. соответствующее ему ограничение неактивно} 
\label{ch1_dr_geom_int2}
\end{figure}

На рис.~\ref{ch1_dr_geom_int1}--\ref{ch1_dr_geom_int2} представлена геометрическая
интерпретация данного факта.
В общем случае множество
$$
H = \left \{  t  \in  \mathbb{R}^{p + 1}  | \,  \exists \, x  \in  \text{dom } f  \cap  \text{dom } h(x), \, f(x) \le  t_0, \, h_i(x) \le  t_i, \,
i = 1, \, \ldots, \, p \right \},
$$
и вектор весов~$(\nu, 1)$ определяет опорную для~$H$ плоскость (см. \nk{определение}~\ref{ch1_def_sup_hyp}).

\subsubsection{Теоремы отделимости}\label{ch1_subs_otdel}
В этом параграфе приведены определения отделимых множеств
и формулировки теорем об отделимости в конечномерном случае. Доказательства теорем см., например, в~\cite{galeev89}.
Теоремы об отделимости применяются
при выводе необходимых условий экстремума (см. с.~\pageref{ch1_subs_teo_kkt}).

\begin{defin}\label{ch1_def_otdel}
Непустые выпуклые множества $A$ и $B$ называются (собственно) \emph{отделимыми}, если существует аффинная гиперплоскость
$$
H = \left \{x \, \in \, \mathbb{R}^n \, | \, \nu^T x = b \right \},\quad
\nu \, \in \, \mathbb{R}^n, \, \nu \neq 0,
$$
такая, что множества $A$, $B$ лежат в противоположных
по отношению к~$H$ замкнутых
полупространствах
и, по крайней мере, одно
из множеств~$A$, $B$ не
содержится в $H$. В этом случае гиперплоскость $H$ \emph{отделяет} множества $A$, $B$.
\end{defin}

Аналитически определение~\ref{ch1_def_otdel}
записывается так:
\begin{equation}
\label{ch1_eq_otdel}
\exists \, H = \left \{x \, \in \, \mathbb{R}^n \, | \, \nu^T x = b \right \} :
\inf_{x \, \in \, A}  \nu^T x  \,
\ge \, b \, \ge \, \sup_{x \, \in \, B}  \nu^T x,
\
\sup_{x \, \in \, A}  \nu^T x  \, > \, \inf_{x \, \in \, B}  \nu^T x.
\end{equation}

\begin{defin}\label{ch1_def_st_otdel}
Непустые выпуклые множества $A$ и $B$ называются \emph{строго отделимыми}, если 
существуют две различные 
параллельные гиперплоскости
$$
H_1 = \left \{x \, \in \, \mathbb{R}^n \, | \, \nu^T x = b_1 \right \},
\quad
H_2 = \left \{x \, \in \, \mathbb{R}^n \, | \, \nu^T x = b_2 \right \}, \quad
\nu \, \in \, \mathbb{R}^n \setminus \{0\},
$$
такие, что
$$
\inf_{x \, \in \, A}  \nu^T x  \,
\ge \, b_1 > b_2 \, \ge \, \sup_{x \, \in \, B}  \nu^T x.
$$
\end{defin}

\begin{figure}[htb]
\begin{center}
\includegraphics[width=0.6\linewidth]{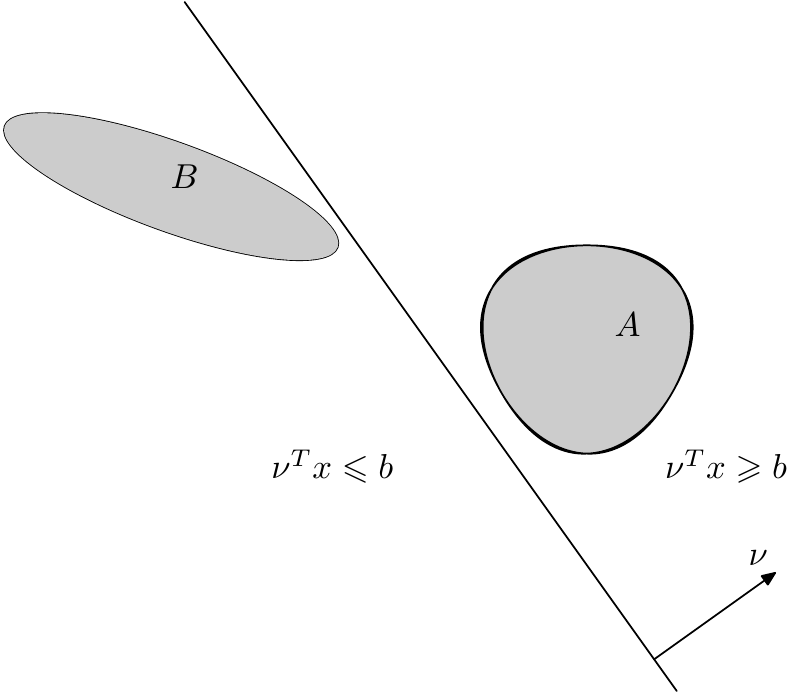}
\end{center}
\caption{Иллюстрация к теореме об отделимости~\ref{ch1_teo_otdel1}. Гиперплоскость $\{ x\!\!\!\in \linebreak \in \mathbb{R}^n \, | \, \nu^T x = b \}$ отделяет
непересекающиеся выпуклые множества $A$ и $B$} \label{ch1_dr_otdel2-19}
\end{figure}


\begin{teo}[первая теорема отделимости в конечномерном случае]\label{ch1_teo_otdel1}
Пусть $A$ и $B$~--- непустые выпуклые множества в $\mathbb{R}^n$, $A \, \cap \, B = \varnothing$.
Тогда множества $A$ и $B$ отделимы.
\end{teo}

\begin{remark}
Можно доказать и более сильную
версию теоремы~\ref{ch1_teo_otdel1}: непустые выпуклые множества
отделимы тогда и только тогда, когда их относительные внутренности не пересекаются. Здесь относительная внутренность множества $A$~--- это его внутренность как подмножества аффинной оболочки множества $A$. Доказательство можно найти в~\cite[теорема~3.16]{Hildebrand_lectures_opt}.
\end{remark}

\begin{figure}[htb]
\begin{center}
\includegraphics[width=0.3\linewidth]{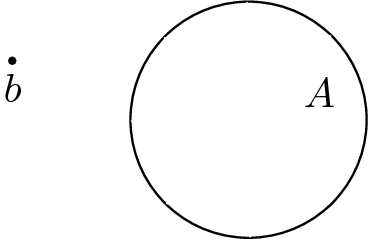}
\end{center}
\caption{Иллюстрация к теореме~\ref{ch1_teo_otdel2}} \label{ch1_dr_otdel}
\end{figure}

\begin{teo}[вторая теорема отделимости в конечномерном случае]\label{ch1_teo_otdel2}
Пусть $A$~--- непустое выпуклое замкнутое множество в~$\mathbb{R}^n$ и~пусть точка~${b  \notin  A}$.
Тогда точку~$b$ можно строго отделить от~$A$ {\rm(}см.~рис.~{\rm\ref{ch1_dr_otdel})}.
\end{teo}

\subsubsection{Теорема Каруша--Куна--Таккера}\label{ch1_subs_teo_kkt}
\begin{teo}\label{ch1_teo_kkt}
Пусть выполняются условия задачи~\eqref{ch1_eq_min_conv_in}.
\begin{itemize}
    \item 
Тогда если $x^*$~--- решение задачи~\eqref{ch1_eq_min_conv_in}, то существует
ненулевой вектор множителей Лагранжа $\nu = (\nu_0, \, \nu_1, \, \ldots, \, \nu_p)$, такой, что для функции Лагранжа
$$
L(x, \, \nu) = \nu_0 f(x) + \sum_{i = 1}^p \nu_i h_i(x)
$$
выполняются:
\begin{enumerate}
\renewcommand\labelenumi{\textup{\arabic{enumi})}}
    \item принцип минимума для функции Лагранжа
    $$
    \min_{x \, \in \, Q} L(x, \, \nu) = L(x^*, \, \nu);
    $$
    \item условия дополняющей нежёсткости
    $$
    \nu_i h_i(x^*) = 0, \, i = 1, \, \ldots, \, p;
    $$
    \item условия неотрицательности
    $$
    \nu_i \, \ge \, 0, \, i = 0, \, \ldots, \, p.
    $$
\end{enumerate}
\item Если $\nu_0 \neq 0$, то условия \textup{1)\,--\,3)} предыдущего утверждения достаточны для того, чтобы
допустимая точка~$x^*$ была решением задачи~\eqref{ch1_eq_min_conv_in}.
\end{itemize}
\end{teo}

Докажем первое утверждение данной теоремы (более подробное доказательство см. в~\cite{galeev89}).
Не ограничивая общности, будем считать, что $f(x^*) = 0$.

1. Пусть 
\begin{align*} \label{ij_quadratic} 
	B = \left \{b = (b_0, \, b_1, \, \ldots, \, b_p) \, \in 
	\, \mathbb{R}^{p+1} \, | \, \exists \, x \, \in \, Q :
	f(x) \, \right.&\le \, b_0, \,
	h_i(x) \, \le \\ &\le \,\left. b_i, \, i = 1, \, \ldots, \, p
	\right \}.
\end{align*}
Докажем, что $B$ является непустым выпуклым множеством.
Если в качестве $x$ взять $x^*$, то любой вектор~$b$ со всеми неотрицательными компонентами
будет принадлежать множеству~$B$, т.е. $B$ непусто.
Для доказательства выпуклости~$B$ возьмём два вектора
$b \, \in \, B$ и $b' \, \in \, B$. Пусть
$x$ и $x'$~--- такие элементы из~$Q$, что 
$f(x) \, \le \, b_0$,
$h_i(x) \, \le \, b_i$, 
$f(x') \, \le \, b'_0$,
$h_i(x') \, \le \, b'_i$ соответственно, $i = 1, \, \ldots, \, p$. Для $\alpha \, \in \, [0, \, 1]$
$x_\alpha = \alpha x + (1 - \alpha) x' \, \in \, Q$
в силу выпуклости~$Q$.
Ввиду выпуклости функций~$f$, $h_i$, $i = 1, \, \ldots,
\, p$,
$$
f(x_\alpha) = f(\alpha x + (1 - \alpha) x') 
\, \le \, \alpha f(x) + (1 - \alpha) f(x') \, \le \,
\alpha b_0 + (1 - \alpha) b'_0
$$
и
$$
h_i(x_\alpha) = h_i(\alpha x + (1 - \alpha) x') 
\, \le \, \alpha h_i(x) + (1 - \alpha) h_i(x') \, \le \,
\alpha b_i + (1 - \alpha) b'_i,
$$
т.е. точка $\alpha b + (1 - \alpha) b' \, \in \, B$.

2. Обозначим $C = \left \{c = (c_0, \, 0, \, \ldots, \, 0) \, \in \,
\mathbb{R}^{p+1} \, | \, c_0 \, < \, 0 \right \}$.

Очевидно, множество $C$ является непустым и выпуклым.
Кроме того, $C \cap B = \varnothing$.
Действительно, предположим противное. Пусть существует
точка~$c=(c_0, \, 0, \, \ldots, \, 0)$, $c_0 \, < \, 0$,
$c \, \in \, B$. Тогда по определению множества~$B$
 существует $\tilde{x} \, \in \, Q$, для которого
 $$
 f(\tilde{x}) \, \le \, c_0 < 0, \quad
 h_i(\tilde{x}) \, \le \, 0, \, i = 1, \, \ldots, \, p.
$$
Отсюда следует, что $\tilde{x}$~--- допустимая точка
задачи~\eqref{ch1_eq_min_conv_in} и
$$
f(\tilde{x}) \, < f(x^*),
$$
что противоречит тому, что $x^*$~--- решение задачи~\eqref{ch1_eq_min_conv_in}. Следовательно,
наше предположение неверно и $C \cap B = \varnothing$.

По первой теореме отделимости в конечномерном случае
(см. п.~\ref{ch1_subs_otdel}) множества~$B$ и $C$
можно отделить, т.е.
существует вектор
$\nu \, \in \, \mathbb{R}^{p+1} \, (\nu \neq 0)$,
такой, что
\begin{equation}
\label{ch1_eq_bc}
\inf_{b \, \in \, B} \sum_{i=0}^p \nu_i b_i \,
\ge \, \sup_{c \, \in \, C} \sum_{i=0}^p \nu_i c_i.
\end{equation}

Так как $0 \, \in \, B$, то из последнего неравенства
следует, что
$$
0 \, \ge \, \sup_{c \, \in \, C}  \sum_{i=0}^p \nu_i c_i = 
\sup_{c \, \in \, C} \nu_0 c_0.
$$
Поскольку $c_0 < 0$, то $\nu_0 \, \ge \, 0$ и, следовательно, $\sup_{c \, \in \, C} \nu_0 c_0 = 0$.
Тогда~\eqref{ch1_eq_bc} можно переписать так:
\begin{equation}
\label{ch1_eq_nub}
\sum_{i=0}^p \nu_i b_i \, \ge \, 0 \quad \forall \, b \, \in \, B.    
\end{equation}

3. Докажем, что множители Лагранжа $\nu_i$, $i = 0, \, 1, \, \ldots, \, p$ удовлетворяют условиям неотрицательности.
Вектор 
$$
\hat{b}^i = (0, \, \ldots, \, 0, \, 1, \, 0, \, \ldots, \, 0),
$$
в котором единица стоит на~$i$-м месте,
принадлежит множеству~$B$, как и любой вектор с неотрицательными компонентами.
Подставив $\hat{b}^i$ в~\eqref{ch1_eq_nub}, получим,
что $\nu_i \, \ge \, 0$, что и требовалось.

4. Докажем, что множители $\nu_i$, $i = 1, \, \ldots, \, p$
удовлетворяют условиям дополняющей нежёсткости. Для
$h_i(x^*) = 0$ равенство
$\nu_i h_i(x^*) = 0$ тривиально. Пусть $h_i(x^*) \neq 0$ для некоторого $i$, тогда $h_i(x^*) < 0$.
Возьмём в определении множества~$B$ в качестве $x$ точку $x^*$. Тогда вектор
$$
{b}_{h_i} = (0, \, \ldots, \, 0, \, h_i(x^*), \, 0, \, \ldots, \, 0),
$$
в котором число $h_i(x^*)$ стоит на~$i$-м месте,
принадлежит множеству~$B$.
Следовательно, можно подставить $b_{h_i}$ в~\eqref{ch1_eq_nub}.
Получаем, что $\nu_i h_i(x^*) \, \ge \, 0$. Тогда $\nu_i$
должно быть неположительным (поскольку $h_i(x^*) < 0$).
Но в п.~3 доказано, что все~$\nu_i$ неотрицательны.
Таким образом, $\nu_i = 0$ и $\nu_i h_i(x^*) = 0$, $i = 1, \, \ldots, \, p$.

5. Докажем принцип минимума для функции Лагранжа\footnote{Обратите внимание, что для задачи
математического программирования в общем виде~\eqref{ch1_lagr_all_problem} это свойство формулируется немного иначе и называется нижней оценкой оптимального значения с помощью двойственной функции Лагранжа, см. лемму~\ref{ch1_lem_low_bound}.}.
Пусть~$x \, \in \, Q$. Тогда вектор
$$
b_{fh} = (f(x), \, h_1(x), \, h_2(x), \, \ldots, \, h_p(x)) 
$$
принадлежит~$B$.
Следовательно, мы можем подставить~$b_{fh}$ в~\eqref{ch1_eq_nub}:
$$
\nu_0 f(x) + \sum_{i = 1}^p \nu_i h_i(x) = L(x, \, \nu) \, \ge \, 0.
$$

{
\renewcommand{\baselinestretch}{0.96}
\selectfont

Поскольку $f(x^*) = 0$ и выполняются условия дополняющей нежёсткости, то для любого $x \, \in \, Q$
$$
L(x, \, \nu) \, \ge \, 0 = \nu_0 f(x^*) + \sum_{i = 1}^p \nu_i h_i(x^*) =  L(x^*, \, \nu),
$$
что и требовалось доказать. $\qed$

\subsection{Сильная и слабая двойственность}\label{ws_duality}
Обозначим оптимальное значение для двойственной задачи~\eqref{ch1_eq_lagr_dual_prob} через~$\varphi^*$.
Тогда можно записать очень важное свойство
оптимального значения целевой функции:
\begin{equation}
\label{ch1_eq_in_phi_f}
\varphi^* = \varphi(\lambda^*, \, \nu^*) \, \le \, f^*=f(x^*)
\end{equation}
(см. лемму~\ref{ch1_lem_low_bound} о нижней оценке оптимального значения).

Неравенство~\eqref{ch1_eq_in_phi_f} выполняется, даже если исходная задача невыпуклая\footnote{Кроме того, неравенство выполняется, даже если $\varphi^* = \infty$
и $f^* = \infty$.}, и называется свойством
{\it слабой двойственности} (weak duality).

\begin{defin}
Разность $f^* - \varphi^*$ называется \textit{зазором двойственности} или зазором между решениями прямой и двойственной задач.
\end{defin}

Ясно, что зазор двойственности всегда неотрицателен.
Если имеет место равенство
\begin{equation}
\label{ch1_eq_phi_f}
\varphi^* = \varphi(\lambda^*, \, \nu^*) = f^*,
\end{equation}
то говорят о \textit{сильной двойственности}.

\subsubsection{Сильная двойственность в невыпуклом случае}
Свойство сильной двойственности в очень редких случаях выполняется
для невыпуклых задач.
Рассмотрим, например, задачу минимизации квадратичной
функции на единичном шаре:
\begin{equation}
\label{ch1_eq_quad_ball_prob}
\min_{\|x \|_2^2 \, \le \, 1} x^T A x,
\end{equation}
где
$A$~--- симметричная матрица размера $n \times n$ и
$A$ не является неотрицательно определённой.
Вследствие последнего утверждения задача~\eqref{ch1_eq_quad_ball_prob}
не~является выпуклой.

Составим функцию Лагранжа для этой задачи:
$$
L(x, \, \nu) = x^T A x + \nu (x^T x - 1) = x^T (A + \nu I) x - \nu,
$$
здесь $I$~--- единичная матрица.

}

Двойственная функция
$$
\varphi(\nu) = \inf_x L(x, \, \nu) =
\left \{
\begin{array}{ll}
-\nu, & \text{если } A + \nu I \succeq 0, \\
-\infty, & \text{иначе}.
\end{array}
\right.
$$
Двойственная задача $\max\limits_{\nu \geq 0} \, \varphi(\nu)$ является выпуклой и одномерной. Зазор двойственности в этой случае нулевой:
$$
\{x^*\}^T A x^* = \varphi(\nu^*).
$$

Решением задачи~\eqref{ch1_eq_quad_ball_prob}
является нормированный собственный вектор матрицы $A$, который соответствует её минимальному собственному значению.
Иногда задачу~\eqref{ch1_eq_quad_ball_prob}
называют задачей построения доверительной
области (trust region problem). Такого рода задачи
могут возникать при необходимости нахождения
оптимальных аппроксимаций 
второго порядка различных функций на единичном шаре.

Имеет место и более общий факт наличия сильной двойственности для общей задачи квадратичного программирования, см. далее п.~\ref{ch1_sec_qp}.

\subsubsection[Пример. Задача разбиения множества на две части. Слабая двойственность]{Пример. Задача разбиения множества на две части.\\ Слабая двойственность}
Рассмотрим следующую невыпуклую оптимизационную задачу:
\begin{equation}
\label{ch1_eq_2-way_prob}
\min_{x_i^2 = 1, \, i = 1, \, \ldots, \, n} x^T W x,
\end{equation}
где
$W$~--- симметричная матрица размера $n \times n$.
Ограничения задачи означают, что компоненты вектора~$x$ могут принимать только значения~$1$ или $-1$. 
Задачи оптимизации, в которых некоторые или все переменные должны быть целыми числами,
называют задачами {\it целочисленного программирования}.

Допустимое множество задачи~\eqref{ch1_eq_2-way_prob}
конечно, поскольку содержит всего~$2^n$ точек,
но зависимость от~$n$ экспоненциальная,
поэтому с помощью простого перебора точек
из допустимого множества задачу~\eqref{ch1_eq_2-way_prob}
можно решить только для совсем небольших~$n$,
не более нескольких десятков. В~общем случае задача~\eqref{ch1_eq_2-way_prob} NP-трудна.

Функция Лагранжа задачи~\eqref{ch1_eq_2-way_prob}
имеет вид
\begin{equation}
\label{ch1_eq_2-way_lagr}
L(x, \, \lambda) = x^T W x + \sum_{i = 1}^n \lambda_i
(x_i^2 - 1) = x^T \left (W + \text{diag}(\lambda)
\right ) x - \sum_{i = 1}^n \lambda_i.
\end{equation}
Двойственная функция Лагранжа равна
\begin{equation}
\label{ch1_eq_2-way_dual}    
\varphi(\lambda) = \inf_x L(x, \, \lambda) =
\left \{
\begin{array}{cl}
     - \sum\limits_{i = 1}^n \lambda_i & \text{при } W + \text{diag}(\lambda) \succeq 0, \\
     - \infty, &  \text{иначе}, 
\end{array}
\right.
\end{equation}
поскольку квадратичная форма~\eqref{ch1_eq_2-way_lagr}  
достигает минимума в точке $x = 0$ при условии,
что её матрица $W + \text{diag}(\lambda)$ неотрицательно определена и не ограничена в остальных случаях.

Функция~\eqref{ch1_eq_2-way_dual} даёт нижнюю оценку
оптимального значения задачи~\eqref{ch1_eq_2-way_prob}.
Например, при $\lambda = -\lambda_{\min}(W) \cdot (1, \, \ldots, \, 1)$
нижняя оценка для $f^*$ следующая:
$$
f^* = 
(x^*)^T W x^* \, \ge \, - \sum_{i = 1}^n \lambda_i = n \lambda_{\min}(W).
$$

Задачу~\eqref{ch1_eq_2-way_prob} можно интерпретировать
как задачу разбиения множества из $n$ элементов
$S = \{ 1, \, 2, \, \ldots, \, n \}$ на две части:
$$
\{ 1, \, 2, \, \ldots, \, n \} = \{i \, | \, x_i = -1 \}
\, \cup \,  \{i \, | \, x_i = 1 \}.
$$
Тогда элемент матрицы $w_{ij}$ можно интерпретировать
как стоимость нахождения элементов с номерами $i$ и $j$
из $S$ в одной части, а $-w_{ij}$~--- как стоимость нахождения элементов с номерами $i$ и $j$
из $S$ в разных классах соответственно.
При этом задача~\eqref{ch1_eq_2-way_prob} соответствует
цели построения разбиения с минимальной общей стоимостью.

Мы вернёмся к описанной задаче в п.~\ref{ch1_relaxations}, где посмотрим на неё с позиции возможной (выпуклой) <<релаксации>> для нахождения приближённого решения.

\subsubsection{Условия Слейтера}\label{ch1_subs_slater}
Часто полезно заранее знать, для каких задач выполняется свойство сильной двойственности (т.е. имеется нулевой разрыв двойственности). Одним из таких условий является
\textit{условие Слейтера}.

Рассмотрим задачу выпуклого программирования следующего вида:
\begin{eqnarray}
\label{ch1_eq_conv_sl_prob}
& \min & f(x), \\  
&  h_i(x) \, \le \, 0, \, i = 1, \, \ldots \, p, & \nonumber \\
& Ax = b \nonumber & 
\end{eqnarray}
в которой функции $f$, $h_1$, $\ldots$, $h_p$ являются выпуклыми.

Условие Слейтера (его ещё иногда называют условием регулярности) требует, чтобы для задачи~\eqref{ch1_eq_conv_sl_prob} существовала такая допустимая
точка~$x'$, 
что неравенства в~\eqref{ch1_eq_conv_sl_prob}
являются строгими:
\begin{equation}
\label{ch1_eq_sl_cond}
h_i(x') \, < \, 0, \, i = 1, \, \ldots \, p, \quad Ax' = b.
\end{equation}
Такую точку $x'$ также иногда называют \textit{строго допустимой}.

Теорема Слейтера~\cite{Slater} утверждает, что свойство сильной
двойственности для задачи выпуклого программирования~\eqref{ch1_eq_conv_sl_prob} выполняется,
если выполняются условия~\eqref{ch1_eq_sl_cond}.

Условия Слейтера можно ослабить, если некоторые из
функций-огра\-ни\-чений $h_i$ являются аффинными
(определение
аффинных функций см. в п.~\ref{ch1_sect_conv}).
Пусть первые $k$ функций-ограничений $h_i(x)$, 
$i = 1, \, \ldots, \, k$ аффинные.
Тогда для существования нулевого зазора двойственности достаточно,
чтобы существовала точка $x'$, такая, что
\begin{equation}
\label{ch1_eq_sl_cond_refine}
h_i(x') \, < \, 0, \, i = k + 1, \, \ldots \, p, \quad Ax' = b.
\end{equation}

\subsection{Задача линейного программирования}\label{subsubsect_lp}

Рассмотрим оптимизационную задачу
с линейной целевой функцией и линейными функциями-ограничениями
\begin{eqnarray}
\label{ch1_eq_lp_all}
& \min & c^T x + d, \\  
& Ax = b &\nonumber \\
& Gx \, \ge \, h &\nonumber 
\end{eqnarray}
где $G \, \in \, \mathbb{R}^{p \times n}$
и $A \, \in \, \mathbb{R}^{m \times n}$.
Постоянный вектор $d$ в целевой функции
часто опускают, поскольку он не влияет
на решение задачи. 

Задача~\eqref{ch1_eq_lp_all}, называемая задачей линейного программирования~(ЛП), стала, наверное,
одной из самых известных в ХХ~веке задач оптимизации. В начале 80-х годов XX~века
считалось, что
решение задачи линейного программирования
находится на первом месте в мире среди всех математических задач\footnote{Но в эти задачи
не входила работа с базами данных (сортировка, поиск данных).} по числу
затраченного компьютерного времени.

\begin{figure}[htb]
\begin{center}
\includegraphics[width=0.5\linewidth]{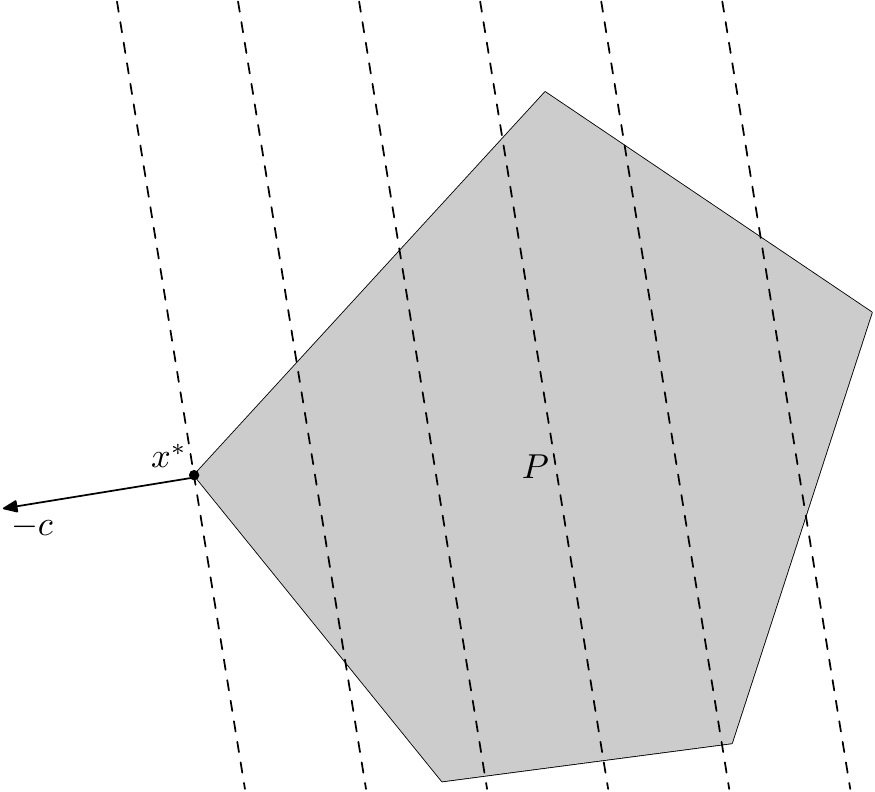}
\end{center}
\caption{Геометрическая интерпретация задачи линейного программирования. $P$~-- допустимое множество (полиэдр\protect\footnotemark), целевая функция~$c^T x$ линейна (её линии уровня показаны пунктиром~-- это гиперплоскости, ортогональные вектору~$c$,
уменьшение значения целевой функции идёт в направлении антиградиента $-c$, поэтому оптимальная точка~$x^*$~-- точка полиэдра~$P$
в направлении $-c$, в которой линия уровня становится касательной к~$P$\nk{)}} \label{ch1_dr_lp_geom}
\end{figure}

Приведём текстовое
описание задачи~ЛП: 
необходимо найти максимум или минимум
линейной функции множества переменных на выпуклом многограннике (задаваемом конечной системой линейных равенств и неравенств).
Задачи~ЛП, конечно, являются подклассом
задач выпуклого программирования. Как отмечено
в текстовом описании, задачи на нахождение
максимума линейной функции при линейных ограничениях также относятся к задачам~ЛП.
Геометрическую иллюстрацию графического
метода решения задачи~ЛП см. на рис.~\ref{ch1_dr_lp_geom}.

С созданием теории линейного программирования начался период активного применения теории и методов оптимизации в экономике и экономических исследованиях. Во-первых, было обнаружено, что целый ряд практических задач организации производства можно формализовать как задачи~ЛП\footnote{Сам Канторович приводит такой список: <<вопросы наилучшего распределения работы станков и механизмов, максимального
уменьшения отходов, наилучшего использования сырья и~местных материалов,
топлива, транспорта и пр.>> (цит. по \cite[с.~43]{Kant_izbr11}).}. 
Во-вторых, к уже
упомянутым экономическим
задачам и
задачам управления
системами, которые позже отнесут
к дисциплине <<исследование
операций>>, возник большой интерес у военных
специалистов. Считается,
что сейчас уже классическая задача~ЛП об оптимальной диете,
описанная и решённая будущим Нобелевским лауреатом
Джорджем Стиглером~\cite{Stigler_diet},
возникла в американской армии в 30--40-е годы
XX века.
Задача о диете ставилась так:
выбрать суточный набор продуктов питания
минимальной стоимости,
который будут обеспечивать
необходимую норму питательных веществ.
Первые версии задачи Стиглер решал эвристическим методом (77~переменных~--- продуктов и 9~ограничений)\footnote{Позже, после появления симплекс-метода Данцига, результаты
были пересчитаны и полностью подтвердились.}.
Оказалось, что для
первой версии задачи оптимальный рацион состоит
из хлеба, молока, капусты, шпината и фасоли.


\footnotetext{Полиэдр (или выпуклый полиэдр)~--- 
множество, которое может быть представлено как пересечение 
конечного числа полупространств в пространстве~$\mathbb{R}^n$.
Множество, являющееся выпуклой оболочкой конечного числа 
точек, называется {\it политопом {\rm(}или выпуклым многогранником{\rm)}}. 
То есть политоп, по сути, ограниченный полиэдр.}

В-третьих, были разработаны и реализованы вычислительные методы, которые сделали реальным решение таких задач 
на ЭВМ, даже относительно больших размеров (как и по числу переменных, так и по числу ограничений).
Что интересно, начальный этап развития~ЛП
на Западе (а это было почти на десятилетие позже, чем в СССР) происходил практически одновременно
с развитием компьютеров, и там это сразу
означало, что
впервые в истории управленцы получили
мощный и практический метод 
поиска оптимального направления
действий из огромного числа возможных направлений.
Так, 
линейное программирование открыло дорогу к практическому выполнению расчётов, направленных на оптимизацию использования ресурсов во многих секторах экономики и военно-промышленном комплексе.

\begin{remark}
Помимо важнейших практических приложений, 
есть ещё одна причина важности
задач линейного программирования:
они часто возникают как один из этапов алгоритмов для решения
задач нелинейной оптимизации.
\end{remark}

Приоритет в создании теории
и моделей задач~ЛП
принадлежит советскому математику
Леониду Витальевичу Канторовичу (1912--1986).
Первая печатная работа 
 <<Математические методы организации
и планирования производства>>
с постановками таких задач
и методом их решения
была опубликована
Канторовичем в 1939 году~\cite{Kant_lp39}.
В качестве метода решения задач~ЛП в
работе был представлен метод разрешающих множителей~--- прообраз симплекс-метода Данцига (см. ниже).

Эта основополагающая работа Канторовича по линейному программированию стала известна и доступна коллегам
из западных стран только через два десятилетия~(!), в 1960~г., когда в журнале <<Management Science>> опубликовали английский
перевод~\cite{Kantor_eng_ManS60}. Инициатива этого перевода и публикации принадлежала
американскому экономисту и математику голландского происхождения
Т.~Купмансу.

В 1940 г. Л.\,В.~Канторович опубликовал
математическую версию брошюры 1939~г.~--- статью <<Об одном эффективном методе peшeния некоторых классов экстремальных проблем>>,
в которой оптимизация впервые 
была рассмотрена с позиции функционального
анализа~\cite{Kant_lp40}. В данной работе был сформулирован
критерий оптимальноcти в задаче максимизации
функционала на компактном множестве.

За свой вклад в теорию~ЛП  и экономику
Канторович (вместе с независимо получившим
те же результаты, но позже,  Т.~Купмансом) был удостоен
Нобелевской премии по экономике
(с формулировкой <<за вклад в теорию оптимального распределения ресурсов>>) в 1975~году. 

Активное развитие теория~ЛП (а позже
и выпуклая оптимизация) получила в Соединённых
Штатах Америки в середине и конце 40-х годов
XX~века. Помимо Канторовича и Купманса,
значительный вклад в стремительное развитие
линейного программирования внесли математики
Джон фон Нейман и Джордж Бернард~Данциг,
а также известный экономист Василий Леонтьев.

Фон Нейман связал основную задачу теории матричных игр с нулевой суммой (теорему о минимаксе)  с линейным программированием. 
Принадлежащая фон Нейману модель расширяющейся экономики 
с возможностью определения максимального темпа роста экономической системы 
содержит элементы, которые до некоторой степени предвосхищают линейное программирование. 

Василий Леонтьев (1905--1999)
в 1930-е~гг. начал работу над
экономико-математической моделью <<затраты--выпуск>>
(модель многоотраслевой экономики).
В модели взаимодействие между различными
отраслями экономики
учитывается таким образом:
объём продукции каждой отрасли является
линейной комбинацией объёмов
продукции остальных отраслей (что
отражает необходимые затраты на обеспечение выпуска) с прибавленным
объёмом конечного потребления по этой отрасли.
В формальной постановке модель описывается
с помощью системы линейных алгебраических
уравнений. 
В~1973~году Василий Леонтьев был удостоен Нобелевской премии по экономике <<за развитие метода «затраты~-- выпуск>> и за его применение к важным экономическим проблемам>>.

Д.~Данциг (1914--2005)
в 1947~г. разработал первый вычислительный алгоритм
для решения задач~ЛП~--- симплекс-метод. Эта работа
выполнялась для американской армии.
Метод последовательно перебирает смежные
вершины графа допустимого множества, так, чтобы
в каждой последующей вершине значение
целевой функции или улучшалось,
или оставалось неизменным.
Результатом работы метода будет либо
множество оптимальных  решений~--- одна из вершин
допустимого полиэдра или его сторона, либо определение
несовместимости множества ограничений,
либо заключение о неограниченности целевой функции.
Благодаря простоте реализации и хорошей практической скорости работы, несмотря на полувековые попытки заменить его, симплекс-метод до сих пор остаётся одним из наиболее
популярных методов решения задач~ЛП на практике. Он
реализован во многих оптимизационных вычислительных пакетах.

В 1950-е годы продолжилось
развитие теоретических основ линейного программирования, его применение
для промышленных задач и создание
первых программных кодов для решения
задач~ЛП~\cite{Todd_LP}. В 1960-е годы стало возможным
решать задачи линейного программирования
большой размерности.
Симплекс-метод хорошо работал на практике:
эксперименты подтверждали,
что число требуемых итераций 
почти линейно зависит
от числа 
ограничений. 
Но в начале 1970-х годов авторитет
 симплекс-метода как универсального
 средства решения задач~ЛП сильно
 упал.
 Появилась теория вычислительной сложности,
 и на специальных примерах было показано, что симплекс-метод (при самом широко используемом правиле выбора
направляющего столбца и направляющей строки) 
требует экспоненциально большого (от размерности пространства) числа итераций\footnote{
Понятия экспоненциальной и полиномиальной разрешимости классов задач 
возникают в теории сложности алгоритмов. 
Если алгоритм решает любую задачу из
рассматриваемого класса задач за время, ограниченное
экспонентой от некоторого многочлена (полинома) от размерности
задачи, то говорят, что алгоритм имеет {\it экспоненциальную
сложность}.
{\it Полиномиальный} по времени алгоритм~--- это алгоритм
решения любой задачи из рассматриваемого класса, в котором число 
арифметических шагов (время работы) ограничено некоторым многочленом от размерности задачи  или, иными словами, полиномиально ограничено. 
Сам класс задач
в этом случае является {\it полиномиально разрешимым}.}.
Первыми такие примеры представили  американские математики В.~Кли и Д.Д.~Минти в 1972~г.~\cite{KleeMin} (см. пример~\ref{ch1_ex_KleeM}). Позже примеры задач~ЛП, на которых
симплекс-метод работает экспоненциально долго\footnote{Это были модифицированные
варианты примера~\ref{ch1_ex_KleeM}.}, были предложены
и для других правил выбора направляющей строки и направляющего столбца (см., например, \cite{Jer_spm}).

Тем не менее в упомянутых примерах речь идёт
о худшем случае работы алгоритма, а 
в целом на практике
симплекс-метод обычно работает
очень быстро, и многочисленные
эксперименты и исследования метода
подтверждали
полиномиальное время работы на практике. 
В 70-е и 80-е годы XX~века эта
разница между теоретическими
и практическими результатами
не давала покоя многим исследователям.
Было
несколько попыток объяснить загадку этих
наблюдений путём вероятностного анализа работы симплекс-метода
на случайных входных данных, которые моделировали
задачи из реальных приложений (см., например,
\cite{Borgwardt, Smale_SM, Todd_SM, Vershik_lp_avg}).
Были получены результаты, что {\bf в среднем} симплекс-метод
полиномиален, поскольку
требует $O(m^3)$ операций (при общем количестве итераций $2m$,
$m$~--- число строк в исходной матрице~$A$ задачи).
Но наиболее успешная попытка
дать объяснение
факту быстрой работы симплекс-метода
на практике была сделана в эпохальной работе начала
XXI~века Спилмана и Тенга~\cite{SpielmanTeng}:
необходимо проводить так называемый 
сглаженный анализ алгоритмов
(smoothed analysis), который
является промежуточным звеном между
анализом худшего случая и анализом
в среднем.  Аргументация такого промежуточного подхода
в следующем: при анализе работы алгоритма
в худшем случае задачи, на которых алгоритм
работает медленно, обычно подбираются
искусственно. \nk{Но они} не являются типичными
и почти никогда не возникают на практике,
при расчётах с реальными данными.
Что касается анализа работы алгоритма в среднем случае,
который был призван решить упомянутую проблему
анализа худшего случая, то и у этого подхода тоже есть свои недостатки.
В частности, при генерации входных данных эти данные
должны соответствовать какому-то распределению,
которое, в свою очередь, моделирует исходные
данные реальных задач. Но распределения данных
реальных задач сильно отличаются друг от друга.
Кроме того, не во всех случаях в принципе для
задачи можно определить распределение исходных
данных. В результате входные данные задач,
которые встречаются на практике, могут
почти не иметь ничего общего со входами,
сгенерированными для анализа работы алгоритма в среднем.

В качестве ответа на перечисленные недостатки анализа
работы алгоритмов Спилман и Тенг предложили альтернативный подход к моделированию
реальных данных. Для начала необходимо определить
основные свойства данных типичных практических
задач, которые позволяет решать алгоритм. Затем
разработать модель генерации таких входных данных.
А далее изучать работу алгоритма на этой модели,
учитывая тот факт, что реальные данные часто имеют
некоторую погрешность (случайный шум в версии Спилмана и~Тенга).

Таким образом, сглаженный анализ является
вероятностным анализом алгоритма: проверяется
работа алгоритма при незначительных случайных 
возмущениях (конкретнее~--- гауссовых возмущениях) конкретных входных данных,
при этом ищется зависимость производительности
алгоритма от размера входа и от 
среднеквадратичного отклонения 
возмущений.

Более формально, идея Спилмана и Тенга в следующем.
Рассматривается задача~ЛП вида 
\begin{equation}
\label{ch1_eq_lp_spil_teng}
\max_{Ax \, \ge \, b} c^Tx.
\end{equation}
Пусть матрица~$A$ представляется в виде 
$$
A = \bar{A} + \|\bar{A}\| G,
$$
где каждый элемент матрицы~$G$~--- нормально распределённая
независимая случайная величина с нулевым математическим ожиданием
и дисперсией~$\sigma^2$. Если в случае анализа алгоритма
в среднем $\sigma^2 \rightarrow \infty$, то для сглаженной
сложности $\sigma^2 \rightarrow 0$. Если
мы ищем время решения задачи~\eqref{ch1_eq_lp_spil_teng}
как функцию $T(A)$ от заданной матрицы~$A$ размерности $m \times n$, то для любого $\bar{A}$
значение $T(A)$ является случайной величиной, и
её математическое ожидание есть полином от $n$, $m$, $\log \frac{1}{\varepsilon}$, $\frac{1}{\sigma}$
в категориях сложности по Блюм--Шубу--Смейлу (см. п.~\ref{ch2_subsectBSSm}).

Спилман и Тенг показали, что симплекс-метод имеет полиномиальную
сглаженную сложность\footnote{В том же 2006~году
оценка симплекс-метода Спилмана и Тенга
была улучшена Р.~Вершининым~\cite{Versh}. Он показал,
что ожидаемое время работы симплекс-метода на незначительно
изменённых входных данных является полиномом от логарифма
количества ограничений~$n$.}. Эти результаты
означают, что хотя и существуют задачи,
на которых симплекс-метод будет работать
экспоненциально долго, но если исходные
данные таких задач (коэффициенты целевой
функции и ограничений) подвергнуть
незначительному изменению, то с достаточно
высокой вероятностью симплекс-метод
на возмущённой задаче уже будет работать
за полиномиальное время.

В конце 1970-х годов обнаружили,
что алгоритмы, в общем случае решающие задачи~ЛП за полиномиальное
время, всё-таки существуют (см.~подробнее
главу~\ref{chapt_num_met}, п.~\ref{ch2_subsect_ellips_meth}).
Все такие полиномиальные алгоритмы~--- 
методы внутренней точки (самый первый метод внутренней точки~--- метод И.\,И.~Дикина~\cite{Dikin}, первый метод с полиномиальной сложностью~--- метод Кармаркара \cite{Karmarkar}, см. также п.~\ref{subs:history_IP}) и метод эллипсоидов~--- отличались от симплекс-метода геометрическим подходом. Эти алгоритмы либо генерировали последовательность точек во внутренности допустимого множества, либо помещали допустимое множество в эллипсоид.
В течение 50~лет оставался открытым вопрос, 
существует ли полиномиальный алгоритм, который работает подобно симплекс-методу, перебирая только вершины (угловые точки) допустимого множества задачи.
Ответ на этот вопрос дали Келнер и Спилман в 2006~г.~\cite{Kelner2006}:
они представили рандомизированный симплекс-метод
с полиномиальным временем работы. Как
и другие известные полиномиальные алгоритмы
для решения задач~ЛП, время работы предложенного 
Келнером и Спилманом алгоритма полиномиально зависит
от числа битов, необходимых для представления входной информации
задачи. Данный алгоритм показывает полиномиальную скорость
работы и на <<плохих>> для симплекс-метода задачах.


\begin{figure}[htb]
\begin{center}
\begin{tikzpicture}[x  = {(0.9cm,-0.076cm)},
                    y  = {(-0.06cm,0.95cm)},
                    z  = {(-0.44cm,-0.29cm)},
                    scale = 2,
                    color = {green}]

  \definecolor{pointcolor_unnamed__1}{rgb}{ 0,0,0 }
  \tikzstyle{pointstyle_unnamed__1} = [fill=pointcolor_unnamed__1]

  \coordinate (v0_unnamed__1) at (2, 0.5, 1.875);
  \coordinate (v1_unnamed__1) at (2, 1.5, 1.625);
  \coordinate (v2_unnamed__1) at (2, 1.5, 0.375);
  \coordinate (v3_unnamed__1) at (2, 0.5, 0.125);
  \coordinate (v4_unnamed__1) at (0, 0, 0);
  \coordinate (v5_unnamed__1) at (0, 2, 0.5);
  \coordinate (v6_unnamed__1) at (0, 2, 1.5);
  \coordinate (v7_unnamed__1) at (0, 0, 2);

  \definecolor{edgecolor_unnamed__1}{rgb}{ 0,0,0 }

  \definecolor{facetcolor_unnamed__1}{rgb}{ 0.8, 0.8, 0.8}

  \tikzstyle{facestyle_unnamed__1} = [fill=facetcolor_unnamed__1, fill opacity=0.75, draw=edgecolor_unnamed__1, line width=1 pt, line cap=round, line join=round]

  \draw[facestyle_unnamed__1] (v0_unnamed__1) -- (v7_unnamed__1) -- (v4_unnamed__1) -- (v3_unnamed__1) -- (v0_unnamed__1) -- cycle;
  \draw[facestyle_unnamed__1] (v4_unnamed__1) -- (v7_unnamed__1) -- (v6_unnamed__1) -- (v5_unnamed__1) -- (v4_unnamed__1) -- cycle;
  \draw[facestyle_unnamed__1] (v4_unnamed__1) -- (v5_unnamed__1) -- (v2_unnamed__1) -- (v3_unnamed__1) -- (v4_unnamed__1) -- cycle;

  \fill[pointcolor_unnamed__1] (v4_unnamed__1) circle (1 pt);

  \draw[facestyle_unnamed__1] (v5_unnamed__1) -- (v6_unnamed__1) -- (v1_unnamed__1) -- (v2_unnamed__1) -- (v5_unnamed__1) -- cycle;

  \fill[pointcolor_unnamed__1] (v5_unnamed__1) circle (1 pt);

  \draw[facestyle_unnamed__1] (v1_unnamed__1) -- (v0_unnamed__1) -- (v3_unnamed__1) -- (v2_unnamed__1) -- (v1_unnamed__1) -- cycle;

  \fill[pointcolor_unnamed__1] (v3_unnamed__1) circle (1 pt);
  \fill[pointcolor_unnamed__1] (v2_unnamed__1) circle (1 pt);

  \draw[facestyle_unnamed__1] (v6_unnamed__1) -- (v7_unnamed__1) -- (v0_unnamed__1) -- (v1_unnamed__1) -- (v6_unnamed__1) -- cycle;

  \fill[pointcolor_unnamed__1] (v6_unnamed__1) circle (1 pt);
  \fill[pointcolor_unnamed__1] (v7_unnamed__1) circle (1 pt);
  \fill[pointcolor_unnamed__1] (v0_unnamed__1) circle (1 pt);
  \fill[pointcolor_unnamed__1] (v1_unnamed__1) circle (1 pt);


\end{tikzpicture}
\end{center}
\caption{Пример куба Кли--Минти при $n = 3$. Рисунок сделан с помощью программы для
проведения исследований в области
полиэдральной геометрии~{\bf polymake}~\cite{Polymake, Polymake_art}} \label{ch1_dr_kleeM}
\end{figure}

 \begin{example}[пример Кли и Минти]\label{ch1_ex_KleeM}
 Рассмотрим задачу 
\begin{equation}
\label{ch1_eq_cubeKM}
  \begin{array}{cl}
  \max  & 2^{n-1} x_1 + 2^{n-2} x_2 + \ldots + 2x_{n-1} + x_n,  \\
   x_1 \, \le \,  5  &  \\
   4x_1 + x_2 \, \le \,  25 &  \\
  8x_1 + 4x_2 + x_3 \, \le \,  125  &  \\
  \ldots &  \\
  2^n x_1 + 2^{n-1}x_2 + 2^{n-2}x_3 + \ldots + & 4 x_{n-1} + x_n \, \le \, 5^n   \\
  x \, \ge \,  0 &  
   \end{array}
\end{equation}
 где $x \, \in \, \mathbb{R}^n$.
Начиная с точки $x_0 = 0$ и следуя симплекс-методу Данцига, необходимо обойти  ${2^n-1}$~вершину, следовательно, алгоритму
требуется экспоненциально большое
число итераций. 
Например, для $n = 71$\linebreak
$2^{n-1} \approx 1.18 \cdot 10^{21}$.
Если обходить 1000~вершин
за 1~секунду, для решения потребуется
примерно 37~миллиардов лет (для сравнения~--- возраст Вселенной 13,8~миллиардов лет).

Закончит работу алгоритм в оптимальной точке~$(0, \, 0, \, \ldots, \, 5^{n})$. Пример получен
в результате небольшой деформации $n$-мерного гиперкуба\protect\footnotemark 
$\ \,0  \le  x_i  \le  1$,
$i = 1, \, \ldots, \, n$, заставляющей симплекс-метод перебирать
все вершины куба.

\begin{remark}
В общем виде примеры Кли и Минти (см. также рис.~\ref{ch1_dr_kleeM}), для которых
симплекс-методу требуется экспоненциальное
число итераций, можно описать таким образом~\cite{KM_exPSZ}:
$$
  \begin{array}{cl}
  \max  & \sum_{j = 1}^n \varepsilon^{n-j} x_j,  \\
   x_1 \, \le \, 1  &  \\
  2 \sum_{j = 1}^{i-1} \varepsilon^{i-j} x_j + x_i \, \le \, 1, \, i = 2, \, 3, \, \ldots, \, n&  \\
  x \, \ge \, 0 &  
   \end{array}
$$
где $0 < \varepsilon \, \le \, 1/3$, $x \, \in \, \mathbb{R}^n$.
\end{remark}
 \end{example}
\footnotetext{Сейчас допустимое множество системы линейных
неравенств из~\eqref{ch1_eq_cubeKM},
предложенной Кли и Минти,
называют кубом, или политопом Кли--Минти.}

\begin{remark}\ 
О главной нерешённой  
проблеме в теории линейного программирования, сформулированной
Стивеном Смейлом, см.~в~п.~\ref{ch2_subsectBSSm}. 
\end{remark}


Существуют два частных случая задачи~ЛП~\eqref{ch1_eq_lp_all},
которые так часто встречаются
в приложениях, что имеют специальные названия.
Один из них~--- {\it задача~ЛП в стандартной форме}:
\begin{eqnarray}
\label{ch1_eq_lp}
& \min & c^T x, \\  
& Ax = b &\nonumber \\
& x \, \ge \, 0 & \nonumber
\end{eqnarray}
т.е. задача~ЛП, в которой ограничениями-неравенствами являются только
условия неотрицательности всех компонент вектора~$x$.
Для приведения~\eqref{ch1_eq_lp} к каноническому
виду~\eqref{ch1_lagr_all_problem} необходимо
учитывать, что $h_i(x) = -x_i$, $i = 1, \, \ldots, \, p$.

Вторым часто встречающимся случаем задачи~ЛП
является задача без ограничений-равенств вида
\begin{eqnarray}
& \min & c^T x. \\  
& A x \, \le \, b & \nonumber
\end{eqnarray}

Перейдём к аналитическому рассмотрению задачи~\eqref{ch1_eq_lp}.
Функция Лагранжа для этой задачи выглядит
следующим образом:
$$
L(x, \, \lambda, \, \nu) = c^T x + \lambda^T (Ax - b) -
\sum_{i = 1}^p \nu_i x_i = -b^T \lambda + \left ( c + A^T \lambda - \nu \right )^T x.
$$

Двойственная функция, соответственно,
\begin{equation*}
\varphi(\lambda, \, \nu) = \inf_x L(x, \, \lambda, \, \nu) = 
-b^T \lambda + \inf_x \left \{ \left ( c + A^T \lambda - \nu \right )^T x \right \}.
\end{equation*}

Задача на минимум в этом выражении легко решается аналитически:

$$
\varphi(\lambda, \, \nu) = \left \{
\begin{array}{cc}
        -b^T \lambda  & \text{при } \left ( c + A^T \lambda - \nu \right ) = 0, \\
-\infty     & \text{иначе}. 
\end{array}
\right.
$$
Таким образом, нижняя оценка оптимального значения $f^*$
является нетривиальной только в случае, когда
$\nu_i \, \ge \, 0$, $i = 1, \, \ldots, \, p$
и $\left ( c + A^T \lambda - \nu \right ) = 0$.
Двойственная к задаче линейного программирования 
записывается так:
\begin{eqnarray}
\label{ch1_eq_lp_dual}
& \max & -b^T \lambda. \\  
&  c + A^T \lambda - \nu  = 0 & \nonumber \\
& \nu \, \ge \, 0 & \nonumber
\end{eqnarray}

\begin{teo}[слабая двойственность]
Пусть $x^*$~--- допустимое решение задачи~\eqref{ch1_eq_lp}, $\lambda^*$~--- допустимое решение задачи~\eqref{ch1_eq_lp_dual}.
Тогда $c^T x^* \geq -b^T \lambda^*$.
\end{teo}

\begin{teo}[сильная двойственность, основная теорема линейного программирования]
Двойственная задача~\eqref{ch1_eq_lp_dual} имеет оптимальное решение тогда и
только тогда, когда оптимальное решение
имеет задача~\eqref{ch1_eq_lp}.
Пусть $x^*$~--- конечное оптимальное решение задачи~\eqref{ch1_eq_lp}, $\lambda^*$~--- конечное оптимальное решение задачи~\eqref{ch1_eq_lp_dual}.
Тогда $c^T x^* = -b^T \lambda^*$.
\end{teo}

Доказательства последних двух теорем см., например, в \cite{AMP_77}.

Обратите внимание: переход к двойственной задаче в случае,
когда число ограничений значительно меньше числа переменных
в прямой задаче, может сильно уменьшить размерность задачи, которую требуется решать.

В современных оптимизационных пакетах задачи ЛП средней размерности (вплоть до сотен тысяч и даже миллиона переменных или ограничений) решаются методами внутренней точки. В частности, методами, описанными в п.~\ref{SB}. Задачи ЛП сверхбольшой размерности требуют разработки более специальных методов, лучше учитывающих структуру задачи (см., например, п.~\ref{truss}).

\subsection{Оценка на множители Лагранжа}
Выполнение условий Слейтера для задачи выпуклого программирования позволяет 
получить верхнюю оценку оптимальных множителей Лагранжа $\nu^* \, \in \, \mathbb{R}_+^p$.
Последняя может быть полезна при использовании методов штрафных функций \cite{Nocedal} или двойственного подхода к решению исходной задачи (см. пример в п.~\ref{comp_prox}).

 
Пусть для  задачи~\eqref{ch1_eq_min_conv_in} выполняются условия
регулярности Слейтера:
$$
\exists \, \, x' \, \in \, Q, \, \, 
i \, \in \, \{1, \, \dots, \, p\}: \quad h_i\left( x' \right) \, < \, 0. 
$$
И пусть число $\gamma$ таково, что
$\gamma = \min\limits_{i=1, \, \ldots, \, p} \left\{ {-h_i \left( x' \right)} \right\}$.

По 
определению двойственной функции
$$
\varphi \left( { {\nu^* }} \right) = \inf_{x \, \in \, Q} \left\{ {f\left( x 
\right)+\sum\limits_{i=1}^p {\nu^*_i h_i \left( x \right)} } \right\} \leq 
f\left( {x'} \right) + \sum_{i=1}^p {\nu^*_i h_i \left( x' \right)} .
$$
Тогда
$$
\gamma \sum_{i=1}^p \nu^*_i \, \le \, -\sum_{i=1}^p \nu^*_i h_i \left( x' \right) \, \le \, f(x') - \varphi({\nu^*}) 
$$
и 
$$
 \| \nu^* \|_1 =  \sum_{i=1}^p | \nu^*_i | 
= 
 \sum_{i=1}^p \nu^*_i \, \le \, \frac{1}{\gamma} \left (f(x') - \varphi({\nu^*}) \right ).
$$

Из последнего неравенства получаем оценку
на оптимальный множитель:
$$
\| \nu^* \|_1 \, \le \,  
\frac{1}{\gamma} \left (f(x') - \varphi({\nu^*}) \right ).
$$

\subsection{Геометрическая интерпретация}\label{geom}
Рассматривается задача с ограничениями типа равенств 
и неравенств вида~\eqref{ch1_lagr_all_problem}.
Область определения $\mathbb{D}$ задачи определяется
в~\eqref{ch1_eq1_dom_all_prob}.

Определим множество значений целевой функции и ограничений задачи:
$$
G = \left \{ \left. \left (f(x), \, g_1(x), \ldots,  g_m(x), \, h_1(x),  \ldots, h_p(x) \right )
 \in \mathbb{R} \times \mathbb{R}^m \times \mathbb{R}^p   \, \right | 
\, x \in \mathbb{D}
\right\}.
$$
Тогда оптимум задачи~\eqref{ch1_lagr_all_problem}
$f^* = f(x^*)$ можно выразить следующим образом:
$$
f^* = \inf \left \{ t \,
\right |  \, 
(t, u, v) \, \in \, G, \left. \, v \, \le \, 0, \, u = 0 \right \}.
$$

Двойственную функцию $\varphi(\lambda, \, \nu)$
можно вычислить как 
\begin{eqnarray*}
\varphi(\lambda, \, \nu) &=& \inf \left \{
\left.
(1, \lambda, \nu)^T (t, u, v) \, \right | \, (t, u, v) \, \in \, G
\right \} = \\
&=& \inf \left \{
t + \sum_{i = 1}^m \lambda_i u_i  
 + \left. \sum_{i = 1}^p \nu_i v_i  
 \, \right |  \, 
 (t, u, v) \, \in \, G
\right \}.
\end{eqnarray*}
В случае, если точная нижняя грань конечна, неравенство
\begin{equation}
\label{ch1_eq_sup_hyp}
(1, \lambda, \nu)^T (t, u, v) \, \ge \, \varphi(\lambda, \, \nu)
\end{equation}
определяет так называемую {\it опорную гиперплоскость} к $G$. 
\begin{defin}\label{ch1_def_sup_hyp}
Гиперплоскость $H = \left \{ x \, | \, a^T x = \alpha \right \}$ называется опорной к множеству $G$ в точке $x^0 \in G$,
если $x^0 \, \in \, H$ и всё множество~$G$ лежит
в полупространстве, задаваемом $H$, т.е. $\forall \, x \, \in \, G \, \rightarrow \, a^T x \, \le \, \alpha$.
\end{defin}

\begin{figure}[htb]
\begin{center}
\includegraphics[width=0.7\linewidth]{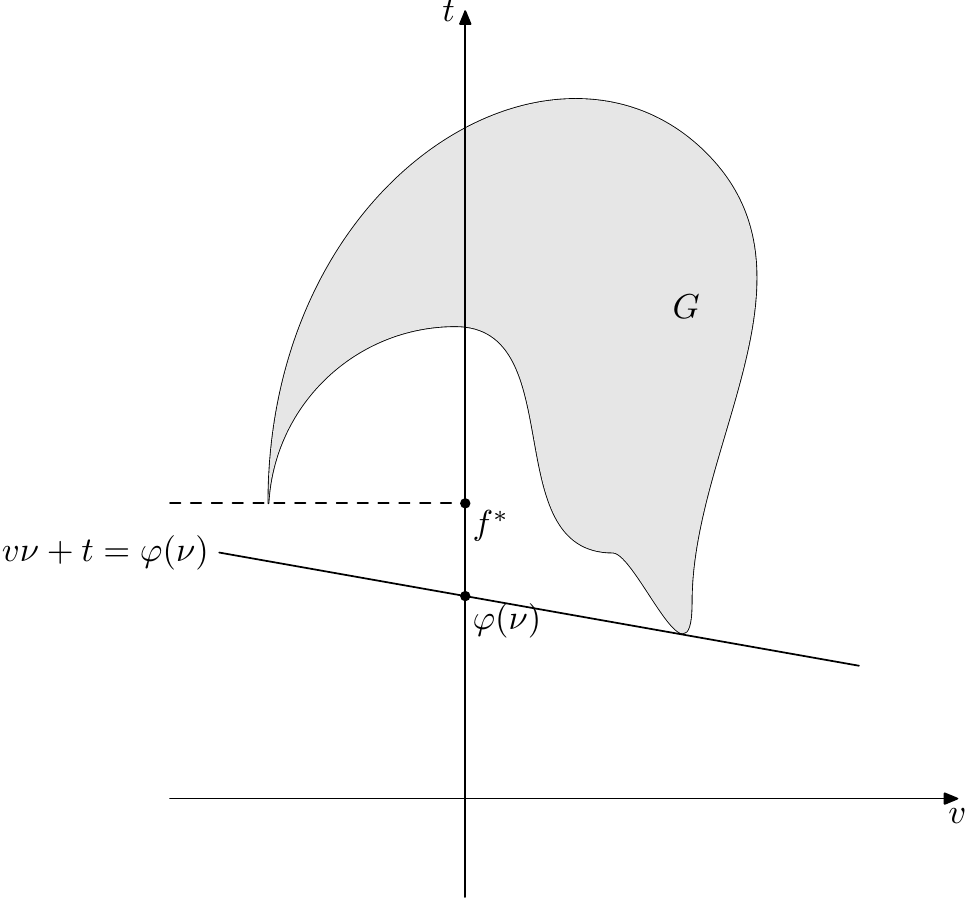}
\end{center}
\caption{Геометрическая иллюстрация
факта выполнения неравенства ${\varphi(\nu)\le  f^*}$. Случай задачи с одним
ограничением-неравенством: $\min\limits_{h(x)  \le  0} f(x)$, $G\!\!\!= \linebreak =\{ (f(x),  h(x)) \, | \, x  \in  \mathbb{D}$.
Уравнение опорной гиперплоскости
(в данном двумерном случае это прямая): $(1, \, \nu)^T (t, \, v) = 
v \nu + t = \varphi(\nu)$ 
(\cite{Boyd_book_2004}, с изменениями обозначений и~формы~$G$)} \label{ch1_dr_5.3}
\end{figure}

Опорную гиперплоскость~\eqref{ch1_eq_sup_hyp} иногда
называют {\it невертикальной}
(см.~рис.~\ref{ch1_dr_5.3}).


При $\nu \, \ge \, 0$ в общем случае в задаче~\eqref{ch1_lagr_all_problem} имеется
ненулевой зазор двойственности (выполняется только свойство
слабой двойственности).
В~самом деле, $t \, \ge \, (1, \lambda, \nu)^T (t, u, v)$
при $v \, \le \, 0$ и $u = 0$. Тогда
\begin{eqnarray*}
f^* & = & \inf \left \{ t \, | \, (t, u, v) \, \in \, G, \, v \, \le \, 0, \, u = 0  \right \} \\
& \ge &  \inf \left \{
(1, \lambda, \nu)^T (t, u, v) \, | \, (t, u, v) \, \in \, G, \, v \, \le \, 0, \, u = 0 
\right \} \\
& \ge & \inf \left \{
(1, \lambda, \nu)^T (t, u, v) \, | \, (t, u, v) \, \in \, G \right \} \\
& = & \varphi(\lambda, \, \nu),
\end{eqnarray*}
т.е. зазор двойственности в общем случае ненулевой (см. рис.~\ref{ch1_dr_5.3}).


\begin{figure}[htb]
\begin{center}
\includegraphics[width=0.7\linewidth]{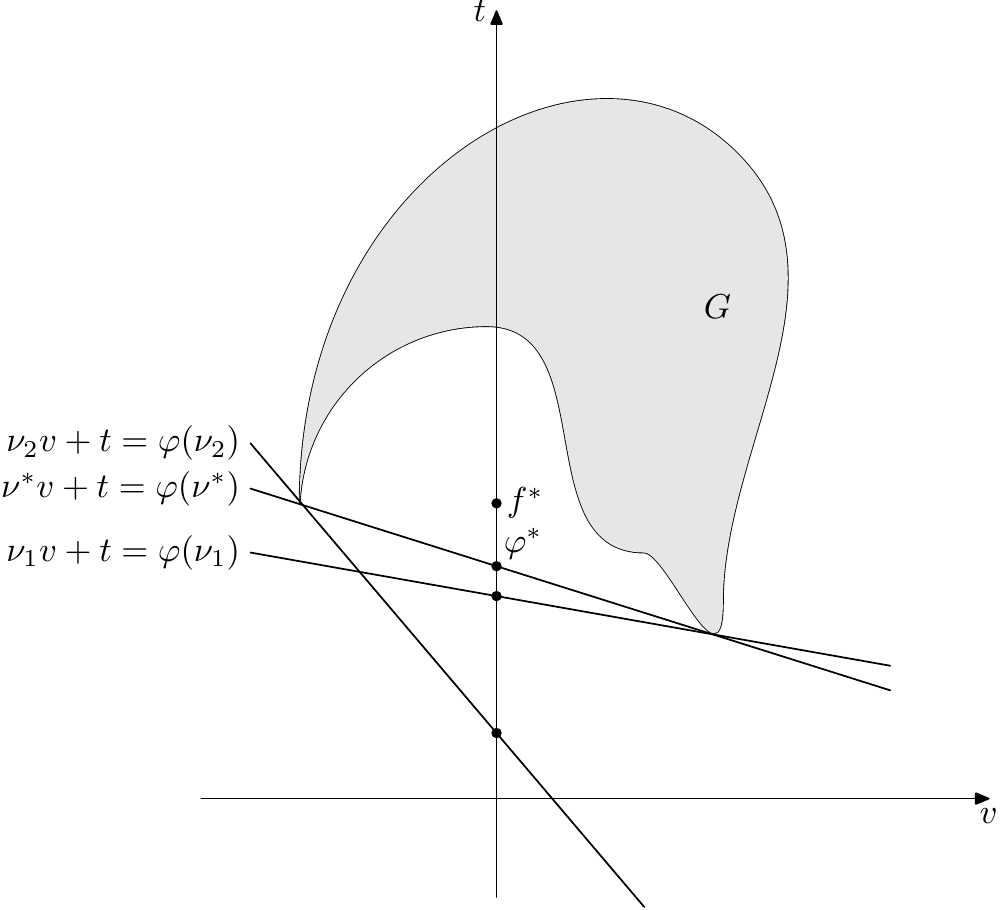}
\end{center}
\caption{Три опорные гиперплоскости, соответствующие трём значениям
двойственной функции.
Случай задачи с одним
ограничением-неравенством: $\min\limits_{h(x) \, \le \, 0} f(x)$, $G = \{ (f(x), \, h(x)) \, | \, x \, \in \, \mathbb{D}\}$.
Для данной задачи не выполняется свойство сильной двойственности: зазор
двойственности $f^* - \varphi^*$ положителен  
(\cite{Boyd_book_2004}, с изменениями обозначений и формы~$G$)} \label{ch1_dr_5.4}
\end{figure}

\subsubsection{Выпуклый случай}
В этом разделе, следуя~\cite{Boyd_book_2004}, 
покажем, что свойство сильной двойственности
часто имеет место для задач выпуклого программирования.
Определим множество $\mathbb{A} \, \subseteq \, \mathbb{R} \times \mathbb{R}^m \times \mathbb{R}^p$
как

\begin{eqnarray}
\mathbb{A} = \left \{
(t, u, v) \, | \, \exists \, x  \in  \mathbb{D}, \,
h_i(x) \le v_i, \, i = 1, \ldots, p, \,
g_i(x) = u_i, \, i = 1, \ldots,  m, \right. \nonumber  \\
 f(x) \le  t \}. \qquad \qquad
 \label{ch1_eq_a}
\end{eqnarray}

{
\renewcommand{\baselinestretch}{0.96}
\selectfont

\noindent Если выразить оптимум через множество $\mathbb{A}$,
то получим
$$
f^* = \inf \left \{
t \, | \, (t, 0, 0) \, \in \, \mathbb{A}
\right \}.
$$
При $\nu \, \ge \, 0$ двойственная функция
\begin{equation}
\label{ch1_eq_dual_conv_geom}
\varphi(\lambda, \, \nu) = \inf \left \{
(1, \lambda, \nu)^T (t, u, v) \, | \, (t, u, v) \, \in \, \mathbb{A}
\right \}. 
\end{equation}

Поскольку $(f^*, 0, 0)$ лежит на границе множества~$\mathbb{A}$, имеем
\begin{equation}
\label{ch1_eq_dual_conv_weak}
f^* = (1, \lambda, \nu)^T (f^*, 0, 0) \, \ge \, \varphi(\lambda, \, \nu)
\end{equation}
(при условии, что точная нижняя грань в~\eqref{ch1_eq_dual_conv_geom} конечна).
В неравенстве~\eqref{ch1_eq_dual_conv_weak}
свойство сильной двойственности будет
выполняться только в том случае, если существует
невертикальная опорная гиперплоскость\footnote{См. определение~\ref{ch1_def_sup_hyp}.} к $\mathbb{A}$
в точке $(f^*, 0, 0)$.

Далее докажем, что выполнение условий Слейтера гарантирует
наличие свойства сильной двойственности у задачи выпуклого программирования.

Рассмотрим выпуклую задачу~\eqref{ch1_eq_conv_sl_prob}.
Пусть в точке $x' \, \in \, \text{int } \mathbb{D} \,$ $h_i(x') \, < \, 0$, $i = 1, \, \ldots, \, p$,
и $Ax' = b$. Здесь область определения задачи $\mathbb{D} =  \textup{dom } f \, \cap \,  \bigcap\limits_{i = 1}^p \,\textup{dom } h_i(x)$, $\text{int } \mathbb{D}$~--- внутренность
множества~$\mathbb{D}$, и для упрощения доказательства сделаны следующие допущения: множество~$\mathbb{D}$
имеет непустую внутренность и ранг матрицы~$A$ (обозначается~$\text{rank } A$) равен $m$, а также допущение, что оптимум~$f^*$ конечен.

Определённое в~\eqref{ch1_eq_a} множество~$\mathbb{A}$
выпукло для выпуклой задачи~\eqref{ch1_eq_conv_sl_prob}.
Добавим ещё одно выпуклое множество~$\mathbb{B}$:
$$
\mathbb{B} = \left \{
(s, \, 0, \, 0) \, \in \, \mathbb{R} \times \mathbb{R}^m 
\times \mathbb{R}^p \, | \, s < f^*
\right \}.
$$

Множества~$\mathbb{A}$ и $\mathbb{B}$ не пересекаются. В самом деле, предположим противное. Пусть существует точка~$(t, u, v)$, принадлежащая $\mathbb{A} \cap \mathbb{B}$. Из того, что $(t, u, v) \, \in \, \mathbb{B}$,
следует, что $u = 0$, $v = 0$, $t < f^*$. Также
$(t, u, v) \, \in \, \mathbb{A}$, поэтому
$\exists \, \tilde{x} \, \in \, \mathbb{D}$, такой, что 
$h_i(\tilde{x}) \, \le \, 0$, $i = 1, \, \ldots, \, p$,
$A \tilde{x} - b = 0$, $f(\tilde{x}) \, \le \, t < f^*$.
Но это противоречит тому, что $f^*$ является оптимальным
значением в задаче~\eqref{ch1_eq_conv_sl_prob}. \qed

Итак, множества $\mathbb{A}$ и $\mathbb{B}$ выпуклы и не пересекаются. Тогда по первой теореме отделимости~\ref{ch1_teo_otdel1} существуют $(\mu, \tilde{\lambda}, \tilde{\nu}) \neq 0$ 
и число $\alpha$, для которых
\begin{equation}
\label{ch1_eq_a_otdel}
(\mu, \tilde{\lambda}, \tilde{\nu})^T (t, u, v) \, \ge \, \alpha
\end{equation}
для всех $(t, u, v) \, \in \, \mathbb{A}$
и
\begin{equation}
\label{ch1_eq_b_otdel}
(\mu, \tilde{\lambda}, \tilde{\nu})^T (t, u, v) \, \le \, \alpha
\end{equation}
для всех $(t, u, v) \, \in \, \mathbb{B}$.

}

Из~\eqref{ch1_eq_a_otdel} следует, что $\mu \, \ge \, 0$ и $\tilde{\nu} \, \ge \, 0$.
Неравенство~\eqref{ch1_eq_b_otdel} сводится к $\mu t \, \le \, \alpha$ для всех $t < f^*$. Отсюда $\mu f^* \, \le \, \alpha$.
Итак, для любого $x \, \in \, \mathbb{D}$
\begin{equation}
\label{ch1_eq_slat_proof_d}
\sum_{i = 1}^p \tilde{\nu}_i h_i(x) + \tilde{\lambda}^T (Ax - b) + \mu f(x) \, \ge \, \alpha \, \ge \, \mu f^*.
\end{equation}

Рассмотрим два случая --- когда $\mu > 0$ и когда $\mu = 0$.

1. Пусть $\mu > 0$. Разделим обе части неравенства~\eqref{ch1_eq_slat_proof_d} на $\mu$:
\begin{equation}
\label{ch1_eq_slat_proof_l}
L(x, \, \tilde{\lambda}/\mu, \, \tilde{\nu}/\mu) \,
\ge \,
f^*
\end{equation}
для всех $x \, \in \, \mathbb{D}$.
Пусть $\lambda = \tilde{\lambda}/\mu$, $\nu = \tilde{\nu}/\mu$.
Из~\eqref{ch1_eq_slat_proof_l} следует, что
$\inf\limits_x L(x, \, \lambda, \, \nu) \,
\ge \, 
f^*$, т.е. $\varphi(\lambda, \, \nu) \, \ge \, f^*$.
Но из условия слабой двойственности
$\varphi(\lambda, \, \nu) \, \le \, f^*$,
поэтому $\varphi(\lambda, \, \nu) = f^*$,
и свойство сильной двойственности имеет место.


\begin{figure}[htb]
\begin{center}
\includegraphics[width=0.6\linewidth]{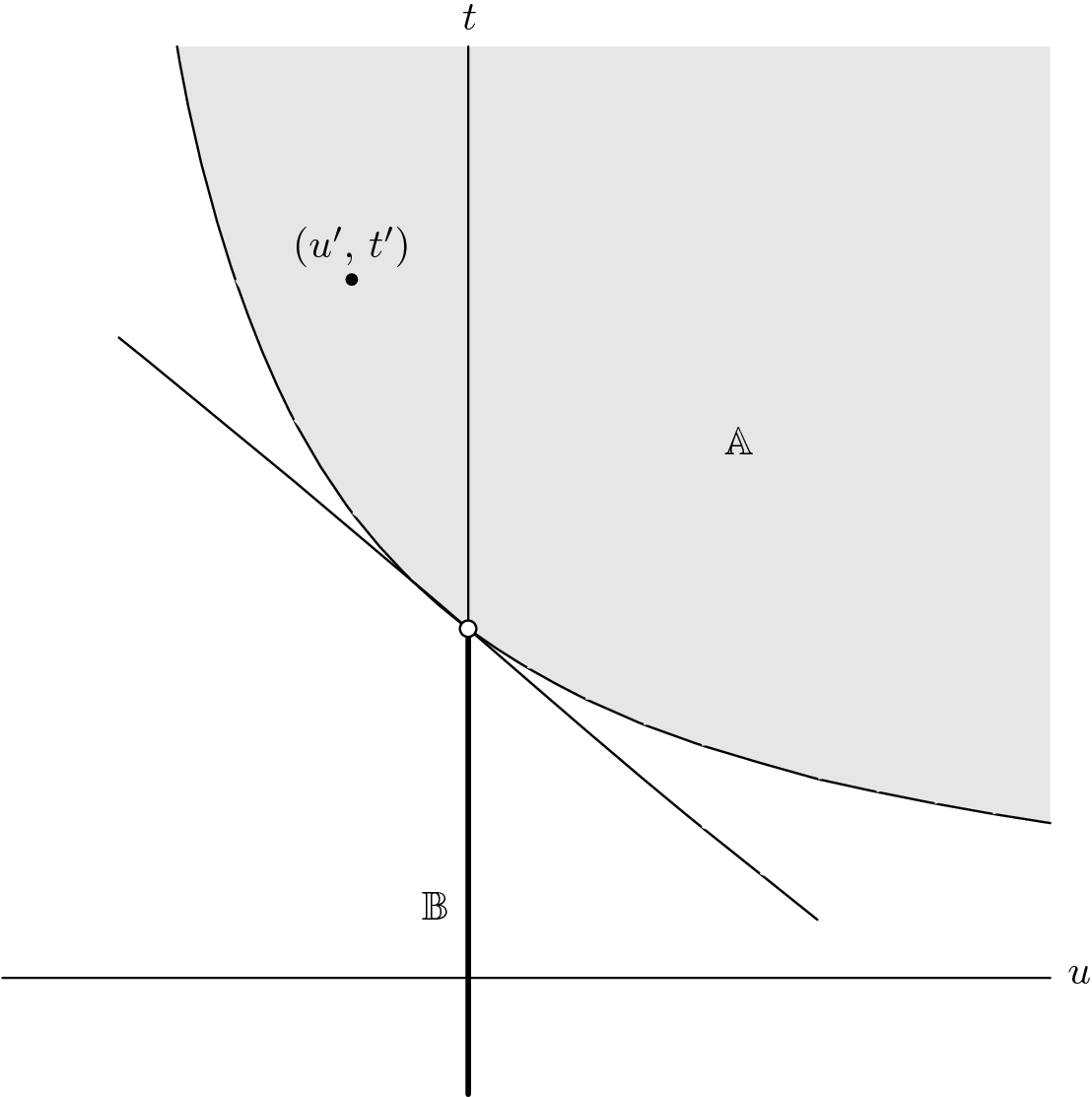}
\end{center}
\caption{Иллюстрация доказательства свойства
сильной двойственности для задачи
выпуклого программирования с одним ограничением (\cite{Boyd_book_2004}, с изменениями обозначений)} \label{ch1_dr_slater5.6}
\end{figure}

2. Пусть $\mu = 0$. Из~\eqref{ch1_eq_slat_proof_d}
следует, что
\begin{equation}
\label{ch1_eq_slat_proof_mu0}
\sum_{i = 1}^p \tilde{\nu}_i h_i(x) + \tilde{\lambda}^T (Ax - b)  \, \ge \, 0.
\end{equation}

Неравенство~\eqref{ch1_eq_slat_proof_mu0} при $x = x'$
принимает вид
\begin{equation}
\label{ch1_eq_slat_proof_vf}
\sum_{i = 1}^p \tilde{\nu}_i h_i(x')  \, \ge \, 0.    
\end{equation}
Поскольку
$h_i(x') < 0$, а $\tilde{\nu_i} \, \ge \, 0$,
из~\eqref{ch1_eq_slat_proof_vf} следует, что~$\tilde{\nu}=0$. Так как $(\mu, \tilde{\lambda}, \tilde{\nu})$ не должно быть равно~$0$, а~$\tilde{\nu}=0$,
$\mu = 0$, остаётся только $\tilde{\lambda}\neq 0$.
Тогда из~\eqref{ch1_eq_slat_proof_mu0} следует, что
для всех $x \, \in \, \mathbb{D}$ выполняется
$\tilde{\lambda}^T (Ax - b)  \, \ge \, 0$.
Но для~$x'$ $\tilde{\lambda}^T (Ax' - b)  = 0$,
и, поскольку $x'$ принадлежит внутренности множества~$\mathbb{D}$, 
должно быть $A^T \tilde{\lambda} = 0$.
Но последнее равенство противоречит допущению 
$\text{rank } A = m$. Противоречие доказывает, что $\mu$
не может быть равно $0$.

Данное доказательство проиллюстрировано
на рис.~\ref{ch1_dr_slater5.6} для случая с 
одним ограничением~($m = 0$, $p = 1$). 
Отделяющая множества~$\mathbb{A}$ и $\mathbb{B}$
гиперплоскость с нормалью $(\nu^*, \, 1)$ является невертикальной и опорной к~$\mathbb{A}$
в точке $(0, f^*)$. Условие Слейтера для точки~$x'$ гарантирует,
что любая отделяющая плоскость будет невертикальной,
поскольку она проходит слева от точки~$(u', \, t') = (h_1(x'), \, f(x'))$.

Итак, условие Слейтера
гарантирует выполнение свойства сильной двойственности
для задачи выпуклого программирования, и доказать данный факт можно с помощью теоремы отделимости.

\paragraph{Возмущения задачи.} 
Рассмотрим задачу с ограни\-че\-ниями-не\-равенст\-вами:
\begin{eqnarray}
\label{ch1_eq_pertub_problem}
& \min & f(x) \\  
& h_i(x) \, \le \, u_i, \, i = 1, \, \ldots, \, p, & \nonumber
\end{eqnarray}
где $ x \, \in \, \mathbb{R}^n$. 
Здесь изменение $u_i$ даёт возможность
ослабить или сделать более жёстким
ограничение $h_i(x) \, \le \, u_i$.
Определим функцию\footnote{См. \eqref{ch1_eq1_dom_all_prob}.} $f^*(u)$:
$$
f^*(u) = \inf \, \left \{ f(x) \, | \, \exists \, x \, \in \, \mathbb{D},
\, 
h_i(x) \, \le \, u_i, \, i = 1, \, \ldots, \, p \right \}.
$$
Ясно, что $f^*(0)$ равно оптимальному
значению задачи
\begin{eqnarray}
\label{ch1_eq_unpertub_problem}
& \min & f(x) \\  
& h_i(x) \, \le \, 0, \, i = 1, \, \ldots, \, p. & \nonumber
\end{eqnarray}

\begin{figure}[htb]
\begin{center}
\includegraphics[width=0.7\linewidth]{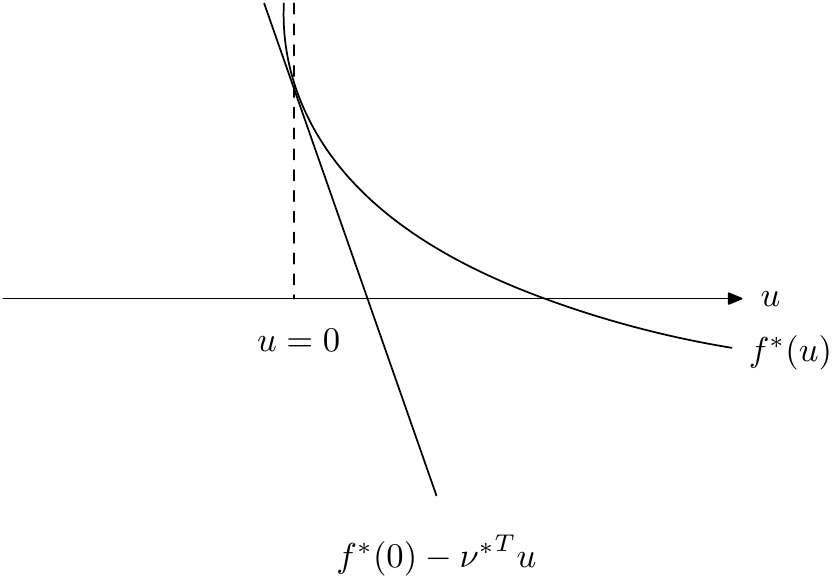}
\end{center}
\caption{Оптимальное значение $f^*(u)$ выпуклой
задачи с одним ограничением $h_1(x) \ge  u$. Функция~$f^*(0) - {\nu^*}^T u$ является
нижней границей $f^*(u)$
(\cite{Boyd_book_2004}, с изменениями обозначений)} \label{ch1_dr_pertub}
\end{figure}

Далее считаем, что задача~\eqref{ch1_eq_pertub_problem} выпукла
и выполняются условия Слейтера.
Пусть $\nu^*$~--- решение двойственной
задачи
$$
\max_{\nu \, \ge \, 0} \varphi(\nu)
$$
к \eqref{ch1_eq_unpertub_problem}.
В силу выполнения свойства сильной двойственности для любого $x$ из
допустимой области задачи~\eqref{ch1_eq_pertub_problem}
(т.е. $h_i(x) \, \le \, u_i$, $i = 1, \, \ldots, \, p$):
$$
f^*(0) = \varphi(\nu^*) \, \le \, f(x) + 
\sum_{i = 1}^p \nu_i^* h_i(x)
\, \le \,
f(x) + {\nu_i^*}^T u.
$$
Следовательно,
$$
f(x) \, \ge \, f^*(0) - {\nu_i^*}^T u.
$$
Тогда для всех $u$ имеет место {\it нижняя оценка
возмущённого оптимума}:
\begin{equation}
\label{ch1_eq_pertub_prop}
f^*(u) \, \ge \, f^*(0) - {\nu^*}^T u.
\end{equation}
Графическая иллюстрация свойства~\eqref{ch1_eq_pertub_prop}
дана на рис.~\ref{ch1_dr_pertub}.


С помощью неравенства~\eqref{ch1_eq_pertub_prop} можно сделать следующие заключения
о чувствительности оптимума по отношению
к изменению ограничений:
\begin{itemize}
    \item Если мы ужесточаем $i$-е ограничение (т.е. выбираем $u_i < 0$) и множитель Лагранжа $\nu_i^*$ имеет достаточно большое значение, то $f^*(u)$ сильно увеличится.
    \item Если мы ослабляем $i$-е ограничение (т.е. выбираем $u_i > 0$) и множитель Лагранжа $\nu_i^*$ имеет достаточно маленькое значение, то $f^*(u)$ не будет сильно уменьшаться.
\end{itemize}

\begin{exercise}
Напишите, как изменится нижняя оценка
и анализ чувствительности (неравенство~\eqref{ch1_eq_pertub_prop})
для возмущённой задачи с ограничениями-неравенствами
и равенствами.
\end{exercise}

\subsection{Условия дополняющей нежёсткости}
Пусть $x^*$~--- оптимальное решение прямой задачи с ограничениями типа равенств и неравенств~\eqref{ch1_lagr_all_problem},
и $(\lambda^*, \, \nu^*)$~--- двойственные множители. 
Предполагаем, что имеет место сильная двойственность.
Тогда
\begin{align*}
f^* &\!=\! f(x^*) \!=\! \varphi^* \!=\! \varphi(\lambda^*, 
\nu^*) \!=\! \inf_x L(x, \lambda^*,  \nu^*)  \le 
 f(x^*) \!+\! \sum_{i = 1}^m \lambda_i^* g_i (x^*) \!+\! \sum_{i = 1}^p \nu_i^* h_i (x^*) \\
&\!\le\!  f(x^*).
\end{align*}
Последнее неравенство выполняется, так как 
$g_i(x^*) = 0, \, i = 1, \, \ldots, \, m$, 
$h_i(x^*) \, \le \, 0, \, i = 1, \, \ldots, \, p$.

Из данных выкладок следует важное заключение:
$$
\sum_{i = 1}^p \nu_i^* h_i (x^*) = 0.
$$
В силу неположительности каждого элемента данной
суммы
\begin{equation}
\label{ch1_eq_com_slack}
\nu_i^* h_i (x^*) = 0, \, i = 1, \, \ldots, \, p.
\end{equation}
Как указано в формулировке теоремы~\ref{ch1_teo_kkt},
равенства~\eqref{ch1_eq_com_slack} называют 
{\it условиями дополняющей нежёсткости}.
По сути, эти условия означают, что двойственный
множитель Лагранжа~$\nu_i$ равен~$0$ при 
неактивном ограничении~$h_i(x^*) < 0$, и
если~$\nu_i$ положителен, то~$h_i(x^*) = 0$.  

\subsection{Седловые точки}\label{minmax}
Рассмотрим задачу с ограничениями-неравенствами следующего вида:
\begin{eqnarray}
\label{ch1_eq_ineq_prob}
& \min & f(x), \\  
&  h_i(x) \, \le \, 0, \, i = 1, \, \ldots \, p, & \nonumber \\
& x \, \in \, Q \nonumber & 
\end{eqnarray}
в которой выпуклость функций $f$, $h_1$, $\ldots$, $h_p$ не предполагается.

Определим для задачи~\eqref{ch1_eq_ineq_prob} функцию~$\bar{L}(x) = \sup\limits_{\nu \ge 0} L(x, \, \nu)$ и рассмотрим задачу
\begin{equation}
\label{ch1_eq_min_bL}
\min_{x \, \in \, Q} \bar{L}(x).
\end{equation}
Значение функции $\bar{L}(x)$ равно~$+\infty$ для любой недопустимой точки $x$
задачи~\eqref{ch1_eq_ineq_prob} и равно~$f(x)$
для любой допустимой точки. Следовательно,
задачи~\eqref{ch1_eq_min_bL} и \eqref{ch1_eq_ineq_prob} эквивалентны.
Отсюда очевидна следующая связь между прямой
и двойственной задачами через функцию Лагранжа:
в прямой задаче происходит \underline{минимизация} по~$Q$
результата \underline{максимизации} функции Лагранжа по~$\nu \ge 0$~\eqref{ch1_eq_min_bL}, 
а в двойственной задаче
$$
\max_{\nu \, \ge \, 0} \varphi(\nu) = 
\max_{\nu \, \ge \, 0} \inf_{x \, \in \, Q} L(x, \, \nu)
$$
происходит \underline{максимизация} по~$\nu \ge 0$ результата
\underline{минимизации} функции Лагранжа по~$x \, \in \, Q$. 

Свойство сильной двойственности, или совпадение
решений прямой и двойственной задач\nk{,} выполняется
в случае, если задача~\eqref{ch1_eq_ineq_prob}
выпукла и имеет место условие регулярности Слейтера
(см. п.~\ref{ch1_subs_slater}). 
Имеется ещё один способ сформулировать условие
равенства решений прямой и двойственной задач~---
они равны, если существует {\it седловая точка
функции Лагранжа}, т.е.
точка~$(x^*, \, \nu^*)$ такая, что~$x^* \, \in \, Q$,
и выполняются неравенства:
$$
L(x^*, \, \nu) \le L(x^*, \, \nu^*) \le L(x, \, \nu^*) \quad \, \forall \, x \, \in \, Q, \, \forall \, \nu \ge 0.
$$

Поскольку для любого $y  \in  Q$ 
$\inf\limits_{x \in  Q} L(x,  \nu)  \leq L(y, \nu)$, то
для любого~${y \in Q}$ имеем
$\max\limits_{\nu \, \ge \, 0} \inf\limits_{x \, \in \, Q} L(x, \, \nu)  \leq \max\limits_{\nu \, \ge \, 0} L(y, \, \nu) = \bar{L}(y)$.
Последнее неравенство означает, что
$\max\limits_{\nu  \ge 0} \inf\limits_{x  \in  Q} L(x, \nu)$
является нижней оценкой функции~$\bar{L}(y)$, ${y \in  Q}$, а~также и нижней оценкой функции~$\min\limits_{y \, \in \, Q} \bar{L}(y)$, 
т.е. всегда выполняется неравенство
$$
\max_{\nu \, \ge \, 0} \inf_{x \, \in \, Q} L(x, \, \nu)
\le
\min_{x \, \in \, Q} \sup_{\nu \, \ge \, 0} L(x, \, \nu),
$$
которое есть не что иное, как условие слабой двойственности, не зависящее, впрочем, от свойств
функции\footnote{В общем случае
$$
\sup_{y \, \in \, A} \inf_{x \, \in \, B} f(x, \, y)
\le
\inf_{x \, \in \, B} \sup_{y \, \in \, A} f(x, \, y),
$$
и данное неравенство верно 
для любой
функции~$f: B \times A \rightarrow \mathbb{R}$ и любых
$A \, \subset \, \mathbb{R}^{n_1}$ и $B \, \subset \, \mathbb{R}^{n_2}$.}
$L$.

Если существует седловая точка функции Лагранжа, то решения прямой и двойственной задачи совпадают. Имеет место и обратное утверждение.
\begin{teo}[о структуре множества седловых точек \cite{Ben-Tal}]\label{ch1_teo_saddle_set}  
Множество седловых точек функции~$L: Q \times A \rightarrow \mathbb{R}$ непусто тогда и только тогда,
когда задачи 
\begin{equation}
\label{ch1_eq_sedl1}
\min_{x \, \in \, Q} \sup_{\nu \, \in \, A} L(x, \, \nu)     
\end{equation}
и
\begin{equation}
\label{ch1_eq_sedl2}
\max_{\nu \, \in \, A} \inf_{x \, \in \, Q} L(x, \, \nu)
\end{equation}
разрешимы и оптимальные значения в данных задачах равны друг другу. В этом случае седловые точки функции~$L$~--- это все такие пары~$(x^*, \, \nu^*)$,
где~$x^*$~--- оптимальное решение задачи~\eqref{ch1_eq_sedl1}, а~~$\nu^*$~--- оптимальное решение задачи~\eqref{ch1_eq_sedl2},
и значение функции~$L(\cdot, \, \cdot)$ в каждой
из этих точек равно оптимальному значению
в задаче~\eqref{ch1_eq_sedl1} и в задаче~\eqref{ch1_eq_sedl2}.
\end{teo}

\subsubsection{Существование седловых точек}\label{subs:kakutani}
Вполне <<безобидная>> функция может не иметь седловых точек. 
\begin{example}
Рассмотрим, например, функцию\protect\footnotemark
$$
L(x, \, \nu) = (x - \nu)^2
$$
на единичном квадрате~$[0, \, 1] \, \times \, [0, \, 1]$.
Тогда 
$$
\sup_{\nu \, \in \, [0, \, 1]} L(x, \, \nu) = 
\sup_{\nu \, \in \, [0, \, 1]} (x - \nu)^2 =
\max \{ x^2, \, (1 - x)^2 \},
$$
и
$$
\inf_{x \, \in \, [0, \, 1]} L(x, \, \nu) = 
\inf_{x \, \in \, [0, \, 1]} (x - \nu)^2 = 0.
$$
Поэтому оптимальное значение в задаче~\eqref{ch1_eq_sedl1} равно~$1/4$
(при $x^* = 1/2$), а оптимальное значение в задаче~\eqref{ch1_eq_sedl2} равно~$0$. Следовательно,
по теореме~\ref{ch1_teo_saddle_set} функция~$L$
не имеет седловых точек.
\end{example}
\footnotetext{Пример взят из книги~\cite{Ben-Tal}.}

Возникает вопрос: можно ли по каким-то свойствам функции определить, что она имеет седловые точки?
Более важны в практическом смысле вопросы: при каких условиях эквивалентны
задачи на минимакс и максимин и можно ли менять местами
минимум и максимум? Ответы на эти вопросы дают теоремы Сиона--Какутани. Приведём формулировку 
теоремы об эквивалентности задач на максимин и минимакс. 

\begin{teo}[теорема Сиона--Какутани] \label{minmax_theorem}
Пусть $X$ и $\Lambda$~--- выпуклые множества в~$\mathbb{R}^n$ и $\mathbb{R}^m$ соответственно,
и $X$~--- компактное множество. Если непрерывная функция
$$
L(x, \, \nu) \, : \, X \, \times \, \Lambda \rightarrow \mathbb{R}
$$
выпукла по~$x \, \in \, X$ для любого фиксированного~$\nu \, \in \, \Lambda$ и
вогнута по~$\nu \, \in \, \Lambda$ для
любого фиксированного~$x \, \in \, X$, то
$$
\sup_{\nu \, \in \, \Lambda} \inf_{x \, \in \, X} L(x, \, \nu)
=
\inf_{x \, \in \, X} \sup_{\nu \, \in \, \Lambda} L(x, \, \nu).
$$
\end{teo}

Доказательство теоремы см., например, в~\cite{Ben-Tal}. Обобщения на бесконечномерные пространства, см. в \cite{Stoer}.

\section{Задача квадратичного программирования}\label{ch1_sec_qp}
Задача минимизации квадратичной формы при 
линейных
ограничениях
\begin{equation}
\label{ch1_eq_quad_prob}
 {x^T Cx}  +  {d^T x}  \to 
\mathop {\min }\limits_{Ax \, \le \, b} 
\end{equation}
называется \textit{задачей квадратичного программирования}.
Здесь $x \, \in \, \mathbb{R}^n$, $C$~--- симметричная матрица размера $n \times n$, $b \, \in \, \mathbb{R}^m$,
$A$~--- матрица размера~$m \times n$. 
Квадратичная форма является выпуклой только при $C \succeq 0$. Здесь и далее будем считать, что это условие выполняется.

\begin{figure}[h!]
\begin{center}
\includegraphics[width=0.5\linewidth]{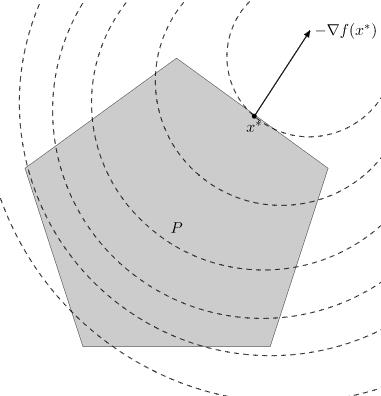}
\end{center}
\caption{Пример графического представления задачи квадратичной оптимизации. Допустимое множество~$P = \{x \, | \, Ax \le b \}$ является полиэдром. Пунктиром обозначены линии уровня целевой функции~$f$. Решение задачи~--- точка~$x^*$ (\cite{Boyd_book_2004}, с изменениями обозначений)}
\label{ch1_quad_prog}
\end{figure}

Функция Лагранжа для задачи~\eqref{ch1_eq_quad_prob}
выглядит следующим образом:
$$
L(x, \, \nu) =  {x^T Cx}  +  {d^T x}  +  \nu^T (Ax - b).
$$
Двойственная функция имеет вид
$$
\varphi(\nu) = \inf_x L(x, \, \nu).
$$
Если квадратичная функция ограничена снизу  на $\mathbb{R}^n$, то она достигает минимума.
Поэтому (см. аналог теоремы Ферма~\ref{ch1_teo_ferma})
выражение для $\varphi(\nu)$
записывается так:
$$
\varphi(\nu) = L(x(\nu), \, \nu),
$$
где $x(\nu)$ находится из условия $L_x'(x, \, \nu) = 0$:
\begin{equation}
\label{ch1_eq_xnu_quad}
2  C x + d + A^T \nu = 0.
\end{equation}
Соотношение~\eqref{ch1_eq_xnu_quad} связывает прямые и двойственные переменные.
Двойственная функция принимает вид
$$
\varphi(\nu) = L(x(\nu), \, \nu) = x^T Cx + (A^T \nu + d)^T x
- b^T \nu = - x^T C x - b^T \nu.
$$
Двойственная к~\eqref{ch1_eq_quad_prob}
задача имеет вид
\begin{equation}
\label{ch1_eq_dual_qp_max}
- x^T C x - b^T \nu  \to 
\mathop {\max }\limits_{2 C x + d + A^T \nu = 0, \, \nu \, \ge \, 0}. 
\end{equation}
Задачу~\eqref{ch1_eq_dual_qp_max} можно переписать
так:
\begin{equation}
\label{ch1_eq_dual_qp_min}
 x^T C x + b^T \nu  \to 
\mathop {\min }\limits_{2 C x + d + A^T \nu = 0, \, \nu \, \ge \, 0}. 
\end{equation}
Важным частным случаем задачи~\eqref{ch1_eq_dual_qp_min}
является случай $C = 0$, тогда задача~\eqref{ch1_eq_quad_prob} превращается
в задачу линейного программирования.
Кроме того, при $C  \succ  0$ 
выражение
для $x$ задаётся явно:
$$
x = \frac{1}{2} C^{-1} (- d - A^T \nu).
$$
Сказанное выше можно подытожить и сформулировать в виде теоремы.
\begin{teo}
Необходимым и достаточным условием экстремума для задачи~\eqref{ch1_eq_quad_prob} при $C \succeq 0$
в точке $x^*$, $Ax^* \, \le \, b$ является
существование $\nu^* \, \in \, \mathbb{R}^m$, $\nu^* \, \ge \, 0$
таких, что
$$
2  C x^* + d + A^T \nu^* = 0, \quad \langle \nu^*, A x^* - b \rangle = 0.
$$
Кроме того, при $C \succ 0$ двойственная задача
для задачи квадратичного программирования~\eqref{ch1_eq_quad_prob}
записывается следующим образом:
\begin{equation}
\label{ch1_eq_dual_quad_prob}    
\min_{\nu \, \ge \, 0} \left \{ \frac{1}{4} \left (
C^{-1} (- d - A^T \nu) \right )^T (-d - A^T \nu) + b^T \nu
\right \}.
\end{equation}

В случае $C \succ 0$ для прямой и двойственной задач выполняется свойство сильной двойственности\protect\footnotemark, т.е. оптимальные значения прямой и двойственной задач квадратичного программирования совпадают.
\end{teo}
\footnotetext{Поскольку прямая задача выпукла, а функ\nk{ции-ог}раничения задачи аффинны.}

Задача квадратичного программирования при $C \succ 0$
сводится к минимизации квадратичной функции~\eqref{ch1_eq_dual_quad_prob} на положительном ортанте. Этот
факт очень существенен, поскольку 
задача квадратичного
программирования и соответственно необходимость её решения часто возникают как отдельный этап
итераций многих современных численных методов оптимизации.

\section{Лемма Неймана--Пирсона}\label{neuman_pirson}
В этом разделе рассматривается применение 
изложенной выше теории (принципа множителей Лагранжа)
для обоснования метода 
отношения правдоподобия\footnote{Впервые был предложен Ю.~Нейманом и Э.~Пирсоном в 1933~г.}
как метода построения  
наиболее мощного критерия в задаче выбора из двух простых статистических гипотез о неизвестном параметре распределения\footnote{Подробнее о фундаментальной лемме Неймана--Пирсона
см. в любом учебнике по математической статистике, например, \cite{Ivch_matstat, Leman_matstat}.}.

Пусть $x=\left( {x_1 , \, \ldots, \, x_n } \right)$~---
реализация {\it выборки} (совокупности наблюдаемых
случайных величин) $X$ и
имеются две простые статистические гипотезы о данной
выборке $H_0$ и $H_1$ с функциями правдоподобия
$L\left(  x \,  | \, H_0  \right)$ и
$L\left(  x \, | \, H_1  \right)$ соответственно.
Простота гипотез означает следующее: 
предполагается, что допустимыми 
распределениями наблюдаемых случайных величин $X$ являются лишь два  
заданных распределения (две функции распределения) $F_0(x)$ и $F_1(x)$, и требуется 
по реализации $x$ проверить гипотезу $H_0$: $F_X(x) = F_0(x)$ 
против альтернативы $H_1$: $F_X(x) = F_1(x)$. 

Принятие или отклонение гипотезы $H_0$ (во втором случае соответственно принятие гипотезы $H_1$) 
осуществляется с зависящей от $x$ вероятностью
$\varphi(x)$
(во втором случае соответственно $1 - \varphi(x)$).
Ясно, что {\it критическая функция} $\varphi(x)$
должна принимать значения только из отрезка $[0, \, 1]$.

Тогда вероятность отклонить верную гипотезу\footnote{Говоря <<гипотеза $H_0$ верна>>, мы имеем в виду, что распределение выборки $F_X(x)$ именно такое, как предполагает гипотеза $H_0$: $F_X(x) = F_0(x)$.}  $H_0$
(т.е. совершить {\it ошибку первого рода}) для выборки $X$ равна
$$
P\left( H_1 \, | \, H_0   \right)=\int_{\mathbb{R}^n} {\varphi \left( x \right)} 
L\left(  x \, | \, H_0  \right)dx,
$$
а вероятность $\beta$ ошибочно принять гипотезу $H_0$
(т.е. совершить {\it ошибку второго рода}) равна
$$
P\left( H_0 \, | \, H_1   \right)= 1 - \int_{\mathbb{R}^n} {\varphi \left( x \right)} 
L\left(  x \, | \, H_1  \right)dx.
$$
Поскольку мы не можем минимизировать обе
вероятности ошибок одновременно при заданном числе испытаний,
выберем число $\alpha$ из отрезка $[0, \, 1]$
({\it уровень значимости}) для
верхней границы вероятности отклонения верной
гипотезы $H_0$
и при этом условии будем минимизировать вероятность совершить ошибку второго рода. 

Почему важно сделать минимальной именно вероятность ошибки второго
рода? Предположим, что мы анализируем работу предприятия,
которое выпускает одинаковые монеты.
Среди этих монет могут оказаться бракованные~--- несимметричные, со смещённым центром тяжести.
Контроль продукции производится таким образом,
что в результате получается какая-то выборка измерений
только что выпущенных монет. И по этой выборке нужно
проверить гипотезу~--- все монеты <<правильные>>,
т.е. соответствуют стандартам. Если гипотеза подтверждается,
монеты прошли контроль. Если нет, то все они будут
забракованы.
Но для предприятия гораздо важнее минимизировать
выпуск бракованных монет. Если даже какие-то <<правильные>>
монеты будут забракованы (мы совершим
ошибку первого рода), это не так критично, как
если какая-то бракованная монета будет признана
настоящей. Именно поэтому обычно 
поступают таким образом: фиксируют ошибку первого рода на достаточно низком безопасном уровне
и 
минимизируют
вероятность ошибки второго рода.

Итак, возникает оптимизационная задача
нахождения оптимальной критической функции
$\varphi(x)$:
\begin{equation}
\label{ch1_eq_min_beta}
\min_{P\left( H_1 \, | \, H_0  \right) \, \le \, \alpha } \beta.
\end{equation}

Для дальнейшего изложения понадобятся
несколько определений.
\begin{defin}
Вероятность $1 - \beta(\varphi)$ называют \emph{мощностью} критерия.
\end{defin}

\begin{defin}
\emph{Критической областью} (областью отклонения
гипотезы $H_0$) называется множество тех 
реализаций $x = \left ( x_1, \, \ldots, \, x_n  \right )$,
для которых принимается гипотеза $H_1$.
\end{defin}

Обычно критическую область задают с помощью
статистики (функции от выборки).

\begin{defin}
В случае двух простых гипотез \emph{наилучшей критической областью} называется область, которая при заданном уровне значимости $\alpha$ обеспечивает минимум ошибки второго
рода.
\end{defin}

Вернёмся к задаче~\eqref{ch1_eq_min_beta}. Применим к её решению
метод множителей Лагранжа:
$$
\min_{0 \, \le \, \varphi \left ( x 
\right) \, \le \, 1} L\left( {\varphi \left( {\,\cdot \,} \right), \, \nu } \right),$$
где
$$L\left( {\varphi \left( {\,\cdot \,} \right), \, \nu } \right)=1-\int_{\mathbb{R}^n} {\varphi 
\left( x \right)} L\left(  x  \, | \, H_1  \right)dx +
\nu 
\left( {\int_{\mathbb{R}^n} {\varphi \left( x \right)} L\left(  x \, | \, H_0  
\right) dx-\alpha } \right). 
$$

После преобразований получаем
\begin{equation}
\label{ch1_eq_min_np}
\min_{0 \, \le \, \varphi \left( x \right) \, \le \,  1} \int_{\mathbb{R}^n} {\varphi \left( x \right)} \left ( \nu L \left (  x \, | \, H_0  
\right ) - L \left (  x \, | \, H_1  \right ) \right ) dx .
\end{equation}

Решение~\eqref{ch1_eq_min_np} можно записать так: 
\begin{equation}
\label{ch1_eq_phi_lnp}
\varphi^* \left( x \right)=\left \{ 
\begin{array}{ll}
 1, & \text{при } \Lambda \left( x \right)>\nu, \\ 
 p\left( x \right), & \text{при } \Lambda \left( x \right)=\nu, \\ 
 0, & \text{при } \Lambda \left( x \right)<\nu, \\ 
 \end{array}
  \right.
\end{equation}
где {\it отношение правдоподобия}
$$
\Lambda \left( x \right)=\frac{L \left( x \, | \,  H_1  
\right ) }
{L\left( x \, | \, H_0  \right)},
$$
функция $p(x)$ и $\nu$ определяются
из условия $\int_{\mathbb{R}^n} {\varphi^* \left( x \right)} 
L\left(  x \, | \, H_0  \right)dx = \alpha$,
а точнее
\begin{equation}
\label{ch1_eq_intlhea}
\Phi(\nu) = \int\limits_{\left \{ x \, | \, \Lambda \left( x \right)>\nu \right \}} {L\left(  x 
\, | \, H_0 \right ) dx} + \int\limits_{\left \{ x \, | \, \Lambda \left( x \right)=\nu \right \} } 
{p\left( x \right)L\left(  x \, | \, H_0 \right ) dx} =\alpha. 
\end{equation}

Множитель $\nu > 0$ определяется
единственным образом (см. ниже), а вероятность ошибки второго рода~$\beta(\varphi^*)$ 
не зависит от выбора функции $p\left( x \right)$. 
Часто функцию~$p(x)$ полагают равной константе.
Возникает вопрос: даже если считать $p(x)$ константой,
почему есть гарантия, что из~\eqref{ch1_eq_intlhea}
можно однозначно определить~$\nu$? Дело в том, что
для непрерывной функции по теореме Брауэра о неподвижной точке в какой-то момент функция~$\Phi(\nu)$
будет иметь пересечение с уровнем~$\alpha$. Если
же функция~$\Phi(\nu)$ не непрерывна, а имеет скачки,
то при попадании~$\alpha$ на этот скачок можно определить~$p(x)$.

Далее приведём формулировку леммы Неймана--Пирсона.
\begin{lemma} [{\emph{Неймана--Пирсона}}]
Пусть $\alpha  \in [0, \, 1]$, тогда при заданном уровне значимости $\alpha$ наиболее мощный среди всех $\varphi$-критериев с уровнем значимости $\alpha$ имеет критическую функцию $\varphi^*(x)$ вида~\eqref{ch1_eq_phi_lnp},
где $0  \le \linebreak \le p(x) \le  1$, функция правдоподобия
$L\left(  x \,  | \, H_0  \right)$
соответствует гипотезе $H_0$,
а $L\left(  x \,  | \, H_1  \right)$~--- конкурирующей гипотезе $H_1$.
\end{lemma}

Таким образом, с помощью метода множителей Лагранжа можно доказать фундаментальную лемму Неймана--Пирсона
о существовании наиболее мощного критерия.
Кроме того, лемма Неймана--Пирсона показывает, как устроена наилучшая критическая область
(в неё включаются те точки, вероятность появления которых
при $H_1$ существенно больше, чем при $H_0$):
$$
\left \{ x \, | \, \Lambda \left( x \right) \, \ge \,  \nu  \right \}.
$$

Более подробно о проверке статистических гипотез
с помощью выпуклой оптимизации можно прочитать в доступной в электронном виде книге~\cite{StatOpt_LN}.

На самом деле выпуклость постановок задач, в которых присутствует интегрирование, имеет место при намного более общих условиях, чем рассмотренные выше. Приведём сначала яркую цитату на эту тему из введения к одной из лучших известных нам книг по выпуклому анализу \cite{MagTih_convan}: <<\textit{Здесь уместно раскрыть некоторые причины большого прикладного значения выпуклого анализа. Одна из главных заключена в тезисе: \textbf{интегрирование по непрерывной мере порождает выпуклость образа интегрального отображения}. Математическим проявлением этого тезиса является теорема Ляпунова о векторных мерах. Существенная доля основополагающих фактов теории классического вариационного исчисления и оптимального управления связана с этим тезисом (условия Лежандра и Вейерштрасса, принцип максимума Понтрягина и т.п.). А через теорию оптимизации легко объясняется широкий спектр приложений выпуклого анализа к геометрии, теории аппроксимаций, теории неравенств и т.п. В экономике типично суммирование многих мелких факторов. Такие суммы естественно аппроксимируются интералами, и это также даёт ключ к пониманию роли выпуклого анализа в математической экономике}>>.

Поясним теперь написанное выше одним простым примером. Рассмотрим два функционала 
$$F(u) = \int_{0}^1 f\left(t,u(t)\right)dt,$$
$$H(u) = \int_{0}^1 h\left(t,u(t)\right)dt,$$
заданных на множестве всех ограниченных измеримых функций $u(~\cdot~)$ на отрезке $[0,1]$.
Предполагается, что интегралы всегда существуют. 
Покажем, что для любых (ограниченных измеримых) $u_1$ и $u_2$ и любого $\alpha \in (0,1)$ найдётся (также ограниченная измеримая)
$u_{\alpha}$ такая, что (эти равенства и выражают ту выпуклость, которую порождает интегрирование):
\begin{equation}
   \label{Funct} 
F\left(u_{\alpha}\right) = \alpha F\left(u_1\right) + (1-\alpha)F\left(u_2\right),
H\left(u_{\alpha}\right) = \alpha H\left(u_1\right) + (1-\alpha)H\left(u_2\right),
\end{equation}
Для простоты дополнительно предположим, что $u_1$ и $u_2$ --- кусочно-постоянные функции (в качестве пределов таких функций можно получать измеримые функции). Более того, будем также считать, что отрезок [0,1] разбит на равные промежутки, на каждом из которых эти функции постоянны. Тогда $u_{\alpha}$ можно выбирать, например, следюущим образом: на первой $\alpha$-части каждого промежутка $u_{\alpha} = u_1$, а на оставшейся $(1 - \alpha)$-части $u_{\alpha} = u_2$. Осталось только перейти к пределу по мелкости разбиения.  

\ag{Теперь нужно сделать главное наблюдение, что свойство \eqref{Funct}  обеспечивает возможность применять к задаче 
$$\min_{H(u)\le 0} F(u)$$
принцип множителей Лагранжа. Это следует из выпуклости множества $\mathbb{A}$ (см. формулу~\eqref{ch1_eq_a}) и, как следствие, возможности применять теорему об отделимости граничной точки этого множества (отвечающей решению задачи) от этого множества, см. п.~\ref{geom}.}

\ag{Написанное выше останется верным, если вместо ограничения вида неравенства (или вместе с ним) присутсвует и ограничение вида равенства
$$\int_{0}^1 g\left(t,u(t)\right)dt = 0.$$}

\section{Лемма Фаркаша}\label{ch1_sect_fark}
В теории математического программирования важную роль играют
теоремы об альтернативных системах
линейных неравенств. Большинство таких
теорем об альтернативах имеют следующий вид: 
даны две альтернативные
системы линейных неравенств,
построенные по определённым правилам. 
Из~этих двух систем будет
иметь решение либо в точности одна (тогда
такие альтернативы называются {\it сильными}), 
либо не \nk{более чем} одна
система (тогда
такие альтернативы называются {\it слабыми}).
Формулировок таких теорем
существует достаточно большое
количество, есть попытки их классифицировать~\cite{Motzkin}.
Иногда все эти теоремы обобщённо называют
{\it теоремами или леммами Фаркаша}, по фамилии венгерского математика Юлиуса Фаркаша,
доказавшего в 1902 г. одну из теорем об альтернативных системах линейных
неравенств и равенств. Но чаще леммой Фаркаша называют только одно из утверждений об альтернативных системах линейных неравенств и равенств. Остальные варианты формулировок связаны с именами многих математиков (Моцкина, Гейла, Фредгольма, Минковского, Гордана и других~\cite{Zork, Motzkin}).
Лемма Фаркаша используется, например, в
следующих случаях:
в одном из способов
доказательства необходимых условий оптимальности
первого порядка для задач с ограничениями типа равенств и неравенств (условий Каруша--Куна--Таккера, см. п.~\ref{ch1_subs_teo_kkt}); при определении принадлежности
вектора конусу; конечно, как результат,
описывающий все препятствия к решению систем линейных неравенств над вещественными числами, и во многих других областях.

Напомним, что непустое множество~$K$ называется {\it конусом}, если для любого $x \, \in \, K$ и  
$\lambda \, \ge \, 0$
выполняется включение $\lambda x \, \in \, K$.
Конусом является, например, множество $C = \{ c \, | \, c = Ax, \, x \ge 0\}$.

Следующее утверждение, очевидное на первый взгляд, но требующее нетривиальных рассуждений, будет использовано при доказательстве леммы Фаркаша.
\begin{lemma}[\cite{Karmanov, Nocedal}]\label{ch1_lem_close_fark}
Конус~$C = \{c \, | \, c = Ax, \, x \, \in \, \mathbb{R}^n, \, 
x \ge 0\}$ является замкнутым множеством.
\end{lemma}

\begin{proof}
Разобьём матрицу~$A$ размерности~$m \times n$ на столбцы:
$$
A = [a_1, \, a_2, \, \ldots, \, a_n].
$$
Докажем факт замкнутости множества~$C$ с помощью метода математической индукции по числу~$n$.

При $n = 1$ множество~$C$ является лучом, следовательно,
$C$~замкнуто.

Предположим, что при $n = k - 1$ конус~$C^{k-1}$, порождённый
векторами $a_1, \, a_2, \, \ldots, \, a_{k-1}$, замкнут.

\begin{enumerate}
    \item Если векторы $-a_1, \, -a_2, \, \ldots, \, -a_{k}$
    принадлежат~$C$, то $C$ является подпространством размерности,
    не превышающей $k$, и, следовательно, замкнутым множеством.
    
    \item Предположим, что хотя бы один из векторов $-a_{k'} \notin C$. Любой вектор $c \, \in \, C$ представим в виде $c = \bar{c} + \alpha a_k$, где $\alpha \ge 0$, $\bar{c}  \, \in \, C^{k-1}$.
    Рассмотрим последовательность~$\{ c_r \} \subset C$, сходящуюся к некоторому вектору $c$.
    Для всех $r$ $c_r = \bar{c}_r + \alpha_r a_k$, $\alpha_r \ge 0$.
    Если последовательность $\{ \alpha_r \}$ ограничена, то,
    без потери общности, можно считать, что $\lim\limits_{r \rightarrow \infty} \alpha_r = \alpha$ и, следовательно,
    $c - \alpha a_k = \lim\limits_{r \rightarrow \infty} (c_r - \alpha_r a_k) = \lim\limits_{r \rightarrow \infty} \bar{c}_r = \bar{c} \, \in \, C^{k-1}$, так как $C^{k-1}$ замкнуто по предположению индукции. Следовательно, $c = \bar{c} + \alpha a_k \, \in \, C$. Итак, доказано, что если последовательность $\{ \alpha_r \}$ ограничена,
    то множество~$C$ замкнуто.
    
    Предположим, что последовательность $\{ \alpha_r \}$ не является ограниченной, т.е. $\{ \alpha_r \} \rightarrow + \infty$
    при $r \rightarrow \infty$. Тогда $\frac{1}{\alpha_r} c_r = \frac{1}{\alpha_r} \bar{c}_r + a_k$ и $\lim\limits_{r \rightarrow \infty} \frac{1}{\alpha_r} \bar{c}_r = - a_k$.
    В силу замкнутости $C^{k-1}$ $\bar{c}_r \, \in \, C^{k-1}$, соответственно, и $-a_k \, \in \, C^{k-1}$,
    что противоречит предположению второго пункта данного доказательства. 
\end{enumerate} 
\end{proof}

\begin{lemma}[{лемма Фаркаша--Минковского}]\label{ch1_lem_farkMin}
Пусть $A$~--- вещественная матрица размера~$m \times n$,
$x$~--- вектор-столбец из $n$ компонент, $y$ и $b$~--- векторы-столбцы из $m$ компонент. Верна 
только одна из двух следующих альтернатив\footnote{Следовательно, альтернативы 1 и 2 \emph{сильные}.}{\rm :}
\begin{enumerate}
\renewcommand\labelenumi{\textup{\arabic{enumi}.}}
    \item Существует $x \, \in \, \mathbb{R}^n$ такой, что
    $x \ge 0$ и $Ax = b$.
    \item Существует $y \, \in \, \mathbb{R}^m$ такой, что
    $y^T A \ge 0$ и $y^T b + 1 = 0$.
\end{enumerate}
\end{lemma}

\noindent{\it Доказательство} леммы~\ref{ch1_lem_farkMin} частично следует из~\cite{Wildon_farkas}. Предположим, что первая
альтернатива не выполняется, т.е. вектор~$b$ не лежит в конусе: 
$$
C = \{c \, | \, c = Ax, \, x \, \in \, \mathbb{R}^n, \, 
x \ge 0\}.
$$
Необходимо доказать, что тогда верна вторая альтернатива.
Пусть $B$~--- замкнутый шар радиуса~$R$ с центром в точке~$b$, где $R \, \ge \, \|b\|$. Тогда $B \cap C \neq \varnothing$. Конус~$C$ является замкнутым множеством по лемме~\ref{ch1_lem_close_fark}.

Поскольку $C$ замкнуто и $B$~--- компактное множество, существует ближайший в евклидовой норме к~$b$  вектор~$w \, \in \, B \cap C$, при этом~$\| b - w \| \le R$.

Возьмём любой вектор $v \, \in \, C$ и рассмотрим
отрезок, соединяющий концы векторов~$v$ и $w$. В силу
выпуклости множества~$C$, данный отрезок принадлежит~$C$. Поэтому для любого~$\alpha \, \in \, [0, \, 1]$ имеем
\begin{eqnarray*}
\| b - w \|^2 = \| w - b \|^2 & \le & \| \alpha v + ( 1 - \alpha ) w - b \|^2 \\
& = & \| \alpha (v - w) + (w - b) \|^2 \\
& = & \alpha^2 \| v - w \|^2 + 2 \alpha (v - w)^T
(w - b) + \| w - b \|^2.
\end{eqnarray*}
Следовательно, при положительном~$\alpha$ выполняется неравенство $(v - w)^T\!\times\linebreak\times\,(w - b) \ge 0$.

Обозначим~$z = w - b$. При $v = 0$ (нулевой вектор также принадлежит~$C$) скалярное произведение~$w^T z \le 0$,
поэтому
$$
b^T z = (b - w)^T z + w^T z = - \| b - w \|^2 + w^T z \, < \, 0,
$$
и существует такое~$\gamma \, < \, 0$, что
$$
v^T z \ge w^T z \, > \, \gamma \, > \, b^T z
$$
для всех $v \, \in \, C$. Зафиксируем какой-то вектор $v \, \in \, C$. Поскольку $\lambda v^T z \, > \, \gamma$
для всех~$\lambda \ge 0$, то после деления на $\lambda \neq 0$ получаем $v ^T z \, > \, \gamma / \lambda$. При достаточно больших $\lambda$ $v^T z \ge 0$. Поэтому $y = z / | z^T b | \, \in \, \mathbb{R}^m$~--- вектор, удовлетворяющий
альтернативе~2. 
$\qed$

Геометрический смысл леммы Фаркаша в формулировке~\ref{ch1_lem_farkMin} состоит в следующем:
{\bf либо} вектор~$b$ принадлежит конусу, натянутому
на векторы\footnote{Т.е.
$b$ можно представить в виде линейной комбинации
столбцов $Ae_1$, $A e_2$, $\ldots$, $A e_n$
с неотрицательными коэффициентами, иными словами,
$Ax = b$ при $x \ge 0$.} $Ae_1$, $A e_2$, $\ldots$, $A e_n$, где
$e_1, \, e_2, \, \ldots, \, e_n$~--- базис в~$\mathbb{R}^n$,
{\bf либо} существует гиперплоскость $H = \{ v \, \in \, \mathbb{R}^n \, | \, y^T v = 0\}$, отделяющая 
вектор~$b$ от вышеупомянутого конуса.
Следствием такой геометрической интерпретации леммы
Фаркаша является следующий факт: проблема
поиска допустимых решений в задаче линейного программирования
имеет такую же сложность, как задача поиска
плоскости, отделяющей точку от конуса. И второе: лемма
Фаркаша--Минковского
особенно полезна, например, в таком случае~--- если существует вектор~$y$, удовлетворяющий второй альтернативе 
леммы~\ref{ch1_lem_farkMin}, то можно сразу заключить, что задача~ЛП с ограничениями
вида $\left \{x \, | \, x \ge 0, \, Ax = b \right \}$ (первая альтернатива)
не имеет решений.

Приведём ещё одну формулировку леммы Фаркаша в терминах сопряжённых конусов.

\begin{lemma}[ лемма Фаркаша, версия~2]\label{ch1_lem_fark_pol_cone}
Пусть $K_i = \{ x \, \in \, \mathbb{R}^n \, | \, (a^i)^T x \ge 0\}$, $i = 1, \, \ldots, \, m$, $K = K_1 \cap \, \ldots \, \cap K_m \neq \varnothing$. Тогда
$$
K^* = K_1^* + \, \ldots \, + K_m^* = \left \{ \sum_{i = 1}^m y_i a^i, \, y_i \ge 0 \right \}.
$$
\end{lemma}
Тот же результат можно записать в более узнаваемой по сравнению с леммой~\ref{ch1_lem_farkMin} форме:

Пусть $A$~--- матрица размера~$m \times n$, строками которой являются~$a^i$. Если 
$K = \{ x \, \in \, \mathbb{R}^n \, | \, Ax  \ge 0\}$,
то
$K^* = \left \{ x \, \in \, \mathbb{R}^n \, | \, x = A^T y, \, y \ge 0, \, y \, \in \, \mathbb{R}^m \right \}$.

Возникает вопрос: существует ли вариант леммы Фаркаша не только для многогранных конусов, как в лемме~\ref{ch1_lem_fark_pol_cone}, но и для конусов
общего вида? Ответ~--- да, существует, для замкнутых выпуклых конусов.

\begin{lemma} [лемма Фаркаша, версия~3]\label{ch1_lem_fark_cone}
Рассмотрим систему 
\begin{equation}
\label{ch1_eq_fark3_sys1}
\left \{ x \, | \, Ax = b, \, x \, \in \, C \right \},    
\end{equation}
где матрица $A$ размера $m \times n$,
$b \, \in \, \mathbb{R}^m$, $C$~--- выпуклый замкнутый конус. Предположим, что существует $\bar{y}$ такой, что
$-A^T \bar{y} \, \in \, \textup{\mbox{int}} \, C^*$. Тогда система~\eqref{ch1_eq_fark3_sys1} имеет допустимое решение~$x$ тогда и только тогда, когда система
$$
\left \{ y \, | \, -A^T y \, \in \, C^*, \, b^T y \, > \, 0 \, \, (b^T y = 1) \right \}
$$
не имеет решений.
\end{lemma}

%
\subsection{Двойственность и лемма Фаркаша}
Для задачи ЛП
\begin{eqnarray}
\label{ch1_eq_lp_fark}
& \min & c^T x \\  
&  Ax \ge b &  \nonumber
\end{eqnarray}
двойственная к ней выглядит следующим образом:
\begin{eqnarray}
\label{ch1_eq_dual_lp_fark}
& \max & b^T \lambda. \\  
&  A^T \lambda - c  = 0 & \nonumber \\
& \lambda \ge 0 \nonumber & 
\end{eqnarray}
Предполагается, что решение двойственной задачи
существует. Обозначим его~$\varphi^*$.
Рассмотрим задачу
\begin{eqnarray}
\label{ch1_eq_dual2_lp_fark}
& \min & 0^Tx. \\  
&  Ax \ge b &  \nonumber \\
& c^T x \le \varphi^* & \nonumber 
\end{eqnarray}

Для применения к данным задачам 
понадобится лемма Фаркаша в следующей
формулировке.
\begin{lemma}[лемма Фаркаша, версия~4]\label{ch1_lem_fark4}
Пусть $A$~--- вещественная матрица размера~$m \times n$,
$x$~--- вектор-столбец из $n$ компонент, $y$ и $b$~--- векторы-столбцы из $m$ компонент. Верна 
только одна из двух следующих альтернатив{\rm :}
\begin{enumerate}
\renewcommand\labelenumi{\textup{\arabic{enumi}.}}
    \item Существует $x \, \in \, \mathbb{R}^n$ такой, что
    $Ax \ge b$.
    \item Существует $y \, \in \, \mathbb{R}^m$ такой, что
    $y \ge 0$, $y^T A = 0$ и $y^T b \, > \, 0$.
\end{enumerate}
\end{lemma}

Если решение задачи~\eqref{ch1_eq_dual2_lp_fark} существует, то и у~\eqref{ch1_eq_lp_fark}
должно существовать решение с оптимальным
значением не больше~$\varphi^*$.

Предположим, что~\eqref{ch1_eq_dual2_lp_fark}
не имеет решения. Тогда
по лемме Фаркаша~4 выполняется альтернатива 2, соответственно существует
$\left ( \begin{array}{c}
    \lambda   \\
     \nu  
\end{array} \right ) \ge 0$ такой, что
$$
(\lambda^T\ \nu) \left ( \begin{array}{c}
    A   \\
     -c^T  
\end{array} \right ) = 0 \mbox{ и }
(\lambda^T\ \nu) \left ( \begin{array}{c}
    b   \\
     -\varphi^*  
\end{array} \right ) \, > \, 0.
$$
Отсюда $\lambda^T A - \nu c^T = 0$ и $\lambda^T b - \nu \varphi^* \, > \, 0$.

1. Если $\nu = 0$, то $\lambda^T A = 0$, и
$\lambda^T b \, > \, 0$.
Выберем решение задачи \eqref{ch1_eq_dual_lp_fark} $z \ge 0$ такое, что $A^T z = c$ и
$b^T z = \varphi^*$. Тогда для~$\varepsilon \, > \, 0$
$$
A^T (z + \varepsilon \lambda) = c,
\quad z + \varepsilon \lambda \ge 0,
$$
поскольку $\lambda \ge 0$, и
$$
b^T ( z + \varepsilon \lambda) = \varphi^* +
\varepsilon b^T \lambda > \varphi^*,
$$
что противоречит условию
оптимальности значения~$\varphi^*$ для
задачи~\eqref{ch1_eq_dual_lp_fark}.

2. Иначе, при $\nu > 0$, не теряя общности,
можно считать, что $\nu = 1$ (нужно просто уменьшить $\lambda \ge 0$ в соответствующее число раз). Тогда $\lambda^T A = c^T$, $\lambda^T b \, > \, \varphi^*$. Но 
тогда $\lambda$ также является решением~\eqref{ch1_eq_dual_lp_fark} со значением, большим $\varphi^*$. Противоречие.

Поэтому решение задачи~\eqref{ch1_eq_dual2_lp_fark} существует,
соответственно существует и решение
задачи~\eqref{ch1_eq_lp_fark} с оптимальным
значением не больше~$\varphi^*$.
\qed

\subsection[Связь теоремы Гильберта о нулях и леммы Фаркаша]{Связь теоремы Гильберта о нулях\\ и леммы Фаркаша}\label{ch1_subsect_Hilb}
{\it Теорема Гильберта о нулях}
(другие её названия~ --- {\it теорема Гильберта о корнях}, часто используют первоначальное немецкое название {\it Nullstellensatz})~---
одна из фундаментальных теорем, устанавливающих взаимосвязь между геометрией и алгеброй
(и лежащей в основе алгебраической геометрии). 
Она
показывает, как выглядят {\it все следствия} данной системы полиномиальных алгебраических уравнений, если мы решаем систему в алгебраически замкнутом поле\footnote{Поле $K$ {\it алгебраически замкнуто}, если каждый отличный от константы многочлен
из $F[x]$ (кольца многочленов от переменной $x$
с коэффициентами из $F$) имеет хотя бы один корень
в $F$
(т.е. хотя бы один линейный множитель в $F[x]$).
Поле~--- это коммутативное кольцо с единицей, в котором каждый ненулевой элемент имеет мультипликативный обратный элемент (т.е. обратный по умножению). Пример алгебраически замкнутого поля~--- поле комплексных чисел~$\mathbb{C}$.}.
Впервые она была доказана одним из великих математиков своего времени  Давидом Гильбертом (1862--1943).

Пусть $F$~--- фиксированное алгебраически замкнутое поле.
Теорему Гильберта о нулях обычно формулируют
в алгебраических терминах идеалов и радикалов идеалов,
но её можно переформулировать на языке систем
полиномиальных уравнений $P_1(x) = \ldots = P_m(x) = 0$ от переменных $x = (x_1, \, \ldots, \, x_d) \, \in \, F^d$ над $F$ (т.е. $P_1, \, \ldots, \, P_m \, \in \, F[x]$~--- многочлены от $d$ переменных).

\begin{teo}[теорема Гильберта о нулях \cite{ttao}] \label{ch1_teo_null}
Пусть $P_1$, $\ldots ,$ $P_m, \, R  \in\linebreak \in F[x]$~--- многочлены. Тогда выполняется \emph{ровно одно}
из следующих утверждений{\rm :}
\begin{enumerate}
\renewcommand\labelenumi{\textup{\arabic{enumi}.}}
    \item Система уравнений $P_1(x) = \ldots = P_m(x) = 0$; $R(x) \neq 0$ имеет решение $x \, \in \, F^d$.
    \item Существуют такие многочлены $Q_1, \, \ldots, \, Q_m \, \in \, F[x]$ и такое неотрицательное целое число $r$, что
    $P_1 Q_1 + \ldots + P_m Q_m = R^r$.
\end{enumerate}
\end{teo}

Предположение об алгебраической замкнутости
поля~$F$ для теоремы~\ref{ch1_teo_null} существенно.
Например, если $F$~--- поле вещественных чисел
$\mathbb{R}$, то уравнение $x^2 + 1 = 0$ не имеет
решения, но не существует такого многочлена $Q_1(x)$, что выполняется $(x^2 + 1) Q_1(x) = 1$.

Эквивалентность формулировки теоремы~\ref{ch1_teo_null}
и теоремы о нулях в терминах идеалов и радикалов
следует из не менее фундаментальной теоремы Гильберта о базисе~\cite{arj_basis, Churkin_pol_eq}, из
которой следует, что всякая бесконечная система полиномиальных уравнений от конечного числа переменных над полем равносильна своей конечной
подсистеме.

В формулировке теоремы~\ref{ch1_teo_null}
приведена ссылка на книгу <<Структура и случайность>>~\cite{ttao} одного из самых
знаменитых современных математиков,
лауреата Филдсовской премии
Теренса Тао\footnote{Книга является
переводом части статей из блога
Тао \url{https://terrytao.wordpress.com},
посвящённого обзору,
популярному изложению и обсуждению
открытых математических проблем и
сюжетов,
а также других, связанных с математикой и физикой, тем.
Например, как пишет Тао в предисловии к книге, самой читаемой и комментируемой записью
за первый месяц существования блога стала
популярная статья <<Квантовая механика и расхитительница гробниц>>.}. В параграфе, посвящённом
теореме Гильберта о нулях, Тао пишет, что
сила этой теоремы состоит в том, что из предположения о несуществовании (решения системы
полиномиальных уравнений) извлекается вывод о существовании
(многочленов $Q_1, \, \ldots, \, Q_m$). Такая
возможность получить <<что-то из ничего>>
должна быть очень нетривиальной и полезной.
И действительно, помимо вышеупомянутой
основополагающей связи между
алгеброй (коммутативной, утверждение~$2$ теоремы~\ref{ch1_teo_null}) и геометрией
(алгебраической, утверждение~$1$ теоремы~\ref{ch1_teo_null}),
оказывается, например,
что разбираемую в данных разделах
лемму Фаркаша можно
сформулировать так, что она будет
имитировать теорему Гильберта о нулях\footnote{
И доказать эти утверждения можно
будет с помощью аналогичных приёмов,
индукцией по размерности, см. доказательства
в~\cite{ttao}.}. И если теорема
Гильберта была проинтерпретирована Т.~Тао~\cite{ttao}
как результат, описывающий все препятствия
к решению системы полиномиальных уравнений и неравенств над алгебраически замкнутым полем,
то лемму Фаркаша он формулирует как результат,
описывающий все препятствия к решению системы линейных неравенств над вещественными числами.
Далее эта лемма связывается с результатами,
близкими теореме Хана--Банаха: теоремой Дьёдонне
об отделимости, теоремой фон Неймана о минимаксе,
теоремой Менгера и теоремой Хелли. Ниже
приведём ещё одну формулировку
леммы Фаркаша и теоремы Дьёдонне и Хелли, следуя Т.~Tао~\cite{ttao}.
Рассматривается только конечномерный случай.

Рассмотрим систему линейных неравенств вида
\begin{equation}
\label{ch1_eq_slineq_tao}
P_1(x), \, \ldots, \, P_m(x) \, \ge \, 0,
\end{equation}
где $x \, \in \, \mathbb{R}^d$ лежит
в конечномерном вещественном векторном пространстве,
а $P_1(x), \, \ldots, \, P_m(x) : \mathbb{R}^d \rightarrow \mathbb{R}$~--- аффинные
функционалы.
Классический вопрос ЛП: имеет ли решение система~\eqref{ch1_eq_slineq_tao}? Очевидное препятствие к разрешимости: система неравенств~\eqref{ch1_eq_slineq_tao} может быть несовместна в том смысле, что из неё можно вывести
неравенство $-1 \, \ge \, 0$. Точнее,
если найдутся такие вещественные числа
$q_1$, $q_2$, $\ldots$, $q_m \, \ge \, 0$,
что
$$
q_1 P_1 + \ldots + q_m P_m = -1,
$$
то система неравенств~\ref{ch1_eq_slineq_tao} неразрешима. Лемма Фаркаша утверждает,
что такое препятствие~--- единственно возможное.

\begin{lemma} [лемма Фаркаша, версия~5~\cite{ttao}]\label{ch1_lem_fark5}
Пусть $P_1(x), \, \ldots, \, P_m(x): \mathbb{R}^d \rightarrow \mathbb{R}$~--- аффинные
функционалы. Тогда выполняется ровно одно из следующих утверждений{\rm :}
\begin{enumerate}
\renewcommand\labelenumi{\textup{\arabic{enumi}.}}
    \item Система неравенств $P_1(x), \, \ldots,
    \, P_m(x) \, \ge \, 0$ имеет решение $x \, \in \, \mathbb{R}^d$.
    \item Существуют такие вещественные числа
$q_1$, $q_2$, $\ldots$, $q_m \, \ge \, 0$,
что
$$
q_1 P_1 + \ldots + q_m P_m = -1.
$$
\end{enumerate}
\end{lemma}
Доказательство леммы с помощью индукции по $d$ см. в~\cite{ttao}.

\begin{teo}[теорема Дьёдонне об отделимости]
Пусть $A$, $B$~--- непересекающиеся выпуклые многогранники в $\mathbb{R}^d$. Тогда их можно отделить гиперплоскостью, т.е. существует такой аффинно-линейный функционал $P : \mathbb{R}^d \rightarrow \mathbb{R}$,
что $P(x) \geq 1$ при $x \, \in \, A$
и $P(x) \leq -1$ при $x \, \in \, B$.

\end{teo}

Приведём в заключение важнейший результат выпуклой конечномерной геометрии, теорему Хелли, которая
имеет многочисленные приложения в различных областях.

\begin{teo}[теорема Хелли, формулировка Т.\;Тао~\cite{ttao}]
Пусть $B_1, \, \ldots, \, B_m$~--- набор выпуклых тел в $\mathbb{R}^d$, 
причём $m > d$. Если любые $d + 1$ из этих тел имеют общую точку, то все $m$~тел имеют общую точку.
\end{teo}
Доказательство основано на лемме Фаркаша и свойствах линейной комбинации векторов.

Возможна и такая формулировка теоремы Хелли.
\begin{teo}[теорема Хелли~\cite{MagTih_convan}]
Пусть $\mathcal{A}$~--- некоторое подмножество индексов и $\{ \mathcal{A}_\alpha\}_{\alpha \, \in \, \mathcal{A}}$~--- семейство замкнутых выпуклых подмножеств~$\mathbb{R}^d$,
по крайней мере одно из которых компактно. Тогда если любое подсемейство из $d + 1$ множества имеет непустое пересечение, то и всё семейство имеет непустое пересечение.
\end{teo}
В такой формулировке эту теорему можно доказать индукцией по $d$ с применением теоремы Радона (любое конечное множество в~$\mathbb{R}^d$, состоящее из не менее чем $d+2$~точек, можно разбить на два непересекающихся подмножества так, что выпуклые оболочки этих подмножеств не будут пересекаться). Ещё одно из многочисленных доказательств теоремы Хелли основано на использовании теорем об очистке (см. п.~\ref{ch1_sect_dem-dan} и \cite{Klimov_matan}).

\subsection{Арбитражная теорема}
В данном параграфе рассматривается модель CRR (Кокса––Росса––\linebreak Рубинштейна, Cox––Ross––Rubinstein \cite{CRR, Ross_arbitrage}), которая играет в финансовой математике роль,
сходную со схемой Бернулли в классической теории вероятностей. С помощью этой модели
можно определить многие финансовые характеристики, например, справедливые цены опционов. Удивительным образом оказывается,
что решение данной задачи 
однозначно связано с одной
из рассмотренных в параграфе~\ref{ch1_sect_fark} лемм Фаркаша.

На фондовом рынке существуют различные типы ценных бумаг. Мы рассмотрим следующие: {\it акции} (stock) и {\it опционы} (option). Цены всех бумаг моделируются случайными процессами на дискретном или непрерывном временном интервале.

Акции~--- это долевые ценные бумаги, выпускаемые компаниями и корпорациями для
увеличения капитала.

Опционы~--- это ценные бумаги, которые дают владельцу право, но не обязанность, купить или продать (у того,
кто предоставил владельцу данный опцион)
фиксированное количество соответствующих финансовых активов по определённой цене в конкретный период.

Опцион-колл (call) даёт право купить, а опцион-пут (put)~--- соответственно право продать финансовый
актив. Существуют два типа опционов по времени исполнения: европейские и американские.

Если опцион можно исполнить только в заранее определённый
момент времени~$N$, то говорят, что $N$~--- момент исполнения, а опцион является {\it опционом европейского типа}.
Если же опцион может быть предъявлен к исполнению в любой (случайный) момент времени, но не позже определённого $N$, то говорят, что это {\it опцион американского типа}.

Если, например, оказалось, что
для опциона европейского типа на
покупку акций эти акции
подорожали к дате исполнения опциона,
то владелец опцион\nk{а-к}олл воспользуется
своим правом покупки и получит
прибыль, купив акции
по цене дешевле рынка
(продавец опциона
обязан продать эти акции по такой цене). Если же цена акций оказалась меньше цены исполнения опциона, то 
владельцу уже невыгодно покупать акции по 
цене дороже рынка,
и он не исполняет свой опцион, т.е. ничего не делает
(ведь опцион даёт право, ни к чему не обязывая).

Пусть на <<идеализированном>> фондовом рынке имеется всего две ценные бумаги, и торговля
осуществляется в два момента времени. Пусть цена первой бумаги~$S$ (акции) известна в первый момент времени.
Цена второй бумаги~$C$ (опцион европейского типа) неизвестна в первый момент. Пусть с ненулевой вероятностью
$p > 0$ ($p$ неизвестно, но от него ничего зависеть в итоге не будет)
к моменту времени~$2$ цена акции вырастет в $u > 1$ (up) раз и с вероятностью $1 - p$ цена акции <<вырастет>> в $d < 1$ (down) раз, т.е. упадёт.
Пусть также известны возможные цены опциона во второй момент времени:
$C_u$, если акция выросла в цене, и $C_d$, если акция упала в цене. Для простоты
считается, что банк работает с нулевым процентом, т.е.
класть деньги в банк в расчёте на проценты бессмысленно.

В данном случае дата исполнения опциона~--- момент времени~$2$, платежи в момент исполнения~--- $C_u$
и $C_d$. Причём эти платежи~--- заранее известные функции от цены
акции в этот момент (введение опционов было мотивировано желанием <<хеджироваться>> (страховаться) от нежелательных изменений цен
акций). Основная задача заключается в установлении <<справедливой>>
цены опциона~$C$ в момент времени~$1$ (см.~\cite{Shiryaev_finn}).

Допустим, что в начальный момент времени~$1$ мы приобретаем $k_S$ единиц акции и $k_C$ единиц опции. При этом $k_S,k_C$ могут быть как положительными, так и отрицательными. В последнем случае мы <<встаём в короткую позицию>>~--- продаём ценные бумаги, приобретённые в долг, с обязательством вернуть их в момент времени~$2$ по цене $S$. Наши затраты на эти транзакции составляют 
\begin{equation}
\label{ch1_eq_arb_0}
X(1) = k_S S + k_C C    
\end{equation}
единиц капитала. При этом $X(1)$ может быть и отрицательным, в случае если выручка от проданных в долг активов превышает цену приобретённых.

Ситуация, в которой вырученная закрытием позиций в момент времени~$2$ сумма $X(2)$ наверняка не меньше затраченной изначально суммы $X(1)$ и при этом с некоторой положительной вероятностью даже строго больше, называется \emph{арбитражем}. Другими словами, арбитраж имеет место, когда
\begin{equation}
\label{ch1_eq_arb_x2}
P\{X(2) \ge X(1)\} = P\{k_S S(2) + k_C C(2) \ge k_S S + k_C C\} = 1,   
\end{equation}
причём 
\begin{equation}
\label{ch1_eq_arb_pol}
P\{X(2) > X(1)\} > 0.  
\end{equation}

Говорят, что {\it рынок безарбитражный}, если не существуют такие $k_S$, $k_C$, что выполняются условия~\eqref{ch1_eq_arb_0}\,--\,\eqref{ch1_eq_arb_pol}.
Найдём справедливую цену опциона~$C$ для безарбитражного рынка.

Условие \eqref{ch1_eq_arb_x2} означает, что при любой реализации цены акции в момент времени~$2$ имеет место неравенство $X(2) \geq X(1)$, т.е.
\begin{eqnarray*}
S(u - 1) k_S + (C_u - C) k_C & \ge & 0,  \\
S(d - 1) k_S + (C_d - C) k_C & \ge & 0.
\end{eqnarray*}
Условие \eqref{ch1_eq_arb_pol} означает, что как минимум одно из этих неравенств~--- строгое.
Система неравенств записывается в виде~$Ax \ge 0$, где
\[
    A = \left ( 
    \begin{array}{cc}
        S(u - 1) & C_u - C \\
        S(d - 1) & C_d - C 
    \end{array}
    \right ), \quad 
x = \left (
\begin{array}{c}
     k_S  \\
     k_C 
\end{array}
\right ).
\]
Поиск арбитража сводится к вопросу о существовании
такого $x$, что $Ax \ge 0$, $Ax \not= 0$. По лемме Фаркаша
об альтернативах, это условие эквивалентно существованию вектора $y \geq 0$, $y \not= 0$ такого, что
$y^T A = 0$. Учитывая то обстоятельство, что $u-1$ и $d-1$ имеют противоположные знаки, это условие просто эквивалентно вырожденности матрицы $A$. В этом случае можно подобрать решение $y = (1-d,u-1)^T > 0$. Следовательно, справедливую цену опциона~$C$ можно найти из условия
\begin{equation}
\label{ch1_eq_sc_syst1}
S(u - 1)(C_d - C) = S(d - 1)(C_u - C) \quad \Leftrightarrow \quad C = \frac{u-1}{u-d}C_d + \frac{1-d}{u-d}C_u.
\end{equation}

Таким образом, справедливую цену опциона для безарбитражного рынка можно вычислить по формуле
$$
C  = p C_u + (1 - p) C_d, 
$$
где $p = \frac{1 - d} {u - d} $.

Числа $p$ и $q = 1 - p$ задают мартингальную
вероятность, или мартингальную меру (см.~\cite{StochAn}). Если существует единственная
мартингальная мера, то рынок называется {\it полным}. На полном рынке
неизвестная цена опциона~$C$ в начальный момент определяется однозначно и может интерпретироваться как <<справедливая цена>>.

Далее для справки приводится арбитражная теорема, также известная по курсу В.Г.~Жадана <<Методы оптимизации>> как теорема Штимке (см.~\cite{Zhadan1}). 

\begin{teo}[арбитражная теорема] 
Пусть имеется $m$ ценных бумаг и $n$
возможных исходов, $r_{ij}$~--- платежи {\rm(}или убытки{\rm)} по $i$-й ценной бумаге,
если произошёл исход~$j$. Под стратегией инвестора понимается вектор $x \, \in \, \mathbb{R}^m$ {\rm(}$x_i$~--– количество ценных бумаг типа~$i${\rm)}. Имеет место одно
и только одно из приводимых ниже утверждений{\rm :}
\begin{enumerate}
\renewcommand\labelenumi{\textup{\arabic{enumi}.}}
    \item $\exists \, p \, > \, 0 \, : \, p^T R = 0$, где
    $R = \| r_{ij} \|_{i, \, j = 1}^{m, \, n}$ {\rm(}мартингальное соотношение{\rm)};
    \item $\exists \, x \, \in \, \mathbb{R}^m \, : \, Rx \ge 0, \, Rx \neq 0$ {\rm(}условие существования арбитража{\rm)}.
\end{enumerate}
\end{teo}

\subsection[Когда выгоден переход к двойственной задаче?]{Когда выгоден переход к двойственной задаче?}
Пусть система $Ax=b$ совместна. Рассмотрим задачу
$$
\min_{Ax = b} \frac{1}{2}\left \| x \right \|_2^2.
$$
В силу выполнения свойства сильной
двойственности можно сделать переход к двойственной задаче:
$$
\mathop {\min }\limits_{Ax=b} \frac{1}{2}\left\| x \right\|_2^2 =
\mathop 
{\min }\limits_x \mathop {\max }\limits_\lambda \left\{ {\frac{1}{2} \left \| 
x \right\|_2^2 +  (b-Ax)^T \lambda  } \right\} = 
$$
$$
=\mathop {\max }\limits_\lambda \mathop {\min }\limits_x \left\{ 
{\frac{1}{2}\left\| x \right\|_2^2 +  (b-Ax)^T \lambda  } \right\}=\mathop {\max }\limits_\lambda \left\{ 
{b^T \lambda  -\frac{1}{2}\left\| {A^T\lambda } 
\right\|_2^2 } \right\}.
$$
Поскольку система $Ax=b$ совместна,
по одной из теорем об альтернативах, рассмотренных 
в данном разделе, 
не существует такого~$\lambda $, что
выполняется $A^T \lambda =0$ и $b^T \lambda \,  > \, 0$.
Следовательно, двойственная задача имеет конечное решение. 
В самом деле, если существует
такой~$x$, что $Ax=b$, тогда для всех
$\lambda $: 
$$
(Ax)^T \lambda  = b^T \lambda. 
$$
Следовательно, $x^T A^T \lambda =b^T \lambda $.
Предположим, что существует такое $\lambda $, что 
$A^T \lambda = 0$ и $b^T \lambda \, > \, 0$. 
Но тогда мы пришли к противоречию:
$$
0 = x^T A^T\lambda = b^T \lambda \,  > \, 0,
$$
следовательно, верна только одна из альтернатив. 

Итак, вместо прямой оптимизационной задачи с ограничениями
можно решать двойственную к ней {\bf безусловную} (!) задачу
$$
\max_\lambda \left\{ 
{b^T \lambda  -\frac{1}{2}\left\| {A^T\lambda } 
\right\|_2^2 } \right\}
$$ 
и затем восстановить решение прямой задачи
по формуле  
\begin{equation}
\label{ch1_eq_prim_dual_fred}
x \left( \lambda 
\right) = A^T \lambda.
\end{equation}

Если матрица~$A$ не является матрицей полного ранга\footnote{
Матрицами полного ранга называются прямоугольные матрицы размера 
$m \, \times \, n$, ранг которых совпадает с минимальным из чисел
$m$ и $n$.},
то двойственная задача имеет не единственное решение
(все решения составляют аффинное многообразие).
Но для любого из них решения прямой задачи, вычисленные 
по формуле~\eqref{ch1_eq_prim_dual_fred},
совпадают.

\section{Формулировка оптимизационных задач} \label{sec:problem_representation}

Задача оптимизации в самой общей форме имеет вид
\[ \min_{x \in X}\, f(x),
\]
где $x$ --- искомая переменная, $X$ --- допустимое множество, $f$ --- целевая функция.

Какой алгоритм применять для решения данной задачи, зависит не только от того, выпукла ли она или нет, но и от алгоритмической доступности определяющих её элементов. Эти данные могут быть определены разными способами.

Функция цены $f$ может быть задана аналитическим выражением. В этом случае часто можно вычислить её значение и производные и работать с \nk{ними везде,} где они существуют, а также получить оценки на константы Липшица и другие параметры, определяющие сложность задачи и класс алгоритмов решения, которые нужно использовать.

Однако информация о функции может быть и более ограниченной. Например, её значение может определяться оптимальным значением другой задачи, параметризованной переменной $x$. Такая ситуация встречается, например, в минимаксных (седловых) задачах. В наименее информативном случае функция цены задаётся т.н. \emph{ чёрным ящиком} или \emph{оракулом}. Это программа, выдающая для каждого данного $x \in X$ значение функции, а может и (суб-)градиент, гессиан и т.д. Такая ситуация может возникнуть, когда значение функции цены является результатом сложного физического процесса, который можно моделировать только численно, или когда фирмы или учреждения, причастные к постановке задачи, не хотят или не способны выдавать слишком подробную информацию. При выборе алгоритма решения в таких случаях нужно учитывать стоимость вызова подпрограммы, считающей эти величины.

Вышесказанное также может относиться к описанию допустимого множества $X$. Последнее может быть задано аналитически, например, равенствами $f_i(x) = 0$, неравенствами $g_i(x) \leq 0$, матричными неравенствами $A_i(x) \succeq 0$, включениями $x_i \in \{0,1\}$, $x_i \in \mathbb Z$ и т.д. Однако информация о множестве может быть и более ограниченной. В~случае выпуклого множества минимальная информация, выдаваемая оракулом, состоит из утвердительного или отрицательного ответа на вопрос о принадлежности данной точки $x$ к множеству $X$, и в случае отрицательного ответа дополнительно из отделяющей плоскости.

В случае аналитического описания задачи её степень сложности определяется характером выражений, задающих функцию цены и допустимое множество. Они могут быть линейными, квадратичными, полиномиальными, тригонометрическими, дискретными (бинарными, целочисленными) и т.д. В зависимости от характера функции цены и ограничений можно выделять разные стандартные классы задач:
\begin{itemize}
\item линейные программы (linear program, LP): линейная функция цены и линейные ограничения;
\item линейные программы с целочисленными ограничениями (mixed-integer linear program, MILP): линейная функция цены, линейные и целочисленные или бинарные ограничения;
\item квадратичные программы (quadratic program, QP): квадратичная функция цены, линейные ограничения;
\item квадратично ограниченные квадратичные программы (quadratically constrained quadratic program, QCQP): квадратичная функция цены, квадратичные ограничения;
\item полиномиальные задачи: полиномиальная функция цены, полиномиальные ограничения;
\item геометрические программы (geometric program, GP): позиномиальная функция цены, позиномиальные ограничения.
\end{itemize}
Для разных классов задач существуют разные алгоритмы решения, имеющие известную сложность и поведение. В первом приближении задача эффективно решается, если она выпукла, т.е. функция цены $f$ и допустимое множество $X$ выпуклы. В п.~\ref{ch3_polynomial_problems} мы, однако, рассмотрим некоторые примеры, идущие вразрез с этим правилом.

Оптимизационные задачи, которые возникают в реальных приложениях, часто имеют необычную и не всегда каноническую форму. Специалистам, решающим такие задачи, необходимо иметь очень важный навык~--- способность узнавать и классифицировать задачи оптимизации, в том числе и <<замаскированные>>. Первый шаг --- это привести конкретную задачу к одной из известных стандартных форм, в которой её можно решить имеющимися в распоряжении методами. Ясно, что для удачного поиска представлений необходимо иметь некий <<каталог>> функций и множеств. Оказывается, что такой <<каталог>> в виде конического представления функций и множеств (см.~раздел \ref{ch1_sect_cone_prog}) существует и активно используется в популярном пакете для численного решения задач оптимизации~CVX (см. п.~\ref{ch1_sect_opt_package}).

Задачи некоторых классов можно преобразовать в задачи из других стандартных классов, поскольку одна и та же задача оптимизации может быть сформулирована разными эквивалентными способами:
\begin{itemize}
\item Функцию цены можно преобразовать в линейную, заменив
\[ \min_x\, f(x)
\]
на
\[ \min_{x,t}\, t\ : \quad t \geq f(x).
\]
\item Бинарные ограничения $x_i \in \{0,1\}$ или $x_i \in \{-1,+1\}$ можно заменить на квадратичные $x_i^2 = x_i$ или $x_i^2 = 1$ соответственно.
\item Полиномиальные ограничения можно преобразовать в квадратичные заменой мономов $x^{\alpha}$ на дополнительные переменные $x_{\alpha}$, удовлетворяющие ограничениям $x_{\alpha+\beta} = x_{\alpha}x_{\beta}$, $x_{e_i} = x_i$, где $\alpha = (\alpha_1, \, \dots, \, \alpha_n) \in \mathbb N^n$ --- мульти-индекс, а $x^{\alpha} := \prod_{i=1}^n x_i^{\alpha_i}$.
\end{itemize}

Таким образом, например, широкий класс полиномиальных задач сводится к QCQP, имеющей общий вид
\begin{equation} \label{ch1_qcqp}
\begin{split}
\min_{x \in \mathbb R^n}\ &\ \left(x^TA_0x + 2\langle b_0, \, x \rangle + c_0\right)\ : \\  x^TA_ix + 2\langle b_i,x \rangle + c_i \leq 0,\ &\ x^TA_jx + 2\langle b_j, \, x \rangle + c_j = 0.
\end{split}
\end{equation}
Здесь $A_0,A_i,A_j$ --- вещественные симметричные матрицы, $b_0$, $b_i$, $b_j$ --- векторы, а $c_0,c_i,c_j$ --- скаляры.

Если задача не эквивалентна задаче стандартного вида или для данного класса не существует эффективных алгоритмов решения, то её можно попытаться аппроксимировать, или \emph{релаксировать}. Релаксация задачи обычно состоит в опущении сложного условия или замене его на более простое, например замена области определения $\{0,1\}$ бинарной переменной на непрерывный интервал $[0,1]$. Иногда упрощение производят после некоторых преобразований исходной задачи, например, после перехода к двойственной. Релаксация имеет смысл, если при этом множество допустимых точек либо сужается, либо расширяется, с тем чтобы оптимальное значение релаксации представляло собой границу (верхнюю или нижнюю) для значения исходной задачи. В разделе \ref{ch1_relaxations} мы рассмотрим примеры, в которых можно оценить ошибку, сделанную при релаксации.

\section{Коническое программирование}\label{ch1_sect_cone_prog}

\subsection{Конические программы}

В этом разделе мы рассмотрим широкий класс выпуклых задач оптимизации в стандартной форме. Эта форма, т.н. \emph{коническая программа}, имеет общий вид
\begin{equation} \label{ch1_conic_indirect}
\min_{x \in K \times \mathbb R^{n'}} \langle c,x \rangle\ : \quad Ax = b,
\end{equation}
где вектор $c \, \in \, \mathbb{R}^n$, $A$~--- матрица размера~$m \, \times \, n$, вектор~$b \, \in \, \mathbb{R}^m$, $K \subset \mathbb R^{n-n'}$ --- замкнутый выпуклый конус с непустой внутренностью, не содержащий прямой. Такие конусы называются \emph{регулярными} или \emph{правильными}. Помимо условия принадлежности к конусу, на $x$ налагаются только линейные условия типа равенства. Функция цены также является линейной. Любую выпуклую задачу оптимизации можно переписать в виде задачи конического программирования.

Таким образом, в задаче конического программирования в качестве множества допустимых точек выступает пересечение регулярного конуса с аффинным подпространством. В формулировке \eqref{ch1_conic_indirect} переменная $x$ (или часть её если $n' > 0$) параметризует точки конуса, и на неё дополнительно накладываются линейные условия, задающие аффинное подпространство. 

Коническую программу над регулярным конусом $K$ можно также записать в виде
\begin{equation}
\label{ch1_eq_con_prog}
\min_{Ax \, \succeq_K \, b}  \langle c,x \rangle.
\end{equation}
Здесь матрица $A$ и векторы $b$, $c$, $x$ имеют те же размерности, что и выше, но конус $K$ определён в пространстве $\mathbb R^m$. В этой формулировке переменная $x$ параметризует непосредственно допустимое множество. В зависимости от ситуации удобнее использовать одну или другую формулировку. Обе формулировки \eqref{ch1_conic_indirect} и \eqref{ch1_eq_con_prog} легко переводятся друг в друга, но при этом размерность и коэффициенты матрицы $A$ и векторов $b,c$ меняются.

Очевидно сложность задачи конического программирования в основном зависит от конуса $K$. Последний может задаваться разными способами, от аналитического описания до  чёрного ящика или через решение вспомогательной задачи. Основным классом конусов, для которых конические программы эффективно решаются, является класс \emph{симметричных} конусов.

\begin{defin} Регулярный конус $K \subset \mathbb R^n$ называется \emph{симметричным}, если выполняются следующие условия:

а) для любых двух внутренних точек $x,y \in K^o$ существует линейный автоморфизм $A$ конуса $K$, который переводит $x$ в $y$ (однородность), и

б) при подходящем выборе базиса в $\mathbb R^n$ и соответствующего двойственного базиса в $\mathbb R_n = (\mathbb R^n)^*$ двойственный конус $K^*$ совпадает с $K$.
\end{defin}

\medskip

Симметрические конусы полностью классифицированы. В оптимизации важную роль играют следующие семейства симметричных конусов:
\begin{itemize}
\item положительный ортант $\mathbb R_+^n$;
\item конус Лоренца $L^n$ и произведения вида $L^{n_1} \times \dots \times L^{n_k}$;
\item вещественный матричный конус $\mathbb{S}_+^n = \{ A \in \mathbb{S}^n \mid A \succeq 0 \}$;
\item комплексно-эрмитовый матричный конус $\mathbb{H}_+^n = \{ A \in \mathbb{H}^n \mid A \succeq 0 \}$.
\end{itemize}
Здесь $\mathbb{S}^n$ и $\mathbb{H}^n$ --- пространства вещественных симметричных и комплексных эрмитовых матриц размера $n \times n$ соответственно. Отметим, что последнее рассматривается как векторное пространство над $\mathbb R$.

Конические программы над симметричными конусами вышеперечисленных семейств эквивалентны задачам линейного программирования (LP), конично-квадратичного программирования (conic-quadratic program, CQP) или, эквивалентно, конического программирования второго порядка (second-order cone program, SOCP), полуопределённого программирования (semi-definite program, SDP) с вещественными и комплексными коэффициентами соответственно.

\ag{Конические программы широко используются в исследовании операций, инженерии \cite{BenTalNemirovski01}, анализе динамических систем \cite{Boyd_rob}, управлении \cite{RaoWrightRawlings98}, маршрутизации \cite{XiaoJohanssonBoyd04}, машинном обучении \cite{FerrisMunson03,Lanckriet_etal04}, методе главных компонент \cite{Aspremont_etal07}, восстановлении сигнала \cite{AspremontGhaoui11}, восстановлении фазы \cite{WaldspurgerAspremontMallat15} и~др. Большое количество невыпуклых задач можно аппроксимировать выпуклыми задачами, представляемыми в виде конической программы, например, такие комбинаторные проблемы, как максимальный разрез графа (MaxCut, см.~п.~\ref{ch1_relaxations}), выполнимость булевых формул (SAT) \cite{GoemansWilliamson}, кликовое число графа (см.~п.~\ref{subs:MaxClique}) \cite{Alizadeh95}, контактное число (kissing number) \cite{BachocVallentin08}, квадратичная задача о назначениях \cite{deKlerkSotirov10}, а также линейные программы с целочисленными ограничениями \cite{Alizadeh95}, квадратичные программы с квадратичными ограничениями \cite{Nesterov98} и задачи полиномиальной оптимизации (см.~п.~\ref{ch3_polynomial_problems}) \cite{NesterovSOS,ParriloThesis,LasserreMoments}.}

Рассмотрим более подробно класс конично-квадратичных программ. Напомним, что для $n \ge 2$ конус Лоренца определяется как (см.\ также с.~\pageref{ch1_ex_loren})
$$
L^n =
\left \{ x = (x_1, \, \ldots, \, x_{n-1}, \, x_{n})^T \, \in \, \mathbb{R}^{n} \, \left| \, x_{n} \, \ge \, \sqrt{x_1^2 + x_2^2 + \ldots + x_{n-1}^2} \right. \, \right \}.
$$
Итак, в конично-квадратичных программах имеем
\begin{eqnarray}
K = L^{n_1} \, \times \, L^{n_2} \, \times \, \ldots \, \times \, L^{n_k} = \nonumber \\
= \left \{
y =
\left ( \left.
\begin{array}{c}
y[1] \\
y[2] \\
\vdots \\
y[k]
\end{array}
\right )
\, \right | \, y[i] \, \in \, L^{n_i}, \quad i = 1, \, \ldots, \, k
\right \}.
\label{ch1_eq_con2}
\end{eqnarray}
Иными словами, задача конического программирования второго порядка является задачей минимизации линейной функции с конечным числом ограничений вида
$$
A_i x - b_i \succeq_{L^{n_i}} 0, \quad i = 1, \, \ldots, \, k,
$$
где
$A = [A_1, \, A_2, \, \ldots, \, A_k]^T$, $b = [b_1, \, b_2, \, \ldots, \, b_k]^T$~--- разбиения на блоки, соответствующие
разбиению матрицы~$y$ в~\eqref{ch1_eq_con2}.
Если составить матрицу~$[A_i, \, b_i]$ и разбить её на блоки,
выделив последнюю компоненту:
$$
[A_i, \, b_i] =
\left [
\begin{array}{cc}
     D_i &  d_i \\
     p_i^T & q_i
\end{array}
\right ],
$$
где $D_i$ имеет размер $(m_i - 1) \, \times \, \dim x$,
а затем вспомнить смысл отношения $\succeq_{L^m}$,
то задачу~\eqref{ch1_eq_con_prog} с~$K$ из \eqref{ch1_eq_con2} можно  переписать так:
$$
\min_{\| D_i x - d_i \|_2 \, \le \, p_i^T x - q_i, \, i = 1, \, \ldots, \, k} c^T x.
$$

Преобразование исходной задачи оптимизации к конической программе над симметричным конусом гарантирует её эффективное решение при условии, что размерность задачи при этом не сильно возрастёт. Если не удаётся произвести эквивалентное преобразование, то можно попытаться \emph{аппроксимировать} задачу конической программой из одного из вышеперечисленных стандартных классов. В следующем разделе мы рассмотрим в качестве примера, как можно подойти к задачам класса QCQP с позиций конического программирования.

\subsection{Преобразование и релаксации квадратичных задач}

Квадратичная задача \eqref{ch1_qcqp} класса QCQP общего вида трудно разрешима. Однако если матрицы $A_0,A_i$ положительно определены, а $A_j = 0$, то задача выпукла и сводится далее к конично-квадратичной программе, которая эффективно решается. В противном случае задачу QCQP можно аппроксимировать выпуклой задачей. Несколько вариантов такой аппроксимации мы увидим ниже в п.~\ref{ch1_relaxations}.

Рассмотрим сначала выпуклый случай. Условия равенства линейные и не требуют дальнейших преобразований. Функцию цены можно заменить на линейную функцию $t$ посредством введения новой переменной и нового ограничения $t \geq x^TA_0x + 2\langle b_0,x \rangle + c_0$. Это ограничение~--- выпуклое квадратичное неравенство, имеющее такой же вид, как и исходные условия типа неравенства. Осталось рассмотреть эти неравенства. Нетрудно проверить, что для неотрицательно определённой матрицы $A$ условие $x^TAx + b^Tx + c \leq 0$ эквивалентно включению
\[ \left(\frac{1+b^Tx+c}{2},A^{1/2}x,\frac{1-b^Tx-c}{2}\right) \in L^{n+2}.
\]
Здесь $A^{1/2}$ обозначает матричный корень, т.е. неотрицательно определённую симметричную матрицу такую, что $A^{1/2} \cdot A^{1/2} = A$. Таким образом, выпуклая задача QCQP переведена в квадратично-коническую программу.

Рассмотрим теперь общую задачу QCQP \eqref{ch1_qcqp}. Введём дополнительную матричную переменную $X$ размера $(n+1) \times (n+1)$ и матрицы коэффициентов:
\[ C = \begin{pmatrix} c_0 & b_0^T \\ b_0 & A_0 \end{pmatrix}, \quad F_i = \begin{pmatrix} c_i & b_i^T \\ b_i & A_i \end{pmatrix}, \quad G_j = \begin{pmatrix} c_j & b_j^T \\ b_j & A_j \end{pmatrix}.
\]
Тогда задачу можно переписать в виде
\[ \min_{x,X} \, \langle C,X \rangle\ :\quad \langle F_i,X \rangle \leq 0,\quad \langle G_j,X \rangle = 0,\quad X = \begin{pmatrix} 1 & x^T \\ x & xx^T \end{pmatrix}.
\]
В этой формулировке функция цены и все ограничения линейны, кроме последнего. Однако последнее ограничение можно переписать в виде
\[ X_{00} = 1,\quad X_{0k} = x_k,\ k = 1, \, \dots, \, n, \quad X \succeq 0, \quad \rk \, X = 1.
\]
Все ограничения линейны или полуопределённые, кроме последнего. Это ограничение на ранг невыпукло и делает задачу сложной.

Стандартная полуопределённая релаксация задачи QCQP состоит в опущении ограничения на ранг. Так как допустимое множество увеличивается, оптимальное значение релаксации ограничивает оптимальное значение исходной задачи снизу. Если среди решений $(x,X)$ полуопределённой релаксации найдётся такое, что $\rk\,X = 1$, то релаксация точна и компонента $x$ является решением исходной задачи. Этот случай имеет место, если в исходной задаче имеется только одно ограничение\footnote{Другие случаи, в которых полуопределённая релаксация квадратичной задачи точна, описаны в работе \cite{WangKilincKarzan21}.}. Для доказательства нам понадобится следующий результат \cite{Dines41}.

\begin{lemma} \label{lem:Dines}
Пусть $A,B$ --- произвольные симметричные матрицы размера $n \times n$. Тогда их \emph{численный образ} {\rm(}numerical range{\rm)}, задающийся множеством
\[ W(A,B) = \left\{ \left(x^TAx,x^TBx\right) \in \mathbb R^2 \mid x \in \mathbb R^n \right\},
\]
является выпуклым.
\end{lemma}

\begin{corr}
Если в задаче QCQP имеется только одно ограничение, то её полуопределённая релаксация точна.
\end{corr}

\begin{proof}
Пусть ограничение имеет вид $\langle F,X \rangle \leq 0$ или $\langle F,X \rangle = 0$.

Пусть $X \in \mathbb{S}_+^{n+1}$ неотрицательно определена. Тогда найдутся векторы $r_1, \, \dots, \, r_k \in \mathbb R^{n+1}$ такие, что $X = \frac{1}{k}\sum_{j=1}^k r_jr_j^T$. Поэтому
\[ (\langle C,X \rangle,\langle F,X \rangle) = \frac{1}{k}\sum_{j=1}^k \left(r_j^TCr_j,r_j^TFr_j\right).
\]
Но пары $\left(r_j^TCr_j,r_j^TFr_j\right)$ являются элементами численного образа $W(C,F)$ для всех $j = 1, \, \dots, \, k$. Так как этот численный образ выпуклый, пара $(\langle C,X \rangle,\langle F,X \rangle)$ также является его элементом. Но тогда существует матрица $\tilde X = rr^T \in \mathbb{S}_+^{n+1}$ ранга 1 такая, что $(\langle C,X \rangle,\langle F,X \rangle) = \left(\langle C,\tilde X \rangle,\langle F,\tilde X \rangle\right)$.

В частности, если $X$ --- решение релаксации, то найдётся другое её решение, имеющее ранг 1. Поэтому релаксация точна. 
\end{proof}

\subsection{Коническое представление множеств и функций} \label{ch1_sect_cone_repr}

Встречающиеся на практике задачи оптимизации изначально не имеют форму \eqref{ch1_conic_indirect} или \eqref{ch1_eq_con_prog} конической программы, к которой можно применить стандартизированные алгоритмы решения. Однако процесс трансформации исходной задачи в коническую программу можно в некоторой степени автоматизировать. Для этого сперва надо понимать, какие типы ограничений и функций цены позволяют переписать проблему в виде конической программы над тем или иным конусом.

Введём следующее определение.

\begin{defin}
Множество $X \subset \mathbb R^n$ называется \emph{представимым} через конус $K \subset \mathbb R^m$, если найдутся матрицы $A_1,A_2$, вектор $b$, и число $l$ такие, что условие включения $x \in X$ эквивалентно условию
\begin{equation} \label{basic_conic_constraint} 
\exists\ y: \ (x,y) \in K \times \mathbb R^l,\ A_1x + A_2y = b.
\end{equation}
Множество $X$ называется \emph{представимым} над семейством конусов ${\cal K}$, если существует конус $K \in {\cal K}$ такой, что $X$ представимо через $K$.
\end{defin}

Другими словами, множество $X$ представимо через $K$, если оно является линейной проекцией некоторого аффинного сечения произведения $K \times \mathbb R^l$. Здесь проекция задаётся $(x,y) \mapsto x$, и точка $(x,y) \in K \times \mathbb R^l$ называется \emph{поднятием} точки $x \in X$. Дополнительные переменные $y$ называются переменными поднятия (lifting variables). Обратим внимание на то, что в представлении \eqref{basic_conic_constraint} не предполагается, что разбиение пары $(x,y)$ соответствует разбиению произведения $K \times \mathbb R^l$ на множители.

В дальнейшем мы будем работать с семействами ${\cal K}$ конусов, замкнутых по отношению к взятию прямых произведений конечного числа конусов. В частности, такое семейство задаётся множеством $\{ \mathbb R_+^n \mid n \in \mathbb N_+\}$ положительных ортантов всех размерностей. Также мы рассмотрим множества, представимые по отношению к семейству, порождённому конусами Лоренца $L^n$ или матричными конусами $\mathbb S_+^n$.

Ясно, что конечное число конических условий \eqref{basic_conic_constraint} над разными конусами $K_1, \, \dots, \, K_r$ можно объединить в одно условие такого вида над прямым произведением $K_1 \times \dots \times K_r$. Поэтому в коническую программу можно интегрировать любое конечное число ограничений вида $x_i \in X_i$, если все встречающиеся множества $X_i$ представимы через подходящие конусы.

Приведём ряд операций, позволяющих построить одни представимые над семейством ${\cal K}$ множества из других, более простых таких множеств.

\emph{Конечные пересечения.} Из вышесказанного следует, что если множества $X_i$ представимы над некоторым семейством конусов ${\cal K}$, то над ${\cal K}$ также представимо их пересечение $x \in \bigcap_{i \in I}X_i$, если индексное множество $I$ конечно.

\emph{Аффинные сечения.} Если $X \in \mathbb R^n$ представимо над ${\cal K}$, то его пересечение $X \cap A$ с аффинным подпространством $A \subset \mathbb R^n$ также представимо над ${\cal K}$. Представление пересечения строится добавлением линейных условий типа равенства $x \in A$ в представление множества $X$.

\emph{Аффинный образ.} Если $X \in \mathbb R^n$ представимо над ${\cal K}$, а $L: \mathbb R^n \to \mathbb R^m$~--- аффинное отображение, то образ $X' = L[X]$ также представим над ${\cal K}$. А именно, пусть представление $X$ задано \eqref{basic_conic_constraint}. Тогда условие $z \in X'$ эквивалентно
\[ \exists\ x,y: \ (x,y,z) \in K \times \mathbb R^{l+m},\ A_1x + A_2y = b,\ Lx = z.
\]

\emph{Аффинный прообраз.} Если $X \in \mathbb R^n$ представимо над ${\cal K}$, а $L: \mathbb R^m \to \mathbb R^n$~--- инъективное аффинное отображение, то прообраз $X' = L^{-1}[X]$ также представим над ${\cal K}$. Пусть $X$ представлено в виде \eqref{basic_conic_constraint}, где $y$ имеет размерность $s$. Определим аффинное отображение $L' = L \times Id_{\mathbb R^s}: \mathbb R^m \times \mathbb R^s \to \mathbb R^n \times \mathbb R^s$. Тогда прообраз множества $K \times \mathbb R^l \subset R^n \times \mathbb R^s$ имеет вид $K \times \mathbb R^{l'}$ для $l' = l + m - n$. Условие $x \in X'$ запишется в виде
\[ \exists\ y: \ (x,y) \in K \times \mathbb R^{l'},\ A_1Lx + A_2y = b.
\]

До сих пор мы описывали конические представления множеств $X$. Эти представления можно использовать для включения ограничений $x \in X$ в коническую программу над рассматриваемыми конусами. Точно так же можно говорить о коническом представлении \emph{функций}. В п.~\ref{sec:problem_representation} мы видели, как преобразовать задачу минимизации произвольной выпуклой функции $f$ в задачу минимизации линейной функции над надграфиком $f$. Таким образом, коническое представление функций можно свести к представлению множеств.

Мы будем называть функцию $f$ представимой через конус $K$, если её надграфик $X$ представим через $K$. Соответственно определим представление функции над семейством ${\cal K}$ конусов. Аналогично операциям над множествами имеем следующую операцию над функциями.

\emph{Максимум семейства функций.} Если функции $f_1, \, \dots, \, f_l$ представимы над семейством конусов ${\cal K}$, то их максимум $f(x) = \max_{i = 1, \, \dots, \, l} f_i(x)$ также представим над ${\cal K}$. Это следует из того, что надграфик $f$ является пересечением надграфиков функций $f_i$.

Чем шире семейство конусов, которым мы разрешаем фигурировать в конической программе, тем богаче класс представимых ограничений и функций. Ниже мы рассмотрим некоторые множества, представимые над произведениями конусов Лоренца или матричных конусов. Для простоты изложения мы приведём непосредственно ограничения, задающие эти множества.

\medskip

Проблема с линейными ограничениями и линейной функцией цены представляется в виде линейной программы.

\medskip

Коническая программа второго порядка способна представить гораздо более широкий спектр проблем. Приведём примеры ограничений, которые можно переписать как условия включения в конус Лоренца.

\begin{itemize}
\item $||x||_2^2 \leq t$, $x \in \mathbb R^n$ эквивалентно условию $\left(x,\frac{t-1}{2},\frac{t+1}{2}\right) \in L^{n+2}$;
\item $\frac{||x||_2^2}{s} \leq t$, $s \geq 0$, $x \in \mathbb R^n$ выражается в виде $\left(x,\frac{t-s}{2},\frac{t+s}{2}\right) \in L^{n+2}$;
\item $ts \geq 1$, $t,s > 0$ эквивалентно условию $\left(1,\frac{t-s}{2},\frac{t+s}{2}\right) \in L^3$;
\item $x^TAx + b^Tx + c \leq t$ при $A \succeq 0$ эквивалентно условию $\left(A^{1/2}x,\frac{t-b^Tx-c-1}{2},\frac{t-b^Tx-c+1}{2}\right) \in L^{n+2}$;
\item $|t| \leq \sqrt{x_1x_2}$, $x_1,x_2 \geq 0$ эквивалентно условию $\left(t,\frac{x_1-x_2}{2},\frac{x_1+x_2}{2}\right) \in L^3$;
\item $t \leq \sqrt{x_1x_2}$, $x_1,x_2 \geq 0$ эквивалентно условиям $t \leq s$, $s \geq 0$, $\left(s,\frac{x_1-x_2}{2},\frac{x_1+x_2}{2}\right) \in L^3$.
\end{itemize}

\medskip

Применим эти соотношения для построения конического представления множеств~\cite{Nemir_ro_tr}.

\begin{example}
Задано множество
$$
X = \left\{x = (x_1, \, x_2, \, x_3, \, x_4) \, \in \, \mathbb{R}_+^4 \, | \, x_1 x_2 x_3 x_4 \, \ge \, 1  \right\}.
$$
Требуется найти его коническое представление.

{\bf Решение.}
Построим поднятие множества $X$, добавив три дополнительные переменные
$u_1$, $u_2$, $u_3$. Представим множество $X$
как проекцию пересечения конусов Лоренца:
$$
X = \left \{ x \, \in \, \mathbb{R}_+^4 \, \left| \,
\exists \, u \, \in \, \mathbb{R}^3 \, :
\left \{
\begin{array}{c}
0 \, \le \,  u_1 \,  \le \, \sqrt{x_1 x_2} \\
0 \, \le \,  u_2 \,  \le \, \sqrt{x_3 x_4} \\
1 \, \le \,  u_3 \,  \le \, \sqrt{u_1 u_2} \\
\end{array}
\right \} \right.
\right \} =
$$
$$
= \left \{
x \, \in \, \mathbb{R}_+^4 \, \left| \,
\exists \, u \, \in \, \mathbb{R}^3 \, :
\left \{
\begin{array}{ccc}
\lbrack x_1; \, x_2; \, x_3; \, x_4; \, u_1; \, u_2; \, u_3 - 1 \rbrack & \in & \mathbb{R}^7_+  \\
\lbrack 2 u_1; \, x_1 - x_2; \, x_1 + x_2 \rbrack & \in & L^3 \\
\lbrack 2 u_2; \, x_3 - x_4; \, x_3 + x_4 \rbrack & \in & L^3 \\
\lbrack 2 u_3; \, u_1 - u_2; \, u_1 + u_2 \rbrack & \in & L^3
\end{array}
\right \} \right.
\right \}.
$$
\end{example}

\begin{exercise}
Проверьте с помощью определения конуса Лоренца, что данное коническое представление действительно эквивалентно исходному множеству~$X$.
\end{exercise}

Обобщая идею, на которой основан пример, можно построить конично-квадратичное представление \nk{подграфика}
\[ X = \left\{ (t,x) \in \mathbb R \times \mathbb R_+^{2^l} \mid t \leq \sqrt[2^l]{x_1x_2\dots x_{2^l}} \right\}
\]
\label{geometric_mean}

\noindent геометрического среднего $2^l$ переменных. Например, для четырёх переменных получаем представление
\[ X = \left \{
(t,x) \, \in \, \mathbb R \times \mathbb{R}^4 \, \left| \,
\exists \, u \, \in \, \mathbb{R}^2 \, :
\left \{
\begin{array}{ccc}
\lbrack x_1; \, x_2; \, x_3; \, x_4; \, u_1; \, u_2 \rbrack & \in & \mathbb{R}^6_+  \\
\lbrack 2 u_1; \, x_1 - x_2; \, x_1 + x_2 \rbrack & \in & L^3 \\
\lbrack 2 u_2; \, x_3 - x_4; \, x_3 + x_4 \rbrack & \in & L^3 \\
\lbrack 2 t; \, u_1 - u_2; \, u_1 + u_2 \rbrack & \in & L^3
\end{array}
\right \} \right.
\right \}.
\]


\medskip

Перейдём к полуопределённым программам. С помощью этого класса можно описать ещё гораздо более широкий спектр проблем. Приведём примеры ограничений, которые можно переписать как условия включения в конус неотрицательно определённых матриц.

\begin{itemize}
\item $\lambda_{\max}(X) \leq t$ эквивалентно условию $tI - X \succeq 0$;
\item $||X||_{\infty} \leq t$ эквивалентно условию $-tI \preceq X \preceq tI$ для симметричных $X$;
\item $\sum_{j=1}^k \lambda_j \leq t$, где $\lambda_1, \, \dots, \, \lambda_n$~--- упорядоченные по низходящей собственные значения $X$, эквивалентно условиям $t \geq ks + \tr\,Z$, $Z \succeq 0$, $Z + sI \succeq X$;
\item $A \succeq BC^{\dagger}B^T$, $C \succeq 0$ эквивалентно условию $\begin{pmatrix} A & B \\ B^T & C \end{pmatrix} \succeq 0$;
\item $||A||_{\infty} \leq t$ эквивалентно условию $\begin{pmatrix} tI & A \\ A^T & tI \end{pmatrix} \succeq 0$ для матриц $A$ общего вида;
\item $(AXB)(AXB)^T + CXD + (CXD)^T + E \preceq Y$ эквивалентно условию $\begin{pmatrix} I & (AXB)^T \\ AXB & Y - E - CXD - (CXD)^T \end{pmatrix} \succeq 0$ (здесь $X,Y$~--- переменные, а $A, \, \dots, \, E$~--- параметры задачи);
\item $X \succeq 0$, $0 \leq \eta \leq \lambda(X)$ эквивалентно условиям
\[ \begin{pmatrix} X & \Delta \\ \Delta^T & \diag\eta \end{pmatrix} \succeq 0,\ \diag\Delta = \eta,\ \Delta_{ij} = 0\ \forall\ i < j,
\]
где $\lambda(X)$ --- спектр матрицы $X$;
\item $x^TAx \geq 0$ для всех $x$ таких, что $x^TBx \geq 0$, и существует $x_0$ такое, что $x_0^TBx_0 > 0$ эквивалентно условиям $A - \lambda B \succeq 0$ и $\lambda \geq 0$ ($S$-лемма).
\end{itemize}
В последнем пункте матрица $A$ и скаляр $\lambda$ являются переменными, а матрица $B$~--- параметром задачи.

Доказательство $S$-леммы основано на лемме \ref{lem:Dines}.

Если $A - \lambda B \succeq 0$ для некоторого $\lambda \geq 0$, то $x^TAx \geq \lambda x^TBx$ для всех $x \in \mathbb R^n$. Это доказывает одно направление эквивалентности.

\begin{figure}[h!]
\begin{center}
\includegraphics[width=8.58cm,height=4.38cm]{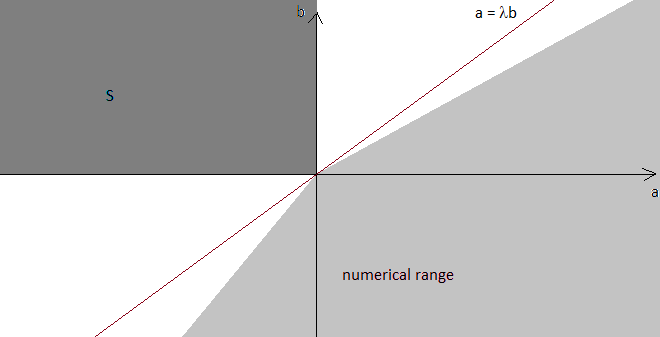}
\end{center}
\caption{Отделение множества $S$ от численного образа}
\label{ch1_fig:slemm}
\end{figure}

Допустим теперь, что $x^TAx \geq 0$ для всех векторов $x$, удовлетворяющих условию $x^TBx \geq 0$. Тогда множество $S = \{ (a,b) \,|\, a < 0,\ b \geq 0\}$ имеет пустое пересечение с численным образом пары $(A,B)$. В силу выпуклости численного образа существует одномерное подпространство $\mathbb R^2$, отделяющее численный образ от $S$. Но численный образ содержит точку $\left(\hat a,\hat b\right)$, удовлетворяющую $\hat a \geq 0$, $\hat b > 0$. Поэтому отделяющее подпространство не может быть горизонтальной осью. Отсюда следует, что оно задаётся уравнением $a = \lambda b$ для некоторого $\lambda \geq 0$ (см. рис.~\ref{ch1_fig:slemm}). Но тогда $a \geq \lambda b$ для всех точек $(a,b)$ в численном образе, что эквивалентно условию $A - \lambda B \succeq 0$. Доказательство $S$-леммы завершено.

\medskip

Применим вышеперечисленные соотношения для построения конического представления множеств~\cite{Nemir_ro_tr}.

\pagebreak

{

\renewcommand{\baselinestretch}{0.94}
\selectfont

\begin{example}
Единичный шар ядерной нормы в пространстве матриц размера $m \times n$. Напомним, что
ядерная норма (nuclear norm, trace norm) матрицы $A \, \in \, \mathbb R^{m \times n}$ задаётся суммой её сингулярных значений:
$$
\| A \|_* = \sum_{i = 1}^{\min \{ m, \, n \}} \sigma_i(A).
$$
Требуется найти коническое представление этого шара.

{\bf Решение.}
$$
\begin{array}{r}
X =  \biggl \{
A \, \in \, \mathbb{R}^{m \times n} \, : \,
\exists \, u = (U \, \in \, \mathbb{S}^m, \, V \, \in \, \mathbb{S}^n) : \\
\text{Tr}(U) + \text{Tr}(V) \, \le \, 2, \quad
\begin{pmatrix}
     U & A \\
     A^T & V
\end{pmatrix} \succeq 0
\biggr \},
\end{array}
$$
где $\mathbb{S}^m$~--- множество вещественных симметричных матриц размера~$m \times m$, $\text{Tr}(A)$~--- след квадратной матрицы~$A$, т.е. сумма диагональных элементов~$A$. 
\end{example}


Более широкий список полуопределённо представимых ограничений можно найти в лекциях А.~Бен-Таля и А.\,С.~Немировского~\cite{Ben-Tal}. В работе~\cite{JuditskyNemirovski21} разработаны методы, с помощью которых некоторые структурированные классы седловых задач и вариационных неравенств сводятся к коническим программам. В этой статье также перечислен ряд операций, с~помощью которых можно строить более сложные такие задачи исходя из более простых.


Отметим, что знать о возможности конического представления функций и множеств очень важно, поскольку для подавляющего большинства приложений, в которых возникают задачи выпуклой оптимизации, в том числе существенно нелинейной, существует возможность представить целевые функции и функции-ограничения подобным образом.  А представленная таким образом задача может быть численно решена стандартными пакетами, например, с помощью оптимизационного пакета CVX.

}

Для многих классов задач выпуклой оптимизации уже удалось найти такие конические представления функций и множеств. Сначала был построен достаточно большой набор
стандартных (библиотечных) выпуклых функций с известными коническими представлениями.
Затем были определены правила сочетания, которые позволяют по известным функциям получать новые\footnote{Подробное описание данных правил конического представления можно найти в~\cite[параграфы~2.3, 2.5, 3.2]{Ben-Tal}.}, и
написана программа, которая организует разумный перебор этих правил для поиска представления новых функций, не присутствующих в библиотеке.
Детали см. в параграфе~\ref{ch1_sect_opt_package} про пакет CVX, \ag{автоматизирующий этот поиск}.

В следующем параграфе мы рассмотрим вопрос представимости выпуклой задачи в виде стандартной конической программы с более общей, теоретической точки зрения.

\subsection{Преобразования конических программ} \label{ch1_lifting}

В этом разделе мы рассмотрим, каким образом можно переписать коническую программу над некоторым конусом $K \subset V_n$ через другую коническую программу над конусом $K' \subset V_m$, где $V_n,V_m$ --- вещественные векторные пространства соответствующих размерностей.

Рассмотрим коническую программу
\begin{equation} \label{original}
\min_{x \in K}\, \langle c,x \rangle\ :\quad x \in A,
\end{equation}
где $A \subset V_n$ --- аффинное подпространство.

Пусть теперь конус $K$ представляется в виде сечения другого конуса $K'$, т.е. существует линейная инъекция $L: V_n \to V_m$ такая, что $L[K] = K' \cap\linebreak \cap\, L[V_n]$. Тогда исходная коническая прог\-рамма \eqref{original} эквивалентна программе
\[ \min_{y \in K'}\, \langle c',y \rangle\ :\quad y \in L[A],
\]
где $c' \in V_m^*$ --- любой линейный функционал, удовлетворяющий соотношению $\langle c',L(x) \rangle = \langle c,x \rangle$ для всех $x \in V_n$. \ag{Эквивалентно, $c'$ должен удовлетворять условию} $L^{\dagger}(c') = c$. Здесь $L^{\dagger}: V_m^* \to V_n^*$ --- сопряжённое к $L$ линейное отображение, определённое соотношением
\[ \left\langle L^{\dagger}(z),y \right\rangle = \langle z,L(y) \rangle\qquad \forall\ y \in V_n,\ z \in V_m^*.
\]

Часто выражения, задающие конус $K$, более или менее явным образом указывают на расширение $K'$, имеющее более простую структуру, чем~$K$. Например, любой полиэдральный конус является сечением ортанта. Линейное сечение матричного конуса называется \emph{спектраэдральным конусом}. Заметим также, что произвольное аффинное сечение матричного конуса называется \emph{спектраэдром}.

\medskip

Менее очевидным является представление $K \in V_n$ линейной проекцией другого конуса $K' \in V_m$. Такое представление задаётся линейной сюръекцией $\Pi: V_m \to V_n$ такой, что $\Pi[K'] = K$. Пусть задана такая сюръекция, также называющаяся \emph{поднятием} (lifting). Тогда исходная коническая программа \eqref{original} эквивалентна программе
\[ \min_{y \in K'}\, \left\langle \Pi^{\dagger}(c),y \right\rangle\ : y \in \Pi^{-1}[A].
\]

Часто поднятие можно реализовать вводом дополнительных переменных. В этом случае $V_m = V_n \times V$, а $\Pi$ определяется как проекция $(x,y) \mapsto x$. Коническая программа над $K'$ принимает вид
\[ \min_{(x,y) \in K'}\, \langle (c,0),(x,y) \rangle\ :\quad (x,y) \in A \times V.
\]

\medskip

Две приведённые выше конструкции двойственны друг к другу. А именно, конус $K$ тогда и только тогда является сечением конуса $K'$, когда двойственный конус $K^*$ является проекцией двойственного конуса $(K')^*$.

\medskip

В общем случае можно комбинировать оба подхода.

\begin{defin} \label{def:representation}
Конус $K \subset V_n$ называется \emph{представимым} через конус $K' \subset V_m$, если существуют векторное пространство $V$, конус~${\tilde K \subset V}$, линейная инъекция $L: V \to V_m$ и линейная сюръекция $\Pi: V \to V_n$ такие, что $\Pi[\tilde K] = K$, $L[\tilde K] = L[V] \cap K'$.
\end{defin}

Иными словами, включение $x \in K$ характеризуется существованием элемента $y \in \tilde K$ такого, что $\Pi(y) = x$, $L(y) \in K'$. Конус $\tilde K$ изоморфен сечению конуса $K'$, а конус $K$ является его проекцией. Коническая программа \eqref{original} перепишется в виде
\[ \min_{z \in K'} \langle c',z \rangle\ :\quad z \in L[\Pi^{-1}[A]].
\]
Здесь $c'$ --- любой линейный функционал, удовлетворяющий соотношению $\langle c',L(y) \rangle = \left\langle \Pi^{\dagger}(c),y \right\rangle$ для всех $y \in V$, или эквивалентно $L^{\dagger}(c') = \Pi^{\dagger}(c)$.

С точностью до линейного изоморфизма конус $K$ можно восстановить из $K'$, если заданы два подпространства $V' = L[V] \subset V_m$ и $V_0 = L[\ker\,\Pi] \subset V'$. Тогда пространство $V_n$ можно отождествить с фактор-пространством $V'/V_0$, а конус $K$~--- с множеством 
\[ \left\{ [x] \in V'/V_0 \mid x \in K' \cap V' \right\}.
\]

\begin{defin} \label{def:sdr}
Конус $K$, представимый через матричный конус $\mathbb{S}_+^n$, называется \emph{полуопределённо представимым} (semi-definite representable, SDr).
\end{defin}

Все симметричные конусы являются спектраэдральными:
\begin{itemize}
\item ортант изоморфен множеству диагональных неотрицательно определённых матриц;
\item конус Лоренца $L^n$ изоморфен множеству неотрицательно определённых матриц вида
\[ \begin{pmatrix} x_n + x_1 & x_2 & \dots & x_{n-1} \\ x_2 & x_n-x_1 & 0 & 0 \\ \vdots & 0 & \ddots & 0 \\ x_{n-1} & 0 & 0 & x_n-x_1 \end{pmatrix};
\]
\item эрмитовый матричный конус, состоящий из матриц вида $S+iA$, где $S$ --- симметричная, а $A$ --- косо-симметричная вещественная матрица, изоморфен множеству неотрицательно определённых матриц вида $\begin{pmatrix} S & A \\ -A & S \end{pmatrix}.$
\end{itemize}

Справедлив следующий, более общий результат \cite{Chua03}.

\begin{teo}
Любой однородный конус является спектраэдральным.
\end{teo}

Если в исходной конической задаче конус является полуопределённо представимым, причём размерность представляющего матричного конуса не слишком большая, то задача сводится к эффективно решаемой полуопределённой программе. Естественно возникает вопрос: каков наименьший размер матричного конуса, или в более общей постановке, конуса из некоторого семейства, который представляет данный конус? Для семейства $\mathbb R_+^n$ ответ на этот вопрос даёт теорема Яннакакиса (см. п.~\ref{ch1_yannakakis}).

\subsection{Программы над несимметричными конусами}

Как говорилось выше, сложность конической программы зависит от свойств конуса $K$, над которым программа определена. Симметрические конусы имеют алгебраическое описание границы и поэтому в принципе не в состоянии представить трансцендентные ограничения. Однако, помимо симметричных конусов, есть целый ряд алгоритмически доступных конусов, в основном маломерных, которые позволяют решать и разные классы выпуклых задач с трансцендентными данными. В качестве примера мы в этом разделе рассмотрим, как представить общую задачу геометрического программирования в виде конической программы.

Геометрическая программа имеет целью минимизацию функции $f_0(x)$ на положительном ортанте $\mathbb R_{++}^n$ при ограничениях $f_i(x) \leq 1$ типа неравенства и ограничениях $h_j(x) = 1$ типа равенства. Здесь $f_i$, $i = 1,\dots,m$~--- \emph{позиномы}, т.е. функции вида $\sum_{k=1}^{k_i} c_k\prod_{l=1}^n x_l^{\alpha_{kl}}$, где $c_k > 0$, а $h_j$, $j = 1,\dots,m'$~--- мономы вида $c\prod_{l=1}^n x_l^{\alpha_l}$ с $c > 0$. Задача формализуется в виде
\begin{gather*} 
\min_{x \in \mathbb R_{++}^n} \, \sum_{k=1}^{k_0} c_{0k}\prod_{l=1}^n x_l^{\alpha_{0kl}}\ :\\
\sum_{k=1}^{k_i} c_{ik}\prod_{l=1}^n x_l^{\alpha_{ikl}} \leq 1,\ i = 1, \, \dots, \, m;\quad c_j\prod_{l=1}^n x_l^{\alpha_{jl}} = 1,\ j = 1, \, \dots, \, m'.
\end{gather*}

Введём новые переменные $y = \log x$ и параметры $b = \log c$, тогда получим эквивалентную форму
\begin{gather*} 
\min_{y \in \mathbb R^n} \, \sum_{k=1}^{k_0} \exp\left(b_{0k} + \sum_{l=1}^n\alpha_{0kl}y_l\right)\ :\\  
\sum_{k=1}^{k_i} \exp\left(b_{ik} + \sum_{l=1}^n\alpha_{ikl}y_l\right) \leq 1,\ i = 1, \, \dots, \, m;\\ b_j + \sum_{l=1}^n\alpha_{jl}y_l = 0,\ j = 1, \, \dots, \, m'.
\end{gather*}
Определим 3-мерный выпуклый конус
\[ K_{\exp} = \{ (x,y,0) \,|\, x \leq 0,\ y \geq 0 \} \cup \left\{ (x,y,z) \,|\, z > 0,\ y \geq ze^{x/z} \right\},
\]
т.н. \emph{экспоненциальный конус}. Конус $K_{\exp}$ определён как замыкание конуса над надграфиком экспоненциальной функции. Он является самодвойственным в более широком смысле, что он линейно изоморфен двойственному конусу $K_{\exp}^*$, но не является однородным.

Вводом дополнительных переменных задачу можно переписать в виде
\begin{gather*} 
\min_{y \in \mathbb R^n} \,\sum_{k=1}^{k_0} d_{0k}\ :\quad \left(b_{ik} + \sum_{l=1}^n\alpha_{ikl}y_l,d_{ik},1\right) \in K_{\exp},\ k = 1,\dots,k_i,\ i = 0,\, \dots, \, m;\\ \sum_{k=1}^{k_i} d_{ik} = 1,\ i = 1,\dots,m;
\quad
b_j + \sum_{l=1}^n\alpha_{jl}y_l = 0,\ j = 1, \, \dots, \, m'.
\end{gather*}
Теперь функция цены линейна, и помимо линейных ограничений типа равенства мы имеем только конические ограничения над конусом $K_{\exp}$. Таким образом, геометрическая программа записана в виде конической программы над произведением конечного числа копий этого 3-мерного конуса.

Отметим, что конический солвер для решения геометрических программ, основанный на этом представлении, предлагается в пакете MOSEK.

Конус $K_{\exp}$ аппроксимируется с точностью $\epsilon$ числом $O\left(\log\epsilon^{-1}\right)$ конично-квадратичных ограничений \cite{Nemir_conrut}. Поэтому геометрическую программу можно приблизить конично-квадратичной.

\subsection{Теорема Яннакакиса и её обобщения} \label{ch1_yannakakis}

Ответ на вопрос, представим ли данный конус $K \subset V_n$ через ортант $V_m = \mathbb R_+^m$ для некоторого $m$, имеет очевидный ответ. Такое представление существует тогда и только тогда, когда $K$ полиэдральный. Однако установить, какова наименьшая размерность $m$, которая для этого потребуется, оказывается нетривиальной задачей.

Этот вопрос был изучен Михалисом Яннакакисом в работе \cite{Yannakakis91} в связи с неоднократным выходом в свет <<доказательств>>\   соотношения $P = NP$ посредством (ошибочного) представления экспоненциально сложных политопов через полиномиально сложные.

Пусть $K \subset V_n$ --- регулярный полиэдральный конус, а $K^* \subset V_n^*$ --- двойственный к нему конус. Пусть $K$ имеет $r$ экстремальных лучей, а $K^*$ --- $s$ экстремальных лучей. Построим матрицы $A,B$ размера $n \times r$ и $n \times s$ соответственно, такие что столбцы $A$ и $B$ генерируют все экстремальные лучи конусов $K,K^*$ соответственно. Матрица $M = B^TA$ размера $s \times r$ называется \emph{матрицой невязок} (slack matrix). Очевидно, элементы этой матрицы неотрицательны. Заметим, что матрица невязок не зависит от выбора системы координат в пространстве $V_n$, а изменение нумерации экстремальных лучей конусов $K,K^*$ приводит к перестановке её столбцов и строк.

\begin{defin}
\emph{Неотрицательным рангом} матрицы $M$ размера $s \times r$ с неотрицательными элементами называется минимальное число $m$ такое, что существуют неотрицательные матрицы $F,G$ размера $m \times r$ и $m \times s$, соответственно, такие что $M = G^TF$.
\end{defin}

Вычисление неотрицательного ранга произвольной неотрицательной матрицы является $NP$-сложной задачей. Понятно, что неотрицательный ранг матрицы невязок зависит только от самого конуса $K$. Мы имеем следующую связь с представлениями конуса $K$ через ортант.

\begin{lemma} \label{lem:fakt1}
Пусть $M = G^TF$ --- неотрицательная факторизация матрицы невязок $M = B^TA$ конуса $K$, в которой матрицы $F,G$ имеют размер $m \times r$ и $m \times s$ соответственно. Тогда $K$ имеет следующее представление через ортант $\mathbb R_+^m${\rm :}
\[ K = \left\{ \left(BB^T\right)^{-1}BG^Ty \mid y \in \mathbb R_+^m \cap Im\,F \right\}.
\]
\end{lemma}

\begin{proof}
Заметим сначала, что $B$ имеет полный ранг и $r \geq n$, поскольку $K$ и вместе с ним и $K^*$ --- регулярный конус. Поэтому $BB^T$ обратимо.

Пусть $a_j$ --- $j$-й столбец матрицы $A$, $y_j$ --- $j$-й столбец матрицы $F$, а $e_j$ --- $j$-й базисный вектор в $\mathbb R^r$. Заметим, что $y_j \in \mathbb R_+^m \cap Im\,F$, поскольку $F$ неотрицательна. Тогда имеем
\begin{align*} 
\left(BB^T\right)^{-1}BG^Ty_j &= \left(BB^T\right)^{-1}BG^TFe_j = \left(BB^T\right)^{-1}BMe_j= \\ &= \left(BB^T\right)^{-1}BB^TAe_j = Ae_j = a_j.
\end{align*}
Из этого следует, что все генераторы экстремальных лучей $K$, а значит, и любые элементы $K$ имеют требуемое представление.

Пусть теперь $y$ --- произвольный вектор из множества $\mathbb R_+^m \cap Im\,F$. Покажем, что $x = \left(BB^T\right)^{-1}BG^Ty$ является элементом $K$. Для этого достаточно показать, что $B^Tx \geq 0$, т.е. $x$ неотрицателен на всех генерирующих элементах двойственного конуса $K^*$. Так как $y \in Im\,F$, то существует вектор $z \in \mathbb R^r$ такой, что $y = Fz$. Из этого следует, что
\begin{align*} 
B^Tx &= B^T\left(BB^T\right)^{-1}BG^Ty = B^T\left(BB^T\right)^{-1}BG^TFz= \\ &= B^T\left(BB^T\right)^{-1}BB^TAz = B^TAz = G^TFz = G^Ty \geq 0
\end{align*}
в силу неотрицательности матрицы $G$. Этого и требовалось доказать. 
\end{proof}

\medskip

С другой стороны, любое представление конуса $K$ через ортант определяет неотрицательную факторизацию матрицы невязок. Для доказательства этого факта нам понадобится следующая вспомогательная лемма.

\begin{lemma} \label{lp_auxiliary}
Пусть $K \subset V_m$ --- выпуклый полиэдральный конус, а $L: V \to V_m$ --- линейная инъекция. Пусть $w \in V^*$ --- линейный функционал такой, что образ $L[H]$ открытого полупространства $H = \{ x \in V \mid \langle w,x \rangle < 0 \}$ имеет пустое пересечение с $K$. Тогда существует элемент $g \in K^*$ такой, что $L^{\dagger}(g) = w$.
\end{lemma}

\begin{proof}
Сперва отметим, что образ $L[H]$ имеет пустое пересечение с $K$ тогда и только тогда, когда образ множества $H_1 = \{ x \in V \mid \langle w,x \rangle \leq -1 \}$ имеет пустое пересечение с $K$.

Так как $L$ инъективно, то сопряжённое отображение $L^{\dagger}$ будет сюръективным. Поэтому найдётся элемент $\tilde g \in V_m^*$ такой, что $L^{\dagger}(\tilde g) = w$. Из этого следует $\langle \tilde g,y \rangle = \langle w,x \rangle$ для всех $y = L(x)$, $x \in V$. Пусть далее $l_1,\, \dots, \, l_k \in V_m^*$ генерируют ортогональное подпространство $(L[V])^{\perp} \subset V_m^*$, т.е. ядро отображения $L^{\dagger}$, а $b_1, \, \dots, \, b_N \in V_m^*$ генерируют двойственный конус $K^*$, т.е. $K = \{ y \in V_m \mid \langle b_i,y \rangle \geq 0\ \forall\ i = 1, \, \dots, \, N \}$. 

Тогда система линейных равенств и неравенств
\[ \langle b_i,y \rangle \geq 0,\ \forall\ i = 1, \, \dots, \, N, \quad \langle l_j,y \rangle = 0,\ \forall\ j = 1, \, \dots, \, k,\quad -\langle \tilde g,y \rangle \geq 1
\]
на вектор $y \in V_m$ несовместима. По лемме об альтернативе существуют $\lambda_1, \, \dots, \, \lambda_k \in \mathbb R$ и $\mu_0,\mu_1, \, \dots, \, \mu_N \geq 0$ такие, что $-\mu_0\tilde g + \sum_{i=1}^N \mu_i b_i + \sum_{j=1}^k \lambda_j l_j = 0$, $\mu_0 > 0$. Положим $g = \tilde g - \mu_0^{-1}\sum_{j=1}^k \lambda_j l_j$. В силу $L^{\dagger}(l_j) = 0$ имеем $L^{\dagger}(g) = w$, а~в~силу $g = \mu_0^{-1}\sum_{i=1}^N \mu_i b_i$ имеем $g \in K^*$, что и требовалось доказать. 
\end{proof}

\medskip

\begin{lemma} \label{lem:fakt2}
Пусть существует представление конуса $K$ через ортант $\mathbb R_+^m$. Тогда неотрицательный ранг матрицы невязок не превосходит $m$.
\end{lemma}

\begin{proof}
Пусть отображения $L: V \to \mathbb R^m$, $\Pi: V \to V_n$ задают представление $K$ через $\mathbb R_+^m$, т.е. существует конус $\tilde K \subset V$ такой, что $\Pi[\tilde K] = \linebreak =K$, $L[\tilde K] = L[V] \cap \mathbb R_+^m$. 

Построим неотрицательную матрицу $F$ размера $m \times r$ следующим образом. В качестве $j$-го столбца $F$ определим вектор $y_j = L(z_j) \in L[V] \cap \mathbb R_+^m$, где $z_j$ --- любой прообраз из пересечения $\tilde K \cap \Pi^{-1}(a_j)$.

Определим векторы $w_i = \Pi^{\dagger}(b_i) \in V^*$. Имеем ${\langle w_i,z \rangle = \langle b_i,\Pi(z) \rangle \geq 0}$ для любого $z \in \Pi^{-1}[K]$, и поэтому $w_i \in (\Pi^{-1}[K])^*$. Вследствие ${K = \Pi[\tilde K]}$ имеем $\tilde K \subset \Pi^{-1}[K]$, $(\Pi^{-1}[K])^* \subset \tilde K^*$, и $w_i \in \tilde K^*$. Поэтому открытое полупространство $H = \{ z \in V \mid \langle w_i,z \rangle < 0 \}$ имеет пустое пересечение с $\tilde K$, а~его образ $L[H]$ имеет пустое пересечение с $\mathbb R_+^m$. В силу леммы \ref{lp_auxiliary} найдётся элемент $g_i \geq 0$ такой, что $L^{\dagger}(g_i) = w_i$. Построим матрицу $G$ размера $m \times s$ из столбцов $g_i$.

По построению имеем $\langle g_i,y_j \rangle = \langle w_i,z_j \rangle = \langle b_i,a_j \rangle$ для всех индексных пар $(i,j)$. Это даёт нам искомую неотрицательную факторизацию $M = G^TF$ матрицы невязок. 
\end{proof}

Комбинируя обе леммы, получаем следующий результат.

\begin{teo} \label{thm:fact_lp}
Пусть $K$ --- регулярный полиэдральный конус. Тогда наименьшее число $m$ такое, что существует представление $K$ через ортант $\mathbb R_+^m$, равно неотрицательному рангу матрицы невязок конуса~$K$. \qed
\end{teo}

\medskip

\begin{example}
Построить представление конуса над правильным шестиугольником через ортант $\mathbb R_+^5$.

{\bf Решение}. Конус $K \subset V = \mathbb R^3$ над правильным шестиугольником имеет шесть экстремальных лучей, которые генерируются векторами $a_k = \left(1,\cos\frac{2\pi k}{6},\sin\frac{2\pi k}{6}\right)^T$, $k = 0, \, \dots, \, 5$. Двойственный к нему конус $K^* \subset V^*$ также имеет шесть экстремальных лучей. Нетрудно увидеть, что эти лучи генерируются векторами $b_0 = a_5 \times a_0$, $b_k = a_{k-1} \times a_k$, $k = 1, \, \dots, \, 5$. Составим матрицы $A,B$ размера $3 \times 6$ из векторов $a_k,b_k$ соответственно. Тогда матрица невязок примет вид
\[ M = B^TA = \frac{\sqrt{3}}{2}\begin{pmatrix} 0 & 1 & 2 & 2 & 1 & 0 \\ 0 & 0 & 1 & 2 & 2 & 1 \\ 1 & 0 & 0 & 1 & 2 & 2 \\ 2 & 1 & 0 & 0 & 1 & 2 \\ 2 & 2 & 1 & 0 & 0 & 1 \\ 1 & 2 & 2 & 1 & 0 & 0 \end{pmatrix}.
\]
Эту матрицу можно разложить в произведение
\[ M = G^TF = \frac{\sqrt{3}}{2} \begin{pmatrix} 0 & 0 & 0 & 1 & 1 & 0 \\ 1 & 1 & 0 & 0 & 0 & 0 \\ 0 & 1 & 2 & 1 & 0 & 0 \\ 1 & 0 & 0 & 0 & 1 & 2 \\ 0 & 0 & 1 & 0 & 0 & 1 \end{pmatrix}^T \cdot \begin{pmatrix} 2 & 1 & 0 & 0 & 0 & 1 \\ 0 & 0 & 1 & 2 & 1 & 0 \\ 0 & 0 & 0 & 0 & 1 & 1 \\ 0 & 1 & 1 & 0 & 0 & 0 \\ 1 & 0 & 0 & 1 & 0 & 0 \end{pmatrix}
\]
неотрицательных матриц. Нетрудно проверить, что это минимальное разложение, и неотрицательный ранг матрицы невязок равен 5.

Отсюда получаем
\[ \Pi = \left(BB^T\right)^{-1}BG^T = \begin{pmatrix} \frac{2\sqrt{3}}{9} & \frac{2\sqrt{3}}{9} & \frac{4\sqrt{3}}{9} & \frac{4\sqrt{3}}{9} & \frac{2\sqrt{3}}{9} \\ \frac{\sqrt{3}}{3} & -\frac{\sqrt{3}}{3} & 0 & 0 & 0 \\ 0 & 0 & -1 & 1 & 0 \end{pmatrix},
\]
\[ Im\,F = \left\{x \in \mathbb R^5 \mid x_1+x_2-x_3-x_4-2x_5 = 0 \right\}.
\]
В итоге конус $K$ представлен в виде
\[ K = \left\{ \Pi x \mid x \in \mathbb R_+^5 \cap Im\,F \right\}.
\]
\end{example}

\medskip

Общий случай был исследован в работе \cite{GouveiaParriloThomas13}. Вместо матрицы невязок регулярному конусу $K$ общего вида соответствует \emph{оператор невязок} (slack operator). Пусть $\ext(K),\ext(K^*)$ --- множества генераторов экстремальных лучей конусов $K,K^*$ соответственно. Тогда оператор невязок конуса $K$ задаётся отображением $M_K: \ext(K^*) \times \ext(K) \to \mathbb R_+$, $M_K: (y,x) \mapsto \langle y,x \rangle$.

\begin{defin}
Пусть $X,Y$ --- произвольные множества, а $M: Y \times\linebreak \times\ X \to \mathbb R_+$ --- неотрицательная функция на их произведении. Пусть $K$ --- произвольный регулярный конус. Функция $M$ \emph{факторизуется через $K$}, если существуют функции $F: X \to K$, $G: Y \to K^*$ такие, что для всех $x \in X$, $y \in Y$ имеет место соотношение $M(y,x) = \langle G(y),F(x) \rangle$.
\end{defin}

Заметим, что разложение неотрицательной матрицы $M = G^TF$ на неотрицательные факторы $F,G$ размера $m \times r$, $m \times s$ является факторизацией $M$ через ортант $\mathbb R_+^m$.
Имеем следующее обобщение леммы \ref{lem:fakt1}.

\begin{lemma} \label{lem:fakt1gen}
Пусть оператор невязок $M_K$ регулярного конуса $K \subset V_n$ факторизуется через регулярный конус $K' \subset V_m$. Тогда $K$ имеет представление через конус $K'$.
\end{lemma}

\begin{proof}
Рассмотрим линейное пространство
\[ V = \left\{ (x,y) \in V_n \times \spa{Im\,F} \mid \langle w,x \rangle = \langle G(w),y \rangle\ \forall\ w \in \ext(K^*) \right\},
\]
его проекцию $\Pi: (x,y) \mapsto x$ на $V_n$, а также отображение $L: (x,y) \mapsto y$. Покажем, что этим задаётся представление $K$ через $K'$. Обозначим $\tilde K = \{ (x,y) \in V \mid y \in K' \}$.

Для любого $x \in \ext(K)$ имеем $(x,F(x)) \in V$ по определению $M_K$. Так как $K$ --- регулярный конус, то $\Pi$ сюръективно. Более того, вследствие $F(x) \in K'$ имеем $(x,F(x)) \in \tilde K$, $x \in \Pi[\tilde K]$, и тем самым $K \subset \Pi[\tilde K]$.

Далее, из $(x,0) \in V$ следует, что $\langle w,x \rangle = 0$ для всех $w \in \ext(K^*)$, что вследствие регулярности $K$ может иметь место только тогда, когда $x = 0$. Поэтому $L$ инъективно.

Осталось показать, что $\Pi[\tilde K] \subset K$. Пусть $y \in K'$, $(x,y) \in \tilde K$, $w \in \ext(K^*)$. Тогда $\langle w,x \rangle = \langle G(w),y \rangle \geq 0$, и $x$ неотрицателен на всех экстремальных лучах двойственного конуса $K^*$. Но тогда $x \in K$, что и требовалось доказать. 
\end{proof}

\medskip

Лемму \ref{lem:fakt2} также можно обобщить, но вследствие отсутствия аналога леммы \ref{lp_auxiliary} необходимо потребовать дополнительное условие.

\begin{lemma} \label{lem:fakt2gen}
Пусть существует представление регулярного конуса $K \subset V_n$ через регулярный конус $K' \subset V_m$ такое, что подпространство $L[V] \subset V_m$ имеет непустое пересечение с внутренностью конуса $K'$. Тогда оператор невязок $M_K$ конуса $K$ факторизуется через конус $K'$.
\end{lemma}

\begin{proof}
Отображение $F$ строится, как в доказательстве леммы~\ref{lem:fakt2}, т.е. $a \in \ext(K)$ отображается в $L(z)$, где $z \in \tilde K \cap \Pi^{-1}(a)$.

Построим отображение $G$. Пусть $b \in \ext(K^*)$, $w = \Pi^{\dagger}(b)$. Пусть функционал $g \in V_m^*$ отделяет конус $K'$ от множества $L[H]$, где ${H = \{ z \in V \mid \langle w,z \rangle < 0 \}}$. По построению образ $L^{\dagger}(g)$ будет пропорционален функционалу $w$, однако надо исключить случай $L^{\dagger}(g) = 0$. Это гарантируется дополнительным условием. Если $L^{\dagger}(g) = 0$, то функционал $g$ равен нулю на всех точках подпространства $L[V]$ \nk{и, таким образом,} на внутренних точках $K'$, что невозможно вследствие регулярности $K'$. Итак, $g$ можно нормализовать так, что $L^{\dagger}(g) = w$. Положим $G(b) = g$.

Тогда имеем $\langle G(b),F(a) \rangle = \langle g,L(z) \rangle = \langle w,z \rangle = \left\langle \Pi^{\dagger}(b),z \right\rangle = \langle b,a \rangle$. Таким образом, $F,G$ задают искомую факторизацию.
\end{proof}

\medskip

В случае полуопределённых представлений особенности матричных конусов позволяют получить полный аналог теоремы \ref{thm:fact_lp}. А именно, если конус $K$ представим через матричный конус $\mathbb{S}_+^m$, но подпространство $L[V]$ не содержит положительно определённых матриц, то $K' = \mathbb{S}_+^m$ можно заменить на минимальный фасад, содержащий пересечение $L[V] \cap \mathbb{S}_+^m$, а объемлющее пространство $V_m$ --- на линейную оболочку этого фасада. Но этот фасад изоморфен матричному конусу $\mathbb{S}_+^{m'}$ для некоторого $m' < m$, а его линейная оболочка --- пространству симметричных матриц $\mathbb{S}^{m'}$ размера $m' \times m'$. Поэтому $K$ представим через $\mathbb{S}_+^{m'}$, но теперь дополнительное условие в формулировке леммы \ref{lem:fakt2gen} уже выполняется. Для формулировки аналога теоремы \ref{thm:fact_lp} нам понадобится следующее обобщение неотрицательного ранга.

\begin{defin}
Пусть $X,Y$ --- произвольные множества, а $M\!\!: Y \times X \to \mathbb R_+$ --- неотрицательная функция на их произведении. \emph{Полуопределённым рангом} функции $M$ называется минимальное число $m$ такое, что существуют отображения $F: X \to \mathbb{S}_+^m$, $G: Y \to \mathbb{S}_+^m$ такие, что $M(y,x) = \langle G(y),F(x) \rangle$ для всех $x \in X$, $y \in Y$.
\end{defin}

В случае, когда $M$ --- матрица, т.е. множества $X,Y$ конечны, свойства полуопределённого ранга были изучены в \cite{FawziEtAl15}.

Мы получаем следующий результат.

\begin{teo} \label{thm:fact_sdp}
Пусть $K$ --- регулярный конус. Тогда наименьшее число $m$ такое, что существует представление $K$ через матричный конус $\mathbb{S}_+^m$, равно полуопределённому рангу оператора невязок конуса $K$. \qed
\end{teo}

В частности, конус $K$ не является полуопределённо представимым, если его оператор невязок не допускает факторизации через матричный конус какой бы то ни было размерности.

Так как нетривиальные фасады конуса Лоренца изоморфны лучу~${\mathbb R_+ \simeq L^1}$, для представления конуса $K$ через \nk{конусы} Лоренца верен аналогичный результат.

\subsection[Аппроксимация полуопределённо представимыми\\ конусами]{Аппроксимация полуопределённо представимыми\\ конусами} \label{ch1_relaxations}

Представлением одного конуса через другой можно порою сильно понизить сложность задачи. Если представление неизвестно или не существует, то можно попытаться сделать аппроксимацию.

Например, конус над правильным многогранником с $N$ гранями можно представить через ортант размерности $O(\log N)$. Конусом над таким многогранником можно аппроксимировать конус Лоренца $L^3$ с точностью $\epsilon = O(N^{-1})$, а конус Лоренца $L^n$ можно представить через прямое произведение $(L^3)^{n-2}$. В итоге конус Лоренца $L^n$ можно с точностью $\epsilon$ аппроксимировать полиэдральным конусом, обладающим представлением через ортант $\mathbb R^m$ с $m \sim O\left(n\log\epsilon^{-1}\right)$. Этот удивительный факт был доказан в 1998~г. А.\,С.~Немировским и А.~Бен-Талем~\cite{Ben-Tal, Nemir_conrut}.

\medskip

Иногда можно гарантировать некоторую точность определённой релаксации, но нет семейства релаксаций со всё возрастающей точностью. Рассмотрим два примера.

\medskip

Рассмотрим задачу максимизации однородной квадратичной формы $A$ на вершинах гиперкуба $[-1,1]^n$. Эта задача записывается в виде квадратичной программы с квадратичными ограничениями
\begin{equation} \label{MaxCube}
\max_{x \in \mathbb R^n}\, x^TAx\ : \quad x_i^2 = 1,\ \forall\ i = 1,\, \dots, \, n.
\end{equation}
Стандартная полуопределённая релаксация этой задачи получается из эквивалентной формулировки
\[ \max_{X \in \mathbb{S}_+^n}\, \langle A,X \rangle\ :\quad \diag\,X = {\bf 1},\quad \rk\,X = 1
\]
опущением условия на ранг.

Матрицы $X \succeq 0$ ранга 1 с единичной диагональю являются вершинами политопа максимального разреза (MaxCut polytope), обозначаемого через ${\cal MC}$. Этот политоп имеет $2^{n-1}$ вершин. Множество допустимых матриц в полуопределённой релаксации, т.е. неотрицательно определённых матриц с единичной диагональю, содержит политоп ${\cal MC}$ и является его внешней аппроксимацией. Это множество обозначается через ${\cal SR}$. Кроме того, введём ещё \emph{тригонометрическую аппроксимацию}: \cite{Hirschfeld}
\[ {\cal TA} = \left\{ \left. \frac{2}{\pi}\arcsin\,X \right| X \in {\cal SR} \right\} = \left\{ \left. \frac{2}{\pi}\arcsin\,X \right| X \succeq 0,\ \diag\,X = {\bf 1} \right\},
\]
где функция $t \mapsto \frac{2}{\pi}\arcsin\,t$ применяется к $X$ поэлементно. Следующий результат показывает, что невыпуклое множество ${\cal TA}$ является внутренней аппроксимацией политопа ${\cal MC}$ \cite{Nesterov98}.

\begin{lemma} \label{lemTA}
Пусть $Z \in {\cal TA}$, $X = \sin\left(\frac{\pi}{2}Z\right) \in {\cal SR}$. Пусть $\xi \sim {\cal N}({\bf 0},X)$ --- нормально распределённый случайный вектор, и пусть $x = \sign\,\xi$. Тогда $Z = \mathbb E_{\xi} xx^T$.
\end{lemma}

\begin{proof}
Для диагональных элементов имеем $Z_{ii} = 1 = \mathbb E_{\xi} x_i^2$.

Пусть теперь $i \not= j$. Вычислим вероятность того, что $x_ix_j = -1$, т.е. что элементы $\xi_i,\xi_j$ имеют разный знак. 

Пусть $\rk\,X = k$ и $X = FF^T$ --- факторизация неотрицательно определённой матрицы $X$, с фактором $F$ размера $n \times k$. Тогда случайный вектор $\xi$ можно представить в виде $\xi = F\psi$, где $\psi \sim {\cal N}({\bf 0},I)$ --- стандартный нормально распределённый вектор в $\mathbb R^k$. Обозначим строки фактора $F$ через $f_1,\dots,f_n \in \mathbb R^k$. Тогда $\xi_i = \langle f_i,\psi \rangle$, $\xi_j = \langle f_j,\psi \rangle$, и элементы $\xi_i,\xi_j$ имеют разный знак тогда и только тогда, когда гиперплоскость, ортогональная вектору $\psi$, разделяет $f_i,f_j$. Но вероятность этого события пропорциональна углу $\varphi_{ij}$ между этими векторами, а именно, $\mathbb P(x_ix_j = -1) = \frac{\varphi_{ij}}{\pi}$. С другой стороны, векторы $f_i$ имеют единичную длину, и поэтому $\cos\varphi_{ij} = \langle f_i,f_j \rangle = X_{ij}$.

В итоге математическое ожидание произведения $x_ix_j$ вычисляется по формуле
\[ \mathbb E_{\xi}x_ix_j = \mathbb P(x_ix_j = 1) - \mathbb P(x_ix_j = -1) = 1 - 2\frac{\arccos X_{ij}}{\pi} = \frac{2}{\pi}\arcsin X_{ij} = Z_{ij}.
\]
Этим утверждение доказано и для недиагональных элементов.
\end{proof}

Так как экстремальные точки политопа ${\cal MC}$ содержатся в ${\cal TA}$, то из включения ${\cal TA} \subset {\cal MC}$ также следует равенство $\conv\,{\cal TA} = {\cal MC}$.

\begin{teo}[$\pi/2$-теорема Нестерова]
Пусть $A \succeq 0$. Тогда имеет место оценка
\[ \frac{2}{\pi}\max_{X \in {\cal SR}}\, \langle A,X \rangle \leq \max_{X \in {\cal TA}}\, \langle A,X \rangle = \max_{X \in {\cal MC}}\, \langle A,X \rangle \leq \max_{X \in {\cal SR}}\, \langle A,X \rangle.
\]
\end{teo}

Константа в этой теореме не может быть улучшена \cite{NemirovskiICM}.

\begin{proof}
Равенство имеет место, поскольку линейная форма на компактном множестве принимает те же экстремальные значения, что и на его выпуклой оболочке. Правое неравенство следует из включения ${\cal MC} \subset {\cal SR}$. Для доказательства левого неравенства рассмотрим матрицу $X \in {\cal SR}$, на которой линейная форма $A$ принимает максимум, и положим $Z = \frac{2}{\pi}\arcsin\,X \in {\cal TA}$. Тогда
\[ \langle A,Z \rangle = \frac{2}{\pi}\left\langle A,\sum_{k=0}^{\infty} \frac{(2k-1)!!}{(2k)!!}\frac{X^{2k+1}}{2k+1} \right\rangle \geq \frac{2}{\pi}\langle A,X \rangle,
\]
поскольку $\langle A,X^k \rangle \geq 0$ для всех $k \in \mathbb N$ вследствие неотрицательной определённости матриц $A$ и $X^k$. 
\end{proof}

Заметим, что лемма \ref{lemTA} даёт нам возможность генерировать суб-оптимальные решения задачи \eqref{MaxCube}. Пусть $X^*$ --- решение полуопределённой релаксации задачи. Генерируя случайные векторы $\xi \sim {\cal N}({\bf 0},X^*)$ и полагая $x = \sign\,\xi$, получаем случайные вершины гиперкуба. Если $A \succeq 0$, т.е. максимизируемая квадратичная форма выпукла, то \nk{математическое ожидание} значения целевой функции на этих вершинах будет не меньше, чем $\frac{2}{\pi}$ умножить на оптимальное значение.

\medskip

Перейдём ко второму примеру, который числится в списке Карпа NP-полных (NP-complete) проблем \cite{Karp72}. Пусть задан граф с $n$ вершинами и неотрицательными весами на рёбрах. Задача нахождения максимального разреза (MaxCut) состоит в разбиении множества вершин графа на два непересекающихся подмножества таким образом, чтобы общий вес рёбер, соединяющих эти два подмножества, был максимален. Задача сводится к программе
\[ \frac14 \max_{X \in {\cal MC}}\, \langle L,X \rangle = \frac14 \max_{X \in {\cal MC}}\, \langle W,{\bf 1} - X \rangle,
\]
где $W$ --- неотрицательная поэлементно матрица весов, а $L = \diag(W{\bf 1}) - W$ --- её лапласиан. Максимальный разрез следующим образом определяется из решения $X^*$ задачи. Так как функция цены линейна, то решение является вершиной политопа ${\cal MC}$, т.е. одноранговой матрицей вида $xx^T$, где вектор $x$ является вершиной гиперкуба в $\mathbb R^n$. Оптимальное разбиение тогда определяется знаками элементов вектора $x$.

Теорема Гёманса--Виллиамсона даёт следующую (см. теорему~\ref{GW}) оценку ошибки 
полуопределённой релаксации, получающейся заменой сложного политопа ${\cal MC}$ на простое полуопределённо представимое множество ${\cal SR}$ \cite{GoemansWilliamson}.

\begin{teo}\label{GW}
Пусть $W$ --- неотрицательная поэлементно матрица. Тогда
\[ \gamma\max_{X \in {\cal SR}}\, \langle W,{\bf 1} - X \rangle \leq \max_{X \in {\cal TA}}\, \langle W,{\bf 1} - X \rangle = \max_{X \in {\cal MC}}\, \langle W,{\bf 1} - X \rangle \leq \max_{X \in {\cal SR}}\, \langle W,{\bf 1} - X \rangle,
\]
где $\gamma = \min_{t \in [-1, \, 1)}\frac{1-\frac{2}{\pi}\arcsin\,t}{1-t} \approx 0.878$.
\end{teo}

\begin{proof}
Правое неравенство и равенство доказываются, как и в предыдущем случае. Для доказательства левого неравенства рассмотрим матрицу $X \in {\cal SR}$, на которой выражение $\langle W,{\bf 1} - X \rangle$ принимает максимум, и положим $Z = \frac{2}{\pi}\arcsin\,X \in {\cal TA}$. Тогда
\[ \langle W,{\bf 1} - Z \rangle = \left\langle W,{\bf 1} - \frac{2}{\pi}\arcsin\,X \right\rangle \geq \gamma \langle W,{\bf 1} - X \rangle,
\]
поскольку $W \geq 0$. 
\end{proof}

Лемма \ref{lemTA} позволяет строить суб-оптимальные случайные разбиения, математическое ожидание веса которых не меньше, чем $\gamma$ умножить на вес максимального разреза. Если \emph{Гипотеза Уникальных Игр} верна, то константу $\gamma$ невозможно улучшить алгоритмом полиномиальной сложности \cite{KhotKindler07}.

\subsection{Двойственные конические программы} \label{ch1_conic_duality}

Пусть дана коническая программа
\[ \min_{x \in K}\, \langle c,x \rangle\ :\qquad Ax = b
\]
над некоторым регулярным выпуклым конусом $K$. Рассмотрим произвольный элемент $s \in K^*$ вида $s = c - A^Ty$, где $y$ -- некоторый вектор соответствующей размерности. Тогда для любого $x$ из допустимого множества исходной конической программы имеем
\[ \langle s,x \rangle = \langle c,x \rangle - \left\langle A^Ty,x \right\rangle = \langle c,x \rangle - \langle y,Ax \rangle = \langle c,x \rangle - \langle y,b \rangle \geq 0.
\]
Таким образом, мы получили \emph{нижнюю} оценку $\langle b,y \rangle$ на оптимальное значение исходной конической программы.

Двойственная коническая программа формулируется как задача \emph{максимизации} этой оценки и записывается в виде
\[ \max_y\ \langle b,y \rangle\ :\qquad c - A^Ty \in K^*.
\]
Она определена над двойственным конусом $K^*$. Очевидно, для каждой допустимой точки $x$ исходной конической программы величина $\langle c,x \rangle$ является \emph{верхней} оценкой оптимального значения двойственной программы.

Разница между оптимальными значениями прямой и двойственной задачи называется \emph{разрывом двойственности}.

\medskip

Двойственную пару конических программ можно записать в более симметричном виде. Аффинную оболочку допустимого множества исходной программы с функцией цены $\langle c,x \rangle$ можно представить в виде суммы $b + L$, где $L \subset \mathbb R^n$ -- некоторое линейное подпространство, а вектор $b \in \mathbb R^n$ выбран так, чтобы имело место соотношение $\langle c,b \rangle = 0$. Таким образом, получаем программу
\[ \min_{x \in K}\, \langle c,x \rangle \ :\qquad x \in b + L.
\]
Тогда двойственную программу можно записать в виде
\[ \max_{s \in K^*}\, -\langle b,s \rangle\ :\qquad s \in c + L^{\perp},
\]
где $L^{\perp} \subset \mathbb R_n$ -- ортогональное к $L$ линейное подпространство. Для любой пары $(x,s)$ допустимых точек будет иметь место соотношение
\[ \langle x,s \rangle = \langle c,x \rangle + \langle b,s \rangle \geq 0.
\]
По этому зазору в ходе итераций можно оценить, насколько близко прямо-двойственная пара $(x,s)$ находится к решению задачи.

\section{Робастная оптимизация}\label{ch1_sect_rob_opt}
В практических задачах реальные
значения параметров, требующихся для решения,
как правило, известны лишь
с какой-то степенью точности, т.е. приближённо.
Иногда удаётся пренебречь небольшими
неточностями в данных и считать, что
параметры заданы точно, в надежде,
что неизвестные небольшие расхождения
между реальными и известными (к примеру,
измеренными с какой-то
погрешностью) данными не повлияют на
оценку оптимального решения.
Но чем сложнее задача по структуре, числу переменных, количеству необходимых для получения решения вычислительных операций, чем выше требования к точности решения, тем выше вероятность
того, что игнорировать неопределённость в исходных данных уже не получится.
Часто бывает и так, что совсем небольшие
 расхождения между
 реальными и известными
данными, буквально в сотые доли процента,
приводят к нарушению ограничений задачи
в сотни тысяч (!) раз (см. примеры в~\cite[Раздел A]{BenNem_Rob9}). Тем самым
найденное решение становится полностью
бесполезным.
Таким образом,
для практического решения задач
оптимизации нужна заранее учитывающая
возможные неточности в данных
методология. Эта методология должна также
позволять получать осмысленные ответы,
согласованные с исходными неточностями.

В настоящее время существует несколько
конкурирующих подходов
к описанию неопределённости в исходных данных.

Приведём следующие три модели описания
неопределённых данных:
\begin{itemize}
    \item вероятностно-статистическая (стохастическая) \cite{Vent};
    \item нечёткая \cite{Bellman};
    \item интервальная \cite{Moore, Shar_book}.
\end{itemize}

В классической
{\it вероятностно-статистической} модели
задаётся вероятностное пространство и измеряемые
исходные данные считаются случайными величинами.
 Недостатком
 модели является
 следующий факт: даже
 если реальные данные и в самом деле выбираются
по некоторому распределению,
что бывает довольно редко,
параметры этого распределения
всё равно остаются неизвестными.

В {\it <<нечёткой>> модели} используется
понятие нечёткого множества вида
$\{ (x, \, \mu(x)) \, | \, x \, \in \, X, \,
0 \le \mu(x) \le 1 \}$,
где $\mu(x)$~--- функция принадлежности
конкретного значения неточно заданного параметра~$x$ данному нечёткому
множеству; $X$~--- заданная область возможных
значений~$x$. Функция принадлежности
обычно задаётся экспертным путём на основании
данных об источниках неточности значений
параметра~$x$.

В {\it интервальной} модели описания данных
неопределённость параметров описывается
только границами их возможных значений,
т.е. для параметра~$x$ задаётся интервал неопределённости $[x_{\min}, \, x_{\max}]$.
Интервалы неопределённости
позволяют описать неоднозначные, вариабельные
и/или неточные исходные данные.
Все значения внутри интервала предполагаются равновозможными.
В ряде прикладных задач интервальная
модель оказывается наиболее предпочтительной
(см., например, \cite{Vosh_2004}).

Особый, но довольно близкий к интервальной
модели подход к решению только
оптимизационных задач с неопределённостями
в исходных данных предлагает
ставшая популярной за последние
два десятилетия {\it робастная оптимизация}.
Развитие данной методологии в какой-то
степени
опиралось на задачи,
возникавшие
в уже существующей к
тому времени теории
робастного управления \cite{Polyak_rob, Jolen_pr_int, Hor_rob, Boyd_rob}.
В таких задачах управления
объектами
 параметры моделей объектов
и даже сами
цели управления
могут быть определены
неточно.

Так же, как и для интервального подхода,
в робастной оптимизации считается,
что данные меняются только в заданной области неопределённости. Решения обязательно
должны быть робастно допустимыми,
т.е. при всех наборах данных из заданной
области неопределённости
решения должны удовлетворять
ограничениям задачи.
Среди робастно допустимых решений выбирается
наилучшее по целевому функционалу.

Стохастическая задача также может быть решена с помощью такого подхода. В этом случае область неопределённости задаётся \emph{множеством доверия}, в котором случайная величина находится с вероятностью, не меньшей некоторого заданного наперёд порога, например, $95\%$ или $99\%$. Если данные распределены нормально, то область неопределённости задаётся соответствующим \emph{эллипсоидом доверия}.

Рассмотрим более подробно парадигму робастной
оптимизации (РО) на примере задачи линейного
программирования:
\begin{equation}
\label{ch1_eq_lp_ro}
\min_{A x \, \le \, b} c^T x.
\end{equation}

В парадигме РО задачу~\eqref{ch1_eq_lp_ro}
необходимо дополнить
условием принадлежности всех коэффициентов
задачи {\it множеству неопределённости}~$\mathbf{U}$:
\begin{equation}
\label{ch1_eq_rob_set}
(c, \, A, \, b) \, \in \, {\mathbf{U}}.
\end{equation}
Тут сразу же возникает вопрос о необходимости
различать типы неопределённости,
которые несут таким образом описанные множества
входных данных. И смысл поставленной задачи,
и результаты её решения сильно различаются
в зависимости от применения разных
логических кванторов, существования
или всеобщности, при определении множеств в~\eqref{ch1_eq_rob_set}. Ответом
на этот вопрос для классической
постановки задачи РО является
следующее уточнение~\eqref{ch1_eq_lp_ro} и~\eqref{ch1_eq_rob_set}:
\begin{equation}
\label{ch1_robust_versionLP}
\min_{x: A x \, \le \, b \ \forall \, (c, \, A, \, b) \, \in \, \mathbf{U}}\,\sup_{\exists (A,b): (c, \, A, \, b) \, \in \, \mathbf{U}}
c^Tx,
\end{equation}
т.е. среди всех возможных значений целевой функции
выбирается наихудшее, максимальное, а затем ищется наилучшее, минимальное значение среди всех <<худших>>.


Впервые задача в форме~\eqref{ch1_eq_lp_ro}--\eqref{ch1_eq_rob_set} рассматривалась в статье Сойстера 1973~г.~\cite{Soyster_rob},
но активные исследования в области линейной~РО
стали проводиться лишь спустя 25 лет, в конце 90-х~гг.,
в том числе и по дискретной
оптимизации~\cite{Kou_rob},
и по непрерывной выпуклой~РО~\cite{BenNem_conv_rob98, El_rob98}.
С этого времени
началось бурное развитие данного направления,
подробности со ссылками см. в обзоре~\cite{BenNem_rob_obz02}.
Методология современной
робастной оптимизации подробно
описана в книге А.С.~Немировского, А.~Бен-Тала и Л.~Эль-Гауи~\cite{BenNem_Rob9} и кратко~--- в
замечательном курсе лекций по современной выпуклой оптимизации А.С.~Немировского и А.~Бен-Тала~\cite{Ben-Tal}.

Можно ли эффективно решать задачи типа \eqref{ch1_robust_versionLP}, т.е. находить оптимальные робастные решения?
Оказывается, в некоторых случаях можно, причём не только линейные, но и более сложные робастные программы.


Если множество неопределённости $\mathbf{U}$ выпукло, то такую задачу можно преобразовать в детерминированную коническую программу над более сложным конусом положительных отображений. Эта программа называется робастной версией исходной задачи. Мы рассмотрим, как сводить робастные версии разных классов стандартных программ снова к стандартной программе. Если такое преобразование невозможно, то можно провести аппроксимацию. Для случая матричного куба мы представим результат, оценивающий ошибку, сделанную при аппроксимации.


\subsection[Робастная версия задачи конического\\ программирования]{Робастная версия задачи конического\\ программирования}

В п.~\ref{ch1_sect_cone_prog} мы представили две разных формы конической программы. В виде \eqref{ch1_conic_indirect} на переменные $x$ накладываются коническое условие и линейные условия типа равенства. В виде \eqref{ch1_eq_con_prog} переменные $x$ прямо параметризуют допустимое множество задачи. При переходе к робастной версии эти представления ведут себя по-разному.

Очевидно, в общем случае невозможно гарантировать выполнение условия типа равенства для всех наборов коэффициентов из множества неопределённости. Поэтому неопределённые коэффициенты могут присутствовать только в условиях типа неравенства, к которым можно отнести также коническое ограничение, и робастная версия имеет смысл только для конической задачи в форме \eqref{ch1_eq_con_prog}. Рассмотрим задачу конического программирования
\[ \min_x \, \langle c,x \rangle\ :\quad Ax - b \in K
\]
с неопределёнными коэффициентами $(A,b) \in \mathbf{U}$. Здесь $K \subset V$ --- регулярный выпуклый конус. Условие, что коэффициенты $c$ функции цены не являются неопределёнными, не ограничивает общности, поскольку в противном случае этого легко можно добиться вводом новой переменной $t$ и ограничением $t \geq \langle c,x \rangle$. Для любого набора коэффициентов $(A,b)$ из множества неопределённости $\mathbf{U}$ мы имеем обычную коническую задачу, но решение, оптимальное для одного набора коэффициентов, может быть недопустимым для другого набора.

Робастная версия задачи записывается в виде
\[ \min_x \, \langle c,x \rangle\ :\quad Ax - b \in K\qquad \forall\ (A,b) \in \mathbf{U},
\]
т.е. мы ищем наилучшее решение, удовлетворяющее коническому ограничению для всех наборов коэффициентов из множества неопределённости.

Мы будем предполагать, что множество $\mathbf{U}$ выпукло и параметризовано линейно некоторыми переменными $u_1, \, \dots, \, u_m$,
\[ \mathbf{U} = \left\{ \left. (A,b) = (A_0,b_0) + \sum_{k=1}^m u_k (A_k,b_k) \right| (u_1, \, \dots, \, u_m) \in U \right\},
\]
где $U$ --- выпуклое компактное множество с непустой внутренностью, характеризующее множество неопределённости коэффициентов задачи. Покажем, что тогда робастная версия задачи конического программирования сама является задачей конического программирования.

Введём однородную версию множества $U$, а именно регулярный конус
\[ K_U = \left\{ u = (u_0,u_1, \, \dots, \, u_m) \in \mathbb R^{m+1} \mid u_0 \geq 0,\ (u_1, \, \dots, \, u_m) \in u_0U \right\}.
\]
Так как коническое ограничение $Ax - b \in K$ не изменится, если умножить $(A,b)$ на положительную константу, робастная версия запишется в виде
\[ \min_x \, \langle c,x \rangle\ :\quad \sum_{k=0}^m u_k (A_kx - b_k) \in K\qquad \forall\ u \in K_U.
\]
Определим линейное отображение:
\[ L_x: \mathbb R^{m+1} \to V, \qquad L_x: u = (u_0, \, \dots, \, u_m) \mapsto \sum_{k=0}^m u_k (A_kx - b_k),
\]
коэффициенты которого аффинны по $x$. Тогда задача примет вид
\[ \min_x \, \langle c,x \rangle\ :\quad L_x[K_U] \subset K.
\]

\begin{defin}
Пусть $K \subset V$, $K' \subset V'$ --- регулярные выпуклые конуса. Линейное отображение $L: V \to V'$, переводящее конус $K$ в подмножество конуса $K'$, называется \emph{положительным} по отношению к этим конусам. Множество положительных отображений обозначается через $\Pos(K,K')$.
\end{defin}

Нетрудно показать, что множество положительных отображений само является регулярным выпуклым конусом в пространстве $\End(V,V')$ всех линейных отображений, или эндоморфизмов $L: V \to V'$.

С этим определением робастная версия принимает вид
\[ \min_x \, \langle c,x \rangle\ :\quad L_x \in \Pos(K_U,K).
\]
Таким образом, она сама является задачей конического программирования, но над конусом положительных отображений $\Pos(K_U,K)$. Здесь переменная $x$ по-прежнему параметризует аффинное подпространство, задающееся линейными ограничениями конической программы в её стандартной формулировке.

Очевидно, сложность робастной версии зависит от того, насколько алгоритмически доступен конус положительных отображений $\Pos(K_U,K)$. Это, в свою очередь, зависит и от конуса $K$, лежащего в основе исходной задачи, и от конуса $K_U$, задаваемого множеством неопределённости.

\subsection{Конусы положительных отображений}

В этом разделе мы рассмотрим конусы положительных отображений, встречающиеся в робастных версиях конических программ над симметричными конусами, их представления и аппроксимации. Сперва, однако, приведём некоторые результаты общего характера.

\begin{lemma}
Пусть $K \subset V$, $K' = K_1' \times \dots \times K_N' \subset V' = V_1' \times \dots \times V_N'$ --- регулярные выпуклые \nk{конусы}. Тогда $\Pos(K,K') = \Pos(K,K_1') \times \dots \times \Pos(K,K_N')$.
\end{lemma}

\begin{proof}
Для линейного отображения $L = (L_1, \, \dots, \, L_N): V \to V'$, $L_i: V \to V_i'$, и точки $x \in K$ имеем $L(x) \in K'$ тогда и только тогда, когда $L_i(x) \in K_i'$ для всех $i$. Поэтому $L \in \Pos(K,K')$ эквивалентно совокупности условий $L_i \in \Pos(K,K_i')$, $i = 1, \, \dots, \, N$. 
\end{proof}

\begin{lemma}
Пусть $K \subset V$, $K' \subset V'$ --- регулярные выпуклые \nk{конусы}, и пусть $\tilde V \subset V'$ --- линейное подпространство. Обозначим пересечение $\tilde V \cap K'$ через $\tilde K$. Тогда конус положительных отображений $\Pos(K,\tilde K)$ канонически изоморфен пересечению конуса $\Pos(K,K')$ с линейным подпространством $\{ L: V \to V' \mid Im\,L \subset \tilde V \}$.
\end{lemma}

\begin{proof}
Линейное отображение $\tilde L: V \to \tilde V$ является положительным тогда и только тогда, когда $L = I \circ \tilde L$ является положительным, где $I: \tilde V \to V'$ --- идентичное вложение. 
\end{proof}

Таким образом, если конус $K$ является линейным сечением конуса $K'$, то конус положительных отображений $\Pos(K_U,K)$ представим через $\Pos(K_U,K')$. Для линейных проекций подобное утверждение не имеет место.

\begin{lemma}
Пусть $K \subset V$, $K' \subset V'$ --- регулярные выпуклые \nk{конусы}, а $L: V \to V'$ --- линейное отображение. Тогда имеем $L \in \Pos(K,K')$ тогда и только тогда, когда $L^{\dagger} \in \Pos((K')^*,K^*)$.
\end{lemma}

\begin{proof}
Пусть $L \in \Pos(K,K')$, $z \in (K')^*$. Тогда для любого $x \in K$ имеем $\left\langle L^{\dagger}(z),x \right\rangle = \langle z,L(x) \rangle \geq 0$, поскольку $L(x) \in K'$. Из этого следует, что $L^{\dagger}(z) \in K^*$. Отсюда вытекает $L^{\dagger} \in \Pos((K')^*,K^*)$.

Обратная импликация доказывается аналогичным образом. 
\end{proof}

В частности, \nk{конусы} положительных отображений $\Pos(K_U,K)$ и $\Pos(K^*,K_U^*)$ изоморфны, и если один из них представим через некоторый конус $K'$, то второй также представим через $K'$.

Рассмотрим теперь более конкретные примеры. Пусть $K_U \subset \mathbb R^{m+1}$, $K \subset V$ --- регулярные выпуклые \nk{конусы}.

\medskip

\begin{lemma}
Если конус $K_U$ полиэдральный с $N$ экстремальными лучами, то конус положительных отображений $\Pos(K_U,K)$ представим через прямое произведение $K^N$.
\end{lemma}

\begin{proof}
Пусть $x_1, \, \dots, \, x_N \, \in \, \mathbb R^{m+1}$ --- генераторы экстремальных лучей конуса $K_U$. Тогда $L \in \Pos(K_U,K)$ тогда и только тогда, когда $L(x_i) \in K$ для всех $i = 1, \dots, N$. Но это условие эквивалентно условию $(L(x_1), \, \dots, \, L(x_N)) \in K^N$, задающим изоморфизм между конусом $\Pos(K_U,K)$ и некоторым линейным сечением произведения $K^N$. 
\end{proof}

\begin{corr}
Если исходная коническая программа класса LP (SOCP, SDP), а множество неопределённости $\mathbf{U}$ является политопом, то робастная версия конической программы также класса LP (SOCP, SDP).
\end{corr}

\begin{proof}
Если $U$ --- политоп, то конус $K_U$ является полиэдральным. Однако семейство конусов, над которыми задаются программы классов LP и SOCP, замкнуты относительно операции построения прямых произведений. Прямое произведение же матричных конусов изоморфно линейному сечению матричного конуса большей размерности, задающегося блочно-диагональными матрицами соответствующего типа. 
\end{proof}

Если число экстремальных точек политопа $\mathbf{U}$ небольшое, например, если $\mathbf{U}$ параметризуется $L_1$-шаром $U = \{ u \mid ||u||_1 \leq r \}$, то сложность робастной версии конической программы сравнима со сложностью самой конической программы. В этом случае экстремальные точки множества $\mathbf{U}$ называются \emph{сценариями}, а решение робастной задачи является наилучшей точкой, удовлетворяющей ограничениям при любом сценарии.

Если политоп $\mathbf{U}$ имеет большое количество экстремальных точек, например, если $\mathbf{U}$ параметризуется $L_{\infty}$-шаром $U = \{ u \mid ||u||_{\infty} \leq r \}$, то робастная версия, хотя и формально является задачей того же класса, может быть значительно сложнее, чем исходная коническая программа.

\medskip

Часто неопределённые параметры задачи нормально распределены, и соответствующее множество неопределённости эллипсоидально. В этом случае конус $K_U$ изоморфен конусу Лоренца $L_{m+1}$.

Если исходная коническая программа из класса SOCP, то её робастная версия с эллипсоидальным множеством неопределённости является задачей класса SDP. Этот факт является следствием полуопределённой представимости конуса $\Pos(L_n,L_m)$. Для того чтобы описать это представление, мы прибегнем к полуопределённому представлению конуса Лоренца, задаваемого инъекцией $f_n: \mathbb R^n \to \mathbb{S}^{n-1}$,
\begin{equation} \label{Lorentz_rep}
f_n: (x_0, \, \dots, \, x_{n-1})^T \mapsto \begin{pmatrix} x_0 + x_1 & x_2 & \dots & x_{n-1} \\ x_2 & x_0 - x_1 & 0 & 0 \\ \vdots & 0 & \ddots & 0 \\ x_{n-1} & 0 & 0 & x_0-x_1 \end{pmatrix}.
\end{equation}
Напомним, что линейное отображение $L: \mathbb R^n \to \mathbb R^m$ можно описать элементом из тензорного произведения $\mathbb R^m \otimes \mathbb R_n$, где $\mathbb R_n$ --- двойственное к $\mathbb R^n$ пространство. Элемент $z \otimes y$, где $y \in \mathbb R_n$, $z \in \mathbb R^m$, действует на вектор $x \in \mathbb R^n$ по правилу $x \mapsto \langle y,x \rangle \cdot z$ и тем самым задаёт соответствующее $z \otimes y$ линейное отображение $L$. Для представления положительного отображения $L \in \Pos(L_n,L_m)$, поэтому напрашивается использовать тензорное произведение $f_m \otimes f_n: \mathbb R^m \otimes \mathbb R_n \to \mathbb{S}^{m-1} \otimes \mathbb{S}^{n-1}$.

Произведение $\mathbb{S}^{m-1} \otimes \mathbb{S}^{n-1}$ можно рассматривать как подпространство пространства матриц $\mathbb{S}^{(m-1)(n-1)}$. При этом элементу $A \otimes B \in \mathbb{S}^{m-1} \otimes \mathbb{S}^{n-1}$, где $A \in \mathbb{S}^{m-1}$, $B \in \mathbb{S}^{n-1}$, соответствует кронекеровское произведение $A \otimes B$. (По этой причине один и тот же знак $\otimes$ обозначает и кронекеровское произведение матриц, и тензорное произведение векторных пространств.) Более конкретно, матричное пространство $\mathbb{S}^{(m-1)(n-1)}$ представляется в виде прямой суммы $(\mathbb{S}^{m-1} \otimes \mathbb{S}^{n-1}) \oplus (\mathbb{A}^{m-1} \otimes \mathbb{A}^{n-1})$, где $\mathbb{A}^n$ --- пространство косо-симметричных матриц размера $n \times n$. Справедлив следующий результат \cite{Hildebrand11e}.

\begin{lemma} \label{lem:L_n_L_m}
Пусть $L: \mathbb R^n \to \mathbb R^m$ --- линейное отображение. Отображение $L$ положительное, $L \in \Pos(L_n,L_m)$, тогда и только тогда, когда
\[ \exists\ A \in \mathbb{A}^{m-1} \otimes \mathbb{A}^{n-1}\ : \quad (f_m \otimes f_n)(L) + A \succeq 0.
\]
Здесь $L$ рассматривается как элемент пространства $\mathbb R^m \otimes \mathbb R_n$, а $f_n$ на $\mathbb R_n$ задаётся той же формулой \eqref{Lorentz_rep}.
\end{lemma}

\medskip

\begin{exercise}
В явном виде построить полуопределённое представление конуса $\Pos(L_3,L_3)$ через матричный конус $\mathbb{S}_+^4$.
\end{exercise}


Перейдём теперь к робастным версиям полуопределённых программ с эллипсоидальным множеством неопределённости. В общем случае проверить включение $L \in \Pos(L_n,\mathbb{S}_+^m)$ является сложной задачей \cite{NesterovMatEllipsoid}. Однако можно использовать представление \eqref{Lorentz_rep} для построения полуопределённо представимой аппроксимации конуса $\Pos(L_n,\mathbb{S}_+^m)$. Линейное отображение $L: \mathbb R^n \to \mathbb{S}^m$ можно представить в виде элемента тензорного произведения $\mathbb{S}^m \otimes \mathbb R_n$. А именно, элементу $S \otimes y \in \mathbb{S}^m \otimes \mathbb R_n$ сопоставим отображение $L$, действующее на элементах $x \in \mathbb R^n$ по правилу $L: x \mapsto \langle y,x \rangle \cdot S$. Справедлив следующий результат \cite{Hildebrand07b}.

\begin{lemma}
Пусть $L: \mathbb R^n \to \mathbb{S}^m$ --- линейное отображение, представленное в виде элемента пространства $\mathbb{S}^m \otimes \mathbb R_n$. Если выполнено условие
\[ \exists\ A \in \mathbb{A}^m \otimes \mathbb{A}^{n-1}\ : \quad (Id \otimes f_n)(L) + A \succeq 0,
\]
то $L$ является положительным, $L \in \Pos(L_n,\mathbb{S}_+^m)$. Здесь $f_n$ на $\mathbb R_n$ задаётся той же формулой \eqref{Lorentz_rep}.
\end{lemma}

Данная полуопределённая аппроксимация конуса $\Pos(L_n,\mathbb{S}_+^m)$ является внутренней. Применяя её, мы сужаем множество допустимых точек, и решение аппроксимации более консервативно, чем робастная версия задачи.

Однако если полуопределённый конус $K$ в исходной конической программе является прямым произведением матричных конусов $\mathbb{S}_+^{m_i}$ с $m_i \leq 3$, то релаксация точна. Этот результат вытекает из следующей леммы \cite{Hildebrand07b}.

\begin{lemma}
Пусть $L: \mathbb R^n \to \mathbb{S}^3$ --- линейное отображение, представленное в виде элемента пространства $\mathbb{S}^3 \otimes \mathbb R_n$. Включение $L \in \Pos(L_n,\mathbb{S}_+^3)$ имеет место тогда и только тогда, когда
выполнено условие
\[ \exists\ A \in \mathbb{A}^3 \otimes \mathbb{A}^{n-1}\ : \quad (Id \otimes f_n)(L) + A \succeq 0.
\]
\end{lemma}

\begin{exercise}
Доказать аналог этой леммы для случая $m = 2$, т.е. построить полуопределённое представление конуса положительных отображений $\Pos(L_n,\mathbb{S}_+^2)$. (Подсказка: использовать лемму \ref{lem:L_n_L_m} и тот факт, что \nk{конусы} $L_3$ и $\mathbb{S}_+^2$ изоморфны.)
\end{exercise}


Второй частный случай, в котором релаксация точна, это когда эллипсоид, задающий множество неопределённости $\mathbf{U}$, двумерный. В этом случае конус $K_U$ изоморфен $L_3$. Справедлив следующий результат.


\begin{lemma}
Пусть $L: \mathbb R^3 \to \mathbb{S}^m$ --- линейное отображение, заданное формулой $L: (x_0,x_1,x_2) \mapsto x_0S_0 + x_1S_1 + x_2S_2$. Тогда $L \in \Pos(L_3,\mathbb{S}_+^m)$ тогда и только тогда, когда выполнено условие
\[ \exists\ A \in \mathbb{A}^m\ : \quad \begin{pmatrix} S_0+S_1 & S_2+A \\ S_2-A & S_0-S_1 \end{pmatrix} \succeq 0.
\]
\end{lemma}

\subsection{Теорема о матричном кубе}

В этом разделе мы представим результат А.~Бен-Таля и А.С.~Немировского о матричном кубе \cite{BenTalNemirovski02}. Этот результат состоит в оценке на точность полуопределённой релаксации конуса положительных отображений $\Pos(K_U,\mathbb{S}_+^n)$, где $K_U \subset \mathbb R^{m+1}$ построен над единичным кубом $U = [-1, \, 1]^m$. Этот конус возникает в робастных версиях полуопределённых программ, в которых множество неопределённости $\mathbf{U}$ аффинно изоморфно $L_{\infty}$-шару.

Линейное отображение $L: \mathbb R^{m+1} \to \mathbb{S}^n$ задаётся матрицами $B_0, \, \dots, \, B_m \in \mathbb{S}^n$ по формуле $x = (x_0, \, \dots, \, x_m)^T \mapsto \sum_{l=0}^m x_lB_l$. Ясно, что $L \in \Pos(K_U,\mathbb{S}_+^n)$ тогда и только тогда, когда $B_0 + \sum_{l=1}^m \epsilon_lB_l \succeq 0$ для всех комбинаций $\epsilon_l \in \{-1, \, +1\}$. Для того чтобы проверить это включение, надо проверить экспоненциальное число --- $2^m$ линейных матричных неравенств. Однако, это условие можно усилить с понижением сложности до линейной по $m$. Рассмотрим полуопределённую релаксацию
\[ {\cal SR} = \left\{ \left. L: x \mapsto \sum_{l=0}^m x_lB_l \,\right|\, \exists\ X_l \in \mathbb{S}^n:\ X_l \succeq \pm B_l,\ l = 1, \, \dots, \, m; B_0 \succeq \sum_{l=1}^m X_l \right\}.
\]
Очевидно, ${\cal SR} \subset \Pos(K_U,\mathbb{S}_+^n)$, и ${\cal SR}$ является внутренней аппроксимацией. Проверка условия $L \in {\cal SR}$ эквивалентна проверке на совместимость системы из $2m+1$ линейных матричных неравенств.

\begin{exercise}
Доказать, что эта релаксация точна для размера матриц $n = 1$.
\end{exercise}


Для того чтобы описать оценку на ошибку релаксации, нам понадобится следующая величина. Для $k \in \mathbb N_+$, определим
\begin{align*}
\eta(k) &= \min_{V \in \mathbb{S}^k} \mathbb E_{\xi}|\xi^TV\xi|\ :\quad ||V||_1 = 1, \\ &= \min_{\lambda \in \mathbb R^k} \mathbb E_{\kappa}\left|\sum_{j=1}^k \lambda_j \kappa_j\right| \ :\quad ||\lambda||_1 = 1,
\end{align*}
где $\xi \in \mathbb R^k$ --- стандартный нормально распределённый случайный вектор, а $\kappa_1, \, \dots, \, \kappa_k$ --- независимые неотрицательные случайные скаляры, распределённые с плотностью $\mu(t) = (2\pi te^t)^{-1/2}$, т.е. по закону $\chi^2(1)$. Первые значения $\eta(k)$ заданы в следующей таблице.

\medskip

\begin{tabular}{c|c|c}
$k$ & $\eta(k)$, точно & $\eta(k)$, численно \\
\hline
1 & 1 & $1.0000$ \\
2 & $2/\pi$ & $0.6366$ \\
3 & Root($t^3 + 9t^2 + 135t - 81$) & $0.5764$ \\
4 & $1/2$ & $0.5000$
\end{tabular}

\medskip

Заметим, что \nk{математическое ожидание} выпукло по $\lambda$ и инвариантно относительно перестановок элементов $\lambda$. 
Поэтому минимум по $\lambda$ на единичной сфере $1$-нормы достигается в точке вида $\lambda = (\lambda_+, \, \dots, \, \lambda_+,\lambda_-, \, \dots, \, \lambda_-)$, где $\lambda_+ > 0$, $\lambda_- < 0$ --- некоторые вещественные числа, встречающиеся $k_+$ и $k_-$ раз соответственно. Соответствующее значение вычисляется по формуле
\[ \frac{2\Gamma(\frac{k}{2}+1)}{\Gamma(\frac{k_+}{2})\Gamma(\frac{k_-}{2})} \int_0^1 |\lambda_+\tau + \lambda_-(1-\tau)| \tau^{k_+/2-1}(1-\tau)^{k_-/2-1}\,d\tau.
\]
Подставляя $k_+ = k_- = \frac{k}{2}$, $\lambda_+ = -\lambda_- = \frac{1}{k}$, получаем значение $\frac{\Gamma(\frac{k}{2}+1)}{2^{k/2}\Gamma(\frac{k}{4}+1)^2} \approx \linebreak \approx\frac{2}{\sqrt{\pi k}}$. При чётных $k$ это значение достигается на допустимом векторе $\lambda$. Кроме того, имеем нижнюю оценку $\eta(k) \geq \frac{2}{\pi\sqrt{k}}$ \cite{BenTalNemirovski02}, и $\eta(k)$ по построению убывает. Отсюда получаем асимптотику $\eta(k) \sim k^{-1/2}$.


\smallskip

Допустим теперь, что $L \not\in {\cal SR}$, где $L$ задаётся набором матриц $(B_0, \, \dots, \, B_m)$. Тогда найдётся элемент $Z = (Z_0, \, \dots, \, Z_m)$ из двойственного конуса ${\cal SR}^*$ такой, что $\langle Z,L \rangle = \sum_{j=0}^m \langle Z_j,B_j \rangle < 0$. По построению ${\cal SR}$ является линейной проекцией произведения $(\mathbb{S}_+^n)^{2m+1}$. Следовательно, двойственный конус ${\cal SR}^*$ является сечением этого произведения, а именно, \ag{он задаётся множеством}
\[ \left\{ (Z_0,Z_1^+-Z_1^-, \, \dots, \, Z_m^+-Z_m^-) \mid Z_j^+,Z_j^- \succeq 0,\ Z_j^+ + Z_j^- = Z_0\ \forall\ j = 1, \, \dots, \, m \right\}.
\]

Нетрудно доказать, что для матриц $B \in \mathbb{S}^n$, $Z \in \mathbb{S}_+^n$ оптимальное значение полуопределённой программы
\[ \min_{Z^+,Z^- \succeq 0} \langle B,Z^+-Z^- \rangle\ :\quad Z^+ + Z^- = Z
\]
равно $-||Z^{1/2}BZ^{1/2}||_1$.


{
\renewcommand{\baselinestretch}{0.95}
\selectfont

Заметим, что для матрицы $V \in \mathbb{S}^n$ ранга не больше $k$ имеет место оценка
\[ \mathbb E_{\xi} |\xi^TV\xi| \geq \eta(k) \cdot ||V||_1,
\]
где $\xi \sim {\cal N}({\bf 0},I)$ --- случайный нормально распределённый вектор.

Отсюда следует, что
\begin{align*} 0 &> \langle Z,L \rangle \geq \tr\,\left(Z_0^{1/2}B_0Z_0^{1/2}\right) - \sum_{j=1}^m \left\|Z_0^{1/2}B_jZ_0^{1/2}\right\|_1 \geq\\ &\geq \mathbb E_{\xi}\left( \xi^TZ_0^{1/2}B_0Z_0^{1/2}\xi - \sum_{j=1}^m \frac{1}{\eta(\rk\,B_j)}\left|\xi^TZ_0^{1/2}B_jZ_0^{1/2}\xi\right| \right).
\end{align*}
Таким образом, найдётся вектор $\zeta \in \mathbb R^n$ такой, что
\[ \zeta^TB_0\zeta - \sum_{j=1}^m \frac{1}{\eta(\rk\,B_j)}|\zeta^TB_j\zeta| < 0.
\]
Полагая $\epsilon_j = -\sign(\zeta^TB_j\zeta)$, получаем ${\zeta^T\left(B_0 + \sum_{j=1}^m \frac{\epsilon_j}{\eta(\rk\,B_j)}B_j\right)\zeta < 0}$. Отсюда вытекает, что отображение $\tilde L$, задающееся набором матриц $\left(B_0,\frac{1}{\eta(\rk\,B_1)}B_1, \, \dots, \, \frac{1}{\eta(\rk\,B_m)}B_m\right)$, не является элементом конуса положительных отображений $\Pos(K_U,\mathbb{S}_+^n)$.

Таким образом, в точках, которые соответствуют матрицам $B_1, \, \dots, \, B_m$ малого ранга, граница внутренней аппроксимации ${\cal SR}$ конуса $\Pos(K_U,\mathbb{S}_+^n)$ находится близко к границе этого конуса. Получаем следующую теорему.

\begin{teo}
	Пусть матричный куб $$\mathbf{U} = \left\{ \left. B_0 + \sum_{j=1}^m u_jB_j \,\right|\, u_j \in [-1, \, 1] \right\}$$ задаётся набором матриц $(B_0, \, \dots, \, B_m)$, не включенным в множес\-тво~${\cal SR}$. Тогда расширенный матричный куб $$\tilde{\mathbf{U}} = \left\{ \left. B_0 + \sum_{j=1}^m u_jB_j \,\right|\, \eta(\rk\,B_j) \cdot |u_j| \leq 1 \right\}$$ не является подмножеством матричного конуса $\mathbb{S}_+^n$.
\end{teo}


}

\chapter[Эффективность численных методов выпуклой оптимизации\dotfill]{Эффективность численных методов выпуклой оптимизации}\label{chapt_num_met}
\chaptermark{Эффективность методов выпуклой оптимизации}
\epigraph{The main lesson of the development of our field in the last few
decades is that efficient optimization methods can be developed only by intelligently employing the structure of particular instances of problems.}{Nesterov Yu. Lectures on Convex Optimization. Springer, 2018.}

В этой главе исследуются возможности численных методов решения задач выпуклой оптимизации.
Вводятся понятия {\it оракула},  {\it чёрного ящика}, {\it функциональной модели оптимизационной задачи}.

Обратим внимание на то, что в главу мы включили совсем свежие результаты об ускоренных градиентных методах. В частности, результаты о том, что ускоренные методы (градиентные и тензорные) для разнообразных задач выпуклой оптимизации можно строить с помощью общей оболочки <<Ускоренного Мета-алгоритма>> (УМ). Достаточно подробно рассматриваются вопросы о накоплении неточности в градиенте в оценках скорости сходимости ускоренных методов. Приводятся примеры, в которых возникают неточные градиенты (безградиентные методы и решения обратных задач). Разбираются методы внутренней точки (глобализирующие сходимость метода Ньютона), которые положены в основу ряда современных пакетов решения задач выпуклой оптимизации умеренных размеров.

%
%
\section{Концепция чёрного ящика}
\label{sect-black-box}
Нет никакой возможности разумно определить наиболее
эффективный метод решения конкретной задачи непрерывной оптимизации.
Действительно, мы можем просто <<подсказать>> методу
найденное другим способом решение, и он будет выдавать
его мгновенно, за константное время~$O(1)$. Но такой
метод будет совершенно бесполезен для решения других задач,
с отличающимся оптимальным решением. 
Поэтому, как правило,  численные
методы разрабатываются для решения многих однотипных задач с
близкими характеристиками (говорят: для задач {\it одного класса}). Следовательно, эффективность метода 
на каком-то классе задач можно считать важной характеристикой качества метода.

\begin{defin}
Та информация о задаче, которая заранее известна численному методу, называется {\it моделью решаемой задачи}.
\end{defin}

Как минимум, в модель включается формализованная постановка задачи, например, для
задачи с ограничениями типа равенств~--- в форме~\eqref{ch1_lagr_eqs_problem}.

Здесь и далее эффективность (оптимальность) метода на классе задач будет пониматься в смысле Бахвалова–Немировского~\cite{NemirYd79}~--- число обращений (по ходу работы метода) к {\it оракулу}.
Оракулом называется подпрограмма расчёта значений
целевой функции (в некоторых случаях~--- градиента, в 
совсем отдельных случаях~--- и производных более высоких порядков) для достижения заданной точности решения в зависимости от параметров, характеризующих класс рассматриваемых задач (от модели задачи). 
Наивысший доступный для вычисления порядок производной
в модели задачи
называется порядком оракула (т.е. если доступны
только значения целевой функции, то говорят, что
имеется {\it оракул нулевого порядка}, соответственно,
если доступен градиент (или его аналог), то имеется {\it оракул
первого порядка} и т.д.).

Итак, аналитическая сложность конкретной оптимизационной
задачи~$P$ для какого-то метода
характеризуется тем, сколько запросов
к оракулу необходимо сделать
для нахождения приближённого решения~$P$ с заранее заданной точностью~$\varepsilon > 0$. 
При этом для обоснования {\it верхних оценок}
сложности необходимо предложить алгоритм
решения задачи~$P$, в то время, как
{\it нижние оценки}, зависящие
от точности~$\varepsilon$, могут быть получены теоретически.

Такой (оракульный) взгляд
на сложность задач выпуклой оптимизации оказался очень удобным и популярным~\cite{Nest, Bubeck15}. Связано это
с тем, что, с одной стороны, существует хорошо разработанная теория оракульной сложности задач выпуклой оптимизации~\cite{NemirYd79, Nemir_comlec95}, хорошо подтверждающаяся на практике, с другой стороны, для большинства методов первого порядка и большинства задач
наиболее вычислительно затратной частью итерации является именно расчёт градиента. Таким образом,
число обращений к оракулу отвечает за число итераций метода, что во многом определяет и общую сложность (время работы) метода.

Монография Немировского--Юдина~\cite{NemirYd79} стала в своё время (с конца \mbox{70-х} годов XX~века) настоящим
прорывом. Эта книга во многом и определила последующее развитие численных методов выпуклой оптимизации. Запас оригинальных идей, заложенных в данной (совсем непростой для чтения) книге, по-прежнему
вдохновляет большое число исследователей по всему миру. 

В последнее время концепция оракульной сложности
и чёрного ящика переживает новую
волну популярности  в связи с появлением большого числа новых постановок задач, приходящих из машинного обучения. Среди прочих, отметим особо работы последних 6-7 лет Натана Сребро (Nathan Srebro) с коллегами О.~Шамир (Ohad Shamir), Ю.~Аржевани (Yossi Arjevani), Б.~Вудворс (Blake Woodworth) и др.


Ещё раз подчеркнём принципы описанной
выше {\it концепции чёрного ящика}~ --- стандартного предположения, которое
необходимо для получения большинства результатов теории сложности задач оптимизации.
\begin{enumerate}
\item Единственное, что может
использовать в ходе своей работы
итерационный метод, это ответы
оракула.
\item Ответы оракула являются локальными: небольшое изменение задачи, произведённое достаточно далеко от тестовой точки~$x$ и согласованное с
описанием данного класса задач, не обязано привести к изменению исходного ответа в точке $x$.
\end{enumerate}
Эта концепция является одной из самых полезных изобретений в
численном анализе. В частности, она позволяет
получить соответствующие верхние и нижние
оценки оракульной сложности решения оптимизационных
задач для различных
подклассов выпуклых функций.


Опишем основные известные результаты об оптимальных вычислительных гарантиях (оценках скорости сходимости) для задач (в общем случае негладкой) выпуклой оптимизации~\cite{Polyak,NemirYd79}.
 
Мы начнём с задач небольшой размерности, когда возможна ситуация $N \geqslant n$, где $n$~--- размерность пространства, а $N$~--- количество вызовов оракула (например, число вычислений субградиента $f$ в текущей точке). Рассмотрим задачу выпуклой оптимизации
\begin{equation}
\label{p0_eq1}
\mathop {\min }\limits_{x \in Q} f\left (x \right),
\end{equation}
где $Q$~--- выпуклое компактное множество простой структуры. С использованием $N$ обращений к оракулу для субградиента ставится задача нахождения такой точки $x_N$, для которой
\begin{equation}
\label{p0_eq1_1}
f\left( {x_N} \right)-f(x^*) \leq \varepsilon,
\end{equation}
где $f(x^*)$~--- искомое минимальное значение функции в (\ref{p0_eq1_1}), $x^*$~--- точное решение задачи (\ref{p0_eq1}). Нижняя и верхняя оценки требуемого количества обращений к оракулу (с точностью до порядка множителя, имеющего логарифмическую зависимость от некоторой характеристики допустимого множества $Q$) равны $N\sim n\log \left( {{\Delta f} \mathord{\left/ {\vphantom {{\Delta f}\varepsilon }} \right. \kern-\nulldelimiterspace} \varepsilon } \right)$, где $\Delta f=\mathop {\sup }\limits_{x,y \in Q} \left\{ {f\left( y\right)-f\left( x \right)} \right\}$. Согласно указанной оценке сходится метод центров тяжести (метод внутренней точки). При $n=1$ этот метод является простым бинарным поиском \cite{Brent}. Однако при $n>1$ реализовать этот метод уже непросто. Это связано с тем, что сложность каждой итерации слишком высока, поскольку на каждой итерации требуется нахождение центра тяжести текущего допустимого множества. Как хорошо известно \cite{NemirYd79,Shor85}, для метода эллипсоидов достижение приемлемого качества решения необходимо $N= \ag{O}\left( {n^2\log \left( {{\Delta f} \mathord{\left/{\vphantom {{\Delta f} \varepsilon }} \right. \kern-\nulldelimiterspace}\varepsilon } \right)} \right)$ вызовов оракула при сложности каждой итерации $O\left( {n^2}\right)$. В \cite{Bubeck15,Vaidya} была предложена специальная версия метода отсекающих гиперплоскостей. Для такого метода \ag{требуется}\footnote{Здесь и всюду далее для всех (больших) $n$: $\widetilde{O}(g(n)) \leqslant C\cdot(\log n)^r g(n)$ с некоторыми константами $C>0$ и $r\geqslant 0$. Как правило, $r = 1$. Если $r=0$, то $\widetilde{O}(\cdot) = O(\cdot)$.} $N=\widetilde{O}\left( {n\log \left( {{\Delta f} \mathord{\left/{\vphantom {{\Delta f} \varepsilon }} \right. \kern-\nulldelimiterspace}\varepsilon } \right)} \right)$ вызовов оракула при сложности итерации $O\left(n^{\log_2 7}\right)$ (сложность обращения матрицы $n$ на $n$, равная по порядку сложности перемножения матриц таких размеров, практическим алгоритмом Штрассена, в теории можно предложить и лучший метод перемножения матриц $\widetilde {O}\left( {n^{2.37}} \right)$, однако на практике такой сложности пока не удалось достичь). В работе \cite{LeeSidford_15} предложен метод с $N=\widetilde{O}\left( {n\log \left( {{\Delta f} \mathord{\left/ {\vphantom {{\Delta f}\varepsilon }} \right. \kern-\nulldelimiterspace} \varepsilon } \right)}\right)$ вызовами оракула и сложностью итерации $\widetilde {O}\left( {n^2}\right)$. Этот метод пока не удалось сделать эффективным на практике, по-видимому, в теории это объясняется наличием больших порядков (степеней) логарифмов в $\widetilde{O}$. Подробнее о методах решения задач выпуклой оптимизации небольшой размерности будет рассказано в разделе \ref{smallD}.

Как видим, упомянутые выше оценки сложности существенно зависят от размерности задачи, что ставит под вопрос их применимость в случае большой размерности задачи $n$. Наиболее известный подход к пониманию эффективности алгоритмических процедур в многомерной оптимизации основан на теории сложности, восходящей к известной монографии А.~С.~Немировского и Д.~Б.~Юдина \cite{NemirYd79}. Разработанная А.~С.~Немировским и Д.~Б.~Юдиным теория нижних оценок возможной эффективности методов выпуклой минимизации градиентного типа для различных классов задач была подкреплена оптимальными методами, которые реализовывали эти оценки. Для класса методов, у которых на каждой итерации разрешается не более чем ${\rm O}\left( 1 \right)$ раз обращаться к оракулу (подпрограмме) для расчёта градиента $\nabla f\left( x \right)$, оценка числа итераций $N$, необходимых для достижения точности решения задачи $\varepsilon $ (по функции)~\eqref{p0_eq1_1}
указана в таблице~\ref{p0_tab001} в зависимости от класса рассматриваемых задач.
Как видим, указанные оценки сложности не зависят от размерности задачи. При этом нижние оценки в негладком случае ($f$ липшицева) достигаются на \textit{субградиентных методах} (см. раздел \ref{subgradient}), а в гладком случае ($f$ имеет липшицев градиент) достигаются на \textit{быстром {\rm(}ускоренном{\rm)} градиентном методе} (см. разделы \ref{gradientTaylor}, \ref{AM}, \ref{model}), предложенным Ю.\,Е.\,Нестеровым в 1983~г.~\cite{Nesterov83}.



\begin{table}[h!]
\centering
\caption{\centering\bf Оптимальные оценки количества обращений $N$ к (суб)градиенту}
\label{p0_tab001}
\begin{tabular}{|c|c|c|}
\hline
$N\leqslant n$&
\begin{tabular}[c]{@{}c@{}}
$|f(y)-f(x)|\leq $\\$\leq M\|y-x\|$\end{tabular}&
\begin{tabular}[c]{@{}c@{}}
${\|\nabla f(y)-\nabla f(x)\|}_\ast\leqslant $\\$\leqslant L\|y-x\|$\end{tabular}\\\hline
$f(x)$ выпукла&
$O\left(\frac{M^2R^2}{\varepsilon ^2}\right)$&
$O\left(\sqrt{\frac{LR^2}{\varepsilon}}\right)$\\\hline
\begin{tabular}[c]{@{}c@{}}
$f(x)$ $\mu$-сильно\\выпукла по $\|\cdot\|$-норме\end{tabular}&
$O\left(\frac{M^2}{\mu\varepsilon}\right)$&
\begin{tabular}[c]{@{}c@{}}
$O\left(\sqrt{\frac{L}{\mu}}\left\lceil\log \left(\frac{\mu R^2}{\varepsilon} \right)\right\rceil \right)\; (\forall N) $\end{tabular}\\\hline
\end{tabular}
\end{table}



После этого на некоторый период градиентные методы были практически забыты. Начиная с работы Н.\,Кармаркара\ag{~\cite{Karmarkar}} (1984 г.) и примерно до 2004~г.
развитие теории и методов оптимизации было в основном связано с прогрессом в теории полиномиальных методов внутренней точки. К примеру, была разработана общая теория самосогласованных функций, которая позволяла строить полиномиальные методы внутренней точки для всех выпуклых задач с явной структурой. Однако, как было отмечено выше, сложность методов внутренней точки сильно зависит (растёт) от размерности решаемой задачи. В связи с указанным обстоятельством в 2000-е годы возник интерес к оптимизационным методам градиентного типа, итерации которых требуют меньших затрат памяти, что интересно ввиду приложений к задачам оптимизации с большими данными. Тогда вновь возник интерес к теории сложности и подкрепляющим её методам градиентного типа. В рамках этой теории установлены оптимальные оценки числа вызовов оракула для задачи выпуклой оптимизации (\ref{p0_eq1}) вне зависимости от размерности задачи $n$, которые приведены в табл.~\ref{p0_tab001}. Отметим, что здесь уже $Q$~--- необязательно компактное множество. В табл.~\ref{p0_tab001} $R$~--- это <<расстояние>> (с точностью до множителя $\log n$) между начальной точкой и ближайшим решением $R=\widetilde {O}\left( \left\| {x_0-x^{*} } \right\| \right).$

\begin{remark}\label{ideas}
Заметим, что оценки в табл.~\ref{p0_tab001} в не сильно выпуклом случае могут быть получены исходя из оценок в сильно выпуклом случае, если положить $\mu\simeq \varepsilon/R^2$. Этому есть простое объяснение. Вместо задачи  \eqref{p0_tab001} предлагается рассматривать регуляризованную задачу (для простоты рассматриваем случай $\|~\cdot~\| = \|~\cdot~\|_2$)
\begin{equation*}
\label{p0_eq1*}
\mathop{\min }\limits_{x \in Q} f_{\mu}\left (x \right):= f\left (x \right) +\frac{\mu}{2}\|x_0 - x^*\|_2^2 .
\end{equation*}
Если $\mu\simeq \varepsilon/R^2$ и
\begin{equation*}
\label{p0_eq1_1*}
f_{\mu}\left(x_N\right)- \min_{x \in Q}f_{\mu}(x) \leqslant \varepsilon/2,
\end{equation*}
то выполняется \eqref{p0_eq1_1}:
\begin{equation*}
\label{p0_eq1_1**}
f\left(x_N\right)- \min_{x \in Q}f(x) = f\left( {x_N} \right)- f(x^*)  \leqslant \varepsilon.
\end{equation*}

Интересно отметить, что существует универсальная конструкция (рестарты) \cite{Aspremont2021}, которая в некотором смысле является  обратной к описанной выше. Продемонстрируем это на примере получения оценки в гладком сильно выпуклом случае из оценки в гладком выпуклом случае. Предположим, что у нас есть метод, который сходится следующим образом:
\begin{equation*}
f\left({x_N} \right) - f(x^*)  \lesssim \frac{L\|x_0 - x^*\|_2^2}{N^2}.
\end{equation*}
Тогда, выбирая $\bar{N} = 2\sqrt{L/\mu}$ и используя сильную выпуклость $f$  из
\begin{equation*}
\frac{\mu}{2}\|x_{\bar{N}} - x^*\|_2^2\le f\left(x_{\bar{N}}\right) - f(x^*)  \lesssim \frac{L\|x_0 - x^*\|_2^2}{N^2},
\end{equation*}
получим, что 
\begin{equation*}
\|x_{\bar{N}} - x^*\|_2^2\le\frac{1}{2}\|x_0 - x^*\|_2^2.
\end{equation*}
Делая логарифмическое по желаемой точности число рестартов (присваиваний $x_0:=x_{\bar{N}}$), получим требуемую оценку. Аналогичные (но чуть более громоздкие) рассуждения можно провести и в негладком случае. Более подробно написанное выше разбирается далее в разделе~\ref{restarts}.

Отметим также, что оценки негладкого случая можно получить из соответствующих оценок гладкого случая, используя результаты о накоплении неточности в рассматриваемых методах. А именно, приведённые в таблице оценки в гладком случае могут быть получены при немного более слабом условии, которое следует из $L$-липшицевости градиента: для всех $x,y\in Q$ выполняется
\begin{equation}\label{smth}
f(y)\le f(x) + \langle \nabla f(x), y-x \rangle + \frac{L}{2}\|y-x\|^2.
\end{equation}
Если допустить, что \eqref{smth} выполняется неточно 
\begin{equation}\label{dlt}
f(y)\le f(x) + \langle \nabla f(x), y-x \rangle + \frac{L}{2}\|y-x\|^2 + \delta,
\end{equation}
то как это изменит скорость сходимости быстрого (ускоренного) градиентного метода оптимального для задач гладкой выпуклой оптимизации? Ответ даёт следующая формула (детали см. в разделе \ref{fastGradMethod}):
\begin{equation}\label{inex}
f(x_N) - f(x_*) \lesssim \frac{LR^2}{N^2} + N\delta.
\end{equation}
Ключевое наблюдение, которое позволяет сводить анализ негладких задач к гладким, заключается в наблюдении, что $M$-липшицева функция (такой функцией, например, будет негладкая скалярная функция $f(x) = M|x|$) для любого $\delta > 0$ удовлетворяет \eqref{dlt} с $L = M^2/(2\delta)$.

Приравнивая каждое слагаемое в правой части \eqref{inex} $\varepsilon/2$, получим, что $N\simeq \sqrt{LR^2/\varepsilon}$, $\delta \simeq \varepsilon/N$. Учитывая далее, что $L = M^2/(2\delta)$, получим, $N \simeq M^2R^2/\varepsilon^2.$ Аналогичные рассуждения можно провести и в сильно выпуклом случае.
\end{remark}

%
%
%
\subsection*{Обозначения}
Как и ранее, будем обозначать через $x^*$ 
такую точку, что $x^* \, \in \, Q$ и 
$x^* = \min\limits_{x \, \in \, Q} f(x)$
(т.е. $x^*$~--- оптимальное
решение задачи). Если не оговорено иное,
будем предполагать, что $x^*$
существует.
Оптимальное значение функции $f$ 
на множестве\footnote{Множество
$Q$~--- область априорной локализации
минимума. Если в практической задаче
$Q$ совпадает со всем пространством $\mathbb{R}^n$,
то можно указать, например, возможный диапазон
изменения по каждому измерению, и этот параллелепипед
принять за~$Q$.} $Q$
обозначим через $f^*$ (т.е. $f^* = f(x^*)$).

%
%
\section[Методы выпуклой оптимизации для
задач небольшой\\ размерности]{Методы выпуклой оптимизации для
задач\\ небольшой размерности%
\sectionmark{Методы для задач небольшой размерности}}
\sectionmark{Методы для задач небольшой размерности}\label{smallD}
Пусть $Q \, \subset \, \mathbb{R}^n$~--- выпуклое тело
(т.е. выпуклое компактное множество
с непустой внутренностью), $f \, : \, Q \, \rightarrow \, \mathbb{R}$~--- выпуклая непрерывная функция.
Далее в этом параграфе описаны методы оракульного типа
для решения задачи
\begin{equation}
\label{ch1_eq_main_prob}
\min_{x \, \in \, Q} f(x).
\end{equation}
Все представленные в этом параграфе методы имеют линейную
или квадратичную по отношению к размерности~$n$ оценку на 
число обращений к оракулу, поэтому их можно
применять только для задач не очень большой размерности.

Другой важной особенностью обсуждаемых в этом параграфе
методов является то, что им для работы
требуется так называемый {\it отсекающий оракул} для множества ограничений~$Q$. 
Для любого~$x \, \in \, \mathbb{R}^n$ отсекающий
оракул для $Q$ либо сообщает, что $x \, \in \, Q$,
либо строит отделяющую гиперплоскость между $x$ и $Q$.
Такие методы в литературе часто объединяют
под названием {\it методов отсекающей гиперплоскости}\label{ch2_idea_separ_hyper}.

\subsection{Методы решения одномерных задач}\label{subs:bisection}

Задача минимизации выпуклой функции от одной вещественной переменной встречается как подзадача поиска оптимальной длины шага в итеративных методах минимизации функции от многих переменных, например, в градиентном спуске. Иногда необходимо отминимизировать функцию цены по скалярному множителю Лагранжа. Здесь мы рассмотрим два метода для решения одномерных задач, в зависимости от того, доступен ли знак производной функции или её значение.

Представленные здесь методы на самом деле применимы для минимизации более широкого класса функций. А именно, достаточно потребовать, чтобы существовала строго монотонно возрастающая функция $g: \mathbb R \to \mathbb R$ такая, что композиция $g \circ f$ выпукла.

\subsubsection{Метод бисекции (деления отрезка пополам)} 

Мы хотим минимизировать выпуклую функцию $f: \mathbb R \supset I \to \mathbb R$, где $I$~--- интервал, не обязательно конечный. Мы предполагаем, что точка глобального минимума $x^*$ существует и что для произвольной точки $x$ существует способ определить, какое из неравенств $x \leq x^*$ и $x \geq x^*$ имеет место, т.е., доступен оракул, отвечающий на этот вопрос. Отметим, что если в точке $x$ верны оба неравенства одновременно, то $x$ является точкой глобального минимума, независимо от единственности этого минимума. Если функция $f$ дифференцируема, то оракул просто возвращает знак производной $f'(x)$.

Метод начинает с некоторой точки $x_0$, которую он передаёт оракулу. Допустим без ограничения общности, что оракул возвращает ответ $x_0 < x^*$. Выберем константу $d > 0$ и итеративно вычисляем последовательность точек $x_{k+1} = x_k + 2^kd$, которые мы также передаём оракулу. Через $N = \left\lceil\frac{\log\left(1 + \frac{|x^*-x_0|}{d}\right)}{\log 2}\right\rceil$ шагов мы получаем точку $x_N$, для которой верны неравенства $x_{N-1} < x^* \leq x_N$ и соотношение $|x_N - x_{N-1}| = 2^{N-1}d$.

Переходим ко второй фазе алгоритма, которая начинает с конечного интервала $I_0 = [x_{N-1},x_N]$ и итеративно строит последовательность вложенных друг в друга интервалов $I_k \subset I_{k-1}$ с центрами $y_k$. На каждом шаге следующий интервал $I_k$ определяется либо как левая, либо как правая половина интервала $I_{k-1}$, в зависимости от того, справедливо ли неравенство $y_{k-1} \geq x^*$ или $y_{k-1} \leq x^*$. Сделав $\left\lceil\log_2\frac{|x_N - x_{N-1}|}{\epsilon}\right\rceil = \left\lceil\log_2\frac{2^{N-1}d}{\epsilon}\right\rceil$ шагов, мы получаем интервал длины не более $\epsilon$, который содержит глобальный минимум функции.

Таким образом, общее количество обращений к оракулу составляет
\[ 1 + N + \left\lceil\log_2\frac{2^{N-1}d}{\epsilon}\right\rceil < 3 + \frac{2\log\left(d + |x^*-x_0|\right) - \log d - \log\epsilon}{\log 2}.
\]
Минимум оценки справа достигается при $d = |x^*-x_0|$. Поэтому начальный шаг $d$ следует выбрать порядка расстояния до минимума.

Количество значимых знаков минимума $x^*$ возрастает пропорционально числу итераций. Поэтому метод бисекции имеет линейную сходимость.

\subsubsection{Метод золотого сечения}

Мы решаем ту же задачу, что и в предыдущем разделе, но в этом случае оракул выдаёт значение $f(x)$ функции в заданной точке $x$. По-прежнему мы полагаем, что глобальный минимум $x^*$ выпуклой функции $f$ существует. Для того чтобы ограничить область, в которой лежит минимум $x^*$, знания значения функции в одной лишь точке недостаточно. Однако это возможно, если мы имеем в распоряжении значения функции в двух точках $x_0 < x_1$. Действительно, допустим без ограничения общности, что $f(x_0) > f(x_1)$. Тогда $x^*$ не может лежать левее $x_0$, и мы можем вывести ограничение $x^* > x_0$.

Первая фаза метода золотого сечения аналогична первой фазе метода бисекции. Определим число золотого сечения $\alpha = \frac{\sqrt{5}-1}{2} \approx 0.618$, являющееся корнем квадратичного уравнения $\alpha^2 + \alpha - 1 = 0$. Положим $\lambda = 1 + \alpha$, $d = x_1 - x_0$ и итеративно построим последовательность точек $x_{k+1} = x_k + d\cdot\lambda^k$, пока после $N = \left\lceil\frac{\log\left(1 + \frac{|x^*-x_0|(\lambda - 1)}{d}\right)}{\log\lambda}\right\rceil$ шагов не получим точку $x_N$, удовлетворяющую условию $f(x_N) \geq f(x_{N-1})$. Тогда три последние точки $x_{N-2} < x_{N-1} < x_N$ последовательности удовлетворяют условию $f(x_{N-1}) < \min(f(x_{N-2}),f(x_N))$. Отсюда следует, что минимум $x^*$ лежит в конечном интервале $I_0 = [x_{N-2},x_N]$ длины $\lambda^2(1+\lambda)d = \lambda^Nd$. Отметим также, что вследствие выбора коэффициента $\lambda$ имеем $\frac{x_N - x_{N-1}}{x_N - x_{N-2}} = \frac{\lambda}{1+\lambda} = \alpha$.

Во второй фазе метода мы строим последовательность вложенных друг в друга интервалов $I_{k+1} \subset I_k$, содержащих минимум $x^*$. На каждом шаге мы уже имеем значения функции $f$ в крайних точках $x_{k-} < x_{k+}$ интервала $I_k$ в нашем распоряжении, а также значение $f$ в некоторой внутренней точке $\tilde x_k \in I_k$, удовлетворяющей условию $\frac{x_{k+} - \tilde x_k}{x_{k+} - x_{k-}} \in \{\alpha,1-\alpha\}$, т.е. $\tilde x_k$ разделяет интервал $I_k$ в пропорции золотого сечения. Более того, значения функций удовлетворяют условию $f(\tilde x_k) \leq \min\left(f(x_{k-}),f(x_{k+})\right)$. Определим следующую точку, которую мы передадим оракулу, как $x'_{k+1} = x_{k-} + x_{k+} - \tilde x_k$, т.е. симметричной к $\tilde x_k$ по отношению к центру интервала $I_k$. Отсюда следует, что точка $x'_{k+1}$ делит интервал $I_k$ также в пропорции золотого сечения. Следующий интервал $I_{k+1}$ строится следующим образом. Если $f(x'_{k+1}) > f(\tilde x_k)$, то $x'_{k+1}$ становится крайней точкой интервала $I_{k+1}$, а точка $\tilde x_k$ становится новой внутренней точкой $\tilde x_{k+1}$ интервала $I_{k+1}$. В противном случае крайней точкой становится $\tilde x_k$, а внутренней точкой $\tilde x_{k+1}$ становится $x'_{k+1}$. Вторая крайняя точка интервала $I_{k+1}$ совпадает с соответствующей крайней точкой интервала $I_k$. Новая внутренняя точка $\tilde x_{k+1}$ снова делит интервал $I_{k+1} = [x_{k+1-},x_{k+1+}]$ в пропорции золотого сечения.

На каждом шаге длина интервала умножается на $\alpha < 1$ и поэтому убывает в геометрической прогрессии. Через $\left\lceil\log_{\alpha^{-1}}\frac{|x_N - x_{N-2}|}{\epsilon}\right\rceil$ шагов мы получаем интервал длины не более $\epsilon$, содержащий минимум $x^*$. Общее количество обращений к оракулу составляет
\[ 1 + N + \left\lceil\log_{\alpha^{-1}}\frac{\lambda^Nd}{\epsilon}\right\rceil < 4 + \frac{2\log\left(d + \alpha|x^*-x_0|\right) - \log d - \log\epsilon}{\log\lambda}.
\]

Метод золотого сечения также имеет линейную сходимость.

%
%
\subsection[Метод центров тяжести]
{Метод центров тяжести%
\subsectionmark{Метод центров тяжести}}
Важной отправной точкой развития методов выпуклой оптимизации в 60-е годы XX~века стал предложенный в 1965~г. А.Ю. Левиным метод центров
тяжести~\cite{Levin_65}. Независимо от А.\,Ю.~Левина, в этом же году, метод центров тяжести в США
разработал Д. Ньюман~\cite{Newman_gravity}.

Метод базируется на теореме Грюнбаума--Хаммера--Митягина\footnote{Данная теорема более известна как теорема Грюнбаума--Хаммера~\cite{Grun} (1960 г.), но
несколько лет спустя она была независимо переоткрыта
Борисом Самуиловичем Митягиным и в 1969~году опубликована в журнале <<Математические заметки>>~\cite{Mityagin}.
}. Согласно этой теореме, любая гиперплоскость,
проведённая через {\it центр тяжести} выпуклого 
ограниченного тела в $\mathbb{R}^n$, разбивает его на два множества,
отношение объёмов которых не меньше
$$
\frac{\left ( \frac{n}{n+1} \right )^n}{1 - \left ( \frac{n}{n+1} \right )^n} \geq \frac{1}{e-1}.
$$
Иными словами, объём части, отсекаемой от $Q$
гиперплоскостью, проходящей через центр тяжести,
составляет всегда не менее чем $e^{-1}$ часть
от всего объёма~$Q$.

Этот факт позволяет перенести метод дихотомии
(деления отрезка пополам)
решения одномерной задачи минимизации 
выпуклой функции на отрезке (с использованием значения
производной или субградиента) на многомерный случай.
Если выпуклая функция~$f$ задана на выпуклом
множестве~$Q$ и в некоторой точке~$x_0$ вычислен её субградиент~$g$, то для всех $x \, \in \, Q$, для
которых 
\begin{equation}
\label{ch2_eq_grun_grad_prop}
g^T (x - x_0) > 0,
\end{equation}
функция минимума не достигает.
По теореме Грюнбаума--Хаммера--Ми\-тя\-гина, основанные
на идее~\eqref{ch2_eq_grun_grad_prop} методы отсекающей
гиперплоскости позволяют приближаться к минимуму со скоростью геометрической прогрессии.

Рассмотрим следующую итерационную схему для
метода центров тяжести (см. алгоритм~\ref{ch2_alg_CGM},  $\text{vol}(Q)$~--- объём множества~$Q$).

\floatname{algorithm}{Алгоритм}
	\begin{algorithm}
		\caption{Метод центров тяжести}\label{ch2_alg_CGM}
		\begin{algorithmic}[1]
			\REQUIRE $f$~--- выпуклая функция; $Q$~--- допустимое множество; $N+1$~--- количество итераций.
			\ENSURE точка $x_N$
\STATE $Q_0 \leftarrow Q$
\FOR{$k=0, \, \dots, \, N$}
\STATE Вычислить 
$$
c_k \leftarrow \frac{1}{\mbox{vol}(Q_k)} \int_{x \, \in \, Q_k} x\,dx
$$
\STATE Вычислить значение какого-либо
из субградиентов $g_k \, \in \, \partial f(c_k)$
\STATE 
$$
Q_{k+1} \leftarrow Q_k \cap \{ x \, \in \, \mathbb{R}^n 
\, | \, g_k^T (x - c_k) \le 0 \}.
$$
\ENDFOR

\STATE $x_N \leftarrow \argmin\limits_{0 \, \le \, r \, \le \, N} f(c_r)$

\STATE \RETURN $x_N$
\end{algorithmic}

\end{algorithm}

\begin{teo}
\label{ch2_teo_cgm}
Для метода центров тяжести выполняется следующая оценка:
$$
f(x_N) - f^* \le C \left ( 1 - \left ( \frac{n}{n+1} \right )^n \right )^{N/n} \le  
C \left ( 1 - \frac{1}{e} \right )^{N/n},
$$
где $C = \max\limits_{x \, \in \, Q} \, (f(x) - f^*)$.
\end{teo}

\begin{proof}
Доказательство теоремы основано на работе~\cite{Polyak}.
По теореме Грюнбаума--Хаммера--Митягина:
$$
\text{vol}(Q_{N}) \le   \left ( 1 - \left ( \frac{n}{n+1} \right )^n \right ) \text{vol}(Q_{N-1}).
$$
Следовательно,
$$
\text{vol}(Q_{N}) \le   \left ( 1 - \left ( \frac{n}{n+1} \right )^n \right )^{N} \text{vol}(Q_{0}).
$$
Возьмём произвольную точку минимума $x^* \, \in \, Q_{N}$
и построим множество~$S$, получающееся из $Q_{N}$
подобным преобразованием с центром в $x^*$ и
коэффициентом растяжения $\alpha = \left ( 1 - \left ( \frac{n}{n+1} \right )^n \right )^{-N/n}$.
Тогда 
$$
\text{vol}(S)=\alpha^n \text{vol}(Q_{N}) \le
\alpha^n \left ( 1 - \left ( \frac{n}{n+1} \right )^n \right )^{N} \text{vol}(Q_{0}) = \text{vol}(Q_{0}),
$$
поэтому множество~$Q_0$ не может помещаться строго внутри $S$,
а, значит, найдётся $z \, \in \, Q_0$, $z \, \notin \, S$.
Отсюда следует, что $u = \left (1 - \alpha^{-1} \right ) x^* + \alpha^{-1} z$ не должно принадлежать $Q_N$ 
(поскольку $z = \alpha (u - x^*) + x^*$, т.е. $z$ получается
из $u$ указанным выше растяжением).
Тогда, по построению $Q_N$ (см. алгоритм~\ref{ch2_alg_CGM}),
найдётся такое $i$, $0 \le i \le N$, что $g_i \, \in \, \partial f(x_i)$
$$
g_i^T (u - x_i) \ge 0.
$$
По определению субградиента и последнему неравенству
$$
f(u) \ge f(x_i) +  g_i^T (u - x_i) \ge f(x_i) \ge f(x_N).
$$
В силу выпуклости $f(x)$ по неравенству Йенсена:
$$
f(x_N) \le f(u) = f \left ( \left (1 - \alpha^{-1} \right ) x^* + \alpha^{-1} z \right ) \le \left (1 - \alpha^{-1} \right ) f^* + \alpha^{-1} f(z) \le f^* + C/\alpha,
$$
где $C = \max\limits_{x \, \in \, Q_0} \, (f(x) - f^*)$,
причём $C < \infty$ в силу непрерывности $f(x)$ и ограниченности $Q$.
Итак, 
$$
f(x_N) - f^* \le C \alpha^{-1} = \left ( 1 - \left ( \frac{n}{n+1} \right )^n \right )^{N/n} \le  
C \left ( 1 - \frac{1}{e} \right )^{N/n},
$$
что и требовалось доказать. 
\end{proof}

Итак, для того, чтобы найти решение задачи~\eqref{ch1_eq_main_prob} при заданной точности
нахождения решения $\varepsilon$ с помощью метода центров
тяжести, необходимо выполнить не более $O(n \log (C / \varepsilon) )$ обращений к оракулам первого и нулевого порядков.

В теореме~\ref{ch2_teo_cgm}
утверждается, что скорость сходимости метода центров тяжести
является {\it линейной}. Это означает,
что оценка ошибки метода на итерации $k+1$
зависит от оценки ошибки метода на предыдущей итерации {\it линейно}.

Важное замечание следует сделать о вычислительной
сложности метода: для $n > 2$ отыскание центра тяжести~$c_k$ множества~$Q$ превращается в чрезвычайно громоздкую задачу~--- в настоящее время не существует эффективной
{\bf нерандомизированной}\footnote{Рандомизацией детерминированного метода называется
модификация метода посредством искусственного введения случайности в алгоритм решения.} вычислительной процедуры решения такой задачи.  

Тем не менее метод центров тяжести представляет большой теоретический интерес. Во-первых, оценка скорости сходимости метода зависит только от размерности пространства и величины <<начальной неопределённости>> целевой функции $\left ( \max\limits_{x \, \in \, Q} \, f(x) - \min\limits_{x \, \in \, Q} \, f(x) \right )$, но не от других характеристик функции типа её обусловленности.
Но самое важное свойство метода заключается в том, что метод центров тяжести является в определённом смысле оптимальным.

%
%
\subsubsection{Оптимальные методы}\label{ch2_subsubsect_optim}
Как уже обсуждалось выше, для задач минимизации
выпуклых функций можно установить потенциальные
возможности любого метода, использующего
лишь оракул определённого порядка ({\it нижние оценки}).
Так, для задач негладкой выпуклой, негладкой сильно
выпуклой и гладкой выпуклой оптимизации
на множествах простой структуры в $\mathbb{R}^n$,
когда число итераций $N \ge n$,
нижние оценки числа обращений к оракулу первого
порядка имеют одинаковый (с точностью до числового множителя~$C$) вид
\begin{equation}
\label{ch2_eq_optim_low_bound}
N \ge C n \log \left ( \frac{1}{\varepsilon} \right ),
\end{equation}
$\varepsilon$~--- относительная точность решения задачи по функции
(доказательство см.~в монографии А.\,С. Немировского и Д.\,Б. Юдина~\cite{NemirYd79}).

Сопоставляя этот факт с полученной выше оценкой скорости
сходимости метода центров тяжести (см. теорему~\ref{ch2_teo_cgm}), получаем, что
{\it скорость сходимости метода центров тяжести не может быть по порядку превзойдена ни для какого метода оптимизации,
решающего задачу выпуклой оптимизации и использующего оракул первого порядка}.
Иными словами, метод центров тяжести является неулучшаемым (т.е. оптимальным) при конечной размерности, и попытки создать более быстро сходящиеся методы будут заведомо неудачными.
Но также здесь нужно пояснить, что данный вывод относится
к широкому классу <<всех выпуклых функций>>.
Для более узких классов функций и/или для функций с конкретной заранее известной структурой (например, $\max\limits_{1 \, \le \, i \, \le \, m} f_i(x)$) могут существовать и более эффективные методы. 

Ещё одно важное замечание, касающееся оценки~\eqref{ch2_eq_optim_low_bound} и оптимальности
метода центров тяжести,
состоит в том, что всё это верно лишь для сходимости по функции
в задачах на выпуклых компактах с оракулом, наделённым дополнительными нетривиальными возможностями.

Хотя метод центров тяжести не получил практического распространения из-за сложности задачи вычисления центра тяжести в многомерном пространстве, его основная
идея применима к построению методов, имеющих большое
практическое значение. Например, к {\it методу эллипсоидов}.

%
%
\subsection{Метод эллипсоидов}\label{ch2_subsect_ellips_meth}
Метод эллипсоидов является модификацией метода центров тяжести (МЦТ), избавленной от главного недостатка МЦТ~--- трудоёмкая операция отыскания центра тяжести выпуклого
множества заменяется на задачу нахождения центра
описанного эллипсоида. Но платой за простоту итерации метода эллипсоидов оказалось значительное увеличение необходимого числа обращений
к оракулу для достижения заданной точности~$\varepsilon$. Методу эллипсоидов требуется $O(n^2 \log (1 / \varepsilon) )$ обращений к оракулу (см. ниже теорему~\ref{ch2_teo_ellips_cong}), в то
время как оптимальному методу (таким
является МЦТ, см. п.~\ref{ch2_subsubsect_optim}) необходимо $O(n \log (1 / \varepsilon) )$  обращений к оракулу. Для больших размерностей пространства проигрыш в скорости сходимости значителен.

Напомним определение эллипсоида.
\begin{defin}
Эллипсоидом называется выпуклое множество следующей формы:
$$
\mathcal{E} = \left\{ x \, \in \, \mathbb{R}^n \, | \, (x-c)^T H^{-1} (x - c) \le 1 \right\},
$$
где $c \, \in \, \mathbb{R}^n$, $H$~--- симметричная положительно определённая матрица.
\end{defin}

С геометрической точки зрения $c$~--- центр эллипсоида,
полуоси~$\mathcal{E}$ являются собственными векторами
$H$,  
а длины полуосей~$\mathcal{E}$ равны квадратным
корням из соответствующих собственных значений матрицы $H$.

Метод эллипсоидов основан на двух идеях~--- идее отсекающей гиперплоскости (см. стр.~\pageref{ch2_idea_separ_hyper})
с добавлением следующего геометрического факта: {\it половину эллипсоида можно поместить в эллипсоид меньшего,
чем объём изначального эллипсоида, объёма}. Центр нового
эллипсоида можно вычислить, затратив порядка $O(n^2)$ операций\footnote{То есть стоимость (трудоёмкость) итерации метода эллипсоидов $O(n^2)$. В эту оценку не входит расчёт субградиента, но обычно (суб-)градиент можно посчитать за $O(n^2)$ (см., например, \cite[\S $\,$ 2]{Gas_Pos18}, поэтому можно считать,
что $O(n^2)$~--- оценка общей сложности итерации.}.

Опишем алгоритм метода эллипсоидов\footnote{Обратите внимание на сходство формул метода
эллипсоидов с формулами квазиньютоновского
метода с матрицей--аппроксимацией гессиана $H_k^{-1}$.
И при этом для метода эллипсоидов не требуется дифференцируемость целевой функции!} (см. алгоритм~\ref{ch2_alg_ellips}).

\floatname{algorithm}{Алгоритм}
	\begin{algorithm}
		\caption{Метод эллипсоидов}\label{ch2_alg_ellips}
		\begin{algorithmic}
			\REQUIRE $f$~--- выпуклая функция; $Q$~--- допустимое множество; $x_0$~--- начальная точка; $N+1$~--- количество итераций.
			\ENSURE точка $x_N$
\STATE $c_0 \leftarrow x_0$			
\STATE $\mathcal{E}_0$ --- евклидов шар радиуса~$R$
с центром $c_0$, который
содержит~$Q$
\STATE $H_0 \leftarrow R^2 I_n$, где $I_n$~--- единичная матрица размера $n \times n$ 
\FOR{$k=0, \, \dots, \, N-1$}
\IF{$c_k \, \notin \, Q$}
\STATE С помощью оракула получить 
отделяющую гиперплоскость $g_k \, \in \, \mathbb{R}^n$,
такую, что $Q \subset \{ x \, | \, (x-c_k)^T g_k  \le 0 \}$
\ELSE
\STATE Вычислить значение какого-либо
из субградиентов $g_k \, \in \, \partial f(c_k)$
\ENDIF

\STATE
Пусть $\mathcal{E}_{k+1}$~--- эллипсоид минимального объёма, содержащий
$\{ x \!\! \in \linebreak \in \, \mathcal{E}_{k} \, | \, (x-c_k)^T g_k  \le 0 \}$, т.е.
$$
\mathcal{E}_{k+1} \leftarrow \{ x \, \in \, \mathbb{R}^n \, | \, (x-c_{k+1})^T H_{k+1}^{-1} (x - c_{k+1}) \le 1 \},
$$
где
$$
c_{k+1} = c_k - \frac{1}{n+1} \frac{H_k g_k} {\sqrt{g_k^T H_k g_k}},
$$
$$
H_{k+1} = \frac{n^2}{n^2 - 1} \left ( H_k - \frac{2}{n+1}
\frac{H_k g_k g_k^T H_k}{g_k^T H_k g_k} \right )
$$
\ENDFOR

\STATE 
\begin{equation}
\label{ch2_eq_ellip_xN}
x_N \leftarrow \argmin\limits_{c \, \in \, \{ c_0, \, \ldots, \, c_N\} \, \cap \, Q} f(c)
\end{equation}

\STATE \RETURN $x_N$
\end{algorithmic}

\end{algorithm}

\begin{exercise}
Проверьте, что при $n = 1$ метод эллипсоидов превращается в метод деления отрезка пополам (дихотомии) на $\mathbb{R}$.
\end{exercise}
 Можно даже сказать, что метод эллипсоидов является обобщением метода дихотомии
на многомерные пространства (что
подтверждает и оценка требуемого
числа итераций для достижения заданной точности $\varepsilon$~--- для метода эллипсоидов она имеет линейную зависимость от $\log(1/\varepsilon)$, см. ниже теорему~\ref{ch2_teo_ellips_cong}).

Оценим эффективность метода эллипсоидов.
\begin{teo}
\label{ch2_teo_ellips_cong}
Пусть функция $f$ липшицева на евклидовом шаре\protect\footnotemark радиуса~$R$, c центром в точке~$x_0$, $B_2(x_0, \, R)$ с некоторой константой $M$
и $x^* \, \in \, B_2(x_0, \, R)$.
Тогда для того, чтобы найти $f^*$ с точностью $\varepsilon$ с помощью
метода эллипсоидов, необходимо выполнить не более $2n^2 \log \frac{MR}{\varepsilon}$ обращений к оракулу.
\end{teo}
\footnotetext{Здесь и далее {\it шар} (или {\it шаровое множество}) вида
$$
\{ x \, \in \, \mathbb{R}^n \, | \, \|x - x_0 \| \le r \}, \quad r \ge 0,
$$
где $r$~--- радиус шара, $x_0 \, \in \, \mathbb{R}^n$~---
центр шара, будем обозначать $B_{\|\cdot\|}(x_0, \, r)$.
Евклидов шар с центром в точке $x_0$ радиуса~$r$ в таких обозначениях будет записываться как $B_2(x_0, \, r)$.}

Доказательство теоремы основано на работах~\cite{Boyd_ellip_method, Nest}.

\begin{proof}
Даже если $\mathcal{E}_{k+1}$ (c центром в точке $c_{k+1}$) будет больше, чем
$\mathcal{E}_{k}$ (с центром в точке $c_k$) по максимальной полуоси (т.е.  $\lambda _{\max } \left( H_{k+1} \right) \, >\linebreak > \, \lambda _{\max } \left( H_{k} \right)$, это возможно), объём $\mathcal{E}_{k+1}$ всё равно будет меньше:
$$
\frac{\text{vol}(\mathcal{E}_{k+1})}{\text{vol}(\mathcal{E}_{k})} = \left ( \frac{\det H_{k+1}}{\det H_k} \right )^{1/2} = \left ( \left ( \frac{n^2}{n^2 - 1} \right )^n \frac{n-1}{n+1} \right )^{1/2}.
$$
Получаем, что
\begin{align*} 
2n\log\frac{\text{vol}(\mathcal{E}_{k+1})}{\text{vol}(\mathcal{E}_{k})} &= -n^2\log\left ( 1 -  \frac{1}{n^2} \right) + n\log\left( 1 - \frac{1}{n} \right) - n\log\left( 1 + \frac{1}{n} \right)= \\ 
=& \left( 1 + \frac{1}{2n^2} + \frac{1}{3n^4} + \frac{1}{4n^6} + \dots \right) - \left( 2 + \frac{2}{3n^2} + \frac{2}{5n^4} + \frac{2}{7n^6} + \dots \right) \\ 
=& - 1 - \frac{1}{2\cdot 3\cdot n^2} - \frac{1}{3 \cdot 5 \cdot n^4} - \frac{1}{4 \cdot 7 \cdot n^6} - \dots < -1.
\end{align*}
Ряд сходится при всех $n \geq 1$, поэтому последнее неравенство также справедливо при всех $n \geq 1$. Из этого следует, что 
\begin{equation}
\label{ch2_eq_ellip_vol_ek1k}
\text{vol}(\mathcal{E}_{k+1}) < e^{-1/2n} \text{vol}(\mathcal{E}_{k}).
\end{equation}

Пусть $f(x_N) - f^* \ge \varepsilon$, где $\varepsilon > 0$. Тогда для $k = 0, \, \ldots, \, N$ $f(c_k) \ge f^* + \varepsilon$ (см. \eqref{ch2_eq_ellip_xN}), и для каждой
точки $x$ из отсекаемой на итерациях $0, \, \ldots, \, N-1$ части эллипсоида выполняется неравенство $f(x) \ge f^* + \varepsilon$.

Константа Липшица $M$ для $f$ на шаре $B_2(x_0, \, R)$ из условия теоремы является
также максимальной из всех норм субградиентов $f$ в $\mathcal{E}_0 = B_2(x_0, \, R)$.
В силу того, что $M$~--- константа Липшица для $f$, для любых $x$ из шара $B_2\left(x^*, \, \frac{\varepsilon}{M}\right)$
выполняется неравенство $f(x) - f^* \le \varepsilon$.
Не теряя общности, предположим, что
$B_2\left(x^*, \, \frac{\varepsilon}{M}\right) \subseteq Q$,
и, следовательно,
$$
B_2\left(x^*, \, \frac{\varepsilon}{M}\right) \subseteq \mathcal{E}_N
$$
и $\text{vol} \left ( B_2\left(x^*, \, \frac{\varepsilon}{M}\right) \right ) \le 
\text{vol}(\mathcal{E}_N)$.
В силу~\eqref{ch2_eq_ellip_vol_ek1k}
$$
\frac{\text{vol}(\mathcal{E}_N)}{\text{vol}(\mathcal{E}_0)} \le 
e^{-N/2n},
$$
поэтому
$$
\frac{\text{vol}\left ( B_2\left(x^*, \, \frac{\varepsilon}{M}\right) \right )}{\text{vol}(\mathcal{E}_0)} \le 
e^{-N/2n},
$$
или
$\frac{\varepsilon^n R^n}{M^n} \le 
e^{-N/2n}$.
Логарифмируя последнее неравенство, получаем
$$
- \frac{N}{2n} \ge n \log \frac{\varepsilon}{M} - n \log R,
$$
откуда
$$
N \le 2n^2 \log \frac{MR}{\varepsilon},
$$
т.е. для того, чтобы найти $f^*$ с точностью $\varepsilon$, необходимо выполнить $O\left(n^2 \log (1 / \varepsilon) \right)$ итераций (или обращений к оракулу)
метода эллипсоидов.
\end{proof}

Метод эллипсоидов был впервые предложен в Москве Д.\,Б.~Юдиным и
A.\,C.~Не\-ми\-ровским~\cite{Ydin_ellipse} в 1974--1976 гг.\footnote{Про историю открытия метода центров тяжести,
метода эллипсоидов, симплекс-метода можно прочитать
в интересной статье В.М.~Тихомирова ''The Evolution of Methods of
Convex Optimization''~\cite{Tikh_history}.}  и независимо переоткрыт в Киеве
Н.\,З.~Шором~\cite{Shor_ellipse} в 1977 г.

Интересно, что Н.З.~Шор пришёл к методу эллипсоидов из других соображений, путём модификации субградиентного метода.
В методе Шора субградиентный алгоритм объединяется с процедурой {\it растяжения пространства} переменных. В первых вариантах метода растяжение пространства проводилось в направлении последнего субградиента,
во втором семействе таких методов, которые получили название {\it $r$-алгоритмы}, растяжение пространства происходит в направлении
разности двух последовательных субградиентов. Тогда 
соответствующие антисубградиенты в растянутом пространстве с достаточным коэффициентом растяжения\footnote{Коэффициент растяжения является параметром $r$-алгоритма, он выбирается эвристически для конкретной задачи.} становятся направлениями
убывания целевой функции.

Имеется также связь и сходство между методом эллипсоидов
и методами {\it переменной метрики} (см. \cite[\S $\,$ 3 гл.~3]{Polyak}).

Существует ещё несколько методов, которые относятся
к семейству методов отсекающей гиперплоскости. В 2015~г. был предложен метод, который с точностью до логарифмического по $n$ множителя работает по оценке~\eqref{ch2_eq_optim_low_bound} в
смысле требуемого числа итераций\footnote{Но с очень большим числовым множителем $10^{24}$.}
с общими затратами на \ag{всех} итерациях порядка
$O \left ( n^3 \log^{O(1)}(n/\varepsilon) \right )$~\cite{LeeSidford_15}. Отметим также {\it метод вписанных эллипсоидов} Л.\,Г.~Хачияна, С.\,П.~Тарасова, И.\,И.~Эрлиха, предложенный в 1988~г.~\cite[стр. 253--259]{Hach_izbr}, проигрывающий методу 2015~года
лишь в оценке стоимости итерации.

Нельзя не упомянуть, что
на основе метода эллипсоидов
Л.\,Г.~Хачияном в 1978~г. был получен и обоснован первый полиномиальный алгоритм решения
задачи линейного программирования в битовой сложности~\cite{Hach}. Так изящно был решён долгое время
\nk{остававшийся} открытым вопрос, является ли задача ЛП NP-сложной или
нет\footnote{Приведём цитату из воспоминаний В.\,А.~Гурвича о предшествовавшем созданию алгоритма Хачияна
разговоре между А.\,С.~Немировским и самим Хачияном:
<<В конце 1977-го (а, может, уже в 1978 году) Аркадий Немировский
рассказывал об информационной сложности задач выпуклого программирования. Оценивалось число итераций метода эллипсоидов.
После доклада Леонид спросил: <<Что из этого следует для ЛП?>> Немировский ответил кратко: <<Ничего>>. Однако меньше чем за месяц Леонид
доказал, что метод эллипсоидов решает ЛП за полиномиальное время>> (цит. по~\cite[с. 468]{Hach_izbr}).
И см. там же поясняющее комментарий Аркадия Семёновича
примечание С.\,П.~Тарасова:
<<В рамках общей теории информационной сложности задач математического и, в частности, выпуклого программирования, которую разрабатывали Немировский и Юдин, метод эллипсоидов был важным, но рядовым результатом. Более
того, применять метод эллипсоидов к задаче~ЛП было вообще неинтересно, поскольку необходимая информация для её решения восстанавливается мгновенно (а
именно: это центральный вопрос теории). Ровно это Немировский и имел в виду,
отвечая Леониду: <<Ничего>>. Битовая же сложность вычислений не является вполне естественной для <<непрерывных>> задач математического программирования и не
была особенно популярна у специалистов>>.}. 

В следующем параграфе приведено описание 
полиномиального алгоритма Леонида Хачияна.

\subsubsection{Полиномиальный алгоритм для ЛП}
В этом параграфе приведено описание 
алгоритма Леонида Хачияна, записанное по его лекциям 1996~г. в Ратгерском университете (США)  Х.~Эльбассиони
и напечатанное в~\cite[стр. 453--461]{Hach_izbr}\footnote{
Приведём цитату из воспоминаний Х.~Эльбассиони о  Хачияне:
<<Леонид владел магией изложения сложных, на первый взгляд, вещей так, что они становились простыми и естественными. Совершенно очевидно, что его собственный <<мир ЛП>>
сильно отличался от стандартного, известного остальным>>
(цит. по~\cite[с. 454]{Hach_izbr}).}.

\paragraph{Битовая модель.}
Рассмотрим общую задачу ЛП 
\begin{eqnarray}
\label{ch2_eq_lp_prob_xach}
& \delta^* = \max c^T x, &  \\
& Ax \le b, \quad x \, \in \, \mathbb{R}^n, & \nonumber
\end{eqnarray}
где $A=(a_{ij})_{i, \, j} \, \in \, \mathbb{Z}^{m \times n}$~--- целочисленная $(m \times n)$-матрица, $b = (b_i)_i \, \in \, \mathbb{Z}^m$ и $c = (c_j)_j
\, \in \, \mathbb{Z}^n$~--- целочисленные векторы. 
Обозначим $a_i^{T}$ $i$-ю строку матрицы~$A$.

В {\it равномерной {\rm(}unit-cost}) модели вычислений считается, что умножение двух целых чисел, $x$ и $y$,
занимает единицу времени: $\text{cost}(x \cdot y) = 1$.
Но в {\it битовой модели} считается, что $\text{cost}(x \cdot y) = \ell(x) + \ell(y)$,
где
$\ell(x) = \log_2(1 + |x|)$, т.е. $\ell(x)$ равно числу символов в двоичном представлении $x$.

Например, рассмотрим следующую программу для вычисления $x = a^{2^k}$ (см. алгоритм~\ref{ch2_alg_a2k}).
В равномерной модели эта программа требует $k$~тактов,
а в битовой модели $\log_2 x = 2^k \ell(a)$.

\begin{algorithm}[htb]
\caption{Вычисление $a^{2^k}$}\label{ch2_alg_a2k}
		\begin{algorithmic}[1]
			\REQUIRE Числа $a$, $k$
			\ENSURE Число $x = a^{2^k}$
\STATE $x \leftarrow a$			
\FOR{$i=1, \, \dots, \, k$}
\STATE $x \leftarrow x^2$			
\ENDFOR
\RETURN $x$
\end{algorithmic}
\end{algorithm}

Пусть 
$\ell = \max \{ \ell(a_{ij}), \, \ell({b_i}), \, \ell({c_j}) \}$, 
$L$~--- полное число битов, необходимых для представления входной информации, т.е.
$$
L = \sum_{i, \, j} \ell(a_{ij}) + \sum_{i} \ell(b_i) + \sum_{j} \ell (c_j).
$$
Говорят, что алгоритм является полиномиальным в битовой модели, если он требует не более чем полиномиального числа $\text{poly}(n, \, m, \, \ell)$ операций ($\{+, \, -, \, *, \, / \}$) и при этом длина двоичной записи используемых в вычислениях чисел не превышает
$\text{poly}(n, \, m, \, \ell)$. 

\paragraph{Метод эллипсоидов.}
Фиксируем $\varepsilon > 0$. Назовём $\tilde{x}$ {\it $\varepsilon$-приближённым} решением~\eqref{ch2_eq_lp_prob_xach}, если $c^T \tilde{x} \ge \delta^* - \varepsilon$ и $a_i^T \tilde{x} \le b_i + \varepsilon$, $i = 1, \, \ldots, \, m$. Соответственно,
назовём $\hat{x}$ {\it $\varepsilon$-допустимым} решением~\eqref{ch2_eq_lp_prob_xach}, если
выполнены ограничения
$a_i^T \hat{x} \le b_i + \varepsilon$, $i = 1, \, \ldots, \, m$.

{\bf Условие 1.} Пусть известен радиус $R \, \in \, \mathbb{R}_{+}$, такой, что \eqref{ch2_eq_lp_prob_xach}
имеет оптимальное решение $x^*$ в шаре $B_2(0, \, R)$.

Пусть $h = \max \{ |a_{ij}|, \, |b_i|, \, |c_j| \}$, т.е.
$h = 2^\ell - 1$. Тогда при выполнении условия~1 можно
с помощью метода эллипсоидов найти $\varepsilon$-приближённое решение~\eqref{ch2_eq_lp_prob_xach} за $O((n+m)n^3 \log(Rhn/\varepsilon))$ арифметических операций над 
$O(\log(Rhn/\varepsilon))$-разрядными двоичными числами.

\begin{teo}[об эллипсоиде и гиперплоскости \cite{Hach_izbr}]
\label{ch2_teo_hach_elbas_fact2}
Пусть $\mathcal{E}$~--- эллипсоид в $\mathbb{R}^n$ с центром в $\eta$ и $a \, \in \, \mathbb{R}^n$~--- ненулевой вектор. Гиперплоскость $\pi = \{ x \, | \, a^T (x - \eta) = 0 \}$ разбивает $\mathcal{E}$ на две части, $\mathcal{E}^+$ и $\mathcal{E}^{-}$. Тогда существует эллипсоид $\mathcal{E}'$ такой, что
$$
\mathcal{E}^- \subseteq \mathcal{E}' \quad \mbox{и} \quad 
\frac{\textup{vol}(\mathcal{E}')}{\textup{vol}(\mathcal{E})} \le e^{-\frac{1}{2(n+1)}} \approx 1 - \frac{1}{2n}.
$$
\end{teo}

Доказательство см. в~\cite{Hach_izbr}, сравните также с~\eqref{ch2_eq_ellip_vol_ek1k}.

Пусть $X_\varepsilon^*= \left \{ x \, | \, c^Tx \ge \delta^* - \varepsilon, \, a_i^T x \le b_i + \varepsilon \mbox{ для } i = 1, \, \ldots, \, m, \, \|x\| \le R \right \}$.
Предполагаем, что точное решение~\eqref{ch2_eq_lp_prob_xach} $x^* \, \in \, B_2(0, R)$ существует.
Тогда из условия~1 следует, что $X_\varepsilon^* \neq \varnothing$, т.к. $x^* \, \in \, X_\varepsilon^*$.

{\bf Факт 1} \cite{Hach_izbr}.
$$
\frac{\text{vol}(X_\varepsilon^*)}{\text{vol}(B_2(0, R))} \ge \left ( \frac{\varepsilon}{h \sqrt{n} R} \right )^n.
$$

\begin{teo}[об эллипсоиде меньшего объёма \cite{Hach_izbr}]
\label{ch2_teo_hach_elbas_fact4}
Пусть дан эллипсоид $\mathcal{E}$ с центром в $\eta$ такой, что $X_\varepsilon^* \subseteq \mathcal{E}$, но
$\eta \, \notin \, X_\varepsilon^*$. Тогда, используя $O((n+m)n)$ операций, можно вычислить такой новый эллипсоид $\mathcal{E}'$, что 
$$
X_\varepsilon^* \subseteq \mathcal{E}'
$$
и
$$
\textup{vol}(\mathcal{E}') \le e^{-\frac{1}{2(n+1)}}
\textup{vol}(\mathcal{E}).
$$
\end{teo}

\begin{proof}
Сначала проверим, является ли $\eta$ $\varepsilon$-допустимым решением (в таком случае должно выполняться
неравенство $a_i^T \eta \le b_i + \varepsilon$
для всех $i \, \in \, \{1, \, \ldots, \, m \}=M$).

\begin{enumerate}
\item Если $\eta$~--- не $\varepsilon$-допустимое решение,
то найдётся такой индекс $i_* \, \in \, M$, 
что $a_{i_*}^T \eta > b_i + \varepsilon$.
Тогда, если бы для какого-то $x$
выполнялось условие $a_{i_*}^T x > a_{i_*}^T \eta$,
то в этом случае
$x \, \notin \, X_\varepsilon^*$.
Следовательно, для любого $x \, \in \, X_\varepsilon^*$
должно выполняться
$a_{i_*}^T x \le a_{i_*}^T \eta$, т.е. 
$a_{i_*}^T (x - \eta) \le 0$.
Значит, в этом случае можно положить $a = a_{i_*}$ в качестве нормали отсекающей гиперплоскости из теоремы~\ref{ch2_teo_hach_elbas_fact2}.

\item Пусть теперь центр $\eta$~---  $\varepsilon$-допустимое решение.
\begin{enumerate}
    \item Если $\| \eta \| > R$, то $\eta^T (x - \eta) \le 0$ для всех $x \, \in \, X_\varepsilon^*$, и
    в этом случае мы полагаем $a = \eta$.
    
    \item Иначе, пусть $\| \eta \| \le R$. Поскольку
    $\eta \, \notin \, X_\varepsilon^*$, то, по определению $X_\varepsilon^*$ $c \eta \, < \, \delta^* - \varepsilon$. Следовательно, любой $x$ такой, что $c^Tx \, < \, c^T \eta$ не принадлежит $X_\varepsilon^*$.
    Тогда для любого $\varepsilon$-приближённого решения
    $c^Tx \ge c^T \eta$, т.е. $-c^T(x- \eta) \le 0$. И можно
    в качестве нормали~$a$ выбрать $-c$ и использовать
    теорему~\ref{ch2_teo_hach_elbas_fact2}. $\qedhere$
\end{enumerate}
\end{enumerate}
\end{proof}

Метод эллипсоидов стартует с $\mathcal{E}_0 = B_2(0, \, R)$ и порождает последовательность эллипсоидов
$\mathcal{E}_0 = B_2(0, \, R)$, $\mathcal{E}_1$,
$\mathcal{E}_2, \, \ldots, \, \mathcal{E}_N, \, \mathcal{E}_{N+1}$. Пока $\eta_k$ (центр эллипсоида $\mathcal{E}_k$) не попал в $X_\varepsilon^*$,
можно построить новый эллипсоид $\mathcal{E}_{k+1}$
такой, что $X_\varepsilon^* \subseteq \mathcal{E}_{k+1}$
и $\frac{\textup{vol}(\mathcal{E}_{k+1})}{\textup{vol}(\mathcal{E}_{k})} \ge \exp\left(-\frac{1}{2(n+1)}\right)$.
В частности,
$$
\textup{vol}(X_\varepsilon^*) \le \textup{vol}(\mathcal{E}_{N}) \le \exp\left(-\frac{N}{2(n+1)}\right) \cdot
\textup{vol}(B_2(0, \, R)).
$$
Из факта~1 следует
$$
\textup{vol}\left ( B_2(0,  R) \right ) \left ( \frac{\varepsilon}{h \sqrt{n} R} \right )^{\!n}\!\le\!\textup{vol}(X_\varepsilon^*)\!\le\!\textup{vol}(\mathcal{E}_{N})\!\le\!\exp\left(-\frac{N}{2(n+1)}\right) \cdot
\textup{vol}(B_2(0,  R)),
$$
откуда получается верхняя оценка числа итераций
$$
N \le 2 n (n+1) \ln \frac{R \sqrt{n} h}{\varepsilon} = O \left ( n^2 \log \frac{Rnh}{\varepsilon} \right ).
$$

Все вычисления требуют $$O \left ((n+m) n^3 \log (Rnh/\varepsilon) \right ) = O \left (m n^3 \log (Rnh/\varepsilon) \right )$$ арифметических операций\footnote{Если
$m \, < \, n$, то можно за время $O(m)$ найти
базис и уменьшить размерность задачи.}
над $O(\log (Rnh/\varepsilon))$-разрядными двоичными числами.

\paragraph{Как восстановить точное решение по $\varepsilon$-приближённому.}
Вернёмся к нахождению точного решения задачи ЛП~\eqref{ch2_eq_lp_prob_xach} с целочисленными $A$ и $b$.

Для целочисленной матрицы $A$ обозначим
$$
\Delta (A) = \max \, \left \{ |\det B | \, | \, B \mbox{ --- квадратная подматрица } A \right \}.
$$
\begin{lemma}[\cite{Hach_izbr}] 
Если задача ЛП~\eqref{ch2_eq_lp_prob_xach} с целочисленными $A$ и $b$ разрешима, то у неё есть такое оптимальное решение $x^*$, что 
$$
\|x^* \|_\infty \le n \| b \|_\infty \Delta(A),
$$
где $\| a \|_\infty = \max\limits_i |a_i|$~--- так называемая \textup{ бесконечная (или максимальная) норма} вектора~$a$.
\end{lemma}

Таким образом, если $R = h n^{3/2} \Delta(A)$, то
гарантированно $\|x^*\|_2 \le R$ (поскольку $\| \cdot \|_2 \le n^{1/2} \| \cdot \|_\infty$).

\begin{lemma}[\cite{Hach_izbr}]
\label{ch2_lem_hach_elbas_lem4}
Рассмотрим систему
\begin{equation}
\label{ch2_eq_hach_elbas_lem4}
    Ax \le b, \quad x \, \in \, \mathbb{R}^n.
\end{equation}
Пусть $\varepsilon_0 = \frac{1}{(n+2)\Delta(A)}$. Если~\eqref{ch2_eq_hach_elbas_lem4} имеет $\varepsilon_0$-допустимое решение $a_i^T \tilde{x} \le\linebreak\le b_i + \varepsilon_0$, $i \, \in \, M = \{1, \, 2, \, \ldots, \, m\}$, то система~\eqref{ch2_eq_hach_elbas_lem4}
совместна.
\end{lemma}

Покажем, следуя~\cite{Hach_izbr}, как, зная $\tilde{x}$~--- $\varepsilon_0$-допустимое решение ~\eqref{ch2_eq_hach_elbas_lem4}, вычислить точное решение
данной системы. Причём необходимые затраты составят $m n^2$ операций. 

Если $\tilde{x}$~--- $\varepsilon$-допустимое решение~\eqref{ch2_eq_hach_elbas_lem4}, то
$a_i^T \tilde{x} \le b_i + \varepsilon_0$, $i \, \in \, M$. Пусть $I_0 = \left \{  i \, \in \, M \, | \, |a_i^T \tilde{x} - b_i | \le \varepsilon_0 \right \}$.
Тогда для любого $j \, \in (M - I_0)$ выполнено $a_j^T \tilde{x} \, < \, - \varepsilon_0 +  b_j$.

Подсистема 
\begin{equation}
\label{ch2_eq_hach_elbas_axebi0}
a_i^T x = b_i, \quad i \, \in \, I_0,
\end{equation}
совместна по лемме~\ref{ch2_lem_hach_elbas_lem4},
так как её можно записать в виде системы линейных неравенств
$$
a_i^T x \le b_i, \, - a_i^T x \le - b_i, \quad i \, \in \, I_0, 
$$
и $\tilde{x}$ является её $\varepsilon_0$-допустимым решением\footnote{$\Delta 
\left [ 
\begin{array}{c}
A_{I_0} \\
-A_{I_0}
\end{array}
\right ] \le \Delta \left (A_{I_0} \right ) \le \Delta \left ( A \right ) 
$,
где $A_{I_0}$~--- это подматрица матрицы~$A$, строки которой пронумерованы $I_0$.
}. 

Пусть $\hat{x}$~--- какое-нибудь точное решение системы~\eqref{ch2_eq_hach_elbas_axebi0}. 
Если $\tilde{x}$ не удовлетворяет исходной системе, то
определим
$$
x(t) = (1 - t) \tilde{x} + t \hat{x}, \quad 0 \le t \le 1.
$$
Тогда для $i \, \in \, I_0$ получаем
\begin{align*}
\left| a_i^T x(t) - b_i \right| &= \left| (1 - t) a_i^T \tilde{x} + t a_i^T \hat{x} - (1-t) b_i - t b_i \right| = \\
&= \left| (1-t) \left[a_i^T \tilde{x} - b_i\right] + t \left[a_i^T \hat{x} - b_i\right] \right| \\
&\le (1 - t) \left| a_i^T \tilde{x} - b_i \right| +
t \left|a_i^T \hat{x} - b_i \right| \le (1 - t) \varepsilon_0 + t \cdot 0 \le \varepsilon_0,
\end{align*}
т.е. $x(t)$ удовлетворяет ограничениям множества~$I_0$ для всех 
$t \, \in \, [0, \, 1]$.

Теперь можно получить новую точку, в которой множество
$\varepsilon_0$-до\-пус\-ти\-мых ограничений будет строго больше, чем $I_0$.
Для этого выберем $t$ так:
$$
t = \min_{j \, \in \, (M - I_0)} \left \{ \left. - \frac{a_j^T \tilde{x} - b_j}{a_j^T \hat{x} - a_j^T \tilde{x}}
\, \right |  \, a_j^T \hat{x} - b_j > 0
\right \}.
$$
Пусть $j_{\min}$~--- индекс, на котором минимум достигается. Тогда 
$$
a_{j_{\min}}^T x(t) - {b_{j_{\min}}} = (1 - t) \left (a_{j_{\min}}^T \tilde{x} - {b_{j_{\min}}} \right ) + t
\left (a_{j_{\min}}^T \hat{x} - {b_{j_{\min}}} \right ) = 0.
$$
Добавим $j_{\min}$ к $I_0$, заменим $\tilde{x}$ на $x(t)$ и повторим процедуру. Таким способом будет проведено не более $n$ итераций, причём каждая итерация (решение системы линейных уравнений) требует не более $O(m n^2)$ операций.  Но и общее время работы алгоритма
равно $O(m n^2)$.

Такой же алгоритм восстановления можно применить и к задаче ЛП. Для этого вместо $\varepsilon_0$-допустимого решения нужно искать $\varepsilon_1 = \frac{1}{4 n^{5/2} \Delta^3(A) \| c \|}$-приближённое решение.

\paragraph{Резюме.} При $\varepsilon_1 = \frac{1}{4 n^{5/2} \Delta^3(A) \| c \|}$ и $R = n^{3/2} h \Delta(A)$ методу эллипсоидов требуется выполнить порядка $O(m n^3 \log(Rhn/\varepsilon_1)) = O \left ( m n^3  \log \left ( hn \Delta(A) \right ) \right )$ операций.

Пусть $L = L(A, \, b, \, c)$~--- битовая длина входа задачи ЛП~\eqref{ch2_eq_lp_prob_xach}, и пусть
$$
L(A) = \sum_{i, \, j} \log_2 \left ( 1 + |a_{ij}| \right ).
$$
Тогда $\Delta(A) \le 2^{L(A)}$, так как
$\Delta(A) \le \prod_{i, \, j} \left ( 1 + |a_{ij}| \right )$.
Кроме того, $n \le 2^{L(A)}$, $h \le 2^{L(A, \, b, \, c)}$. Таким образом, для решения задачи ЛП методу эллипсоидов требуется $O(m n^3 L)$ операций. Причём для
вычислений достаточно точности $O(n \Delta) = O(L)$ двоичных разрядов.

\section[Субградиентные методы]{Субградиентные методы\\ для задач выпуклой оптимизации\sectionmark{Субградиентные методы}}
\sectionmark{Субградиентные методы}\label{subgradient}

Опишем оптимальный метод на классе выпуклых липшицевых (возможно, негладких) целевых функций в самом простом случае: $Q=\mathbb{R}^n$ и $\left\|~\cdot~ \right\|=\left\|~\cdot~ \right\|_2$. На классе задач минимизации выпуклых негладких и липшицевых функций оптимальным будет обычный субградиентный метод, впервые предложенный \ag{Н.\,З.~Шором~\cite{Shor_book}}. 
Пусть $B_2^n(x^*, R)=\left\{ {x\in \mathbb{R}^n:\;\;\left\|{x-x^*} \right\|_2 \leqslant R} \right\}.$ 
Будем всюду далее обозначать субградиент $f$ (некоторый элемент субдифференциала $\partial f\left( x \right)$) в точке $x$ как $\nabla f\left( x \right)$. Если функция $f$ дифференцируема в точке $x$, то $\nabla f(x)$~--- её градиент. Итерационный процесс субградиентного метода имеет вид
\begin{equation}\label{p0_eq2}
x_{k+1}=\;x_k-h\nabla f(x_k).
\end{equation}
Предположим, что
\begin{equation}\label{p0_eq3}
\left\|\nabla f(x) \right\|_2 \leqslant M \text{ для всякого } x\in B_2^n(x^*, R\sqrt{2}).
\end{equation}

Если допустить (базируясь на условии~\eqref{p0_eq3}, которое мы впоследствии проверим), что при всяком $k \geqslant 0$ верно $\left\|\nabla f(x_k) \right\|_2 \leqslant M$, то из \eqref{p0_eq2} следует
$$\left\| {x-x_{k+1}} \right\|_2^2 \mbox{=}\left\|x-x_k+h\nabla f\left({x_k} \right) \right\|_2^2 =$$
$$=\;\left\| {x-x_k} \right\|_2^2 +2h\left\langle {\nabla f\left( {x_k}\right),x-x_k} \right\rangle +h^2\left\| {\nabla f(x_k)}\right\|_2^2 \leqslant$$
\begin{equation}
\label{p0_KI}
\leqslant \;\left\| {x-x_k} \right\|_2^2 +2h\left\langle {\nabla f(x_k),x-x_k} \right\rangle + h^2M^2.
\end{equation}

Выберем теперь $x=x^*$ (если $x^*$ не единственно, то выбирается ближайшее решение $x^*$ к точке старта $x_0)$. Имеем:
$$f\left( {\frac{1}{N}\sum\limits_{k=0}^{N-1} {x_k} } \right)-f(x^*) \leqslant$$
$$\leqslant\frac{1}{N}\sum\limits_{k=0}^{N-1} {f\left( {x_k} \right)} -f\left( {x^*}\right)\leqslant\frac{1}{N}\sum\limits_{k=0}^{N-1} {\left\langle {\nabla f\left( {x_k} \right),x_k-x^* } \right\rangle } \leqslant$$
$$\leqslant \frac{1}{2hN}\sum\limits_{k=0}^{N-1} {\left( {\left\| {x^* - x_k}\right\|_2^2 -\left\| {x^* - x_{k+1}} \right\|_2^2 } \right)}+\frac{hM^2}{2}=$$
$$=\frac{1}{2hN}\left( {\left\| {x^* - x_0} \right\|_2^2 -\left\| {x^*-x_N} \right\|_2^2 } \right)+\frac{hM^2}{2}.$$

Если выбрать шаг субградиентного метода и точку выхода алгоритма следующим образом (здесь $R=\left\| x_0-x^*  \right\|_2 $):
\begin{equation}
\label{p0_eq4}
h=\frac{R}{M\sqrt N },
\quad
\bar {x}_N=\frac{1}{N}\sum\limits_{k=0}^{N-1} {x_k},
\end{equation}
то будет выполняться неравенство
\begin{equation}
\label{p0_eq5}
f\left( {\bar {x}_N} \right)-f(x^*) \leqslant \frac{MR}{\sqrt N}.
\end{equation}
Заметим, что точная нижняя оценка на классе задач выпуклой оптимизации с условием (\ref{p0_eq3}) для методов первого порядка вида
$$x_{k+1} \in x_0 + \text{span}\left\{\nabla f(x_0),...,\nabla f(x_k)\right\},$$
где для всех $j = 0, 1, ..., k$ верно $\nabla f(x_j) \in \partial f(x_j)$,
имеет вид~\cite{DroriTeboulle16}:
$$
f\left({x_N} \right)-f(x^*) \geqslant \frac{MR}{\sqrt {N+1} }.$$
Неравенство (\ref{p0_eq5}) означает, что при выборе
$$N=\frac{M^2R^2}{\varepsilon ^2},\quad h=\frac{\varepsilon }{M^2}$$
будет достигаться оценка $f\left( {\bar{x}_N} \right)-f(x^*) \leqslant \varepsilon$, которая оптимальна с точностью до умножения на константу. Далее ввиду (\ref{p0_eq2}) при выборе в (суб)градиентном методе переменной величины шага
\begin{equation}\label{p0_eq6}
x_{k+1}\mbox{=}\;x_k-h_k \nabla f(x_k),\quad h_k =\;\frac{\varepsilon }{\left\| {\nabla f\left( {x_k} \right)}\right\|_2^2 }
\end{equation}
оценка (\ref{p0_eq5}) сохранится для точки выхода вида 
$$\bar{x}_N=\frac{1}{\sum\limits_{k=0}^{N-1} {h_k } }\sum\limits_{k=0}^{N-1}{h_k x_k}.$$
Если же в~(\ref{p0_eq6}) выбрать шаги $h_k$ по аналогии с (\ref{p0_eq4}):
$$h_k =\frac{R}{\left\|\nabla f(x_k) \right\|_2 \sqrt N },$$
то будет верна следующая оценка, похожая на (\ref{p0_eq5}):
$$\mathop{ \min\limits_{k=0,...,N-1} f\left( {x_k} \right)-f(x^*) \leqslant\frac{MR}{\sqrt N }}.$$

Заметим, что при этом для всякого $k \geqslant 0$ верны неравенства
$$
0\leq\frac{1}{2kh_k}\biggl(\left\|x^{*}-x_{0}\right\|^{2}_{2}-\left\|x^{*}-x_{k}\right\|^{2}_{2}\biggr)+\frac{h_kM^{2}}{2}.
$$
Поэтому при любом $k=0,\ldots,N$
$$
\left\|x^{*}-x_{k}\right\|^{2}_{2}\leq\left\|x^{*}-x_{0}\right\|^{2}_{2}+h_k^{2}M^{2}k\leq 2\left\|x^{*}-x_{0}\right\|^{2}_{2},
$$
т.е.
$$
\left\|x_{k}-x^{*}\right\|_{2}\leq \sqrt{2}\left\|x_{0}-x^{*}\right\|_{2},\quad k=0,\ldots,N.
$$
Это означает, что в случае предлагаемых способов подбора шагов и в силу \eqref{p0_eq3} можно считать
$\left\|\nabla f(x_k) \right\|_2 \leqslant M$ для произвольного $ k=0, \, \ldots, \, N$. То~есть предположение, которые было сделано в самом начале раздела, подтвердилось.

\begin{remark}
Возможно пояснение изложенного подхода с точки зрения теории размерностей физических величин (подробнее см.~\cite{Zorich_mathan}, гл.~I). Напомним сначала \textit{$\Pi$-теорему теории размерностей} (\cite{Zorich_mathan}, параграф~3.1) на простом примере колебания маятника. 

Пусть к потолку на нити длиной $l$ [м] подвешен груз массой $m$ [кг]. Груз находится в поле силы тяжести $g$ [м/с$^2$]. Нужно определить период колебания груза $T$ [с], считая, что он может зависеть только от перечисленных выше величин. Идея решения физических задач из соображений теории размерностей заключается в  том, чтобы составить независимые (функционально не выражающиеся через друг друга) дробно рациональные выражения (вида произведения доступных переменных, возведенных в, вообще говоря, различные рациональные, в том числе и отрицательные, степени), которые были бы  безразмерными. Утверждение $\Pi$-теоремы заключается в том, что существует такая функция, число аргументов которой равно числу независимых безразмерных комбинаций, которая тождественно равна константе при подстановке в качестве ее аргументов отмеченных безразмерных комбинаций параметров (переменных). В случае примера с маятником существует всего одна такая безразмерная комбинация $Tl^{-1/2}g^{1/2}$. Отметим, что информация о $m$ оказалась ненужной. Из~$\Pi$-теоремы следует, что существует такая функция $\Pi$, что $\Pi\left(Tl^{-1/2}g^{1/2}\right)\equiv \text{const}$, то есть $T = \Pi^{-1}\left(\text{const}\right)\sqrt{l/g}$. При этом такая константа $\Pi^{-1}\left(\text{const}\right)$ универсальна (одинакова) для всех маятников. То есть можно один раз провести эксперимент (в кораблестроении такие эксперименты называют <<экспериментами в ванной>> \cite{Sedov}), определив эту константу (в нашем случае это будет $\Pi^{-1}\left(\text{const}\right) = 2\pi$), и получить, таким образом, готовую формулу для периода колебания математического маятника, практически не привлекая никакие физические соображения, кроме того, что в самом начале бы (с запасом) обрисовали набор переменных (параметров) задачи, от которых период колебания может зависеть.

Можно ввести размерности и в численные методы оптимизации (и, на самом деле, часто бывает весьма полезно это делать для контроля правильности получаемых формул). А именно, предположим, что целевая функция измеряется в рублях [руб]. То есть $\varepsilon$ [руб], а аргумент функции рассчитывается в килограммах [кг], т.е. $x$ [кг]. Тогда $R$ [кг], $M$ [руб/кг], количество итераций $N$~--- безразмерная величина. Шаг субградиентного метода $h$ [кг$^2$/руб]~--- это следует из \eqref{p0_eq2}.

Считая, что $N$ может зависеть только от $R$, $M$,  $\varepsilon$, из $\Pi$-теоремы сразу можно заключить, что существуют такие числа $\text{const}$ и $\alpha$, что $N = \text{const}\cdot\left(MR/\varepsilon\right)^{\alpha}$. Это, конечно, не готовый ответ, но хорошо согласуется с полученным выше результатом $N = M^2R^2/\varepsilon^2$.

Аналогично можно поступить и с зависимостью шага метода $h$ от параметров задачи $R$, $M$, $\varepsilon$. Возможные варианты $h = \text{const}_1\cdot\varepsilon/M^2$ и $h = \text{const}_2\cdot R/M$. Собственно, выше у нас уже получалось, что $h = \varepsilon/M^2$ и $h = R/\left(M\sqrt{N}\right)$.  \ag{Первая} формула получается из второй после подстановки $N = M^2R^2/\varepsilon^2$.
\end{remark}

Похожие методы можно применять и к задачам с функциональными ограничениями вида
\begin{equation}\label{equv_state_problem}
f(x) \rightarrow \min_{x \in Q}, \quad g(x) \leqslant 0, 
\end{equation}
где $f,g:~~\mathbb{R}^n \to \mathbb{R}$. \ag{Субградиентные методы для задач с ограничениями-неравенствами впервые были исследованы Б.\,Т.~Поляком в \cite{Polyak_subgrad}.}

Здесь уже $Q$ --- замкнутое подмножество $\mathbb{R}^n$, а $Pr_Q$ --- оператор проектирования на множество $Q$:
$$
Pr_Q(y) = \argmin_{x\in Q} \|y-x\|_2.
$$

Точнее говоря, к задачам типа \eqref{equv_state_problem} возможно применять субградиентные схемы с переключениями по продуктивным и непродуктивным шагам следующего вида. 

\floatname{algorithm}{Алгоритм}
	\begin{algorithm}
		\caption{Субградиентный метод с переключениями для задач выпуклого программирования \eqref{equv_state_problem}}\label{Switch_Subgrad}
\begin{algorithmic}[1]
\REQUIRE $x_0$~--- начальная точка, $h_f > 0$, $h_g > 0$, количество шагов $N$.
\STATE {\bf If} $g(x_k)\leqslant\varepsilon$
\STATE
{\bf then}
$k \rightarrow I$ (<<продуктивный шаг>>)
$$
x_{k+1}=Pr_Q\{x_k-h_f\nabla f(x_k)\},
$$
\STATE  {\bf else}
$k \rightarrow J$ (<<непродуктивный шаг>>)
$$
x_{k+1}=Pr_Q\{x_k-h_g\nabla g(x_k)\},
$$
\RETURN $\widehat{x}= \frac{1}{|I|}\sum_{k \in I} x_k.$
\end{algorithmic}
\end{algorithm}

Обозначим множество продуктивных шагов (для которых $g(x_k)\leqslant\varepsilon$) через $I$, а множество других (непродуктивных) шагов через $J$. Предположим, что целевой функционал $f$ удовлетворяет условию Липшица с константой $M_f > 0$ на множестве $Q$, а функционал ограничения $g$ удовлетворяет условию Липшица с константой $M_g > 0$ на этом же множестве. При этом на продуктивных шагах будем выбирать субградиенты целевого функционала такие, что $\|\nabla f(x_k)\|_2 \leqslant M_f$, а на непродуктивных --- субградиенты функционала ограничения такие, что $\|\nabla g(x_k)\|_2 \leqslant M_g$. 
\begin{teo} \label{swithc_subgrad_thm} Если выбрать 
$$h_f = \frac{\varepsilon}{M_f^2} \quad \text{и} \quad h_g = \frac{\varepsilon}{M_g^2},$$
то после 
$$
N = \left\lceil\frac{\|x_0 - x^*\|_2^2\max\{M_f^2,\,M_g^2\}}{\varepsilon^2}\right\rceil
$$
шагов указанного выше алгоритма {\rm \ref{Switch_Subgrad}} гарантированно $I \neq \varnothing$ и будут верны неравенства
$$
f(\widehat{x}) - f(x^*) \leqslant \varepsilon, \quad g(\hat{x}) \leqslant \varepsilon,
$$
где $\widehat{x} = \frac{1}{|I|} \sum\limits_{k \in I} x_k$. \end{teo}
\begin{proof}
Рассмотрим отдельно необходимые неравенства для продуктивных и непродуктивных шагов.

1) Если $k$-й шаг алгоритма продуктивен ($k\in I$), то $h_k = h_f$ и верны неравенства
\begin{equation}\label{p5_eq2a}
\begin{split}
h_k(f(x_k)-f(x^*))\leqslant h_k\langle\nabla f(x_k), x_k-x^*\rangle\leqslant\\
\leqslant\frac{h_k^2}{2}||\nabla f(x_k)||_2^2+\frac{1}{2}\|x^* - x_k\|_2^2 - \frac{1}{2}\|x^* - x_{k+1}\|_2^2\leqslant\\
\leqslant\frac{\varepsilon^2}{2M_f^2}+ \frac{1}{2}\|x^*- x_k\|_2^2 - \frac{1}{2}\|x^*-x_{k+1}\|_2^2.
\end{split}
\end{equation}
Первое неравенство следует из выпуклости $f$, второе получается из свойств  оператора проектирования: 
$$\frac{1}{2}\|x^* - x_{k+1}\|_2^2 =  \frac{1}{2}\left\|x^* - Pr_Q\{x_k-h_f\nabla f(x_k)\}\right\|_2^2 =
$$
$$=\frac{1}{2}\left\|Pr_Q\{x^* - (x_k-h_f\nabla f(x_k))\}\right\|_2^2\leqslant \frac{1}{2}\left\|x^* - (x_k-h_f\nabla f(x_k))\right\|_2^2. 
$$

2) Если $k$-й шаг алгоритма не продуктивен ($k \, \not\in \, I$), то $h_k = h_g$. Поэтому $g(x_k)>\varepsilon$ и $g(x_k)-g(x^*) \geqslant g(x_k)>\varepsilon$. Это означает, что верны неравенства
\begin{equation}\label{p5_eq3}
\begin{split}
\frac{\varepsilon^2}{M_g^2} < h_k (g(x_k)-g(x^*))\leqslant\frac{h_k^2}{2}||\nabla g(x_k)||_2^2+ \frac{1}{2}\|x^*- x_k\|_2^2-\\-\frac{1}{2}\|x^*-x_{k+1}\|_2^2 \leqslant \frac{\varepsilon^2}{2M_g^2}+\frac{1}{2}\|x^*- x_k\|_2^2- \frac{1}{2}\|x^*-x_{k+1}\|_2^2,\text{ или}\\
\frac{\varepsilon^2}{2M_g^2} < \frac{1}{2}\|x^*- x_k\|_2^2- \frac{1}{2}\|x^*-x_{k+1}\|_2^2.
\end{split}
\end{equation}

3) После суммирования неравенств~\eqref{p5_eq2a} и~\eqref{p5_eq3} имеем
$$\sum_{k\in I}h_k(f(x_k)-f(x^*))\leqslant$$
$$\leqslant\sum_{k\in I}\frac{\varepsilon^2}{2M_f^2}-\frac{\varepsilon^2|J|}{2M_g^2}+\frac{1}{2}\|x^*- x_0\|_2^2-\frac{1}{2}\|x^*-x_{N+1}\|_2^2\leqslant$$
$$\leqslant\frac{\varepsilon}{2}\sum_{k\in I}h_k-\frac{\varepsilon^2|J|}{2M_g^2}+\frac{1}{2}\|x^*- x_0\|_2^2=$$
$$=\varepsilon\sum_{k\in I}h_k-\frac{|I|\varepsilon^2}{2M_f^2}-\frac{|J|\varepsilon^2}{2M_g^2}+ \frac{1}{2}\|x^*- x_0\|_2^2.$$

Это означает, что после выполнения 
$$
N = \left\lceil\frac{\|x_0 - x^*\|_2^2\max\{M_f^2,\,M_g^2\}}{\varepsilon^2} \right\rceil
$$
итераций алгоритма \ref{Switch_Subgrad} будет заведомо верно неравенство
$$\sum_{k\in I}h_k\left(f(x_k)-f(x^*)\right)\leqslant\varepsilon\sum_{k\in I}h_k,$$
откуда для выхода $\widehat{x}$ алгоритма \ref{Switch_Subgrad} имеем
$$f(\widehat{x})-f(x^*)\leqslant\varepsilon.$$
При этом для всякого $ k\in I$ верно неравенство $g(x_k)\leqslant\varepsilon$ и поэтому ввиду выпуклости $g$ имеем
$$g(\widehat{x})\leqslant\frac{1}{|I|}\sum_{k\in I}g(x_k)\leqslant\varepsilon.$$

Остаётся лишь показать, что множество продуктивных шагов $I$ не пусто. Если $I=\emptyset$, то $|J|=N$ и это означает, что $N\geqslant\frac{\|x^*- x_0\|_2^2M_g^2}{\varepsilon^2}$. С другой стороны, из~\eqref{p5_eq3} имеем
$$\frac{\varepsilon^2N}{2M_g^2}<\frac{1}{2}\|x^* - x_0\|_2^2.$$
Таким образом, получили противоречие и $I\neq\emptyset$. Теорема доказана.
\end{proof}

Весьма важна {\bf прямо-двойственность} рассматриваемого алгоритма \ref{Switch_Subgrad} для задач выпуклого программирования вида \eqref{equv_state_problem} на компактном допустимом множестве $Q$. Опишем соответствующий теоретический результат. Для этого уточним задание функционала ограничения $g(x)=\displaystyle\max_{l = 1, \, \ldots, \, m}g_{l}(x)$. Рассмотрим следующую двойственную задачу к  \eqref{equv_state_problem}.
\begin{equation}\label{f_8}
\varphi(\lambda) = \min_{x\in Q}\left\{f(x)+\sum_{l=1}^{m}\lambda_{l}g_{l}(x)\right\}\rightarrow\max_{\lambda\geqslant0}.
\end{equation}
Всегда имеет место следующее неравенство (слабая двойственность, см.~п.~\ref{ws_duality}):
$$
0\leqslant f(x)-\varphi(\lambda) := \Delta(x,\lambda), \quad x\in Q,\, g(x) \leqslant0,\, \lambda \geqslant0.
$$
Обозначим решение задачи \eqref{f_8} через $\lambda^{*}$. Сделаем предположение о том, что для задачи \eqref{equv_state_problem} выполняются условия Слейтера (см. п.~\ref{ch1_subs_slater}), т.е. существует такой $\widetilde{x}\in Q$, что $g(\widetilde{x})<0$. Тогда
$$
f(x^{*})=\varphi(\lambda^{*}) := \varphi^{*}.
$$
В этом случае <<качество>> приближённого решения (пары $(x^{N},\lambda^{N})$) вполне естественно оценивать величиной зазора двойственности $\Delta(x^{N},\lambda^{N})$. Чем этот зазор двойственности меньше, тем лучше качество достигнутого приближения искомого решения (связи решения прямой и двойственной задач). Будет полезным также такое обозначение, уточняющее выбор рассматриваемого субградиента функционала $g$: 
$$
g(x_{k})=g_{l(k)}(x_{k}),\quad \nabla g(x_{k})=\nabla g_{l(k)}(x_{k}),\, k\in J,
$$
т.е. для непродуктивной итерации ($k \in J$) под $l(k)$ мы понимаем индекс функционала ограничения, для которого достигается максимум $\max_{l = 1,..., m} g_l(x_k) = g_{l(k)}(x_k) = g(x_k)$. Теперь уточним правило построения искомых приближений двойственных множителей по работе алгоритма \ref{Switch_Subgrad}. Положим
$$
\lambda_{l}^{N}=\frac{1}{h_{f}|I|}\sum_{k\in J}h_{g}I[l(k)=l],
$$
$$
I[\text{predicat}] =
\begin{cases}
1, & \text{predicat}\, = true\\
0, & \text{predicat}\, = false.
\end{cases}
$$

По-прежнему на итерациях алгоритма \ref{Switch_Subgrad} мы используем такие субградиенты $f$ и $g$, для которых
$$
\|\nabla f(x_k)\|_{2}\leqslant M_{f},\quad \|\nabla g(x_k)\|_{2}\leqslant M_{g}.
$$
Это возможно в силу сделанных выше предположений о липшицевости функционалов $f$ и $g$.

\begin{teo}\label{thm_subgr_primal-dual}
Пусть $Q$ --- компакт, причём для некоторого $\overline{R} >0$ верно $ \| x- y\|_2^2 \leqslant  2\overline{R}^2$ для \ag{всех} $x, y \in Q$. Тогда после 
$$
N = 
\left\lceil\frac{2\max\{M^{2}_{f}, M_g^2\}\overline{R}^{2}}{\varepsilon^{2}}\right\rceil
$$
итераций алгоритма {\rm \ref{Switch_Subgrad}} заведомо верны следующие неравенства
$$
\Delta(\widehat{x},\lambda^{N})\leqslant\varepsilon,\quad g(\widehat{x})\leqslant\varepsilon.
$$
\end{teo}
\begin{proof} В доказательстве предыдущей теоремы, по сути, получены следующие неравенства:
$$
h_{f}|I|f(\widehat{x})\leqslant \min_{x\in Q}\left\{h_{f}|I|f(x)+h_{f}\sum_{k\in I}\Bigl\langle\nabla f(x_{k}), x_{k}-x\Bigr\rangle\right\} \leqslant
$$
$$
\leqslant\min_{x\in Q}\Biggr\{h_{f}|I|f(x)+\frac{h_{f}^{2}}{2}\sum_{k\in I}\left\|\nabla f(x_{k})\right\|^{2}_{2} \nk{-}
$$
$$
\nk{-}h_{g}\sum_{k\in J}\underbrace{\Bigl\langle \nabla g(x_{k}), x_{k}-x\Bigr\rangle}_{\geqslant g_{l(k)}(x_{k})-g_{l(k)}(x)} +\frac{h_{g}^{2}}{2}\sum_{k\in I}\left\|\nabla g(x_{k})\right\|^{2}_{2}  +
$$
$$
+\frac{1}{2}\sum_{k = 0}^{N-1}\left( \|x - x_{k}\|_2^2 - \| x - x_{k+1}\|_2^2\right)\Biggr\}.
$$
Далее, имеем
$$
h_{f}|I|f(\widehat{x}) \leqslant\frac{1}{2}h^{2}_{f}|I|M^{2}_{f}-\frac{1}{2M^{2}_{g}}\varepsilon^{2}|J|+\overline{R}^{2}+
$$
$$
+ h_{f}|I|\min_{x\in Q}\left\{f(x)+\sum_{l=1}^{m}\lambda^{N}_{l}g_{l}(x)\right\}=
$$
$$
= \varepsilon h_{f}|I|+\left(\overline{R}^{2}- \frac{|I|\varepsilon^{2}}{2M^{2}_{f}}-\frac{|J|\varepsilon^{2}}{2M^{2}_{g}}\right)+h_{f}|I|\varphi(\lambda^{N}) \leqslant \varepsilon h_{f}|I|+ h_{f}|I|\varphi(\lambda^{N}),
$$
откуда и следует доказываемое утверждение.
\end{proof}

Можно обобщить приведённые выше результаты на случай, когда выбирается не евклидова норма \cite{Bayandina}. Как это делается и зачем это нужно, обсуждается также далее в более общем контексте в разделе \ref{model}.

В заключении этого раздела отметим, что описанные в нём субградиентные методы при правильном выборе шага оптимальны с точки зрения нижних оценок на число вызовов оракула~\cite{NemirYd79}, \cite{Nemir_comlec95} (вычислений субградиента, для класса задач выпуклой оптимизации с ограниченной константой Липшица целевого функционала (функционального ограничения). Для класса задач с функционалом, имеющим липшицев градиент, ни при каком выборе шага в классе методов градиентного спуска вида~\eqref{p0_eq2} не удаётся найти оптимальный метод. Построению оптимальных методов в этом случае посвящены следующие разделы~\ref{gradientTaylor}, \ref{AM} и \ref{fastGradMethod}.

\section{Методы типа градиентного спуска}
\sectionmark{Методы типа градиентного спуска}\label{gradientTaylor}

Метод градиентного спуска, безусловно, является <<краеугольным камнем>> численных методов оптимизации. С некоторой натяжкой можно даже сказать, что этот метод порождает все остальные методы \cite{Gas_Pos18}. 

В данном пособии мы изложим градиентный спуск, следуя \cite{Aspremont2021}. А именно, сначала мы рассмотрим задачи квадратичной оптимизации, для которых градиентный спуск и его производные можно достаточно просто изучать. Уже на примере квадратичных задач будет подмечено, что градиентный спуск не является наилучшим возможным (оптимальным) методом для этого класса задач в классе градиентных методов (методов \textit{первого порядка}). Будет приведён оптимальный \textit{метод Чебышёва}. От этого  метода, оптимального для квадратичных задач, мы перейдём к построению его аналогов (оптимальных / ускоренных градиентных методов) для решения ($\mu$-сильно) выпуклых гладких (у целевой функции есть $L$-липшицев градиент) задач оптимизации. 

Везде в дальнейшем в этом разделе под $\|\,\cdot\,\|$ будем понимать обычную евклидову норму в пространстве $\R^n$ ($\|\,\cdot\,\| = \|\,\cdot\,\|_2$).

Итак, начнём с рассмотрения \textit{задачи квадратичной оптимизации}
\begin{equation}\label{quadratic}
    \min_{x\in\R^n} f(x):=\frac{1}{2}\langle x,Ax \rangle - \langle b,x \rangle.
\end{equation}
Решение этой задачи эквивалентно решению системы линейных уравнений $Ax = b$. Будем предполагать, что $\mu I_n \preceq A \preceq L I_n$,  где $I_n$ -- единичная матрица $n\times n$ (с  единицами на диагонали и нулями вне диагонали). Другими словами, спектр симметричной матрицы $A$ лежит в отрезке $[\mu,L]$. Не ограничивая общности можно считать, что $\lambda_{\min}(A) =\mu\ge 0$, $\lambda_{\max}(A) = L$. Константа $\mu$ -- есть константа сильной выпуклости $f$, а константа $L$ -- константа Липшица градиента $f$.

В качестве численного метода решения данной задачи возьмём однопараметрическое семейство (с параметром $h$ -- размер шага) методов \textit{градиентного спуска}:
\begin{equation}\label{GD}
x_{k+1} = x_k - h\nabla f(x_k) = \left(I_n - hA\right)x_k + hb.
\end{equation}
Обозначая через $x^*$ решение задачи \eqref{quadratic} ($Ax^* = b$), можно переписать \eqref{GD} следующим образом:
\begin{equation*}
x_{k+1} - x^* = \left(I_n - hA\right)\left(x_{k} - x^*\right).
\end{equation*}
Следовательно,
\begin{align}\label{argument}
\|x_{N} - x^*\| = \|\left(I_n - hA\right)^N \left(x_{0} - x^*\right)\| \le \|I_n - hA\|^N\|x_0 - x^*\|, 
\end{align}
\begin{align}\label{FunVal}
f(x_{N}) - f(x^*) = \langle \left(x_{0} - x^*\right), A\left(I_n - hA\right)^{2N} \left(x_{0} - x^*\right)\rangle.
\end{align}
Из \eqref{argument} следует, что для получения оценки скорости сходимости во всём классе квадратичных функций (с константой сильной выпуклости $\mu$ и константой Липшица градиента $L$) надо оценить
\begin{align*}
\max_{\mu I_n \preceq A \preceq L I_n}\|I_n - hA\| \le \max_{\mu\le\lambda\le L} \left| 1 - h\lambda \right| &\le \max_{\mu\le\lambda\le L} \left\{h\lambda - 1, 1 - h\lambda \right\} \le\\ &\le \max \left \{ hL - 1, 1 - h\mu \right \}.
\end{align*}
Первое неравенство следует из того, что 2-норма симметричной матрицы ($\|B\| = \max_{\|x\| = 1} \|Bx\|$, где, напомним, что $\|\,\cdot\,\| = \|\,\cdot\,\|_2$) равняется её максимальному (по модулю) собственному значению. Максимум в выражении $\max \left \{ hL - 1, 1 - h\mu \right \}$ достигается при 
$$h = \frac{2}{L+\mu}.$$
В этом случае 
$$\max \left\{hL - 1, 1 - h\mu \right\} = \frac{L-\mu}{L+\mu}.$$
Поэтому оценка \eqref{argument} примет вид
\begin{align}\label{argument_h_opt}
\|x_{N} - x^*\| \le \left(\frac{L-\mu}{L+\mu}\right)^N \|x_{0} - x^*\|. 
\end{align}
Недостатком выбора шага $h = 2/(L+\mu)$ является необходимость знать $\mu$. Если выбирать шаг по-другому:
$$h = \frac{1}{L},$$
то в этом случае
$$\max \left\{hL - 1, 1 - h\mu \right\} = 1 - \frac{\mu}{L}.$$
Оценка \eqref{argument} примет вид
\begin{align}\label{argument_L}
\|x_{N} - x^*\| \le \left(1 - \frac{\mu}{L}\right)^N \|x_{0} - x^*\|\le \exp\left(-\frac{\mu}{L}N\right)\|x_{0} - x^*\|. 
\end{align}
Так будет сходиться \textit{метод градиентного спуска}:
\begin{equation}\label{GDL}
    x_{k+1} = x_k - \frac{1}{L}\nabla f(x_k)
\end{equation}
в случае, если квадратичная задача \eqref{quadratic} будет $\mu$-сильно выпуклая и иметь $L$-липшицев градиент. 

В связи с полученным результатом возникает много вопросов.
\begin{itemize}
    \item Можно ли улучшить оценку \eqref{argument_L}?
    
    \textbf{Ответ.} Нет, для метода \eqref{GDL} на рассматриваемом классе задач оценку в общем случае (без дополнительных предположений) принципиально улучшить нельзя. Можно улучшить до оценки \eqref{argument_h_opt} за счёт лучшего выбора шага, но не более.
    \item (\textit{вырожденный случай}) Можно ли улучшить оценку \eqref{argument_L} если $\mu$ достаточно мало (в частности, $\mu = 0$)? 
    
    \textbf{Ответ.} Да, если поменять критерий сходимости. В этом случае уже нельзя рассчитывать на сходимость по аргументу \eqref{argument_L}. Однако можно рассчитывать на сходимость по функции. Действительно, рассуждая аналогично, из \eqref{FunVal} можно получить для любого $h \le 1/L$ (разумно выбирать $h = 1/L$, см. ниже):  
\begin{align*}
f(x_{N}) - f(x^*) &= \langle \left(x_{0} - x^*\right), A\left(I_n - hA\right)^{2N} \left(x_{0} - x^*\right)\rangle \\ &\le \max_{\mu\le\lambda\le L} \lambda\left|\left(1 - h\lambda\right)\right|^{2N}\cdot\|x_0 - x^*\|^2;
\end{align*}
\begin{align*}
\max_{\mu\le\lambda\le L} \lambda\left|\left(1 - h\lambda\right)\right|^{2N} &\le \max_{\mu\le\lambda\le L} \lambda \exp\left(-2Nh\lambda\right) = \\ &= \frac{1}{2Nh}\max_{\mu\le\lambda\le L} 2Nh\lambda\exp\left(-2Nh\lambda\right) \le \\ &\le 
\frac{1}{2Nh}\max_{\chi \ge 0} \chi\exp\left(-\chi\right) = \frac{1}{2eNh}\le \frac{1}{4Nh}.
\end{align*}
Таким образом,
 \begin{align*}
f(x_{N}) - f(x^*)\le \frac{LR^2}{4N},
\end{align*}
где $R = \|x_0 - x^*\|$.
\item Можно ли улучшить оценку \eqref{argument_L}, если выйти из класса методов вида \eqref{GD}, перейдя к более общему классу методов вида
$$x_{k+1} \in x_0 + \text{span}\left\{\nabla f(x_0),...,\nabla f(x_k)\right\},$$
где $\text{span}\left\{\nabla f(x_0),...,\nabla f(x_k)\right\}$ -- линейное пространство, натянутое на векторы $\nabla f(x_0),...,\nabla f(x_k)$?  

\textbf{Ответ.} Да, \textit{метод Чебышёва} будет давать оптимальные оценки\footnote{Что можно понимать следующим образом. Какой бы другой метод из рассматриваемого класса методов мы не взяли, существует такая квадратичная функция ($\mu$-сильно выпуклая с константой Липшица градиента $L$), на которой рассматриваемый метод будет сходиться хуже (не лучше), чем сходится метод Чебышёва на любой функции из рассматриваемого класса функций ($\mu$-сильно выпуклых с константой Липшица градиента $L$). В описанной ситуации говорят также о минимаксной оптимальности метода. Метод оптимален на рассматриваемом классе целевых функций и алгоритмов, если для наиболее плохой (для себя) функции из рассматриваемого класса он сходится не хуже, чем все остальные методы сходятся на своих самых плохих функциях. В основном в пособии мы стараемся строить именно такие методы -- минимаксно оптимальные. Говоря более строго, минимаксно оптимальные с точностью до числового множителя. Ввиду сложности подбора минимаксно оптимального метода на практике бывает достаточно найти метод, который в худшем случае проигрывает не более нескольких раз (например, 2--4 раза) в скорости сходимости минимаксно оптимальному методу (на его худшем случае). Отметим, что хотя минимаксная оптимальность метода совсем не означает, что метод хорошо работает на практике, как правило, всё же такой метод, и правда, будет хорошо работать на практике. Именно это небольшое замечание является одним из основных аргументов в современной прикладной оптимизации в пользу поиска практически эффективных методов среди  минимаксно оптимальных методов или близких к ним методов.}, которые в общем случае (без дополнительных предположений) не могут быть улучшены (ограничимся случаем, когда $\mu$ не мало -- см. выше):
\begin{equation}\label{Chebyshev}
x_{k+1} = x_k - \frac{4\delta_k}{L-\mu}\nabla f(x_k) + \left( \frac{2\delta_k\left(L+\mu\right)}{L-\mu} - 1\right)\left(x_{k} - x_{k-1}\right),
\end{equation}
$$x_1 = x_0 - \frac{2}{L+\mu}\nabla f(x_0),$$
$$\delta_{k+1} = \frac{1}{2\frac{L+\mu}{L-\mu} - \delta_k}, \quad \delta_1 = \frac{1}{2\frac{L+\mu}{L-\mu} + 1}.$$
Метод Чебышёва сходится следующим образом:
\begin{align*}
\|x_N - x^*\| \le \frac{2}{\xi^N+\xi^{-N}}\|x_0 - x^*\| \le \\ 2\left(\frac{\sqrt{L} - \sqrt{\mu}}{\sqrt{L} + \sqrt{\mu}}\right)^N\|x_0 - x^*\| \lesssim  2\exp\left(-2\sqrt{\frac{\mu}{L}}N\right)\|x_0 - x^*\|,
\end{align*}
где $$\xi=\frac{\sqrt{L}+\sqrt{\mu}}{\sqrt{L}-\sqrt{\mu}}.$$
Внимательный читатель, знакомый с основами численных методов оптимизации, в этом месте по идее должен был поставить под сомнение написанное выше об оптимальности метода Чебышёва. Ведь есть \textit{метод сопряжённых градиентов}, имеющий следующий вид \cite{Polyak} (глава 3, \S~2, п.~2):
$$x_{k+1} = x_k - \alpha_k \nabla f(x_k) + \beta_k\cdot \left(x_k - x_{k-1}\right),$$
где 
$$\left(\alpha_k,\beta_k\right) = \argmin_{\left(\alpha,\beta\right)} f\left(x_k - \alpha \nabla f(x_k) + \beta\cdot \left(x_k - x_{k-1}\right)\right),$$
что означает
$$\alpha_k = \frac{\|r_k\|^2\langle p_k,Ap_k\rangle - \langle r_k,p_k\rangle \langle p_k,Ar_k\rangle}{\langle r_k,Ar_k\rangle \langle p_k,Ap_k\rangle - \langle p_k,Ar_k\rangle^2},$$
$$\beta_k = \frac{\|r_k\|^2\langle p_k,Ar_k\rangle - \langle r_k,p_k\rangle \langle r_k,Ar_k\rangle}{\langle r_k,Ar_k\rangle \langle p_k,Ap_k\rangle - \langle p_k,Ar_k\rangle^2},$$
$$r_k = \nabla f(x_k) = Ax_k - b,\quad p_k = x_k - x_{k-1},$$
который также является оптимальным методом решения задачи выпуклой (не обязательно сильно выпуклой, т.е. допускается $\mu = 0$) квадратичной оптимизации. Метод сходится следующим образом:  
$$f(x_N)-f(x_*) \lesssim \min\left\{\frac{LR^2}{2(2N+1)^2},2LR^2\exp\left(-2\sqrt{\frac{\mu}{L}}N\right)\right\},$$
где $R = \|x_0 - x^*\|$. Однако никакого противоречия тут нет. Метод Чебышёва строится на основе \textit{многочленов Чебышёва} так же, как и метод сопряжённых градиентов. В некотором смысле можно просто говорить о разных формах записи идейно одного и того же подхода. Форма записи в виде сопряжённых градиентов обладает очевидными преимуществами -- отсутствием необходимости знать параметры $\mu$, $L$ и равномерностью по $\mu$ (не требуется оговорок о том, что $\mu$ не мал, см. выше). В таких случаях говорят о полностью \textit{адаптивных методах}. С другой стороны, метод Чебышёва имеет максимально приближённый вид к градиентному спуску, поскольку содержит только градиент и коэффициенты (шаг и коэффициент при <<моментном члене>> $x_k - x_{k-1}$), не зависящие от вычисленных градиентов и генерируемых методом точек, в отличие от метода сопряжённых градиентов, в котором шаг и коэффициент при моментном члене, очевидным образом, зависят и от градиентов и от генерируемых точек. На первый взгляд, это не является недостатком (в вычислительном плане это всё практически одинаково трудозатратно). Однако на базе методов типа Чебышёва оказалось намного проще строить простые аналоги (методы без вспомогательных маломерных оптимизаций) описываемых здесь оптимальных методов для более общего класса задач гладкой (сильно) выпуклой оптимизации. О чём пойдёт речь далее.     
\item Собственно, заключительный вопрос, к которому мы естественным образом подошли в повествовании, может быть сформулирован таким образом: Какие методы, в классе методов вида
$$x_{k+1} \in x_0 + \text{span}\left\{\nabla f(x_0),...,\nabla f(x_k)\right\},$$
будут оптимальными для задач выпуклой оптимизации с $L$-лип\-ши\-це\-вым градиентом и для задач $\mu$-сильно выпуклой оптимизации с $L$-лип\-ши\-це\-вым градиентом?

Ответу на этот вопрос посвящена оставшаяся часть данного раздела. 
\end{itemize}

Мы не будем приводить в данном разделе технику получения оптимальных методов для (сильно) выпуклых задач, отметим лишь, что эта техника существенным образом использует полуопределённое программирование \cite{TaylorHendrickxGlineur17}. Немного подробнее мы поговорим об этом в следующем разделе на примере получения точной оценки скорости сходимости градиентного спуска на классе задач гладкой выпуклой оптимизации. 

В данном разделе мы лишь приведём соответствующие оптимальные методы и оценки их скорости сходимости.

Но начнём мы не с оптимальных методов, а с важного (не только в историческом плане, но и в практическом) вырожденного случая метода Чебышёва. Речь идет о \textit{методе тяжёлого шарика Б.\,Т.~Поляка} \cite{Polyak}:
$$x_{k+1} = x_{k} - \frac{4}{\left(\sqrt{L}+\sqrt{\mu}\right)^2}\nabla f(x_k) + \frac{\left(\sqrt{L}-\sqrt{\mu}\right)^2}{\left(\sqrt{L}+\sqrt{\mu}\right)^2}\left(x_k - x_{k-1}\right).$$
Этот метод получается в асимптотике $k\to\infty$ из метода Чебышёва, поскольку 
$$\delta_{\infty} = \frac{1}{2\frac{L-\mu}{L+\mu}-\delta_{\infty}},$$
следовательно,
$$\delta_{\infty} = \frac{\sqrt{L} - \sqrt{\mu}}{\sqrt{L} + \sqrt{\mu}}.$$
Метод тяжёлого шарика (также часто можно встретить названия \textit{метод Поляка}, \textit{импульсный метод}) был первым методом (1964 г.), который для невырожденных задач гладкой выпуклой оптимизации локально (в окрестности минимума) сходился  как оптимальные методы сходятся (глобально) для квадратичных задач. Вопрос о построении методов, сходящихся глобально подобно методу тяжёлого шарика, был закрыт лишь спустя 15 лет А.\,С.~Немировским \cite{NemirYd79} в классе градиентных методов, допускающих вспомогательную маломерную оптимизацию. Избавиться от маломерной оптимизации (для квадратичных задач см. выше связь метода сопряжённых градиентов к методу Чебышёва) удалось спустя ещё несколько лет (в 1983 году). В~статье \nk{из журнала <<Доклады академии наук СССР>>} (статью представлял в доклады акад. Л.\,В.~Канторович) Ю.\,Е.~Нестеров сумел предложить такую вариацию метода тяжёлого шарика, которая уже глобально сходилась подобно тому, как метод тяжёлого шарика сходился локально. Далее будут приведены современные варианты методов, предложенных Ю.\,Е.~Нестеровым. Их~отличие от исходных \textit{ускоренных} (также говорят \textit{быстрых}, \textit{моментных}) методов Ю.\,Е.~Нестерова заключается в том, что они доказуемо оптимальные (без оговорок типа <<с точностью до числового множителя>>).

В алгоритме~\ref{Taylor} описывается оптимальный метод Тейлора--Дрори \cite{TaylorDrori21prep} для класса $\mu$-сильно выпуклых задач с $L$-липшицевым градиентом. 

\floatname{algorithm}{Алгоритм}
	\begin{algorithm}
		\caption{Ускоренный метод Тейлора--Дрори}\label{Taylor}
		\begin{algorithmic}[1]
			\REQUIRE $f$~--- $\mu$-сильно выпуклая функция с $L$-липшицевым градиентом, точка старта $x_0$
\STATE \textbf{SET:} $z_0=x_0$, $A_{0} = 0$, $q=\mu/L$
\FOR{$k=0, \, \dots, \, N-1$}
\STATE $A_{k+1}=\frac{(1+q)A_k+2\left(1 + \sqrt{(1+A_k)(1+qA_k)}\right)}{(1-q)^2}$
\STATE $\tau_k = 1 - \frac{A_k}{(1-q)A_{k+1}}$ и $\delta_k = \frac{1}{2}\frac{(1-q)^2A_{k+1} - (1+q)A_k}{1+q+qA_k}$
\STATE $y_k = x_k + \tau_k\left(z_k - x_k\right)$
\STATE $x_{k+1} = y_{k} - \frac{1}{L}\nabla f(y_k)$
\STATE $z_{k+1} = (1 - q\delta_k)z_k + q\delta_k y_k - \frac{\delta_k}{L}\nabla f(y_k)$
\ENDFOR
\RETURN $z_N$
\end{algorithmic}
\end{algorithm}

Ускоренный метод Тейлора--Дрори сходится следующим образом ($x^*$ -- решение задачи):
$$\|z_N - x^*\|\le\frac{1}{\sqrt{1+qA_N}}\|z_0 - x^*\|\le \sqrt{\frac{L}{\mu}}\exp\left(-\sqrt{\frac{\mu}{L}}N\right)\|z_0 - x^*\|.$$
При этом теми же авторами была установлена и нижняя оценка при $2N+1 \le n$ ($n$ -- размерность пространства) \cite{DroriTaylor21prep}:
$$\|x_N - x^*\|\ge\frac{1}{\sqrt{1+qA_N}}\|x_0 - x^*\|,$$
которая справедлива для любого метода вида
$$x_{k+1} \in x_0 + \text{span}\left\{\nabla f(x_0),...,\nabla f(x_k)\right\}$$ 
на классе задач $\mu$-сильно выпуклой оптимизации с $L$-липшицевым градиентом.

В алгоритме~\ref{KimFessler} описывается оптимальный метод Кима--Фесслера для класса выпуклых задач с $L$-липшицевым градиентом \cite{KimFessler16}.

\floatname{algorithm}{Алгоритм}
	\begin{algorithm}
		\caption{Ускоренный метод Кима--Фесслера}\label{KimFessler}
		\begin{algorithmic}[1]
			\REQUIRE $f$~--- выпуклая функция с $L$-липшицевым градиентом, точка старта $x_0$
\STATE \textbf{SET:} $y_0=x_0$, $\theta_{0} = 1$
\FOR{$k=0, \, \dots, \, N-1$}
\STATE $\theta_{k+1}=\frac{1+\sqrt{4\theta_{k}^2 + 1}}{2}$,~~$k\le N-2$;\quad $\theta_{k+1}=\frac{1+\sqrt{8\theta_{k}^2 + 1}}{2}$,~~$k = N-1$
\STATE $x_{k+1} = y_{k} - \frac{1}{L}\nabla f(y_k)$
\STATE $y_{k+1} = x_{k+1} +\frac{\theta_{k} - 1}{\theta_{k+1}}\left(x_{k+1} - x_k\right) + \frac{\theta_k}{\theta_{k+1}}\left(x_{k+1} - y_k\right)$
\ENDFOR
 \RETURN $y_N$
\end{algorithmic}
\end{algorithm}

Ускоренный метод Кима--Фесслера сходится следующим образом ($x^*$ -- решение задачи):
$$f(y_N) - f(x^*)\le\frac{LR^2}{2\theta_N^2}\le \frac{
LR^2}{\left(N+1\right)^2},$$
где $R = \|x_0 - x^*\|$.

При этом Ю. Дрори \cite{Drori17} была установлена и нижняя оценка при $N+1 \le n$ ($n$ -- размерность пространства):
$$f(\ag{x}_N) - f(x^*)\ge\frac{LR^2}{2\theta_N^2},$$
которая справедлива для любого метода вида
$$x_{k+1} \in x_0 + \text{span}\left\{\nabla f(x_0),...,\nabla f(x_k)\right\}$$ 
на классе задач выпуклой оптимизации с $L$-липшицевым градиентом.


\section[Оценка скорости сходимости градиентного спуска]{Оценка скорости сходимости градиентного спуска для задач гладкой выпуклой\\ оптимизации%
\sectionmark{Оценка скорости сходимости градиентного спуска}
}
\sectionmark{Оценка скорости сходимости градиентного спуска}
\label{sect:gradient_speed}

Итеративные методы оптимизации порядка $m$ работают с локальной информацией о минимизируемой функции, т.е. следующая точка $x_{k+1}$, которую выдаёт метод, вычисляется на основе значений самой функции и её производных до порядка $m$ в предыдущих точках $x_0, \, \dots, \, x_k$. Таким образом, каждая новая точка зависит только от конечного числа переменных. На двух функциях $f_1,f_2$ метод будет выдавать одинаковую последовательность итераций, пока значения $f_1^{(j)}(x_k),f_2^{(j)}(x_k)$ функций и их производных до порядка $m$ будут совпадать в точках $x_k$.

В некоторых случаях это обстоятельство позволяет установить точную грань для скорости сходимости метода на некотором классе функций, путём редукции этой проблемы к конечномерной задаче оптимизации. В этом разделе мы рассмотрим конкретный пример, а именно: анализ градиентного спуска с постоянным шагом на классе сильно выпуклых функций с липшицевым градиентом. В работе \cite{TaylorHendrickxGlineur17} эта проблема была сведена к полуопределённой программе.

Пусть ${\cal F}_{\mu,L}$~--- класс функций $f \in C^1(\mathbb R^d)$ с константой сильной выпуклости $\mu$ и с $L$-липшицевым градиентом, $0 < \mu < L < \infty$. В качестве граничных случаев мы получаем классы ${\cal F}_{0,L}$ выпуклых функций с $L$-липшицевым градиентом, ${\cal F}_{\mu,\infty}$ $\mu$-сильно выпуклых функций, и ${\cal F}_{0,\infty}$ выпуклых функций. Пусть $x^*$~--- минимум функции $f$. Это условие можно записать в виде $\nabla f(x^*) = 0$. Рассмотрим метод градиентного спуска:
\begin{equation} \label{gamma_grad} 
x_{j+1} = x_j - \gamma \nabla f(x_j),
\end{equation}
где $\gamma > 0$~--- фиксированная константа. Предположим, что начальная точка $x_0$ находится в шаре $B(x^*,R)$ радиуса $R$ с центром в оптимуме, т.е. $||x_0 - x^*|| \leq R$ (здесь, как и в предыдущем разделе, норма евклидова). Мы хотим найти точную верхнюю грань расстояния до оптимума после $k$ итераций, т.е. решить бесконечномерную невыпуклую проблему оптимизации:
\[ \sup_{x_0,x^*,f \in {\cal F}_{\mu,L}} ||x_k - x^*||\ :\quad ||x_0 - x^*|| \leq R,\ \nabla f(x^*) = 0,
\]
где точки $x_1, \, \dots, \, x_k$ задаются по формуле \eqref{gamma_grad}. Отметим, что в их определении участвует только конечное количество величин\footnote{С похожей ситуацией мы столкнёмся в п.~\ref{ch3_subs_moments}, где рассмотрим моментные релаксации для полиномиальных задач оптимизации. В этом контексте неотрицательная мера $\mu$ входит в формулировку задачи только посредством конечного количества своих линейных проекций~--- моментов. Отметим, что в нашей ситуации мы не имеем дело с линейными проекциями переменных $x_0,f$ задачи. Дело в том, что хотя условия \eqref{sect:gradient_speed} линейны по слагаемым, выражение $\nabla f(x_j)$ линейно только по $f,$ но не по $x_j$.}, а именно, $x_0, \, \nabla f(x_0), \, \dots, \, \nabla f(x_{k-1})$. Поэтому задачу можно переписать в виде
\[ \sup_{x_0,\dots,x_k,x^*,g_0,\dots,g_{k-1}} ||x_k - x^*||\ :
\]
\[ \exists\ f \in {\cal F}_{\mu,L}:\ ||x_0 - x^*|| \leq R,\ g_j = \nabla f(x_j),\ \nabla f(x^*) = 0,
\]
\[x_{j+1} = x_j - \gamma g_j, j = 0,\dots,k-1.
\]

Ключевое наблюдение, ведущее к основному результату работы \cite{TaylorHendrickxGlineur17}, заключается в том, что условие $\exists\ f \in {\cal F}_{\mu,L}$ на бесконечномерный объект $f$ можно выразить квадратичными условиями на конечное число скаляров $f(x_0),\dots,f(x_{k-1})$ и векторов $x_0,\dots,x_{k-1},\nabla f(x_0),\dots,\nabla f(x_{k-1})$. Начнём с граничного класса ${\cal F}_{0,\infty}$ выпуклых функций, для которого характеризация вполне очевидна.

\begin{lemma}
Пусть даны точки $x_0,\dots,x_{k-1} \in \mathbb R^d$, скаляры $c_0,\dots,c_{k-1} \in \mathbb R$, и векторы $g_0,\dots,g_{k-1} \in \mathbb R^d$. Следующие условия эквивалентны:
\begin{itemize}
\item Существует функция $f \in {\cal F}_{0,\infty}$ такая, что $f(x_j) = c_j$, $g_j \in \partial f(x_j)$ для всех $j = 0,\dots,k-1$.
\item Имеют место неравенства $c_j \geq c_i + \langle g_i,x_j - x_i \rangle$ для всех $i,\linebreak j = 0,\dots,k-1$.
\end{itemize}
\end{lemma}

\begin{proof}
Необходимость неравенств очевидна вследствие выпуклости функции $f$.

Достаточность доказывается построением конкретной выпуклой функции из данных скаляров $c_j$ и векторов $x_j,g_j$, а именно максимума линейных функций:
\[ f(x) = \max_{j=0, \, \dots, \, k-1} \, (c_j + \langle g_j,x-x_j \rangle). \qedhere
\]
\end{proof}

Построение условий для более сложных классов функций сводится к вышеописанному простому случаю. А именно:
\[ f \in {\cal F}_{\mu,L} \ \Leftrightarrow\ f - \frac{\mu^2}{2}||x||^2 \in {\cal F}_{0,L - \mu};
\]
\[ f \in {\cal F}_{0,L}\ \Leftrightarrow\ f^* \in {\cal F}_{L^{-1},\infty};
\]
\[ f \in {\cal F}_{\mu,\infty} \ \Leftrightarrow\ f - \frac{\mu^2}{2}||x||^2 \in {\cal F}_{0,\infty}.
\]
Здесь $f^*(p) = \sup_{x \in \mathbb R^d} (\langle x,p \rangle - f(x))$ является сопряжённой к $f$ функцией. Комбинируя эти условия, приходим к следующему результату \cite[Теорема~4]{TaylorHendrickxGlineur17}.

\begin{teo}
Пусть даны точки $x_0, \, \dots, \, x_{k-1} \in \mathbb R^d$, скаляры $c_0, \, \dots, \, c_{k-1} \, \in \, \mathbb R$, и векторы $g_0, \, \dots, \, g_{k-1} \, \in \,  \mathbb R^d$. Следующие условия эквивалентны:
\begin{itemize}
\item Существует функция $f \in {\cal F}_{\mu,L}$ такая, что $f(x_j) = c_j$, $g_j = \linebreak =\nabla f(x_j)$ для всех $j = 0, \, \dots, \, k-1$.
\item Имеют место неравенства 
\begin{align} \label{ij_quadratic} 
c_j \geq c_i &+ \langle g_i,x_j - x_i \rangle + \frac{1}{2L}\|g_j - g_i\|^2 +\nonumber\\ &+\frac{\mu}{2(1-\mu/L)}\left\|x_j - x_i - L^{-1}(g_j - g_i)\right\|^2
\end{align}
для всех $i, \, j = 0, \, \dots, \, k-1$.
\end{itemize}
\end{teo}

Таким образом, проблему можно переписать в виде
\[ \sup_{x_0, \, \dots, \, x_k, \, x_*, \,  c_0, \, \dots, \, c_{k-1}, \, c_*, \, g_0, \, \dots, \, g_{k-1}, \, g_*} ||x_k - x^*|| : 
\]
\[
||x_0 - x^*|| \leq R,\, g_* = 0,\, x_{j+1} = x_j - \gamma g_j
\]
при условиях \eqref{ij_quadratic} для всех $i,j = 0, \, \dots, \, k-1,*$. 

Это невыпуклая задача квадратичной оптимизации. Однако для неё легко выписать полуопределённую релаксацию. Сформируем грамиан $\Gamma$ векторов $x_0-x_*, \, \dots, \, x_k-x_*, \, g_0, \, \dots, \, g_{k-1}, \, g_*$. Тогда все выражения, квадратичные по этим векторам, станут линейными по элементам $\Gamma$. Линейные по векторам условия равенства перейдут в такие же линейные условия равенства на соответствующие столбцы $\Gamma$. В частности, последний столбец~$\Gamma$ равен нулю вследствие условия $g_* = 0$, и его можно опустить. Условие, что все векторы принадлежат пространству $\mathbb R^d$, перейдёт в условие на ранг $\rk\,\Gamma \leq d$. Это условие на ранг является единственным невыпуклым условием в новой формулировке задачи, и его опущением мы получаем искомую полуопределённую релаксацию.

В случае, когда размерность $d$ пространства не меньше размера $2k+1$ грамиана $\Gamma$, условие на ранг тривиально выполнено, и релаксация точна. В~частности, точность имеет место, когда мы изначально не фиксируем размерность пространства, а изучаем скорость сходимости метода на задачах произвольной размерности.

Решая получившуюся полуопределённую программу для разных значений $\gamma$, мы можем численно определить оптимальную длину шага для фиксированного числа итераций $k$. Легко видеть, что все приведённые выше рассуждения остаются в силе, если длина шага зависит от номера итерации, т.е. метод вычисляет следующую точку по формуле
\[ x_{j+1} = x_j - \gamma_j \nabla f(x_j)
\]
с различными константами $\gamma_j$.

\medskip

Схожие рассуждения можно провести и для других итеративных методов и критериев сходимости, например, сходимости по значению функции, когда максимизируется выражение $f(x_k) - f(x^*)$, или по норме градиента $\|\nabla f(x_k)\|$. Заметим, что скорость сходимости градиентного спуска на классе функций ${\cal F}_{0,L}$ была рассмотрена в работе \cite{DroriTeboulle14}, в которой представленный в этом разделе метод анализа был опубликован впервые.

Численный эксперимент показывает \cite{TaylorHendrickxGlineur17}, что для критерия сходимости по значению функции наихудшая функция является одномерной, $f: \mathbb R \to \mathbb R$, а именно вида
\[ f(x) = \left\{ \begin{array}{rcl} \frac{\mu}{2}|x - x^*|^2 + (L-\mu)\tau|x - x^*| - \frac{L-\mu}{2}\tau^2, &\quad& |x - x^*| \geq \tau, \\ \frac{L}{2}|x - x^*|^2,&&|x - x^*| \leq \tau, \end{array} \right.
\]
где либо $\tau = \frac{R\mu}{L(1 - \mu\gamma)^{-2k} - (L - \mu)}$, либо $\tau = +\infty$, в зависимости от того, при каком значении достигается максимум разницы $f(x_k) - f(x_*)$. Это приводит к следующей оценке скорости сходимости:
\[ f(x_k) - f(x^*) \leq \frac{LR^2}{2}\max\left\{\frac{\mu}{L(1 - \mu\gamma)^{-2k} - (L - \mu)},(1-\gamma L)^{2k}\right\}.
\]
Таким образом, оптимальная длина шага $\gamma$, которая минимизирует максимум на правой стороне неравенства, является некоторой функцией от числа итераций $k$. Нетрудно проверить, что максимум минимизируется при значении $\gamma$, при котором оба выражения в фигурных скобках равны. При $k \to \infty$ оптимальная длина шага стремится к $\gamma = \frac{2}{L + \mu}$. А именно, с точностью до более высоких по $k^{-1}$ членов мы имеем
\[ \gamma_{opt} \approx \frac{1 + \kappa^{1/2k}}{1 + \kappa^{1 + 1/2k}},
\]
где $\kappa = \frac{\mu}{L}$. Используя оптимальный шаг вместо стандартного суб-оптимального $\gamma = \frac{2}{L + \mu}$, получаем скорость сходимости
\[ f(x_k) - f(x^*) \lesssim \frac{LR^2}{2}\kappa^{\frac{L}{\mu + L}}\left( \frac{L - \mu}{L + \mu} \right)^{2k}
\]
вместо $\frac{LR^2}{2}\left( \frac{L - \mu}{L + \mu} \right)^{2k}$, что даёт выигрыш в множитель $\left( \frac{\mu}{L} \right)^{\frac{L}{\mu + L}}$.

%
%

\section[Метод условного градиента, или алгоритм Франк--Вульфа]{Метод условного градиента, или алгоритм Франк--Вульфа} 
\label{ch2_sect_fw}

Рассматривается задача с ограничениями
\begin{equation}
\label{ch2_q_cond_grad_prob}
\min_{x \, \in \, Q} f(x),
\end{equation}
где $f(x) : \mathbb{R}^n \rightarrow \mathbb{R}$~--- выпуклая, непрерывно дифференцируемая функция, а множество~$Q$~--- компактное и выпуклое. При этом мы предполагаем, что оптимизировать \emph{линейную} функцию на $Q$~--- это более простая задача, чем исходная \eqref{ch2_q_cond_grad_prob}. 

В основе метода лежит та же идея линеаризации целевой функции, что и в основе градиентного метода: в очередной точке $x_k$ линеаризуем функцию~$f(x)$, затем
решаем задачу минимизации линейной аппроксимации на~$Q$, после этого
найденную точку используем для выбора направления движения (см. алгоритм~\ref{ch2_alg_cond_grad}). 

Применимость метода условного градиента зависит от того, насколько просто решить задачу минимизации линейной функции на~$Q$. В случае, когда $Q$~--- политоп, эта подзадача является линейной программой. В ряде других случаев подзадача решается в явном виде, например, когда $Q$ является единичным шаром в какой-либо простой норме, матричной или векторной \cite{Jaggi_Rev_FW} или проекцией такого шара, например, {\it зонотопом}\footnote{Зонотоп~--- это политоп вида $Z = \left\{ \left. \sum_{i=1}^m \sigma_i u_i \,\right|\, \sigma_i \in [-1,1] \right\}$, где $u_i \in \mathbb R^n$~--- фиксированные векторы. То есть зонотоп, это линейная проекция $m$-мерного куба. Минимум линейного функционала $c$ по $Z$ равен $-\sum_{i=1}^m |\langle c,u_i \rangle|$.}. 

Впервые метод условного градиента был представлен в 1956~г. в работе~\cite{FrankWolfe}.
В последние годы продолжается активное развитие этого направления
и появляются новые современные модификации метода Франк--Вульфа и
его обобщения на другие классы задач (см., например, \cite{Jaggi_Rev_FW, DoikovNest_contr_p_meth_2020}). 

У итерации метода условного градиента, который применяется для решения~\eqref{ch2_q_cond_grad_prob}, есть три замечательных свойства, которые позволяют эффективно применять его в важных приложениях (см., например, п.~\ref{ch3_sect_ls_regr}): 1)~отсутствие необходимости делать настоящую проекцию;
2)~независимость от типа нормы; 3)~возможность использования разреженности. 

\floatname{algorithm}{Алгоритм}
	\begin{algorithm}
		\caption{Метод условного градиента}\label{ch2_alg_cond_grad}
		\begin{algorithmic}[1]
			\REQUIRE $f$~--- выпуклая, непрерывно дифференцируемая функция; $Q$~--- допустимое множество, выпуклое и компактное; $x_0$~--- начальная точка; $N$~--- количество итераций.
			\ENSURE точка $x_N$

\FOR{$k=0,\, \dots, \, N-1$}
\STATE $\gamma_k = \frac{2}{k+1}$, $0 \leq \gamma_k \leq 1$

\STATE $y_{k} = \argmin_{y \, \in \, Q} \langle \nabla f(x_k), \, y \rangle$ 

\STATE $x_{k + 1} = (1 - \gamma_k) x_k + \gamma_k y_k$
\ENDFOR

\RETURN $x_N$
		\end{algorithmic}

\end{algorithm}

Отметим, что на каждом шаге величина $\mu_k = f(x_k) + \langle \nabla f(x_k),y_k \rangle$ является нижней границей для оптимального значения задачи \eqref{ch2_q_cond_grad_prob}, поскольку линейная аппроксимация $f(x) + \langle \nabla f(x_k),y \rangle$ является нижней границей для $f(y)$. Таким образом, разница $f(x_k) - \mu_k = -\langle \nabla f(x_k),y_k \rangle$ ограничивает зазор двойственности.
 
Отметим также, что метод аффинно инвариантен, т.е. после аффинной замены координат метод генерирует ту же самую последовательность точек, если только начать с той же точки $x_0$ \cite{Jaggi_Rev_FW}.

\begin{teo}[\cite{Jaggi_Rev_FW}]\label{thmuslFW}
Пусть $f$~--- выпуклая функция, градиент
которой на $Q$ удовлетворяет условию Липшица с константой~$L$ по отношению к некоторой норме $\| \cdot \| :$ для всех $x,y \in Q$, выполняется
$$\|\nabla f(y) - \nabla f(x)\|_* \le L\|y-x\|,$$
$R = \sup_{x,\, y \, \in \, Q} \| x - y \|$, $\gamma_k = \frac{2}{k+1}$ для $k \geq 1$.
Тогда для любого $k \geq 2$ выполняется
\begin{equation}
\label{ch2_eq_rate_FW}
f(x_k) - f(x^*) \leq \frac{2 L R^2}{k + 2}.
\end{equation}
\end{teo}

\begin{proof} (см.~\cite{Bubeck15}) Из условия Липшица на градиент $f$ и выпуклости $f$ имеем оценку
\[ 0 \leq f(x) - f(y) - \langle \nabla f(y),x-y \rangle \leq \frac{L}{2}\|x - y\|^2
\]
для всех $x,y \in Q$. Получаем цепочку неравенств:
\begin{align*}
\lefteqn{f(x_{k+1}) - f(x_k) \leq \langle \nabla f(x_k),x_{k+1}-x_k \rangle + \frac{L}{2}\|x_{k+1} - x_k\|^2} \phantom{space} \\ 
=& \,\gamma_k \langle \nabla f(x_k),y_k-x_k \rangle + \frac{L\gamma_k^2}{2}\|y_k - x_k\|^2\\ \leq& \, \gamma_k \langle \nabla f(x_k),x^*-x_k \rangle + \frac{L\gamma_k^2}{2}R^2  \leq \,\gamma_k(f(x^*) - f(x_k)) + \gamma_k^2\frac{LR^2}{2}.
\end{align*}
Введём обозначение $\delta_k = \frac{f(x_k) - f(x^*)}{LR^2}$. Тогда неравенство перепишется в виде
\[ \delta_{k+1} \leq (1-\gamma_k)\delta_k + \frac{\gamma_k^2}{2} = \left( 1 - \frac{2}{k+1} \right)\delta_k + \frac{2}{(k+1)^2}.
\]
Начиная с неравенства $\delta_2 \leq \frac12$, применением индукции по $k$ получаем желаемый результат. \end{proof}

Поясним третью важную особенность метода условного градиента~--- использование разреженности
(см. также~\cite{Bubeck15}).
Пусть допустимое множество $Q \, \in \, \mathbb{R}^n$ является политопом, т.е. выпуклой оболочкой конечного множества точек. Тогда по теореме Каратеодори любую точку $q \, \in \, Q$ можно представить
как выпуклую комбинацию не более чем $n + 1$ вершины $Q$. В то же время на $k$-й итерации метода условного градиента $x_{k+1}$ может быть записана как выпуклая комбинация $k$~вершин (при условии, что начальная точка $x_0$~--- тоже вершина политопа~$Q$). Скорость сходимости
метода не зависит от размерности (см.~\eqref{ch2_eq_rate_FW}), поэтому
в случае $k <\!\!< n$ результирующая точка итерации есть комбинация
числа вершин допустимого политопа, которое в разы меньше, чем размерность пространства.

\section{Ускоренный Мета-алгоритм}\label{AM} 

В последние 15 лет в численных методах гладкой выпуклой оптимизации преобладают так называемые ускоренные методы. Прообразом таких методов является метод тяжёлого шарика Б.Т.~Поляка и моментный метод Ю.\,Е.~Нестерова \cite{Gas_Pos18, NesterovLectures}.
Оказалось, что для многих задач гладкой выпуклой оптимизации оптимальные методы (с точки зрения числа вычислений градиента функции; в общем случае, старших производных) могут быть найдены среди ускоренных методов \cite{Gas_Pos18, NesterovLectures, Lan2019}. 
Появилось огромное число работ, в которых предлагаются различные варианты ускоренных методов для разных классов задач, см., например, обзор литературы в \cite{Gas_Pos18, Lan2019}. Каждый раз процедура ускорения принимала свою причудливую форму. Естественно, возникло желание как-то унифицировать всё это. 
В 2015 году это было сделано для широкого класса (рандомизированных) градиентных методов с помощью проксимальной ускоренной оболочки, названной Каталист
\cite{Catalyst}.
С 2013 года данные результаты стали активно переноситься на тензорные методы (использующие старшие производные) \cite{Doikov2019, Gasnikov2019COLT, Monteiro2013, Nesterov2020}. В самое последнее время предпринимаются попытки унификации процедур ускорения  для седловых задач и задач со структурой (композитных задач) \cite{Alkousa2019, Ivanova2020, Kamzolov2020, Lin2020}. Во всех этих направлениях по-прежнему использовалось значительное разнообразие  ускоренных проксимальных оболочек \cite{Gas_Pos18, Catalyst, Doikov2019, Gasnikov2019COLT, Monteiro2013, Nesterov2020, Ivanova2020, Kamzolov2020,  Gasnikov2019COLT, Bubeck2020, Jiang2020,  Ivanova2019}. 
Метод из данного раздела взят из работы \cite{Gasnikov2021} и
будет во многом базироваться на схеме из \cite{Bubeck2020}.

В данном разделе показывается, что достаточно изучить всего одну ускоренную проксимальную оболочку, которая позволяет получать все известные нам ускоренные методы для задач гладкой выпуклой безусловной оптимизации. Причём в ряде случаев предложенный Ускоренный Мета-алгоритм позволяет убирать логарифмические зазоры в оценках сложности (по сравнению с нижними оценками), имевшие место в предыдущих подходах.

\subsection{Основные результаты}\label{section_2}

Рассмотрим следующую задачу ($x^*$ -- решение задачи):
\begin{equation}
\label{eq1}
\min\limits_{x \in \R^n } \, \{ F\left( x \right):=f\left( x \right)+g\left( x \right)\} ,
\end{equation}
где $f$ и $g$ --- выпуклые функции.

Везде в дальнейшем в разделах~\ref{AM},~\ref{AM_Applications} под $\|\,\cdot\,\|$ будем понимать обычную евклидову норму в пространстве $\R^n$ ($\|\,\cdot\,\| = \|\,\cdot\,\|_2$), $$D^k f(x)[h]^k = \sum_{i_1,...,i_d \ge 0:\,\, \sum_{j=1}^d i_j = k} \frac{\partial^k f(x)}{\partial x_1^{i_1}...\partial x_d^{i_d}}h_1^{i_1} \cdot...\cdot h_d^{i_d},$$  $$\|D^k f(x)\| = \max_{\|h\|\le 1} \left\|D^k f(x)[h]^k\right\|.$$
Будем считать, что $f$ имеет липшицевы производные порядка $p$ ($p \in \NN$):
\begin{equation}
    \|D^p f(x)- D^p f(y)\|\leq L_{p,f}\|x-y\|.
    \label{def_lipshitz}
\end{equation}
Здесь и далее (см., например, \eqref{unif_conv}) можно считать, что $x,y \in \R^n$ принадлежат евклидову шару с центром в точке $x^{*}$ и радиусом $O(\|x_0 - x^*\|)$, где $x_0$ -- точка старта \cite{Bubeck2020}.

Введём аппроксимацию рядом Тейлора функции $f$:
\begin{equation}
    \Omega_{p}(f,x;y)=f(x)+\sum_{k=1}^{p}\frac{1}{k!}D^{k}f(x)\left[ y-x \right]^k, y\in \R^n.
    \label{eq_taylor}
\end{equation}
Заметим, что из \eqref{def_lipshitz} следует \cite{NesterovImplementable}
\begin{equation}
   \left|f(y)-\Omega_{p}(f,x;y)\right| \leq \frac{L_{p,f}}{(p+1)!}\|y-x\|^{p+1}.
    \label{eq_sumup}
\end{equation}

Заметим, что при $p=1$ (методов первого порядка, т.е. градиентных методов) условие $$\frac{1}{2} \leq \lambda_{k+1} \frac{H \|y_{k+1} - \tilde{x}_k\|^{p-1}}{p!}  \leq \frac{p}{p+1}$$ может быть переписано следующим образом: $$H = \frac{1}{2\lambda_{k+1}}.$$

\begin{algorithm} [h!]
\caption{Ускоренный Мета-алгоритм (УМ) (УМ($x_0$,$f$,$g$,$p$,$H$,$K$))}\label{alg:highorder}
	\begin{algorithmic}[1]
		\STATE \textbf{Input:} $p \in \NN$, $f : \R^n \rightarrow \R$, 
		$g : \R^n \rightarrow \R$, $H > 0$. 
		\STATE $A_0 = 0, y_0 = x_0$.
		\FOR{ $k = 0$ \TO $k = K- 1$}
		\STATE Определить пару $\lambda_{k+1} > 0$ и $y_{k+1}\in \R^n$ из условий
		\[
\frac{1}{2} \leq \lambda_{k+1} \frac{H \|y_{k+1} - \tilde{x}_k\|^{p-1}}{p!}  \leq \frac{p}{p+1} \,,
\]
где
\begin{equation}
\label{prox_step}
\hspace{-2em} y_{k+1} = \argmin_{y\in \R^n} \left\{\widetilde{\Omega}^k(y):=
\Omega_{p}(f,\tilde{x}_k;y)+g(y) +\frac{H}{(p+1)!}\|y-\tilde{x}_k\|^{p+1} \right\} \,,
\end{equation}
		\[
a_{k+1} = \frac{\lambda_{k+1}+\sqrt{\lambda_{k+1}^2+4\lambda_{k+1}A_k}}{2} 
\text{ , } 
A_{k+1} = A_k+a_{k+1}
\text{ , } 
\]
\[
\tilde{x}_k = \frac{A_k}{A_{k + 1}}y_k + \frac{a_{k+1}}{A_{k+1}} x_k. 
		\]
		\STATE $x_{k+1} := x_k-a_{k+1} \nabla f(y_{k+1}) - a_{k+1}\nabla g(y_{k+1})$.
		\ENDFOR
		\RETURN $y_{K}$ 
	\end{algorithmic}	
\end{algorithm}

Доказательство следующей теоремы см. в \cite{Gasnikov2021}. 
\begin{teo} \label{theoremCATD}
Пусть $y_k$ -- выход Алгоритма~{\rm \ref{alg:highorder}} УМ{\rm(}$x_0$,$f$,$g$,$p$,$H$,$k${\rm)} после $k$ итераций при $p\geq 1$ и $H\ge (p+1)L_{p,f}$. Тогда
 \begin{equation} \label{speedCATD}
 F(y_k) - F(x^{\ast}) \leq \frac{c_p H R^{p+1}}{k^{\frac{3p +1}{2}}} \,,
 \end{equation}
где $c_p = 2^{p-1} (p+1)^{\frac{3p+1}{2}} / p!$,
$R=\|x_0 - x^{\ast}\|$. 

Более того, при $p \ge 2$ для достижения точности $\varepsilon$: $$F(y_k) - F(x_{*}) \leq \varepsilon$$ на каждой итерации УМ вспомогательную задачу \eqref{prox_step} придётся решать {\rm(}процедурой типа бинарного поиска{\rm)} для подбора пары $(\lambda_{k+1},y_{k+1})$ не более чем $O\left(\ln\left(\varepsilon^{-1}\right)\right)$ раз.
\end{teo}
Заметим, что приведённая выше теорема будет справедлива и при условии $H\ge 2L_{p,f}$ (независимо от $p \in \NN$). Это выводится из \eqref{eq_sumup}. Условие $H\ge (p+1)L_{p,f}$ было использовано, поскольку оно гарантирует выпуклость вспомогательной подзадачи \eqref{prox_step} \cite{NesterovImplementable}.  При этом условии и $g \equiv 0$ для $p = 1, 2, 3$ существуют эффективные способы решения вспомогательной задачи \eqref{prox_step} \cite{NesterovImplementable}. Для $p = 1$ существует явная формула для решения \eqref{prox_step}, для $p = 2, 3$ сложность \eqref{prox_step} такая же (с точностью до логарифмического по $\varepsilon$ множителя), как у итерации метода Ньютона \cite{NesterovImplementable}. 

Важно отметить, что вспомогательную задачу \eqref{prox_step} не обязательно решать точно: достаточно \cite{Kamzolov2020,Kamzolov2020Hyperfast} найти точку $\tilde{y}_{k+1}$, удовлетворяющую
\begin{equation}
\label{inexact1}
    \left\|\nabla \widetilde{\Omega}^k(\tilde{y}_{k+1}) \right\| \le \frac{1}{4p(p+1)}\|\nabla F(\tilde{y}_{k+1})\|.
\end{equation}
Такая модификация приведёт лишь к появлению множителя $12/5$ в правой части \eqref{speedCATD}.

Также заметим, что приведённый Алгоритм~\ref{alg:highorder} и Теорему~\ref{theoremCATD}, описывающую его сходимость, можно распространить на задачи условной оптимизации на (выпуклых) множествах простой структуры ($Q$):
\begin{equation}
\label{ch2_eq_eq1_q}
\min\limits_{x \in Q }\{ F\left( x \right):=f\left( x \right)+g\left( x \right)\} .
\end{equation}
В этом случае следует заменить формулу \eqref{prox_step} на
\begin{equation*}
y_{k+1} = \argmin_{y\in Q} \left\{\widetilde{\Omega}^k(y):=
\Omega_{p}(f,\tilde{x}_k;y)+g(y) +\frac{H}{(p+1)!}\|y-\tilde{x}_k\|^{p+1} \right\}
\end{equation*}
и формулу из строчки 5 на 
$$x_{k+1} := \argmin_{y\in Q} \left\{a_{k+1}\left(F(x_k) + \langle \nabla f(y_{k+1}), y - x_k \rangle + g(y)\right) + \frac{1}{2}\|y - x_k^2\|^2\right\}.$$

Отмеченные задачи можно решать неточно, подобно \eqref{inexact1} (см. также \cite{Doikov2019}). Однако в этом случае (а также в случае, когда $g$ --- негладкая) выполнить условие подобное \eqref{inexact1} может быть совсем не просто. В этом случае существуют более удобные критерии точности решения вспомогательной задачи, см., например, раздел~\ref{model}.

Будем говорить, что функция $F$ является $r$-равномерно выпуклой ($p+1 \geq r \geq 2$) с константой $\sigma_r > 0$, если
\begin{equation}\label{unif_conv}
    F(y)\geq F(x) + \langle \nabla F(x), y-x \rangle + \frac{\sigma_r}{r} \|y-x\|^r, \quad  x,y \in \R^n.
\end{equation}

В этом случае, используя \cite{Grapiglia2019}
\begin{equation}\label{inexact2}
  F(\tilde{y}_{k+1}) -   F(x^*) \le \frac{r-1}{r}\left(\frac{1}{\sigma_r}\right)^{\frac{1}{r-1}} \|\nabla F(\tilde{y}_{k+1})\|^{\frac{r}{r-1}},
\end{equation}
можно завязать критерий \eqref{inexact1} на желаемую точность (по функции) решения исходной задачи $\varepsilon$ \cite{Kamzolov2020}:  $\left\|\nabla \widetilde{\Omega}^k(\tilde{y}_{k+1}) \right\|  = O\left(\left(\epsilon^{r-1}\sigma_r\right)^{\frac{1}{r}}\right)$. 

Более того, для $p = 1$ приведённые здесь выкладки можно уточнить, подчеркнув тем самым, что сложность решения вспомогательной задачи может даже не зависеть от $\varepsilon$. Оказывается (см. \cite{Ivanova2019}), что условие
\begin{equation}
\label{inexact3}
\|\tilde{y}_{k+1} - y^*_{k+1}\| \le   \frac{H}{3H + 2L^g_1} \|\tilde{x}_{k} - y^*_{k+1}\|,
\end{equation}
где $y^*_{k+1}$ -- точное решение задачи \eqref{prox_step}, а $L^g_1$ -- константа Липшица градиента $\nabla g$, в теоретическом плане гарантирует то же, что и условие \eqref{inexact1} при $p=1$. А именно теорема~\ref{theoremCATD} останется верной с добавлением в правую часть \eqref{speedCATD} множителя $12/5$. 

Отметим, что оценка скорости сходимости
$\eqref{speedCATD}$ с точностью до числового множителя $c_p$ не может быть улучшена для класса выпуклых задач \eqref{eq1} с липшицевой $p$-й производной и для широкого класса тензорных методов порядка $p$, описанном в \cite{NesterovImplementable}. 
При дополнительном предположении равномерной выпуклости $F$
оптимальный метод можно построить на базе УМ с помощью процедуры рестартов \cite{Kamzolov2020} -- см. Алгоритм~\ref{alg:restarts}.

\begin{algorithm} [h!]
\caption{Рестартованный УМ($x_0$,$f$,$g$,$p$,$r$,$\sigma_r$,$H$,$K$)}\label{alg:restarts}
	\begin{algorithmic}[1]
		\STATE \textbf{Input:} $r$-равномерно выпуклая функция $F = f + g : \R^n \rightarrow \R$ с константой $\sigma_r$ и УМ($x_0$,$f$,$g$,$p$,$H$,$K$).
		\STATE $z_0=x_0$.
		\FOR{$k = 0 $ \TO $K$}
		\STATE $R_k=R_0\cdot 2^{-k}$, 
		\begin{equation}
		\label{numberofrestarts}    
		N_k=\max \left\{ \left\lceil \left(\frac{r c_p H 2^r}{\sigma_r} R_k^{p+1-r}\right)^{\frac{2}{3p+1}}  \right\rceil, 1\right\}.
		\end{equation}
		\STATE  $z_{k+1} := y_{N_k}$, где $y_{N_k}$ -- выход УМ($z_k$,$f$,$g$,$p$,$H$,$N_k$).
		\ENDFOR
		\RETURN $z_{K}$ 
	\end{algorithmic}	
\end{algorithm}

\begin{teo} \label{theoremRestartCATD}
Пусть $y_k$ -- выход алгоритма~{\rm \ref{alg:restarts}}  после $k$ итераций. Тогда если $H\ge (p+1)L_{p,f}$, $\sigma_r > 0$, то общее число вычислений \eqref{prox_step} для достижения $$F(y_k) - F(x^*) \le \varepsilon$$ будет
\begin{equation}
    N = \tilde{O} \left(\left( \frac{H R^{p+1-r}}{\sigma_r} \right)^{\frac{2}{3p+1}}\right),
\end{equation}
где $\tilde{O}(\,)$ -- означает то же самое, что $O(\,)$ с точностью до множителя
$\ln\left({\varepsilon^{-1}}\right)$.
\end{teo}
Всё, что было сказано после теоремы~\ref{theoremCATD}, можно отметить и в данном случае.

\section{Приложения Ускоренного Мета-алгоритма}\label{AM_Applications} 
В данном разделе будут продемонстрированы некоторые приложения Ускоренного Мета-алгоритма.  

\subsection{Ускоренные методы композитной оптимизации}\label{3.1} 

Если не думать о сложности решения подзадачи \eqref{prox_step}, например, считать, что $g$ какая-то простая функция и \eqref{prox_step} решается по явным формулам (как, например, для задачи LASSO), то
УМ описывает класс ускоренных методов (1, 2, 3, ... порядка) композитной оптимизации \cite{NesterovLectures}. При этом функция $g$ не обязана быть гладкой. В общем случае в строчке 5 Алгоритма~\ref{alg:highorder} под $\nabla g(y_{k+1})$ следует понимать такой субградиент функции $g$ в точке $y_{k+1}$, с которым субградиент правой части \eqref{prox_step} равен (близок) к нулю (немного переписав метод, от последнего ограничения можно отказаться). Как уже отмечалось, при $p = 1$ необходимость в поиске параметра $\lambda_{k+1}$ исчезает, что делает метод заметно проще.

\subsection{Ускоренные проксимальные методы. Каталист}\label{3.2} 

Если считать $p = 1$, a $f\equiv 0$, $H > 0$, то получится ускоренный проксимальный метод. Отличительная особенность такого метода (см. также \cite{Ivanova2019}) от других известных ускоренных проксимальных методов заключается в том, что не требуется очень точно решать вспомогательную задачу. Критерий \eqref{inexact3} и сильная (2-равномерная) выпуклость вспомогательной подзадачи \eqref{prox_step} указывают на то, что сложность решения \eqref{inexact3} может не зависеть от желаемой точности решения исходной задачи $\varepsilon$. Таким образом, не теряется логарифмический множитель при использовании такой проксимальной оболочки для ускорения различных неускоренных процедур. Собственно, последнее направление получило название Каталист \cite{Catalyst}. До настоящего момента идея (Каталист) использования ускоренной проксимальной оболочки для <<обёртывания>> неускоренных методов, решающих вспомогательную задачу \eqref{prox_step} на каждой итерации (при должном выборе параметра $H$), являлась наиболее общей идеей разработки ускоренных методов для разных задач. Мы получаем Каталист просто как частный случай УМ. Примеры использования Каталист будут приведены в п.~\ref{3.4}.

\subsection{Разделение оракульных сложностей} \label{3.3} 

Если считать, что для $g$ имеем $L_{p,g} < \infty$ (см. \eqref{def_lipshitz}) и  на вспомогательную задачу \eqref{prox_step} смотреть как на равномерно выпуклую достаточно гладкую задачу (с $f: = g$, $g(x) := 
\Omega_p\left(f,\tilde{x}_k;x\right) + 
\frac{(p+1)L_{p,f}}{\left( p+1\right)!}
\left\| x - \tilde{x}_k \right\|^{p+1}$), то для решения \eqref{prox_step} в свою очередь можно использовать Рестартованный УМ c $H \simeq (p+1)L_{p,g}$. В случае, когда $L_{p,f} \le L_{p,g}$, удаётся получить такие оценки сложности \cite{Kamzolov2020, Lan2019} (см.~теорему~\ref{theoremCATD}):
\begin{center}
$N_f = \tilde{O} \left(\left( \frac{L_{p,f} R^{p+1}}{\varepsilon} \right)^{\frac{2}{3p+1}}\right)$ -- число вызовов оракула для функции $f$,
\end{center}
\begin{center}
$N_g = \tilde{O} \left(\left( \frac{L_{p,g} R^{p+1}}{\varepsilon} \right)^{\frac{2}{3p+1}}\right)$ -- число вызовов оракула для функции $g$.
\end{center}
Вызов оракула подразумевает вычисление (старших) производных до порядка $p$ включительно. Таким образом, число вызовов оракула для каждой из функций $f$, $g$ является квазиоптимальным, т.е. оптимальным с точностью до логарифмического (от желаемой точности по функции) множителя. Аналогичные оценки можно получить и в $r$-равномерно ($r \ge 2$) выпуклом случае, см. п.~\ref{3.4}.

Заметим, что при $p = 1$ внутреннюю задачу \eqref{prox_step} не обязательно решать Рестартованным УМ. Можно использовать (ускоренные) покомпонентные и безградиентные методы, методы редукции дисперсии \cite{Gas_Pos18, Lan2019}. Причём ускорение можно получить из базовых неускоренных вариантов этих методов с помощью 
УМ (см. п.~\ref{3.2}). По сравнению с оболочкой, использованной в \cite{Ivanova2020}, УМ даёт оценку сложности на логарифмический множитель лучше.
Это следует из теоретического анализа и было подтверждено в экспериментах \cite{Gasnikov2021}. 

\subsection{Ускоренные методы для седловых задач} \label{3.4} 

Следуя, например \cite{Alkousa2019, Lin2020}, рассмотрим  выпукло-вогнутую седловую задачу:
\begin{equation}\label{eq:3F}
\min _{x \in  \R^{d_x}}\, \left\{ F(x):=f(x)+\underbrace{\max _{y \in \R^{d_y}}\{G(x, y)-h(y)\}}_{g(x)=G(x, y^*(x))-h(y^*(x))} \;\; \right\},
\end{equation}
где $y^*(x) = \arg\max_{y \in \R^{d_y}}\{G(x, y)-h(y)\}$. Будем считать, что $\nabla f, \nabla G, \nabla h$ являются соответственно $L_f,L_G,L_h$-липшицевыми. 
Также будем считать, что $f(x) + G(x,y)$  является $\mu_x$-сильно (2-равномерно) выпуклой по $x$, а $G(x, y)-h(y)$
-- $\mu_y$-сильно (2-равномерно) вогнутой по $y$. Тогда $F(x)$ будет $\mu_x$-сильно выпуклой, а $\nabla g$ будет $L_g = \left(L_{G} + 2L_G^2/\mu_y\right)$-липшицевым \cite{Alkousa2019, Lin2020}. 

Если считать, что доступен $\nabla g$, то внешнюю задачу \eqref{eq:3F} можно решать ускоренным слайдингом (например, в варианте УМ с $p = 1$, см. п.~\ref{3.3})  за 
$\tilde{O}\left(\sqrt{L_f/\mu_x}\right)$ вычислений $\nabla f$ и $\tilde{O}\left(\sqrt{L_g/\mu_x}\right)$ вычислений $\nabla g$.  

Чтобы приближённо посчитать $\nabla g(x) = \nabla_x G(x,y^*(x))$, надо решить (с достаточной точностью) вспомогательную задачу в \eqref{eq:3F}, т.е. найти с нужной точностью $y^*(x)$. Это в свою очередь также можно сделать с помощью слайдинга (УМ с $p = 1$) за 
$\tilde{O}\left(\sqrt{L_h/\mu_y}\right)$ вычислений $\nabla h$ и $\tilde{O}\left(\sqrt{L_G/\mu_y}\right)$ вычислений $\nabla_y G$.

Резюмируя написанное, получаем, что исходную задачу \eqref{eq:3F} можно решить за $\tilde{O}\left(\sqrt{L_f/\mu_x}\right)$ вычислений $\nabla f$, $\tilde{O}\left(\sqrt{L_g/\mu_x}\right) \simeq \tilde{O}\left(\sqrt{L_G^2/(\mu_x\mu_y)}\right)$ вычислений $\nabla_x G$, $\tilde{O}\left(\sqrt{L_G^3/(\mu_x\mu_y^2)}\right)$ вычислений $\nabla_y G$, $\tilde{O}\left(\sqrt{L_h L_G^2/(\mu_x\mu_y^2)}\right)$ вычислений $\nabla h$. Поменяв порядок взятия $\min$ и $\max$ аналогичным образом, можно прийти к оценкам $\tilde{O}\left(\sqrt{L_h/\mu_y}\right)$ вычислений $\nabla h$, $\tilde{O}\left(\sqrt{L_G^2/(\mu_x\mu_y)}\right)$ вычислений $\nabla_y G$, $\tilde{O}\left(\sqrt{L_G^3/(\mu_x^2\mu_y)}\right)$ вычислений $\nabla_x G$, $\tilde{O}\left(\sqrt{L_f L_G^2/(\mu_x^2\mu_y)}\right)$ вычислений $\nabla f$.

Оценки, полученные на число вычислений $\nabla_x G$ и  $\nabla f$, в последнем случае не являются оптимальными \cite{Lin2020}. Чтобы улучшить данные оценки (сделать их оптимальными с точностью до логарифмических множителей \cite{Lin2020}), воспользуемся Каталистом, см. п.~\ref{3.2} (УМ, с $p = 1$, $H \gg \mu_x$, $f \equiv 0$, $g = F$, где $F$ определяется \eqref{eq:3F}). Если параметр метода $H$, то число итераций метода будет $\tilde{O}\left(\sqrt{H/\mu_x}\right)$, см. теорему~\ref{theoremRestartCATD}. На каждой итерации необходимо будет решать с должной точностью задачу вида \eqref{eq:3F}, в которой $L_f := L_f + H$, $\mu_x := \mu_x + H \simeq H$. Таким образом, для решения внутренней седловой задачи потребуется $\tilde{O}\left(\sqrt{L_h/\mu_y}\right)$ вычислений $\nabla h$, $\tilde{O}\left(\sqrt{L_G^2/(H\mu_y)}\right)$ вычислений $\nabla_y G$, $\tilde{O}\left(\sqrt{L_G^3/(H^2\mu_y)}\right)$ вычислений $\nabla_x G$, $\tilde{O}\left(\sqrt{(L_f + H) L_G^2/(H^2\mu_y)}\right)$ вычислений $\nabla f$. Считая для наглядности $L_f \ge L_G$, выберем $H = L_G$. Тогда итоговые оценки на число вычислений соответствующих градиентов будут такие: $\tilde{O}\left(\sqrt{L_h L_G/(\mu_x\mu_y)}\right)$ вычислений $\nabla h$, $\tilde{O}\left(\sqrt{L_G^2/(\mu_x\mu_y)}\right)$ вычислений $\nabla_y G$, $\tilde{O}\left(\sqrt{L_G^2/(\mu_x\mu_y)}\right)$ вычислений $\nabla_x G$, $\tilde{O}\left(\sqrt{L_f L_G/(\mu_x\mu_y)}\right)$ вычислений~$\nabla f$.

За счёт использования 
УМ приведённая выше схема улучшает похожую схему рассуждений из \cite{Lin2020} на логарифмический (по желаемой точности решения задачи) множитель и обобщает её на случай отличных от тождественного нуля функций $f$ и $h$, см. также \cite{wang2020,Yang2020}. Заметим также, что приведённая выше схема для седловых задач имеет обобщения и на невыпукло- (сильно) вогнутые седловые задачи \cite{Zhang}.

Отметим, что приведённый в этом разделе результат означает, что оптимальные градиентные методы для гладких выпукло-вогнутых седловых задач можно строить на базе оптимальных (ускоренных) градиентных методов для гладких задач выпуклой оптимизации. До недавнего времени в возможность такого построения было трудно поверить. Ещё более удивительным является (также недавно обнаруженная) возможность обратного перехода \cite{cohen2020}: на базе оптимального градиентного метода (Mirror Prox) для гладких выпукло-вогнутых седловых задач построить оптимальный (ускоренный) градиентный метод для задач гладкой выпуклой оптимизации. 

Приведённая здесь схема рассуждений наглядно демонстрирует, как из одной универсальной схемы ускорения удаётся получить (<<собрать как в конструкторе>>) 
оптимальный метод с точностью до логарифмического по желаемой точности множителя (при $f \equiv 0$ и $h \equiv 0$ -- только при этих условиях известны нижние оценки \cite{Lin2020}).

%
%
\section{Методы внутренней точки}\label{ch2_subsectIPM}

В этом разделе мы рассмотрим класс эффективных методов, которые применяются для решения конических задач (см.~п.~\ref{ch1_sect_cone_prog}). Методы внутренней точки отличаются тем, что генерируют последовательность итераций во \emph{внутренности} выпуклого конуса, лежащего в основе задачи. Это отличает их, например, от симплекс-метода и его аналогов, генерирующих последовательность экстремальных точек допустимого множества.

Применимость метода внутренней точки к данной конической задаче определяется наличием эффективно вычислимого т.н. самосогласованного барьера для соответствующего конуса. С помощью такого барьера задачу можно решить методом \emph{короткого шага}, имеющего полиномиальную сложность. В некоторых случаях, в частности, если конус симметричный, барьер обладает дополнительной структурой. Это позволяет существенно ускорить метод на практике, переходя к \emph{длинному шагу}. Стоит отметить, что в теории сложность многих методов с длинным шагом выше сложности методов с коротким шагом, что указывает либо на изъяны в самой теории, либо на существование не встречающихся на практике патологических ситуаций, в которых методы с длинным шагом дают сбой.

Методы внутренней точки можно также разделить на \emph{прямые} и \emph{прямо-двойственные}, в зависимости от того, генерируются ли точки только в прямом конусе или парами в прямом и двойственном. Если в большинстве методов все итерации удовлетворяют линейным ограничениям задачи, то в классе \emph{недопустимых} (infeasible) методов соответствующие невязки ненулевые и стремятся к нулю экспоненциально в ходе прогресса метода. Это избавляет пользователя от необходимости находить начальную точку для основной фазы метода.

Ниже мы рассмотрим базовый механизм методов внутренней точки, начиная с определения самосогласованных функций, и опишем особенности вышеперечисленных вариантов метода.

\subsection{Самосогласованные функции}

Самосогласованные функции были введены Ю.\,Е.~Нестеровым и А.\,С.~Немировским при изучении поведения метода Ньютона. Связывающим звеном здесь является аффинная инвариантность. Последовательность итераций метода Ньютона на данной функции переходит в себя при аффинном преобразовании координат, т.е. метод обладает аффинной инвариантностью. Поэтому естественно изучать поведение метода на классе функций, также являющимся аффинно инвариантным. Самосогласованные функции естественным образом возникают как аффинно инвариантный аналог функций с липшицевым гессианом.

В этом разделе мы рассмотрим поведение метода Ньютона на классе самосогласованных функций. Материал почерпнут из книги Ю.\,Е.~Нестерова и А.\,С.~Немировского \cite{NesNem94}.

Рассмотрим задачу минимизации выпуклой функции $f: D \to \mathbb R$ класса $C^3$, определённой на открытом выпуклом множестве $D \subset A$, где $A$ --- некое аффинное пространство. Для простоты предположим, что гессиан $f''$ функции положительно определён всюду на $D$. Исходя из данной итерации $x_k \in D$, \emph{метод Ньютона} выдаёт следующую итерацию по формуле:
\[ x_{k+1} = x_k - (f''(x_k))^{-1}f'(x_k).
\]
В точке $x_{k+1}$ строго выпуклый полином Тейлора функции $f$ второго порядка вокруг точки $x_k$ достигает минимума,

\begin{align*}
x_{k+1} &= \arg\min_x q_k(x) = \\
&=\arg\min_x \left( f(x_k) + \langle f'(x_k), x-x_k \rangle + \frac12(x-x_k)^Tf''(x_k)(x-x_k) \right).
\end{align*}

При анализе шага Ньютона естественно использовать евклидову норму $||\cdot||_{x_k}$, задаваемую гессианом $f''(x_k)$ на касательном пространстве в точке $x_k$. В этой локальной норме подмножества уровня $\{ x \in A \,|\, q_k(x) \leq c \}$ квадратичной аппроксимации $q_k$ являются шарами с центром в точке минимума $x_{k+1}$. Длина шага в этой норме равна
\begin{equation} \label{gradient_norm}
\rho = \sqrt{(x_{k+1}-x_k)^Tf''(x_k)(x_{k+1}-x_k)} = \sqrt{f'(x_k)^T(f''(x_k))^{-1}f'(x_k)}.
\end{equation}
Это выражение называется \emph{ньютоновским декрементом}. Оно также задаёт длину градиента $f'(x_k)$ в локальной норме.

Очевидно, без дополнительных условий на функцию $f$ невозможно сделать какие-либо утверждения о поведении этой функции в окрестности новой точки $x_{k+1}$, в частности, насколько уменьшилось или уменьшилось ли вообще значение функции $f$ в новой точке по сравнению с предыдущей. Для этого необходимо, чтобы аппроксимация $q_k$ функции $f$, построенная в точке $x_k$, не слишком сильно отличалась от $f$ в окрестности новой точки $x_{k+1}$. Другими словами, локальные нормы $||\cdot||_{x_k},||\cdot||_{x_{k+1}}$ не должны слишком сильно отличаться друг от друга. Это требование соответствует некоему условию типа Липшица на гессиан $f''(x)$. С другой стороны, аффинная инвариантность предполагает, что сравнивать изменение гессиана в окрестности данной точки $x$ можно только с длиной вариации, измеренной в локальной норме $||\cdot||_x$. Это приводит нас к следующему определению.

\begin{defin} \label{def_self_concordant} Выпуклая функция $f: D \to \mathbb R$ класса $C^3$, определённая на выпуклой области $D \subset A$ некоего аффинного пространства называется \emph{самосогласованной}, если она удовлетворяет неравенству
\[ |f'''(x)[h,h,h]| \leq 2(f''(x)[h,h])^{3/2}
\]
для всех $x \in D$ и всех касательных векторов $h \in T_xD$.

Функция называется \emph{сильно} самосогласованной, если дополнительно выполняется условие
\[ \lim_{x \to \partial D}f(x) = +\infty.
\] 
\end{defin}

Напомним, что производные функции $f$ интерпретируются как мульти-линейные формы на касательном пространстве $T_xD$, т.е.
\[ f''(x)[h,h] = \sum_{i,j=1}^n \frac{\partial^2f(x)}{\partial x_i\partial x_j}h_ih_j,\ \  f'''(x)[h,h,h] = \sum_{i,j,k=1}^n \frac{\partial^3f(x)}{\partial x_i\partial x_j\partial x_k}h_ih_jh_k,
\]
где $n$ -- размерность пространства. Показатель $3/2$
на правой стороне неравенства необходим для того, чтобы неравенство было однородным по $h$.

Граничное условие означает, что $f$ стремится к $+\infty$ на любой последовательности точек в $D$, стремящейся к точке на границе $\partial D$. Это условие обеспечивает замкнутость подмножеств уровня $\{ x \in D \,|\, f(x) \leq c \}$ для любых констант $c \in \mathbb R$.

Самосогласованные функции обладают следующими свойствами. Пусть $F$ --- (сильно) самосогласованная функция на области $D$.
\begin{itemize}
\item Если $A$ --- аффинное подпространство, то ограничение $F$ на пересечение $A \cap D$ является (сильно) самосогласованным на этом пересечении.
\item Если $l$ --- линейная функция на $D$, то $F + l$ является (сильно) самосогласованной.
\item Если $\tilde F$ --- (сильно) самосогласованная функция на области $\tilde D$, то $G(x,y) = F(x) + \tilde F(y)$ --- (сильно) самосогласованная функция на прямом произведении $D \times \tilde D$.
\item Для всех $\alpha \geq 1$ произведение $\alpha F$ является (сильно) самосогласованным.
\end{itemize}

\begin{remark} Вообще говоря, условия определения~\ref{def_self_concordant} можно несколько ослабить. Достаточно потребовать, чтобы функция $f$ была класса $C^2$ с положительно определённым гессианом, такая что $\lim_{x \to \partial D}f(x) = +\infty$, открытый эллипсоид Дикина радиуса 1 вокруг произвольной точки $x \in D$ является подмножеством $D$, и что для всех $x,y \in D$, удовлетворяющих $||y-x||_x < 1$, и любого ненулевого вектора $u$ выполнялось неравенство
\[ \frac{||u||_y}{||u||_x} \leq \frac{1}{1 - ||y-x||_x}.
\]
Более подробно разные определения самосогласованности и отношения между ними разобраны в \cite[п.~2.5]{Renegar01}.
\end{remark}

\bigskip

Следующие результаты гарантируют, что метод Ньютона в применении к сильно самосогласованной функции $f$ может безопасно делать шаги конечной длины, где <<безопасно>> означает, что все итерации являются элементами $D$ и значения декремента $\rho$ на последовательности итераций монотонно убывают.

{\lemma Пусть $f: D \to \mathbb R$ -- сильно самосогласованная функция. Тогда для любого $x^* \in D$ центрированный на $x^*$ {\itshape {\bfseries эллипсоид Дикина}}, задающийся
\[ E_{x^*} = \{ x \,|\, (x-x^*)^Tf''(x^*)(x-x^*) < 1 \},
\]
содержится в области $D$. }

Из этого следует, что шаг Ньютона не выводит из области $D$, если градиент $f'(x_k)$ в текущей точке имеет локальную норму, не превышающую единицу. Более того, в этом случае ньютоновский декремент в точке $x_{k+1}$ можно ограничить функцией от декремента в точке $x_k$:
\begin{equation} \label{decrement_decrease} 
\rho(x_{k+1}) \leq \left( \frac{\rho(x_k)}{1+\rho(x_k)} \right)^2.
\end{equation}
В частности, справедлива импликация
\begin{equation} \label{values_iterate}
\rho(x_k) \leq \frac14 \quad \Rightarrow\quad \rho(x_{k+1}) \leq \frac19.
\end{equation}
Следует отметить, что оценка \eqref{decrement_decrease} не оптимальна, и численные значения в \eqref{values_iterate}, которые будут использоваться ниже, можно существенно увеличить.

Оценку можно также улучшить, если делать не полный, а укороченный шаг Ньютона, задающийся формулой
\[ x_{k+1} = x_k - \gamma_k(f''(x_k))^{-1}f'(x_k).
\]
Здесь $\gamma_k \in (0,1]$ --- коэффициент затухания (damping coefficient). Имеем следующую оценку \cite[Теорема 5.2.2.3]{NesterovLectures}.

\begin{lemma}
Пусть $F$ --- сильно самосогласованная функция на $D$. Если $\rho(x_k) < \lambda^* = \mbox{roots}(\lambda^3+\lambda^2+\lambda-1) \approx 0.5437$, то после укороченного шага Ньютона с коэффициентом $\gamma_k = \frac{1+\rho(x_k)}{1+\rho(x_k)+\rho(x_k)^2}$ имеем
\[ \rho(x_{k+1}) \leq \rho(x_k)^2\left( 1 + \rho(x_k) + \frac{\rho(x_k)}{1+\rho(x_k)+\rho(x_k)^2} \right) < \rho(x_k).
\]
\end{lemma}

\subsubsection{Геометрическая интерпретация} \label{ch1_geom_Newton}

Дадим геометрическую интерпретацию метода Ньютона. Точка минимума $x^*$ функции $f$ характеризуется условием оптимальности первого порядка $f'(x^*) = 0$. Это позволяет нам игнорировать значения $f$ и ограничиться поиском корня этого векторного уравнения.

Пусть $V$~--- векторное пространство, лежащее в основе аффинного пространства $A$. Тогда градиент $f'(x_k)$ является элементом двойственного векторного пространства $V^*$. Граф $M$ градиентного отображения $\nabla f: x \mapsto f'(x)$ тогда является $n$-мерным подмногообразием $2n$-мерного произведения $A \times V^*$. Точка минимума $x^*$ соответствует паре $(x^*, \, 0) \in A \times V^*$, которая одновременно является пересечением многообразия $M$ с горизонтальным подпространством $H_0 = \{ (x,0) \,|\, x \in A \}$.

\begin{figure}[ht]
\begin{center}
\includegraphics[width=10.78cm,height=6.32cm]{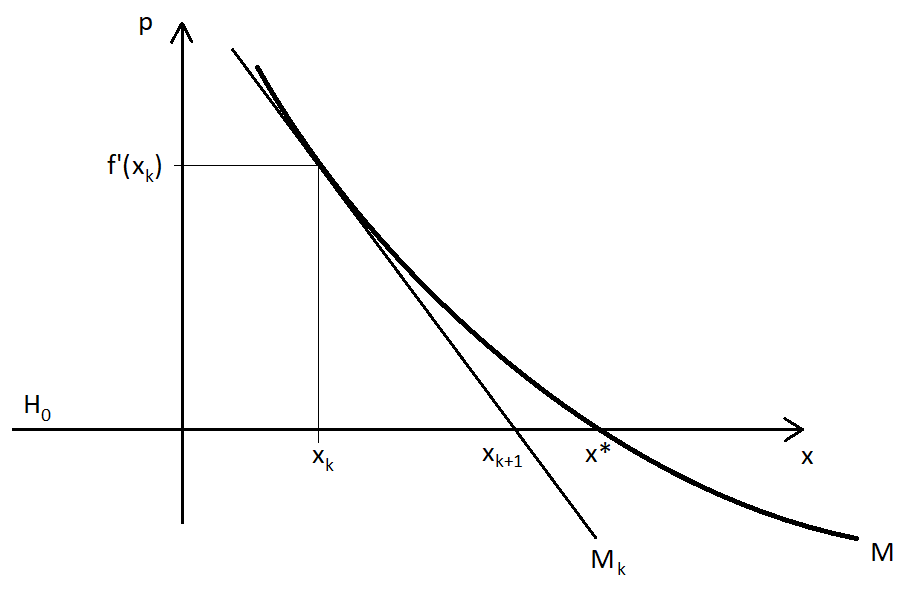}
\end{center}
\caption{Геометрическая интерпретация метода Ньютона}
\label{fig:newton_geometric}
\end{figure}

Тогда итерация метода Ньютона может быть интерпретирована следующим образом. Напомним, что в текущей точке $x_k$ мы строим квадратичную аппроксимацию $q_k$ функции $f$. Граф $M_k \subset A \times V^*$ градиентного отображения $\nabla q_k$ также является $n$-мерным подмногообразием. Более того, градиент $\nabla q_k$ \emph{линейный}, и поэтому $M_k$ является даже аффинным подпространством произведения $A \times V^*$. Так как градиент и гессиан $q_k$ в точке $x_k$ совпадают с градиентом и гессианом $f$, соответственно, это подпространство не что иное, как \emph{касательная плоскость} к $M$ в точке $(x_k,f'(x_k))$. Следующая точка $x_{k+1}$ определяется пересечением $(x_{k+1},0)$ подпространств $M_k$ и $H_0$ (см.~рис.~\ref{fig:newton_geometric}).

\subsection{Прямой метод внутренней точки с коротким шагом}

В этом разделе мы рассмотрим принцип методов внутренней точки на примере самой простой их версии, прямом методе следования центральному пути с коротким шагом.

Вся сложность конической программы заключена в нелинейном ограничении $x \in K$. Основная идея методов внутренней точки заключается в устранении нелинейного условия включения посредством добавления сильно самосогласованной \emph{барьерной функции} $F: K^o \to \mathbb R$ к линейной функции цены. Для удобства положим $F(x) = +\infty$ для всех $x \not\in \interior K$. Напомним, что сильная самосогласованность остаётся в силе, если прибавить к $F$ произвольную линейную функцию.

Пусть $X \subset \mathbb R^n$ -- произвольное замкнутое выпуклое множество с непустой внутренностью $D = X^o$, не содержащее прямой, а $F: X^o \to \mathbb R$ -- сильно самосогласованная функция на $X$. В случае конической задачи роль множества $X$ выполняет допустимое множество, являющееся пересечением конуса $K$ с аффинным подпространством, задающимся линейными условиями. Как правило, в этом случае сильно самосогласованная функция задана на внутренности конуса $K$, а $F$ определяется как её ограничение на область $D$. Поэтому в данной ситуации необходимо потребовать, чтобы существовала допустимая точка во внутренности конуса $K$.

Наша цель --- решить выпуклую задачу оптимизации
\[ \inf_{x \in X} c^Tx,
\]
существование решения $x^* \in \partial X$ которой мы предполагаем. Вместо исходной задачи мы рассмотрим семейство задач
\begin{equation} \label{composite_function}
\min_x\, (\tau \cdot c^Tx + F(x)),
\end{equation}
параметризованное вещественным параметром $\tau \geq 0$. В этих вспомогательных задачах нелинейное условие $x \in X$ отсутствует, но зато мы вместо линейной функции цены минимизируем сильно самосогласованные функции. В предыдущем разделе мы увидели, что для этого можно использовать метод Ньютона. 

Для достаточно больших $\tau$ точки минимума $x^*(\tau)$ вспомогательной задачи существуют и единственны. Более того, они являются дифференцируемыми по параметру $\tau$ и образуют кривую, называемую \emph{центральным путём}.  При $\tau \to +\infty$ центральный путь стремится к решению $x^* = x^*(+\infty)$ исходной задачи. Если у исходной задачи решение не единственно, то пределом центрального пути будет служить точка в относительной внутренности множества решений.

Опишем схему прямого метода следования центральному пути. Для краткости обозначим композитную функцию цены во вспомогательной задаче \eqref{composite_function} через $F_{\tau}$, а ньютоновский декремент этой функции --- через $\rho_{\tau}$. Метод генерирует последовательность точек $x_k \in D$ в окрестности центрального пути, которая одновременно перемещается вдоль этого пути (следует ему). При этом все итерации удовлетворяют линейным ограничениям задачи. Окрестность определена таким образом, что для каждой итерации $x_k$ существует значение $\tau_{k-1}$ такое, что декремент функции $F_{\tau_{k-1}}$ в точке $x_k$ удовлетворяет условию $\rho_{\tau_{k-1}}(x_k) \leq 1/9$. Таким образом, центральный путь имеет непустое пересечение с эллипсоидом Дикина радиуса $1/9$ вокруг текущей итерации.

Следующая точка $x_{k+1}$ получается из текущей $x_k$ посредством шага Ньютона по направлению к некоторой точке $x^*(\tau_k)$ на центральном пути. Другими словами, мы делаем один шаг Ньютона для минимизации вспомогательной функции $F_{\tau_k}$. Параметр $\tau_k$ при этом, с одной стороны, нужно выбрать как можно более большим, чтобы продвинуться вдоль центрального пути по направлению к решению, а с другой стороны, нужно оставаться в окрестности центрального пути, т.е. обеспечить выполнение неравенства $\rho_{\tau_k}(x_{k+1}) \leq 1/9$. В силу \eqref{values_iterate} это условие будет выполнено, если $\rho_{\tau_k}(x_k) \leq 1/4$.

\begin{figure}[ht]
\begin{center}
\includegraphics[width=9.19cm,height=2.83cm]{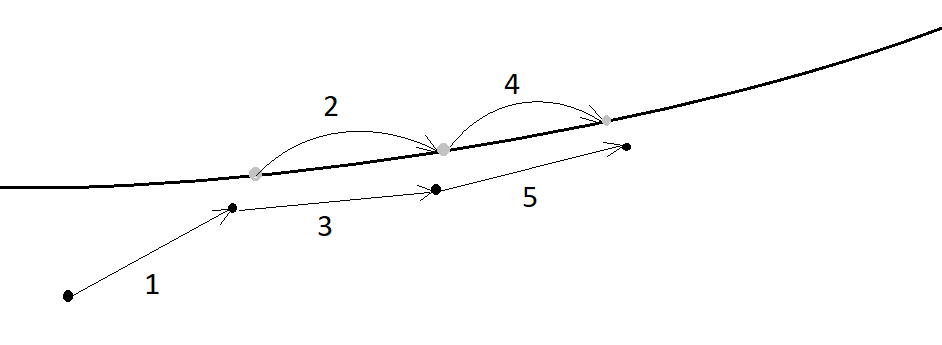}
\end{center}
\caption{Схема прямого метода следования центральному пути: обновление целевой точки (минимума вспомогательной функции цены) на центральном пути перемежается с шагом Ньютона по направлению к текущей целевой точке}
\label{fig:short_step}
\end{figure}

Таким образом, метод попеременно делает шаг Ньютона по направлению к некоторой целевой точке на центральном пути и отодвигает эту точку в сторону решения задачи (см.~рис.~\ref{fig:short_step}). Шаг Ньютона приближает текущую точку к минимуму и \emph{уменьшает} декремент. Обновление параметра $\tau$ удаляет минимум от текущей точки и \emph{увеличивает} декремент. При этом позволительно увеличивать $\tau$ настолько, насколько результирующий проигрыш в декременте компенсируется выигрышем, сделанным на предыдущем шаге Ньютона.

Оценим прогресс в параметре $\tau$, сделанный в ходе итерации. Рассмотрим зависимость величины 
\[ \rho_{\tau}(x_k) = \sqrt{(F'_{\tau}(x_k))^T(F_{\tau}''(x_k))^{-1}F'_{\tau}(x_k)}
\]
от $\tau$. Имеем $F_{\tau}'' = F''$, $F'_{\tau} = F' + \tau c$, и под корнем стоит неотрицательная, квадратичная по $\tau$ функция с лидирующим коэффициентом $c^T(F''(x_k))^{-1}c$. Так как $\rho_{\tau_{k-1}}(x_k) \leq 1/9$, то уравнение $\rho_{\tau}(x_k) = 1/4$ имеет два вещественных корня. Очевидно, \nk{следующее} значение $\tau_k$ параметра $\tau$ \nk{следует} положить равным большему из этих корней.

Граф $\rho_{\tau}(x_k)$ как функции от $\tau$ --- гипербола с асимптотическим наклоном, равным $\sqrt{c^T(F''(x_k))^{-1}c}$. Поэтому $\left|\frac{d\rho_{\tau}(x_k)}{d\tau}\right| \leq \sqrt{c^T(F''(x_k))^{-1}c}$ для~всех $\tau$, откуда получаем 
\[ \tau_k - \tau_{k-1} \geq \frac{\rho_{\tau_k}(x_k) - \rho_{\tau_{k-1}}(x_k)}{\sqrt{c^T(F''(x_k))^{-1}c}} \geq \frac{\frac14 - \frac19}{\sqrt{c^T(F''(x_k))^{-1}c}}.
\]
Оценим величину $||c||_{x_k} = \sqrt{c^T(F''(x_k))^{-1}c}$. Обозначая $\delta_k = F'(x_k) + \tau_kc$ и учитывая $||\delta_k||_{x_k} = \rho_{\tau_k}(x_k) \leq 1/4$, получим 
\[ ||c||_{x_k} = \tau^{-1}_k||-F'(x_k) + \delta_k||_{x_k} \leq \tau^{-1}_k\left(||F'(x_k)||_{x_k} + \frac14\right).
\]
В итоге находим оценку
\[ \log\tau_k - \log\tau_{k-1} \geq \frac{\tau_k - \tau_{k-1}}{\tau_k} \geq \frac{\frac{5}{36}}{||F'(x_k)||_{x_k} + \frac14}
\]
и получаем, что на каждом шаге $\log\tau$ увеличивается на величину порядка $\frac{1}{||F'(x_k)||_{x_k}}$. Чтобы была гарантирована линейная скорость сходимости, необходимо ограничить локальную норму $||F'(x)||_x$ сверху. Это мотивирует следующее определение.

\begin{defin} Сильно самосогласованная функция $F: D \to \mathbb R$, определённая на некоторой выпуклой области $D$, называется \emph{самосогласованным барьером} на $D$ с параметром $\nu$, если для всех $x \in D$ выполнено соотношение
\[ ||F'(x)||_x^2 \leq \nu.
\] 
\end{defin}

Таким образом, если функция $F$ является самосогласованным барьером с параметром $\nu$, то вышеописанный прямой метод следования центральному пути позволяет на каждом шаге увеличивать $\log\tau$ на величину порядка $\nu^{-1/2}$ или, что эквивалентно, умножать $\tau$ на константу $\theta = 1 + O(\nu^{-1/2})$. Чем больше параметр барьера $\nu$, тем меньше множитель $\theta$ и тем медленнее метод будет сходиться. Поэтому на каждом данном выпуклом множестве желательно иметь в наличии барьеры с как можно меньшим параметром. Сложность задачи зависит от наличия эффективно вычислимого самосогласованного барьера с небольшим значением параметра.

\medskip

В коническом программировании барьер определяется на конусе $K$, допускающем действие группы растяжений. Разумно потребовать также инвариантность барьера на $K$ по отношению к этой группе симметрий. Постоянство функции $F$ на внутренних лучах конуса приводит к противоречию с остальными условиями. Однако в методе Ньютона используются только производные минимизируемой функции. Поэтому достаточно потребовать инвариантность производных $F$ по отношению к растяжениям. Это естественным образом приводит к условию логарифмичной однородности. Следующее определение было введено Ю.\,Е.~Нестеровым и А.\,С.~Немировским в~\cite{NesNem94}.

\begin{defin}
Пусть $K \subset \mathbb R^n$ -- выпуклый регулярный конус. Выпуклая функция $F: K^o \to \mathbb R$ класса $C^3$ с положительно определённым гессианом называется \emph{логарифмично однородным самосогласованным барьером} с \emph{параметром} $\nu$, если $F$ является сильно самосогласованной функцией на $K^o$ и $F$ логарифмично однородна порядка $-\nu$, т.е.
\begin{equation} \label{loghom}
F(\alpha x) = -\nu\log\alpha + F(x)
\end{equation}
для всех $x \in K^o$ и $\alpha > 0$. 
\end{defin}

Таким образом, растяжения действуют прибавлением констант к логарифмично однородной функции $F$.

\medskip

Данное для случая конусов определение самосогласованного барьера совместимо с приведённым выше общим определением, поскольку логарифмичная однородность также ограничивает ньютоновский декремент, и локальная норма градиента барьера с параметром $\nu$ равна константе $\sqrt{\nu}$. Действительно, дифференцируя соотношение \eqref{loghom} по $x$, мы получим $\alpha F'(\alpha x) = F'(x)$. Дифференцируя это соотношение и \eqref{loghom} по $\alpha$ при $\alpha = 1$, мы получим
\begin{equation} \label{gradient_rels}
F'(x) + F''(x)\cdot x = 0,\qquad \langle F'(x),x \rangle = -\nu.
\end{equation}
Из этого следует, что $(F''(x))^{-1}F'(x) = -x$ и, следовательно, 
\[ (F'(x))^T(F''(x))^{-1}F'(x) = -\langle F'(x),x \rangle = \nu.
\]


\medskip

Подытожим результаты этого раздела. Для решения конической программы с помощью метода внутренней точки необходимо наличие эффективно вычислимого логарифмично однородного самосогласованного барьера с небольшим значением параметра $\nu$ на соответствующем выпуклом конусе. В методах следования центральному пути нелинейное коническое ограничение устраняется переходом к однопараметрическому семейству вспомогательных задач, каждая из которых состоит в минимизации самосогласованной функции. Метод короткого шага чередует итерацию Ньютона для решения вспомогательной задачи с умножением параметра семейства на величину $1 + O(\nu^{-1/2})$, при этом последовательность итераций сходится с линейной скоростью к решению исходной задачи. Скорость сходимости увеличивается с уменьшением параметра барьера $\nu$. В методах с коротким шагом итерации $x_k$ находятся в окрестности центрального пути. Точнее, центральный путь пересекается с эллипсоидами Дикина некоторого порогового радиуса $\rho < 1$ вокруг точек $x_k$.

\subsection{Двойственные барьеры}

В п.~\ref{ch1_conic_duality} мы видели, что к каждой конической программе можно определить двойственную программу над двойственным конусом. Двойственность также распространяется на самосогласованные барьеры.

\medskip

Для любого логарифмично однородного самосогласованного барьера $F: \interior K \to \mathbb R$ на регулярном выпуклом конусе $K$ с параметром $\nu$ можно определить \emph{двойственный} барьер $F_*: \interior K^* \to \mathbb R$ на двойственном к $K$ конусе посредством преобразования Лежандра:
\[ F_*(s) = \sup_{x \in K}\, (-\langle s,x \rangle - F(x)).
\]
Супремум для данного $s \in \interior K^*$ достигается в точке $x \in K^o$, удовлетворяющей уравнению $F'(x) = -s$. В силу строгой выпуклости $F$ и свойства $F|_{\partial K} = +\infty$ такая точка $x$ существует и единственна для любого $s \in \interior K^*$. Действительно, рассмотрим аффинную гиперплоскость $A = \{ x \in \mathbb R^n \,|\, \langle x,s \rangle = \nu \}$. Так как $s$ находится во внутренности двойственного конуса $K^*$, пересечение $X = A \cap K$ является компактным. На внутренности $X$ функция $F$ строго выпукла, а на границе имеем $F|_{\partial X} = +\infty$. Поэтому минимум функции $F$ на $X$ существует и единственен. Обозначим этот минимум через $x^*$. Условие оптимальности первого порядка в точке $x^*$ сводится к существованию $\lambda \in \mathbb R$ такого, что $F'(x^*) = \lambda s$. В силу \eqref{gradient_rels} из этого следует, что
\[ -\nu = \langle F'(x^*),x^* \rangle = \lambda \langle s,x^* \rangle = \lambda\nu.
\]
Поэтому $\lambda = -1$ и $F'(x^*) = -s$, \nk{что} и требовалось доказать.

\begin{teo} \label{thm:isometry} Функция $F_*$ является логарифмично однородным самосогласованным барьером с параметром $\nu$ на $K^*$ {\rm \cite[Теорема 2.4.4]{NesNem94}}. Отображение ${\cal L}: x \mapsto -F'(x)$ является биекцией между внутренностями конусов $K$ и $K^*$. Если рассмотреть эти внутренности как римановы многообразия, оснащённые метриками $F''(x)$ и $F_*''(s)$ соответственно, то отображение ${\cal L}$ является изометрией {\rm (\cite[стр.~45]{NesNem94},} см.~также {\rm \cite{NesTodd02}).} Более того, отображение ${\cal L}$ переводит тензор $F'''(x)$ третьих производных в точке $x$ в тензор $-F_*'''(s)$ в образе $s = -F'(x)$ {\rm \cite[стр.~45]{NesNem94}.}
\end{teo}

\begin{remark} По определению преобразование Лежандра $f^*$ функции $f$ задаётся формулой $f^*(p) = \sup_{x \in K}\, (\langle p,x \rangle - f(x))$. Знак минус при аргументе двойственного барьера был введён для совместимости с определением двойственного конуса, $F_*(p) = F^*(-p)$. 
\end{remark}

\subsection{Примеры самосогласованных барьеров}

Приведём явные выражения для самосогласованных барьеров на конусах, лежащих в основе ходовых классов конических программ, приведённых выше.

Барьером на прямом произведении $K = \prod_{i=1}^m K_i$ конусов, оснащённых барьерами $F_i$ с параметрами $\nu_i$, может служить сумма
\[ F(x_1, \, \dots, \, x_m) = \sum_{i=1}^m F_i(x_i).
\]
Она обладает параметром $\nu = \sum_{i=1}^m \nu_i$.

Рассмотрим симметричные конусы, лежащие в основе классов программ LP, SOCP и SDP.

\bigskip

\begin{center}
\begin{tabular}{c||c|c|c}
конус & $\mathbb R_+^n$ & $L^n$ & $\mathbb{S}_+^n$ \\ 
\hline
\hline
барьер & $-\sum_{j=1}^n \log x_j$ & $-\log\left(x_n^2 - x_1^2 - \dots - x_{n-1}^2\right)$ & $-\log\det A$ \\ 
\hline
параметр $\nu$ & $n$ & 2 & $n$ \\ 
\end{tabular}
\end{center}

\bigskip

Для барьеров, приведённых в вышестоящей таблице, двойственный барьер идентичен исходному. Параметр этих барьеров оптимальный, т.е. на этих конусах нет барьера со строго меньшим параметром. Но они обладают ещё одним полезным свойством, они являются \emph{авто-шкалированными}. Это свойство допускает особенно эффективные методы внутренней точки для решения конических программ над этими конусами, а именно методы с длинным шагом \cite{NesterovTodd97}, \cite{NesterovTodd98}.

\medskip

\emph{Экспоненциальный конус.} Рассмотрим несимметричный конус
\begin{equation} \label{definition_Kexp} 
K_{\exp} = \left\{ (x,y,0) \,|\, x \leq 0,\ y \geq 0 \right\} \cup \left\{ (x,y,z) \,|\, z > 0,\ y \geq ze^{x/z} \right\},
\end{equation}
встречающийся в геометрическом программировании. На этим конусе существует барьер
\[ F(x,y,z) = -\log\left(z\log\frac{y}{z}-x\right)-\log y - \log z,
\]
значение параметра которого равно $\nu = 3$. Это значение параметра также оптимально.

\subsection{Универсальные конструкции барьеров}

Наличие эффективно вычислимого логарифмично однородного самомогласованного барьера на выпуклом конусе позволяет решать конические программы над этим конусом. Естественно, возникает вопрос о существовании барьеров на произвольных конусах. Ответ на этот вопрос утвердительный, и в этом разделе мы приведём два способа построения таких барьеров.

\medskip

{\bf Универсальный барьер.} Пусть $K \subset \mathbb R^n$ -- регулярный выпуклый конус. Рассмотрим его \emph{характеристическую функцию}
\[ \varphi(x) = \int_{K^*}e^{-\langle x,y \rangle}\,dy,
\]
определённую для всех внутренних точек $x \in K$. Функция $F(x) = \log\varphi(x)$ является логарифмично однородным самосогласованным барьером на $K$ со значением параметра $\nu = n$. Этот барьер называется \emph{универсальным}. Он был введён в \cite{NesNem94} Ю.\,Е. Нестеровым и А.С. Немировским и является первой универсальной конструкцией барьера. Изначально было известно, что параметр универсального барьера ограничен сверху величиной $C \cdot n$, где $C \geq 1$ -- некоторая постоянная. Позже в работе \cite{BubeckEldan19} было установлено, что можно выбрать $C = 1$.

Если $A$ -- линейный автоморфизм конуса $K$, то имеет место соотношение $F(Ax) = F(x) + \log\det A$ для всех $x \in K^o$. Это позволяет вычислить универсальный барьер на однородных конусах, т.е. на внутренности которых группа линейных автоморфизмов действует транзитивно \cite{Guler96,GulerTuncel98}. 

В общем случае универсальный барьер задаётся многомерным интегралом по двойственному конусу и его трудно вычислить даже для сравнительно простых неоднородных конусов.

\medskip

{\bf Энтропический барьер.} Этот барьер двойственный к универсальному. Его параметр также равен размерности конуса. На ограниченных выпуклых множествах энтропический барьер был рассмотрен в работе \cite{BubeckEldan19}.

\medskip

{\bf Канонический барьер.} Построение этого барьера основано на глубокой теореме в теории уравнений в частных производных.

\begin{teo} Пусть $D \subset \mathbb R^n$ -- выпуклая область, не содержащая прямой. Тогда существует единственное выпуклое решение $F: D \to \mathbb R$ уравнения в частных производных $\log\det F'' = 2F$ с граничным условием $\lim_{x \to \partial D}F(x) = +\infty$. 
\end{teo}

В случае когда $D$ -- внутренность выпуклого регулярного конуса $K$, это решение является логарифмично однородным самосогласованным барьером на $K$ со значением параметра $\nu = n$ \cite{Hildebrand14d}. Этот барьер называется \emph{каноническим}. Двойственный барьер к каноническому барьеру на конусе $K$ совпадает с каноническим барьером на двойственном конусе $K^*$.

Канонический барьер обладает тем же свойством инвариантности по отношению к действию группы автоморфизмов конуса $K$, что и универсальный, и поэтому совпадает с ним на однородных конусах.

Для ортанта, конуса Лоренца и матричного конуса все три конструкции приводят к барьерам, пропорциональным авто-шкалированным барьерам, приведённым в вышестоящей таблице.

Примером неоднородного конуса, для которого канонический барьер вычислим аналитически, может послужить экспоненциальный конус \eqref{definition_Kexp}. На этом конусе канонический барьер задаётся формулой
\[ F(x,y,z) = -\log y -2\log z+\phi\left(\log\frac{y}{z}-\frac{x}{z}\right),
\]
где граф функции $\phi: \mathbb R_{++} \to \mathbb R$ определяется кривой
\[ \left\{ \left. \begin{pmatrix} t \\ \phi \end{pmatrix} = \frac12\begin{pmatrix} \log(1+\kappa)+2\kappa \\ \log(1+\kappa)-3\log\kappa \end{pmatrix} \,\right|\, \kappa \in \mathbb R_{++} \right\}.
\]

\subsection{Другие варианты методов внутренней точки}

Как мы видели выше, метод следования центральному пути с коротким шагом генерирует последовательность точек, расположенных в некоторой узкой окрестности центрального пути. При этом параметр $\tau$, определяющий прогресс вдоль центрального пути, на каждом шаге умножается на величину $1+O\left(1/\sqrt{\nu}\right)$, где $\nu$ --- параметр барьера. Поэтому заданная точность $\epsilon$ достигается за $O(\sqrt{\nu}\log\epsilon)$ итераций. На данный момент это самая лучшая скорость сходимости, достигнутая в теории. Однако практика показывает, что скорость можно значительно повысить, отступив от этих правил.

\medskip

Методы с \emph{длинным шагом} увеличивают параметр $\tau$ более агрессивно. В этом случае шаг метода Ньютона, как правило, выводит текущую точку из узкой окрестности центрального пути. Поэтому нужно либо расширить эту окрестность, либо сделать несколько шагов, чтобы вернуться в неё. В первом случае необходимо наличие дополнительной структуры, позволяющей измерять расстояние до центрального пути другим способом, чем ньютоновский декремент. Мы рассмотрим соответствующие алгоритмы, работающие на симметричных конусах, чуть ниже. Во втором случае можно показать, что для возвращения в узкую окрестность центрального пути достаточно сделать порядка $O(\nu)$ шагов метода Ньютона. Таким образом, в теории сложность возрастает до $O(\nu\log\epsilon)$ итераций.

\medskip

\emph{Прямо-двойственные методы} решают одновременно и прямую, и двойственную коническую программу. Для каждой из них генерируется последовательность точек $x_k \in \interior K$ и $s_k \in \interior K^*$, которые аппроксимируют точки $x^*(\tau_k)$, $s^*(\tau_k)$ на центральных путях прямой и двойственной программы соответственно. Эти точки можно объединить в пары $(x_k,s_k)$ в произведении прямого и двойственных пространств, которые находятся в окрестности точки $(x^*(\tau_k),s^*(\tau_k))$ прямо-двойственного центрального пути. Если конус обладает дополнительной структурой, в частности, если он симметричный, то существуют очень эффективные прямо-двойственные методы, которые мы рассмотрим в следующем разделе.

\medskip

\emph{Методы редукции потенциала} производят последовательность допустимых прямо-двойственных пар $(x_k,s_k)$ точек, не привязанных к каким-либо точкам на центральном пути. На этой последовательности монотонно убывают значения потенциала
\[ \Phi(x,s) = (\nu+\sqrt{\nu})\log\langle x,s \rangle + F(x) + F_*(s),
\]
определённого на произведении внутренностей прямого и двойственного конуса. Отметим, что на множестве пар точек $(x,s)$, которые удовлетворяют линейным ограничениям прямой и двойственной задачи, произведение $\langle x,s \rangle$ является \emph{линейным}. Это следует из того, что для любых двух таких пар $(x,s)$ и $(x',s')$ справедливо соотношение $\langle x - x',s - s' \rangle = 0$. Поэтому функция $\Phi$ <<почти>> выпукла, в смысле что она выпукла на любой аффинной гиперплоскости $H_{\mu} = \{ (x,s) \,|\, \langle x,s \rangle = \mu \}$ аффинного подпространства, задаваемого линейными ограничениями прямой и двойственной задач. На каждом шаге потенциал убывает как минимум на некоторую постоянную. Любая последовательность пар точек, на которой значение потенциала стремится к $-\infty$, обязана стремиться к решениям прямой и двойственной задач.

\subsection[Прямо-двойственные методы на симметричных\\ конусах]{Прямо-двойственные методы на симметричных\\ конусах}\label{SB}

Самыми успешными на практике являются алгоритмы решения конических задач над симметричными конусами, т.е. линейных, квадратично-коничных и полуопределённых программ. Причиной тому является богатая структура симметричных конусов, которая допускает существование барьеров с особыми свойствами, так называемых \emph{авто-шкалированных} барьеров. 

Определение авто-шкалированного барьера $F$ на симметричном конусе $K \in V$ довольно техническое. Главное его свойство --- это существование для любой прямо-двойственной пары точек $(x,s) \in \interior K \times \interior K^*$ единственной так называемой \emph{точки шкалировки} $w \in \interior K$, которая определяется как минимум функции $\Xi(z) = \langle F'(z),x \rangle - \langle z,s \rangle$ по $z$. Дополнительно она должна удовлетворять условию $F_*(s) = F(x) - 2F(w) - \nu$, где $F_*$ --- двойственный барьер. Ей также соответствует двойственная точка шкалировки $w_* = -F'(w)$, минимизирующая функцию $\langle F_*'(t),s \rangle - \langle x,t \rangle$ по $t$.

Покажем, что пара точек $(w,w_*)$ имеет простую геометрическую интерпретацию. На прямом произведении $V \times V^*$ прямого и двойственного векторного пространств существует каноническое псевдо-скалярное произведение, определяемое выражением
\[ \langle (u,v),(u',v') \rangle := \frac{\langle u,v' \rangle + \langle u',v \rangle}{2}.
\]
В правой части этого уравнения угловые скобки обозначают обычную свёртку прямого с двойственным вектором. Это произведение порождает псевдо-римановую метрику, в которой квадрат расстояния задаётся выражением
\[ d^2((u,v),(u',v')) = \langle (u-u',v-v'),(u-u',v-v') \rangle = \langle u-u',v-v' \rangle.
\]
Заметим, что это выражение может быть как положительным, так и отрицательным или нулевым.

В произведении $V \times V^*$ можно определить нелинейное многообразие
\[ M = \{ (x, \, -F'(x)) \mid x \in \interior K \},
\]
являющееся графом изометрии ${\cal L}$ между внутренностями прямого и двойственного конусов (см.~теорему \ref{thm:isometry}). Тогда точка $(w,w_*) \in M$ является ближайшей точкой на $M$ к текущей итерации $(x,s)$ в определённой выше псевдо-римановой метрике. Действительно, квадрат расстояния от точки $(x,s)$ до $(z, \, -F'(z)) \in M$ определяется выражением
\[ \langle x - z,s + F'(z) \rangle.
\]
Произведение $\langle F'(z),z \rangle$ является постоянным на $M$ вследствие логарифмической однородности функции $F$, а произведение $\langle x,s \rangle$ также не зависит от $z$. Остаётся в точности выражение $\Xi(z)$, экстремум которого определяет точку шкалировки.

В п.~\ref{ch1_geom_Newton} мы видели, что в методе Ньютона нелинейный граф градиента минимизируемой функции аппроксимируется аффинным подпространством, с помощью которого вычисляется следующая итерация. При этом аппроксимация строится в текущей точке.

В случае прямо-двойственного метода также строится аффинная аппроксимация нелинейного градиентного графа $M$. Однако текущая точка $(x,s)$ в общем случае не лежит на $M$. Поэтому логично построить аппроксимацию на базе ближайшей к текущей точке точки на $M$, которой и является прямо-двойственная пара $(w,w_*)$ точек шкалировки. 

Пара $(w,w_*)$ удовлетворяет условию $\langle w,w_* \rangle = \nu$, где $\nu$~--- параметр барьера. Однако по мере приближения к решению свёртка $\langle x,s \rangle$ стремится к нулю. Поэтому вводят фактор $\mu = \langle x,s \rangle/\nu$, и аппроксимируют не само многообразие $M$, а его <<сжатый>> образ:
\begin{align*} 
\sqrt{\mu}M &= \left\{ (\sqrt{\mu}x,\sqrt{\mu}s) \mid (x,s) \in M \right\} \\ &= \left\{ (x,-\mu F'(x)) \mid x \in \interior K \right\} = \left\{ (-\mu F_*'(s),s) \mid s \in \interior K^* \right\}.
\end{align*}

Пусть $x_*(\tau),s_*(\tau)$ --- точки центральных путей для прямой и двойственной \nk{задач}, отвечающие параметру $\tau$. Оказывается, что точка $(x_*(\tau),s_*(\tau))$ \emph{прямо-двойственного центрального пути} лежит на многообразии $\sqrt{\mu}M$, где $\mu = \tau^{-1}$. Поэтому прямо-двойственной паре $(x,s)$ можно сопоставить точку прямо-двойственного центрального пути с параметром $\tau = \mu^{-1} = \frac{\nu}{\langle x,s \rangle}$ . Аппроксимацией этой точки служит пересечение аффинной аппроксимации нелинейного многообразия $\sqrt{\mu}M$ с аффинным подпространством, задаваемым линейными ограничениями прямой и двойственной \nk{задач}.

\begin{figure}[ht]
\begin{center}
\includegraphics[width=11.05cm,height=6.32cm]{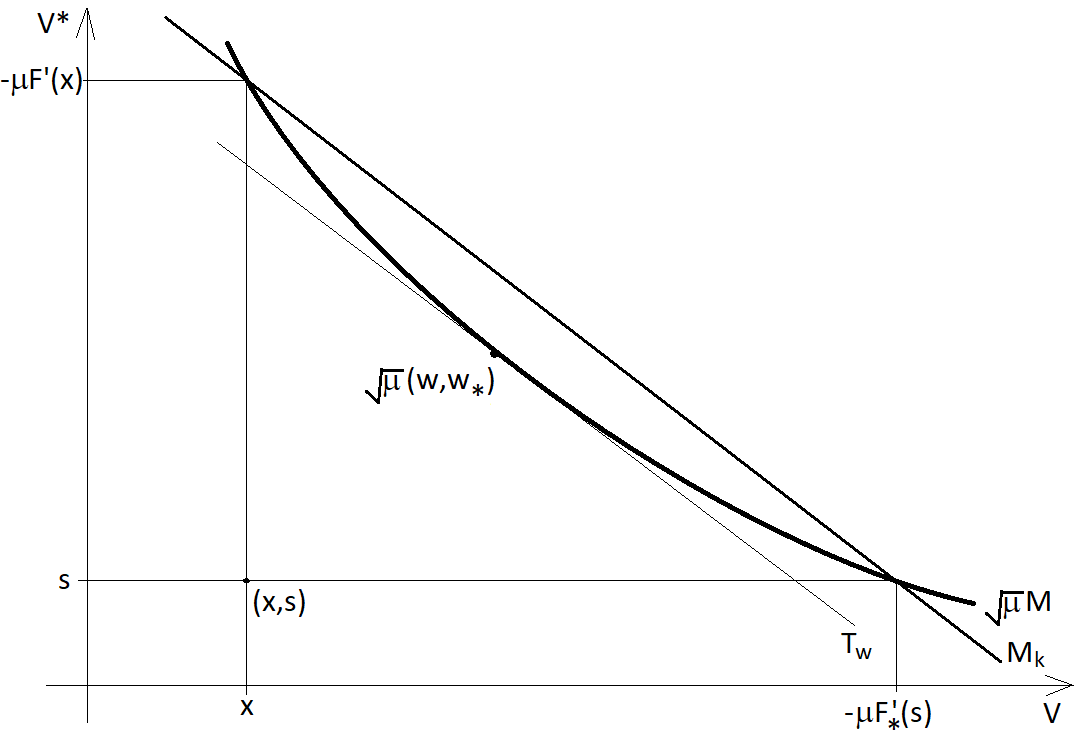}
\end{center}
\caption{Аффинные аппроксимации нелинейного графа градиента}
\label{fig:NTdirection}
\end{figure}

Можно аппроксимировать $\sqrt{\mu}M$ касательной плоскостью $T_w$ в точке $\sqrt{\mu}(w,w_*)$. Однако, как правило, используется аффинная плоскость $M_k$, параллельная к $T_w$ и пролегающая через точки $(x, \, -\mu F'(x)),(-\mu F'_*(s), \, s) \in\linebreak\in \sqrt{\mu}M$ (см.~рис.~\ref{fig:NTdirection}). Такая плоскость существует вследствие специальной структуры барьера на симметричных конусах.

Плоскость $M_k$ соответствует некоторой квадратичной аппроксимации барьерной функции $F$. На её базе строится аппроксимация центрального пути, являющаяся просто прямой в произведении $V \times V^*$, а именно --- пересечением линейного подпространства $\mathbb R \cdot M_k$ с аффинным подпространством, задающимся линейными ограничениями прямой и двойственной задач. Перпендикуляр из текущей точки на эту прямую называется \emph{центрирующим направлением} (centering direction). Направление на точку прямой, соответствующую значению параметра $\mu = 0$, т.е. решению задачи, называют \emph{направлением аффинной шкалировки} (affine scaling direction). Их выпуклые комбинации называются \emph{направлениями Нестерова--Тодда}. Делая шаг по направлению аффинной шкалировки, мы улучшаем значения функций цены прямой и двойственной задач в текущей точке, в то время как шаг вдоль центрирующего направления уменьшает расстояние до центрального пути.

В линейном программировании точка шкалировки для прямо-двойственной пары $(x,s)$ находится за $O(n)$ операций и задаётся вектором $w = \sqrt{x/s}$, где корень и дробь считаются по-элементно. Сложность итерации в линейном программировании в основном определяется сложностью решения линейной системы для нахождения направления шага.

В полуопределённом программировании для прямо-двойственной пары матриц $(X,S)$ точка шкалировки задаётся матрицей $W\!\! =\linebreak= X^{1/2}(X^{1/2}SX^{1/2})^{-1/2}X^{1/2}$. Последняя вычисляется за количество арифметических операций порядка $O(n^4)$, где $n$ --- размер матриц \cite{ToddTohTutuncu98}. Матрица $W$ имеет простую интерпретацию. В базисе пространства $\mathbb R^n$, в котором положительно определённые матрицы $X,S$ одновременно диагонализуются, матрица $W$ также диагональна, и её диагональные элементы равны геометрическому среднему соответствующих диагональных элементов матриц $X,S^{-1}$.

Теория авто-шкалированных барьеров позволяет делать шаги по направлению аффинной шкалировки не порядка $O(\nu^{-1/2})$, а порядка расстояния до границы допустимого множества. Поэтому на практике количество итераций, нужных для достижения определённой точности, слабо зависит от размерности задачи и составляет, как правило, всего несколько десятков. Длина шага в прямо-двойственных методах с длинным шагом ограничивается некоторой <<большой>> окрестностью центрального пути. 

Рассмотрим сначала случай линейного программирования. Ортант $\mathbb R_+^n$ обладает непрерывной группой автоморфизмов, а именно линейных отображений вида $x \mapsto Dx$, где $D$ является положительно определённой диагональной матрицей. Такое преобразование координат индуцирует автоморфизм $s \mapsto D^{-1}s$ двойственного конуса. Преобразование $(x,s) \mapsto (Dx,D^{-1}s)$ сохраняет вектор $y = x \cdot s \in \mathbb R_+^n$, в котором произведение определено по-элементно. Отметим, что многообразие $M$ в линейном программировании задаётся условием $s = x^{-1}$, интерпретируемым по-элементно. Поэтому для $(x,s) \in M$ имеем $y = {\bf 1} = (1,\dots,1)^T$. Таким образом, для любой пары $(x,s)$, лежащей на прямо-двойственном центральном пути задачи, вектор $y$ пропорционален вектору ${\bf 1}$, независимо от вида линейных ограничений, и <<число обусловленности>> $\frac{\max_i y_i}{\min_i y_i}$ вектора $y$ можно использовать для измерения расстояния до центрального пути. <<Большая>> окрестность центрального пути определяется как множество точек
\[ {\cal N}_{\gamma} = \left\{ (x,s) \left| \frac{\max_i x_is_i}{\min_i x_is_i} \leq \gamma \right. \right\},
\]
где типичное значение порога $\gamma$ порядка $10^3$.

В случае полуопределённого программирования, в котором ортант заменяется на матричный конус ${\cal S}_+^n$, окрестность ${\cal N}_{\gamma}$ определяется аналогично. Автоморфизм $X \mapsto AXA^T$ конуса ${\cal S}_+^n$ индуцирует автоморфизм $S \mapsto A^{-T}SA^{-1}$ двойственного конуса, где $A$ --- обратимая матрица размера $n \times n$. Эти преобразования оставляют спектр $y$ произведения $XS$ неизменным. С другой стороны, многообразие $M$ определяется соотношением $y = {\bf 1}$, и расстояние до центрального пути можно измерять с помощью числа обусловленности матрицы $XS$. Отсюда получаем определение
\[ {\cal N}_{\gamma} = \left\{ (X,S) \left| \frac{\lambda_{\max}(XS)}{\lambda_{\min}(XS)} \leq \gamma \right. \right\}.
\]

Вариант типа предиктор-корректор прямо-двойственного метода с длинным шагом идёт до границы окрестности ${\cal N}_{\gamma}$ по направлению аффинной шкалировки (фаза предиктор), после чего снова приближается к центральному пути с помощью шагов вдоль центрирующего направления (фаза корректор). Возможно также ходить по направлению комбинации аффинной шкалировки и центрирующего направления, при этом не выходя за пределы окрестности ${\cal N}_{\gamma}$. На практике используются разные эвристики, управляющие длиной шага и весами направлений в комбинациях \cite{SDPT3}.

\subsection{История методов внутренней точки} \label{subs:history_IP}

\ag{Методы внутренней точки вначале развивались как альтернатива симплекс-методу для решения задач линейного программирования. В отличие от последнего они генерируют последовательность точек, для которых ограничения типа неравенства строго выполнены. Вообще говоря, методы, генерирующие последовательность внутренних точек допустимого множества, применялись уже давно к решению нелинейных задач (см., например, п.~\ref{ch2_sect_fw}). Под методами внутренней точки мы здесь будем понимать методы, которые в явном или неявном виде используют барьер на допустимом множестве или на конусе, лежащем в основе задачи. В случае линейного программирования в роли этого конуса выступает неотрицательный ортант.}

\ag{Первым такой метод предложил И.\,И.~Дикин \cite{Dikin}. По сути, этот метод определял следующую точку как аргумент минимума линейной функции цены на эллипсоиде Дикина с центром в предыдущей итерации. Однако этот метод не обладал гарантиями на скорость сходимости, и ему требовалась строго допустимая начальная точка. В 1984 г. Н.~Кармаркар построил первый метод внутренней точки для линейного программирования, который доказуемо имел полиномиальную сложность \cite{Karmarkar}. Этот метод на каждом шаге применял проективное преобразование допустимого множества перед тем, как построить вписанный в это множество эллипсоид. Так как класс линейных функций не инвариантен по отношению к проективным преобразованиям, пришлось расширить этот класс до дробно-линейных функций, которые и минимизировались на построенных эллипсоидах. Этот недостаток метода Кармаркара послужил мотивацией для разработки вариантов, использующих аффинные преобразования, что в итоге привело к переоткрытию метода Дикина на Западе \cite{VanderbeiMeketonFreedman86, Barnes86}. Все последующие варианты метода Кармаркара могли конкурировать с симплекс-методом на практике, но не имели теоретических гарантий скорости сходимости. Заметим, что сам Кармаркар уже указывал на то, что, с одной стороны, его метод можно интерпретировать как барьерный, т.е. минимизирующий взвешенную сумму исходной функции цены и стандартного логарифмичного барьера на ортанте, а с другой стороны, как метод редукции потенциала.}

\ag{В серии работ Д. Байер и Дж. Лагариас (D.\,A.~Bayer, J.\,C.~Lagarias) подробно исследовали связь методов типа Кармаркара с барьерными методами. В частности, методы были классифицированы на проективно- и аффинно-масштабирующие, в зависимости от того, какие преобразования применяются к допустимому множеству. В дальнейшем развитие получили в основном аффинно-масштабирующие методы. Было также введено понятие центрального пути. В работе \cite{Renegar88}  Дж. Ренегар построил метод следования центральному пути, который выдавал последовательность точек, близких к центральному пути, и который доказуемо обладал полиномиальной сложностью. Прямо-двойственные методы были предложены в работах \cite{KojimaMizunoYoshise89, MonteiroAdler89}, а в \cite{KojimaMegiddoMizuno93, MizunoKojimaTodd95}~--- недопустимые прямо-двойственные методы, на итерациях которых ограничения типа равенства могли не выполняться. В работах \cite{Tanabe88, ToddYe90} были предложены методы редукции потенциала для линейного программирования. Предложенный в последних работах \emph{потенциал Танабэ--Тодда--Йе} (Tanabe--Todd--Ye potential) лёг в основу более поздних, продвинутых методов редукции потенциала для полуопределённого программирования.}

\ag{В конце 80-х гг. XX~в. Ю.\,Е.~Нестеров и А.\,С.~Немировский обобщили методы внутренней точки для линейного программирования на задачи оптимизации с произвольными выпуклыми коническими ограничениями и ввели понятия конической программы и самосогласованного барьера. В частности, они предложили методы с полиномиальной сложностью для решения полуопределённых программ. Также была построена теория конической двойственности, обобщающая двойственность в линейных программах \cite{NesNem92}. Полная теория изложена в монографии \cite{NesNem94}, в которой описаны как методы следования центральному пути, так и методы, понижающие потенциал. Позже Ю.\,Е.~Нестеров установил, что методы редукции потенциала делают шаги, похожие на длинные шаги в методах следования центральному пути \cite{NesLongStep}.}

\ag{Методы решения линейных программ независимо были обобщены на случай полуопределённого программирования в цикле работ Ф.~Ализадэ \cite{Alizadeh95}. Однако его методы использовали масштабирующие автоморфизмы конуса и в принципе не могли быть использованы для решения задач над несимметричными конусами.}

\ag{Следующим этапом было развитие теории автошкалированных барьеров и базирующихся на этой теории методов с длинным шагом, которые использовали богатую структуру симметричных конусов \cite{NesterovTodd97, NesterovTodd98, Monteiro97, MonteiroZhang98, Tuncel00}, см.~также монографию \cite{Renegar01}. Именно этот класс методов показал себя наиболее эффективным на практике и до сих пор используется в солверах \cite{ToddTohTutuncu98, SDPT3}.}

\ag{Заметим, что особая роль симметричных конусов как наиболее естественного обобщения неотрицательного ортанта в коническом программировании впервые была отмечена Л.~Файбузовичем, который в явном виде использовал их алгебраическую структуру для построения алгоритмов \cite{Faybusovich97a, Faybusovich97b}. В~то время как свойство автошкалированности присуще барьеру, симметричность является свойством конуса. В начале 2000-х гг. было обнаружено, что эти два понятия тесно взаимосвязаны и автошкалированные барьеры существуют на всех симметричных конусах и только на них. Обзор этих разработок представлен в \cite{HauserGuler02}.}

\ag{В работе \cite{Guler97} методы, опирающиеся на свойство автошка\-ли\-ро\-ван\-нос\-ти, были обобщены на класс \emph{гиперболических конусов}, а в \cite{Chua09} --- на класс однородных конусов. В работе \cite{Nesterov12} понятие точки шкалировки было определено также для произвольных самосогласованных барьеров.}

\ag{В последнем десятилетии было опубликовано несколько новых универсальных конструкций самосогласованных барьеров для конусов и выпуклых множеств \cite{Hildebrand14d,BubeckEldan19,Lee18}, а также для линейного программирования \cite[п.~6.3]{LeeThesis}. В работе \cite{AbernethyHazan} была найдена связь между методами следования центральному пути с использованием энтропического барьера и алгоритмом имитации отжига (simulated annealing). В работе \cite{NesterovTuncel16} вместо свойства авто-шкалированности для построения алгоритмов использовалось более слабое свойство \emph{отрицательной кривизны}, которое выполнено, в частности, для стандартных барьеров на гиперболических конусах. В последние годы наблюдается повышенный интерес к барьерам на несимметричных конусах, например, на экспоненциальном конусе \cite{DahlAndersen21}.}

\ag{Отметим, что наряду с методами внутренней точки получили развитие обобщения симплекс-метода на полуопределённое программирование, в~частности, в цикле работ В.\,Г.~Жадана.}



%
%
\section[Концепция неточной модели функции]{Концепция неточной модели функции\\ и методы градиентного типа для задач,\\ допускающих\,существование\,таких\,моделей\sectionmark{Концепция неточной модели функции}
}
\sectionmark{Концепция неточной модели функции}
\label{model}

Многие методы оптимизации в своей основе содержат идею замены оптимизируемой функции на некоторую более простую функцию-модель, которая достаточно хорошо аппроксимирует исходную функцию. После этой замены предлагается решать задачу оптимизации уже не для исходной функции, а для её модели. Это можно сказать, например, о методах градиентного типа, популярных сейчас в приложениях к задачам анализа больших данных. Помимо методов с простым градиентным шагом выделяют ещё и так называемые {\it ускоренные методы} (см.~п.~\ref{gradientTaylor} и \ref{AM}), \ag{позволяющие} получать оптимальные оценки сложности для гладких выпуклых задач. В основу как ускоренных, так и неускоренных градиентных методов положена идея аппроксимации функции в исходной точке (текущем положении метода) мажорирующим её параболоидом вращения и выбора точки минимума параболоида вращения в качестве нового положения метода. Такими методами исходно сложная задача заменяется совокупностью более простых задач.

В данном разделе пособия мы постараемся проиллюстрировать важность общего структурного подхода к методам первого порядка на базе абстрактной концепции модели функции \cite{Gas_Pos18}. Использование этого понятия позволяет обобщить концепцию неточного оракула (или ($\delta, \, L$)-оракула) из статьи Деволдера--Глинёра--Нестерова \cite{DevolderGlineurNesterov14}, которая в последние годы получила широкую известность в  оптимизационном сообществе. Понятия неточного оракула и неточной модели дают возможность описать многие естественно возникающие в приложениях \ag{типы} задач: постановки с предположениями о доступности на итерациях неточных значениях целевой функции, градиента или субградиента в негладком случае, задачи с функциями пониженного уровня гладкости с точки зрения свойственной для <<гладкого>> случая оптимизационной модели, задачи со специальной структурой (например, интересен класс задач композитной оптимизации, когда целевая функция есть сумма гладкой и негладкой функций простой структуры). Как оказалось, для задач, которые допускают существование неточного оракула или неточной оптимизационной модели, возможно предложить методы градиентного типа (в том числе ускоренные) с приемлемыми оценками скорости сходимости, не зависящими от размерности задачи. В случае отсутствия погрешностей (на итерациях используется точная информация о значениях функции и/или градиента) такие оценки оптимальны (см. таблицу \ref{p0_tab001} и \cite{NemirYd79}). При этом в разделе \ref{gradMethod} детально рассматривается и неускоренный градиентный метод, поскольку оценки его скорости сходимости оптимальны \cite{DragomirRelativeSmooth} на классе выпуклых оптимизационных задач с {\it относительно гладкими функционалами} \cite{Bauschke_Bolte_2016}, \cite{LuNesterov}. Что же касается ситуации, когда возможны погрешности оракула (доступно неточное значение целевой функции или её градиента), то тут важен следующий момент. \ag{Использование ускоренных градиентных методов может приводить к накоплению этих погрешностей в оценках скорости сходимости. Как следствие, ускоренный метод может давать худшие оценки скорости сходимости по сравнению с обычным неускоренным градиентным методом, что будет пояснено далее (см. \eqref{estim_mist_or1} и \eqref{estim_mist_or2}).}

В разных практических ситуациях (например, в задаче оптимизации газотранспортных систем \cite{Nesterov15b}) возникают задачи с целевыми функционалами пониженного уровня гладкости (промежуточной гладкости): градиент (или субградиент при $\nu = 0$) удовлетворяет условию Гёльдера при некотором $\nu \in [0, 1]$. При этом $\nu = 1$ означает, что градиент удовлетворяет условию Липшица, а $\nu = 0$~--- сама функция удовлетворяет свойству Липшица. Для задач минимизации выпуклых функций такого уровня гладкости Ю.\,Е.~Нестеровым предложены {\it универсальные градиентные методы}, являющиеся развитием методов менее известной работы Немировского--Нестерова~1985 г. \cite{NemirovskiNesterov85}. Они основаны на возможности построить единообразную оптимизационную квадратичную модель (свойственную для функций с липшицевым градиентом) для задач с гёльдеровыми градиентами за счёт введения искусственной неточности (степени близости функции и модели). При этом чем адекватнее такая модель описывает задачу, тем больше соответствующий параметр гладкости модели $L$ (аналог константы Липшица градиента). Увеличение $L$ ухудшает оценку качества решения. Поэтому очень интересной здесь будет идея <<адаптивной>> настройки метода на параметр $L$ гладкости модели, позволяющая в ходе практической реализации метода заменять потенциально большое глобальное значение $L$ на его локальные аналоги. Так появляется <<универсальность>> методов, под которой понимается возможность адаптивной настройки при работе метода на оптимальный в некотором смысле уровень гладкости задачи, а также на величину соответствующей константы Липшица (Гёльдера) градиента целевого функционала. Оказывается, что возможность такой настройки может позволить для некоторых задач экспериментально повысить скорость сходимости даже по сравнению с нижними теоретическими оценками \cite{Nesterov15b}. 

Отметим также, что использование концепции неточной абстрактной модели целевой функции \cite{Gas_Pos18, GasTurin} оптимизационной задачи позволяет расширить класс применимости градиентных методов на задачи со специальной структурой, а также задач с относительно гладкими функционалами \cite{Bauschke_Bolte_2016}, \cite{LuNesterov}. Однако платой за такое расширение рассматриваемого класса задач может быть усложнение вспомогательных подзадач на итерациях численных методов. \nk{По этой причине важно уметь решать такие задачи  с необходимой точностью, а также иметь} чёткое представление о влиянии возникающих при этом погрешностей на итоговые оценки скорости сходимости метода. Мы будем рассматривать подход \cite{BenTalNemirovski01} к описанию погрешностей решения вспомогательных подзадач на итерациях градиентных методов (см. пункт \ref{sect_inexact_sol}). Важно, что рассматриваемый подход к описанию таких погрешностей может приводить к теоретическим результатам об отсутствии их накопления и для ускоренных градиентных методов (см. теорему~\ref{mainTheoremDL}). Это обстоятельство позволяет, в частности, найти интересную связь методов градиентного типа в модельной общности с методами Франк--Вульфа \cite{GasTurin} (см. далее пункты \ref{model_uslovn}, \ref{universal_uslovn}).

В пункте \ref{restarts} для сильно выпуклых задач мы рассмотрим методику рестартов (перезапусков) методов для выпуклых задач. Такой подход позволяет существенно улучшить оценки скорости сходимости градиентных методов по сравнению с выпуклым случаем.

Напомним общую постановку задачи выпуклой оптимизации.
Пусть имеется выпуклая функция $F: Q \longrightarrow \mathbb{R}$ и дана произвольная норма $\|\cdot\|$ в банаховом пространстве $E = \mathbb{R}^n$ (можно рассматривать общие рефлексивные пространства, в частности, евклидовы). Для работы с (суб)градиентами функций важно использовать норму в сопряжённом пространстве $E^*$ (или двойственную к $\|\cdot\|$ норму), которая определяется следующим образом:
\begin{gather*}
\|\lambda\|_*=\max\limits_{\|\nu\| \leq 1;\nu \in \mathbb{R}^n}\langle \lambda,\nu\rangle,\,\text{ для всех  } \, \lambda \in \mathbb{R}^n.
\end{gather*}

Будем полагать, что
\begin{enumerate}
	\item $Q$~--- выпуклое и замкнутое подмножество $\mathbb{R}^n$;
	\item $F(x)$~--- непрерывная и выпуклая функция на $Q$;
	\item $F(x)$ ограничена снизу на $Q$ и достигает своего минимума в некоторой точке (необязательно единственной) $x^* \in Q$.
\end{enumerate}

Рассмотрим следующую задачу выпуклой оптимизации:
\begin{equation}
\label{mainTask3}
\min_{x \in Q} F(x).
\end{equation}

Как правило, при использовании стандартных методов первого порядка для задачи \eqref{mainTask3} требуется, чтобы функция $F$ имела $L$-липшицев градиент на $Q$ (см. разделы \ref{gradientTaylor}, \ref{sect:gradient_speed}). При таком допущении можно гарантировать выполнение неравенства
\begin{equation}\label{eq_lip_grad}
F(x) \leq F(y) + \langle \nabla F(y), x - y \rangle + \frac{L}{2} \|x - y\|^2 
\end{equation}
для произвольных $x, y \in Q$. Неравенство \eqref{eq_lip_grad}, как правило, позволяет обосновать оценки скорости сходимости некоторых градиентных методов. 

Однако на практике известно немало выпуклых дифференцируемых функций, градиент которых не удовлетворяет условию Липшица. Например, можно рассмотреть функцию $F(x)=-\log (x)+x^{2}$ на множестве положительных чисел $Q=\mathbb{R}_{++}$. Конечно, если предположить, что итерации алгоритма гарантируют монотонное убывание значений целевого функционала, то этого вполне достаточно для того, чтобы $F$ имела $L$-липшицев градиент на некотором множестве уровня. Как правило, такие рассуждения применимы к большинству методов первого порядка на классах гладких функций. Однако даже в таком случае постоянная $L$ может быть очень большой. Например, для $F(x)=-\log (x)+x^{2}$ на $Q=\mathbb{R}_{++}$ можно рассмотреть множество уровня $\{x:\,F(x)\leqslant 10\}$, для которого $L\approx e^{20}$. Такое большое значение параметра $L$ представляется нецелесообразным для практического использования. 

В этой связи интересен выделенный несколько лет назад в \cite{Bauschke_Bolte_2016} класс оптимизационных задач, связанный с условием {\it относительной гладкости}. Это условие обобщает неравенство \eqref{eq_lip_grad} путём его замены на неравенство
\begin{equation}\label{eq_rel_smooth}
F(x) \leq F(y) + \langle \nabla F(y), x - y \rangle + LV(x, y)  \text{  при всех  } x, y \in Q,
\end{equation}
где $V(x, y)$ --- некоторая специальная функция, которую называют {\it дивергенцией} или {\it расхождением Брэгмана}. Такой подход позволяет построить методы градиентного типа для некоторых важных типов задач выпуклой оптимизации, теряющих обычные свойства гладкости (липшицевости градиента). Для того чтобы работать с относительно гладкими задачами, необходимо ввести понятия прокс-функции и дивергенции Брэгмана \cite{Bregman67} (см. также \cite[гл.~5]{Ben-Tal}). 

\begin{defin}
Выпуклую функцию $d:Q \rightarrow \mathbb{R}$ будем называть прокс-функцией, если $d$ непрерывно дифференцируема в любой внутренней точке множества $Q$.
\end{defin}
\begin{defin}\label{defin_bregman}
Дивергенцией (или расхождением) Брэгмана называют функцию
\begin{align}\label{div-breg}
V(x,y) := d(x) - d(y) - \langle\nabla d(y), x - y\rangle,
\end{align}
где $d$~--- произвольная прокс-функция.
\end{defin}

\ag{В случае евклидовой нормы для вектора $x = (x_1, x_2, ..., x_n) \in \mathbb{R}^{n}$:
$$
\|x\|_2 = \sqrt{\sum\limits_{i = 1}^{n} x_i^2}
$$
можно положить
$$
d(x) = \frac{1}{2}\|x\|_2^2, \quad V(x,y) = \frac{1}{2}\|x - y\|_2^2.$$
Во многих ситуациях разумно работать и с неевклидовыми  прокс-стуктурами, некоторые их примеры рассмотрены далее в п. \ref{subsect_prox_strongly}. 
}
Много интересных примеров прикладных задач выпуклой оптимизации с относительной гладкостью приведены в \cite{LuNesterov}. Один из таких примеров (двойственная задача о нахождении эллипсоида минимального объёма, покрывающем заданный набор точек) описан далее в разделе \ref{rel_smooth}. \ag{Для этой задачи можно использовать прокс-функцию вида}
$$
d(x):=-\sum^{n}_{j=1}\log (x_{j}),
$$ 
которая выпукла, но не сильно выпукла относительно евклидовой нормы. Отметим также \cite{LuNesterov}, относительно гладкими при подходящем выборе прокс-структуры будут целевые функции, для которых матричная норма $\nabla^2 f(x)$ полиномиально зависит от $\|x\|_2$. Простейший пример такой функции $f(x) = x^4$ ($x \in \mathbb{R}$). Ясно, что такая функция не удовлетворяет условию Липшица градиента, т.е. стандартному предположению о гладкости. 

Если по смыслу решаемой оптимизационной задачи  целесообразно работать с обычной евклидовой нормой $\|\cdot\|_2$, то шаг градиентного спуска (см. \ref{gradientTaylor}) можно записать в виде вспомогательной задачи минимизации:
\begin{equation}\label{eq_euclid_step}
x_{k+1} = \argmin_{x \in Q}\left\{h\langle \nabla F(x_k), x - x_k\rangle + \frac{1}{2}\|x - x_k\|_2^2\right\}.
\end{equation}
На классе же относительно гладких задач ситуация усложняется и шаг градиентного спуска уже нужно записывать в виде
\begin{equation}\label{eq_bregman_step}
x_{k+1} = \argmin_{x \in Q}\{ h\langle \nabla F(x_k), x - x_k\rangle + V(x, x_k)\}.
\end{equation}
Однако при этом усложняются вспомогательные минимизационные подзадачи \eqref{eq_bregman_step}. Например, в случае $Q = \mathbb{R}^n$ задача вида \eqref{eq_bregman_step} сводится к решению относительно $x_{k+1}$ нелинейного уравнения $\nabla d(x_{k+1}) = \nabla d(x_k) - h \nabla f(x_k)$.

Оказывается, выпуклость прокс-функции $d$ и неравенство \eqref{eq_rel_smooth} позволяют получить оценки скорости сходимости неускоренного градиентного метода (см. п.~\ref{gradMethod}). Недавно было доказано \cite{DragomirRelativeSmooth}, что на достаточно широком классе выпуклых относительно гладких оптимизационных задач оценка сложности неускоренного метода оптимальна.

Однако на классе выпуклых и гладких в обычном смысле оптимизационных задач \eqref{eq_lip_grad} оценка сложности неускоренного градиентного метода уже не оптимальна (см.~п.~\ref{sect-black-box}). Для построения ускоренных методов, которые оптимальны на классе выпуклых гладких задач (см.~п.~\ref{fastGradMethod}), уже существенно предположение о $1$-сильной выпуклости прокс-функции $d$. Подробнее обсудим такие прокс-функции далее в п. \ref{subsect_prox_strongly}.

Отметим одно важное, но достаточно простое наблюдение, отчасти объясняющее возникновение условия $1$-сильной выпуклости $d$ относительно нормы $\|\cdot\|$. Если 
$$\|\nabla F(y) - \nabla F(x)\|_* \le L\|y-x\| \text{  при всех  } x, y \in Q,$$
то выполняется неравенство \eqref{eq_lip_grad}. А если ещё при этом и $V(y,x) \ge \frac{1}{2}\|y-x\|^2$ для всех $x, y \in Q$ (для этого как раз достаточно $1$-сильной выпуклости \ag{$d$ на множестве $Q$} относительно нормы $\|\cdot\|$), то выполняется  и неравенство \eqref{eq_rel_smooth}. Формулируя по-другому написанное выше, можно заметить, что для построения теории неускоренных градиентных методов на самом деле достаточно только последнего (более слабого) неравенства \ag{\eqref{eq_rel_smooth}}.

\subsection[Неточный оракул и неточная модель целевой\\ функции: понятия и базовые свойства]{Неточный оракул и неточная модель целевой\\ функции: понятия и базовые свойства}
\label{subsect_inexact_model}

Многие методы первого порядка (методы градиентного типа) основаны на использовании квадратичной функции (параболоида вращения) в качестве оптимизации модели. Однако необязательно рассматривать именно параболоиды вращения. Можно использовать и другие функции, которые позволяют построить более точную локальную модель целевой функции задачи в рассматриваемой точке, что приводит в итоге к увеличению скорости сходимости метода. Здесь можно выделить два направления. Первое направление связано с использованием старших производных в модели функции (см. п.~\ref{AM}), а второе~--- с занесением в модель части постановки задачи. Например, если оптимизируемая функция есть сумма гладкой и негладкой функций, то первую можно заменять параболоидом вращения в модели, а вторую оставить как есть. Насколько нам известно, по такому подходу до последнего момента имелись лишь разрозненные результаты. Стоит отметить, что в случае с параболоидом вращения решение вспомогательной задачи, как правило, не представляет большого труда и часто может быть сделано по явным формулам, то есть без ошибок. Совсем другая ситуация возникает при рассмотрении второго подхода, в котором типична обратная ситуация~--- вспомогательную задачу можно решить только приближённо. В этой связи важно проработать возможность неточного решения вспомогательной задачи. Недавно в \cite{Gas_Pos18} предпринята попытка оформить этот подход на базе абстрактного понятия ($\delta, L$)-модели функции, которое есть прямое обобщение известного понятия \textit{$(\delta, L)$-оракула} Деволдера--Глинёра--Нестерова. Начнём с необходимых определений. Отметим, что в приводимых далее определениях вместо $\frac{L}{2} \|x - y\|^2$ можно писать $LV(x,y)$ и использовать неравенство \eqref{eq_rel_smooth}, о чем мы уже немного упоминали ранее. Однако в этом \ag{пункте} мы ограничимся общностью~\eqref{eq_lip_grad}. 

\begin{defin}
\leavevmode
\label{delta_L_oracle}
Будем говорить, что пара $(F_{\delta}(y), G_{\delta}(y)) \in \mathbb{R} \times E^*$ есть $(\delta, L)$-оракул функции $F(x)$ в точке $y$ при $\delta, L>0$, если для любого $x \in Q$ справедливо неравенство
\begin{gather}
\label{DLoracleineq}
0 \leq F(x) - F_{\delta}(y) - \langle G_{\delta}(y), x - y \rangle \leq \frac{L}{2} \|x - y\|^2 + \delta.
\end{gather}
\end{defin}

В качестве упражнений приведём здесь несколько важных свойств $(\delta, \, L)$-оракула функции в точке.
\begin{exercise}
$F_{\delta}(y)$ из~\eqref{DLoracleineq} есть некоторое $\delta$-приближение точного значения значения функции $F$. Действительно, положив $x=y$ в~\eqref{DLoracleineq}, получим
\begin{equation}\label{4}
F_{\delta}(y)\leqslant F(y)\leqslant F_{\delta}(y)+\delta.
\end{equation}
\end{exercise}

\begin{exercise} Покажите, что $G_{\delta}(y)$ из~\eqref{DLoracleineq} есть $\delta$-субградиент $F$ на $y\in Q$, то есть
$$G_{\delta}(y)\in \partial_{\delta}F(y)=\{z\in E^{*}: \, F(x)\geqslant F(y)+\langle z,x-y\rangle -\delta \text{ при всех } x\in Q\}.$$ 
\end{exercise}
В качестве указания можно отметить, что действительно, используя первое неравенство в~\eqref{DLoracleineq} и (\ref{4}), мы имеем для всех $x,y\in Q$
\begin{equation}\label{5}
F(x)\geqslant F_{\delta}(y)+\langle G_{\delta}(y),x-y\rangle\geqslant F(y)+\langle G_{\delta}(y),x-y\rangle - \delta.
\end{equation}

\begin{remark}
Ясно, что обычный субградиент выпуклой функции также может удовлетворять второму неравенству в~\eqref{DLoracleineq}, что указывает на возможность использования концепции неточного оракула для задач негладкой выпуклой оптимизации.

Отметим, что методы негладкой выпуклой оптимизации с использованием $\delta$-субградиентов довольно хорошо известны и активно используются в самых разных типах задач.
\end{remark}

\begin{exercise}\label{ex_oracle_start} Если $F$ имеет $(\delta,L)$-оракул, то для всякой положительной константы $c>0$ функция $cF$ имеет $(c\delta,cL)$-оракул.
\end{exercise}

\begin{exercise} Если $F_{i}$ имеет $(\delta_{i},L_{i})$-оракулы при $i = 1, 2$, то $F_{1} + F_{2}$ допускает $(\delta_{1}+\delta_{2},L_{1}+L_{2})$-оракул.
\end{exercise}

\begin{exercise}\label{ex_oracle_fin} В случае, когда $Q$ совпадает со всем пространством $\mathbb{R}^n$, разность между $G_{\delta}$ и любым субградиентом $G_{y}\in\partial F(y)$ ограничена следующим образом:
\begin{equation}\label{6}
\|G_{y}-G_{\delta}(y)\|_{*}\leq \sqrt{2\delta L}.
\end{equation}
Отметим, что в случае, когда $Q$ не совпадает со всем пространством $\mathbb{R}^n$ известны более нетривиальные соотношения между $G_y$ и $G_{\delta}(y)$ (\cite{DevolderGlineurNesterov14}, с.~6).
\end{exercise}

Интересно, что концепция Деволдера--Глинёра--Нестерова может допускать следующее обобщение \cite{GasTurin} путём замены линейной по переменной $x$ (при фиксированном $y$) функции $\langle G_{\delta}(y), x - y \rangle$ в~\eqref{DLoracleineq} на некоторую абстрактную выпуклую функцию.
\begin{defin}
\leavevmode
\label{gen_delta_L_oracle}
Будем говорить, что существует $(\delta, L)$-модель функции $F$ для некоторых $\delta >0$ и $L>0$ в точке $y$, и обозначать эту модель $(F_{\delta}(y), \psi_{\delta}(x,y))$, если для любого $x \in Q$ справедливо неравенство
\begin{gather}
\label{exitLDLOrig}
0 \leq F(x) - F_{\delta}(y) - \psi_{\delta}(x,y) \leq \frac{L}{2} \|x - y\|^2 + \delta,
\end{gather}
\begin{gather}
\label{exitLDLOrig_eqaul}
\psi_{\delta}(x,x)=0 \text{ при всех } x \in Q
\end{gather}
и $\psi_{\delta}(x,y)$~--- выпуклая функция по $x$ для  всякого $ y \in Q$.
\end{defin}
Будем считать, что для $F$ найдутся такие $\delta>0$ и $L>0$, что существует $(\delta, L)$-модель в любой точке $x \in Q$.

Если выбрать $x=y$ в (\ref{exitLDLOrig}) и воспользоваться (\ref{exitLDLOrig_eqaul}), то получим 
\begin{gather}
\label{exitLDLOrig2}
F_{\delta}(y) \leq F(y) \leq F_{\delta}(y) + \delta \text{ для произвольных } y \in Q.
\end{gather}

\begin{exercise} Сформулируйте аналоги упражнений \ref{ex_oracle_start} -- \ref{ex_oracle_fin} $(\delta, L)$-оракула для абстрактной концепции $(\delta, L)$-модели функции.
\end{exercise}

Далее (преимущественно в качестве упражнений) обсудим некоторые примеры постановок задач, которые можно решать методами градиентного типа на базе рассмотренной концепции модели оптимизируемой функции.

\subsubsection{Композитная оптимизация}

Оказывается, что абстрактная концепция модели функции позволяет распространить ускоренный вариант градиентного метода на задачи композитной выпуклой оптимизации \cite{Nesterov13}:
\begin{align}
\label{composite}
\min_{x \in Q} F(x) := f(x) + \varphi(x),
\end{align}
где $f$~--- гладкая выпуклая функция с $L$-липшицевым градиентом относительно нормы $\|\cdot\|$ и $\varphi$~--- выпуклая функция (вообще говоря, негладкая). Функция $\varphi$ называется {\it композитом}. Для данной задачи верно следующее неравенство:
\begin{gather}
0 \leq F(x) - F(y) - \langle\nabla f(y), x - y \rangle - \varphi(x) + \varphi(y) \leq \frac{L}{2} \|x - y\|^2\\ \text{ при всех } x,y \in Q. \nonumber
\end{gather}
Это означает, что возможно выбрать
$$
\psi_{\delta}(x,y)=\langle\nabla f(y), x - y \rangle + \varphi(x) - \varphi(y), \; F_{\delta}(y)=F(y) \text{ при } \delta=0.
$$
Такое неравенство позволяет получать для задач выпуклой композитной оптимизации оценки скорости сходимости градиентных методов (см. пункты \ref{gradMethod} и \ref{fastGradMethod}), аналогичные оценкам на обычном классе гладких выпуклых задач.

Задачи композитной оптимизации нередко возникают в приложениях. К одному из наиболее известных примеров стоит отнести так называемую задачу LASSO (Least Absolute Shrinkage and Selection Operator), мотивированную задачами статистики. По сути, это задача линейной регрессии согласно методу наименьших квадратов (подробнее про эту классическую задачу можно почитать, например, в \cite{Boyd_book_2004}, п. 1.2.1), но с добавлением $\ell_1$-нормы как регуляризатора (именно это слагаемое и рассматривается как композит). Использование такой регуляризации~--- один из подходов к проблеме неустойчивости стандартной регрессионной модели, связанной с методом наименьших квадратов. Суть проблемы в том, что метод наименьших квадратов, как известно, может быть сильно чувствительным к аномальным данным измерений (см.,~например,~\cite{Friedman2009}).

Также в качестве конкретного примера можно рассмотреть задачу восстановления матрицы корреспонденций по замерам потоков на линках (рёбрах) в большой компьютерной сети (Minimal Mutual Information Model), которая сводится к задаче композитной оптимизации \cite{Anikin2015} вида

$$\min \limits_{x \in S_n(1)} F(x) = \dfrac{1}{2} \|Ax - b\|_2^2+ \mu \sum \limits_{k=1}^n x_k \ln x_k,$$
где $S_n(1)$~--- единичный симплекс в $n$-мерном пространстве. Этот пример немного подробнее рассматривается в разделе~\ref{restarts}.

\subsubsection{Проксимальный метод для негладких задач выпуклой оптимизации}

Рассмотрим задачу
\begin{align}
\label{main_prox_method}
\min_{x \in Q} F(x),
\end{align}
где $F(x)$~--- вообще говоря, негладкая выпуклая функция. В качестве \ag{($\delta, L$)-модели в таком случае} можно выбрать $$\psi_{\delta}(x,y) = F(x) - F(y), \;\;\; F_{\delta}(y)=F(y) \text{ и }\delta=0 .$$
Условие (\ref{exitLDLOrig}) будет выполнено при любом выборе $L \geq 0$. 

На первый взгляд может показаться, что описанный здесь способ формирования модели не является практичным. Однако в последнее время такой подход оказался достаточно популярным при разработке универсальных ускоренных оболочек (см., например, \cite{ParikhBoyd14}, \cite{LinMairalHarchaoui15}, \cite{Doikov2019} и пп. \ref{AM}, \ref{AM_Applications}).

\subsubsection{Суперпозиция функций}

Рассмотрим следующую задачу \cite{NemirovskiNesterov85}:
\begin{align}
 \min_{x \in Q} F(x) := f(f_1(x), \, \dots, \, f_m(x)),
\end{align}
где $f_i$~--- гладкая выпуклая функция с $L_i$-липшицевым градиентом в норме $\|\cdot\|$ для всякого $i$, $f$~--- $M$-липшицева выпуклая функция относительно $L_1$-нормы, неубывающая по каждому из своих аргументов. Поэтому $F$ также выпуклая функция, причём имеют место неравенства
\begin{gather*}
0 \leq F(x) - f(f_1(y) + \langle\nabla f_1(y), x - y \rangle, \, \dots, \, f_m(y)+\langle\nabla f_m(y), x - y \rangle) \leq\\\leq M\frac{\sum_{i=1}^{m}L_i}{2} \|x - y\|^2 \text{   для произвольных  } x,y \in Q.
\end{gather*}
Далее, 
\begin{gather*}
0 \leq F(x) - F(y) - f(f_1(y) + \langle\nabla f_1(y), x - y \rangle, \, \dots, \, f_m(y)+\langle\nabla f_m(y), x - y \rangle) +\\ +\, F(y) \leq M\frac{\sum_{i=1}^{m}L_i}{2} \|x - y\|^2 \text{  для произвольных  } x,y \in Q.
\end{gather*}
Поэтому можно выбрать
$$\psi_{\delta}(x,y)=f(f_1(y) + \langle\nabla f_1(y), x - y \rangle, \, \dots, \, f_m(y)+\langle\nabla f_m(y), x - y \rangle) - F(y),
$$
$$ F_{\delta}(y)=F(y), L = M\cdot\left(\sum_{i=1}^{m}L_i\right)\text{ и }\delta=0.
$$

Рассмотрим без подробного объяснения ещё некоторые примеры постановок задач, в которых может быть актуальной рассматриваемая нами концепция модели функции.

\subsubsection{Минмин задача}

Рассмотрим следующую задачу:
\begin{align}
\min_{x \in \mathbb R^n} f(x) := \min_{y \in Q}F(y,x).
\label{main_minmin}
\end{align}

Пусть $F(y,x)$~--- дифференцируемая функция и для произвольных ${y,y' \in Q}$,\  $ x,x' \in \mathbb R^n$ выполнено неравенство
\begin{gather*}
\|\nabla F(y',x') -\nabla F(y,x)\|_2 \leq L \|(y',x') -(y,x)\|_2.
\end{gather*}
Тогда если можно найти такую $\widetilde{y}_\delta(x) \in Q$, что
\begin{gather*}
\langle\nabla_y F(\widetilde{y}_\delta(x), x), y - \widetilde{y}_\delta(x)\rangle \geq -\delta \,\, \text{ для всякого } y \in Q,
\end{gather*}
то верны неравенства (см. \cite{GasnikovDvurechenskyNesterov16})
\begin{gather*}
F(\widetilde{y}_\delta(x), x) - f(x) \leq \delta, \;\;\; \|\nabla f(x') -\nabla f(x)\|_2 \leq L \|x' -x\|_2.
\end{gather*}
Это означает, что
\begin{gather*}
\left(F_{\delta}(x)= F(\widetilde{y}_\delta(x), x) - 2\delta, \; \psi_{\delta}(z,x)=\langle\nabla_y F(\widetilde{y}_\delta(x), x), z - x\rangle\right)
\end{gather*}
есть $(6\delta, 2L)$-модель для функции $f$ в точке $x$.

\begin{exercise} Проведите описанные выкладки более подробно. \end{exercise}

\subsubsection{Седловая задача}

Рассмотрим задачу нахождения седловой точки следующего вида:
\begin{align}
\min_{x \in \mathbb R^n} f(x) := \max_{y \in Q}\left[\langle x, b - Ay\rangle - \varphi(y)\right],
\label{main_saddle}
\end{align}
где $\varphi(y)$~--- $\mu$-сильно выпуклая относительно $p$-нормы ($1\leq p\leq2$). Тогда $f$~--- гладкая функция с константой Липшица градиента в $2$-норме \cite{DevolderGlineurNesterov14}
\begin{gather*}
L=\frac{1}{\mu}\max_{\|y\|_p\leq1}\|Ay\|_2^2.
\end{gather*}
\begin{exercise} Докажите это равенство.
\end{exercise}

Если $y_\delta(x)$~--- решение вспомогательной задачи максимизации с точностью по функции $\delta$, то
\begin{gather*}
\left(F_{\delta}(x)=\langle x, b - Ay_\delta(x)\rangle - \varphi(y_\delta(x)), \; \psi_{\delta}(z,x) =\langle b - Ay_\delta(x), z - x\rangle\right)
\end{gather*}
есть $(\delta, 2L)$-модель для функции $f$ в точке $x$.

\subsection[О подходах к описанию погрешностей при решении\\ вспомогательных подзадач, возникающих\\ на итерациях численных методов]{О подходах к описанию погрешностей при решении вспомогательных подзадач, возникающих\\ на итерациях численных методов}
\label{sect_inexact_sol}

Как правило, на практике невозможно достичь точного решения как задачи оптимизации в целом, так и при решении вспомогательных подзадач на каждой итерации численного метода. Если на итерациях накапливаются погрешности решения вспомогательных подзадач, то естественно ожидать, что они окажут влияние на качество итогового результата работы метода. Важно уметь оценивать степень такого влияния, поскольку теоретические результаты о сходимости методов формулируются, как правило, в предположении о точном решении вспомогательных подзадач. Особенно это актуально для задач с неевклидовой прокс-структурой и относительной гладкостью. Если ещё рассмотреть такие предположения в модельной общности, то вспомогательная подзадача примет вид
\begin{equation}\label{eq_bregman_step_model}
x_{k+1} = \argmin_{x \in Q} \left \{ h\psi_{\delta}(x, x_k) + V(x, x_k) \right \}
\end{equation}
для некоторой выпуклой по первой переменной функции $\psi_{\delta}$. В общем случае не известны методы нахождения точного решения задач вида \eqref{eq_bregman_step_model}. Поэтому мы уделим внимание вопросу учёта в оценках скорости сходимости градиентных методов погрешностей, возникающих за счёт неточностей решения вспомогательных подзадач. Рассмотрим следующий подход к определению неточного решения задачи выпуклой минимизации \cite{BenTalNemirovski01}.
\begin{defin}
\label{solNemirovskiy}
Пусть имеется задача
\begin{gather*}
\min_{x \in Q} \psi(x),
\end{gather*}
где $\psi(x)$~--- выпуклая функция. Тогда $ \textnormal{Arg}\min_{x \in Q}^{\widetilde{\delta}}\psi(x)$~--- множество таких $\widetilde{x}$, что
существует вектор $\tilde{g}$ из субдифференциала $\partial\psi(\widetilde{x})$ такой, что
\begin{gather}\label{eqv_inex_sol}
\langle \tilde{g}, x - \widetilde{x} \rangle \geq -\widetilde{\delta} \text{  для всякого  } x \in Q.
\end{gather}
Некоторый выбранный элемент из $\textnormal{Arg}\min_{x \in Q}^{\widetilde{\delta}}\psi(x)$ будем обозначать как $\arg\min_{x \in Q}^{\widetilde{\delta}}\psi(x)$.
\end{defin}

Ясно, что при $\widetilde{\delta} = 0$ неравенство~\eqref{eqv_inex_sol} означает, что $\widetilde{x}$~--- точное решение задачи минимизации $\psi$ на множестве $Q$ и это поясняет смысл концепции неточности из определения~\ref{solNemirovskiy}. Однако, как правило, методы оптимизации гарантируют достижение достаточно хорошей аппроксимации решения задачи не по аргументу, а по функции. Сравним эти подходы к описанию неточности решения задачи выпуклой минимизации.

Пусть $\widetilde{x} \in \textnormal{Arg}\min_{x \in Q}^{\widetilde{\delta}}\psi(x)$, тогда выпуклость $\psi$ означает, что: $$\psi(x) \geq \psi(\widetilde{x}) + \langle \tilde{g}, x - \widetilde{x} \rangle \geq \psi(\widetilde{x})-\widetilde{\delta},$$
где $\tilde{g} \in \partial\psi(\widetilde{x})$. Если $x=x_*$ --- одно из точных решений вспомогательной подзадачи минимизации $\psi$ по $x$, то верно неравенство $\psi(\widetilde{x}) - \psi(x_*) \leq \widetilde{\delta}$. Поэтому $\widetilde{x} \in \textnormal{Arg}\min_{x \in Q}^{\widetilde{\delta}}\psi(x)$ означает, что $\widetilde{x}$~--- $\widetilde{\delta}$-решение задачи минимизации $\psi(x)$ по функции. Обратное утверждение, вообще говоря, неверно. Но мы попробуем представить довольно общие примеры, когда обратное утверждение всё же можно считать верным\footnote{Тривиальный случай, когда $\widetilde{\delta} = 0$. Тогда из критерия оптимальности первого порядка будет следовать, что данные два определения $\widetilde{\delta}$-решения будут эквивалентны.}. Это важно, поскольку далее в обосновании оценок скорости сходимости алгоритмов \ref{Grad_Model} и \ref{Acceler_Grad_Model} (лемма \ref{lemma_maxmin_2} и теоремы \ref{mainTheoremDL_G} и \ref{mainTheoremDL}) существенно используется именно точность решения вспомогательных подзадач в смысле определения ~\ref{solNemirovskiy}, а не просто точность их решения по функции. 

Предположим, что вспомогательная задача, возникающая на итерации метода градиентного типа \eqref{eq_bregman_step_model}, имеет вид
\begin{equation}\label{helpTask}
\min_{x \in Q} \left \{ \psi(x) := h\psi_{\delta}(x, \, x_k) + V(x, \, x_k) \right \},
 \end{equation}
где $\psi_{\delta}(x, x_k)$~--- выпуклая функция по $x$.

Отметим, что вспомогательная подзадача в итерациях методов оптимизации довольно часто имеет такой вид (см. далее алгоритмы \ref{Grad_Model} и \ref{Acceler_Grad_Model}). Конечно, существуют случаи, когда задачу вида \eqref{helpTask} можно решить аналитически. Например, это так, когда в основной задаче решается гладкая задача оптимизации без ограничений с евклидовой прокс-структурой $V(x,y) = (1/2)\|x-y\|_2^2$. В случаях, когда задача (\ref{helpTask}) может быть решена только численно, используют различные подходы в зависимости от постановки задачи.

\ag{Например, рассмотрим} случай сепарабельной задачи в (\ref{helpTask}): 
$$
h\psi_{\delta}(x, \, x_k) + V(x, \, x_k) = \sum_{i=1}^{n}\left[h\psi_i(x_i) + V_i(x_i)\right],
$$
где $\psi_i$ и $V_i$~--- некоторые функционалы от $x_i$ ($x = (x_1,\dots,x_n)$). В таком случае достаточно решить $n$ одномерных задач, каждую из которых можно решать методом деления отрезка пополам (или золотого сечения, см.~п.~\ref{subs:bisection}) за время (необходимое количество итераций)
$$O\left(\log\left(\frac{1}{\varepsilon}\right)\right),$$ где $\varepsilon$ --- точность найденного решения по функции. Заметим также, что условие сепарабельности прокс-функции можно ослабить, сохраняя эффективность решения задачи~\eqref{helpTask}. Важный пример разобран в п.~\ref{comp_prox}.

Часто в приложениях прокс-функцию $d$ выбирают $1$-сильно выпуклой относительно некоторой нормы $\|\cdot\|$, см. п.~\ref{subsect_prox_strongly}. В этом случае появляются дополнительные степени свободы для эффективного решения задачи \eqref{helpTask}. 
Например, если дополнительно предположить, что $\psi$ имеет $L$-липшицев градиент в норме $\|\cdot\|$. Если $V(x,x_k)$ также имеет $L$-липшицев градиент по $x$ в норме $\|\cdot\|$, то для задачи \eqref{helpTask} можно предложить численный метод c линейной скоростью сходимости.
Если же $V(x, \, x_k)$ не имеет $L$-липшицева градиента в норме $\|\cdot\|$, можно в задаче \eqref{helpTask} рассматривать слагаемое $V(x,x_k)$ как композит. Чтобы получить численный метод для таких задач с линейной \nk{скоростью} сходимости, можно воспользоваться, например, техникой рестартов градиентных методов (см.~далее п.~\ref{restarts}).

Итак, при довольно общих предположениях можно добиться линейной скорости сходимости для вспомогательной подзадачи. Предположим, что $\psi(x)$ --- $M$-липшицева функция по норме $\|\cdot\|$. Тогда, если $x_{\tilde{\delta}}$ -- $\varepsilon$-решение по функции задачи минимизации $\psi(x)$, то верно неравенство
\begin{gather}
\langle\nabla \psi(x_{\tilde{\delta}}), x_{\tilde{\delta}} - x_*\rangle \leq \|\nabla \psi (x_{\tilde{\delta}})\|_* \|x_{\tilde{\delta}} - x_*\| \leq M \sqrt{2 \varepsilon}.
\label{bound_1}
\end{gather}
\begin{exercise}
Последнее неравенство следует из $1$-сильной выпуклости $V(x,x_k)$ в норме $\|\cdot\|$. Покажите это, проведя подробные выкладки.
\end{exercise}

Если $\psi$ имеет $L$-липшицев градиент относительно нормы $\|\cdot\|$, то при $Q=\mathbb R^n$ оценку (\ref{bound_1}) можно улучшить:
\begin{gather}
\langle\nabla  \psi(x_{\tilde{\delta}}), x_{\tilde{\delta}} - x_*\rangle \leq 2 \varepsilon
\sqrt{L}.
\label{bound_2}
\end{gather}

Таким образом, при сделанных предположениях всякое \ag{$\varepsilon$}-решение задачи минимизации \eqref{helpTask} по функции будет ($M \sqrt{2 \varepsilon}$)-решением или ($2\varepsilon\sqrt{L}$)-решением в смысле определения~\ref{solNemirovskiy}. 

Таким образом, рассматриваемая нами концепция приближённого решения задачи минимизации выпуклой функции определения~\ref{solNemirovskiy} может считаться вполне приемлемой для описания влияния на итоговый результат возникающих погрешностей на вспомогательных для градиентных методов подзадачах вида \eqref{helpTask}.

\subsection[Градиентный метод с оракулом, использующим на итерациях неточную модель функции]{Градиентный метод с оракулом, использующим\\ на итерациях неточную модель функции}
\label{gradMethod}


В данном пункте мы рассмотрим адаптивный вариант неускоренного градиентного метода (см. раздел \ref{gradientTaylor}) для задачи (\ref{mainTask3}), который применим к весьма общему классу задач, допускающих существование неточной модели целевой функции (см. п. \ref{subsect_inexact_model}) в произвольной точке области определения. Оказывается, что для обоснования оценок скорости сходимости неускоренного градиентного метода не требуется $1$-сильная выпуклость используемой прокс-функции, что позволяет применять неускоренные градиентные методы к широкому классу {\it относительно гладких} выпуклых оптимизационных прикладных задач \cite{Bauschke_Bolte_2016,LuNesterov}. Напомним \cite{Bauschke_Bolte_2016}, что функция $F$ называется относительно гладкой, если для произвольных $x, y \in Q$ верно неравенство \eqref{eq_rel_smooth}.

Введём аналог понятия ($\delta, L$)-модели функции (см. определение \ref{gen_delta_L_oracle}), который применим и к относительно гладким задачам.
\begin{defin}
\leavevmode
\label{gen_delta_L_model}
\nk{Будем} говорить, что существует $(\delta, L, V)$-модель функции $F$ относительно дивергенции $V(x, y)$ в точке $y$ при некоторых $\delta>0$ и $L>0$, и обозначать эту модель $(F_{\delta}(y), \psi_{\delta}(x,y))$, если для любого $x \in Q$ справедливо неравенство
\begin{gather}
\label{exitLDLOrigV}
0 \leq F(x) - F_{\delta}(y) - \psi_{\delta}(x,y) \leq LV(x, y) + \delta,
\end{gather}
\begin{gather}
\label{exitLDLOrig_eqaulV}
\psi_{\delta}(x,x)=0 \text{ при всех } x \in Q
\end{gather}
и $\psi_{\delta}(x,y)$~--- выпуклая функция по $x$ для всякого $ y \in Q$.
\end{defin}
\nk{Предположим}, что для $F$ найдутся такие $\delta>0$ и $L>0$, что существует $(\delta, L,V)$-модель в любой точке $x \in Q$ относительно $V(x, y)$.

Будем \ag{считать}, что дана начальная точка $x_0$, $N$~--- количество шагов метода, $L_0$~--- константа, которая имеет смысл предположительной <<локальной>> константы Липшица градиента в точке $x_0$. Также на вход алгоритму \ref{Grad_Model} (как и ускоренному алгоритму \ref{Acceler_Grad_Model} далее) подаются константы $\delta$ и $\widetilde{\delta}$: $\delta$~--- положительное число такое, что существует $(\delta, L, V)$-модель $F$ в любой точке $x \in Q$, $\widetilde{\delta}$~--- величина, соответствующая точности решения вспомогательных подзадач согласно определению~\ref{solNemirovskiy}. При этом в зависимости от задачи эти величины могут быть как равными нулю, так и иметь постоянное положительное значение. 

В этом пункте \ag{мы рассмотрим вариант неускоренного градиентного метода с {\it адаптивной настройкой} на величину константы Липшица градиента целевого функционала (алгоритм \ref{Grad_Model})}. Точнее говоря, явно в методе не будет использоваться константа $L$. Она будет выбираться специальным образом на каждой итерации и может как увеличиваться, так и уменьшаться. Для \ag{повышения качества} работы метода важно постараться выбрать как можно меньшую константу $L_k$ на каждой итерации, но при этом для гарантированного достижения требуемого качества решения должно выполняться неравенство-критерий выхода из итераций~\eqref{exitLDL_G}. Таким образом, осуществляется балансировка стремления повысить скорость работы метода и необходимости иметь гарантированную в теоретическом плане точность решения. Адаптивная настройка на величину параметра $L$ представляется особо важной для идеологии {\it универсальных методов} решения задач с пониженными требованиями (условие Гёльдера градиента) к гладкости функционалов (см. далее пункты \ref{Univers_Meth} и \ref{universal_uslovn}). Этот подход основан на построении похожей на <<гладкий>> случай оптимизационной ($\delta, L, V$)-модели для задачи с пониженной гладкостью, но за счёт введения в эту модель искусственной неточности $\delta > 0$. При этом на такую неточность завязана также итоговая оценка качества выдаваемого решения. Однако при этом константа $L = L(\delta)$ обратно пропорционально зависит от этой величины $\delta > 0$. Чем меньше $\delta > 0$, тем лучшего качества решения задачи можно достичь. Однако при этом величина $L = L(\delta)$ растёт, что и приводит к идее замены в оценке скорости сходимости глобального параметра $L$ его локально подбираемыми на итерациях методов аналогами $L_k$, которые потенциально могут оказаться меньшими \ag{$L$}. Аналогичный подход для ускоренного метода будет описан далее в п. \ref{fastGradMethod} (см. алгоритм \ref{Acceler_Grad_Model}).

Будем считать, что $L_0 \leq 2L$. Если не делать такого предположения, то в теоретических оценках \ag{качества выдаваемого приближённого решения} из теорем \ref{mainTheoremDL_G} и \ref{mainTheoremDL} для алгоритмов \ref{Grad_Model} и \ref{Acceler_Grad_Model} далее следует полагать $L := \max\left(L_0, L\right)$.

\floatname{algorithm}{Алгоритм}
	\begin{algorithm}
		\caption{Неускоренный градиентный метод с оракулом, использующим $(\delta, L, V)$-модель}\label{Grad_Model}
\begin{algorithmic}[1]
\REQUIRE $x_0$~--- начальная точка, $N$~--- количество шагов, а также постоянные $\delta$, $\widetilde{\delta}$ и $L_0 > 0$.
\STATE $L_{1} := {L_0}/{2}$
\FOR{$k=0, \, \dots, \, N-1$}
\STATE $\alpha_{k+1} := \frac{1}{L_{k+1}}$
\STATE \begin{equation}\label{equmir2DL_G}
\varphi_{k+1}(x)= \alpha_{k+1}\psi_{\delta}(x, x_{k}) + V(x, x_k),
\end{equation}
$$
x_{k+1} := \argmin_{x \in Q}^{\widetilde{\delta}}\varphi_{k+1}(x).
$$
\STATE {\bf If} \begin{equation}\label{exitLDL_G}
F_{\delta}(x_{k+1}) \leq F_{\delta}(x_{k}) + \psi_{\delta}(x_{k+1}, x_{k}) + L_{k+1}V(x_{k+1},x_{k}) + \delta,
\end{equation}
\STATE
{\bf then}
$
L_{k+2} := {L_{k+1}}/{2}, \; k:= k+1
$
and goto item 3.
\STATE  {\bf else} $L_{k+1} := 2\cdot L_{k+1}$ and goto item 5.
\ENDFOR
\RETURN $\bar{x}_N= \frac{1}{N}\sum_{k=0}^{N-1}x_{k+1}.$
\end{algorithmic}
\end{algorithm}

\begin{remark}
В простейшем случае евклидовой нормы и прокс-структуры, $Q = \mathbb{R}^n$, модели $\psi_{\delta} (x, y) = \langle \nabla F(y), x - y \rangle$ и $\tilde{\delta} = 0$ вспомогательная задача минимизации \eqref{equmir2DL_G} в листинге алгоритма \ref{Grad_Model} приводит к обычному шагу градиентного спуска:
$$
x_{k+1} := x_k - \alpha_{k+1} \nabla F(x_k) = x_k - \frac{1}{L_{k+1}} \nabla F(x_k).
$$
Если же $Q \neq \mathbb{R}^n$, \nk{то вспомогательная} задача минимизации \eqref{equmir2DL_G} сводится к задаче проектирования на множество $Q$:
$$
x_{k+1} := Pr_Q\{x_k - \alpha_{k+1} \nabla F(x_k)\} = Pr_Q \left\{x_k - \frac{1}{L_{k+1}} \nabla F(x_k)\right\}.
$$
\end{remark}

\begin{remark}
\leavevmode
\label{remark_maxmin}
Для всех $k \geq 0$ выполнено $	L_{k} \leq 2L$. При $k=0$ это верно ввиду предположения $L_0 \leq 2L$. При $k \geq 1$ ясно, \ag{что} выход из внутреннего цикла (на котором подбирается $L_k$) гарантирован при условии $L_{k} \geq 2L$. Выход из цикла гарантируется тем, что по условию существует $(\delta, L, V)$-модель для $F$ в любой точке $x \in Q$.
\end{remark}

\begin{exercise} Докажите, что количество обращений к вспомогательным задачам~\eqref{equmir2DL_G} в связи с изменением констант $L_k$ на итерациях будет увеличиваться по сравнению с неадаптивным вариантом градиентного метода в среднем \nk{не более чем} в постоянное число раз (меньшее трёх). Это означает, что адаптивность не приводит к существенному повышению стоимости итерации метода. В качестве указания см.~статью Ю.\,Е.~Нестерова \cite{Nesterov15b}.
\end{exercise}

Докажем важную лемму, которая необходима далее для вывода оценки скорости сходимости описываемого метода.

\begin{lemma}
	Пусть $\psi(x)$ --- выпуклая функция и при некотором $z \in Q$ верно 
	\begin{gather*}
	y=\argmin_{x \in Q}^{\widetilde{\delta}}\{\psi(x) + V(x,z)\}.
	\end{gather*}
	Тогда для всякого $x \in Q$ справедливо неравенство
	\begin{equation*}
	\psi(x) + V(x,z) \geq \psi(y) + V(y,z) + V(x,y) - \widetilde{\delta}.
	\end{equation*}
	\label{lemma_maxmin_2}
\end{lemma}
\begin{proof}
	Согласно определению~\ref{solNemirovskiy} существует такой вектор\\ $g \in \partial\psi(y)$, что:
	\begin{gather*}
	\langle g + \nabla_y V(y, z), x - y \rangle \geq -\widetilde{\delta} \text{  для всякого  } x \in Q.
	\end{gather*}
Тогда верно неравенство
	\begin{gather*}
		\psi(x) - \psi(y) \geq \langle g, x - y \rangle \geq \langle \nabla_y V(y, z), y - x \rangle - \widetilde{\delta},
	\end{gather*}
а также следующие равенства:
	\begin{gather*}
	\langle \nabla_y V(y, z), y - x \rangle=\langle \nabla d(y) - \nabla d(z), y - x \rangle=d(y) - d(z) - \langle \nabla d(z), y - z \rangle +\\ + d(x) - d(y) - \langle \nabla d(y), x - y \rangle - d(x) + d(z) + \langle \nabla d(z), x - z \rangle=\\=
	V(y,z) + V(x,y) - V(x,z),
	\end{gather*}
которые и завершают доказательство.
\end{proof}

Приведём ещё одно вспомогательное утверждение, описывающее поведение целевой функции для генерируемой алгоритмом \ref{Grad_Model} последовательности точек $\{x_k\}_{k=1}^{N}$.

\begin{lemma}
\label{lemma_maxmin_3DL_G}
Пусть $F$ допускает $(\delta, L, V)$-модель в произвольной точке множества $Q$ в смысле определения~{\rm \ref{gen_delta_L_model}}. Тогда для алгоритма~{\rm \ref{Grad_Model}} при любом $ x \in Q$ для всякого целого неотрицательного числа $k$ справедливо неравенство
	\begin{gather*}
		F(x_{k+1}) - F(x)\leq 2L(V(x, x_k) - V(x, x_{k+1}) + \widetilde{\delta}) + 2\delta.
	\end{gather*}
\end{lemma}
\begin{proof} Имеем следующие неравенства:
	\begin{gather*}
	F(x_{k+1}) \leqarg{(\ref{exitLDL_G}), (\ref{exitLDLOrig2})} F_{\delta}(x_{k}) + \psi_{\delta}(x_{k+1}, x_{k}) + L_{k+1}V(x_{k+1}, x_{k}) + 2\delta \leq \\ \leq
	F_{\delta}(x_{k}) + \psi_{\delta}(x_{k+1}, x_{k}) + \frac{1}{\alpha_{k+1}}V(x_{k+1}, x_{k}) + 2\delta
	\leq_{{\tiny \circled{1}}} \\\leq
	 F_{\delta}(x_{k}) + \psi_{\delta}(x,x_{k})
	 	 + \frac{1}{\alpha_{k+1}}V(x, x_k) - \frac{1}{\alpha_{k+1}}V(x, x_{k+1}) + \frac{\widetilde{\delta}}{\alpha_{k+1}} + 2\delta \leqarg{(\ref{exitLDLOrigV})} \\ \leq
	 F(x) + \frac{1}{\alpha_{k+1}}V(x, x_k) - \frac{1}{\alpha_{k+1}}V(x, x_{k+1}) + \frac{\widetilde{\delta}}{\alpha_{k+1}}+ 2\delta \leq_{{\tiny \circled{2}}} \\\leq
	 F(x) + 2LV(x, x_k) - 2LV(x, x_{k+1}) + 2L\widetilde{\delta} + 2\delta.
	\end{gather*}

{\small \circled{1}}~--- вытекает из леммы~\ref{lemma_maxmin_2} для
$\psi(x)=\alpha_{k+1} \psi_{\delta}(x, x_{k})$ c учётом левого неравенства в (\ref{exitLDLOrigV}).

{\small \circled{2}}~--- ввиду того, что $1/\alpha_{k+1}=L_{k+1} \leq 2L$.
\end{proof}

\begin{teo}
	\label{mainTheoremDL_G}
Пусть $F$ допускает $(\delta, L, V)$-модель в произвольной точке множества $Q$ в смысле определения~{\rm \ref{gen_delta_L_model}}, а также $V(x^*, x_0) \leq R^2$, где $x_0$~--- начальная точка, а $x^*$~--- ближайшая к $x_0$ точка минимума задачи \eqref{mainTask3} в смысле дивергенции Брэгмана, а также $L_0 \leq 2L$. Тогда для алгоритма~{\rm \ref{Grad_Model}} верно следующее неравенство:
	\begin{equation}\label{grad_estim_thm}
	F(\bar{x}_N) - F(x^*) \leq \frac{2LR^2}{N} + 2L\widetilde{\delta} + 2\delta .
	\end{equation}
\end{teo}
\begin{proof}
	Просуммируем неравенства из леммы\;\ref{lemma_maxmin_3DL_G} по ${k=0,\dots,}$ $N-1$: 
	\begin{gather*}
		\sum_{k=0}^{N-1}F(x_{k+1}) - N F(x) \leq 2LV(x, x_0) - 2LV(x, x_N) + 2NL\widetilde{\delta} + 2N\delta.
	\end{gather*}
	Выберем $x=x^*$. Тогда
	\begin{gather*}
		\sum_{k=0}^{N-1}F(x_{k+1}) - N F(x^*) \leq 2LR^2 - 2LV(x^*, x_N) + 2NL\widetilde{\delta} + 2N\delta.
	\end{gather*}
	Поскольку $V(x^*, x_N) \geq 0$, то верно неравенство
	\begin{gather*}
		\sum_{k=0}^{N-1}F(x_{k+1}) - N F(x^*) \leq 2LR^2 + 2NL\widetilde{\delta} + 2N\delta.
	\end{gather*}
	Разделим обе части последнего неравенства на $N$. Получим
	\begin{gather*}
		\frac{1}{N}\sum_{k=0}^{N-1}F(x_{k+1}) - F(x^*) \leq \frac{2LR^2}{N} + 2L\widetilde{\delta} + 2\delta.
	\end{gather*}
	Учитывая выпуклость $F$, для выхода $\bar{x}$ алгоритма \ref{Grad_Model} в итоге получаем, что
	\begin{gather*}
		F(\bar{x}_N) - F(x^*) \leq \frac{2LR^2}{N} + 2L\widetilde{\delta} + 2\delta.
	\end{gather*}
\end{proof}

\begin{exercise}
Выведите аналог предыдущей теоремы с сохранением оценки~\eqref{grad_estim_thm} в случае, если $F$ допускает ($\delta, L$)-модель в произвольной точке $Q$ согласно определению~{\rm \ref{gen_delta_L_oracle}}. \end{exercise}

\begin{exercise} Как будет выглядеть оценка~\eqref{grad_estim_thm}, если в ней вместо константы Липшица градиента $L$ использовать <<локальные>> константы $L_k$, полученные на итерациях работы метода?
\end{exercise}

\begin{exercise}
Покажите, что в случае выбора постоянного шага $\alpha_{k+1} = 1/L$ в алгоритме \ref{Grad_Model} для произвольного $k = 0, 1, 2, ...$ оценку \eqref{grad_estim_thm} можно уточнить:
	\begin{equation}\label{grad_estim_thm_const_step}
	F(\bar{x}_N) - F(x^*) \leq \frac{LR^2}{N} + L\widetilde{\delta} + 2\delta.
	\end{equation}
\end{exercise}

\begin{exercise}\label{Exer_Grad_Const_Step}
Покажите, что если $\delta = \tilde{\delta} = 0$ в алгоритме \ref{Grad_Model}, то
$F(x_{k+1}) \leq F(x_{k})$ для всякого $k \geq 0$. Покажите также, что в неравенствах \eqref{grad_estim_thm} и \eqref{grad_estim_thm_const_step} вместо $x^*$ можно выбрать произвольное $x \in Q$. Это означает, что при сделанных допущениях оценку \eqref{grad_estim_thm_const_step} можно заменить на следующую:
\begin{equation}\label{grad_estim_thm_const_exact}
	F(x_N) - F(x) \leq \frac{LV(x,x_0)}{N}
	\end{equation}
для произвольного $x \in Q$. Отметим, что далее оценка \eqref{grad_estim_thm_const_exact} будет использована в разделе \ref{rel_smooth}. 	
\end{exercise}

\subsection{Примеры сильно выпуклых прокс-функций}\label{subsect_prox_strongly}
Доказанная в предыдущем пункте оценка скорости сходимости неускоренного градиентного метода оптимальна на классе выпуклых относительно гладких задач. Если же рассмотреть класс выпуклых гладких задач (\ag{напомним, что} под гладкостью мы понимаем условие Липшица градиента), то оптимальные оценки дают ускоренные методы (см. выше пункты \ref{gradientTaylor}, \ref{AM}). Оказывается, что ускоренный градиентный метод можно ввести как для случая с необязательно евклидовой прокс-структурой, так и в модельной общности \cite{GasTurin}. В следующем пункте \ref{fastGradMethod} мы представим адаптивный метод такого типа (алгоритм \ref{Acceler_Grad_Model}). В отличие от неускоренного случая, здесь уже будет существенным предположение о $1$-сильной выпуклости прокс-функции $d$ относительно нормы $\|\cdot\|$. Поэтому мы скажем несколько слов о таких прокс-функциях.

Требование $1$-сильной выпуклости \ag{$d$} относительно нормы $\|\cdot\|$ влечёт по определению \ref{defin_bregman} неравенство
\begin{gather*}V(x,y) \geq \frac{1}{2}\|x - y\|^2.\end{gather*}
В случае евклидовой нормы для вектора $x = (x_1, x_2, ..., x_n) \in \mathbb{R}^{n}$
$$
\|x\|_2 = \sqrt{\sum\limits_{i = 1}^{n} x_i^2}
$$
$$
d(x) = \frac{1}{2}\|x\|_2^2, \quad V(x,y) = \frac{1}{2}\|x - y\|_2^2.$$

Известно \cite{Ben-Tal}, что дивергенцию Брэгмана, связанную с $p$-нормой при $1\leqslant p \leqslant 2$ (1-сильно выпуклую в $p$-норме) можно выбрать согласно табл. \ref{p0_tab1}. Например, для $1$-нормы можно положить 
\begin{equation}\label{eq_prox_a}
\widehat{d}(x) = \frac{e}{2(a-1)}\|x\|_a^2,
\end{equation}
где $a = 1 + \frac{O(1)}{\log n}$.

\medskip
\begin{table}[h!]
\centering
\caption{\centering\bf Примеры прокс-функций и дивергенций Брэгмана}
\label{p0_tab1}
\begin{tabular}{|c|c|c|c|}
\hline
$Q = B_p^n(1)$ & $1 \leq p \leq a$ & $a \leq p \leq 2$ & $2 \leq p \leq \infty$ \\
\hline
$\|\cdot\|$ & $\|\cdot\|_1$ & $\|\cdot\|_p$ & $\|\cdot\|_2$ \\
\hline
$d(x)$ & $d(x) = \frac{1}{2(a-1)}\|x\|_a^2$ & $d(x) = \frac{1}{2(p-1)}\|x\|_p^2$ & $\frac12 \|x\|_2^2$ \\
\hline
\end{tabular}
\end{table}
\medskip

Докажем, что $\hat{d}$ из \eqref{eq_prox_a} $1$-сильно выпукла относительно нормы $\|\cdot\|_1$. Для этого покажем, что функция
$$
d(x) = \frac{1}{2(a-1)}\|x\|_a^2 = \frac{\widehat{d}(x)}{e}
$$
$1$-сильно выпукла относительно нормы $\|\cdot\|_a$.
Будем использовать следующий известный факт (\cite{Rocaf_lemma}, предложения 12.54 и 12.60). {\it Непрерывно-дифференцируемая выпуклая функция $q: \mathbb{R}^n \rightarrow \mathbb{R}$ $\mu$-сильно выпукла относительно нормы $\|\cdot\|$ \ag{при некотором $\mu>0$} тогда и только тогда, когда двойственная ей функция {\rm(}результат преобразования Лежандра--Фенхеля{\rm):}
$$
W(y) = \max\limits_{x \in \mathbb{R}^n} \{\langle y, x \rangle  - q(x)\}
$$
непрерывно-дифференцируема и удовлетворяет следующему неравенству}{\rm :}
\begin{equation}\label{eq_V_eample_prox}
W(z) \leq W(y) + \langle \nabla W(y), z - y\rangle + \frac{1}{2\mu}{\|z - y\|_*^2}\; \text{  при любых }\; y, z \in \mathbb{R}^n.
\end{equation}
Сначала покажем сильную выпуклость $q(x) = \frac{1}{2}\|x\|_a^2$ относительно $\|\cdot\|_a$. При таком выборе $q$ имеем $W(y) = \frac{1}{2}\|y\|_b^2$, где $b = a/(a-1)$. Поскольку $W$ дважды непрерывно дифференцируема в любой внутренней точке области определения, то для доказательства \eqref{eq_V_eample_prox} достаточно проверить, что для всякого ненулевого вектора $r$ верно неравенство 
\begin{equation}\label{eq_ineq_prox_a}
\langle \nabla W(y) r, r \rangle \leq \frac{1}{a-1} \|r\|_b^2.
\end{equation}
В силу однородности $W$ имеем $\nabla^2 W(tr) = \nabla^2 W(r)$ при всяком $t > 0$. Поэтому, не уменьшая общности рассуждений, предположим $\|y\|_b = 1$. Обозначим через $\sign(t)$ знак числа $t$ ($\sign(t) = t/|t|$  при $t \neq 0$ и $\sign(0) = 0$). Имеем
\begin{equation*}
\langle \nabla W(y), r \rangle = \|y\|_b^{2-b} \sum_{i=1}^n |y|_i^{b-1} \sign(y_i) r_i,
\end{equation*}
$$
\langle \nabla^2 W(y) r, r \rangle = (2-b) \|y\|_b^{2-2b} \left ( \sum_{i=1}^n |y|_i^{b-1} \sign(y_i) r_i \right )^2 +
$$
$$
+(b-1) \|y\|_b^{2-b}  \sum_{i=1}^n |y|_i^{b-2} r_i^2 \leq (b-1) \sum_{i=1}^n |y|_i^{b-2} r_i^2,
$$
поскольку $b  \geq 2$. Наконец, в силу неравенства Гёльдера и $ \|y \|_b = 1$ получаем, что
$$
(b-1) \sum_{i=1}^n |y|_i^{b-2} r_i^2 \leq (b-1) \left ( \sum_{i} |y_i|^b \right )^{\frac{b-2}{b}} \left ( \sum_{i} |r|_i^b \right )^{\frac{2}{b}} \leq (b-1) \|r\|_b^2,
$$
откуда и следует \eqref{eq_ineq_prox_a} ввиду $b-1 = \frac{1}{a-1}$. Итак, для произвольных $x, y \in \mathbb{R}^n$ верно
\begin{equation}\label{eq_prox_a_1norm}
d(x) \geq d(y) + \langle \nabla d(y), x - y \rangle + \frac{1}{2} \|x - y\|_a^2.
\end{equation}
Согласно неравенству Гёльдера, 
$$\|x\|_1 \leq \|x\|_{a}n^{\frac{a-1}{a}} \text{ при всяком } x \in \mathbb{R}^n.$$
Если $a = 1 + \frac{C}{\log n}$ для некоторого $C > 0$, то 
$$\log \left(n^{\frac{a-1}{a}}\right) \leq \frac{C\log n}{C\log n + 1} \leq 1.$$

Поэтому при $a = 1 + \frac{O(1)}{\log n}$ для всякого $x$ верно $\|x\|_a \geq \frac{\|x\|_1}{e}$.  Это означает, что из неравенства \eqref{eq_prox_a_1norm} следует ($1/e$)-сильная выпуклость $d$ относительно уже $1$-нормы. Поэтому функция $\widehat{d}(x) = ed(x)$ из \eqref{eq_prox_a} $1$-сильно выпукла относительно $1$-нормы.

Отметим, что далее в п.~\ref{comp_prox} подробно разбирается пример задачи, часто возникающей на итерациях многих современных численных методов, для которой естественен выбор прокс-структуры указанного типа при $a = \frac{2\log n}{2\log n - 1}$ .

Отметим\ag{, что} в связи с приложениями в машинном обучении часто используют следующую прокс-структуру (дивергенция или расхождение Кульбака--Лейблера) 
\begin{equation}\label{eq21}
d(x)=\sum_{i=1}^nx_i\ln x_i,\;V(x,y)=\sum_{i=1}^nx_i\log\left(\frac{x_i}{y_i}\right).
\end{equation}
Такая прокс-структура отвечает $1$-норме на единичном (вероятностном) симплексе 
$$
Q=S_n \left( 1 \right)=\left\{ {x\in \mathbb{R}_+^n:\;\;\sum\limits_{i=1}^n{x_i } =1} \right\}.
$$
Использование такого прокса позволяет выписывать явные формулы проектирования на симплекс, что удобно для задач на симплексах. Действительно, при выборе указанной прокс-структуры шаг \eqref{eq_bregman_step} имеет вид ($i = 1,2, ..., n$): 
$$
x^i_{k+1} =\frac{\exp \left( {-h\sum\limits_{r=0}^k {\nabla_i f\left( {x_r}\right)} } \right)}{\sum\limits_{l=1}^n {\exp \left( {-h\sum\limits_{r=0}^k{\nabla_l f\left( {x_r} \right)}}\right)} }=\frac{x^i_k \exp \left({-h\nabla_i f\left( {x_k} \right)} \right)}{\sum\limits_{l=1}^n {x^l_k \exp\left( {-h\nabla _l f\left( {x_k} \right)} \right)}},
$$
где $x_k = (x^1_k, x^2_k, ..., x^n_k) \in \mathbb{R}^n$, $x_{k+1} = (x^1_{k+1}, x^2_{k+1}, ..., x^n_{k+1}) \in \mathbb{R}^n$. При этом принято выбирать начальную точку $(1/n, 1/n, ..., 1/n)$, чтобы минимизировать максимально возможное удаление от начальной точки до точного решения задачи. А такое удаление в случае неевклидовых норм определяется через введённую дивергенцию (или расхождение) Брэгмана, которая в общем случае заменяет квадрат <<евклидова>> расстояния. Как видно по формулировкам теорем \ref{mainTheoremDL_G} и \ref{mainTheoremDL}, такое <<расстояние>> вида $V(x^*, x_0)$ важно для оценок скорости сходимости численных методов градиентного типа. 

Все приведённые в табл. \ref{p0_tab1} прокс-функции $1$-сильно выпуклы в указанных нормах на всём пространстве. Поэтому их можно использовать с неускоренными и ускоренными градиентными методами и в том случае, когда мы заранее не знаем, где локализовано решение. Например, если стартовать из нулевой точки и заранее знать, что решение разреженно (большинство компонент равно нулю), то оказывается, что естественно выбирать $1$-норму и соответствующую ей прокс-функцию в табл. \ref{p0_tab1}. В таком случае можно ожидать, что $R^2 = V(x^*, x_0)$ не существенно будет зависеть от выбора нормы ($1$-нормы вместо евклидовой). При этом от перехода от $2$-номы к $1$-норме может быть
существенная выгода в оценке константы $L$ в формулировке теорем \ref{mainTheoremDL_G} и \ref{mainTheoremDL}. Действительно, константу
Липшица градиента 
$$
\|\nabla F(x) - \nabla F(y)\|_* \leq L\|x - y\| \; \text{ при всех } x, y \in Q,
$$
если $F$ имеет не только градиент, но и гессиан $\nabla^2 F(x)$, можно определять следующим образом: 
$$
L:= L_p = \sup\limits_{x \in Q, \; \|h\|_p \leq 1} \langle h, \nabla^2 F(x)h \rangle.
$$
Несложно привести пример, когда $L$ будет существенно зависеть от выбора $1$-нормы или $2$-нормы. Он связан с тем, что $L_2/n \leq L_1 \leq L_2$. Например, можно рассмотреть конкретный пример:   
$$F(x) = \frac{1}{2} \langle Ax, x \rangle - \langle b, x\rangle$$
для фиксированной квадратной матрицы $A =(a_{ij})$ размера $n \times n$. 
Тогда
$L_1 = \max_{i, j = 1, 2, ..., n} |a_{ij}|$, а $L_2 = \lambda_{\max}(A)$ --- максимальное собственное значение~$A$. Это означает, что использование градиентных процедур с <<неевклидовыми>> дивергенциями Брэгмана и итерациями вида \eqref{eq_bregman_step} может быть целесообразно не только для относительно гладких задач, но и в случае $1$-сильно выпуклой прокс-функции. При этом $1$-сильная выпуклость прокс-функции позволяет строить ускоренные методы (см. п. \ref{fastGradMethod}).
\begin{exercise}
Приведите пример матрицы $A$, для которой верно равенство $L_2 = n\cdot L_1$. 
\end{exercise}

\subsection[Ускоренный градиентный метод с оракулом, который использует $(\delta, L)$-модель в запрашиваемой точке]{Ускоренный градиентный метод с оракулом,\\ который использует $(\delta, L)$-модель\\ в запрашиваемой точке}
\label{fastGradMethod}

В данном пункте рассмотрим ускоренный вариант градиентного метода в модельной общности, который базируется на так называемом {\it методе подобных треугольников} \cite{GasNestSimTr}. Этот метод получил название по соответствующей геометрической интерпретации (см. рис. \ref{FigureSimTriangle}). Оптимальная оценка скорости сходимости для задачи минимизации
$$
\min_{x\in Q} F(x)
$$
выпуклой функции $F$ с $L$-липшицевым градиентом 
$$
\|\nabla F(x) - \nabla F(y)\|_* \leq L\|x - y\| \; \text{ при всех } x, y \in Q
$$
\noindent для метода подобных треугольников  получается за счёт более тонкой организации итерации по сравнению с алгоритмом~\ref{Grad_Model}. Выпишем листинг метода подобных треугольников в случае $Q = \mathbb{R}^n$ и евклидовой нормы (см.~алгоритм~\ref{Sim_Triangle}). 

\floatname{algorithm}{Алгоритм}
	\begin{algorithm}
		\caption{Метод подобных треугольников}\label{Sim_Triangle}
\begin{algorithmic}[1]
\REQUIRE $x_0$~--- начальная точка, $N$~--- количество шагов, константа Липшица градиента целевой функции $L > 0$.
\STATE $y_0 := x_0.$
\STATE $u_0 := x_0.$
\STATE $\alpha_0 := 0,\,
A_0 := \alpha_0$
\FOR{$k=0, \, \dots, \, N-1$}
\STATE 	
$$\alpha_{k+1}:= \frac{1 + \sqrt{1+4A_kL}}{2L} \text{ --- наибольший корень квадратного уравнения} $$ 
$$L\alpha^2_{k+1} -  \alpha_{k+1} - A_k = 0,$$
$$A_{k+1}=A_k +\alpha_{k+1}.$$
\STATE
$$
y_{k+1} := \frac{\alpha_{k+1} u_k + A_k x_k}{ A_{k+1}},
$$
$$
u_{k+1} := u_k - \alpha_{k+1} \nabla F(y_{k+1}),
$$
$$
x_{k+1} := \frac{\alpha_{k+1}u_{k+1}+A_kx_k}{A_{k+1}}.
$$
\ENDFOR
\RETURN $x_N.$
\end{algorithmic}
\end{algorithm}

Как видим, каждая итерация алгоритма \ref{Sim_Triangle} состоит из нескольких шагов, связанных со вспомогательными последовательностями $u_k$ и $y_k$. Для наглядности на рис.~\ref{FigureSimTriangle} приведём некоторую геометрическую иллюстрацию метода подобных треугольников (аналогичную картинку можно найти в обзорной статье \cite{DvurechenskyFirstOrder2021}). На рис.~\ref{FigureSimTriangle} каждой из точек $A, B, C, D, E$ соответствует используемый на $k$-й итерации метода подобных треугольников элемент основной последовательности $x_k$, либо вспомогательных последовательностей $y_k$ и $u_k$. Эти последовательности нужны для описания <<составного>> шага ускоренного метода, позволяющего получить оптимальную оценку скорости сходимости для выпуклых гладких минимизационных задач. Отрезкам $EC$, $CA$, $ED$ и $BD$ соответствуют векторы $y_{k+1} - x_k$, $u_k - y_{k+1}$, $x_{k+1} - x_k$ и $u_{k+1} - x_{k+1}$.
Можно записать равенства
$$\frac{EC}{AC} = \frac{ED}{BD} = \frac{A_k}{\alpha_{k+1}},$$


\noindent поскольку (см. алгоритм \ref{Sim_Triangle}) верны соотношения
$$y_{k+1}= \frac{\alpha_{k+1}}{A_{k+1}} u_k + \frac{A_k}{A_{k+1}} x_k, \text{  и  } x_{k+1} = \frac{\alpha_{k+1}}{A_{k+1}} u_{k+1} + \frac{A_k}{A_{k+1}} x_k.$$
Поэтому получаем подобные треугольники $\triangle DEC$ и $\triangle BEA$, что и поясняет название {\it метод подобных треугольников}. Из этого подобия получаем соотношение
$$x_{k+1} - y_{k+1} = \frac{\alpha_{k+1}}{A_{k+1}} (u_{k+1} - u_k),$$
поскольку на рис.~\ref{FigureSimTriangle} $\overline{CD} \parallel \overline{AB}$.

\begin{figure}[h!]
\begin{center}
\begin{tikzpicture}[line width=1pt]
\draw (1,4)--(7,2)--(7,6)--(1,4);
\draw (4, 3)--(4,5);
\node[below left] at (1,4) {$x_k (E)$};
\node[below right] at (7,2) {$u_{k+1} (B)$};
\node[below right] at (7,6) {$u_k (A)$};
\node[below] at (4,6) {$y_{k+1} (C)$};
\node[below] at (4,2.8) {$x_{k+1} (D)$};
\end{tikzpicture}
\end{center}

\caption{Геометрическая иллюстрация алгоритма \ref{Sim_Triangle}}\label{FigureSimTriangle}
\end{figure}

Теперь переходим непосредственно к описанию адаптивного ускоренного градиентного метода (алгоритм \ref{Acceler_Grad_Model}) в модельной общности, который использует на итерациях $(\delta, L)$-модель (см. определение \ref{gen_delta_L_model}) в произвольной запрашиваемой точке. В отличие от неускоренного алгоритма \ref{Grad_Model} для ускоренного метода с итерациями вида \eqref{eq_bregman_step} уже существенно предположение о $1$-сильной выпуклости используемой прокс-функции $d$:
$$
V(y, x) \geqslant \frac{1}{2} \|y - x\|^2 \text{  для произвольных  } x, y \in Q.
$$
Напомним, что сильно выпуклые прокс-структуры относительно нормы обсуждались выше в п.~\ref{subsect_prox_strongly}.

Для большей наглядности  доказательства и лучших оценок (см. \cite{DragomirRelativeSmooth}) мы ограничимся далее в этом пункте общностью определения \ref{gen_delta_L_oracle}, а не определения \ref{gen_delta_L_model}.


\floatname{algorithm}{Алгоритм}
	\begin{algorithm}
		\caption{Ускоренный градиентный метод (аналог метода подобных треугольников) с оракулом, использующим на итерациях $(\delta, L)$-модель}\label{Acceler_Grad_Model}
\begin{algorithmic}[1]
\REQUIRE $x_0$~--- начальная точка, $N$~--- количество шагов, а также постоянные $\delta$, $\widetilde{\delta}$, $L_0 > 0$.
\STATE $y_0 := x_0.$
\STATE $u_0 := x_0.$
\STATE $L_1 := {L_0}/{2},$
\STATE $\alpha_0 := 0,\,
A_0 := \alpha_0$
\FOR{$k=0, \, \dots, \, N-1$}
\STATE  	$$\alpha_{k+1}:= \frac{1 + \sqrt{1+4A_kL_{k+1}}}{2L_{k+1}} \text{ --- наибольший корень} $$ 

$$\text{ квадратного уравнения } L_{k+1}\alpha^2_{k+1} -  \alpha_{k+1} - A_k = 0.$$ 

\STATE $$A_{k+1} := A_k + \alpha_{k+1},$$
\begin{equation}\label{eqymir2DL}
y_{k+1} := \frac{\alpha_{k+1}u_k + A_k x_k}{A_{k+1}},
\end{equation}
\begin{equation}\label{equmir2DL}
\varphi_{k+1}(x)=\alpha_{k+1}\psi_{\delta}(x, y_{k+1}) + V(x, u_k), \; 
u_{k+1} := \argmin_{x \in Q}^{\widetilde{\delta}}\varphi_{k+1}(x),
\end{equation}
$$
x_{k+1} := \frac{\alpha_{k+1}u_{k+1}+A_kx_k}{A_{k+1}}.
$$
\STATE {\bf If} \begin{equation}\label{exitLDL}
F_{\delta}(x_{k+1}) \leq F_{\delta}(y_{k+1}) + \psi_{\delta}(x_{k+1}, y_{k+1}) +\frac{L_{k+1}}{2}\|x_{k+1} - y_{k+1}\|^2 + \delta,
\end{equation}
\STATE {\bf then} $L_{k+2} := {L_{k+1}}/{2}, \; k:= k+1$
and goto item 6.
\STATE {\bf else} $L_{k+1} := 2\cdot L_{k+1}$ and goto item 8.
\ENDFOR
\RETURN $x_N.$
\end{algorithmic}
\end{algorithm}


Описание вывода оценки скорости сходимости алгоритма \ref{Acceler_Grad_Model} начнём со следующей вспомогательной леммы.

\begin{lemma}
	\label{lemma_maxmin_1}
	Пусть для используемой на итерациях алгоритма~{\rm \ref{Acceler_Grad_Model}} последовательности $\alpha_k$ имеют место равенства
	\begin{align*}
	\alpha_0=0,\,\,\,
	A_k=\sum_{i=0}^{k}\alpha_i,\,\,\,
	A_k=L_{k}\alpha_k^2,\,\,\,
	\end{align*}
	где $L_k \leq 2L$ для всякого $ k\geq0$ {\rm(}см.~замечание~{\rm \ref{remark_maxmin})}.
	Тогда для произвольного $k \geq 1$ верно неравенство
	 \begin{align}
	 \label{lemma_maxmin_1_1}
	 A_k \geq \frac{(k+1)^2}{8L}.
	 \end{align}
\end{lemma}
\begin{proof}
Проведём доказательство методом математической индукции по $k$. Если $k=1$, то имеем
	\begin{equation*}
	\alpha_1=L_{1}\alpha_1^2
	\end{equation*}
и поэтому
	\begin{equation*}
	A_1=\alpha_1=\frac{1}{L_1} \geq \frac{1}{2L}.
	\end{equation*}
Если же $k \geq 2$, то
	\begin{equation*}
	L_{k+1}\alpha^2_{k+1}=A_{k+1} \;\Leftrightarrow \;
	L_{k+1}\alpha^2_{k+1}=A_{k} + \alpha_{k+1} \;\Leftrightarrow\;
	L_{k+1}\alpha^2_{k+1} - \alpha_{k+1} - A_{k}=0.
	\end{equation*}
Решая данное квадратное уравнение, выбираем наибольший его корень, т.е.
	\begin{equation*}
	\alpha_{k+1}=\frac{1 + \sqrt{1 + 4L_{k+1}A_{k}}}{2L_{k+1}}.
	\end{equation*}
По индукции, допустим, что неравенство (\ref{lemma_maxmin_1_1}) верно для некоторого натурального числа $k$. Тогда получим неравенство (\ref{lemma_maxmin_1_1}) для $k+1$:
	\begin{gather*}
	\alpha_{k+1}=\frac{1}{2L_{k+1}} + \sqrt{\frac{1}{4L_{k+1}^2} + \frac{A_{k}}{L_{k+1}}} \geq
	\frac{1}{2L_{k+1}} + \sqrt{\frac{A_{k}}{L_{k+1}}} \geq \\\geq
	\frac{1}{4L} + \frac{1}{\sqrt{2L}}\frac{k+1}{2\sqrt{2L}} =
	\frac{k+2}{4L}.
	\end{gather*}
Отметим, что последнее неравенство следует из индуктивного предположения. В итоге получаем, что
	\begin{equation*}
	\alpha_{k+1} \geq \frac{k+2}{4L}, \text{ а также }
	A_{k+1}=A_k + \alpha_{k+1}=\frac{(k+1)^2}{8L} + \frac{k+2}{4L} \geq \frac{(k+2)^2}{8L},
	\end{equation*}
что и завершает доказательство леммы.
\end{proof}

Докажем базовую лемму, необходимую для вывода оценки скорости сходимости алгоритма \ref{Acceler_Grad_Model}.

\begin{lemma}
\label{lemma_maxmin_3DL} Пусть $F$ допускает $(\delta, L)$-модель в произвольной точке множества $Q$ в смысле определения~{\rm \ref{gen_delta_L_oracle}}. Тогда
для алгоритма~{\rm \ref{Acceler_Grad_Model}} при всяком $x \in Q$ справедливо неравенство
	\begin{equation*}
		A_{k+1} F(x_{k+1}) - A_{k} F(x_{k}) + V(x, u_{k+1}) - V(x, u_{k}) \leq \alpha_{k+1}F(x) + 2\delta A_{k+1} + \widetilde{\delta}.
	\end{equation*}
\end{lemma}
\begin{proof} Имеем следующую цепочку неравенств:
	\begin{gather*}
	F(x_{k+1}) \leqarg{(\ref{exitLDL}), (\ref{exitLDLOrig2})} F_{\delta}(y_{k+1}) + \psi_{\delta}(x_{k+1},y_{k+1}) + \frac{L_{k+1}}{2} \left\|x_{k+1} - y_{k+1}\right\|^2 + 2\delta= 
	\\=
	F_{\delta}(y_{k+1}) + \psi_{\delta}\left(\frac{\alpha_{k+1}u_{k+1} + A_k x_k}{A_{k+1}},y_{k+1}\right) +\\+ \frac{L_{k+1}}{2}\left\|\frac{\alpha_{k+1}u_{k+1} + A_k x_k}{A_{k+1}} - y_{k+1}\right\|^2 + 2\delta \leqarg{Опр.~\ref{gen_delta_L_oracle}, (\ref{eqymir2DL})} \\
	\leq
	F_{\delta}(y_{k+1}) +
	\frac{\alpha_{k+1}}{A_{k+1}}\psi_{\delta}(u_{k+1}, y_{k+1}) + \\+
	 \frac{A_k}{A_{k+1}}\psi_{\delta}(x_k, y_{k+1}) + \frac{L_{k+1} \alpha^2_{k+1}}{2 A^2_{k+1}}\|u_{k+1} - u_k\|^2 + 2\delta= \\=
	 \frac{A_k}{A_{k+1}}(F_{\delta}(y_{k+1}) + \psi_{\delta}(x_k, y_{k+1}))
	 + \\+
	 \frac{\alpha_{k+1}}{A_{k+1}}(F_{\delta}(y_{k+1}) +
	 \psi_{\delta}(u_{k+1}, y_{k+1}))+ \\
	 + \frac{L_{k+1} \alpha^2_{k+1}}{2 A^2_{k+1}}\|u_{k+1} - u_k\|^2 + 2\delta=_{{\tiny \circled{1}}} \\ =
	 \frac{A_k}{A_{k+1}}(F_{\delta}(y_{k+1}) + \psi_{\delta}(x_k,y_{k+1}))
	 + \\+
	 \frac{\alpha_{k+1}}{A_{k+1}}(F_{\delta}(y_{k+1}) + \psi_{\delta}(u_{k+1},y_{k+1})
	 + \frac{1}{2 \alpha_{k+1}}\|u_{k+1} - u_k\|^2) + 2\delta\leq_{{\tiny \circled{2}}}  \\\leq
	 \frac{A_k}{A_{k+1}}(F_{\delta}(y_{k+1}) + \psi_{\delta}(x_k,y_{k+1}))
	 + \\+
	 \frac{\alpha_{k+1}}{A_{k+1}}\left(F_{\delta}(y_{k+1}) + \psi_{\delta}(u_{k+1},y_{k+1})
	 + \frac{1}{\alpha_{k+1}}V(u_{k+1}, u_k)\right) + 2\delta\leq_{{\tiny \circled{3}}} \\\leq
	 \frac{A_k}{A_{k+1}} F(x_k)\, + \\+
	 \frac{\alpha_{k+1}}{A_{k+1}}\!\left(F_{\delta}(y_{k+1})\!+\!\psi_{\delta}(x,y_{k+1})
	 \!+\!\frac{1}{\alpha_{k+1}}V(x, u_k)\!-\!\frac{1}{\alpha_{k+1}}V(x, u_{k+1})\!+\! \frac{\widetilde{\delta}}{\alpha_{k+1}}\right)\!+\!2\delta\!\!\leqarg{(\ref{exitLDLOrigV})} \\ \leq
	 \frac{A_k}{A_{k+1}} F(x_k) +
	 \frac{\alpha_{k+1}}{A_{k+1}} F(x)
	 + \frac{1}{A_{k+1}}V(x, u_k) - \frac{1}{A_{k+1}}V(x, u_{k+1}) + 2\delta + \frac{\widetilde{\delta}}{A_{k+1}}.
	\end{gather*}
	
	{\small \circled{1}}~--- вытекает из равенства $A_k=L_{k}\alpha^2_k$;
	
	{\small \circled{2}~--- вытекает из $V(u_{k+1}, u_k) \geqslant \frac{1}{2} \|u_{k+1} - u_k\|^2$;}
	
	{\small \circled{3}}~--- следствие леммы~\ref{lemma_maxmin_2} для функции
	$\psi(x)=\alpha_{k+1} \psi_{\delta}(x, y_{k+1})$ и неравенства слева в  (\ref{exitLDLOrigV}).
\end{proof}

Перейдём теперь к основному результату настоящего пункта. 

\begin{teo}
\label{mainTheoremDL}
Пусть $F$ допускает $(\delta, L)$-модель в произвольной точке множества $Q$ в смысле определения~{\rm\ref{gen_delta_L_oracle}}, причём $V(x^*, x_0) \leq R^2$, где $x^*$~--- ближайшая к начальной точке $x_0$ в смысле дивергенции Брэгмана точка минимума $F$, а также $L_0 \leq 2L$. Тогда для алгоритма~{\rm\ref{Acceler_Grad_Model}} выполнены следующие неравенства:
	\begin{gather*}	
		F(x_N) - F(x^*) \leq \frac{R^2}{A_{N}} + \frac{2\delta \sum\limits_{k=0}^{N-1} A_{k+1}}{A_{N}} + \frac{N \widetilde{\delta}}{A_{N}} \leq \\
	\leq \frac{8LR^2}{(N+1)^2} + \frac{8L\widetilde{\delta}}{N+1} + 2N\delta.
	\end{gather*}
\end{teo}
\begin{proof}
	Просуммируем неравенства из леммы~\ref{lemma_maxmin_3DL} по $k=0,\dots,$ $N - 1$:
	\begin{gather*}
		A_{N} F(x_N) - A_{0} F(x_0) + V(x, u_N) - V(x, u_0) \leq (A_N - A_0)F(x) + 2\delta \sum\limits_{k=0}^{N-1} A_{k+1}+ \\ +\, N \widetilde{\delta} \Leftrightarrow
		A_{N} F(x_N) + V(x, u_N) - V(x, u_0) \leq A_NF(x) + 2\delta \sum\limits_{k=0}^{N-1} A_{k+1} + N \widetilde{\delta}.
	\end{gather*}
	Выберем теперь $x=x^*$. В таком случае
	\begin{gather*}
		A_{N} (F(x_N) - F(x^*) \leq R^2 + 2\delta \sum\limits_{k=0}^{N-1} A_{k+1} + N \widetilde{\delta}.
	\end{gather*}
	Разделив обе части предыдущего неравенства на $A_N$, получаем, что
	\begin{gather*}	
		F(x_N) - F(x^*) \leq \frac{R^2}{A_{N}} + \frac{2\delta \sum\limits_{k=0}^{N-1} A_{k+1}}{A_{N}} + \frac{N \widetilde{\delta}}{A_{N}} \leq_{{\tiny \circled{1}}}\\\leq \frac{8LR^2}{(N+1)^2} + \frac{2\delta \sum\limits_{k=0}^{N-1} A_{k+1}}{A_{N}} + \frac{8L\widetilde{\delta}}{N+1} \leq \\
	\leq \frac{8LR^2}{(N+1)^2} + \frac{8L\widetilde{\delta}}{N+1} + 2N\delta.
	\end{gather*}
	{\small \circled{1}}~--- следствие леммы~\ref{lemma_maxmin_1}.
\end{proof}

\subsection[Некоторые следствия из описанных теоретических результатов о сходимости градиентных методов\\ в модельной общности]{Некоторые следствия из описанных теоретических\\ результатов о сходимости градиентных методов\\ в~модельной общности}
\label{sledviya}

\subsubsection[Оценка сложности для ускоренного градиентного\\ метода\,c\,использованием\,точной\,информации на итерациях]{Оценка сложности для ускоренного градиентного метода\\ c использованием точной информации на итерациях}
Предположим, что $F$~--- выпуклая функция с $L$-липшицевым градиентом относительно нормы $\|\cdot\|$. Тогда верно следующее неравенство:
\begin{gather*}
0 \leq F(x) - F(y) - \langle\nabla F(y), x - y \rangle \leq \frac{L}{2}\|x - y\|^2 \,\,\, \text{  для всяких  } x,y \in Q.
\end{gather*}
Рассмотрим ($\delta, L$)-модель минимизириуемой функции $F$ вида
\begin{equation}\label{eq_stand_model}
\psi_{\delta}(x,y)=\langle\nabla F(y), x - y \rangle, \quad F_{\delta}(y)=F(y)
\end{equation}
и $\delta=0$, причём мы предполагаем, что промежуточную задачу можно решать точно, то есть $\widetilde{\delta}=0$. Получаем следующую оценку скорости сходимости для ускоренного варианта градиентного метода (п.~\ref{fastGradMethod} и теорема~\ref{mainTheoremDL}):
\begin{equation}\label{OracleEstimAcceler}
F(x_N) - F(x^*) \leq \frac{8LR^2}{(N+1)^2}.
\end{equation}
Данная оценка скорости сходимости является оптимальной с точностью до умножения на константу. В случае точного оракула ($\delta = 0$)  ускоренные градиентные методы (в нашем случае это алгоритм \ref{Acceler_Grad_Model}) оптимальны на классе задач минимизации выпуклых функций с $L$-липшицевым градиентом (см. пп. \ref{sect-black-box}, \ref{gradientTaylor}). При этом известно, что лучше константы $\simeq 1$ в такой оценке получить нельзя \cite{TaylorHendrickxGlineur17b}, в то время для алгоритма \ref{Acceler_Grad_Model} мы имеем константу 8. На оптимальность ускоренных методов для задач выпуклой минимизации указывает то, что точность $\varepsilon$ по функции \eqref{p0_eq1_1} приближённого решения задачи выпуклой гладкой минимизации достигается за $O\left(\sqrt{\frac{LR^2}{\varepsilon}}\right)$ итераций (см. \eqref{OracleEstimAcceler}). В то же время для неускоренных градиентных методов необходимо $O\left(\frac{LR^2}{\varepsilon}\right)$ итераций для достижения точности $\varepsilon$ по функции приближённого решения задачи выпуклой гладкой минимизации.

Итак, в случае отсутствия погрешностей оракула ($\delta = 0$) и решения вспомогательных подзадач на итерациях ($\tilde{\delta} = 0$) ускоренный градиентный метод на классе задач гладкой выпуклой минимизации даёт лучшие оценки скорости сходимости по сравнению с неускоренным методом. Интересно, что ситуация усложняется при рассмотрении неточного оракула. В отличие от неускоренного градиентного метода, ускоренный вариант может страдать от накопления ошибок. 

\subsubsection[Сравнение оценок для ускоренного и неускоренного методов при наличии погрешностей]{Сравнение оценок для ускоренного и неускоренного методов\\ при наличии погрешностей}\label{subsec_inexact} Предположим, что для целевая функция задачи (\ref{mainTask3}) допускает $(\delta, \, L)$-модель в произвольной точке. Будем считать, что промежуточную задачу в смысле определения~\ref{solNemirovskiy} мы можем решать на каждом шаге с точностью $\widetilde{\delta}$, тогда ввиду теорем~\ref{mainTheoremDL_G} и~\ref{mainTheoremDL} верны следующие неравенства:
\begin{equation}\label{estim_mist_or1}
F(\bar{x}_N) - F(x^*) \leq \frac{2LR^2}{N} + 2L\widetilde{\delta} + 2\delta,
\end{equation}
\begin{equation}\label{estim_mist_or2}
F(x_N) - F(x^*) \leq \frac{8LR^2}{(N+1)^2} + \frac{8L\widetilde{\delta}}{N+1} + 2N\delta,
\end{equation}
где $\bar{x}_N$~--- выход работы неускоренного метода (см. теорему~\ref{mainTheoremDL_G}), а $x_N$~--- выход работы ускоренного градиентного метода (теорема~\ref{mainTheoremDL}). Таким образом, мы видим, что ускоренный метод более устойчив к ошибкам $\widetilde{\delta}$, возникающих при решении промежуточных задач (на самом деле, это является некоторым артефактом определения $\widetilde{\delta}$, поскольку решить вспомогательную задачу \eqref{equmir2DL} для ускоренного метода с той же точностью $\widetilde{\delta}$, что и для неускоренного, всё сложнее и сложнее по мере роста номера итерации, поскольку обусловленность этой задачи растёт вместе с ростом $\alpha_{k+1}$), но при этом его оценка скорости сходимости может накапливать соответствующую $\delta$ величину, т.е. возникающую погрешность при вызове $(\delta, \, L)$-модели. Неускоренный же метод гарантирует отсутствие накопления этого параметра в оценке скорости сходимости, но при $\delta = 0$ (или малом $\delta$) оценка скорости сходимости лучше для ускоренного метода.

\subsubsection{Универсальный ускоренный градиентный метод}
\label{Univers_Meth}
Теперь покажем, как с использованием градиентных методов в модельной общности можно пояснять идеологию так называемых {\it универсальных} методов для задач выпуклой оптимизации, предложенную недавно Ю.\,Е.\,Нестеровым в \cite{Nesterov15b}. Здесь имеется в виду, что для градиента целевой функции (субградиента при $\nu = 0$) выполнено условие Гёльдера: существует $\nu\in[0,1]$ такое, что
\begin{equation}\label{eq_hold_condit}
\|\nabla F(x) - \nabla F(y)\|_* \leq L_\nu\|x - y\|^\nu \text{ для всяких } x,y \in Q.
\end{equation}
В частности, при $\nu = 0$ таким свойством могут обладать также и негладкие функции. Суть метода Ю.\,Е.~Нестерова заключается как на довольно естественной идее адаптивной настройки на константу Липшица (Гёльдера) (суб)градиента функции $F$, так и на довольно нетривиальной возможности настройки метода на подходящий уровень гладкости задачи. В частности, как показывают эксперименты, универсальные варианты Ю.\,Е.~Нестерова обычного и быстрого градиентного методов за счёт этих двух трюков могут реально работать лучше оптимальных теоретических оценок для некоторых задач. Например, для задачи Ферма--Торричелли--Штейнера, целевая функция которой удовлетворяет условию \eqref{eq_hold_condit} при $\nu = 0$, реальная скорость сходимости \ag{при практической реализации метода} может оказаться $O(1/\varepsilon)$~\cite{Nesterov15b}. Это существенно лучше гарантированной общей теорией \ag{оценкой} $O\left(1/\varepsilon^2\right)$ (см.~раздел~\ref{sect-black-box}). \ag{Про задачу Ферма--Торричелли--Штейнера и её связь с проблемой построения оптимальных сетей дорог можно почитать, например, в брошюре~\cite{Protasov}. Известны также и обобщения постановки задачи Ферма--Торричелли--Штейнера о минимизации суммы расстояний от точки до нескольких выпуклых замкнутых множеств (в частности, шаров) \cite{Mordukh}. Такие задачи приводят уже к бесконечному количеству точек недифференцируемости целевой функции, но универсальный метод также может экспериментально работать с оценкой скорости сходимости $O\left(1/\varepsilon\right)$.} Аналогичная ситуация имеет место в задаче поиска равновесного распределения потоков в транспортных сетях~\cite{Baimurzina}.

К выделенному классу задач выпуклой оптимизации можно применять рассмотренные выше в пунктах \ref{gradMethod} и \ref{fastGradMethod} адаптивные градиентные методы в модельной общности. Действительно, для всякой выпуклой функции $F$ с $\nu$-гёльдеровым градиентом (или субградиентом при $\nu = 0$) верно неравенство \cite{Nesterov15b}:
\begin{gather}
0 \leq F(x) - F(y) - \langle\nabla F(y), x - y \rangle \leq \frac{L(\delta)}{2}\|x - y\|^2 + \delta \,\,\, \text{ при всех } x,y \in Q,
\end{gather}
где \begin{gather*}L(\delta) = L_\nu\left[\frac{L_\nu}{2\delta}\frac{1-\nu}{1+\nu}\right]^\frac{1-\nu}{1+\nu}\end{gather*} и $\delta > 0$~--- фиксированный свободный параметр (искусственно вводимая неточность в оптимизационную модель).

Таким образом, возможно использовать градиентные методы при условии выбора ($\delta, L$)-модели вида \eqref{eq_stand_model}. Мы рассмотрим только оценку для ускоренного метода на базе  теоремы~\ref{mainTheoremDL}, поскольку именно ускоренный метод гарантирует оптимальность скорости сходимости при $\nu > 0$. 

\begin{remark}
Поскольку адаптивный алгоритм \ref{Acceler_Grad_Model} реализует возможность настройки на константу $L = L(\delta)$
в ходе работы метода, то для его практической реализации нет необходимости вычислять значения $L(\delta)$. При этом стоит отметить некоторый недостаток универсального метода~--- наличие искусственной погрешности $\delta>0$, от которой можно избавиться лишь в частном случае $\nu =1$. Без введений такой погрешности возникает проблема обобснования выполнимости критерия выхода из итерации \eqref{exitLDL}.
\end{remark}

Будем предполагать, что вспомогательные подзадачи мы можем решать точно, т.е. для всех \ag{$k \geqslant 0$} верно $\widetilde{\delta}=0$. Предположим, что ускоренный градиентный метод отработал ровно $N$ итераций и выберем
\begin{gather}\delta =\varepsilon \min_{k = \overline{0, N-1}}\frac{\alpha_{k+1}}{4A_{k+1}} ,\label{delta_k_univ}\end{gather} где $\varepsilon$~--- необходимая точность решения по функции. Чтобы выполнилось \begin{gather*}\frac{8LR^2}{N^2} \leq \frac{\varepsilon}{2},\end{gather*} достаточно сделать не больше $8\sqrt{L}R/\sqrt{\varepsilon}$ шагов внешнего цикла, где $L =  L(\delta)$. Выбор искусственной неточности (\ref{delta_k_univ}) гарантирует, что количество обращений ко вспомогательным подзадачам на итерациях методов будет конечным. Параметр $L$ можно оценить следующим образом:
\begin{gather*}
L =  L(\delta)=L_\nu\left[\frac{L_\nu}{2\delta}\frac{1-\nu}{1+\nu}\right]^\frac{1-\nu}{1+\nu} \leqslant L_\nu\left[\frac{4L_\nu A_N}{2\varepsilon\alpha_N}\frac{1-\nu}{1+\nu}\right]^\frac{1-\nu}{1+\nu}.
\end{gather*}
Так как
\begin{equation*}
\frac{A_N^2}{\alpha_N^2}=A_N L_N \leq N \alpha_N L_N=N \frac{A_N}{\alpha_N},
\end{equation*}
то $A_N / \alpha_N \leq N$ и поэтому
\begin{gather*}
L \leq L_\nu\left[\frac{2L_\nu N}{\varepsilon}\frac{1-\nu}{1+\nu}\right]^\frac{1-\nu}{1+\nu}.
\end{gather*}
Поскольку при достаточном малом $\varepsilon>0$ достаточное количество итераций $N$ можно оценить сверху $N \leq 8\sqrt{L}R / \sqrt{\varepsilon}$, то имеем
\begin{gather*}
L \leq L_\nu\left[\frac{16L_\nu \sqrt{L}R}{\varepsilon^\frac{3}{2}}\frac{1-\nu}{1+\nu}\right]^\frac{1-\nu}{1+\nu}.
\end{gather*}
Если перенести все выражения с $L$ в левую часть неравенства, то получим
\begin{gather*}
L^\frac{1+3\nu}{2+2\nu} \leq L_\nu\left[\frac{16L_\nu R}{\varepsilon^\frac{3}{2}}\frac{1-\nu}{1+\nu}\right]^\frac{1-\nu}{1+\nu}.
\end{gather*}
Отсюда имеем, что величина параметра $L$ ограничена.

Оценим скорость сходимости ускоренного градиентного метода (п.~\ref{fastGradMethod} и теорема~\ref{mainTheoremDL}) при указанных предположениях:
\begin{gather*}
F(x_N) - F(x^*) \leq \frac{8LR^2}{(N+1)^2}+ \frac{2\sum_{k=0}^{N-1}\delta A_{k+1}}{A_{N}} \leq \frac{8LR^2}{(N+1)^2}\, +  \\
+ \, \frac{2\varepsilon\sum_{k=0}^{N-1}\frac{\alpha_{k+1}}{2A_{k+1}}A_{k+1}}{A_{N}} = \frac{8LR^2}{(N+1)^2} + \frac{\varepsilon}{2} \leq \varepsilon.
\end{gather*}
Оценим достаточное для этого количество итераций $N$:
\begin{gather*}
N \leq \frac{8\sqrt{L}R}{\sqrt{\varepsilon}}=\frac{8\sqrt{L_\nu\left[\frac{L_\nu}{2\delta}\frac{1-\nu}{1+\nu}\right]^\frac{1-\nu}{1+\nu}}R}{\sqrt{\varepsilon}}= \frac{8\sqrt{L_\nu\left[\frac{4L_\nu A_N}{2\varepsilon\alpha_N}\frac{1-\nu}{1+\nu}\right]^\frac{1-\nu}{1+\nu}}R}{\sqrt{\varepsilon}}\\ \leq \frac{8\sqrt{L_\nu\left[\frac{2L_\nu N}{\varepsilon}\frac{1-\nu}{1+\nu}\right]^\frac{1-\nu}{1+\nu}}R}{\sqrt{\varepsilon}},\\
N^2 \leq \frac{64L_\nu\left[\frac{2L_\nu N}{\varepsilon}\frac{1-\nu}{1+\nu}\right]^\frac{1-\nu}{1+\nu}R^2}{\varepsilon}=\frac{64L_\nu^{\frac{2}{1+\nu}}\left[2 N\frac{1-\nu}{1+\nu}\right]^\frac{1-\nu}{1+\nu}R^2}{\varepsilon^\frac{2}{1+\nu}}.
\end{gather*}
Отсюда
\begin{gather*}
N^{\frac{1+3\nu}{1+\nu}} \leq \frac{64L_\nu^{\frac{2}{1+\nu}}\left[ \frac{2-2\nu}{1+\nu}\right]^\frac{1-\nu}{1+\nu}R^2}{\varepsilon^\frac{2}{1+\nu}}.
\end{gather*}
Поэтому получаем оценку числа итераций
\begin{gather*}
N = \inf_{\nu\in[0,1]}\left[64^\frac{1+\nu}{1+3\nu}\left(\frac{2-2\nu}{1+\nu}\right)^\frac{1-\nu}{1+3\nu}\left(\frac{L_\nu R^{1+\nu}}{\varepsilon}\right)^\frac{2}{1+3\nu}\right],
\end{gather*}
\ag{которая оптимальна с точностью до умножения на числовую константу \cite{NemirYd79}}.

\begin{remark}
Отметим, что на практике в ходе работы метода нет необходимости решать задачу минимизации для нахождения точного числа шагов $N$, поскольку \ag{алгоритм \ref{Acceler_Grad_Model}} адаптивен и предполагает <<настройку>> на оптимальный уровень и константу гладкости при выполнении самих итераций.
\end{remark}

\subsubsection{Метод условного градиента с точки зрения концепции модели функции}\label{model_uslovn}

Отметим, что при реализации на практике часто возникающие на итерациях градиентных методов вспомогательные подзадачи (\ref{equmir2DL}) не могут быть решены за разумное время. В этой связи в работе \cite{Jaggi_Rev_FW} показано, что весьма эффективным для некоторых классов задач может быть метод условного градиента Франк--Вульфа (см. раздел \ref{ch2_sect_fw}). Оказывается, что описанный выше результат (теорема \ref{mainTheoremDL}) об оценке скорости сходимости варианта метода подобных треугольников в модельной общности (алгоритм \ref{Acceler_Grad_Model}) с учётом неточности решения вспомогательных задач может привести к интересному обобщению метода условного градиента. В рамках этого обобщения вместо используемой алгоритмом \ref{Acceler_Grad_Model} оптимизационной модели вида
\begin{equation}\label{abstractmodelGrad}
\varphi_{k+1}(x)= \alpha_{k+1}\psi_{\delta}(x, y_{k+1}) + V(x, u_k)
\end{equation}
в (\ref{equmir2DL}) используется более простой её вариант
\begin{equation}\label{abstractmodelFW}
\widetilde{\varphi}_{k+1}(x)=\alpha_{k+1}\psi_{\delta}(x, y_{k+1}).
\end{equation}

Такой подход указывает на интересную связь метода Франк--Вульфа с ускоренным градиентным методом (пункт \ref{fastGradMethod}). В случае $\psi_{\delta}(x, y) = \langle \nabla F(y), x - y \rangle $ использование при реализации градиентных методов оптимизационных моделей вида \eqref{abstractmodelFW} вместо \eqref{abstractmodelGrad} как раз приводит к методу Франк--Вульфа. Чтобы получить интерпретацию результатов из пункта \ref{fastGradMethod} с точки зрения этой связи, рассмотрим замену \eqref{abstractmodelGrad} на \eqref{abstractmodelFW} с точки зрения параметра погрешности $\widetilde{\delta}$. 
Далее будем предполагать, что $F$~--- выпуклая функция с $L$-липшицевым градиентом по норме $\|\cdot\|$ по аналогии с предположениями теоремы \ref{thmuslFW}
$$V(x,y) \leq R_Q^2 \text{  для произвольных   } x,y \in Q.
$$
Пусть 
$$u_{k+1}=\argmin^{\widetilde{\delta}}_{x \in Q}\varphi_{k+1}(x) = \argmin_{x \in Q} \widetilde{\varphi}_{k+1}(x).$$
Тогда существуют субградиенты $\hat{g} \in \partial\varphi_{k+1}(u_{k+1})$ и $\tilde{g} \in \partial\widetilde{\varphi}_{k+1}(u_{k+1})$ такие, что
\begin{gather*}
\langle \hat{g}, x - u_{k+1} \rangle =\\= \langle \tilde{g}, x - u_{k+1} \rangle + \langle \nabla_{u_{k+1}} V(u_{k+1}, u_k), x - u_{k+1}\rangle \geq \langle \nabla_{u_{k+1}} V(u_{k+1}, u_k), x - u_{k+1}\rangle=\\ =
-V(u_{k+1},u_{k}) - V(x,u_{k+1}) + V(x,u_{k}) \geq -2R^2_Q.
\end{gather*}

Поэтому в алгоритме \ref{Acceler_Grad_Model} будем предполагать, что $\widetilde{\delta}=2R^2_Q$ и $\delta = 0$. Воспользуемся результатом пункта \ref{fastGradMethod} для ускоренного градиентного метода для минимизации функций $F$ с $L$-липшицевым градиентом относительно нормы $\|\cdot\|$ (см. теорему~\ref{mainTheoremDL}):
\begin{equation*}
F(x_N) - F(x_*) \leq \frac{8LR^2}{(N+1)^2} + \frac{16LR^2_Q}{N+1}.
\end{equation*}
Известно \cite{Lan2019}, что данная оценка с точностью до числового множителя не может быть улучшена для методов типа Франк--Вульфа.

Описанная в этом пункте конструкция была позаимствована из книги~\cite{Ben-Tal}.

\subsection{Универсальный метод условного градиента}
\label{universal_uslovn}

Снова будем отталкиваться от ускоренного алгоритма \ref{Acceler_Grad_Model} (см. п.~\ref{fastGradMethod}). Рассмотрим случай, когда $F$ --- выпуклая функция с $\nu$-гёльдеровым градиентом ($\nu \in (0; 1]$) \eqref{eq_hold_condit}, причём ставится цель использовать подход с помощью метода типа условного градиента. Положим $\widetilde{\delta}=2R^2_Q$, зафиксируем количество шагов $N$, а также искусственно введём неточность \begin{gather*}\delta = \min_{i = \overline{1, N}} \varepsilon\frac{\alpha_i}{4A_{i}}.\end{gather*} 
Тогда при соответствующем $L>0$ для алгоритма \ref{Acceler_Grad_Model} по теореме~\ref{mainTheoremDL} верно неравенство:
\begin{equation*}
F(x_N) - F(x_*) \leq \frac{8LR^2}{(N+1)^2} + \frac{16LR^2_Q}{N+1} + \frac{\varepsilon}{2} \leq \frac{24LR^2_Q}{N+1} + \frac{\varepsilon}{2}.
\end{equation*}
Поэтому при достаточно малом $\varepsilon > 0$ выполнение не более $96LR_Q^2/\varepsilon$ внешних шагов ускоренного метода гарантирует справедливость неравенства $24LR^2_Q/N \leq \varepsilon / 2$, т.е.
\begin{gather*}
N \leq \frac{96LR_Q^2}{\varepsilon} \leq \frac{96L_\nu\left[\frac{2L_\nu N}{\varepsilon}\frac{1-\nu}{1+\nu}\right]^\frac{1-\nu}{1+\nu}R_Q^2}{\varepsilon}.
\end{gather*}
После переноса $N$ в левую часть последнего неравенства получим, что
\begin{gather*}
N^\frac{2\nu}{1+\nu} \leq \frac{96L_\nu^\frac{2}{1+\nu}\left[\frac{2-2\nu}{1+\nu}\right]^\frac{1-\nu}{1+\nu}R_Q^2}{\varepsilon^\frac{2}{1+\nu}}.
\end{gather*}
Отсюда получаем важное замечание, что $\nu$ должно быть больше нуля, что полностью согласуется с примером из \cite{Nesterov18}, когда метод условного градиента может не сходиться в случае $\nu=0$. Далее, получим оценку числа итераций
\begin{gather*}
N = \inf_{\nu\in(0,1]}\left[96^\frac{1+\nu}{2\nu} \left[\frac{2-2\nu}{1+\nu}\right]^\frac{1-\nu}{2\nu}\left(\frac{L_\nu R_Q^{1+\nu}}{\varepsilon}\right)^\frac{1}{\nu}\right].
\end{gather*}

\begin{remark}
Как и для обычного универсального метода (см. п. \ref{Univers_Meth}), заметим, что на практике в ходе работы метода нет необходимости решать задачу минимизации для нахождения $N$, поскольку \ag{алгоритм \ref{Acceler_Grad_Model}} адаптивен и предполагает <<настройку>> на оптимальный уровень и константу гладкости при выполнении самих итераций.
\end{remark}


\subsection{Регуляризация и рестарты
}\label{restarts}

В данном разделе мы вернемся к изложению приемов \textit{регуляризации} и \textit{рестартов}, упомянутых в замечании~\ref{ideas}. А именно, в этом разделе будет продемонстрировано как, имея 
эффективные алгоритмы решения выпук\-лых\,/\,сильно выпуклых задач вида
\begin{equation}\label{E1}
F(x)\rightarrow\min_{x\in Q},
\end{equation}
можно предложить на их базе алгоритмы решения сильно выпуклых / выпуклых задач аналогичного вида. 

Начнем с изложения перехода [сильно выпуклые задачи] $\rightarrow$ [выпуклые задачи].

Введём семейство $\mu$-сильно выпуклых в норме $\|\cdot\|$ задач ($\mu>0$)
\begin{equation}\label{E22}
\min_{x\in Q} F^{\mu}(x):=F(x)+\mu V(x,y_{0}),
\end{equation}
где дивергенция Брэгмана $V(x,y)$ (см. определение~\ref{div-breg}) -- 1-сильно выпуклая функция по $x$ относительно $\|\cdot\|$-нормы.

\begin{teo}\label{th 4} Пусть
\begin{equation}\label{E23}
\mu\leq\frac{\varepsilon}{2V(x^*,y_{0})}=\frac{\varepsilon}{2R^{2}}
\end{equation}
и удалось найти $\varepsilon/2$-решение задачи~\eqref{E22}, т.е. такое $x^{N}\in Q$, что
$$F^{\mu}(x_{N})-\min_{x\in Q}F^{\mu}(x)
\leq\varepsilon/2.$$
Тогда
$$F(x_{N})-\min_{x\in Q}F(x)=F(x_{N})-F(x^{*})=F(x_{N})-F^{*}\leq\varepsilon.$$
\end{teo}
\begin{proof}
Действительно,
$$F(x_{N})-F^{*}\leq F^{\mu}(x_{N})-F^{*}\leq F^{\mu}(x_{N})-
\min_{x\in Q}F^{\mu}(x)
+\varepsilon/2\leq\varepsilon.$$
Здесь использовалось определение $F^{\mu}(x)$ и формула~\eqref{E23}:
$$
\min_{x\in Q}F^{\mu}(x)=\min_{x\in Q} \left\{F(x)+\gamma V(x,y_{0}) \right\}\leq F(x^{*})+\gamma V(x^{*},y_{0})\leq F^{*}+\varepsilon/2.
$$
\end{proof}

Таким образом, имея метод, решающий $\mu$-сильно выпуклую задачу в $\|~\cdot~\|$-норме, можно использовать его для решения не сильно выпуклой задачи $\eqref{E1}$ путём перехода к регуляризованной задаче $\eqref{E22}$ c $\mu = \varepsilon/(2R^2)$.

Приведём на примере алгоритма~\ref{Acceler_Grad_Model} в некотором смысле обратную конструкцию (рестартов): [выпуклые задачи] $\rightarrow$ [сильно выпуклые задачи].

\begin{teo}\label{th 5}
Пусть функция $F(x)$~--- $\mu$-сильно выпуклая в норме $\|\cdot\|$. Пусть точка $x_{\overline{N}}(y_{0})$ выдаётся ускоренным алгоритмом~\ref{Acceler_Grad_Model} при $\delta = \tilde{\delta} = 0$ (см. раздел \ref{fastGradMethod}; важно отметить, что алгоритм~\ref{Acceler_Grad_Model} был взят для определённости: описываемая далее конструкция позволяет рестартовать любой метод, в оценку скорости сходимости которого входит начальное расстояние / [дивергенция Брэгмана] от точки старта до решения или невязка по функции или норма градиента в точке старта при условии, что в решении норма градиента нулевая \cite{Aspremont2021}; условие $\delta = \tilde{\delta} = 0$ сделано лишь для упрощения формулировок), стартующим из точки $y_{0}$, после
\begin{equation}\label{E24}
\overline{N}=\sqrt{\frac{8L\omega_{n}}{\mu}}
\end{equation}
итераций, где
$$\omega_{n}=\sup_{x\in Q}\frac{2V(x,y_{0})}{\|x-y_{0}\|^{2}}.$$
Положим
$$\left[x_{\overline{N}}(y_{0})\right]^{1}=x_{\overline{N}}(y_{0})$$
и определим по индукции
$$
\left[x_{\overline{N}}(y_{0})\right]^{k+1}:=x_{\overline{N}}\left(\left[x_{\overline{N}}(y_{0})\right]^{k}\right),\quad k=1,\, 2,\, \ldots
$$
При этом на $(k+1)$-м перезапуске (рестарте) алгоритма \ref{Acceler_Grad_Model} также корректируется прокс-функция (считаем, что определённая таким образом функция корректно определена на $Q$ с сохранением свойства сильной выпуклости)
$$d^{k+1}(x):=d\left(x-\left[x_{\overline{N}}(y_{0})\right]^{k}+y_{0} \right)\geq0,$$
чтобы
$$d^{k+1}\left(\left[x_{\overline{N}}(y_{0})\right]^{k}\right)=0.$$
Тогда справедливо неравенство
\begin{equation}\label{E25}
F\left(\left[x_{\overline{N}}(y_{0})\right]^{k}\right)-F^{*}\leq\frac{\mu\|y_{0}-x^{*}\|^{2}}{2^{k+1}}.
\end{equation}
\end{teo}
\begin{proof}
Согласно теореме \ref{mainTheoremDL} при $\delta = \tilde{\delta} = 0$ после $\overline{N}$ итераций алгоритм \ref{Acceler_Grad_Model} выдаёт такой $x_{\overline{N}}$, что
$$\frac{\mu}{2}\|x_{\overline{N}}-x^{*}\|^{2}\leq F(x_{\overline{N}})-F^{*}\leq\frac{8LV(x^{*},y_{0})}{\overline{N}^{2}}.$$
Отсюда имеем
$$\|x_{\overline{N}}-x^{*}\|^{2}\leq\frac{16LV(x^{*},y_{0})}{\mu\overline{N}^{2}}\leq \|y_{0}-x^{*}\|^{2}\frac{8L\omega_{n}}{\mu\overline{N}^{2}}.$$
Поскольку
$$\overline{N}=\sqrt{\frac{16L}{\mu}\omega_{n}},$$
то
$$\|x_{\overline{N}}-x^{*}\|^{2}\leq\frac{1}{2}\|y_{0}-x^{*}\|^{2}.$$
Повторяя эти рассуждения по индукции, получаем
$$F\left(\left[x_{\overline{N}}(y_{0})\right]^{k}\right)-F^{*}\leq\left(\frac{1}{2}\right)^{k}\frac{8L\omega_{n}}{\overline{N}^{2}}\|y_{0}-x^{*}\|^{2}=\frac{\mu\|y_{0}-x^{*}\|^{2}}{2^{k+1}}.$$
\end{proof}

Предыдущая теорема гарантирует достижение описанной процедурой рестартов ускоренного градиентного метода $\varepsilon$-точного по функции решения задачи минимизации \eqref{E1} за
\begin{equation}\label{restart-compl}
O\left(\sqrt{\frac{L\omega_n}{\mu}} \left\lceil \log\left(\frac{\mu \|y_{0}-x^{*}\|^{2}}{\varepsilon}\right) \right\rceil\right).
\end{equation}
обращений к оракулу (подпрорамме) для вычисления градиента целевой функции.  

Из написанного выше представляется, что в процедуре рестартов (теорема~\ref{th 5}) можно использовать вместо предписанного числа итераций $\overline{N}$ на каждом рестарте какой-нибудь
критерий остановки. В частности, дожидаться, когда норма (или квадрат нормы) градиента (а в общем случае, когда минимум достигается не в точке экстремума,~--- норма градиентного отображения) уменьшится вдвое. С таким критерием остановки нет необходимости делать предписанное число итераций на каждом рестарте. Однако, пока неизвестен способ рассуждений, который позволял бы показать, что такая процедура сохраняет при переносе оптимальность оценок (метод, работающий оптимально не в сильно выпуклом случае, порождает оптимальный метод и в сильно выпуклом случае). Впрочем, имеются различные эффективные на практике способы более раннего выхода с каждого рестарта (см., например, \cite{DonoghueCandes15}), позволяющие (в случае задач безусловной оптимизации с евклидовой прокс-структурой) ускорить описанную выше конструкцию на порядок.

Как приложение, рассмотрим конкретный пример задачи сильно выпуклой композитной оптимизации
\begin{equation}
\label{eq15}
\mathop {\min
}\limits_{\sum\limits_{i=1}^n {x_i } =1,\;x\ge 0} F\left( x \right):=\frac{1}{2}\left\| {Ax-b} \right\|_2^2 +\mu
\sum\limits_{i=1}^n {x_i \log x_i }  
\end{equation}
при условии достаточно большого $\mu \gg \varepsilon \mathord{\left/ {\vphantom
{\varepsilon {\left( {2\log n} \right)}}} \right. \kern-\nulldelimiterspace}
{\left( {2\log n} \right)}$ (сильную выпуклость
композита в 1-норме необходимо учитывать).

Для указанной постановки можно выбрать норму в прямом пространстве $\left\|~\cdot~ \right\|=\left\|~\cdot~\right\|_1 $.
Положим в~\eqref{composite}
\[
f\left( x \right)=\frac{1}{2}\left\| {Ax-b} \right\|_2^2 ,
\quad
h\left( x \right)=\mu \sum\limits_{i=1}^n {x_i \log x_i } ,
\]
\[
Q=S_n \left( 1 \right)=\left\{ {x\ge 0:\;\sum\limits_{i=1}^n {x_i } =1}
\right\},
\]
где $L=\mathop {\max }\limits_{j=1,...,n} \left\| {A^{\left\langle j
\right\rangle }} \right\|_2^2$ -- константа Липшица градиента функции $f$ в 1-норме. Здесь  $A^{\left\langle j \right\rangle }$ обозначает $j$-й столбец матрицы $A$.


Заметим, что для случая достаточно малого $\mu$ ($0<\mu \ll\varepsilon \mathord{\left/ {\vphantom
{\varepsilon {\left( {2\log n} \right)}}} \right. \kern-\nulldelimiterspace}
{\left( {2\log n} \right)}$)
можно выбрать
\[
d(x) = \log n+\sum\limits_{i=1}^n {x_i \log x_i }, 
\]
и тогда
\[
V\left( {x,y} \right)=\sum\limits_{i=1}^n {x_i \log \left( {{x_i }
\mathord{\left/ {\vphantom {{x_i } {y_i }}} \right.
\kern-\nulldelimiterspace} {y_i }} \right)} ,
\]
\[
R^2 = V(x^*,x_0)\le \log n,\text{ где } x_0 = (1/n,...,1/n).
\]

Однако в ситуации достаточно большого $\mu$ для
выбранной функции $V\left( {x,y} \right)$ (расстояние Кульбака--Лейблера)  уже получается $\omega _n =\infty $. При этом существует
другой способ выбора прокс-функции ($d(x)$ будет 1-сильно выпуклой относительно 1-нормы на неотрицательном ортанте):
\begin{equation}
\label{eq16}
d\left( x \right)=\frac{1}{2\left( {a-1} \right)}\left\| x \right\|_a^2 ,
\quad
a=\frac{2\log n}{2\log n-1}.
\end{equation}
В таком случае также имеем $R^2={\rm O}\left( {\log n} \right)$. Однако
теперь и $\omega _n ={\rm O}\left( {\log n} \right)$.
Вот в такой ситуации возможно применить описанную выше схему рестартов ускоренного градиентного метода и из оценки сложности \eqref{restart-compl} получить следующий результат.

\begin{teo}\label{tomogravity}
Для задачи выпуклой композитной оптимизации~\eqref{eq15},
где 
$\mu \gg \varepsilon / (2 \log n)$
(т.е. сильную выпуклость композитного члена в $1$-норме необходимо учитывать) алгоритм~\ref{Acceler_Grad_Model} с рестартами
и с прокс-функцией 
$$
d(x) = \frac{1}{2(a-1)} ||x||_a^2,
$$
где $a = \frac{2\log n}{2\log n - 1}$
 приводит к необходимости на каждой итерации наряду с расчётом градиента гладкой части функционала (${\rm O}\left( {nnz\left( A \right)} \right)$ операций, где $nnz\left( A \right)$~-- число  ненулевых  элементов в матрице $A$) 
 решать вспомогательную подзадачу композитной оптимизации вида \eqref{ch1_eq_cx_log} 
 с помощью перехода к двойственной задаче и её решения с помощью прямо-двойственной версии метода эллипсоидов (${\rm
O}\left( {n\log \left( {C \mathord{\left/ {\vphantom {C \varepsilon }}
\right. \kern-\nulldelimiterspace} \varepsilon } \right)\log \left( {{nC}
\mathord{\left/ {\vphantom {{nC} \varepsilon }} \right.
\kern-\nulldelimiterspace} \varepsilon } \right)} \right)$ арифметических операций, см. раздел~\ref{comp_prox}). При этом {\rm(}$\tilde{x}_N$ --- то, что выдаёт алгоритм на последнем рестарте{\rm):}
\[
F\left(\tilde{x}_N \right)-F^* \le \varepsilon ,
\]
если общее число итераций (обращений к оракулу за градиентом)
\begin{equation}
\label{ch3_eq_5515_estim}
N={\rm O}\left( {\sqrt {\frac{L}{\mu }\log n} \left\lceil {\log \left(
{\frac{\mu }{\varepsilon }} \right)} \right\rceil } \right)={\rm O}\left(
{\sqrt {\frac{\mathop {\max }\limits_{j=1,...,n} \left\| {A^{\left\langle j
\right\rangle }} \right\|_2^2 \log n}{\mu }} \left\lceil {\log \left(
{\frac{\mu }{\varepsilon }} \right)} \right\rceil } \right).
\end{equation}
\end{teo}

Заметим, что в <<пороговой>> ситуации, отвечающей регуляризации (см. выше теорему~\ref{th 4}),
$\mu \simeq \varepsilon \mathord{\left/ {\vphantom {\varepsilon {\left(
{2\log n} \right)}}} \right. \kern-\nulldelimiterspace} {\left( {2\log n}
\right)}$. В этом случае формула (\ref{ch3_eq_5515_estim}) примет вид
\[
N={\rm O}\left( {\sqrt {\frac{\mathop {\max }\limits_{k=1, \, \dots, \,n} \left\|
{A^{\left\langle k \right\rangle }} \right\|_2^2 \log ^2n}{\varepsilon }} }
\right),
\]
что с точностью до $\sim\sqrt {\log n}$ соответствует оценке в случае
$0 \, < \, \mu \ll \varepsilon / (2 \log n)$
(когда сильную выпуклость композитного члена в $1$-норме можно не учитывать). Отличие этих двух случаев ещё и в том, что в случае, когда $0 \, < \, \mu \ll \varepsilon / (2 \log n)$,  существует простой способ
добиться стоимости итерации ${\rm O}\left( {nnz\left( A \right)} \right)$ (прямое использование алгоритма~\ref{Acceler_Grad_Model} c $V\left( {x,y} \right)=\sum\limits_{i=1}^n x_i \log \left(x_i/y_i \right)$), а
в случае $\mu \gg \varepsilon / (2 \log n)$ нам не известно более эффективного способа, чем способ
(описанный выше) со стоимостью итерации
\[
{\rm O}\left( {nnz\left( A \right)+n\log \left( {C \mathord{\left/ {\vphantom
{C \varepsilon }} \right. \kern-\nulldelimiterspace} \varepsilon }
\right)\log \left( {{nC} \mathord{\left/ {\vphantom {{nC} \varepsilon }}
\right. \kern-\nulldelimiterspace} \varepsilon } \right)} \right).
\]
Поскольку в типичных приложениях первое слагаемое заметно доминирует второе,
то можно не учитывать (что часто и делают на практике) 
плату за невозможность выполнения <<проектирования>> по явным формулам и не
сильно задумываться, с какой точностью решать вспомогательную задачу, делая
это с точностью длины мантиссы (описанный подход позволяет решать её с очень
хорошей точностью, и сложность решения вспомогательной задачи практически не
чувствительна к этой точности). По-видимому, это
утверждение имеет достаточно
широкий спектр практических приложений. 

В данном разделе на конкретном примере мы
продемонстрировали более подробно, чем это принято на практике, что для большого класса задач наличие явных
формул для шага итерации не есть сколько-нибудь сдерживающее обстоятельство
для использования метода. Используемая при этом техника и способ рассуждений
характерным образом (на наш взгляд) демонстрируют современный арсенал
средств (описанных в этой главе) решения задач выпуклой оптимизации в
пространствах больших размернотей.

Отметим также, что нам не известно никакого другого способа (кроме описанного в этом разделе приема рестартов), позволяющего решать задачу~\eqref{eq15} ускоренным методом, учитывающим её сильную выпуклость в 1-норме.


\section{Другие концепции неточного градиента}\label{inexact_grad}

В действительности, две наиболее популярные на практике концепции неточного градиента отличаются от описанной ранее концепции Деволдера--Глинёра--Нестерова (см. определение \ref{delta_L_oracle} 
и упражнение \ref{ex_oracle_fin}). А именно, для задачи
\begin{equation*}
    \min_{x \in Q} f(x) 
\end{equation*}
считаем, что доступен такой неточный градиент $\tilde{\nabla} f(x)$, что \cite{Polyak}: для всех $x \in Q$
\begin{equation}\label{inexact}
    \|\tilde{\nabla} f(x) - \nabla f(x)\|_2 \le \delta
\end{equation}
или
\begin{equation}\label{relative_inexact}
    \|\tilde{\nabla} f(x) - \nabla f(x)\|_2 \le \alpha\|\nabla f(x)\|_2, \quad \alpha\in[0,1).
\end{equation}
Далее в изложении мы будем в основном следовать работе \cite{Vasin}. Отметим, что многие вопросы, обсуждаемые ниже, возникли ещё в начале 80-х годов XX века. Особо отметим работы Б.Т. Поляка того периода, в которых можно найти частичные ответы на эти вопросы.

\subsection{Основные результаты}
Прежде всего заметим, что для концепции \eqref{inexact} 
\begin{equation*}
    \|\tilde{\nabla} f(x) - \nabla f(x)\|_2 \le \delta
\end{equation*}
даже для гладкой сильно выпуклой задачи и обычного градиентного метода нельзя гарантировать в общем случае (без дополнительных приёмов, типа регуляризации задачи или ранней остановки используемого метода) какие-то хорошие оценки, показывающие, что шумом можно пренебречь в случае его малости. 

Действительно, рассмотрим такой пример:
\begin{equation}\label{example}
 \min_{x\in\mathbb R^n} f(x):=\frac{1}{2}\sum_{i=1}^n \lambda_i (x^i)^2,   
\end{equation}
где $0\le\mu = \lambda_1 \le \lambda_2 \le ...\le \lambda_n = L$, $L\ge2\mu$. Решением задачи \eqref{example} будет $x^* = 0$. Предположим, что неточность имеется только в вычислении первой компоненты градиента. То есть, вместо $\partial f(x) / \partial x^1 = \mu x^1$ нам доступно только $\tilde{\partial} f(x) / \partial x^1 = \mu x^1 - \delta$. Для простейшей динамики градиентного спуска (в этом разделе мы обозначаем номер компоненты вектора верхним индексом, потому что далее будет использоваться одновременно номер итерации и номер компоненты \eqref{lowerbound}):
\ag{$$x_{k+1} = x_{k} - \frac{1}{L}\tilde{\nabla}f(x_{k}),$$}
можно получить, что при \ag{$x^1_0 \ge 0$ и достаточно больших} $k\in\mathbb N$ ($k \gg L/\mu$)
\begin{equation}\label{lowerbound}
 x_k^1 \ge \frac{\delta}{L}\frac{1 - (1 -\mu/L)^k}{1-(1 -\mu/L)} \ag{\simeq} \frac{\delta}{\ag{\mu}}.   
\end{equation}
Следовательно, \ag{при $x^1_0 \ge 0$ и достаточно больших $k\in\mathbb N$}
$$f(x_k) - f(x_*) \ag{\gtrsim} \frac{\delta^2}{2\mu}.$$

Таким образом, имеется проблема с нижней оценкой \eqref{lowerbound}, поскольку $\mu$ может быть малым ($\mu\lesssim \varepsilon$ -- вырожденный режим, где $\varepsilon$ -- желаемая точность решения исходной задачи по функции) в знаменателе правой части неравенства \eqref{lowerbound}. Для ускоренных методов можно ожидать ещё больших проблем.

Чтобы получить наиболее тонкие (из известных сейчас) оценок накопления неточности, немного специфицируем концепцию неточного градиента Деволдера--Глинёра--Нестерова: 
\begin{align}
\label{inexact_model_2}
&f(x)+\langle\tilde{\nabla} f(x),y-x\rangle +\frac{\mu}{2}\|y - x\|_2^2- \delta_1\|y-x\|_2 \le f(y) \le \nonumber \\ 
&\hspace{3cm}\le f(x)+\langle\tilde{\nabla} f(x),y-x\rangle + \frac{L}{2}\|y-x\|_2^2 +\delta_2.
\end{align}

С такой концепцией имеем для неускоренного градиентного метода 
\begin{align}\label{convergenceDNA} 
&f(x_k) - f(x^*) = \nonumber \\ 
&\hspace{0cm} =O\left(\min\left\{ \frac{LR^2}{k} + \tilde{R}\delta_1 + \delta_2,  LR^2\exp\left(-\frac{\mu}{L}k\right) + \tilde{R}\delta_1 + \delta_2\right\}\right), 
\end{align}
для ускоренного (можно брать разные варианты, описанные ранее в пособии в п.~\ref{gradientTaylor}, \ref{AM}, \ref{model}): 
\begin{align}\label{convergenceDA} 
&f(x_k) - f(x^*) = \nonumber \\ 
&\hspace{0cm}
=O\left(\min\left\{ \frac{LR^2}{k^2} + \tilde{R}\delta_1 + k\delta_2,  LR^2\exp\left(-\sqrt{\frac{\mu}{L}}\frac{k}{2}\right) +  \tilde{R}\delta_1 + \sqrt{\frac{L}{\mu}}\delta_2\right\}\right), 
\end{align}
где $\tilde{R}$ -- максимальное расстояние (в 2-норме) между генерируемыми методом точками и решением, а $R$ -- расстояние (в 2-норме) от точки старта до решения (если решений много, то до ближайшего к точке старта). По определению $R\le\tilde{R}$.
 
Таким образом из \eqref{convergenceDNA}, \eqref{convergenceDA} можно получить, что если $\tilde{R}$ ограничен\footnote{Во многих случаях это предположение справедливо. Например, когда $Q$ ограничено или когда  $\mu\gg\varepsilon$.}, то выбирая\footnote{Здесь также используются неравенства
$$
\langle\tilde{\nabla} f(x) - \nabla f(x),y-x\rangle \le \frac{1}{2L}\|\tilde{\nabla} f(x) - \nabla f(x)\|_2^2 + \frac{L}{2}\|y-x\|_2^2,
$$
$$
\frac{\varepsilon^2}{2LR^2}\max\left\{\sqrt{\frac{LR^2}{\varepsilon}},\sqrt{\frac{L}{\mu}}\right\}\le \varepsilon.
$$} 
 \begin{equation*}
    \delta_1 = \delta_{\eqref{inexact}}, \delta_2 = \frac{\delta_{\eqref{inexact}}^2}{2L}, 
 \end{equation*}
  где $\delta_{\eqref{inexact}}$ определено в  \eqref{inexact} как $\delta$, будем иметь следующий результат: 
  
 \textit{Для неускоренного и ускоренного методов возможно достичь уровня $f(x_k) - f(x_*) = \varepsilon$ с $\delta_{\eqref{inexact}}\simeq\varepsilon/\tilde{R}$ за такое же по порядку число итераций, как если бы $\delta_{\eqref{inexact}}=0$.}  
 
Однако в общем случае с предположением об ограниченности $\tilde{R}$ есть проблемы. Как мы уже видели из примера~\eqref{example}, в общем случае для вырожденной задачи (константа сильно выпуклости $\mu$ сильно меньше желаемой точности решения задачи по функции $\varepsilon$) имеет место оценка: $$\tilde{R}\simeq R + \frac{\delta_{\eqref{inexact}}}{\mu}.$$
Отмеченную проблему больших $\tilde{R}$ можно решить, используя регуляризацию задачи с параметром регуляризации $\mu\simeq\varepsilon/R^2$. В этом случае $\delta_{\eqref{inexact}}\simeq\varepsilon/R$ и мы имеем $\tilde{R}\simeq R$. Альтернативным решением является <<ранняя остановка>> метода \cite{Vasin}. Эта терминология используется также в машинном обучении, где <<ранняя остановка>> также используется как альтернатива к регуляризации для препятствования переобучению  \cite{good}.

Для концепции~\eqref{relative_inexact}: \begin{equation}
    \|\tilde{\nabla} f(x) - \nabla f(x)\|_2 \le \alpha\|\nabla f(x)\|_2, \quad \alpha\in[0,1)
\end{equation}
неускоренные методы будут сходиться к решению в $(1-\alpha)^{-1}$ раз медленнее, а про ускоренные на данный момент доказано только то, что они будут сходиться с такой же по порядку скоростью к решению, если $\alpha \lesssim \left(\mu/L\right)^{3/4}$ --- в сильно выпуклом случае и на $k$-й итерации метода $\alpha_k \lesssim k^{-3/2}$~--- в вырожденном случае ($\mu\ge 0$ мало), \ag{см. также \cite{Gannot}.} Заметим, что в численных экспериментах, в которых шум был случайный (не враждебный), эти условия можно существенно ослабить \cite{Vaswani}.

В следующем разделе мы приведём характерный пример того, как в задачах выпуклой оптимизации могут возникать неточные градиенты.

\subsection{Обратная задача для эллиптического уравнения\\ на~квадрате}\label{subsec:inverse}

Одним из наиболее естественных приложений численных методов оптимизации является теория некорректных обратных задач \cite{kab}.

Характерный пример обратной задачи будет рассмотрен ниже. Поскольку эти задачи часто связаны с решением линейных краевых задач для уравнений в частных производных, то вычисление градиента целевого функционала возможно только приближённо.

Для описания способа вычисления градиента целевого функционала в таких задачах сначала напомним (см.~\cite{evt}, а также п.~\ref{ch1_sect_dem-dan}), 
как считать градиент функции, 
$$
J(q):=J(q,u(q)),
$$ 
где $u(q)$ определяется единственным образом из  $$G(q,u)=0.$$ Предположим, что линейный оператор $G_{q}(q,u)$ (для конечномерных задач матрица, сформированная по принципу, что $(i,j)$ элемент равен $\partial G_i/\partial q_j$), имеет обратный $[G_{q}(q,u)]^{-1}$. Тогда 
 $$G_q(q,u) + G_u(q,u)\nabla u(q) = 0,$$
 следовательно,
 $$\nabla u(q) = - \left[G_u(q,u)\right]^{-1}G_q(q,u).$$
Таким образом,
 $$\nabla J(q):=J_q(q,u) + J_u(q,u)\nabla u(q) = J_q(q,u) - J_u(q,u)\left[G_u(q,u)\right]^{-1}G_q(q,u).$$
Аналогичный результат может быть получен с использованием формализма Лагранжа 
$$L(q,u;\psi) = J(q,u(q)) + \langle \psi, G(q,u) \rangle,$$ где $$L_u(q,u;\psi)=0, G_q(q,u)=0$$ и $$\nabla J(q) = L_q(q,u;\psi).$$
Действительно, можно связать эти два подхода соотношением 
$$\psi(q,u) = - \left[G_u(q,u)^T\right]^{-1}J_u(q,u)^T.$$

Теперь можно приступить к демонстрации описанной выше техники вычисления градиента (на базе принципа множителей Лагранжа) для обратной начально-краевой задачи для эллиптического уравнения на квадрате.
	
Пусть $u$ является решением следующей задачи (P):
\[
u_{xx} +u_{yy} =0,
\quad
x,y\in \left( {0,1} \right),
\]
\[
u\left( {1,y} \right)=q\left( y \right),
\quad
y\in \left( {0,1} \right),
\]
\[
u_x \left( {0,y} \right)=0,
\quad
y\in \left( {0,1} \right),
\]
\[
u\left( {x,0} \right)=u\left( {x,1} \right)=0,
\quad
x\in \left( {0,1} \right).
\]

Первые два уравнения
\[
-u_{xx} -u_{yy} =0,
\quad
x,y\in \left( {0,1} \right),
\]
\[
q\left( y \right)-u\left( {1,y} \right)=0,
\quad
y\in \left( {0,1} \right)
\]
обозначим через $G(q,u) = \bar{G}\cdot(q,u)= 0$, а последние два уравнения --- через $u\in Q$.

\textbf{Обратная задача:} \textit{Предположим, что мы хотим оценить $q(y)\in L_2(0,1)$, наблюдая $b(y) = u(0,y)\in L_2(0,1)$, где $u(x,y)\in L_2\left((0,1)\times(0,1)\right)$ -- единственное решение задачи (P).} 

Эта обратная задача, естественным образом, сводится к задаче квадратичной оптимизации:
\begin{equation}\label{inverse}
\min_{q} \mathfrak{J}(q) := \min_{u:\text{ }\bar{G}\cdot(q,u)= 0, u\in Q} J(q,u):= J(u) = \int_0^1 |u(0,y) - b(y)|^2 dy,
\end{equation}
которую  можно решать численно, если удастся посчитать $\nabla \mathfrak{J}(q)$.
Применим для вычисления $\nabla \mathfrak{J}(q)$ в задаче \eqref{inverse} описанный выше принцип множителей Лагранжа:
\begin{multline*}
    L\left(q,u;\psi:= \left(\psi(x,y),\lambda(y)\right)\right) = J(u) + \langle \psi, \bar{G}\cdot(q,u) \rangle = \\
    = \int_0^1 |u(0,y) - b(y)|^2\,dy - \int_0^1\int_0^1\left(u_{xx}+u_{yy}\right)\psi(x,y)\,dx\,dy \,+ \\
    + \int_0^1 \left(q(y) - u(1,y) \right)\lambda(y)\,dy.
\end{multline*}
Чтобы получить \textit{сопряжённую} задачу на  $\psi$, следует проварьировать $L\left(q,u;\psi\right)$ по $\delta u$ при условии, что $u \in Q$:
\begin{multline}\label{du}
 \delta_u L\left(q,u;\psi\right) =  2\int_0^1 \left(u(0,y) - b(y)\right)\delta u(0,y)\,dy -  \\
   - \int_0^1\int_0^1\left(\delta u_{xx}+\delta u_{yy}\right)\psi(x,y)\,dx\,dy  - \int_0^1 \delta u(1,y)\lambda(y)\,dy,
\end{multline}
где
\[
\delta u_x \left( {0,y} \right)=0,
\quad
y\in \left( {0,1} \right),
\]
\[
\delta  u\left( {x,0} \right)=\delta u\left( {x,1} \right)=0,
\quad
x\in \left( {0,1} \right).
\]
Интегрируя по частям, из \eqref{du} получим
\begin{multline*}
 \delta_u L\left(q,u;\psi\right) =  \int_0^1 \left(2\left(u(0,y) - b(y)\right) - \psi_x(0,y)\right)\delta u(0,y)\,dy -  \\
 - \int_0^1 \psi(1,y)\delta u_x(1,y)\,dy -  \int_0^1 \psi(x,1)\delta u_y(x,1)\,dx\, + 
  + \int_0^1 \psi(x,0)\delta u_y(x,0)\,dy +
 \\
    \int_0^1\int_0^1\left(\psi_{xx}+  \psi_{yy}\right)\delta u(x,y)\,dx\,dy  +
    \int_0^1 \left(\psi_x(1,y) - \lambda(y)\right)\delta u(1,y)\,dy.
\end{multline*}

Рассмотрим сопряжённую задачу (D):
\[
\psi _{xx} +\psi _{yy} =0,
\quad
x,y\in \left( {0,1} \right),
\]
\[
\psi _x \left( {0,y} \right)=2\left(u(0,y) - b(y)\right),
\quad
y\in \left( {0,1} \right),
\]
\[
\psi \left( {1,y} \right)=0,
\quad
y\in \left( {0,1} \right),
\]
\[
\psi \left( {x,0} \right)=\psi \left( {x,1} \right)=0,
\quad
x\in \left( {0,1} \right)
\]
и дополнительные условия между множителями Лагранжа: 
\begin{equation}\label{psi}
\lambda(y) = \psi_x(1,y),\quad y\in \left( {0,1} \right).
\end{equation}
Эти соотношения выводятся из того, что $\delta_u L\left(q,u;\psi\right) = 0$ и $\delta u(0,y),\delta u_x(1,y),\delta u(1,y) \in L_2(0,1)$; $\delta u_y(x,1),\delta u_y(x,0) \in L_2(0,1)$; $\delta u(x,y) \in L_2\left((0,1)\times(0,1)\right)$ произвольные.    

Поскольку (запись $u: (q,u)\in (P)$ и $\psi\in(D)$ следует понимать так, что $u$ является решением задачи (P), а $\psi$ -- решением задачи (D)):
\[
\mathfrak{J}(q) = \min_{u: (q,u)\in (P)} J(u) = \min_{u:\text{ }\bar{G}\cdot(q,u)= 0, u\in Q} J(u) =  \min_{u\in Q}\max_{\psi\in(D)} L(q,u;\psi),
\]
то по формуле Демьянова--Данскина--Рубинова (см.~п.~\ref{ch1_sect_dem-dan}):
\[
\nabla \mathfrak{J}(q) = \nabla_q \min_{u\in Q}\max_{\psi\in(D)} L(q,u;\psi) = L_q(q,u(q);\psi(q)),
\]
где $u(q)$ -- решение задачи (P), а $\psi(q)$ -- решение задачи (D), в которой 
\[\psi _x \left( {0,y} \right)=2\left(u(0,y) - b(y)\right),
\quad
y\in \left( {0,1} \right)
\]
и $u(0,y)$ зависит от $q(y)$ через задачу (P). Причём, в то же самое время $\left(u(q),\psi(q)\right)$ является решением седловой задачи $$\min_{u\in Q}\max_{\psi\in(D)} L(q,u;\psi).$$ 
Действительно, из $\delta_{\psi} L(q,u;\psi) = 0$ следует $\bar{G}\cdot(q,u)= 0$, что соответствует задаче (P), если добавить ограничение $u\in Q$. Из $\delta_u L_(q,u;\psi) = 0$ при $u\in Q$ получается задача (D), как было показано выше. 

Отметим также, что
\[L_q(q,u(q);\psi(q))(y) = \lambda(y), \quad
y\in \left( {0,1} \right).\]
Следовательно, согласно \eqref{psi} 
\[
\nabla \mathfrak{J}(q)(y) = \psi_x(1,y),\quad y\in \left( {0,1} \right).
\]
Таким образом, мы свели вычисление $\nabla \mathfrak{J}(q)(y)$ к решению двух корректных начально-краевых задач ((P) и (D)) для уравнений эллиптического типа на квадрате.

Полученный результат также можно проинтерпретировать немного по-другому.
Введём линейный оператор
$$
A:\quad q(y):=u(1,y)\mapsto u(0,y),
$$
где $u(x,y)$ -- решение задачи (P). В работе \cite{kab} было показано, что $A$ -- ограниченный оператор: 
$$
A: L_2(0,1) \rightarrow  L_2(0,1).
$$

Сопряжённый оператор (также ограниченный \cite{kab}) будет иметь вид
$$
A^\ast:\quad p(y):=\psi_x(0,y)\mapsto \psi_x(1,y), \qquad 
A^\ast : L_2(0,1)\rightarrow   L_2(0,1),
$$
где $\psi(x,y)$ -- решение сопряжённой задачи (D). 
Таким образом, считая, что
\[
\mathfrak{J}(q)(y) = \|Aq - b\|_2^2, 
\]
можно получить
\[
\nabla \mathfrak{J}(q)(y) = A^\ast\left( 2\left(Aq - b\right)\right), 
\]
что полностью соответствует описанной выше схеме:

\textbf{1.} \textit{Зная $q(y)$, решить задачу {\rm (}P{\rm )} и получить $u(0,y) = Aq(y)$. Затем определить $p(y) = 2\left(u(0,y) - b(y)\right)$.}

\textbf{2.} \textit{Зная $p(y)$ {\rm (}см.~п.~{\rm 1),} решить сопряжённую задачу {\rm (}D{\rm )} и вычислить $\nabla \mathfrak{J}(q)(y) = A^\ast p(y) = \psi_x(1,y)$.}

Таким образом, неточность при вычислении градиента $\nabla \mathfrak{J}(q)$ неизбежно будет возникать, поскольку мы можем решить краевые задачи (P) и (D) только численно.

Описанная выше техника может быть использована при решении большого числа обратных задач \cite{kab} и задач оптимального управления  \cite{Vas_met_opt2}. Отметим, что для задач оптимального управления на практике часто используют другую стратегию \cite{evt}. Сначала дискретизируют саму постановку задачи, а потом для дискретизированной задачи с помощью описанной выше техники считают точно градиент. В теории оптимального управления таким образом получают, например, дискретный принцип максимума. Ну, а сама процедура вычисления градиента оказывается аналогичной процедуре \textit{автоматического дифференцирования}  
(\textit{обратного распространения} по терминологии глубокого обучения \cite{good}), см.~п.~\ref{subsect_autodiff}.

\ag{Пионерами в приложении градиентных методов к задачам оптимизации в гильбертовых пространствах являются Л.\,В.~Канторович (40-е годы XX~века) и Б.\,Т.~Поляк (60-е годы XX века).} 



\chapter[Примеры задач и их практическое решение\dotfill]{Примеры задач\\ и их практическое решение}\label{chapt_examples}
\chaptermark{Примеры задач и их практическое решение}

В этой главе собраны примеры задач выпуклой оптимизации, которые позволяют продемонстрировать работу численных методов, описанных в предыдущей главе. Глава начинается с задачи расчёта барицентра Васерштейна. Оказывается, вместо того, чтобы решать исходную (прямую задачу), выгодно её специальным образом сгладить (энтропийное сглаживание) и построить двойственную. Двойственная задача будет обладать заметно лучшими свойствами (в частности, гладкостью). Её и предлагается решать ускоренным (быстрым) градиентным методом. Примеры на метод Франк--Вульфа хорошо демонстрируют, что оптимальный по числу обращений к оракулу (итераций) метод совсем не обязан оставаться оптимальным и по времени работы. Особенность метода Франк--Вульфа хорошо проявляется при решении задач выпуклой оптимизации на симплексе. В этом случае часто удаётся сильно сэкономить на стоимости каждой итерации, что в конечном итоге может приводить к лучшим оценкам времени работы метода Франк--Вульфа по сравнению с оптимальными (ускоренными) методами для гладких выпуклых задач. На вариациях задачи оптимального проектирования механических конструкций будет продемонстрирована работа метода субградиентного спуска для задач с ограничениями, а также теория полуопределённого программирования (semidefinite programming). Эта важная теория демонстрируется ещё на нескольких примерах из анализа сигналов. 

\section{Расстояние и барицентры Васерштейна}
\epigraph{...more or less everything in
this subject should be called after Kantorovich.}{\textit{C.~Villani. Optimal transport, old and new}}
Описываемая в этом разделе задача (Монжа) является ярким примером задачи с большой историей (берущей своё начало в XVIII веке). Однако за последние 10 лет задача Монжа нашла много новых интересных приложений. В частности,  полученные здесь в последние годы достижения разные специалисты отмечают в числе наиболее важных событий в области численных методов оптимизации последнего времени, см., например, \url{https://blogs.princeton.edu/imabandit/2019/12/30/a-decade-of-fun-and-learning/}.

\subsection{Примеры. Задача Монжа}
Для того, чтобы сравнить друг с другом объекты,
которые отличаются 
вероятностными распределениями
или гистограммами
и моделируют объекты реального
мира~--- изображения, видео, тексты и т.д.,
необходимо чётко определить,
как именно измерять
степень <<похожести>> или <<близости>>
между такими объектами.
Например, документы можно описать с помощью
встречающихся в них ключевых слов~\cite{EMD_techrep_yaleu}.
В этом случае математическим аналогом
документа является вероятностное распределение
в пространстве слов. Тогда
для сравнения документов нужно
искать отличия в их
вероятностных распределениях.
Таким образом, часто возникает необходимость
вычислить расстояние между точками
с неизвестными заранее местоположениями,
подчиняющимися только каким-то вероятностным
законам\footnote{Это могут быть
какие-то выборки или результаты наблюдений и т.п.}.

\begin{figure}[htb]
\begin{center}
\includegraphics[width=0.3\linewidth]{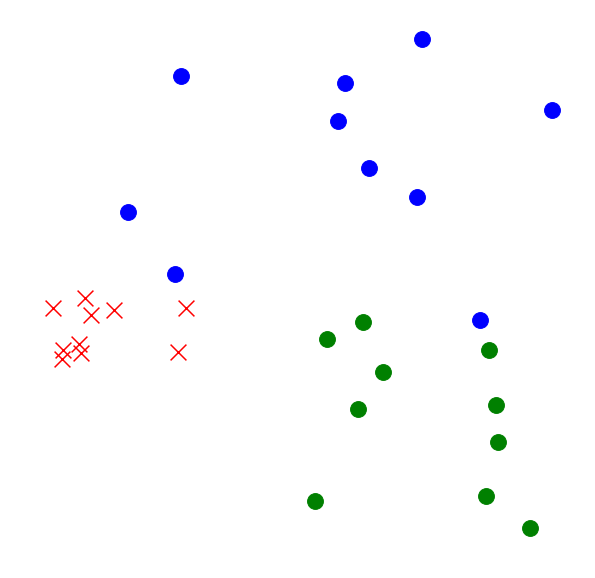}
\end{center}
\caption{Какие кластеры из точек расположены ближе друг к другу~-- голубые кружки к зелёным кружкам или голубые кружки к красным крестикам?} 
\label{ch1_dr_wasser_1}
\end{figure}

\begin{example}\label{ch1_exam_bar1}
На плоскости лежит множество точек, принадлежащих
трём различным кластерам. Соответственно, точки раскрашены в три цвета
(см. рис.~\ref{ch1_dr_wasser_1}). 
Как определить, какие из всех пар кластеров 
расположены ближе друг к другу?

Существуют несколько неуниверсальных способов это сделать~--- считать расстояние от ближайшей точки одного кластера до ближайшей точки другого (не подходит, если точки кластеров перемешаны),
или найти центры масс каждого кластера,
и считать расстояние между этими центрами
(не подходит для кластеров с совпадающими
центрами масс\protect\footnotemark).

Будем рассматривать множество точек одного цвета
(назовём его~$A$)
как выборку дискретной случайной величины, 
т.е.
для любой ${\xi \in A}$ вероятность ${p_{\xi} = 1/|A|}$, где $|A|$ обозначает количество элементов в~$A$.
Расстояние между $A$ и множеством точек другого цвета $B$
можно найти как решение задачи линейного программирования.
 Прежде чем пояснить детали,
рассмотрим следующий пример.
\end{example}
\footnotetext{Если мы хотим, чтобы
расстояние между такими кластерами
было отлично от нуля.}

\begin{figure}[htb]
\begin{minipage}[h]{0.49\linewidth}
\center{\includegraphics[width=0.2\linewidth]{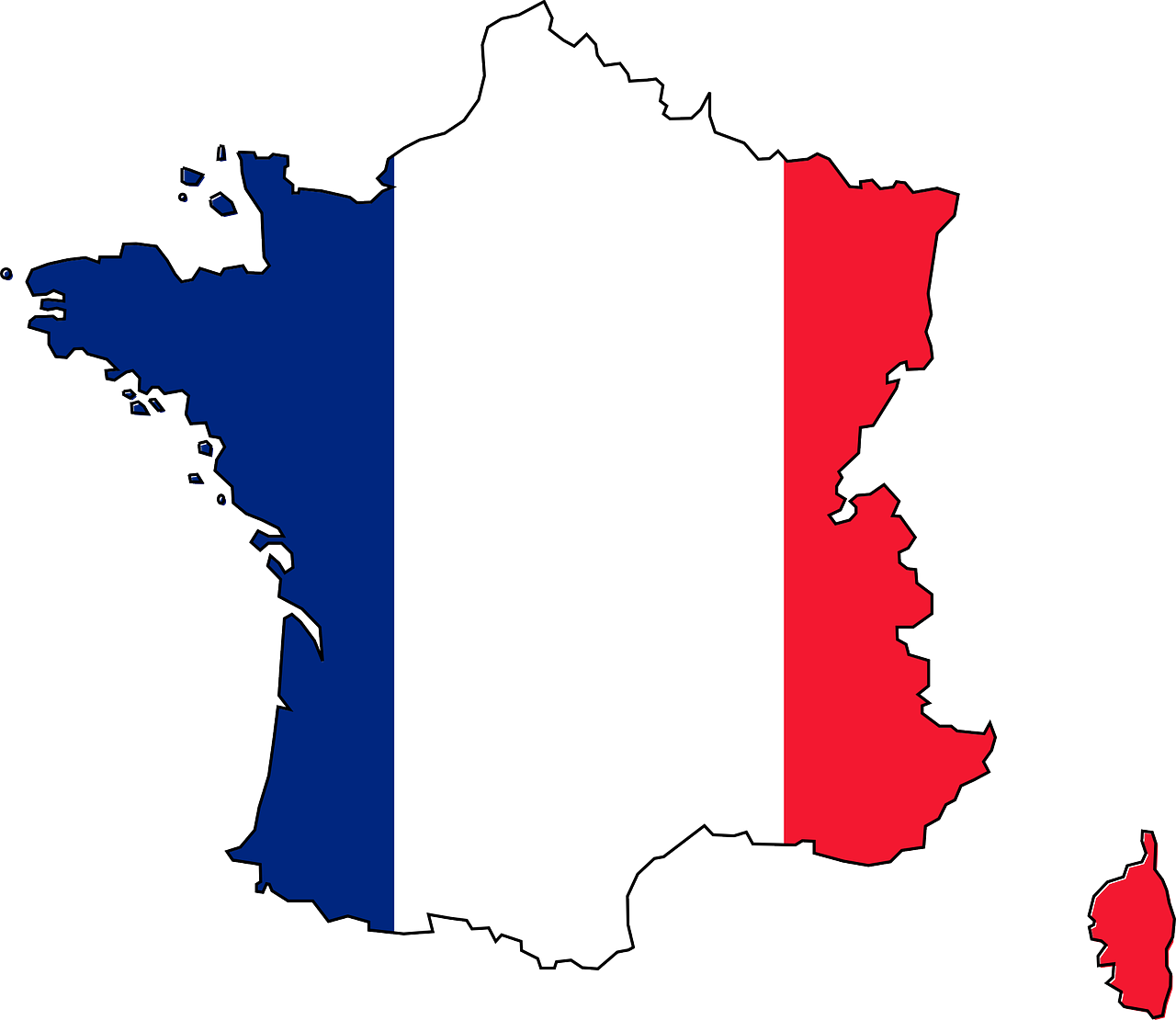} }
\end{minipage}
\hfill
\begin{minipage}[h]{0.49\linewidth}
\center{\includegraphics[width=0.99\linewidth]{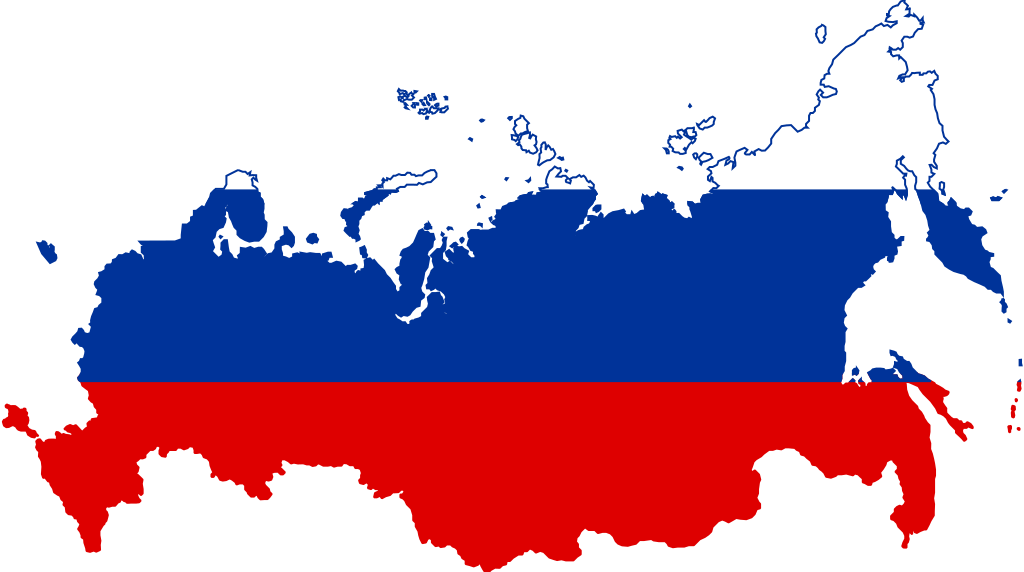}}
\end{minipage}
\caption[Как посчитать расстояние от России до Франции?]%
{Как посчитать расстояние от России до Франции?
\protect\footnotemark
}
\label{ch1_dr_wasser_2}
\end{figure}
\footnotetext{Отношение площадей стран на картинке не совсем соответствует реальным данным, на самом деле площадь Франции в 26.7 раза меньше площади Российской Федерации.}

\begin{example}\label{ch1_exam_bar2}
Необходимо вычислить
расстояние между реальными
географическими областями, например,
между двумя странами~--- Россией и Францией
\protect\footnotemark
(см. рис.~\ref{ch1_dr_wasser_2}).
\end{example}
\footnotetext{Франция выбрана не случайно: во-первых,
рассматриваемая в данном параграфе
тема барицентров Васерштейна
связана с теорией {\it оптимального
транспорта}, и далее приводится ссылка на основополагающую книгу ''Topics in Optimal Transportation'' французского математика,
лауреата Филдсовской
премии Седрика Виллани, 
во-вторых, в настоящее время (2020 год) Франция  на втором месте по количеству обладателей Филдсовской премии за всю историю
вручения премии с 1936~года, в-третьих, Франция, как и Россия, дала миру много замечательных оптимизаторов.
}

В дополнение к рассмотренным в примере~\ref{ch1_exam_bar1}
способам вычислить это расстояние
можно ещё предложить, например, такой способ~---
разделить обе области на маленькие одинаковые
квадраты и считать расстояния между такими
квадратами, затем усреднить полученные значения.
Понятно, что все эти способы дадут совершенно
разные ответы. 

Задачи из примеров~\ref{ch1_exam_bar1} и \ref{ch1_exam_bar2} можно объединить следующим
образом: как наиболее рационально переместить
кучу земли из одного участка (там образуется ров)
на другой участок (там нужно построить земляное
защитное укрепление~--- вал). Или наоборот~---
имеется куча земли и яма одинаковых объемов,
как с наименьшими усилиями засыпать яму землей?
Иногда эту задачу называют задачей землекопа,
или экскаваторщика.
Более точно задача формулируется так: 
требуется разбить две равновеликие плоские
или пространственные области  на бесконечно
малые участки. В примере~$1$: каждой $\xi \, \in \, A$
сопоставляется куча земли весом $p_\xi$,
каждой $\zeta \, \in \, B$ сопоставляется
ров глубиной $p_\zeta = 1 / |B|$.
Требуется сопоставить участки одной области
с участками другой так, чтобы сумма произведений
расстояний между участками (здесь можно взять любую
метрику) на площади или объёмы перемещающихся 
участков (т.е. затраченная работа, или усилия по перемещению,
или стоимость перемещения) была
бы минимальной. 
Ясно, что решение задачи определяется геометрией
рассматриваемых участков.
Строгую математическую
постановку задачи в терминах метрических пространств
и мер см., например, в \cite[с.~159]{MagTih_convan},
также см.~\cite{Kol_TranspConcent}.

В дискретной постановке рассмотренная
проблема превращается в классическую
{\it транспортную задачу}.

\begin{example}[{\it Транспортная задача}]
Имеются два набора точек~${x_i}$ (пункты
производства) и ${y_i}$ (пункты потребления), $1 \le i \le N$. Задана функция стоимости перевозки~$c(x, \, y)$ (например, какая-то
зависимость от расстояния между точками). 
Требуется построить взаимно однозначное отображение $T$,
сопоставляющее каждой точке из первого набора~${x_i}$ точку из второго набора ${y_i}$
так, чтобы суммарная стоимость перевозок $\sum_{i=1}^N c(x_i, \,  y_i)$ была минимальной.
\end{example}

Впервые транспортная задача (точнее,
её непрерывный аналог из примера~\ref{ch1_exam_bar2}~---
<<задача о рвах и насыпях>>) была поставлена
и исследована Гаспаром Монжем в 1781~году~\cite{Monge}.
Задача Монжа может не иметь решения.

\subsection{Релаксация Канторовича}
Существенным ограничением 
транспортной задачи в формулировке Монжа 
является возможность сравнения
гистограмм {\it только} одного и того же размера.
Транспортная задача о перемещении
масс в
релаксационной по сравнению с задачей Монжа постановке,
в которой допускается {\it расщепление
масс} (т.е. перемещение массы
из одного источника не в один,
а в несколько пунктов
назначения),
с целью нахождения оптимальных транспортных планов, была сформулирована и решена Л.В.~Канторовичем
в 40-х годах XX~века (см. \cite{Kant_mass42, Kant_mass48}) и получила в дальнейшем название {\it задачи Монжа--Канторовича}. 
Таким образом, Канторович перешёл от 
детерминированной постановки задачи
перемещения масс к вероятностной~--- когда
имеется возможность распределить массу из точки~${x_i}$ по нескольким пунктам
назначения~${y_i}$.
По
сути, это была первая
бесконечномерная задача линейного программирования.
\begin{defin}
Пусть $X$ и $Y$~--- два компактных
метрических пространства,
$\mu$ и $\nu$~--- 
борелевские
вероятностные меры,
т.е. меры,
определённые на всех открытых множествах,
с одной и той же вероятностной массой.
Задана непрерывная функция стоимости перевозки
$c \, : \, X \times Y \rightarrow \mathbb{R}$.
Оптимизационная задача Монжа--Канторовича
формулируется следующим образом:
\begin{equation}
\label{ch1_eq_MongeKant}
\inf_{\gamma \, \in \, \Pi(\mu, \, \nu)}
\int_{X \times Y} c(x, y) d \gamma (x, y)
\end{equation}
и заключается в нахождении оптимального
транспортного плана среди множества
всех транспортных планов $\Pi(\mu, \, \nu)$
(т.е. $\Pi(\mu, \, \nu)$~--- множество
борелевских вероятностных мер на $X \times Y$
с маргинальными для
$X$ и $Y$ соответственно функциями\protect\footnotemark $\mu$ и $\nu$).
\end{defin} 
\footnotetext{Это означает, что 
для любого борелевского множества $A_X \subset X$
и любого борелевского множества $B_Y \subset Y$
$$
\gamma(A_X \times Y) = \mu (A_X),
\gamma(X \times B_Y) = \nu(B_Y).
$$
}

При весьма общих предположениях задача Монжа--Канторовича имеет
решение.

Приведём практический пример задачи распределения
ресурсов, которая сводится к задаче Монжа--Канторовича.

\begin{example}[{\it Шахты и заводы} \cite{Comp_OT, Hitchcock}\protect\footnotemark]\label{ch1_ex_mines}
Предположим, что Центр владеет $m$~заводами и $n$~складами. На складах хранится сырьё, которое
требуется заводам для их работы. Конкретнее,
на $i$-м складе хранится $a_i$ условных единиц
(например, тонн) сырья, $i = 1, \, \ldots, \, n$. Заводу под
номером $j$ требуется $b_j$ тонн сырья,
$j = 1, \, \ldots, \, m$. Чтобы перевезти
сырье от складов на фабрики, Центру необходимо
заплатить транспортной компании, у которой
есть матрица стоимостей перевозок $C$,
где $C_{i, j}$~--- стоимость перевозки одной
тонны ресурса (сырья) от $i$-го склада до $j$-го завода. Зависимость стоимости перевозки от
количества перевозимого материала предполагается
линейной. Оптимальный
транспортный план перевозок~--- матрицу
$P^* \, \in \, \mathbb{R}_{+}^{n \times m}$,
где $P_{i, j}$ обозначает количество
перевозимого материала из пункта $i$ на завод $j$, ---
можно определить,
решив задачу~ЛП
$$
\min_{\sum\limits_j P_{i, j} = a, \, \sum\limits_i P_{i, j} = b}
\sum_{i, j} C_{i, j} P_{i, j}.
$$
\end{example}
\footnotetext{С.~Виллани в~\cite[Введение, часть~3]{Villani_2009} приводит аналогичный пример
с пекарнями и кафе, продающими выпечку.}

\subsection{Оптимальный транспорт}
Решением задач, подобных рассмотренным
выше, и их усложнённых аналогов, 
занимается  теория {\it оптимального 
транспорта} (ОТ) (op\-ti\-mal trans\-por\-ta\-tion или
op\-ti\-mal trans\-port\footnote{До начала 2000-х годов~--- mass trans\-por\-ta\-tion.}).
Уникальность теории~ОТ видна уже из того факта,
что её можно отнести сразу к трём математическим
дисциплинам~--- теории вероятностей, математическому
анализу и оптимизации. Основной целью~ОТ
является нахождение способов сравнить вероятностные
распределения с учётом
геометрии пространств,
на которых рассматриваются меры.

Наиболее подробное руководство
по данной теории, с учётом
последних достижений, написал уже
упомянутый выше математик Седрик Виллани~\cite{Villani_2003, Villani_2009} (в настоящее время вышли два издания, с разницей в 5~лет).
В последние десятилетия эта теория
стала очень популярной вследствие её успешных применений для решения целого ряда прикладных задач большой размерности,
в особенности, задач машинного
обучения (см., например, \cite{Wass_GAN}) и задач 
сравнительного анализа и обработки изображений~\cite{EMD_image, Wasser_DDGUN_0618, Conv_wasser15}. Популярность теории оптимального
транспорта продолжает расти в том числе и из-за
неожиданно открывающихся связей этой теории
с задачами в физике, астрономии, геометрии,
теории дифференциальных уравнений в частных производных и т.д.\footnote{И даже в политологии~\cite{emd_polit}.}
В начале 2018~года
вышла книга известных специалистов
в области~ОТ M.~Cuturi и G.~Peyr\'{e} ''Computational Optimal Transport''~\cite{Comp_OT}, акцентирующая внимание
на численных методах теории~ОТ и актуальных
приложениях к задачам анализа данных (книга
доступна для скачивания на arxiv.org).
Кстати, они же принимали участие в создании
свободно распространяемой 
библиотеки <<POT: Python Optimal Transport>>~\cite{POT_Python},
на языке программирования~Python, в которой
реализованы алгоритмы решения задач обработки
сигналов и изображений, машинного обучения,
относящиеся к теории оптимального транспорта
и приведён ряд наглядных примеров решения
учебных задач в Jupyter Notebook.

В теории оптимального транспорта
одной из важных задач является
задача вычисления {\it расстояния (метрики) Канторовича, или Васерштейна {\rm\cite{Was_69}} (в~англоязычной литературе <<Wasserstein distance>>)}\footnote{Название
<<расстояние Васерштейна>> уже закрепилось
в научной литературе, хотя вопрос о
приоритете открытия решается неоднозначно~\cite{Villani_2009}.
Здесь всё в соответствии с принципом Арнольда~\cite{arno_prin} и законом Стиглера:
\begin{enumerate}
    \item Такие
расстояния и их частные случаи в течение XX~века рассматривались
в работах разных авторов независимо 
друг от друга, в их числе и Канторович, и Хёфдинг,
и Фреше, и Рубинштейн, и, конечно, сам Леонид Васерштейн.
\item Название <<метрика Васерштейна>> было предложено Р.Л.~Добрушиным в 1970 году.
\item Сам Васерштейн рассматривал только случай $p = 1$~\cite{Was_69}.
\item Этот же случай рассматривали и Канторович совместно с Рубинштейном, поэтому частный случай метрики Васерштейна $W_1$ часто называют метрикой Канторовича--Рубинштейна.
\item Расстояние в задаче землекопа, Earth Mover’s distance (EMD), эквивалентно $W_1$, в случае
расчёта EMD между вероятностными распределениями.
\end{enumerate}}
~--- именно такого расстояния
(точнее, характеристики близости
двух вероятностных мер), которое
учитывает все упомянутые выше особенности
и трудности вычисления расстояния в пространстве
вероятностных мер. Одна из главных особенностей
расстояния Васерштейна~--- учёт геометрии метрического пространства
при вычислениях.
\begin{defin}
Пусть $\Omega$~--- компактное метрическое пространство
с метрикой~$c$, $M_1^+(\Omega)$~--- множество борелевских
вероятностных мер на $\Omega$, т.е. множество мер,
определённых на борелевской $\sigma$-алгебре множества~$\Omega$. Для $p \, \in \, [1, \, \infty)$
$p$-метрикой (расстоянием) Васерштейна
между $\mu \, \in \, M_1^+(\Omega)$ и $\nu \, \in \, M_1^+(\Omega)$ называется следующая величина
\protect\footnotemark:
$$
W_p(\mu, \, \nu) = \left ( \inf_{\gamma \, \in \, \Pi (\mu, \, \nu)} \int_{\Omega \times \Omega} c(x, \, y)^p d\gamma(x, \, y) \right )^{1/p}.
$$
\end{defin}
\footnotetext{Сравните с~\eqref{ch1_eq_MongeKant}.}

\subsubsection{Васерштейн и транспортная
задача}
Рассмотрим два множества точек
$X = (x_1, \, \ldots, \, x_n)$ и 
$Y = (y_1, \, \ldots, \, y_m)$, 
$X \, \in \, \Omega^n$, $Y \, \in \, \Omega^m$.
Пусть $\mu$, $\nu$~--- эмпирические
вероятностные меры, т.е. $\mu = \sum_{i=1}^n a_i \delta_{x_i}$,
$\nu = \sum_{i=1}^m b_i \delta_{y_i}$.
Для любого $x \, \in \, \Omega$ $\delta_x$~---
мера Дирака, сосредоточенная в точке~$x$.
Кроме того, набор $a = (a_1, \, \ldots, \, a_n)$
принадлежит единичному симплексу~$S_n(1)$: 
\begin{equation}
\label{ch1_eq_simplex_sn1}
S_n(1) \overset{\text{def}}{=} \left \{ a \, \in \, \mathbb{R}_+^n \, \left |  \, \sum_{i=1}^n a_i = 1 \right. \right \}.
\end{equation}
Тогда расстояние Васерштейна $W_p(\mu, \, \nu)$,
по сути, является оптимумом для
транспортной задачи (в степени $1/p$)
(см. пример~\ref{ch1_ex_mines} на стр.~\pageref{ch1_ex_mines}).

Важным достижением (активно использующимся сейчас на практике \cite{Comp_OT}) в разработке эффективных алгоритмов приближенного вычисления расстояния Васерштейна является энтропийная регуляризация (сглаживание), о которой подробнее будет написано в следующем разделе.


\subsection{Барицентры}
В задачах машинного обучения без учителя часто
возникает проблема вычисления некоторого <<среднего
арифметического>>~--- {\it барицентра}\footnote{Барицентром геометрической фигуры называют её геометрический центр, т.е. усреднённое положение всех точек фигуры по всем координатным направлениям. 
В физике термин <<барицентр>> обозначает
либо центр масс, либо центр тяжести физического объекта.} нескольких многомерных рядов данных \cite{Dvinskikh_OMS}.
В евклидовом случае барицентром точек~$(x_1, \, \ldots, \, x_s)$ с барицентрическими координатами
$(\lambda_1, \, \ldots, \, \lambda_s)$
является точка
\begin{equation}
\label{ch1_eq_barycentre}
\min_x \sum_{i = 1}^s \lambda_i d(x, \, x_i)^p,
\end{equation}
где
расстояние $d(x, \, x_i) = \| x - x_i \|$, $p = 2$.
В общем случае задача нахождения барицентра
обычно является задачей невыпуклой оптимизации.
Но в теории~ОТ, когда метрика~$d$ становится
расстоянием Васерштейна, задачу нахождения
барицентра Васерштейна можно сформулировать
как задачу выпуклой оптимизации,
для которой доказано существование
решения и разработаны эффективные численные
методы~\cite[глава~9]{Comp_OT}.

\begin{defin} Пусть
$\nu_i$, $i = 1, \, \dots, \, s$ и $\mu$~--- вероятностные меры
на симплексе~$S_n(1)$
(см.~\eqref{ch1_eq_simplex_sn1}).
Барицентр Васерштейна~\cite{agueh_bar_was_11} мер $\nu_i$ с барицентрическими координатами
$(\lambda_1, \, \ldots, \, \lambda_s)$ является
решением задачи оптимизации:
\begin{equation}
\label{ch1_eq_bar_wasser}
\min_{\mu \, \in \, S_n(1)} \sum_{i = 1}^s \lambda_i W_p(\mu, \, \nu_i)^p.
\end{equation}
\end{defin}

Итак, барицентр Васерштейна~--- это усреднённый представитель всех рассматриваемых объектов
(элемент пространства мер) по минимальной сумме расстояний Васерштейна от каждого
из объектов выборки. Геометрическую структуру всех
объектов барицентр Васерштейна обобщает
гораздо лучше усреднённых по другим 
метрикам представителей~\cite{Cuturi_bar_2014, Comp_OT}.

\begin{exercise}\label{ch1_exlist_was_entr}
Рассмотрим частный случай конструкции~\eqref{ch1_eq_bar_wasser} при
$\lambda_i \equiv 1$, $i=1, \, \ldots, \, s$,
$p=2$. Расстояние Васерштейна возьмём энтропийно-сглаженное~\cite{Cuturi_bar_2014}: 
\begin{equation}
\label{ch1_eq_entr_wasser}
W_{2,r}(\mu, \, \nu)^2 = \mathop {\min }\limits_{\begin{array}{c}
 \sum\limits_{j=1}^n {x_{ij} } = \mu_i , \; \sum\limits_{i=1}^n {x_{ij} } = \nu_j \\ 
 x_{ij} \ge 0,\;i,j=1,...,n \\ 
 \end{array}} \!\!\!\left\{ {
 \sum\limits_{i,j=1}^n {c_{ij} x_{ij} } + 
 r \sum\limits_{i,j=1}^n {x_{ij} \log x_{ij} } 
 } \right\}.
\end{equation}
Здесь матрица стоимостей~$C=\{c_{ij}\}_{i, \, j = 1}^{n, \, n}$ сформирована из квадратов
попарных расстояний от носителя меры~$i$ до носителя
меры~$j$.
Перейдём от~\eqref{ch1_eq_entr_wasser} к двойственной
задаче:
$$
W_{2,r}(\mu, \, \nu)^2 = \mathop {\max }\limits_{\lambda, \, \beta } \left\{ {
{\lambda^T \mu} 
 +  \beta^T \nu  - r 
\sum\limits_{i,j=1}^n {\exp \left( {\frac{-c_{ij} +\lambda_i +\beta_j 
}{r }-1} \right)} } \right \} =
$$
\begin{equation}
\label{ch1_eq_was_dv_prob}
=\mathop {\max }\limits_\lambda \left \{ {  \lambda^T \mu  }  -
 r \sum\limits_{j=1}^n \nu_j 
\log \left( \frac{1}{\nu_j } \sum\limits_{i=1}^n \exp \left( \frac{-c_{ij} 
+\lambda _i }{r } \right )  \right )  
\right \}. 
\end{equation} 
Определим при $\mu \, \in \, S_n(1)$
функцию $H_\nu(\mu) = W_{2,r}(\mu, \, \nu)^2$.

По формуле Демьянова--Данскина--Рубинова (см. п.~\ref{ch1_sect_dem-dan})
для $\mu \, \in \, S_n( 1)$ 
функция $H_\nu(\mu)$~--- гладкая с градиентом 
$$
\nabla H_\nu(\mu) =\lambda ^\ast ,
$$
где $\lambda ^\ast $~--- единственное
решение~\eqref{ch1_eq_was_dv_prob}, удовлетворяющее условию $\sum_{i=1}^n \lambda^{\ast}_i = 0$.

Найдём сопряжённую функцию по Юнгу--Фенхелю:
\begin{equation}\label{conj}
H_\nu^*(\lambda) = \mathop {\max }\limits_{\mu \, \in \, S_n(1) } \left\{ { \lambda^T \mu  - H_\nu ( \mu )} \right\} =
r \sum\limits_{j=1}^n {\nu_j \log \left( {\frac{1}{\nu_j 
}\sum\limits_{i=1}^n {\exp \left( {\frac{-c_{ij} +\lambda _i }{r }} 
\right)} } \right)} 
\end{equation}
(см.~\eqref{ch1_eq_was_dv_prob}). Отметим, что эта функция будет иметь константу Липшица градиента в 2-норме, равную $r^{-1}$.

Перепишите задачу вычисления барицентра Васерштейна 
$$
\min_{\mu \, \in \, S_n(1)} \sum_{i = 1}^s  H_{\nu_i}(\mu) 
$$ 
следующим образом:
\begin{align*}
 \mathop {\min }\limits_{\begin{array}{c}
\mu_1 = \mu_s,...,\mu_{s-1} = \mu_s; \\ \mu_1,...,\mu_s \in S_n(1)
 \end{array}}
\sum_{i = 1}^s  H_{\nu_i}(\mu_i) = \\ = -\mathop {\max }\limits_{\begin{array}{c}
\mu_1 = \mu_s,...,\mu_{s-1} = \mu_s; \\ \mu_1,...,\mu_s \in S_n(1)
 \end{array}}
\sum_{i = 1}^s - H_{\nu_i}(\mu_i). 
\end{align*}
Последнюю задачу удобно рассматривать без знака минус
\begin{equation*}
\mathop {\max }\limits_{\begin{array}{c}
\mu_1 = \mu_s,...,\mu_{s-1} = \mu_s; \\ \mu_1,...,\mu_s \in S_n(1)
 \end{array}}
\sum_{i = 1}^s - H_{\nu_i}(\mu_i). 
\end{equation*}
Воспользуйтесь принципом множителей Лагранжа и перенесите ограничения в функционал
$$
 \mathop {\max }\limits_{\begin{array}{c}
\mu_1 = \mu_s,...,\mu_{s-1} = \mu_s; \\ \mu_1,...,\mu_s \in S_n(1)
 \end{array}}
\sum_{i = 1}^s  -H_{\nu_i}(\mu_i) = 
$$
$$
=\max_{\mu_1,...,\mu_s \in S_n(1)} \left[\sum_{i = 1}^{s-1} \min_{\lambda_i}\left\{\lambda_i^T\mu_i - H_{\nu_i}(\mu_i)\right\} +\right.
$$
$$
\left. + \min_{\lambda_1,...,\lambda_{s-1}}\left\{\left(-\sum_{j=1}^{s-1}\lambda_j\right)^T\mu_s- H_{\nu_s}(\mu_s)\right\}\right].
$$
Последнюю задачу можно переписать следующим образом (\textit{минмакс} теорема фон-Неймана--Сиона--Какутани о возможности перестановочности минимума и максимума для выпукло-вогнутых задач с компактными ограничениями по одной из групп переменных, см.~п.~\ref{subs:kakutani}):
\begin{align*}
\min_{\lambda_1,...,\lambda_{s-1}} \left[\sum_{i = 1}^{s-1} \max_{\mu_i\in S_n(1)}\left\{\lambda_i^T\mu_i - H_{\nu_i}(\mu_i)\right\} + \right.\\ \left. + \max_{\mu_s\in S_n(1)}\left\{\left(-\sum_{j=1}^{s-1}\lambda_j\right)^T\mu_s- H_{\nu_s}(\mu_s)\right\}\right].
\end{align*}
Используя \eqref{conj} и формулу Демьянова--Данскина--Рубинова (см. п.~\ref{ch1_sect_dem-dan}), можно заметить, что для $i=1,...,s-1$:
\begin{equation}\label{WB}
\ag{\mu_i=\nabla_{\lambda_i}\max_{\mu_i\in S_n(1)}\left\{\lambda_i^T\mu_i - H_{\nu_i}(\mu_i)\right\} = \nabla H_{\nu_i}^{\ast}\left(\lambda_i\right)\nk{.}}
\end{equation}
\nk{Поскольку} в решении $\mu_{\ag{1}}^{\ast} = ... = \mu_s^{\ast} = \mu^{\ast},$
где $\mu^{\ast}$ -- искомый барицентр, \ag{то $\mu^{\ast} = \nabla H_{\nu_i}^{\ast}\left(\lambda_i^{\ast}\right)$ для любого $i=1,...,s-1$.}
Итак, мы пришли к задаче
$$
\max_{\lambda_1,...,\lambda_{s-1}} \sum_{i = 1}^{s-1} H_{\nu_i}^{\ast}(\lambda_i) + H_{\nu_s}^{\ast}\left(-\sum_{j=1}^{s-1}\lambda_j\right), 
$$
решение которой позволяет восстанавливать барицентр согласно \eqref{WB}.
Последняя задача является задачей гладкой выпуклой безусловной оптимизации, которую можно решать, например, ускоренными градиентными методами (см. п.~\ref{AM}), являющимися прямо-двойственными \cite{Gas_Pos18}, т.е. по решению двойственной задачи можно восстанавливать решение прямой (исходной) задачи с той же точностью, используя всю последовательность, генерируемую методом в двойственном пространстве и формулу связи прямых и двойственных переменных \eqref{WB}. На таком пути строятся эффективные (в том числе и распределённые) алгоритмы поиска барицентра Васерштейна ~\cite{GasDvu_TranspOptima1516}, \cite{Wasser_DDGUN_0618}. Наилучший на данный момент в теоретическом плане алгоритм поиска настоящего (нерегуляризованного) барицентра Васерштейна  был получен недавно выпускницей школы ПМИ МФТИ Дариной Двинских и студентом ФКН ВШЭ Даниилом Тяпкиным \cite{Dvinskikh_AISTATS}. Сложность этого алгоритма (число арифметических операций) для нахождения $\varepsilon$-барицентра Васерштейна (по функции) будет $\tilde{O}\left(sn^2/\varepsilon\right)$, что в $\sqrt{n}$ раз лучше описанного выше способа, при  выборе параметра регуляризации $r\sim\varepsilon$, что необходимо для аппроксимации исходной \nk{задачи регуляризованной}. По-видимому, оценку Двинских--Тяпкина уже нельзя в дальнейшем будет улучшить более чем на логарифмические множители.
\end{exercise}

\section{Релаксации полиномиальных задач} \label{ch3_polynomial_problems}

В этом разделе мы рассмотрим полуопределённые релаксации оптимизационных задач с полиномиальными данными. Здесь можно выделить два концептуально разных, но тесно связанных между собой подхода, а именно релаксации типа сумм квадратов (sum of squares --- SOS) и моментные релаксации.

Тематика сумм квадратов восходит к Гильберту, но в применении к оптимизации путём полуопределённых релаксаций, насколько нам известно, впервые упомянута в статье Ю.Е.~Нестерова \cite{NesterovSOS} и диссертации П.~Паррило \cite{ParriloThesis} в 2000~г. В рассматриваемых в этом контексте задачах искомым объектом является \emph{полином}, представленный вектором своих коэффициентов, а в качестве конических ограничений выступают условия неотрицательности этого полинома либо на пространстве $\mathbb R^n$, либо на множествах, задаваемых полиномиальными ограничениями. По сути, задача ставится в некотором конечномерном пространстве полиномиальных функций.

Во втором подходе искомым объектом является вектор \emph{моментов} неотрицательных мер на некотором множестве $X \subset \mathbb R^n$, задающемся полиномиальными ограничениями. Этот подход более гибкий с точки зрения возможностей \nk{построения полуопределённых релаксаций}. Его теория и приложения подробно описаны в являющейся ныне стандартом в этой области монографии Ж.-Б.~Лассера \cite{LasserreMoments}.

\subsection[Выпуклая формулировка общей задачи оптимизации]{Выпуклая формулировка общей задачи\\ оптимизации}

Чтобы пояснить идею, стоящую за моментным подходом к релаксации полиномиальных задач, сперва рассмотрим, как \emph{произвольную} задачу оптимизации можно эквивалентно записать в виде выпуклой задачи. Рассмотрим общую проблему оптимизации
\[ \min_{x \in X}\ f(x)
\]
с допустимым множеством $X \subset \mathbb R^n$ и функцией цены $f$. В качестве требований мы предъявим лишь измеримость множества $X$ и непрерывность функции $f$. Задачу можно переписать в виде
\[ \min_{\mu \in {\cal M}}\ \int_X f(x)\,d\mu(x)\ :\quad \int_X d\mu(x) = 1,
\]
где ${\cal M}$ --- конус неотрицательных мер на множестве $X$. На этом бесконечномерном конусе мы минимизируем линейный функционал с одним линейным ограничением. Другими словами, мы минимизируем \nk{математическое} ожидание значения функции цены $f$ по всем вероятностным мерам на допустимом множестве~$X$.

Отметим концептуальное сходство этой формулировки с полуопределённой релаксацией квадратичных задач. В обоих случаях допустимое множество представляется в виде невыпуклого множества в некотором векторном пространстве, причем таким образом, что функция цены описывается неким линейным функционалом. В случае полуопределённой релаксации точки $x \in \mathbb R^n$ допустимого множества представляются матрицами $xx^T \in \mathbb{S}_+^n$ ранга 1, а в рассматриваемом выше случае --- дельта-функциями с носителем в точке $x$. Далее допустимое множество расширяется до некоторого выпуклого множества, содержащего выпуклую оболочку образа исходного множества $X$. В случае релаксации квадратичных задач получаем матричный конус $\mathbb{S}_+^n$, пересеченный с линейными ограничениями исходной задачи, а в вышеописанном случае --- множество вероятностных мер. Отметим, что дельта-функции являются в точности экстремальными точками множества вероятностных мер.

Решение выпуклой формулировки в пространстве мер наталкивается на ряд принципиальных трудностей. Во-первых, это пространство и определённый в нем конус неотрицательных мер бесконечномерны. Для дизайна численных алгоритмов нужно свести задачу к конечномерной. Если множество $X$ и функция $f$ обладают полиномиальным описанием, то этот шаг удаётся сделать без дополнительных аппроксимаций. В этом случае все возникающие выражения, зависящие от меры, описываются через конечное число моментов этой меры. Этим объясняется, что данный подход на практике весьма успешно применяется к полиномиальным задачам.

Вторая трудность заключается в том, что даже конечномерный образ конуса неотрицательных мер в общем случае не имеет численно эффективно доступного описания, и приходится прибегать к релаксациям. Тем не менее в ряде важных случаев эти релаксации оказываются точными, что позволяет эффективно решать соответствующие классы задач.

Перед тем, как перейти к описанию метода моментов для полиномиальных задач, мы представим интуитивно более доступную и исторически гораздо более раннюю концепцию сумм квадратов.

\subsection{Положительные полиномы и суммы квадратов}

Теория, изложенная в этом разделе, применима к коническим программам над конусами неотрицательных полиномов. Мы рассмотрим полуопределённые аппроксимации, а в некоторых случаях и полуопределённые представления этих конусов.

Обозначим конус неотрицательных однородных полиномов степени $d$ от $n$ вещественных переменных через $P_{d,n}$. Ясно, что $P_{d,n}$ содержит ненулевой элемент тогда и только тогда, когда $d$ чётно. В этом случае $P_{d,n}$ является регулярным выпуклым конусом. Полином $p \in P_{d,n}$ является внутренней точкой конуса тогда и только тогда, когда он строго положителен на единичной сфере $S^{n-1} \subset \mathbb R^n$. В противном случае он должен обращаться в ноль хотя бы в одной точке $x \in S^{n-1}$.

\begin{exercise} \label{example_copositive}
Пусть $p \in \partial P_{d,n}$ --- полином на границе конуса неотрицательных полиномов. Построить подпирающую гиперплоскость к конусу $P_{d,n}$ в точке $p$, т.е., гиперплоскость $H$ такую, что $p \in H$, а конус $P_{d,n}$ целиком лежит в одном из двух замкнутых полупространств, определяемых $H$.
\end{exercise}

Рассмотрим в качестве примера конус ${\cal COP}^n$ коположительных матриц. Симметрическая матрица $A \in \mathbb{S}^n$ называется \emph{коположительной}, если для всех неотрицательных векторов $x \in \mathbb R_+^n$ выполнено неравенство $x^TAx \geq 0$. Ясно, что матрица $A$ коположительна тогда и только тогда, когда полином четвёртой степени
\[ p_A(x) = \sum_{i,j = 1}^n A_{ij}x_i^2x_j^2
\]
является неотрицательным, т.е., элементом конуса $P_{4,n}$. Таким образом, конус ${\cal COP}^n$ обладает представлением через линейное сечение конуса неотрицательных полиномов $P_{4,n}$. Однако установить принадлежность произвольной матрицы $A \in \mathbb{S}^n$ к конусу ${\cal COP}^n$ является ко-NP-полной задачей \cite{MurtyKabadi87}.

Поэтому решить, принадлежит ли данный полином $p(x)$ к конусу $P_{2d,n}$ или нет, в общем случае является сложной проблемой. Основной алгоритмически доступной аппроксимацией конуса неотрицательных полиномов является конус сумм квадратов.

\begin{defin}
Конусом сумм квадратов $\Sigma_{2d,n}$ называется множество однородных полиномов $p(x)$ чётной степени $2d$ от $n$ вещественных переменных, которые представляются в виде конечной суммы $p(x) = \sum_{k=1}^m q_k^2(x)$, где $q_k(x)$, $k = 1,\dots,m$~--- однородные полиномы степени $d$.
\end{defin}

Ясно, что $\Sigma_{d,n} \subset P_{d,n}$, и $\Sigma_{d,n}$ является \emph{внутренней} аппроксимацией $P_{d,n}$.

\medskip

Построим полуопределённое представление конуса $\Sigma_{2d,n}$. Для этого введём вектор ${\bf x}$ всех мономов вида $x^{\alpha} := \prod_{k=1}^n x_k^{\alpha_k}$ степени $d$. Здесь $\alpha = (\alpha_1, \, \dots, \, \alpha_n) \in \mathbb N^n$ --- мульти-индекс, в котором сумма $|\alpha| := \sum_{k=1}^n \alpha_k$ экспонент равна $d$. Размерность $\begin{pmatrix} n+d-1 \\ d \end{pmatrix}$ вектора ${\bf x}$ обозначим через $N$.

\begin{lemma}
Однородный полином $p$ степени $2d$ от $n$ переменных $x_1, \, \dots, \, x_n$ является элементом конуса $\Sigma_{2d,n}$ тогда и только тогда, когда существует неотрицательно определённая матрица $A \in \mathbb{S}_+^N$ такая, что $p(x) = {\bf x}^TA{\bf x}$.
\end{lemma}

\begin{proof}
Пусть существует матрица $A \in \mathbb{S}_+^N$ такая, что $p(x) = {\bf x}^TA{\bf x}$. Представим $A$ в виде произведения $A = B^TB$, где $B \in \mathbb R^{k \times N}$. Обозначим строки фактора $B$ через $b_j$, $j = 1, \, \dots, \, k$. Тогда имеем
\[ p(x) = {\bf x}^TB^TB{\bf x} = \sum_{j=1}^k \langle b_j,{\bf x} \rangle^2,
\]
и $p$ представлено в виде суммы квадратов $k$ полиномов $q_j(x) = \langle b_j,{\bf x} \rangle$ степени $d$. Отсюда следует включение $p \in \Sigma_{2d,n}$.

Теперь допустим, что $p \in \Sigma_{2d,n}$. Тогда существуют $k$ однородных полиномов $q_1(x), \, \dots, \, q_k(x)$ степени $d$ таких, что $p(x) = \sum_{j=1}^k q_j^2(x)$. Полином $q_j(x)$ можно записать как скалярное произведение $\langle c_j,{\bf x} \rangle$, где $c_j \in \mathbb R^N$ --- его вектор коэффициентов. Составим матрицу $C \in \mathbb R^{k \times N}$ таким образом, что векторы $c_j$ являются её строками. Тогда получим
\[ p(x) = \sum_{j=1}^k \langle c_j,{\bf x} \rangle^2 = {\bf x}^TC^TC{\bf x},
\]
и полином $p(x)$ представлен в виде произведения ${\bf x}^TA{\bf x}$ с неотрицательно определённой матрицей $A = C^TC$. 
\end{proof}

Из леммы вытекает, что конус сумм квадратов $\Sigma_{2d,n}$ линейно изоморфен проекции матричного конуса $\mathbb{S}_+^N$. Зададим эту проекцию в явном виде. Напомним, что коэффициенты полиномов $p \in \Sigma_{2d,n}$ индексируются мульти-индексами $\alpha \in \mathbb N^n$ с суммой $2d$, в то время как элементы вектора ${\bf x}$ индексируются мульти-индексами $\beta \in \mathbb N^n$ с суммой $d$. Определим линейную проекцию $\Pi$, сопоставляющую каждой матрице $A \in \mathbb{S}^N$ полином $p(x)$ с коэффициентами
\begin{equation} \label{lin_SOS}
c_{\alpha} = \sum_{\beta,\gamma: \beta+\gamma = \alpha} A_{\beta\gamma}.
\end{equation}
Тогда $\Pi(A) = {\bf x}^TA{\bf x}$, и из леммы следует представление $\Sigma_{2d,n} = \Pi\left[\mathbb{S}_+^N\right]$.

Технически для описания условия $p \in \Sigma_{2d,n}$ достаточно ввести дополнительную матричную переменную $A \in \mathbb{S}_+^N$ и связать её элементы с коэффициентами $p$ через линейные условия \eqref{lin_SOS}. Применимость релаксаций типа сумм квадратов ограничена сложностью полуопределённого условия $A \succeq 0$, поскольку порядок $N$ матрицы очень быстро увеличивается с ростом $n$ и $d$. С другой стороны, линейные условия \eqref{lin_SOS} на элементы $A$ задаются разреженными матрицами. В последние годы было предложено аппроксимировать условие $A \succeq 0$ более простыми линейными или конично-квадратичными ограничениями \cite{AhmadiMajumbar19}.

\medskip

Конуса неотрицательных полиномов и сумм квадратов можно аналогичным образом определить в произвольных конечномерных линейных пространствах полиномов, в частности, в пространствах неоднородных полиномов от $n$ переменных степени, не превосходящей $d$.

\medskip

В ряде случаев релаксация конуса положительных полиномов суммами квадратов точна. Когда именно равенство $P_{2d,n} = \Sigma_{2d,n}$ имеет место, было известно ещё в 19-ом веке.

\begin{teo}
Пусть $n,d \in \mathbb N_+$, $d$ чётное. Равенство $\Sigma_{d,n} = P_{d,n}$ справедливо тогда и только тогда, когда либо $\min(d,n) \leq 2$, либо $(d,n) = (4,3)$.
\end{teo}

Во всех остальных случаях конус $P_{2d,n}$ не является полуопределённо представимым \cite{Scheiderer18a}.

\medskip

Для того, чтобы проверить существование представления в виде суммы квадратов для конкретного полинома, не обязательно работать в пространстве матриц $\mathbb{S}^N$. Введём следующее понятие.

\begin{defin}
Пусть $p(x) = \sum_{\alpha} c_{\alpha}x^{\alpha}$ --- полином от $n$ переменных. \emph{Политопом Ньютона} полинома $p$ называют выпуклую оболочку множества мульти-индексов $\{ \alpha \in \mathbb N^n \mid c_{\alpha} \not= 0 \}$.
\end{defin}

\begin{lemma}
Пусть $p(x) = \sum_{\alpha} c_{\alpha}x^{\alpha}$ --- полином от $n$ переменных, представимый в виде суммы квадратов $\sum_j q_j(x)^2$, и пусть $P$ --- его политоп Ньютона. Тогда справедливы следующие утверждения:
\begin{itemize}
\item экстремальные точки политопа $P$ имеют только чётные компоненты,
\item коэффициенты $c_{\alpha}$, соответствующие экстремальным точкам $\alpha$ политопа $P$, строго положительны,
\item мульти-индексы $\beta$, соответствующие ненулевым коэффициентам полиномов $q_j$, являются элементами политопа $\frac12P$.
\end{itemize}
\end{lemma}

В частности, если полином $p$ --- разреженный, то часто для проверки включения $p \in \Sigma_{2d,n}$ возможно обойтись рассмотрением матриц размера, гораздо меньшего, чем $N$.

Стандартным примером неотрицательного полинома, не представимого в виде суммы квадратов других полиномов, является \emph{полином Моцкина} $p_M(x,y,z) = x^4y^2 + x^2y^4 + z^6 - 3x^2y^2z^2$. Этот полином является элементом конуса $P_{6,3}$ вследствие неравенства между алгебраическим и геометрическим средним.

\begin{exercise}
Доказать, что полином Моцкина не является элементом конуса $\Sigma_{6,3}$.
\end{exercise}

Тем не менее, неотрицательность полинома Моцкина можно все же доказать с помощью представления его в виде суммы квадратов. Но для этого нужно расширить или изменить класс функций, которые могут выступать в роли факторов $q_j$. А именно, имеем представление
\[ p_M(x,y,z) = (u^2+v^2+z^2) \cdot \begin{pmatrix} u^2 \\ v^2 \\ z^2 \end{pmatrix}^T \begin{pmatrix} 1 & -0.5 & -0.5 \\ -0.5 & 1 & -0.5 \\ -0.5 & -0.5 & 1 \end{pmatrix} \begin{pmatrix} u^2 \\ v^2 \\ z^2 \end{pmatrix},
\]
где $u = x^{2/3}y^{1/3}$, $v = x^{1/3}y^{2/3}$. То, что правая часть является суммой квадратов, вытекает из соотношения $J = \left(\sqrt{\frac23}J\right)^2$, которому удовлетворяет $3 \times 3$ матрица коэффициентов $J$, фигурирующая в формуле.

Сертификат неотрицательности полинома Моцкина можно также получить, разложив произведение $(x^2 + y^2 + z^2) \cdot p_M(x,y,z)$ в сумму квадратов полиномов 4-й степени. Но тогда и произведение $(x^2 + y^2 + z^2)^2 \cdot p_M(x,y,z)$ является суммой квадратов полиномов, и $p_M$ представляемо в виде суммы квадратов от рациональных функций.

Справедлив следующий общий результат.

\medskip

\emph{17-я проблема Гильберта:} Любой неотрицательный полином может быть представлен в виде суммы квадратов \emph{рациональных} функций.

\medskip

Утверждение было доказано Э.~Артином в 20-х гг. \nk{XX века.}

\begin{example}
Рассмотрим коположительный конус ${\cal COP}^n$ всех квадратичных форм $A \in \mathbb{S}^n$, неотрицательных на неотрицательном ортанте. Выше мы охарактеризовали включение $A \in {\cal COP}^n$ условием $p_A(x) = \sum_{i,j=1}^n A_{ij}x_i^2x_j^2 \in P_{4,n}$. Отсюда вытекает, что множество
\[ {\cal K}_0 = \{ A \in \mathbb{S}^n \,|\, p_A \in \Sigma_{4,n} \}
\]
является внутренней полуопределённой релаксацией конуса ${\cal COP}^n$. Вычислим эту релаксацию в явном виде.

Сформируем вектор ${{\bf x} = (x_1^2,  \dots,  x_n^2,x_1x_2,x_1x_3,  \dots,  x_{n-1}x_n)^T \in \mathbb R^N}$, где $N = \frac{n(n+1)}{2}$. Тогда условие $A \in {\cal K}_0$ эквивалентно существованию матрицы ${\bf A} \in \mathbb{S}_+^N$ такой, что $p_A(x) = {\bf x}^T{\bf A}{\bf x}$. Разобъём ${\bf A}$ на 4 блока, отвечающих разбиению ${\bf x}$ на подвектор ${\bf x}^1$ длины $n$ и подвектор ${\bf x}^2$ длины $\frac{n(n-1)}{2}$. Тогда коэффициенты при мономах $x_i^4$ и $x_i^2x_j^2$ в полиноме ${\bf x}^T{\bf A}{\bf x}$ зависят только от элементов блока ${\bf A}_{11}$ и диагональных элементов блока ${\bf A}_{22}$. Поэтому достаточно ограничиться неотрицательной определённостью блочно-диагональной подматрицы $\diag(B,c_{12},c_{13}, \, \dots, \, c_{n-1,n})$, составленной из этих элементов, а все остальные элементы положить равными нулю. Таким образом, условие $A \in {\cal K}_0$ эквивалентно существованию матрицы $B \in \mathbb{S}_+^n$ и скаляров $c_{ij} \geq 0$, $1 \leq i < j \leq n$, таких что
\[ \sum_{i,j=1}^n A_{ij}x_i^2x_j^2 = \sum_{i,j=1}^n B_{ij}x_i^2x_j^2 + \sum_{i < j}c_{ij}x_i^2x_j^2.
\]
Сравнивая коэффициенты, мы получаем условия $\diag A = \diag B$ и $A_{ij} = B_{ij} + c_{ij}$ для всех $i < j$. Отсюда получаем явное представление ${\cal K}_0 = \mathbb{S}_+^n + {\cal N}^n$,
где ${\cal N}^n$ --- конус поэлементно неотрицательных матриц с нулевой диагональю. Легко видеть, что условие на диагональ можно опустить. 

Можно показать, что равенство ${\cal COP}^n = \mathbb{S}_+^n + {\cal N}^n$ имеет место тогда и только тогда, когда $n \leq 4$.

\medskip

Внутренняя релаксация ${\cal K}_0$ может быть усилена. Определим иерархию конусов ${\cal K}_0 \subset {\cal K}_1 \dots$, параметризованную целым числом $r \geq 0$:
\[ {\cal K}_r = \left\{ A \in \mathbb{S}^n \,\left|\, \left(\sum_{j=1}^n x_j^2\right)^r \cdot p_A(x) \in \Sigma_{4+2r,n} \right. \right\}.
\]
Это так называемая {\it иерархия релаксаций Паррило} для коположительного конуса. Сложность релаксации быстро увеличивается с ростом $r$. Справедлив следующий результат \cite{ParriloThesis}:

\begin{teo}
Пусть $A \in \interior {\cal COP}^n$. Тогда существует $r \geq 0$ такое, что $A \in {\cal K}_{r'}$ для всех $r' \geq r$.
\end{teo}
\end{example}

Теперь мы перейдём к важному для вычислительной практики случая \emph{матричных} полиномов от \emph{одной} переменной, для которого релаксация положительных полиномов суммами квадратов точна.

\begin{teo}
Пусть $A_0, \, \dots, \, A_{2m} \in \mathbb{S}^n$ --- симметричные матрицы такие, что матричный полином $A(t) = \sum_{j=0}^{2m} A_jt^j$ степени $2m$ положительный, т.е., $x^TA(t)x \geq 0$ для всех $t \in \mathbb R$ и $x \in \mathbb R^n$. Тогда найдутся матрицы $B_0, \, \dots, \, B_m \in \mathbb R^{n \times k}$ такие, что $A(t) = B(t)B(t)^T$, где $B(t) = \sum_{j=0}^m B_jt^j$ --- матричный полином степени $m$.
\end{teo}

Доказательство теоремы разобьём на несколько этапов.

\begin{lemma} \label{lem:blockHankel}
Пусть $H \in \mathbb{S}_+^{n(m+1)}$ --- неотрицательно определённая блочно-ханкелевая матрица с блоками размера $n$. Тогда $H$ представляется в виде суммы неотрицательно определённых блочно-ханкелевых матриц ранга 1.
\end{lemma}

\begin{proof}
Профакторизуем $H = FF^T$ так, что $F \in \mathbb R^{n(m+1) \times l}$ состоит из блоков $F_0, \, \dots, \, F_m \in \mathbb R^{n \times l}$. Умножением $F$ справа на ортогональную матрицу можно добиться того, что $F_i$ разделится на подблоки $F_i = (F_{i1},F_{i2})$ размера $n \times l_1$ и $n \times l_2$, соответственно, так чтобы подматрица $F_I \in \mathbb R^{nm \times l_2}$ матрицы $F$, состоящая из блоков $F_{02}, \, \dots, \, F_{m-1,2}$, имела полный ранг, а $F_{01} = \dots = F_{m-1,1} = 0$. Определим также подматрицу $F_{II} \in \mathbb R^{nm \times l_2}$ матрицы $F$, состоящую из блоков $F_{12}, \, \dots, \, F_{m2}$.

Так как $H$ --- блочно-ханкелевая, её верхний правый угол размера $mn \times mn$ совпадает с нижним левым углом того же размера. Из этого следует соотношение $F_IF_{II}^T = F_{II}F_I^T$. Так как $F_I$ имеет полный ранг, то образ $F_{II}$ содержится в образе $F_I$, и найдётся матрица $\Lambda \in \mathbb R^{l_2 \times l_2}$ такая, что $F_{II} = F_I\Lambda$. Далее из соотношения $F_IF_{II}^T = F_{II}F_I^T$ следует, что $\Lambda$ симметричная. Пусть $\Lambda = UDU^T$ --- спектральное разложение матрицы $\Lambda$, где $U$ --- ортогональная, а $D$ --- диагональная матрица. Заменив $F_I,F_{II}$ на $F_IU,F_{II}U$, можно считать, что $\Lambda = D = \diag(d_{l_1+1}, \, \dots, \, d_l)$.

Таким образом, первые $l_1$ столбцов $f_j$ фактора $F$ имеют вид $(0; \, 0; \, \dots; \, 0; \, x_j)$, $j = 1,\, \dots, \, l_1$, а последние $l_2$ столбцов имеют вид $(x_j; \, d_jx_j; \, d_j^2x_j; \, \dots; \, d_j^mx_j)$, $j = l_1+1, \, \dots, \, l$. Из этого следует, что матрицы $f_jf_j^T$ ранга 1 сами блочно-ханкелевы. 
\end{proof}

\begin{lemma}
Пусть $A(t) = \sum_{j=0}^{2m} A_jt^j$ --- положительный матричный полином, а $H \in \mathbb{S}_+^{n(m+1)}$ --- блочно-ханкелевая матрица с блоками ${H_0,  \dots,  H_{2m}  \in  \mathbb{S}^n}$. Тогда $\sum_{i=0}^{2m} \langle A_i,H_i \rangle \geq 0$. В частности, существует такая матрица ${\bf A}\!\in\linebreak\in  \mathbb{S}_+^{n(m+1)}$, состоящая из $n \times n$ блоков $A_{ij}$, $i, \, j = 0, \, \dots, \, m$, что $A_k = \sum_{i+j=k} A_{ij}$.
\end{lemma}

\begin{proof}
Ввиду предыдущей леммы можно без ограничения общности считать, что ранг $H$ равен 1. Но тогда либо у $H$ ненулевой всего лишь правый нижний угол размера $n$, либо найдётся вектор $x \in \mathbb R^n$ и число $\lambda \in \mathbb R$ такое, что $H_i = \lambda^ixx^T$ для всех $i = 0, \, \dots, \, 2m$. В обоих случаях выражение $\sum_{i=0}^{2m} \langle A_i,H_i \rangle$ будет неотрицательно определённой одноранговой матрицей вследствие условия на $A$. В первом случае, потому что старший коэффициент $A_{2m}$ неотрицательно определённый, а во втором, потому что выражение имеет вид $\left(\sum_{i=0}^{2m} A_i\lambda^i\right)xx^T$.

Второе утверждение получается применением конической двойственности. Раз набор симметричных матриц $(A_0, \, \dots, \, A_{2m})$ неотрицателен на любом наборе симметричных матриц $(H_0, \, \dots, \, H_{2m})$, представляющих блочно-ханкелевую неотрицательно определённую матрицу, то $(A_0, \, \dots, \, A_{2m})$ является элементом конуса, двойственного к конусу блочно-ханкелевых неотрицательно определённых матриц. Но этот двойственный конус характеризуется в точности существованием неотрицательно определённой матрицы ${\bf A}$ с требуемыми свойствами. 
\end{proof}

Теперь теорема доказывается факторизацией ${\bf A} = {\bf B}{\bf B}^T$ матрицы ${\bf A}$, где ${\bf B} \in \mathbb R^{n(m+1) \times k}$. Искомые коэффициенты $B_j$ получаются разбиением ${\bf B}$ на $m+1$ блоков. \qed

\medskip

В этом разделе мы рассматривали полуопределённые аппроксимации отдельных конусов положительных полиномов. Далее мы перейдём к релаксации задач оптимизации с полиномиальными данными.

\subsection[Релаксации типа сумм квадратов\\ для полиномиальных задач]{Релаксации типа сумм квадратов\\ для полиномиальных задач}

Задача полиномиальной оптимизации характеризуется полиномиальной функцией цены и полиномиальными ограничениями. Введём следующий класс множеств.

\newpage

\begin{defin}
Множество $K \subset \mathbb R^n$ называется \emph{базовым полуалгебраическим} если оно задаётся
\begin{equation} \label{basicSA}
K = \{ x \,|\, f_i(x) = 0,\ g_j(x) \leq 0 \}
\end{equation}
для некоторых полиномов $f_i,g_j: \mathbb R^n \to \mathbb R$.
\end{defin}

Мы рассмотрим проблему минимизации полинома на базовом полуалгебраическом множестве $K$:
\begin{equation} \label{polynomial_problem}
\min_{x \in K} f_0(x).
\end{equation}

Для того чтобы переписать проблему в виде конической программы, введём конус $P_{d,K}$, состоящий из полиномов степени, не превосходящей $d$, которые неотрицательны на базовом полуалгебраическом множестве $K$. Этот конус конечномерный, замкнутый и выпуклый. Задача оптимизации запишется в виде
\begin{equation} \label{max_poly_equivalent} 
\max \tau:\quad f_0(x) - \tau \in P_{d,K},
\end{equation}
где $d$ --- не меньше степени $f_0$.

В общем случае конус $P_{d,K}$ не поддаётся эффективному описанию. Следуя философии предыдущего раздела, аппроксимируем этот конус полуопределённо представимым конусом $\Sigma_{d,K}$, состоящим из всех полиномов $p(x)$ степени, не превосходящей $d$, которые представимы в виде суммы
\begin{equation} \label{p_sos}
p(x) = \sigma_0(x) + \sum_i p_i(x)f_i(x) - \sum_j \sigma_j(x)g_j(x),
\end{equation}
где $p_i(x)$ --- произвольные полиномы, а $\sigma_0(x),\sigma_j(x)$ --- суммы квадратов полиномов. Здесь полиномы $f_i$ и $g_j$ определяют множество $K$ по формуле \eqref{basicSA}. Тогда любой полином из конуса $\Sigma_{d,K}$ неотрицательный на $K$, и, следовательно, $\Sigma_{d,K} \subset P_{d,K}$. Таким образом, конус $\Sigma_{d,K}$ является внутренней аппроксимацией $P_{d,K}$.

Полуопределённая представимость $\Sigma_{d,K}$ следует из того, что разложение \eqref{p_sos} приводит к условиям типа равенства, линейным по коэффициентам полинома $p$ и неизвестных полиномов $\sigma_0,p_i,\sigma_j$. Поэтому включение $p \in \Sigma_{d,K}$ эквивалентно конечному набору полуопределённых и линейных ограничений.

В итоге мы аппроксимируем исходную полиномиальную задачу оптимизации полуопределённой программой
\[ \max \tau:\quad f_0(x) - \tau \in \Sigma_{d,K}.
\]

\noindent Эта релаксация может быть усилена, если в основу определения конуса $\Sigma_{d,K}$ вместо условия \eqref{p_sos} положить условие
\[ p(x) \cdot \left( \sum_{i=1}^n x_i^2 \right)^r = \sigma_0(x) + \sum_i p_i(x)f_i(x) - \sum_j \sigma_j(x)g_j(x).
\]

\begin{example}
Рассмотрим полиномиальную проблему оптимизации
\begin{equation} \label{ex_pr}
\min x+y:\quad x \geq 0,\ x^2 + y^2 = 1.
\end{equation}

\noindent Множество $K = \{ (x,y) \,|\, x \geq 0,\ x^2 + y^2 = 1\}$ является базовым полуалгебраическим. Выберем максимальную степень $d = 3$. Аппроксимируем конус полиномов $P_{3,K}$ конусом $\Sigma_{3,K}$ полиномов, представимых в виде
\[ p(x,y) = \sigma_0(x,y) + l(x,y)(x^2+y^2-1) + \sigma_1(x,y)x,
\]

\noindent где $\sigma_0,\sigma_1$ --- суммы квадратов линейных полиномов, а $l$ --- линейный полином. Введём вектор мономов ${\bf x} = (x,y,1)^T$ степени, не превосходящей 1. Тогда имеем $p \in \Sigma_{3,K}$ тогда и только тогда, когда $p$ выражается в виде

\begin{align*}
p(x,y) =& {\bf x}^TA^0{\bf x} + {\bf l}^T{\bf x}\cdot (x^2+y^2-1) + ({\bf x}^TA^1{\bf x})\cdot x= \\
=& (A^1_{11}\!+\!l_x)x^3\!+\!(2A^1_{12}\!+\!l_y)x^2y\!+\!(A^0_{11}\!+\!2A^1_{13}\!+\!l_1)x^2\!+\!(A^1_{22}\!+\!l_x)xy^2+\\ &+ (2A^0_{12} + 2A^1_{23})xy + (2A^0_{13} + A^1_{33} - l_x)x + l_yy^3 + (A^0_{22} + l_1)y^2+\\ &+ (2A^0_{23} - l_y)y + A^0_{33} - l_1,
\end{align*}
где ${\bf l} = (l_x,l_y,l_1)^T \in \mathbb R^3$ и $A^0,A^1 \in \mathbb{S}_+^3$.

Полуопределённая релаксация задачи принимает вид
\[ \max_{A^0,A^1 \in \mathbb{S}_+^3} \tau:\quad A^1_{11} + l_x = 2A^1_{12} + l_y = A^0_{11} + 2A^1_{13} + l_1 = A^1_{22} + l_x = 0,
\]
\[ 2A^0_{12} + 2A^1_{23} = l_y = A^0_{22} + l_1 = 0,
\]
\[ 2A^0_{13} + A^1_{33} - l_x = 2A^0_{23} - l_y = 1,\ A^0_{33} - l_1 = -\tau.
\]

\noindent Используем линейные соотношения типа равенства для того, чтобы элиминировать часть переменных. Тогда получим эквивалентную полуопределённую программу
\[ \max\, -(A^0_{33} + A^0_{11} + 2A^1_{13}):\quad \begin{pmatrix} A^0_{11} &  A^0_{12} &  A^0_{13} \\  A^0_{12}  &  A^0_{11}+2A^1_{13} &  \frac12 \\ A^0_{13} &  \frac12 & A^0_{33} \end{pmatrix} \succeq 0,
\]
\[ \begin{pmatrix} A^1_{11} &  0  & A^1_{13} \\ 0 &  A^1_{11}  & -A^0_{12} \\  A^1_{13}  & -A^0_{12} &  1-A^1_{11}-2A^0_{13} \end{pmatrix} \succeq 0.
\]

\noindent Её решение приводит к оптимальному значению $-1$, которое является оптимальным также для исходной задачи.
\end{example}

\subsection{Моментные релаксации} \label{ch3_subs_moments}

Множество неотрицательных мер с носителем, содержащемся в некотором множестве $K \subset \mathbb R^n$, является выпуклым конусом. Если $K$ состоит из бесконечного числа точек, то этот конус бесконечномерный. Экстремальные меры этого конуса генерируются $\delta$-функциями $\mu(x) = \delta(x-\hat x)$, где $\hat x \in K$. Мера $\delta(x-\hat x)$ имеет носитель $\{ \hat x \}$ и интеграл от функции $f$ по этой мере равен
\[ \int_{\mathbb R^n} f(x)\delta(x-\hat x)\,dx = f(\hat x).
\]

Пусть $\mu$ --- неотрицательная мера на $\mathbb R^n$ с носителем $\Supp \mu$. Введём в пространстве $\mathbb R^n$ координаты $x_1, \, \dots, \, x_n$.

\begin{defin}
Пусть $\alpha = (\alpha_1, \, \dots, \, \alpha_n) \in \mathbb N^n$ --- мульти-индекс. \emph{Моментом} $m_{\alpha}$ меры $\mu$ будем называть значение интеграла
\[ m_{\alpha}(\mu) = \int_{\mathbb R^n} x_1^{\alpha_1}\cdot x_2^{\alpha_2}\cdot \dots \cdot x_n^{\alpha_n}\mu(x)\,dx = \int_{\mathbb R^n} x^{\alpha}\mu(x)\,dx.
\]
\end{defin}

Для данного $\alpha$ момент $m_{\alpha}$ является линейным функционалом на пространстве мер. Область определения этого функционала не совпадает со всем пространством, поскольку интеграл, задающий значение функционала, может расходиться. Для $\delta$-функции $\mu(x) = \delta(x - \hat x)$, однако, все моменты существуют и равны $m_{\alpha}(\mu) = \hat x^{\alpha}$.

Так как нам необходимо работать с конечномерными объектами, мы фиксируем натуральное число $d$ и рассматриваем только моменты $m_{\alpha}$, для которых степень $|\alpha| = \sum_{i=1}^n \alpha_i$ не превосходит $d$. Таких моментов конечное число, равное $N = \begin{pmatrix} n + d \\ d \end{pmatrix}$, и они образуют $N$-мерный \emph{вектор моментов} $m(\mu) = (m_{\alpha}(\mu))_{\alpha:|\alpha| \leq d}$.

\emph{Конус моментов} $M_d \subset \mathbb R^N$ определяется как множество всех векторов, которые представляются как вектор моментов некоторой неотрицательной меры $\mu$. Для подмножества $K \subset \mathbb R^n$ мы также определим конус $M_{d,K} \subset \mathbb R^N$, состоящий из векторов моментов неотрицательных мер $\mu$ с носителем в $K$. Таким образом, конусы моментов являются конечномерными проекциями бесконечномерного конуса неотрицательных мер.

\begin{exercise}
Описать конус моментов $M_d$ неотрицательных мер на $\mathbb R$. (Использовать лемму \ref{lem:blockHankel} с $n = 1$.)
\end{exercise}

Конусы моментов $M_{d,K}$ в общем случае не обладают алгоритмически эффективным описанием. Мы рассмотрим необходимые условия для того, чтобы данный вектор являлся вектором моментов некоторой неотрицательной меры. Множество векторов, удовлетворяющих этим необходимым условиям, образует \emph{внешнюю} аппроксимацию конуса моментов.

Пусть ${\bf x} = \left(1, \, x_1, \, \dots, \, x_n, \, x_1^2, \, x_1x_2, \, \dots, \, x_n^{[d/2]}\right)^T$ --- вектор мономов $x^{\alpha}$, для которых степень $|\alpha|$ не превосходит целую часть от $d/2$. Тогда все элементы одноранговой матрицы ${\bf x}{\bf x}^T$ являются мономами со степенью, не превосходящей $d$. Рассмотрим матрично-значный интеграл
\[ \int_{\mathbb R^n} {\bf x}{\bf x}^T\mu(x)\,dx,
\]
где $\mu$ --- неотрицательная мера. Этот интеграл является неотрицательно определённой матрицей, элементы которой равны неким элементам вектора моментов $m(\mu)$. Отсюда вытекает полуопределённое ограничение на вектор моментов, а именно, что матрица, составленная соответствующим образом, является неотрицательно определённой.

Пусть теперь $K = \{ x \in \mathbb R^n \,|\, f_i(x) = 0,\ g_j(x) \leq 0 \}$ --- базовое полуалгебраическое множество, и пусть $\mu$ --- неотрицательная мера с носителем в~$K$.

Обозначим степень полинома $f_i$ через $d_i$. Тогда для любого полинома $p$ со степенью, не превосходящей $d - d_i$, имеем
\[ \int_K p(x)f_i(x)\mu(x)\,dx = 0.
\]
С другой стороны, данный интеграл выражается через линейную комбинацию элементов вектора моментов $m(\mu)$. Таким образом, мы получаем линейное ограничение типа равенства на вектор моментов $m(\mu)$. Если полином $p(x)$ пробегает все мономы $x^{\beta}$ степени $|\beta| \leq d - d_i$, то мы получим максимальный линейно независимый набор таких линейных ограничений.

Образуем вектор ${\bf x}'$ всех мономов степени, не превосходящей целую часть от $\frac{d-d_j}{2}$, и рассмотрим матрично-значный интеграл
\[ -\int_K {\bf x}'({\bf x}')^Tg_j(x)\mu(x)\,dx.
\]
Этот интеграл является неотрицательно определённой матрицей, каждый элемент которой представляется в виде линейной комбинации неких элементов вектора моментов $m(\mu)$. Отсюда получаем полуопределённое ограничение на $m(\mu)$.

\bigskip

Рассмотрим снова проблему полиномиальной оптимизации \eqref{polynomial_problem}, где $f_0 = \sum_{\alpha} c_{\alpha}x^{\alpha}$ --- полином степени, не превосходящей $d$, а $K$ --- базовое полуалгебраическое множество. Перепишем эту проблему в виде задачи минимизации
\begin{equation} \label{moment_min_equivalent} 
\min_{\mu \geq 0: \Supp\mu \subset K} \int_{\mathbb R^n} f_0(x)\mu(x)\,dx:\quad \int_{\mathbb R^n} \mu(x)\,dx = 1
\end{equation}
на множестве всех вероятностных мер с носителем в $K$.

Условие типа равенства на $\mu$ можно записать в виде $m_0(\mu) = 1$, а интеграл, задающий функцию цены, в виде линейной комбинации $\sum_{\alpha} c_{\alpha}m_{\alpha}(\mu)$ элементов вектора моментов $m(\mu)$. Таким образом, проблема предстанет в виде
\[ \min_{m \in M_{d,K}} \sum_{\alpha} c_{\alpha}m_{\alpha}:\quad m_0 = 1.
\]

Заменив не поддающееся эффективному описанию включение $m \in M_{d,K}$ набором полуопределённых и линейных необходимых условий, построенных выше, мы получаем полуопределённую релаксацию исходной проблемы.

\begin{example}
Рассмотрим снова проблему \eqref{ex_pr}. Положим $d = 3$, тогда вектор моментов будет 10-мерным. Описанная выше полуопределённая релаксация запишется в виде
\[ \min m_{10} + m_{01}:\quad \begin{pmatrix} m_{00} & m_{10} & m_{01} \\  m_{10} & m_{20} & m_{11} \\  m_{01} & m_{11} & m_{02} \end{pmatrix} \succeq 0,
\]
\[ m_{20}+m_{02}-m_{00} = m_{30}+m_{12}-m_{10} = m_{21}+m_{03}-m_{01} = 0,
\]
\[ \begin{pmatrix} m_{10} & m_{20} & m_{11} \\  m_{20} & m_{30} & m_{21} \\  m_{11} & m_{21} & m_{12} \end{pmatrix} \succeq 0,\ m_{00} = 1.
\]
Её решение также даёт оптимальное значение $-1$.
\end{example}

\subsection{Релаксации типа сумм квадратов для верхних оценок}

Рассмотрим снова проблему полиномиальной оптимизации \eqref{polynomial_problem}. Релаксации, основанные на замене конуса положительных полиномов $P_{d,K}$ в эквивалентной формулировке \eqref{max_poly_equivalent} на конус сумм квадратов $\Sigma_{d,K}$, сужают множество допустимых точек, и поэтому дают нижнюю оценку на оптимальное значение этих задач. Моментные релаксации, рассмотренные в п.~\ref{ch3_subs_moments}, также приводят к нижним оценкам на оптимальное значение, поскольку расширяется допустимое множество в эквивалентной формулировке \eqref{moment_min_equivalent}.

В этом разделе мы рассмотрим релаксации, приводящие к \emph{верхним} оценкам на оптимальное значение для задачи \eqref{polynomial_problem}. Идея состоит в том, чтобы сузить допустимое множество в формулировке \eqref{moment_min_equivalent}, рассматривая только положительные меры на базовом полуалгебраическом множестве $K$, имеющие плотность, выражающуюся в виде суммы квадратов. Этот подход впервые был предложен в работе \cite{Lasserre11}. 

Пусть $\mu$~--- \emph{фиксированная} положительная мера с носителем $K$. Пусть далее $h \in \Sigma_{d,K}$~--- полином, представимый как сумма вида \eqref{p_sos}, и поэтому неотрицательный на $K$. Тогда мера $\mu_h(x) = \mu(x) \cdot h(x)$ также неотрицательна на $K$. Моменты этой меры являются линейными комбинациями моментов меры $\mu$, с коэффициентами, равными коэффициентам полинома $h$. А именно, если полином $h$ представляется в виде $h(x) = \sum_{|\alpha| \leq d} c_{\alpha}x^{\alpha}$, то моменты меры $\mu_h$ задаются по формуле
\[ m_{\beta}(\mu_h) = \int_K x^{\beta}\mu_h(x)\,dx = \int_K \sum_{|\alpha| \leq d} c_{\alpha}x^{\alpha+\beta}\mu(x)\,dx = \sum_{|\alpha| \leq d} c_{\alpha}m_{\alpha+\beta}(\mu).
\]
Коэффициенты полинома $h \in \Sigma_{d,K}$, в свою очередь, задаются полуопределёнными условиями.

Релаксация исходной полиномиальной проблемы оптимизации \eqref{polynomial_problem} запишется в виде полуопределённой программы:
\[ \min_{h \in \Sigma_{d,K}} \int_K h(x) \cdot f_0(x) \cdot \mu(x)\,dx: \quad \int_K h(x) \cdot \mu(x)\,dx = 1.
\]
Переменными в этой программе являются коэффициенты полинома $h \in \Sigma_{d,K}$, которые линейно входят и в функцию цены, и в ограничение программы. 

Точность релаксации возрастает с максимальной степенью $d$ полинома~$h$. Для компактных множеств $K$ иерархия этих релаксаций асимптотически точна, т.е., получаемые верхние оценки на оптимальное значение исходной задачи стремятся к этому значению при $d \to \infty$. Однако, так как носитель меры $\mu_h$ всегда будет совпадать с $K$, релаксация любой степени может быть точной только в том случае, когда минимизируемая функция $f_0$ постоянна на~$K$. 

Описанный подход предполагает наличие положительной меры $\mu$ с носителем $K$ такой, что моменты этой меры легко вычислимы или задаются аналитическим выражением. Поэтому данный подход применим только для относительно простых множеств $K$. 

\subsection{Тригонометрические полиномы} \label{ch3_subs_trig_moments}

В этом разделе мы рассмотрим тригонометрические полиномы, скалярные или матрично-значные, от одной переменной. В роли этой переменной выступает угол $\phi \in [-\pi,\pi]$, и рассматриваемые функции $2\pi$-периодически продолжаются на всю вещественную ось. Эквивалентный подход состоит в том, чтобы рассматривать величину $z = e^{i\phi}$, принимающую значения на единичном круге $\mathbb T \subset \mathbb C$, как независимую переменную. В этом случае тригонометрический полином соответствует конечному ряду Лорана по переменной $z$. Тригонометрические полиномы сводятся к обычным заменой независимой переменной $\cos\phi = \frac{z + \bar z}{2} = \frac{1-t^2}{1+t^2}$, $\sin\phi = \frac{z - \bar z}{2i} = \frac{2t}{1+t^2}$, $t \in \mathbb R \cup \{\infty\}$.

Тригонометрические полиномы широко используются в обработке сигналов, управлении и идентификации линейных динамических систем. Рациональные по переменной $z$ функции, скалярные или матрично-значные, выступают в качестве передаточных функций линейных динамических систем в дискретном времени с конечной памятью. Неотрицательные меры на~$\mathbb T$, скалярные или матрично-значные, играют роль спектров стационарных сигналов. Здесь оператор умножения на $z$ соответствует сдвигу по времени на одну единицу.

Тригонометрический полином степени $d$ имеет общий вид $p(\phi)\!= \linebreak=\sum_{k = -d}^d A_ke^{ik\phi} = \sum_{k = -d}^d A_kz^k$, где $A_{-k} = A_k^*$ для всех $k = 0, \, \dots, \, d$. Здесь $A_k$ --- матрицы размера $n \times n$, а $A^*$ --- эрмитово сопряжённая (комплексно сопряжённая транспонированная) к $A$ матрица. Полином принимает значения в пространстве $\mathbb{H}^n$ эрмитовых матриц размера $n$. В случае $n = 1$ мы имеем дело со скалярными полиномами, принимающими значения в $\mathbb R$. Множество неотрицательных полиномов степени, не превышающей $d$, является выпуклым конусом. В случае матрично-значных полиномов под неотрицательностью подразумевается, что значения полинома являются неотрицательно определёнными эрмитовыми матрицами.

Как и в случае обычных полиномов от одной переменной, конус неотрицательных полиномов имеет полуопределённое представление.

\begin{teo}
Пусть $p(\phi) = \sum_{k = -m}^m A_ke^{ik\phi}$ --- тригонометрический полином. Полином $p$ является неотрицательным тогда и только тогда, когда существует эрмитовая неотрицательно определённая матрица ${\bf A} \in \mathbb{H}_+^{n(m+1)}$ такая, что $A_k = \sum_{i-j = k} A_{ij}$, где $A_{ij}$, $i,j = 0, \, \dots, \, m$ --- блоки матрицы ${\bf A}$ размера $n$.
\end{teo}

Доказательство следует той же схеме, что и в случае обычных неотрицательных матрично-значных полиномов от одной переменной.

\begin{lemma}
Пусть $T \in \mathbb{H}_+^{n(m+1)}$ --- неотрицательно определённая блочно-тёплицевая эрмитова матрица с блоками размера $n$. Тогда $T$ представляется в виде суммы неотрицательно определённых блочно-тёплицевых матриц ранга 1.
\end{lemma}

%
%

\begin{lemma}
Пусть $p(\phi) = \sum_{j = -m}^m A_je^{ij\phi}$ --- положительный тригонометрический матричный полином, а $T \in \mathbb{H}_+^{n(m+1)}$ --- блочно-тёплицевая матрица с блоками $T_{-m}, \, \dots, \, T_m \, \in \,  \mathbb{H}^n$. Тогда $\sum_{j=-m}^{m} \langle A_j,T_j \rangle \geq 0$. В частности, существует такая матрица ${\bf A} \in \mathbb{H}_+^{n(m+1)}$, состоящая из $n \times n$ блоков $A_{ij}$, $i,j = 0, \, \dots, \, m$, что $A_k = \sum_{i-j=k} A_{ij}$.
\end{lemma}

%

Факторизуя матрицу ${\bf A} = {\bf B}{\bf B}^*$ и разбивая фактор ${\bf B} \in \mathbb C^{n(m+1) \times l}$ на блоки $B_0, \, \dots, \, B_m \, \in \, \mathbb C^{n \times l}$, мы получаем, что положительный тригонометрический полином $p$ представляется в виде квадрата $p(\phi) = \left(\sum_{k=0}^m B_ke^{ik\phi}\right)\left(\sum_{k=0}^m B_ke^{ik\phi}\right)^*$.


\section{Оптимизация топологии фермы}

\subsection[Оптимизация топологии фермы: подход на базе\\ полуопределённого программирования]{Оптимизация топологии фермы: подход на базе\\ полуопределённого программирования} \label{subs:ch3_truss}

\emph{Ферма} --- это стержневая система в строительной механике, остающаяся геометрически неизменяемой после замены её жёстких узлов шарнирными. При нагрузке узлов в элементах фермы возникают только усилия растяжения--сжатия. В этом разделе мы рассмотрим, как проектировать фермы, оптимальные по отношению к заданным сценариям нагрузок, с помощью полуопределённого программирования. Этот подход был разработан А.\,С.~Немировским и А.~Бен-Талем в работе \cite{NemirovskiTruss}.

Ферму можно моделировать множеством узлов, соединённых стержнями. При нагрузке, т.е. приложении силы к некоторому подмножеству узлов, конструкция деформируется. Эта деформация сопровождается растяжением или сжатием стержней, которые при этом запасают потенциальную энергию и вызывают силу реакции. Ферма приходит в новое состояние равновесия, в котором силы реакции компенсируют приложенную нагрузку. Задача оптимизации состоит в том, чтобы минимизировать суммарную потенциальную энергию, соответствующую состоянию равновесия, или, что эквивалентно, максимизировать жёсткость конструкции. Искомыми переменными в этой задаче выступают массы стержней между узлами. При этом масса данного стержня может быть равной нулю, что на практике соответствует отсутствию этого стержня. Поэтому, по сути, решается задача оптимизации топологии фермы, т.е. расположения узлов и стержней, в ней участвующих. Суммарная масса $M$ стержней при этом ограничена.

Чтобы формализовать задачу, нам нужно знать, какую силу $f_{ik}$ вызывает смещение $x_i,x_j \in \mathbb R^3$ узлов $i,j$, связанных стержнем $k$ массы $m_k$, на узел $i$. Здесь и далее через $k$ обозначаем индекс стержня, а под суммой по индексу $k$ подразумеваем сумму по всем стержням. Пусть $v_i,v_j \in \mathbb R^3$ --- векторs позиций узлов $i,j$. Так как стержень реагирует только на растяжение или сжатие, сила будет пропорциональна скалярному произведению $\left\langle \frac{v_i-v_j}{\|v_i-v_j\|^2},x_i-x_j \right\rangle$ (здесь и далее $\|~\cdot~\| = \|~\cdot~\|_2$), которое задаёт относительное изменение его длины. Далее, сила, с которой стержень воздействует на узел $i$, будет \nk{направлена по направлению} $\frac{v_j-v_i}{\|v_i-v_j\|}$ в случае его растяжения и по направлению $\frac{v_i-v_j}{\|v_i-v_j\|}$ в случае сжатия. Кроме того, она будет пропорциональна сечению, т.е. массе стержня на единицу его длины $\frac{m_k}{\|v_i - v_j\|}$. В конце концов, сумма сил, действующих на узлы $i$ и $j$, равна нулю. Отсюда получаем $f_{ik} = -f_{jk} = -c_km_k\left\langle \frac{v_i-v_j}{\|v_i-v_j\|^2},x_i-x_j \right\rangle\frac{v_i-v_j}{\|v_i-v_j\|^2}$, где $c_k > 0$ --- некоторая независящая от смещений $x$ константа, задаваемая модулем Юнга стержня.

Соберем все смещения $x_i$ и силы $f_{ik}$ в большие векторы $x,f_k \in \mathbb R^{3n}$, где $n$ --- количество узлов. Тогда сила $f_k$, генерируемая стержнем $k$, записывается в виде $f_k = -m_kA_kx$, где $A_k$ --- симметричная неотрицательно определённая одноранговая матрица $a_ka_k^T$. В векторе $a_k \in \mathbb R^{3n}$ все компоненты нулевые, кроме $a_{ik},a_{jk}$, соответствующих узлам $i,j$ и задающихся $a_{ik} = -a_{jk} = \sqrt{c_k}\frac{v_i-v_j}{\|v_i-v_j\|^2}$. Если $f$ --- сила, вызываемая нагрузкой на ферму, то уравнения равновесия $f + \sum_k f_k = 0$ приводят к соотношению $\sum_k m_kA_kx = Ax = f$, где $A = \sum_k m_ka_ka_k^T$.

Пусть $V \subset \mathbb R^{3n}$ --- линейное подпространство возможных смещений узлов фермы. Оно задаётся условиями, накладываемыми точками опоры. Мы предположим, что на подпространстве $V$ матрица $A$ положительно определена. Это условие эквивалентно тому, что фиксация фермы на опорах не даёт ей двигаться как твёрдое тело. В этом случае деформация $x$ однозначно вычисляется из нагрузки $f$, а энергия деформации задаётся выражением
\[ c_f(m) = \langle f,x \rangle = \sup_{u \in V} \ \left( 2\langle f,u \rangle - u^TAu \right).
\]

Мы будем предполагать, что нагрузка $f$ может варьировать в некотором эллипсоиде
\[ E = \{ Qe \mid \|e\| \leq 1 \},
\]
где $Q \in \mathbb R^{3n \times q}$ --- некоторая матрица, задающая форму эллипсоида, а $q$ --- размерность эллипсоида. Данный подход позволяет одновременно учесть конечное число сценариев, т.е. заданных наперёд нагрузок, и получить робастное решение, жёсткое по отношению к произвольным малым нагрузкам, прилагаемым к активным (т.е. участвующим в конструкции) узлам. Выбор эллипсоида в качестве множества возможных нагрузок приводит к некоторой консервативности решения. Оно обусловлено тем, что менее консервативные или точные формулировки задачи не могут быть решены эффективно.

Таким образом, мы решаем задачу
\begin{equation} \label{truss_problem}
\inf_{m \in \mathbb R_+^d}\,\sup_{u \in V,\|e\| \leq 1}\, \left(2u^TQe - u^TAu\right): \quad A = \sum_{k = 1}^d m_kA_k,\ \sum_{k=1}^d m_k = M,
\end{equation}
где $d$ --- количество потенциальных стержней, $m = (m_1,\dots,m_d)^T$ --- вектор их масс, а $M$ --- их суммарная масса.

Перепишем задачу в виде полуопределённой программы. Определим две квадратичные функции ${\cal Q}_1(e,u) = \|e\|^2$, ${\cal Q}_2(e,u) = 2u^TQe - u^TAu$. Ясно, что супремум в \eqref{truss_problem} не больше $\tau$ тогда и только тогда, когда квадратичная функция $\tau {\cal Q}_1 - {\cal Q}_2$ неотрицательно определена, т.е.
\[ \begin{pmatrix} \tau I & -Q^T \\ -Q & A \end{pmatrix} \succeq 0.
\]
В итоге получаем программу
\begin{equation} \label{primal_truss_problem} 
\min_{\tau,m}\tau: \quad m_k \geq 0,\ \sum_{k=1}^d m_k = M,\ \begin{pmatrix} \tau I & -Q^T \\ -Q & \sum_{k = 1}^d m_ka_ka_k^T \end{pmatrix} \succeq 0.
\end{equation}
Конус, лежащий в основе этой программы, задаётся прямым произведением $\mathbb R_+^d \times \mathbb{S}_+^{3n+q}$. Его размерность равна $d + \frac{(3n+q)(3n+q+1)}{2}$, а количество переменных равно $d$ ($d$ переменных $m_k$, связанных одним линейным условием типа равенства, и переменная $\tau$).

Оценим сложность решения программы прямым методом внутренней точки в предположении $q \ll n$, $d \approx \frac{n^2}{2}$. Вследствие разреженности векторов $a_k$ в цене шага Ньютона доминирует стоимость решения линейной системы уравнений на прирост переменных задачи. Это порядка куба числа переменных, т.е. $O(n^6)$. Число итераций пропорционально корню из параметра барьера, который задаётся суммой количества линейных ограничений типа неравенства и размера матрицы, т.е. $d + 3n + q = O(n^2)$. Таким образом, сложность метода пропорциональна $n^7$.

Рассмотрим теперь двойственную программу. Она также определена над конусом $\mathbb R_+^d \times \mathbb{S}_+^{3n+q}$, поскольку этот конус самодвойственный. Число переменных в двойственной задаче равно размерности конуса за вычетом количества переменных в исходной задаче, т.е. $\frac{(3n+q)(3n+q+1)}{2}$. Нетрудно проверить, что двойственная задача имеет вид
\[ \max_{\Lambda,X,Y,\rho}(2\langle Q,X \rangle-M\rho):\quad \begin{pmatrix} \Lambda & X^T \\ X & Y \end{pmatrix} \succeq 0,\ \mu_k = \rho - \langle a_ka_k^T,Y \rangle \geq 0,\ \tr\,\Lambda = 1.
\]
Элемент двойственного конуса задаётся вектором $\mu = (\mu_1, \, \dots, \, \mu_d)^T \, \in \, \mathbb R_+^d$ и матрицей $\begin{pmatrix} \Lambda & X^T \\ X & Y \end{pmatrix}$. 

Пусть $\left(\hat\Lambda,\hat X,\hat Y,\hat\rho\right)$ --- решение задачи. Тогда и $\left(\hat\Lambda,\hat X,\hat X\hat\Lambda^{\dagger}\hat X^T,\hat\rho\right)$ --- решение, поскольку $\hat Y \succeq \hat X\hat\Lambda^{\dagger}\hat X^T$ в силу матричного неравенства (см. п.~\ref{ch1_section_schur}) и вследствие этого $\langle a_ka_k^T,\hat Y \rangle \geq \langle a_ka_k^T,\hat X\hat\Lambda^{\dagger}\hat X^T \rangle$. Здесь $\hat\Lambda^{\dagger}$ --- псевдо-обратная матрицы $\hat\Lambda$. Поэтому можно положить $Y = X\Lambda^{\dagger}X^T$, и скалярные неравенства примут вид 
\[ \rho - \langle a_ka_k^T,X\Lambda^{\dagger}X^T \rangle = \rho - a_k^TX\Lambda^{\dagger}X^Ta_k \geq 0.
\]
Дополнение по Шуру позволяет записать эти неравенства в виде $\begin{pmatrix} \Lambda & X^Ta_k \\ a_k^TX & \rho \end{pmatrix} \succeq 0$. Таким образом мы избавляемся от переменной $Y$, и программа принимает вид
\begin{equation} \label{dual_truss_problem}
\max_{\Lambda,X,\rho}\,(2\langle Q,X \rangle-M\rho):\ \begin{pmatrix} \Lambda & X^Ta_k \\ a_k^TX & \rho \end{pmatrix} \succeq 0,\ k = 1, \, \dots, \, d,\ \tr\,\Lambda = 1.
\end{equation}
Количество переменных в этой формулировке равно $\frac{q(q+1)}{2} + 3qn = O(qn)$, а конус, над которым определена программа, задаётся прямым произведением $d$ копий матричного конуса $\mathbb{S}_+^{q+1}$. Параметр барьера этого конуса равен $d(q+1) = O(n^2q)$. Отсюда получаем сложность $O\left((qn)^3(n^2q)^{1/2}\right) = O(q^{7/2}n^4)$ решения двойственной задачи в формулировке \eqref{dual_truss_problem}.

Эта сложность значительно ниже сложности $O(n^7)$ решения исходной программы \eqref{primal_truss_problem}. Дело в том, что при переформулировке двойственной задачи повышение в сложности конуса (т.е., параметра барьера) более чем компенсировалось уменьшением количества переменных. Перейдя от \eqref{dual_truss_problem} снова к двойственной задаче, получаем программу
\begin{equation} \label{truss_second_primal}
\min_{P_k}\,\tau: \ P_k = \begin{pmatrix} L_k & -d_k \\ -d_k^T & m_k \end{pmatrix} \succeq 0,\ \sum_{k=1}^d L_k = \tau I,\ \sum_{k=1}^d a_kd_k^T = Q,\ \sum_{k=1}^d m_k = M.
\end{equation}
Нетрудно видеть, что переменные $m_k$ в этой формулировке соответствуют искомым массам стержней. Отслеживая двойственные переменные в ходе решения программы \eqref{dual_truss_problem} методом следования центральному пути, можно получить сколь угодно точные оценки этих масс.

То, что задача \eqref{truss_second_primal} эквивалентна задаче \eqref{primal_truss_problem}, вытекает также из следующей леммы.

\begin{lemma} \label{lem:matrix_decomp}
Пусть $P = \begin{pmatrix} A & B^T \\ B & \sum_{k=1}^s C_k \end{pmatrix}$ --- неотрицательно определённая симметричная матрица, при этом слагаемые $C_k$, $k = 1,\dots,s$~--- также симметричные и неотрицательно определённые. Тогда существует разбиение матрицы $P$ на неотрицательно определённые слагаемые $P_k = \begin{pmatrix} A_k & B_k^T \\ B_k & C_k \end{pmatrix}$.
\end{lemma}

\begin{exercise}
Доказать лемму \ref{lem:matrix_decomp}.
\end{exercise}

Из леммы следует, что матрицу в задаче \eqref{primal_truss_problem} можно разложить на неотрицательно определённые слагаемые $\begin{pmatrix} L_k & -Q_k^T \\ -Q_k & m_ka_ka_k^T \end{pmatrix}$. Так как подматрица в правом нижнем углу имеет ранг 1 с образом, генерируемым вектором $a_k$, то блок $Q_k$ должен иметь вид $a_kd_k^T$ для некоторого вектора $d_k$. Но неотрицательная определённость матрицы $\begin{pmatrix} L_k & -d_ka_k^T \\ -a_kd_k^T & m_ka_ka_k^T \end{pmatrix}$ эквивалентна неотрицательной определённости матрицы $\begin{pmatrix} L_k & -d_k \\ -d_k^T & m_k \end{pmatrix}$. Отсюда получаем формулировку \eqref{truss_second_primal}.

\medskip

На примере этого раздела видно, что сложность полуопределённой задачи иногда можно понизить преобразованиями матриц или переходом к двойственной задаче. Некоторые солверы в зависимости от количества переменных самостоятельно принимают такое решение.

\subsection[Оптимизация топологии фермы: подход на базе\\ прямо-двойственного субградиентного метода]{Оптимизация топологии фермы: подход на базе\\ прямо-двойственного субградиентного метода}\label{truss}

Задача оптимизации топологии фермы формулируется таким образом, что позиции потенциальных узлов и стержней в пространстве заданы наперёд. Чтобы получить хорошее приближение оптимального расположения узлов и стержней, нужно покрыть область достаточно густой сеткой потенциальных позиций этих объектов. В типичной ситуации решение задачи таково, что большинство стержней имеют нулевую массу, а большинство узлов не соединено с другими узлами стержнями положительной массы, т.е. эти узлы и стержни попросту отсутствуют в соответствующей оптимальной физической конструкции. Поэтому, несмотря на то, что в конечной конструкции количество узлов и стержней обозримо, число потенциальных объектов, и тем самым и размерность задачи, могут быть настолько велики, что решение задачи с помощью конического программирования (см. п.~\ref{subs:ch3_truss}) представляется невозможным. 

В этом разделе мы опишем альтернативный подход к задаче оптимизации топологии фермы на базе алгоритма \ref{Switch_Subgrad} \cite{NesterovShpirko14}. Отметим, что метод обобщается на линейные программы общего вида с однородно разреженной матрицей коэффициентов в линейных ограничениях.

Как сказано выше в п.~\ref{subs:ch3_truss}, задача состоит в нахождении наиболее жёсткой механической структуры, способной выдержать внешнюю нагрузку, но при этом связанной с ограниченными затратами материала на саму конструкцию. Математическую формулировку возможно выписать так (см.~также п.~\ref{subs:ch3_truss}):
\begin{equation}\label{p5_4.1}
\min_{m,x}\,\langle f,x\rangle:\quad A(m)x=f,\, m \geq 0,\,\langle {\bf  e},m\rangle=M,
\end{equation}
где $f\in \mathbb R^{3n}$~--- вектор внешних сил, $x\in \mathbb R^{3n}$~--- вектор виртуального смещения $n$ узлов в $\mathbb R^{3}$, $m\in \mathbb R^{d}_{+}$~--- вектор весов $d$ балок и $M$~--- общий вес конструкции, ${\bf e} = (1, 1, ..., 1)$. Матрица соответствия $A(m)$ имеет следующий вид:
$$A(m)=\sum_{i=1}^{d}m_{i}a_{i}a_{i}^{T},$$
где $a_{i}\in \mathbb R^{3n}$ --- вектор, описывающий взаимодействие двух узлов, соединённых балкой с индексом $i$. 

Существенно, что такие векторы можно считать разреженными: на каждый вектор $a_i$ приходится не более 6 ненулевых элементов. Мы будем также считать что вектор внешней силы $f$ разреженный, что эквивалентно тому, что сила прилагается только к одному или нескольким узлам. Кроме того, мы сделаем ещё одно важное допущение. А именно, мы полагаем, что в каждом узле встречаются не больше, чем $s = O(1)$ стержней. Это условие выполнено, например, если каждый узел соединён только с соседними узлами.

Для простоты предположим, что все смещения $x \in \mathbb R^{3n}$ допустимы, и для $m > 0$ матрица $A$ положительно определена. Эти условия выполнены, если некоторое подмножество узлов жёстко зафиксировано, конструкция не имеет степеней свободы движения, а $x$ содержит только позиции оставшихся узлов.

Опишем подход \cite{NesterovShpirko14} к указанной задаче~\eqref{p5_4.1}. Обозначим множество допустимых векторов масс через $\triangle(M)=\left\{\,m\geqslant0,\,\langle {\bf e},m\rangle=M\right\}$. Оказывается, что проблему можно переписать в виде задачи линейного программирования
\begin{equation}\label{p5_4.2}
\begin{split}
&\min_{x,m}\left\{\langle f,x\rangle:\,A(m)x=f,\,m\in \triangle(M)\right\}=\\
&=\min_{m\in \triangle(M)}\max\limits_{x}\left\{2\langle f,x\rangle-\langle A(m)x,x\rangle\right\} = \max\limits_{x}\min_{m\in \triangle(M)}\left\{2\langle f,x\rangle-\langle A(m)x,x\rangle\right\}=\\
&=\max\limits_{x}\left\{2\langle f,x\rangle- M\max\limits_{1\leqslant i\leqslant d}\langle a_{i},x\rangle^{2}\right\}=\max\limits_{\alpha,y}\left\{2\alpha\langle f,y\rangle- \alpha^{2}M\max\limits_{1\leqslant i\leqslant d}\langle a_{i},y\rangle^{2}\right\}=\\
&=\max\limits_{y}\frac{\langle f,y\rangle^{2}}{M\max\limits_{1\leqslant i\leqslant d}\langle a_{i},y\rangle^{2}}=\frac{1}{M}\left(\max\limits_{y}\left\{\langle f,y\rangle:\,\max\limits_{1\leqslant i\leqslant d}|\langle a_{i},y\rangle|\leqslant1\right\}\right)^{2}.
\end{split}
\end{equation}
Равенство во второй строке следует из теоремы Сиона--Какутани \ref{minmax_theorem} (см. п.~\ref{subs:kakutani}), поскольку функция $2\langle f,x\rangle-\langle A(m)x,x\rangle$ вогнута по $x$, выпукла (даже линейна) по $m$, и множество $\triangle(M)$ компактно. Обозначим через $y^{*}$ решение задачи максимизации
\begin{equation}\label{p5_4.7}
\max\limits_{y }\,\langle f,y\rangle:\quad g(y) := \max\limits_{1\leqslant i\leqslant d}\left( \pm\langle a_{i},y\rangle - 1 \right) = \max_{1 \leq i \leq 2d} g_i(y) \leqslant 0,
\end{equation}
где мы определили
\[ g_i(y) = \langle a_{i},y\rangle - 1, \quad g_{i+d}(y) = -\langle a_{i},y\rangle - 1,\qquad i = 1,\dots,d.
\]

Из решения задачи \eqref{p5_4.7} невозможно напрямую восстановить искомые массы $m_i$ стержней. Поэтому рассмотрим двойственную к этой задаче, которая имеет вид
\begin{equation} \label{truss_dual} 
\min_{\mu \in \mathbb R_+^{2d}} \max_y \left( \langle f,y \rangle - \sum_{i=1}^{2d} \mu_i g_i(y) \right).
\end{equation}
Из оптимального решения $\mu^*$ этой задачи получаем множители Лагранжа $\lambda^{*}\in \mathbb R^{d}_{+}$ по формуле $\lambda_i^* = \mu_i^* + \mu_{i+d}^*$, удовлетворяющие
\begin{equation}\label{p5_4.3}
f=\sum_{i\in J_{+}}a_{i}\lambda_{i}^{*}-\sum_{i\in J_{-}}a_{i}\lambda_{i}^{*}\quad \text{и } \lambda_{i}^{*}=0\text{ при }i\not\in J_{+}\bigcup J_{-},
\end{equation}
где $J_{+}=\{i:\,\langle a_{i},y^{*}\rangle=1\}$ и $J_{-}=\{i:\,\langle a_{i},y^{*}\rangle=-1\}$. Умножив первое уравнение в~\eqref{p5_4.3} на $y^{*}$, получим
\[
\langle f,y^{*}\rangle = \langle {\bf e},\lambda^{*}\rangle.
\]
Заметим, что первое уравнение в~\eqref{p5_4.3} может быть записано как
\begin{equation}\label{p5_4.5}
f=A(\lambda^{*})y^{*}.
\end{equation}
Покажем, как можно по решению прямо-двойственной пары задач \eqref{p5_4.7},\eqref{truss_dual} восстановить решение исходной задачи \eqref{p5_4.1}. Обозначим
\begin{equation}\label{p5_4.6}
m^{*}=\frac{M}{\langle {\bf e},\lambda^{*}\rangle}\cdot \lambda^{*},\qquad x^{*}=\frac{\langle {\bf e},\lambda^{*}\rangle}{M}\cdot y^{*}.
\end{equation}
Тогда с учётом~\eqref{p5_4.5} имеем $f=A(m^{*})x^{*}$ и $m^{*}\in \triangle(M)$. Таким образом, пара в~\eqref{p5_4.6} допустима для задачи \eqref{p5_4.1}. С другой стороны, 
$$\langle f,x^{*}\rangle=\left\langle f,\frac{\langle {\bf e},\lambda^{*}\rangle}{M}\cdot y^{*}\right\rangle =\frac{1}{M}\cdot\left\langle {\bf e},\lambda^{*}\right\rangle\cdot\langle f,y^{*}\rangle=\frac{1}{M}\cdot\langle f,y^{*}\rangle^{2}.$$
Поэтому пара векторов $(m^{*}, x^{*})$, определяемая~\eqref{p5_4.6}, есть решение задачи \eqref{p5_4.1}.


Итак, приведённые выше рассуждения позволяют заменить исходную задачу \eqref{p5_4.1} на прямо-двойственную пару задач \eqref{p5_4.7},\eqref{truss_dual}. 
Эти задачи можно решать алгоритмом \ref{Switch_Subgrad}. Заметим, что в задаче \eqref{p5_4.7} мы ищем максимум, а не минимум, как в \eqref{equv_state_problem}, и соответственно \eqref{truss_dual} является задачей на минимум, в то время как двойственная задача \eqref{f_8} есть задача на максимум по двойственной переменной.

Обратим внимание на то, что функция цены и ограничение в задаче \eqref{p5_4.7} удовлетворяют условиям теоремы \ref{swithc_subgrad_thm}. А именно, 2-норма градиента $f$ функции цены ограничена константой $M_f = \|f\|_2$, а 2-норма субградиента ограничения $g$ ограничена константой $M_g = \max_{1 \leq i \leq d}\|a_i\|_2$. Последнее утверждение следует из того, что в качестве субградиента функции $g$ в точке $y$ можно выбрать вектор $\pm a_i$, в зависимости от того, на каком индексе $i$ достигается максимум в определении $g$.

Поэтому после выполнения необходимого количества итераций алгоритма \ref{Switch_Subgrad} по теоремам \ref{swithc_subgrad_thm}, \ref{thm_subgr_primal-dual} можно гарантировать достижение некоторого $\varepsilon$-приближённого решения задачи \eqref{p5_4.7}, а также сокращения зазора двойственности между достигнутыми значениями задач \eqref{p5_4.7},\eqref{truss_dual} до $\varepsilon$, если допустить возможность локализации значений $y$ в пределах некоторого компактного подмножества $Q \subset \mathbb{R}^{3n}$. Отсюда получаем, что  решение задачи \eqref{truss_dual}, а таким образом и решение исходной задачи \eqref{p5_4.1}  $2\varepsilon$-при\-бли\-жён\-ное. 

Напомним правило построения искомых приближений $\widehat{\mu}$ переменной двойственной задачи \eqref{truss_dual} по работе алгоритма \ref{Switch_Subgrad} для задачи \eqref{p5_4.7}. Положим
$$
\widehat{\mu}_{i}=\frac{1}{h_{f}|I|}\sum_{k\in J}h_{g}{\bf I}[i(k)=i],\ i = 1,\dots,2d, 
$$
где
$$
{\bf I}[\text{predicat}] =
\begin{cases}
1, & \text{predicat}\, = true\\
0, & \text{predicat}\, = false,
\end{cases}
$$
а также $h_f$ и $h_g$ --- величины продуктивных и непродуктивных шагов алгоритма \ref{Switch_Subgrad} соответственно, $I$ и $J$ --- множества номеров продуктивных и непродуктивных итераций алгоритма \ref{Switch_Subgrad}.


Основное преимущество применения алгоритма \ref{Switch_Subgrad} заключается в том, что сложность одной итерации сильно понижается за счёт разреженности векторов $a_i,f$. Из этого следует, что на каждом шаге 
\[
y_{k+1} = y_k - h_g \cdot \nabla g(y_k) \qquad \mbox{или} \qquad y_{k+1} = y_k + h_f \cdot f
\]
в векторе $y$ обновляется не больше, чем $r = O(1)$ элементов. Из того, что в каждом узле встречаются не более, чем $s$ стержней, следует, что обновление $y$ влечёт за собой обновление максимум $sr$ скалярных произведений $\langle a_i,y \rangle$, каждое из которых можно выполнить за $O(r)$ операций. Это сильно упрощает поиск индекса, доставляющего максимум в выражениях \eqref{p5_4.7} для функции $g$ на следующей итерации. В конце концов, положив начальный вектор $y_0$ равным нулю, мы инициализируем и вектор $y$, и скалярные произведения $\langle a_i,y \rangle$ нулями, и они остаются разреженными в течение всей работы алгоритма.

Описанная в предыдущем параграфе идея в работе \cite{Nesterov15} была применена к целому классу разреженных задач оптимизации большой размерности, которых выгодно сводить к негладкой задаче и использовать субградиентные методы. В её основе лежит эффективная процедура для пересчёта максимума $2d$ чисел, если на каждом шаге обновляется только $t = O(1)$ из них, а~начальные значения всех чисел равны нулю.


{\fontsize{10}{10}\selectfont

\begin{table}[h!]
\centering
\caption{\centering\bf Бинарное дерево вычислений для нахождения максимума $2d$ чисел}
\label{tab-truss}
\begin{center}
\begin{tabular}{|*{8}{c |}}\hline
\multicolumn{8}{|c|}{$\max \{p_1, p_2, \ldots, p_{2d} \}$} \\ \hline
\multicolumn{4}{|c|}{$\max \{p_1, p_2, \ldots, p_{2^{q-1}} \}$} & \multicolumn{4}{c|}{$\max \{p_{2^{q-1} +1 }, p_{2^{q-1} +2}, \ldots, p_{2d} \}$} \\ \hline
\multicolumn{8}{|c|}{$\ldots$} \\ \hline

\multicolumn{2}{|c|}{$\max \{ p_1, p_2\}$} & \multicolumn{2}{c|}{$\max \{ p_3, p_4\}$} & \multicolumn{2}{c|}{$\ldots$} & \multicolumn{2}{c|}{$\max \{p_{2d-1}, p_{2d} \}$} \\ \hline
$p_1$ & $p_2$ & $p_3$ & $p_4$ & $\ldots$ & $\ldots$ & $p_{2d-1}$ & $p_{2d}$  \\ \hline
\end{tabular}
\end{center}
\end{table}

}

Опишем вариант организации процедуры нахож\-де\-ния\linebreak  $\max \{p_1, p_2, \ldots, p_{2d} \}$ для обновления $g(y_{k+1})$. Для наглядности выберем $d = 2^{q-1}$ при некотором $q > 0$. В этом случае удобно использовать бинарное дерево вычислений \cite[рис.~1]{Nesterov15}, приведённое в табл.~\ref{tab-truss}.

На первом шаге (нижний уровень дерева) сравниваются пары чисел $(p_1, p_2), \; (p_3, p_4), \ldots (p_{2d-1}, p_{2d})$. На следующем уровне уже найденные на первом шаге максимумы ($\max  \{p_1, p_2 \}, \; \max  \{p_3, p_4 \}$). Далее, рекурсивно процесс продолжается $q$ шагов. На последнем шаге сравниваются два числа $\max \{p_1, p_2, \ldots, p_{2^{q-1}} \}$ и $\max \{p_{2^{q-1} +1 }, p_{2^{q-1} +2}, \ldots, p_{2d} \}$ и после этого находится искомый максимум чисел $p_1, p_2, \ldots, p_{2d}$. Если обновляется одно из чисел $p_i$, то обновить всё дерево можно за всего $q$ операций, идя вверх от данного элемента по направлению к корню дерева. Обновление $s$ чисел $p_i$ приводит к трудоёмкости в $sq$ операций для обновления всего дерева.

Подчеркнём, что в ходе работы алгоритма никогда не хранится полное дерево, имеющее порядка $O(d)$ элементов. Так как начальные значения скалярных произведений равны нулю, то дерево изначально пустое, и только в ходе работы алгоритма начинает заполнятся. Одно из преимуществ использования именно структуры дерева заключается в том, что местонахождение каждого обновляемого элемента можно найти за $q$ операций. Кроме того, обновление одного из элементов $p_i$ влечёт за собой только $O(q)$ операций для обновления максимума. В итоге общее количество операций за всё время работы алгоритма может быть сильно меньше, чем размерность $d$ задачи. 




\section{Идентификация систем} \label{ch3_identification}

В этом разделе мы рассмотрим приложение полуопределённого программирования к обработке сигналов, а именно к оптимальному дизайну эксперимента для оценки параметров линейной динамической системы по данным входа и выхода. Здесь под <<системой>> могут подразумеваться самые различные объекты, от канала связи в мобильной телефонии до самолёта, получающего на вход сигналы управления, и выдающий на выходе такие характеристики как скорость, наклон и т.п. Задача идентификации состоит в том, чтобы получить информацию о системе посредством измерения сигнала на входе и выходе. Оценка параметров является стандартной задачей в статистике и её сложность зависит в основном от самой параметризации системы. Нас здесь интересует другой вопрос, а именно, как подобрать такой сигнал для подачи на вход, чтобы собранные данные стали по возможности более информативными. Это задача оптимального дизайна эксперимента. Стандартным пособием по идентификации систем является книга \cite{Ljung99}.

\subsection{Сигналы, спектры, передаточные функции}

Мы рассматриваем каузальные линейные стационарные динамические системы в дискретном времени. Систему $G$ можно представить в виде линейного оператора, переводящего входной сигнал $\{u_t\}$ в выходной сигнал $\{y_t\}$, $t \in \mathbb Z$. Условие каузальности означает, что значение $y_t$ сигнала на выходе зависит только от предыдущих значений $u_{t-1}, \, u_{t-2}, \, \dots$ сигнала на входе. Стационарность системы эквивалентна тому, что оператор $G$ коммутирует с оператором сдвига по времени $z$, т.е. сдвиг последовательности $u$ приводит к сдвигу последовательности $y$. Линейность оператора означает, что выходной сигнал представляется в виде ряда
\begin{equation} \label{noiseless_system_equation} 
y_t = \sum_{k=1}^{+\infty} g_ku_{t-k}\quad \forall\ t \in \mathbb Z
\end{equation}
с коэффициентами $g_k$, зависящими только от системы. Функция $G(z) = \sum_{k=1}^{+\infty} g_kz^{-k}$ называется \emph{передаточной функцией} системы. С её помощью соотношение \eqref{noiseless_system_equation} можно формально записать в виде $y = Gu$. Если величины $u_t,y_t$ являются векторами, то коэффициенты $g_k$ --- матрицы соответствующего размера. В расчётах оператор $z$ заменяется комплексной переменной, а передаточная функция представляется голоморфной функцией.

Мы будем предполагать, что сигналы $u,y$ являются \emph{квазистационарными}. В этом случае для всех $\tau \in \mathbb Z$ существуют пределы $R_u(\tau) = \lim_{N \to +\infty}\frac{1}{N}\sum_{k=1}^Nu_{t+k}u_{t+k-\tau}^T$, независимые от $t$. Функция $R_u$ называется \emph{функцией корреляции} сигнала $u$. Для пары сигналов $u,e$ можно определить \emph{функцию кросс-корреляции} $R_{ue}(\tau) = \lim_{N \to +\infty}\frac{1}{N}\sum_{k=1}^Nu_{t+k}e_{t+k-\tau}^T$. С помощью этих функций определяются \emph{спектры} сигналов
\[ \Phi_u(\omega) = \sum_{\tau \in \mathbb Z} R_u(\tau) e^{-i\tau\omega},\ \Phi_{ue}(\omega) = \sum_{\tau \in \mathbb Z} R_{ue}(\tau) e^{-i\tau\omega} = \Phi_{eu}(\omega)^*.
\]
Спектр $\Phi_u$ --- комплексно-эрмитовая матрица, и $\Phi_u(-\omega) = \overline{\Phi_u(\omega)}$, $\Phi_{ue}(-\omega) = \overline{\Phi_{ue}(\omega)}$. Блочную матрицу $\begin{pmatrix} \Phi_u & \Phi_{ue} \\ \Phi_{eu} & \Phi_e \end{pmatrix}$ называют \emph{совместным спектром} сигналов $u,e$.

Нетрудно проверить, что спектр сигнала $Gu$ имеет вид $$\Phi_{Gu}(\omega) = G(e^{i\omega})\Phi_u(\omega)G(e^{i\omega})^*.$$ Таким образом, спектр $Gu$ определяется спектром $u$ и значениями передаточной функции $G$ на единичном круге. С помощью спектра также легко посчитать \nk{математическое} ожидание:
\[ \mathbb E uu^T = \frac{1}{2\pi} \int_{-\pi}^{\pi} \Phi_u(\omega)\,d\omega.
\]

В реальном мире сигналы подвержены различного рода возмущениям. Поэтому к идеальному выходному сигналу $Gu$ прибавляется \emph{шум}, моделируемый выражением $He$, где $H(z) = I + \sum_{k=1}^{+\infty} h_kz^{-k}$ --- передаточная функция модели шума, а $e$ --- последовательность независимых одинаково распределённых случайных величин с нулевым средним и вариацией $\lambda_0I$ (обычно гауссовых). При этом функция $H$ предполагается аналитической и обратимой вне единичного диска на комплексной плоскости (включая бесконечно удалённую точку). Таким образом, выходной сигнал получается из входного по формуле $y = Gu + He$ и сам является случайной величиной. Отметим, что спектр $\Phi_e$ не зависит от частоты $\omega$ и равен $\lambda_0I$.

В системах с \emph{открытым контуром} входной сигнал $u$ никак не зависит от выхода $y$, и вследствие и от шума, и поэтому $\Phi_{ue} = 0$. Часто необходимо систему \emph{регулировать}, например, потому что она неустойчива (функция $G$ имеет полюса вне единичного круга). Для этого вводится \emph{обратная связь} с помощью некоторого \emph{регулятора} $K$. Конкретно, полагается $y = r - Ky$, где $r$ --- внешнее возмущение, независимое от шума. В этом случае имеем
\begin{equation} \label{spectra_formula} 
\Phi_u = \lambda_0 SKH (SKH)^* + S\Phi_rS^*,\ \Phi_{ue} = -\lambda_0 SKH,
\end{equation}
где $S = (I + KG)^{-1}$.

\subsection{Оценка параметров}

Задача \emph{идентификации} системы состоит в том, чтобы с помощью конечных отрезков длины $N$ последовательностей $u,y$ получить оценку параметров $g_k,h_k$, или, что эквивалентно, передаточных функций $G,H$ системы. Чтобы ограничить количество оцениваемых величин, обычно предполагают, что $(G,H)$ --- пара рациональных функций, параметризованная конечным числом параметров, собранных в некотором векторе $\theta$. При этом динамика системы задаётся парой $(G,H)$, соответствующей некоторому <<истинному>> значению $\theta_0$ вектора параметров. Это допущение, конечно, является идеализацией, но необходимо для дальнейшего анализа.

Каждому значению $\theta$ параметров системы можно сопоставить \emph{предиктор}
\[ \hat y = y - e = (I - H(\theta)^{-1})y + H(\theta)^{-1}G(\theta)u,
\]
задающий наиболее точный прогноз выхода $y_t$ при известных значениях $y_{t-k},u_{t-k}$, $k \geq 1$ сигналов в прошлом. Исходя из заданных значений сигналов $u_t, \, y_t$, $t = 1, \, \dots, \, N$, оценка $\hat\theta$ параметра вычисляется минимизацией нормы \emph{ошибки предсказания} $\|y - \hat y\|$ по $\theta$. Так как $y$ (а часто и $u$) --- случайная величина, то оценка $\hat\theta$ также будет случайной величиной. При достаточно слабых условиях, накладываемых на сигнал входа $u$ и параметризацию $(G(\theta),H(\theta))$, \nk{математическое} ожидание оценки $\hat\theta$ будет совпадать с истинным значением параметра~$\theta_0$. Её~ковариация будет обратной от \emph{матрицы информации} Фишера
\[ M = \frac{N}{\lambda_0} \mathbb E \psi(\theta_0)\psi(\theta_0)^T,
\]
где $\psi = \frac{\partial \hat y}{\partial\theta}^T = H^{-1}\frac{\partial G}{\partial\theta}u + H^{-1}\frac{\partial H}{\partial\theta}e$ --- градиент предиктора. Отсюда получаем для элементов матрицы информации
\begin{equation} \label{information_matrix} 
M_{kl} = \frac{N}{2\pi\lambda_0} \int_{-\pi}^{\pi} \tr\, F_k(e^{i\omega}) \begin{pmatrix} \Phi_u & \Phi_{ue} \\ \Phi_{ue}^* & \lambda_0I \end{pmatrix} F_l(e^{i\omega})^* \,d\omega,
\end{equation}
где $F_k(z) = \left[H^{-1}(z)\frac{\partial G(z)}{\partial\theta_k}, H^{-1}(z)\frac{\partial H(z)}{\partial\theta_k}\right]$. Заметим, что элементы $F_k$ также являются рациональными функциями от $z$.

\begin{example} \label{ident_example}
Рассмотрим модель авто-регрессии с входом (ARX)
\[ G(z) = \frac{bz^{-1}}{1+az^{-1}},\ H(z) = \frac{1}{1+az^{-1}}.
\]
Параметры системы задаются вектором $\theta = (a,b)$ с <<истинным>> значением $\theta_0 = (a_0,b_0)$, $|a_0| < 1$. Выходной сигнал системы с параметрами $(a,b)$ удовлетворяет соотношению $(1+az^{-1})y = bz^{-1}u + e$, т.е. $y_t + ay_{t-1} = bu_{t-1} + e_t$. Предиктор имеет вид
\[ \hat y_t = y_t - e_t = -ay_{t-1}+bu_{t-1}
\]
и линейный по оцениваемым параметрам $(a,b)$. Поэтому задача оценки сводится к задаче о наименьших квадратах.

Градиент предиктора равен
\[ \psi = \begin{pmatrix} -z^{-1}y \\ z^{-1}u \end{pmatrix} = \begin{pmatrix} \frac{-b_0z^{-2}}{1+a_0z^{-1}} & \frac{-z^{-1}}{1+a_0z^{-1}} \\ z^{-1} & 0 \end{pmatrix}\begin{pmatrix} u \\ e \end{pmatrix}.
\]
Матрицы $F_k$ задаются строками матрицы коэффициентов в выражении справа. Отсюда получаем матрицу информации
\[ M = \frac{N}{2\pi\lambda_0} \int_{-\pi}^{\pi} F \begin{pmatrix} \Phi_u & \Phi_{ue} \\ \Phi_{ue}^* & \lambda_0 \end{pmatrix} F^* \,d\omega,
\]
где $F = \begin{pmatrix} \frac{-b_0e^{-2i\omega}}{1+a_0e^{-i\omega}} & \frac{-e^{-i\omega}}{1+a_0e^{-i\omega}} \\ e^{-i\omega} & 0 \end{pmatrix}$.
\end{example}

\subsection{Дизайн эксперимента}

Теперь допустим, что у нас есть возможность выбора входного сигнала $u$. Мы хотим оценить вектор $\theta$ с как можно меньшей ошибкой. Для этого нужно подобрать спектры $\Phi_u,\Phi_{ue}$ так, чтобы матрица информации \eqref{information_matrix} была большой. Понятно, что есть разные способы измерить величину матрицы информации. По возможности максимизируемую величину выбирают с учётом дальнейшего использования оценки. Например, если на базе оценки будет создаваться новый регулятор для стабилизации системы, то критерием будет служить уровень шума на выходе системы вкупе с этим регулятором.

Для простоты будем полагать, что мы хотим отмаксимизировать след матрицы информации при ограничении на общую энергию сигнала входа. Тогда получим задачу оптимизации
\[ \max_u \sum_k \frac{N}{2\pi\lambda_0} \int_{-\pi}^{\pi} \tr\, F_k(e^{i\omega}) \begin{pmatrix} \Phi_u(\omega) & \Phi_{ue}(\omega) \\ \Phi_{ue}(\omega)^* & \lambda_0I \end{pmatrix} F_k(e^{i\omega})^* \,d\omega:
\]
\[\frac{N}{2\pi} \int_{-\pi}^{\pi} \tr\,\Phi_u(\omega)\,d\omega = c.
\]
Здесь $c$ --- суммарная энергия сигнала входа. Рассмотрим, как преобразовать эту задачу в полуопределённую программу. 

Во-первых, в исходной формулировке задача бесконечномерная, поскольку спектры $\Phi_u,\Phi_{ue}$ являются элементами функциональных пространств. Эту проблему можно решить с помощью введения моментов (см.~п.~\ref{ch3_subs_moments}), используя то обстоятельство, что функции $F_k$ --- рациональные. Представим эти функции в виде $F_k(z) = \frac{p_k(z^{-1})}{q(z^{-1})}$, где $p_k(z) = \sum_{l=1}^{n_k} p_{k,l}z^l$ --- матрично-значные полиномы, а $q(z) = \sum_{k=0}^{n_q} q_kz^k$ --- скалярный полином, такой что $|q(z)|^2$ является общим знаменателем всех используемых в постановке задачи рациональных функций. Введём также (обобщённые) тригонометрические моменты
\begin{equation} \label{moments_definition} 
m_n = \begin{pmatrix} m_{u,n} & m_{ue,n} \\ m_{ue,n}^T & m_{e,n} \end{pmatrix} = \frac{1}{2\pi} \int_{-\pi}^{\pi} \frac{e^{i\omega n}}{|q(e^{i\omega})|^2} \begin{pmatrix} \Phi_u(\omega) & \Phi_{ue}(\omega) \\ \Phi_{ue}(\omega)^* & \lambda_0I \end{pmatrix}\,d\omega
\end{equation}
для всех $n \in \mathbb Z$. Отметим, что моменты $m_n$ вещественны, и $m_{-n} = m_n^T$. Тогда задача примет вид
\[ \max_u \,\sum_k \frac{N}{2\pi\lambda_0} \int_{-\pi}^{\pi} \tr\, \sum_{l,l'=1}^{n_k} p_{k,l} \begin{pmatrix} \Phi_u(\omega) & \Phi_{ue}(\omega) \\ \Phi_{ue}(\omega)^* & \lambda_0I \end{pmatrix} p_{k,l'}^*\frac{e^{-i(l-l')\omega}}{|q(e^{i\omega})|^2} \,d\omega:
\]
\[ \frac{N}{2\pi} \int_{-\pi}^{\pi} \frac{\sum_{k,k'=0}^{n_q} q_kq_{k'} e^{-i(k-k')\omega}}{|q(e^{i\omega})|^2}\tr\,\Phi_u(\omega)\,d\omega = c.
\]
Выражая интегралы через моменты, получаем
\begin{equation} \label{moment_problem1}
\begin{split}
\max_u & \sum_k \frac{N}{\lambda_0} \tr\, \sum_{l,l'=1}^{n_k} p_{k,l} m_{l'-l} p_{k,l'}^*: \\
& N \sum_{k,k'=0}^{n_q} q_kq_{k'} \tr\,m_{u,k'-k} = c.
\end{split}
\end{equation}
Так как в формулировке участвует только конечное число моментов, задача сведена к конечномерной.

В случае системы с открытым контуром внедиагональные блоки матриц $m_n$ равны нулю. Моменты спектра $\Phi_e$ --- числа, зависящие от коэффициентов $q_k$. Таким образом, в качестве переменных выступают только моменты $m_{u,n}$ неотрицательной матрично-значной функции $\Phi_u$ на интервале $[-\pi, \, \pi]$. Дальнейшее преобразование в полуопределённую программу следует стандартным рецептам, используя критерий Каратеодори--Фейера. Этот результат даёт необходимое и достаточное условие на последовательность матриц $m_{-n}, \, \dots, \, m_n$ являться последовательностью тригонометрических моментов некоторой неотрицательной меры на единичном круге. А именно, блочно-тёплицевая матрица
\[ T(m_k) = \begin{pmatrix} m_0 & m_1 & \ddots & m_{n-1} & m_n \\ m_{-1} & m_0 & m_1 & \ddots & m_{n-1} \\ \ddots & \ddots & \ddots & \ddots & \ddots \\ m_{-n} & \ddots & m_{-2} & m_{-1} & m_0 \end{pmatrix},
\]
составленная из этой последовательности, должна быть неотрицательно определённой \cite[Гл.~6, теорема 4.1]{KarlinStudden} (см.~также п.~\ref{ch3_subs_trig_moments}). Таким образом, к линейному условию в \eqref{moment_problem1} следует прибавить ещё полуопределённое условие $T(m_{u,k}) \succeq 0$.

Если в системе присутствует фиксированный регулятор $K$, то функция $\Phi_{ue}$ не зависит от искомого сигнала $r$. Поэтому моменты $m_{ue,n}$ также фиксированы, и условие $T(m_k) \succeq 0$ на оставшиеся переменные $m_{u,n}$ тоже необходимо и достаточно для того, чтобы они были тригонометрическими моментами спектра $\Phi_u$, который можно получить по формуле \eqref{spectra_formula} для некоторого сигнала $r$.

Более сложным является случай, когда оптимизация дизайна эксперимента проводится и по внешнему возмущению $r$, и по регулятору $K$. Дополнительная степень свободы, которую мы получаем выбором регулятора, возводит моменты $m_{ue,n}$ спектра $\Phi_{ue}$ в ранг переменных. Отметим, что $R_{ue}(\tau) = 0$ для всех $\tau < 0$, поскольку $u_t$ зависит только от значений $e$ до момента $t$ включительно. Поэтому функция $f_{ue}(z) = \sum_{\tau \geq 0} R_{ue}(\tau) z^{-\tau}$ аналитична вне единичного диска, и мы получаем линейные условия
\[ \sum_{k=0}^{n_q} q_k m_{ue,s+k} = \frac{1}{2\pi} \int_{-\pi}^{\pi} \frac{e^{i\omega s}\Phi_{ue}(\omega)}{q(e^{-i\omega})} \,d\omega = \frac{1}{2\pi i} \int_{\mathbb T} \frac{f_{ue}(z)z^{s-1}}{q(z^{-1})}\,dz = 0
\]
на переменные $m_{ue,k}$ для всех $s < 0$. Здесь последний интеграл по единичному кругу равен нулю по теореме Коши. Оказывается, что эти линейные условия совместно с полуопределённым условием $T(m_k) \succeq 0$ необходимы и достаточны для того, чтобы последовательность пар $(m_{u,k},m_{ue,k})$ задавала искомые тригонометрические моменты \cite{HildebrandGeversSolari15}.

\begin{example}
Рассмотрим систему из примера \ref{ident_example}. В качестве полинома $q$ в определении моментов \eqref{moments_definition} следует выбрать $1 + a_0z$. Тогда матрица информации $M$ примет вид
\[ \frac{N}{2\pi\lambda_0}\!\int\limits_{-\pi}^{\pi}\!\frac{1}{|q(e^{i\omega})|^2}\!\begin{pmatrix} -b_0e^{-i\omega} & -1 \\ 1+a_0e^{-i\omega} & 0 \end{pmatrix}\!\begin{pmatrix} \Phi_u & \Phi_{ue} \\ \Phi_{ue}^* & \lambda_0 \end{pmatrix}\!\begin{pmatrix} -b_0e^{-i\omega} & -1 \\ 1+a_0e^{-i\omega} & 0 \end{pmatrix}^*\!d\omega 
\]
с диагональными элементами
\[ \frac{N}{\lambda_0} \left(b_0^2m_{u,0} + b_0m_{ue, \, -1} + b_0m_{eu, \, 1} + m_{e, \, 0}\right),
\]
\[ 
\frac{N}{\lambda_0} \left((1+a_0^2)m_{u, \, 0} + a_0m_{u, \, 1} + a_0m_{u, \, -1}\right),
\]
а энергия сигнала входа запишется в виде
\[ N \left( (1+a_0^2)m_{u, \, 0} + a_0m_{u, \, -1} + a_0m_{u, \, 1} \right).
\]
Далее,
\[ m_{e, \, 0} = \frac{1}{2\pi} \int_{-\pi}^{\pi} \frac{\lambda_0}{1+a_0^2+2a_0\cos\omega}\,d\omega = \frac{\lambda_0}{1 - a_0^2},
\]
\[ m_{e,1} = m_{e, \, -1} = \frac{1}{2\pi} \int_{-\pi}^{\pi} \frac{\lambda_0e^{i\omega}}{1+a_0^2+2a_0\cos\omega}\,d\omega = -\frac{a_0\lambda_0}{1 - a_0^2},
\]
\[ m_{u,1} = m_{u, \, -1},\quad m_{ue, \, \mp1} = m_{eu, \, \pm1},\quad m_{ue, \, 0} = m_{eu,\, 0}.
\]
Переменные $m_{ue, \, k}$ равны нулю в случае отсутствия регулятора, и удовлетворяют соотношению $m_{ue, \, -1} + a_0m_{ue, \, 0} = 0$ в общем случае.

В итоге получаем полуопределённую программу
\[ \max \frac{N}{\lambda_0} \left( (1+a_0^2+b_0^2)m_{u, \, 0} + 2a_0m_{u, \, 1} + 2b_0m_{ue, \, -1} + \frac{\lambda_0}{|1 - a_0^2|} \right):
\]
\[ N \left( (1+a_0^2)m_{u, \, 0} + 2a_0m_{u, \, 1} \right) = c,\ m_{ue, \, -1} + a_0m_{ue, \, 0} = 0,
\]
\[ \begin{pmatrix} m_{u, \, 0} & m_{ue, \, 0} & m_{u, \, 1} & m_{ue, \, 1} \\ m_{ue, \, 0} & m_{e, \, 0} & m_{ue, \, -1} & m_{e, \, 1} \\ m_{u, \, 1} & m_{ue, \, -1} & m_{u, \, 0} & m_{ue, \, 0} \\ m_{ue, \, 1} & m_{e, \, 1} & m_{ue, \, 0} & m_{e, \, 0} \end{pmatrix} \succeq 0,
\]
а в случае отсутствия регулятора имеем дополнительные условия $m_{ue,k} = 0$.

В этом случае матрица в полуопределённом условии распадается на два блока. Блок, содержащий переменные $m_{e, \, k}$, автоматически неотрицательно определён, а второй блок даёт эквивалентное линейное условие $|m_{u, \, 1}| \leq m_{u_0}$. В итоге программа становится линейной, и её решение имеет вид
\[ m_{u, \, 0} = -\sign a_0 \cdot m_{u, \, 1} = \frac{c}{(1-|a_0|)^2N}.
\]
Нетрудно проверить, что это решение соответствует постоянному сигналу входа $u_t \equiv \sqrt{\frac{c}{N}}$ при $a_0 < 0$ и альтернирующему $u_t = (-1)^t\sqrt{\frac{c}{N}}$ при $a_0 > 0$.
\end{example}

Подчеркнём, что ключевым элементом вышеизложенного подхода является полуопределённое представление конуса моментов спектров $\Phi_u,\Phi_{ue}$. Любая задача на оптимальный дизайн эксперимента, в которой функция цены и ограничения записываются в виде линейных комбинаций этих моментов, представляется в виде полуопределённой программы. 

\section[Приложения полуопределённого программирования]{Приложения полуопределённого\\ программирования}

В этом разделе мы рассмотрим некоторые задачи, которые можно эффективно решать, сформулировав их в виде полуопределённой программы. Примеры почерпнуты из лекций А.~Бен-Таля и А.С.~Немировского \cite{Ben-Tal}, в которых описаны и другие приложения полуопределённого программирования.

\subsection{Максимальный вписанный в политоп эллипсоид}

Пусть набором $m$ линейных неравенств задан политоп $P = \{ x \in \mathbb R^n \mid Ax \leq b \}$ с непустой внутренностью. Задача состоит в том, чтобы найти вписанный в $P$ эллипсоид
\[ E = \{ x = Cu + c \mid u \in \mathbb R^n,\ \|u\|_2 \leq 1 \}
\]
максимального объёма. 

Эллипсоид представляется в виде аффинного образа единичного шара $B_1$ в пространстве $\mathbb R^n$. Его центр задаётся вектором $c$, а его форма --- положительно определённой матрицей $C$. Объём эллипсоида вычисляется по формуле
\[ Vol(E) = \det\,C \cdot Vol(B_1).
\]

Так как объём единичного шара --- константа, задача сводится к максимизации определителя матрицы $C$ при ограничениях
\[ A(Cu + c) \leq b \qquad \forall\ u \in B_1.
\]
Каждое из этих неравенств имеет вид 
\[ \langle a_i,Cu + c \rangle \leq b_i \quad \Leftrightarrow\quad \langle C^Ta_i,u \rangle \leq b_i - a_i^Tc,
\]
где $a_i$, $i = 1,\dots,m$, --- строки матрицы $A$. Поэтому условие запишется в виде конично-квадратичного ограничения
\[ \left( C^Ta_i, b_i - a_i^Tc \right) \in L^{n+1} \quad \forall\ i = 1,\dots,m.
\]
Определитель $\det\,C$ можно заменить на геометрическое среднее $2^l$ величин $\lambda_1(C),\dots,\lambda_n(C),1,\dots,1$, где $l = \lceil \log_2n \rceil$, а $\lambda_i(C)$ --- собственные значения $C$. Конично-квадратичное представление подграфика геометрического среднего было описано на стр.~\pageref{geometric_mean}, там же полуопределённое представление подграфика спектра.

Для размерности $n = 3$, например, получаем следующую полуопределённую программу:
\[ \max_{t,u,\eta,C,c,\Delta}\,t: \quad \begin{pmatrix} C & \Delta \\ \Delta^T & \diag\,\eta \end{pmatrix} \succeq 0,\quad \diag\Delta = \eta,\quad \Delta_{ij} = 0\quad \forall\ 1 \leq i < j \leq 3;
\]
\[ \eta_1,\eta_2,\eta_3,u_1,u_2 \geq 0,\] 
\[( 2 u_1, \, \eta_1 - \eta_2, \, \eta_1 + \eta_2 ),
( 2 u_2, \, \eta_3 - 1, \, \eta_3 + 1 ),
( 2 t, \, u_1 - u_2, \, u_1 + u_2 ) \in L^3;
\]
\[ \left( C^Ta_i, b_i - a_i^Tc \right) \in L^4 \quad \forall\ i = 1,\dots,m.
\]

\subsection[Минимальный описанный вокруг политопа эллипсоид]{Минимальный описанный вокруг политопа\\ эллипсоид} \label{subs:min_ellipsoid}

Пусть набором $m$ вершин задан политоп $P = \conv\{ x_1,\dots,x_m \} \subset \mathbb R^n$ с~непустой внутренностью. Требуется найти описанный вокруг $P$ эллипсоид
\[ E = \{ x \in \mathbb R^n \mid (x-D^{-1}d)^TD(x-D^{-1}d) \leq 1 \}
\]
минимального объёма.

Эллипсоид представляется выпуклым квадратичным неравенством, задающим единичный шар вокруг некоторой центральной точки $D^{-1}d$ в евклидовой норме, заданной положительно определённой матрицей $D$. Объём эллипсоида вычисляется по формуле
\[ Vol(E) = \frac{Vol(B_1)}{\det\,D}.
\]

Таким образом, задача сводится к максимизации определителя матрицы $D$ при ограничениях
\[ (x_i-D^{-1}d)^TD(x_i-D^{-1}d) \leq 1\quad \forall\ i = 1,\dots,m.
\]
Каждое из этих неравенств переписывается в виде
\[ x_i^TDx_i - 2d^Tx_i + d^TD^{-1}d \leq 1 \ \ \Leftrightarrow\ \ \exists\ s:\ \ \begin{pmatrix} s & d^T \\ d & D \end{pmatrix} \succeq 0,\ x_i^TDx_i - 2d^Tx_i + s \leq 1.
\]
Полуопределённая программа достраивается как и в предыдущем пункте с помощью представлений \nk{подграфиков} геометрического среднего и спектра.

Для размерности $n = 3$ получаем полуопределённую программу
\[ \max_{t,u,\eta,D,d,\Delta,s}\,t: \quad \begin{pmatrix} D & \Delta \\ \Delta^T & \diag\,\eta \end{pmatrix} \succeq 0,\ \ \diag\Delta = \eta,\ \ \Delta_{ij} = 0\ \ \forall\ 1 \leq i < j \leq 3;
\]
\[ \eta_1,\eta_2,\eta_3,u_1,u_2 \geq 0,\]
\[( 2 u_1, \, \eta_1 - \eta_2, \, \eta_1 + \eta_2 ),
( 2 u_2, \, \eta_3 - 1, \, \eta_3 + 1 ),
( 2 t, \, u_1 - u_2, \, u_1 + u_2 ) \in L^3;
\]
\[ \begin{pmatrix} s & d^T \\ d & D \end{pmatrix} \succeq 0,\ x_i^TDx_i - 2d^Tx_i + s \leq 1\quad \forall\ i = 1,\dots,m.
\]

\subsection{Поиск функции Ляпунова} \label{subs:Lyapunov}

Рассмотрим линейную динамическую систему
\[ \dot x(t) = A(t)x(t),
\]
про которую известно, что матрица системы $A(t)$ в любой момент времени принадлежит некоторому множеству ${\cal U}$.

Квадратичная функция $L(x) = x^TXx$, $X \succ 0$, называется \emph{функцией Ляпунова} системы, если существует $s > 0$ такое, что \ag{производная этой функции в силу динамической системы удовлетворяет условию}
\[ \frac{d}{dt}L = x^T(A^TX + XA)x \leq -sL
\]
для всех $t$. Наличие функции Ляпунова гарантирует, что система устойчива и любая её траектория стремится к началу координат.

Достаточным для существования квадратичной функции Ляпунова очевидно является условие
\begin{equation} \label{Lyapunov}
\exists\ s > 0,\ X \succ 0:\quad A^TX + XA \preceq -sX \quad \forall\ A \in {\cal U}.
\end{equation}
Для данной матрицы $A \in {\cal U}$ это условие не является полуопределённым, так как оно не линейно по искомым переменным $s,X$.

Однако можно воспользоваться тем фактом, что матричные неравенства однородные по $X$. Поэтому условие $X \succ 0$ можно заменить на условие $X \succeq I$. С другой стороны, в силу конечной обусловленности положительно определённой матрицы $X$ существование константы $s > 0$, удовлетворяющей неравенству $A^TX + XA \preceq -sI$, гарантирует существование константы $s' > 0$ такой, что имеет место неравенство $A^TX + XA \preceq -s'X$. Поэтому \eqref{Lyapunov} эквивалентно условиям
\[ \exists\ s > 0,\ X \succeq I: \quad A^TX + XA \preceq -sI\quad \forall\ A \in {\cal U}.
\]

Можно ли эту проблему решить с помощью полуопределённой программы, зависит от множества ${\cal U}$. Так как условия выпуклы по $A$, задача решается в случае, когда ${\cal U}$ --- политоп, т.е. выпуклая оболочка конечного числа матриц $A_1,\dots,A_m$.

Тогда поиск функции Ляпунова сводится к полуопределённой программе
\[ \max_{s,X}\,s:\ X \succeq I,\ -sI - A_i^TX - XA_i \succeq 0\ \forall\ i = 1,\dots,m.
\]
Если существует строго положительное допустимое значение переменной $s$, то соответствующая матрица $X$ задаёт функцию Ляпунова. В противном случае достаточное условие \eqref{Lyapunov} не выполняется.

\subsection{Поиск линейного стабилизирующего регулятора}

Рассмотрим линейную управляемую систему
\[ \dot x = Ax + Bu,
\]
где $x$ --- вектор состояния системы и $u$ --- вектор управления. Задача состоит в построении линейного закона управления $u = Kx$, стабилизирующего систему. 

Достаточным условием устойчивости управляемой регулятором $K$ системы является существование квадратичной функции Ляпунова $L(x) = x^TXx$, $X \succ 0$, удовлетворяющей условию
\begin{equation} \label{Lyapunov_control} 
\exists\ s > 0:\quad \frac{d}{dt}L = x^T((A+BK)^TX + X(A+BK))x \leq -sL = -sx^TXx
\end{equation}
для всех $x$. Условие эквивалентно нелинейному матричному неравенству
\[ \exists\ s > 0,\ X \succ 0:\quad (A+BK)^TX + X(A+BK) \preceq -sX.
\]
В этом случае нелинейность происходит не только от правой части, но и от произведения $XBK$ в левой части неравенства, в которое входят обе матричные переменные $K,X$.

Если произведение $sX$ можно устранить как в п.~\ref{subs:Lyapunov}, то от произведения $XBK$ можно избавиться умножением неравенства на $Y = X^{-1}$ слева и справа. Получаем неравенство
\[ \exists\ s > 0,\ Y \succ 0:\quad YA^T + Z^TB^T + AY + BZ \preceq -sY,
\]
которое линейно по $Y$ и $Z = KY$. В итоге получаем полуопределённую программу
\[ \max_{s,Y,Z}\,s:\ Y \succeq I,\ -sI - (AY + BZ) - (AY + BZ)^T \succeq 0.
\]
Если существует строго положительное допустимое значение переменной $s$, то из соответствующей матрицы $Y$ можно восстановить функцию Ляпунова и управление по формулам $X = Y^{-1}$, $K = ZX$. В противном случае достаточное условие \eqref{Lyapunov_control} не выполняется.

\subsection{Задача о максимальной клике} \label{subs:MaxClique}

Задачей о максимальной клике (MaxClique) называют проблему вычисления кликового числа данного графа $G$. Эта комбинаторная задача входит в список Карпа NP-полных проблем \cite{Karp72}.

\begin{defin}
\emph{Кликой} графа $G$ называют подмножество $S$ вершин такое, что любые две вершины из $S$ соединены ребром. \emph{Максимальной} кликой называется клика, которая перестаёт быть кликой при добавлении любой дополнительной вершины. \emph{Кликовым числом} $\alpha(G)$ графа $G$ называется мощность наибольшей клики.
\end{defin}

Верхней оценкой кликового числа является \emph{$\vartheta$-функция Ловаша}, которую можно вычислить полуопределённой программой
\[ \max_{X \succeq 0}\,\langle X,{\bf 1} \rangle:\quad X \bullet A_{\bar G} = 0,\ \tr\,X = 1.
\]
Здесь $A_G$ --- матрица смежности графа, $\bar G$ --- комплементарный к $G$ граф, а ${\bf 1}$ обозначает матрицу, все элементы которой равны единице. Умножение $\bullet$ по-элементное, так называемое {\it произведение Адамара}.

Действительно, пусть $V$ --- множество вершин графа $G$, а $S \subset V$ --- наибольшая его клика, мощности $k$. Определим матрицу $X = (X_{ij})$ посредством
\[ X_{ij} = \left\{ \begin{array}{rcl} \frac{1}{k},&\quad& i,j \in S, \\
0,& & \{i,j\} \not\subset S. \end{array} \right.
\]
Тогда имеем
\[ \tr\,X = 1,\quad X \succeq 0,\quad X \bullet A_{\bar G} = 0,\quad \langle X,{\bf 1} \rangle = k.
\]
Из этого следует, что матрица $X$ допустима, и поэтому её цена $k$ не превышает оптимальное значение, т.е. $\vartheta(G)$.

Кликовое число графа $G$ можно представить в виде оптимального значения коположительной программы
\[ \min_{Z \in {\cal COP}^n}\,\alpha:\quad Z = \alpha(I + A_{\bar G}) - {\bf 1} = (\alpha - 1){\bf 1} - \alpha A_G.
\]
Заменяя коположительный конус ${\cal COP}^n$ на его внутреннюю полуопределённую аппроксимацию ${\cal K}_0$ (см. пример \ref{example_copositive}), получаем верхнюю оценку кликового числа $\alpha(G)$ оптимальным значением полуопределённой релаксации
\[ \min_{Z \succeq 0}\,\lambda:\quad Z \leq \lambda(I + A_{\bar G}) - {\bf 1} = (\lambda - 1){\bf 1} - \lambda A_G.
\]
Эта релаксация сложнее, чем $\vartheta$-функция Ловаша, поскольку помимо матричного неравенства порядка $n$ она содержит ещё и $O(n^2)$ линейных неравенств, но она и сильнее.

\section{Приложения метода условного градиента}

\subsection{Аппроксимация сигналов} \label{ch3_sect_ls_regr}
В данном разделе изложение во многом следует~\cite{Bubeck15}. Рассмотрим проблему
аппроксимации сигнала~$Y \, \in \, \mathbb{R}^n$  линейной комбинацией
некоторых векторов $s_1, \, \ldots, \, s_N \, \in \, \mathbb{R}^n$:
$\sum_{i = 1}^N x^i s_i$, $x = (x^1, \, \ldots, \, x^N) \, \in \, \mathbb{R}^N$.
Момент времени $t$ фиксирован.
Зафиксируем также некоторое вещественное число~$\lambda$. Для того, чтобы обобщённый полином $\sum_{i = 1}^N x^i s_i$ аппроксимировал сигнал~$Y$,
коэффициенты полинома $x^i$ можно выбирать, например,
из условия минимизации нормы ошибки с некоторой $L_1$~-регуляризацией:
$$
\min_{x \, \in \, \mathbb{R}^N} \left\| Y - \sum_{i = 1}^N x^i s_i \right\|_2^2 + \lambda \| x \|_1.
$$
Обозначим через $S$ вещественную матрицу размера $\mathbb{R}^{n \times N}$,
столбцы которой образованы векторами $s_i$. Тогда регуляризованную задачу
можно заменить на задачу с ограничением
(где $r \, \in \, \mathbb{R}$ фиксировано):
\begin{equation}
    \min_{\|x \|_1 \leq r} \| Y - S x \|^2_2,
\end{equation}
что эквивалентно
\begin{equation}
\label{ch3_ls_regr_main2}
    \min_{\|x \|_1 \leq 1} \| Y/r - S x \|^2_2.
\end{equation}

Предположим, что количество векторов $s_i$ очень большое,
настолько большое, что число~$N$ экспоненциально
зависит от $n$. Можно ли в таком случае создать алгоритм решения задачи~\eqref{ch3_ls_regr_main2}, сложность которого полиномиально зависит от размерности~$n$? В общем случае это невозможно.
Поэтому наложим следующее условие на структуру
векторов~$s_i$:
для любого $y \, \in \, \mathbb{R}^n$ за полином~$p(n)$
от размерности $n$ можно
решить задачу линейного программирования
$$
\min_{1 \leq i \leq N} y^T s_i.
$$

Кроме того, пусть все $s_i$ ограничены сверху по $2$-норме некоторым положительным числом~$m$:
$\| s_i \|_2 \leq m$, $i = 1, \, 2, \, \ldots, \, N$.

Итак, необходимо решить задачу минимизации функции $f(x) = \frac{1}{2} \| Y - S x \|_2^2$
на $l_1$-шаре в $\mathbb{R}^N$ за полиномиальное время от $n$. На первый взгляд, это невозможно,
ведь одно только вычисление всех компонент $x \, \in \, \mathbb{R}^N$ займёт линейное время от $N$. С другой стороны, частичное вычисление компонент вектора, т.е. использование
свойства разреженности, может спасти ситуацию. Применим для решения задачи
метод условного градиента (см. п.~\ref{ch2_sect_fw}).

Заметим, что
$$
\nabla f(x) =  S^T ( Sx - Y ).
$$
Пусть $z_k = S x_k - Y \, \in \, \mathbb{R}^n$ уже вычислено, тогда для выполнения шага~$3$ алгоритма Франк--Вульфа на стр.~\pageref{ch2_alg_cond_grad} ($y_{k} = \argmin_{y \, \in \, S} \langle \nabla f(x_k), \, y \rangle$) нужно найти координату $i_k \in \{ 1, \, 2, \, \ldots, \, N \}$, которая максимизирует  $| \nabla f(x_k)^i |$ (это можно сделать, решив две задачи
максимизации: $\max_i \, \langle s_i, \, z_k\rangle$ и $\max_i \, - \langle s_i, \, z_k\rangle$).
Таким образом, шаг $3$ алгоритма Франк--Вульфа выполняется за полином $p(n)$ от $n$.

Следующая точка $x_{k+1}$ будет линейной комбинацией точки $x_k$ и единичного вектора $e^{i_k} \in \mathbb R^N$. Таким образом, $x_{k+1}$ вычисляется за $O(k)$ операций. Вектор $z_{k+1}$ является линейной комбинацией $z_k$, $Y$ и $s_{i_k}$ и вычисляется за $O(n)$ операций.

Общая сложность первых $k$ итераций, таким образом, составляет\linebreak  $O(k \cdot p(n) + k^2)$ операций. Для оценки скорости сходимости установим степень гладкости целевой функции. Для любых $x,x'$ в $L_1$-шаре имеем
\[ \|\nabla f(x) - \nabla f(x')\|_{\infty} = \| S^TS(x-x') \|_{\infty}
\leq m^2 \|x - x'\|_1,
\]
так как все элементы матрицы $S^TS$ ограничены константой $m^2$ по модулю. Таким образом, целевая функция имеет константу Липшица $m^2$, а допустимое множество~--- диаметр~2. Отсюда получаем оценку (см.~\eqref{ch2_eq_rate_FW}):
\[ f(x_k) - f(x_*) \leq \frac{8m^2}{k+2}.
\]

\subsection{Выявление аномалий}

Допустим, что даны точки $x^1, \, \dots, \, x^m \, \in \, \mathbb R^n$, представляющие собой  реализации некоторого распределения. Нам необходимо установить, может ли некая новая точка $x \in \mathbb R^n$ быть реализацией того же распределения, или она \emph{аномальная} и представляет собой результат возникновения потенциальной проблемы. Один из возможных подходов к этой задаче, это построение некоторого простого множества, содержащего точки $x^1, \, \dots, \, x^m$ и характеризующего нормальные реализации. Если новая точка окажется вне этого множества, она классифицируется как аномальная.

В качестве такого множества можно взять, например, шар наименьшего радиуса, содержащего все точки $x^1, \, \dots, \, x^m$ \cite{Lan2020slides}. Данная задача формулируется в виде проблемы оптимизации
\[ \min_{c,r}\,r:\quad \langle x^i - c,x^i - c \rangle \leq r,\ \forall\ i = 1, \, \dots, \, m.
\]
Перейдём к двойственной задаче. Введём множители Лагранжа $\lambda_i$. Получаем
\[ \max_{\lambda \geq 0}\min_{c,r}\,\left(r + \sum_{i=1}^m \lambda_i (\langle x^i - c,x^i - c \rangle - r)\right).
\]
Минимизация по $r \in \mathbb R$ приводит к условию $\sum_{i=1}^m \lambda_i = 1$, минимизация по $c \in \mathbb R^n$ к условию $c = \sum_{i=1}^m \lambda_i x^i$. Собрав точки $x^i$ в матрицу $X \in \mathbb R^{n \times m}$, получим двойственную задачу
\[ \max_{\lambda \in \Delta_m}\, ( \langle \lambda,\diag(X^TX) \rangle - \lambda^TX^TX \lambda ),
\]
где $\Delta_m = S_m(1)= \left\{ \lambda \geq 0 \,\left|\, \sum_{i=1}^m \lambda_i = 1 \right. \right\}$~--- стандартный (вероятностный) симплекс.

Эту задачу можно решать методом Франк--Вульфа (см. п.~\ref{ch2_sect_fw}). На шаге $k$ алгоритма нам нужно максимизировать градиент целевой функции в текущей точке $\lambda^k$ (в данном разделе индекс номера итерации у $\lambda$ поставлен сверху, поскольку есть нижний индекс, отвечающий номеру компоненты), т.е. линейную функцию $\diag(X^TX) - 2X^TX\lambda^k$ по симплексу, т.е. найти максимальную компоненту $i_k$ этого градиента. Следующий двойственный вектор $\lambda^{k+1}$ вычисляется как линейная комбинация предыдущего вектора $\lambda^k$ и единичного вектора $e^{i_k} \in \Delta_m$. Отсюда следует, что вектор $\lambda^k$ имеет не более $k+1$ ненулевых компонент.

При условии ограничения на 2-норму точек $x^i$ скорость сходимости получается аналогичной полученной в предыдущем п.~\ref{ch3_sect_ls_regr}.

\section[Пример использования метода эллипсоидов для решения двойственной задачи малой размерности]{Пример использования метода эллипсоидов\\ для решения двойственной задачи малой\\ размерности%
\sectionmark{Пример использования метода эллипсоидов}
}
\sectionmark{Пример использования метода эллипсоидов}
\label{comp_prox}

Рассмотрим ещё один пример оптимизационной задачи~\cite{GasKamMen_2016} (см. также теорему~\ref{tomogravity}), для которой оказывается выгоднее решать двойственную задачу, а не прямую задачу. Задачи такого вида часто необходимо решать на итерациях
некоторых современных численных методов оптимизации~\cite{Ben-Tal} (см. также раздел~\ref{restarts}): 
\begin{equation}
\label{ch1_eq_cx_log}
\min_{x\in S_n \left( 1 \right)} c^T x +\left\| x 
\right\|_a^2 +\gamma \sum\limits_{i=1}^n {x_i \log x_i },
\end{equation}
где $S_n(1) = \left \{ x \ge 0 \, \left | \, \sum\limits_{i = 1}^n x_i = 1 \right. \right \}$, $1 < a <2$.

Задачу~\eqref{ch1_eq_cx_log} можно переписать в следующем,
почти <<сепарабельном>> виде:
$$
\min_{\begin{array}{c}
 x\in S_n \left( 1 \right),\;\left\| x \right\|_a^a \le t^{a \mathord{\left/ 
{\vphantom {a 2}} \right. \kern-\nulldelimiterspace} 2}, \\ 
 0\le t\le n^{2 \mathord{\left/ {\vphantom {2 a}} \right. 
\kern-\nulldelimiterspace} a},\;0\le x_i \le 1, \; i=1, \, \ldots, \, n \\  \end{array}} 
c^T x +t+\gamma \sum\limits_{i=1}^n {x_i \log x_i }. $$
С помощью метода множителей Лагранжа получаем, что необходимо решить задачу
\begin{equation}
\label{ch1_eq_dual_maxglam}
\max_{\lambda_1 \, \in \, \mathbb{R}, \; \lambda_2 \ge 0} 
 \tilde{G} (\lambda),
\end{equation}
где $\lambda = (\lambda_1, \, \lambda_2)$ и
$$
\tilde {G}\left( \lambda_1, \, \lambda_2 \right)=\mathop {\min }\limits_{\begin{array}{c}
 0\le t\le n^{2 \mathord{\left/ {\vphantom {2 a}} \right. 
\kern-\nulldelimiterspace} a}, \\ 
 0 \le x_i \le 1, \; i=1, \, \ldots, \, n\; \\ 
 \end{array}} \left\{ \sum\limits_{i=1}^n {c_i x_i } +t+\lambda _1 \cdot 
\left( {\sum\limits_{i=1}^n {x_i } -1} \right) + \right.
$$
$$
\left.
+\, \lambda _2 \cdot \left( 
{\sum\limits_{i=1}^n {x_i^a } -t^{a \mathord{\left/ {\vphantom {a 2}} 
\right. \kern-\nulldelimiterspace} 2}} \right)+\gamma \sum\limits_{i=1}^n 
{x_i \log x_i }  \right\} .
$$

Прямые переменные $x$, $t$ и двойственные $\lambda_1$, $\lambda_2$ связаны соотношением для $t(\lambda)$:
$$
t\left( \lambda \right)=\min \left\{ {\left( {\frac{\lambda _2 a}{2}} 
\right)^{\frac{2}{2-a}}, \, \, n^{\frac{2}{a}}} \right\},
$$
а $x(\lambda)$ можно найти численным способом,
решив $n$ задач одномерной выпуклой оптимизации на отрезке~$[0, \, 1]$, например, с помощью метода деления отрезка пополам,
или метода Фибоначчи / золотого сечения (см.~\cite{Vas_met_opt} и раздел~\ref{subs:bisection}).
Таким образом, если задаться некоторой точностью $\sigma \, > \, 0$, то за время $O(n \log (n/\sigma))$
можно найти решение прямой задачи~--- такой $x^{\sigma} (\lambda)$, что
\begin{equation}
\label{ch3_eq_exB_x_eps}
\| x^{\sigma}(\lambda) - x(\lambda) \|_1 = O(\sigma).
\end{equation}

Попробуем теперь оценить время решения двойственной задачи~\eqref{ch1_eq_dual_maxglam},
следуя~\cite{GasKamMen_2016}.
По формуле Демьянова--Данскина--Рубинова (см. п.~\ref{ch1_sect_dem-dan}):
\[
\frac{\partial \tilde {G}}{\partial \lambda _1 }=\sum\limits_{i=1}^n {x_i 
\left( \lambda \right)} -1,
\quad
\frac{\partial \tilde {G}}{\partial \lambda _2 }=\sum\limits_{i=1}^n {x_i 
\left( \lambda \right)^a} -t\left( \lambda \right)^{a \mathord{\left/ 
{\vphantom {a 2}} \right. \kern-\nulldelimiterspace} 2}.
\]
Поскольку~\eqref{ch1_eq_dual_maxglam}
является
задачей оптимизации на двумерной плоскости (т.е. в пространстве малой
размерности), то её можно решить, например, методом эллипсоидов (см. п.~\ref{ch2_subsect_ellips_meth}).

Далее попробуем оценить
сверху двойственные
переменные,
используя \textit{слейтеровский подход}.
Как обычно, $\tilde {G}^{\ast} = \tilde {G} (\lambda^{\ast})
= \max_{\lambda_1 \, \in \, \mathbb{R}, \; \lambda_2 \ge 0} 
 \tilde{G} (\lambda)$
(в приводимой
далее выкладке, приводящей к формуле~\eqref{ch1_eq_lamleC}), для упрощения записи мы опускаем нижний
индекс * у $\lambda$). 
Из сильной двойственности (см.~раздел~\ref{ws_duality})
\[
-\left\| c \right\|_\infty -\gamma \log n\le \tilde {G}^{\ast} 
\le \sum\limits_{i=1}^n {c_i \bar {x}_i } +\bar {t}+\lambda _1 \cdot \left( 
{\sum\limits_{i=1}^n {\bar {x}_i } -1} \right)+
\]
\[
+\lambda _2 \cdot \left( {\sum\limits_{i=1}^n {\bar {x}_i^a } -\bar {t}^{a 
\mathord{\left/ {\vphantom {a 2}} \right. \kern-\nulldelimiterspace} 2}} 
\right) + \gamma \sum\limits_{i=1}^n {\bar {x}_i \log \bar {x}_i }.
\]
Если $\lambda _1 \ge 0$, положим $\bar {t}=1$, $\bar {x}_i =1 \mathord{\left/ 
{\vphantom {1 {\left( {2n} \right)}}} \right. \kern-\nulldelimiterspace} 
{\left( {2n} \right)}$, $i=1, \, \ldots, \, n$. 
Тогда
\[
\frac{1}{2}\lambda _1 +\frac{1}{2}\lambda _2 \le 2\left\| c \right\|_\infty 
+2\gamma \log \left( n \right)+1.
\]
Если $\lambda _1 < 0$, положим $\bar {t}=8$, $\bar {x}_i =2 \mathord{\left/ 
{\vphantom {2 n}} \right. \kern-\nulldelimiterspace} n$, $i=1, \, \ldots, \, n$. 
Тогда
\[
\left| {\lambda _1 } \right|+\frac{1}{2}\lambda _2 \le 3\left\| c 
\right\|_\infty +2\gamma \log \left( {2n} \right)+8.
\]
В любом случае, можно гарантировать 
\begin{equation}
\label{ch1_eq_lamleC}
\left\| {\lambda^{\ast} } \right\|_1 \le 6\left\| c \right\|_\infty +4\gamma 
\log \left( {2n} \right)+16\mathop =\limits^{\text{def}} C.
\end{equation}
Следовательно, задачу~\eqref{ch1_eq_dual_maxglam} можно
переписать:
\begin{equation}
\label{ch1_eq_dual_minglamC}
\min_{\lambda_1 \, \in \,
{\mathbb{R}},\;\lambda _2 \, \ge \, 0,\;\left\| \lambda \right\|_1 \, \le \, C}  -\tilde {G}\left( \lambda \right).
\end{equation}

К сожалению, мы можем найти только приближённое решение данной задачи. Для этого найдём приближённое значение градиента.  Точнее говоря, в силу (\ref{ch3_eq_exB_x_eps}), (\ref{ch1_eq_lamleC}) мы можем
найти на множестве $\left\{(\lambda_1,\lambda_2)\in\mathbb{R}\otimes\mathbb{R}_{+}:~\|\lambda\|_1\le C\right\}$ для задачи (\ref{ch1_eq_dual_minglamC}) $\delta ={\rm O}\left( {C\sigma } \right)$-градиент\footnote{Вектор $g$ называется $\delta$-(суб)градиентом выпуклой функции $f$ на множестве $Q$, если для всех $x,y\in Q$ верно $f(y) \ge f(x) +\langle g, y-x \rangle - \delta$. Отметим, что из выпуклости функции при этом следует, что $x,y\in Q$ верно $f(y) \ge f(x) +\langle \nabla f(x), y-x \rangle - \delta$ для всех $\nabla f(x) \in \partial f(x)$, см. раздел~\ref{ch2_subs_nonsmooth}.}
$\nabla_{\delta} \left(-\tilde{G}\left( \lambda \right)\right)$ (см., например,
\cite{Polyak}). Если использовать метод эллипсоидов в
пространстве размерности $r$ (в нашем случае $r=2$), где вместо градиента используется $\delta$-градиент (чаще говорят $\delta$-субградиент, но в нашем случае можно
говорить о градиенте), то при $\delta = {\rm O}\left( \varepsilon \right)$ имеют место следующие оценки~\cite{NemirYd79}, \cite{Gladin2020}:
\begin{equation}
\label{ch3_ex_comp_estim}
\tilde{G}(\lambda^{*}) - \tilde{G}(\lambda^N)
\le \varepsilon ,
\quad
N={\rm O}\left( {r^2\log \left( {C \mathord{\left/ {\vphantom {C \varepsilon
}} \right. \kern-\nulldelimiterspace} \varepsilon } \right)} \right).
\end{equation}
При этом стоимость одной (без учета стоимости вычисления субградиента) итерации будет ${\rm O}\left( {r^2} \right)$. Число итераций можно сократить в $\sim r$ раз, немного увеличив сложность итерации (см., например, \cite{Bubeck15} и \cite{Gladin2021}).

Резюмируя написанное выше, получим, что в нашем случае ($r = 2$) число итераций метода эллипсоидов равно ${\rm O}\left(\log\left(C/\varepsilon\right)\right)$, а полная стоимость каждой
итерации равна ${\rm O}\left( {n\log \left( {{nC} \mathord{\left/ {\vphantom
{{nC} \varepsilon }} \right. \kern-\nulldelimiterspace} \varepsilon }
\right)} \right)$.

Однако решение задачи (\ref{ch1_eq_dual_maxglam}) (или (\ref{ch1_eq_dual_minglamC})), в смысле (\ref{ch3_ex_comp_estim}), ещё не гарантирует
возможность точного восстановления решения задачи (\ref{ch1_eq_cx_log}). Для того чтобы
показать, что метод эллипсоидов с той же по порядку точностью $\varepsilon$
позволяет восстанавливать (без каких бы то ни было существенных
дополнительных затрат) решение задачи (\ref{ch1_eq_cx_log}), нужно воспользоваться
\textbf{\textit{прямо-двойственностью}} этого метода~\cite{Nemir_acc_cert2010}. Ввиду компактности множества (единичный
симплекс), на котором ведётся оптимизация, в прямом пространстве и сильной
выпуклости функционала прямой задачи (\ref{ch1_eq_cx_log}) мы не просто восстанавливаем из
прямо-двойственной процедуры метода эллипсоидов решение задачи (\ref{ch1_eq_cx_log}) с
точностью по функционалу (прямой задачи) порядка $\varepsilon$, но и делаем
это в нужном нам более сильном смысле (сравните также с \cite[п.~4.6]{Devolder_thesis}  и \cite[п.~2]{GasKamMen_2016}). Формула~5.5.15
из \cite{Ben-Tal} и результаты раздела~\ref{restarts} гарантируют при этом справедливость теоремы~\ref{th 5}.
\section[Пример относительно гладкой оптимизационной задачи]
{Пример относительно гладкой\\ оптимизационной задачи%
\sectionmark{Относительно гладкая задача}
}
\sectionmark{Относительно гладкая задача}
\label{rel_smooth}

В качестве примера применимости свойства~\eqref{eq_rel_smooth} относительной гладкости функции рассмотрим задачу {\it $D$-оптимального плана эксперимента} ($D$-optimal design). В~\cite{LuNesterov} на базе концепции относительной гладкости обоснована возможность использования  для этой задачи неускоренного метода градиентного типа (см. п.~\ref{gradMethod}). До появления статьи~\cite{LuNesterov} для такой задачи были известны оценки сложности для методов внутренней точки (эти методы описываются выше в разделе~\ref{ch2_subsectIPM}) или метода Хачияна и его аналогов. Метод Хачияна несколько похож на метод условного градиента (см. раздел~\ref{ch2_sect_fw}). 

Приведём постановку задачи $D$-оптимального плана эксперимента. Пусть $H \in \mathbb{R}^{m \times n}$ --- матрица ранга $m$, $n \geqslant m+1$. Рассмотрим задачу
\begin{equation}\label{p0_a23}
\min_{x} \, F(x) := -\log \det \left( HXH^{T} \right):\quad \langle {\bf e}, \, x \rangle = 
\sum_{k = 1}^{n} x_k = 1, \ x_k \geq 0,\ k = 1, \dots, n,
\end{equation}
где $X := \diag(x)$~--- диагональная матрица, ${\bf e} = (1, \, 1, \, \ldots, \, 1)$. Задача $D$-оптимального плана эксперимента возникает как лагранжева двойственная задача к хорошо известной в вычислительной геометрии задаче о нахождении эллипсоида минимального объёма, покрывающего заданный набор точек (см. п.~\ref{subs:min_ellipsoid}). Она имеет приложения в вычислительной  статистике~\cite{CrouxHaesbroekRousseesuv} и анализе данных~\cite{Edwin}. Исследования по задаче об эллипсоиде минимального объёма, который покрывает заданный набор точек, начались с работы~\cite{John_ell}. Это означает, что разные подходы к данной задаче исследуются уже более 60 лет. Относительно современный обзор достижений по ней можно найти в~\cite{ToddEllipsoid}.

По задаче $D$-оптимального плана эксперимента известно довольно много работ, основанных на применении методов внутренней точки (см., например, \cite{KhachiyanToddEll, ZhangEll}), а также метода Хачияна или его аналогов (см., например \cite{KhachiyanEllipsoid, AhipasaogluTodd}). В работе~\cite{LuNesterov} показана применимость к данной задаче  неускоренного градиентного метода с относительной гладкостью. В настоящем пособии приведён аналог этого метода в модельной общности~--- алгоритм~\ref{Grad_Model}. Для задачи~\eqref{p0_a23} за счёт подходящего выбора прокс-функции \begin{equation}\label{proxoptdesigh}
d(x) := -\sum^{n}_{j=1}\log\,x_{j} = -\log\det\,X
\end{equation}
удаётся вывести из оценки скорости сходимости для неускоренного градиентного алгоритма~\ref{Grad_Model} результат \cite[Теорема 4.1]{LuNesterov} о линейной скорости сходимости. 
В качестве информации о гладкости задачи \eqref{p0_a23} используем условие относительной гладкости \eqref{eq_rel_smooth}, в котором в качестве специальной функции $V(x,y)$ выступает дивергенция Брэгмана \eqref{div-breg} с прокс-функцией \eqref{proxoptdesigh}. Покажем, что функция $F$ из \eqref{p0_a23} $L$-гладкая на $\mathbb{R}^{n}_{++}$ относительно выбранной прокс-функции $d$ с константой 
\begin{equation} \label{Leq1}
L = 1.
\end{equation}

Так как $F$ и $d$ дважды дифференцируемы на положительном ортанте, а следовательно, во всех точках внутренности допустимого множества 
\[ Q = \Delta_n = \{ x \geq 0 \mid \langle {\bf e},x \rangle = 1 \}
\]
задачи \eqref{p0_a23}, то условие $L$-гладкости $F$ относительно $d$ эквивалентно условию $\nabla^2 F(x)\preceq L\nabla^2 d(x)$. Имеем $\nabla F(x) = \diag (-C)$, а также $\nabla^2 F(x) = C \circ C$, где  $C := H^T (HXH^T)^{-1} H$ --- вспомогательная матрица, а $C\circ C$ --- произведение Шура матрицы $C$ саму на себя. Если $U = H X^{\frac{1}{2}}$, то $U^T (UU^T)^{-1} U \preceq I$, поскольку левая часть этого матричного неравенства является оператором проектирования. Поэтому $X^{\frac{1}{2}} H^T (HXH^T)^{-1} H X^{\frac{1}{2}} \preceq I$. Умножив обе части данного матричного неравенства $X^{-\frac{1}{2}}$, получим $C \preceq X^{-1}$. Поэтому 
\begin{equation}\label{eqfunctrel-smooth}
\nabla^2 F(x) = C\circ C \preceq C\circ X^{-1} \preceq X^{-1} \circ X^{-1} = X^{-2} = \nabla^2 d(x),
\end{equation}
где первое и второе матричные неравенства вытекают из $C \preceq X^{-1}$ c учётом того, что произведение Шура двух симметричных положительно полуопределённых матриц есть снова симметричная положительно полуопределённая матрица. Соотношения \eqref{eqfunctrel-smooth} означают, что $F$ из \eqref{p0_a23} $1$-гладкая относительно дивергенции Брэгмана прокс-функции $d$ из \eqref{proxoptdesigh} на $\mathbb{R}^{n}_{++}$.
 
Напомним, что стандартная схема итерации методов первого порядка с общей дивергенцией Брэгмана имеет следующий вид (см. раздел \ref{model}):
\begin{equation}\label{p0_a10}
x_{k+1}:=\argmin\limits_{x\in Q}\left\{f(x_{k})+\langle\nabla F(x_{k}), x-x_{k}\rangle+LV(x,x_{k})\right\}.
\end{equation}
Итерация вида~\eqref{p0_a10} приводит к необходимости решать вспомогательные подзадачи вида
\begin{equation}\label{p0_a11}
\min_{x\in Q}\left\{\langle c,x\rangle+d(x)\right\}
\end{equation}
для соответствующего однозначно определяемого линейного функционала $c$. Если, как имеет место в рассматриваемом случае, функция $d$ обладает барьерным свойством по отношению к множеству $Q$, т.е. $\lim_{x \to \partial Q} d(x) = +\infty$, то шаг \eqref{p0_a10} можно записать в виде 
$$
\nabla d(x_{k+1}) = \nabla d(x_{k}) - \frac{1}{L}\nabla F(x_k). 
$$
Для функций $d$ общего вида разрешение этого уравнения по отношению к $x_{k+1}$ приводит к нелинейной системе уравнений. Ниже мы покажем, однако, что в случае прокс-функции \eqref{proxoptdesigh} вычисление следующей точки $x_{k+1}$ сводится к выпуклой одномерной задаче.

Приведём формулировку и доказательство результата \cite[теорема~4.1]{LuNesterov} о скорости сходимости градиентного метода \eqref{p0_a10} для задачи \eqref{p0_a23}. Этот результат аналогичен выводам пункта \ref{gradMethod} для алгоритма \ref{Grad_Model}, но в качестве точки-выхода используется $x_N$ вместо усреднения, выписанного в листинге алгоритма \ref{Grad_Model}. Напомним (см. упражнение \ref{Exer_Grad_Const_Step}), что для неускоренного градиентного метода с постоянным шагом \eqref{p0_a10} можно получить оценку
\begin{equation}\label{Nest_Rel_Grad}
F(x_N) - F(x) \leq \frac{LV(x,x_0)}{N}
\end{equation}
для всякого $x$ из допустимого множества задачи, где $x_0$ --- начальная точка. Напомним также, что для рассматриваемой нами функции $F$ из \eqref{p0_a23} при выборе прокс-функции $d$ из \eqref{proxoptdesigh} можно положить $L = 1$ (см. \eqref{Leq1}).  

\begin{teo}
Пусть для задачи \eqref{p0_a23} применяется градиентный метод \eqref{p0_a10} с начальной точкой $x_0 = {\bf e}/n = (1/n, 1/n, ..., 1/n)$ и прокс-функцией \eqref{proxoptdesigh}. Если при $\varepsilon\leq F(x_0)-F^*$ выбрать количество итераций алгоритма \eqref{p0_a10} согласно условию
\begin{equation}\label{EqNestDoptN}
N = \left\lceil\frac{2n\log \left(\frac{2(F(x_0)-F^*)}{\varepsilon}\right)}{\varepsilon}\right\rceil,
\end{equation}
то справедливо неравенство $F(x_N)-F^* \leq\varepsilon$.
\end{teo}
\begin{proof}
Введём вспомогательную величину $\Delta=\frac{\varepsilon}{2(F(x_0)-F^*)}$. Тогда $\Delta \leq 1/2$, поскольку $\varepsilon\leq F(x_0)-F^*$. Пусть $\hat{x}:=(1-\Delta)x^*+\Delta x_0$. Выпуклость функции $F(\cdot)$ означает, что 
$$
F(\hat{x})\leq(1-\Delta)F^*+\Delta F(x_0),
$$
откуда
\begin{equation}\label{EqNest40}
F(\hat{x})-F^*\leq \Delta(F(x_0)-F^*).
\end{equation}

Далее, ввиду \eqref{proxoptdesigh} с учётом выбора начальной точки $x_0$ имеем
$$
V(\hat{x},x_0) = d(\hat{x}) - d(x_0)-\langle \nabla d(x_0),\hat{x}-x_0\rangle = d(\hat{x})-d(x_0)\leq
$$
\begin{equation}\label{EqNest41}
\leq -n \log\left(\frac{\Delta}{n}\right)+n\log\left(\frac{1}{n}\right) = n\log\left(\frac{1}{\Delta}\right), \end{equation}
поскольку $\nabla d(x_0)=-n {\bf e} = (-n, -n, ..., -n)$ и 
$\langle \nabla d(x_0),\hat{x}- x_0\rangle = 0$, а также $\hat{x}\geq(\Delta / n) {\bf e}$. Следовательно, после $N$ итераций алгоритма \eqref{p0_a10} из условия теоремы получаем
$$
F(x_N)-F^* = F(x_N) - F(\hat{x}) + F(\hat{x}) - F^* 
\leq\frac{V(\hat{x},x_0)}{N}+\Delta(F(x_0)-F^*) \leq
$$
\begin{equation}\label{EqNest42}
\leq\frac{n\log \left(1/\Delta \right)}{N}+\frac{\varepsilon}{2}
\leq\varepsilon,
\end{equation}
где первое неравенство в \eqref{EqNest42} следует из \eqref{Nest_Rel_Grad} при $x=\hat{x}$ и неравенства \eqref{EqNest40}, второе неравенство следует из \eqref{EqNest41} и выбора вспомогательной величины $\Delta$, а третье неравенство --- из \eqref{EqNestDoptN}. Цепочка неравенств \eqref{EqNest42} завершает доказательство теоремы.
\end{proof}

Концепция относительной гладкости позволяет расширить класс задач, к которым применимы неускоренные градиентные методы с сохранением оценки скорости сходимости по функции вида $O(1/N)$ ($N$ --- число итераций). Однако некоторая расплата за общность --- необходимость решать подзадачи типа ~\eqref{p0_a11}, которые в случае $Q = \mathbb{R}^n$, или когда $\lim_{x \to \partial Q} d(x) = +\infty$, сводятся к нелинейным уравнениям относительно $x$ вида  $$\nabla d(x) = \nabla d(x_k) - \frac{1}{L} \nabla f(x_k).$$ В общем случае такие задачи не могут быть решены точно. В этом контексте интересен приведённый ранее в пункте \ref{sect_inexact_sol} подход, который позволяет описывать влияние на оценки скорости сходимости градиентных методов погрешностей решения вспомогательных подзадач на их итерациях. 

Однако для прокс-структуры \eqref{proxoptdesigh}, соответствующей задаче $D$-оп\-ти\-маль\-ного плана эксперимента, задача \eqref{p0_a11} может быть записана в виде
\[ \min_{x\in\Delta_{n}}\,\left( \langle c,x\rangle-\sum^{n}_{j=1}\log\,x_{j} \right)
\]
и сводится к одномерной задаче. Условия оптимальности первого порядка принимают вид 
\[ x_1, x_2, ..., x_n >0,\ \langle {\bf e},x\rangle = x_1 + x_2 +... + x_n = 1,\ c-X^{-1}{\bf e}=-\theta {\bf e}
\]
для некоторого скалярного множителя $\theta$. Очевидно, что для всех $j = 1,\dots,n$ имеем $x_{j}=1/(c_{j}+\theta)$, где $c = (c_1, c_2, ..., c_n )$, и
остаётся лишь определить значение параметра $\theta$. Обратим внимание, что $\theta$ должно удовлетворять условию
\begin{equation}\label{p0_a24}
\varphi(\theta) := \sum^{n}_{j=1}\frac{1}{c_{j}+\theta}-1=0
\end{equation}
и является элементом открытого интервала $\mathcal{F}:=\left(-\min_{j}\{c_{j}\},\infty\right)$. Ясно, что $\varphi$ строго убывает на $\mathcal{F}$ и
$\varphi(\theta)\rightarrow +\infty$ при $\theta\searrow-\min_{j}\{c_{j}\}$, $\varphi(\theta)\rightarrow -1$ при $\theta\rightarrow +\infty$. Поэтому уравнение~\eqref{p0_a24} имеет единственное решение на интервале $\mathcal{F}$. Для нахождения этого решения можно применить, например, метод Ньютона.

\begin{table}[h!]
\centering
\caption{{\hspace{1.1cm}\bf Сравнение вычислительных гарантий для задачи \eqref{p0_a23} метода\\ Хачияна и градиентного метода вида \eqref{p0_a11}. Все константы в оценках опущены}}
\label{tabDoptDesign}
{\fontsize{8}{8.5}\selectfont
\begin{tabular}{|c|l|l|l|}
\hline
 Метод &
  \multicolumn{1}{c|}{\begin{tabular}[c]{@{}c@{}}Оценки \nk{количества}\\ итераций метода\end{tabular}} &
  \multicolumn{1}{c|}{\begin{tabular}[c]{@{}c@{}}Сложность\\ итерации\end{tabular}} &
  \multicolumn{1}{c|}{\begin{tabular}[c]{@{}c@{}}Общая сложность\\ подхода\end{tabular}} \\ \hline
Метод Хачияна &
  \begin{tabular}[c]{@{}l@{}}$m\log (F(x_0)-F^{*})\,+$ \\ $+ m^2 \varepsilon^{-1}$\end{tabular} &
  $mn$ &
  \begin{tabular}[c]{@{}l@{}}$m^2 n\log(F(x_0)-F^{*})+$\\ $+ \, m^3 n \varepsilon^{-1}$\end{tabular} \\ \hline
\begin{tabular}[c]{@{}c@{}}Градиентный\\ метод с \\ \nk{относительной}\\ гладкостью\end{tabular} &
  \begin{tabular}[c]{@{}l@{}}$n\;\varepsilon^{-1} \log(F(x_0)-F^{*})\, +$\\ $ + \, n \varepsilon^{-1}\log\varepsilon^{-1}$\end{tabular} &
  $m^2 n$ &
  \begin{tabular}[c]{@{}l@{}}$ m^2 n^2 \varepsilon^{-1} \log(F(x_0)-F^{*}) + $\\ $+ \, m^2 n^2 \varepsilon^{-1}\log\varepsilon^{-1}$\end{tabular} \\ \hline
\end{tabular}

}
\end{table}

В табл.~\ref{tabDoptDesign} приведено сравнение известных результатов об оценках сложности задачи $D$-оптимального плана эксперимента для градиентного метода с относительной гладкостью \cite{LuNesterov} и метода типа условного градиента (метод Хачияна) \cite{KhachiyanEllipsoid}. В обоих методах используется начальная точка $x_0 = (1/n, 1/n, ..., 1/n)$.

Как видим из табл.~\ref{tabDoptDesign}, хотя количество итераций градиентного метода с использованием относительной гладкости может быть меньшим по сравнению с методом Хачияна (некоторый аналог метода условного градиента), последний всё равно выигрывает за счёт менее дорогой итерации. Однако стоит отметить, что такие оценки для метода Хачияна \cite{KhachiyanEllipsoid} получены за счёт специфики конкретной постановки задачи \eqref{p0_a23}. Метод Хачияна \cite{KhachiyanEllipsoid} и подходы других смежных работ заточены на геометрию задачи об эллипсоиде минимального объёма ($D$-оптимального плана эксперимента) и не являются частью общей теории методов условного градиента. В~то~же время для градиентного метода с относительной гладкостью известны доказанные результаты \cite{LuNesterov} об оценках скорости сходимости, которые оптимальны на классе выпуклых относительно гладких задач \cite{DragomirRelativeSmooth}. Кроме того, гладкость задачи \eqref{p0_a23} относительно подходящей прокс-функции указывает на применимость всех подходов и результатов п.~\ref{gradMethod}, т.е. возможность адаптивной настройки в ходе работы метода на константу (относительной) гладкости $L > 0$ с возможностями использования на итерациях неточной информации (значений целевой функции и/или градиента). Применение неускоренных методов гарантирует отсутствие накопления в оценках сложности таких погрешностей. Похожих результатов для методов условного градиента не известно. Отметим, однако, в этом плане описанную нами в~п.~\ref{model_uslovn} связь метода условного градиента с ускоренными градиентными методами. По-видимому, это означает, что возможно выписать оценки для метода условного градиента в модельной общности, с учётом неточностей используемой информации. Но при этом, как и для ускоренных градиентных методов, возможно накопление в таких оценках параметра, соответствующего погрешностям (см.~п.~\ref{subsec_inexact}).

Из табл.~\ref{tabDoptDesign} видно, что при близких значениях $m$ и $n$ оценки сложности для метода Хачияна и градиентного метода с относительной гладкостью сопоставимы при не слишком малом значении $\varepsilon$. При этом проведённые Д.\,А.~Пасечнюком эксперименты в \cite[Appendix F]{ModelOMSpap} показали, что применение адаптивного варианта градиентного метода (алгоритм  \ref{Grad_Model}) с относительной гладкостью может позволить существенно улучшить качество выдаваемого методом решения.


\section{Вычислительные аспекты}
%
В этом разделе собраны сюжеты, не имеющие прямого отношения к выпуклой оптимизации и к оптимизации в целом, однако часто и активно использующиеся специалистами по численным методам оптимизации в своей работе. Было решено собрать такие сюжеты вместе и привести их в конце пособия.

\subsection{Автоматическое дифференцирование}
\label{subsect_autodiff}
Как было продемонстрировано ранее на различных примерах, для решения многих практических задач оптимизации необходимо иметь возможность численно вычислять производные целевых функций (градиенты, гессианы, в отдельных случаях производные высших порядков) в точках. В частности, большая часть современных
алгоритмов машинного обучения основана на обязательном
использовании градиентов для решения задачи оптимизации. В этом разделе будут рассмотрены возможные варианты вычисления производных и подробно разобран наиболее перспективный и часто использующийся в машинном обучении вариант~--- так называемое {\it автоматическое дифференцирование}~(АД)\footnote{В зарубежной литературе АД иногда называют {\it алгоритмическим дифференцированием} или даже просто <<autodiff>>.}
\cite{AD_survey2015}.

Все существующие на данный момент методы вычисления производных для использования в дальнейших вычислениях на компьютере можно разделить на четыре группы:
\begin{enumerate}
    \item  Ручной аналитический вывод формулы для производной и
    последующее использование формулы в программе.
    \item Численное дифференцирование с помощью метода конечных разностей. В самой простой форме метода используется
    определение производной через предел. В одном
    из таких вариантов частную производную   
    $\frac{\partial f}{\partial x_i}$
    функции
    $f(x) \, : \, \mathbb{R}^n \rightarrow \mathbb{R}$
 можно вычислить по приближённой формуле
 $$
 \frac{\partial f}{\partial x_i} \approx 
 \frac{f(x + h e_i) - f(x - h e_i)}{2h},
 $$ 
 где $e_i$~--- $i$-й единичный вектор и $h \, > \, 0$~ --- шаг.
    \item Символьное дифференцирование. Реализовано в системах компьютерной алгебры, таких, как Maple, Mathematica, Maxima.
    \item Автоматическое дифференцирование.
\end{enumerate}
Главным недостатком первого способа вычисления производных является необходимость иметь аналитическую формулу целевой функции. Во многих приложениях, включая нейронные сети, это невозможно. Что касается численного дифференцирования,
то у этого способа есть несколько существенных недостатков: невысокая точность и численная неустойчивость из-за ошибок округления и аппроксимации, необходимость тщательно выбирать шаг~$h$, а также $2n$  обращений к оракулу целевой функции для вычисления градиента в $\mathbb{R}^n$.

Далее рассмотрим подробнее последний из перечисленных способов вычисления производных, который магическим\footnote{На самом деле никакой магии, конечно, а обычная математика.}
образом лишён главных недостатков предыдущих способов,
но сохраняет их достоинства.

\begin{defin}
Автоматическое дифференцирование~--- это общее название
техник вычисления аналитических значений производных функции с использованием компьютерного представления
данной функции.
\end{defin}

\begin{figure}[htb]
\begin{center}
\includegraphics[width=0.9\linewidth]{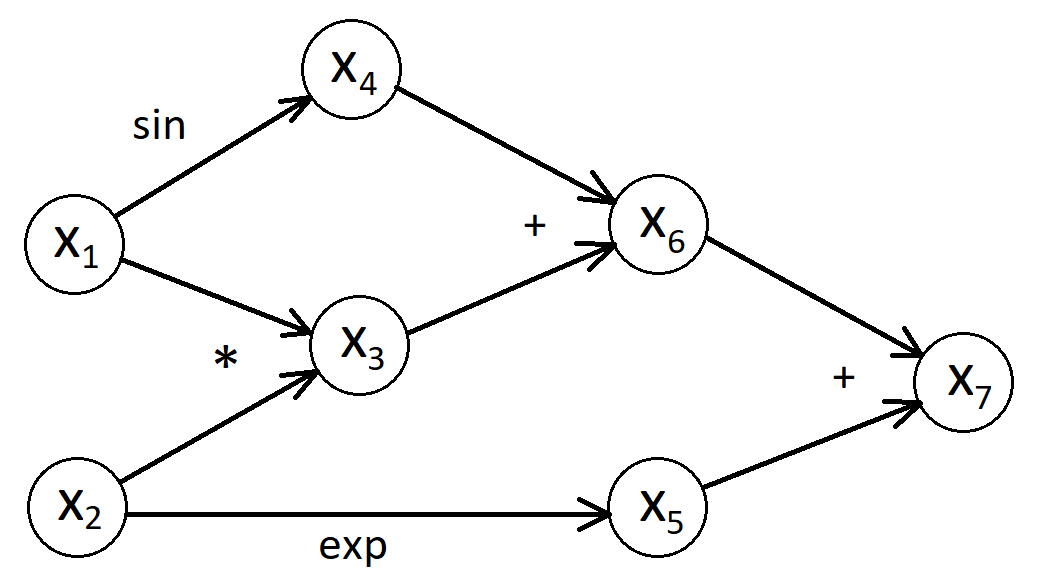}
\end{center}
\caption{Вычислительный граф функции $f(x_1, \, x_2) = x_1 x_2 + \sin x_1 + e^{x_2}$} 
\label{ch3_dr_graph_autodiff_1}
\end{figure}

Возможность создания методов~АД основана на простом наблюдении: вычисление значения любой сколь угодно сложной функции от аргумента можно
представить как иерархическую последовательность
элементарных в текущей модели операций над одним или двумя аргументами на текущем уровне иерархии.
Такое представление в виде композиции называется {\it вычислительным графом}\label{page_comp_graph}, или {\it графом вычислений}
функции (см. пример вычислительного графа на рис.~\ref{ch3_dr_graph_autodiff_1}). При этом вершины графа~--- это элементарные функции, а рёбра служат для передачи в  эти функции их аргументов.
В классической модели~АД элементарные бинарные операции 
включают сложение, умножение, деление и возведение в степень~$x^y$. Унарными операциями являются, например,
тригонометрические функции. 
Следующим общим компонентом всех методов~АД 
является использование цепного правила дифференцирования
(chain rule). Согласно этому правилу, если
$f$ является функцией вектора $g = (g_1, \, \ldots, \, g_m) \, \in \, \mathbb{R}^m$,
который, в свою очередь, зависит от вектора $x \, \in \, \mathbb{R}^n$, то производную $f$ по $k$-й компоненте~$x$ можно определить
по следующей формуле:
$$
\frac{ \partial f(g(x)) }{ \partial x_k}= \sum_{i = 1}^m \frac{\partial f}{\partial g_i}  \frac{ \partial g_i (x)}{ \partial x_k}  
= \left\langle \nabla f, \,
\frac{\partial g}{ \partial x_k} \right\rangle.
$$

В АД есть два основных режима: прямой и обратный. Проиллюстрируем разницу между этими режимами на простом примере.

\begin{example}
Рассмотрим функцию двух переменных
$$
f(x_1, \, x_2) = x_1 x_2 + \sin x_1 + e^{x_2}.
$$
На рис.~\ref{ch3_dr_graph_autodiff_1} показано, что процесс вычисления данной функции может быть представлен
как последовательность элементарных операций, выполняемых в определённом порядке. Например, мы не можем выполнить сложение $\sin x_1$ и $x_1 x_2$ до тех пор, пока не посчитаем результат произведения $x_1 x_2$. 

В прямом режиме вычисление градиента происходит параллельно с вычислением самой функции. Ассоциируем с каждой вершиной графа на рис.~\ref{ch3_dr_graph_autodiff_1} некоторую векторно-значную переменную $v_i$, соответствующую градиенту $\nabla x_i$ переменной, относящейся к данной вершине. Одновременно с независимыми переменными $x_1,x_2$ инициализируем векторы $v_1,v_2$, соответствующие этим переменным, единичными векторами $v_1 = e_1$, $v_2 = e_2$, в соответствии с соотношением $\frac{\partial x_i}{\partial x_j} = \delta_{ij}$, $i,j = 1,2$, где $\delta_{ij}$ --- символ Кронекера. Далее на каждом шаге вектор $v_i$, соответствующий вычисляемой текущей переменной $x_i$, считается в соответствии с правилом дифференцирования, соответствующим текущей вершине, т.е.,
\[ v_3 = v_1x_2 + x_1v_2,\ v_4 = v_1 \cdot \cos x_1,\ v_5 = v_2 \cdot x_5,\ v_6 = v_3 + v_4,\ v_7 = v_5 + v_6.
\]
Таким образом, градиент $v_i$ считается параллельно со значениями $x_i$. Отметим, что гораздо дешевле, чем весь градиент, считать производную по направлению. Для этого векторно-значная переменная $v_i$ заменяется скалярной и инициализируется элементами вектора направления. Прямой режим автоматического дифференцирования можно осуществить, например, перегрузкой функций или процедур, вычисляющих значения промежуточных переменных $x_i$. 

В обратном режиме вычисление \emph{всего} градиента скалярной функции является таким же затратным, как вычисление одной производной по направлению в прямом режиме. В обратном режиме каждой вершине сопоставляется скаляр $v_i$, соответствующий коэффициенту чувствительности выходной функции $f$ к изменениям значения переменной $x_i$, т.е., частной производной $\frac{\partial f}{\partial x_i}$ (хотя это обозначение не совсем корректно, поскольку промежуточные переменные $x_i$ не являются независимыми). Однако вычисление значений $v_i$ происходит не одновременно с вычислением $x_i$, а во время отдельной второй фазы, которая запускается уже после того, как все $x_i$ вычислены. Эта фаза начинается инициализацией скаляра $v_7$, соответствующего выходной вершине, единицей, $v_7 = 1$. Затем остальные значения $v_i$ заполняются от вершины к вершине в \emph{обратном} направлении в соответствии с правилом дифференцирования, соответствующим текущей вершине, т.е.,
\[ v_6 = v_7,\ v_5 = v_7,\ v_4 = v_6,\ v_3 = v_6,\ v_2 = v_3 \cdot x_1 + v_5 \cdot x_5,\ v_1 = v_3 \cdot x_2 + v_4 \cdot \cos x_1.
\]
Объясним этот принцип на примере вычисления значения $v_2$. Значение переменной $x_2$ двояким образом влияет на конечный результат $f = x_7$, а именно посредством зависимости переменных $x_3$ и $x_5$ от $x_2$. Так как эти эффекты складываются в первом порядке малости, значение $v_2$ вычисляется как сумма двух соответствующих членов. Чувствительность значения $f$ к изменениям значений $x_3,x_5$ определяется уже посчитанными значениями $v_3,v_5$. В итоге получаем
\[ v_2 = v_3\frac{\partial x_3}{\partial x_2} + v_5\frac{\partial x_5}{\partial x_2} = v_3 \cdot x_1 + v_5 \cdot x_5.
\]
Заметим, что частные производные $\frac{\partial x_3}{\partial x_2}$, $\frac{\partial x_5}{\partial x_2}$ в этом примере ассоциируются с рёбрами вычислительного графа. Эти производные можно уже вычислить в прямом проходе параллельно со значениями переменных $x_i$, перед проходом в обратном режиме, вычисляющим значения $v_i$.
\end{example}

Вычисление градиента скалярной функции в обратном режиме требует выполнения числа операций в 4--5 раз большего, чем для самой скалярной функции \cite{Nocedal}. Однако, в отличие от прямого режима, требуется сохранять значения $x_i$ всего вычислительного графа. Если выходных переменных несколько, то для вычисления градиента каждой из них необходим одни проход по графу в обратном режиме. Таким образом, один проход в прямом режиме считает один столбец матрицы Якоби частных производных, в то время как \nk{один} проход в обратном режиме считает одну строку матрицы Якоби. Если необходимо вычислить все частные производные, то в первом приближении прямой режим менее затратный в том случае, когда больше выходных переменных, чем входных, а обратный режим менее затратный в случае, когда перевешивает число входных переменных.

В ситуации, когда нужно минимизировать некоторую функцию невязки от выходных переменных, например, в обучении нейронных сетей, можно ограничиться одним проходом в обратном режиме, инициализировав скаляры $v_i$, соответствующие выходным переменным, частными производными функции невязки по этим переменным.

\begin{remark}
В глубоком обучении (нейронных сетей) обратный режим АД часто называют методом обратного распространения ошибки (backpropagation). Именно АД является одной из ключевых технологий,
которые определили успех <<революции>> глубокого обучения.
\end{remark}

\begin{remark}
На данный момент существует несколько удобных программных реализаций
для создания вычислительного графа и выполнения операций с ним. Эти
реализации можно поделить на два типа (подробнее см. п.~\ref{ch1_sect_opt_package}): 
\begin{enumerate}
    \item Вычислительный граф полностью строится заранее, до выполнения
    программы. Пример пакета: Tensorflow.
    \item Вычислительный граф строится динамически, во время выполнения 
программы (JIT, компиляция ''just-in-time''). Примеры пакетов: Tensorflow Eager, JAX, PyTorch.\qedhere
\end{enumerate}
\end{remark}

\subsection{Оптимизационные пакеты} \label{ch1_sect_opt_package}
Понимать, что стоящая перед вами задача оптимизации является выпуклой, а также уметь сводить задачи оптимизации к выпуклым задачам (в некоторых случаях это удаётся), очень важно.
После этого для численного решения задачи выпуклого программирования можно воспользоваться существующими оптимизационными пакетами. Однако в настоящее время размерность задач, которые можно решить с помощью таких пакетов, ограничена примерно $10^5$--$10^6$ переменными.

Очень часто, для того, чтобы можно было воспользоваться возможностями современных пакетов, выпуклая оптимизационная задача должна быть хорошо структурирована, 
то есть переформулирована в виде конической программы над
симметричными конусами одного из <<популярных>> семейств (положительный ортант,
конусы Лоренца, матричные конусы). Конические программы над данными конусами
эквивалентны, соответственно, задачам линейного программирования, конического программирования второго порядка и полуопределённого программирования 
(см.~подробнее гл.~1, раздел~\ref{ch1_sect_cone_prog}). 
В свою очередь, все основные приложения выпуклой оптимизации в
инженерии, статистике, принятии решений и других областях могут
быть сведены к решению задач одного из таких типов.

\subsubsection{Система моделирования CVX}
Наиболее известным и популярным из оптимизационных пакетов для
выпуклого программирования является пакет {\bf CVX} (\url{http://cvxr.com/cvx/}), являющийся расширением Matlab. После установки данного расширения в среде Matlab можно использовать специальный язык моделирования
для решения задач выпуклой оптимизации. При этом целевая
функция и ограничения задачи указываются с помощью стандартных конструкций Matlab.
Обычный код на Matlab можно легко 
комбинировать с кодом для пакета CVX.
Это позволяет выполнить в Matlab вычисления, необходимые для формулировки оптимизационной
проблемы, и сразу же решить её с помощью CVX,
и затем обработать полученные результаты также
в Matlab.

\begin{example}
Для следующей задачи
\begin{eqnarray*}
& \min & \| Ax - b \|_2 \\
& Cx = d &\\
& \| x \|_\infty \, \le \, e &
\end{eqnarray*}
код на CVX, с помощью которого можно получить решение частного случая этой задачи в Matlab, выглядит так:

\begin{lstlisting}[language=Matlab]
m = 10; n = 5; p = 4;
A = randn(m, n); b = randn(m, 1);
C = randn(p, n); d = randn(p, 1); 
e = rand;
cvx_begin 
    variable x(n)
    minimize( norm(A * x - b, 2) )
    subject to
        C * x == d
        norm(x, Inf) <= e
cvx_end
\end{lstlisting}
\end{example}

В пакете CVX реализован особый подход к решению задач выпуклого программирования, который авторы пакета Стивен Бойд и Майкл Грант назвали 
<<управляемое выпуклое программирование>>
(disciplined convex programming, DCP) 
\cite{DCP_06}.
Управляемое выпуклое программирование~--- это методология, предназначенная для
конструирования выпуклых задач по особому
набору правил\footnote{См. также \cite{JuditskyNemirovski21, BenTalNemirovski01, Ben-Tal}. В \cite{JuditskyNemirovski21} указано, что DCP по сути выполняется с помощью специального типа <<анализа>>, который состоит из двух частей:
1) <<строительных материалов>>, <<кирпичиков>>~--- набора
множеств и функций, которые допускают явное представление через
упомянутые выше семейства симметричных конусов;
2) набора операций (<<правил>>) над этими коническими представлениями, которые
сохраняют выпуклость (например, неотрицательные линейные комбинации выпуклых функций,
конечное пересечение выпуклых множеств,
см. подробнее пп.~\ref{ch1_sect_conv}, \ref{ch1_sect_preserve_conv_func}).} (\url{http://cvxr.com/dcp/}). 
Этот набор правил представляет собой достаточные,
но не необходимые условия выпуклости целевой функции
и ограничений. Выпуклые функции и множества
в концепции DCP получаются путём применения
к стартовой библиотеке функций и множеств,
про которые точно известно, что они выпуклые,
некоторого набора преобразований,
сохраняющих выпуклость.
Полная алгоритмичность этого набора преобразований позволяет
легко реализовать эту концепцию на языке программирования,
что и было крайне успешно выполнено в пакете~CVX.

\begin{figure}[htb]
\begin{center}
\includegraphics[width=0.78\linewidth]{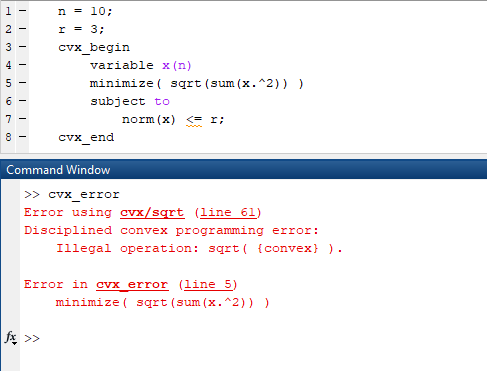}
\end{center}
\caption{Пример неправильной формулировки задачи для пакета CVX} 
\label{ch1_dr2_cvx_error}
\end{figure}

Если пользователь пакета CVX
знает, что его задача выпукла, и может
представить её согласно набору правил DCP,
то задача легко и быстро решается в пакете CVX.
Следует отметить, что если выпуклая задача сформулирована
с нарушением правил DCP, то пакет не сможет её решить в такой форме. Но если переформулировать задачу, то решение будет найдено. Так, например,
если стоит задача минимизировать 2-норму $\| x \|_2$
при каких-то выпуклых ограничениях, которая, как известно, является выпуклой функцией, но для
CVX задача сформулирована как \mcode{minimize sqrt(sum(x.^2))}, то пакет выдаст ошибку, так как 
выпуклость данной функции не следует из набора правил DCP (в DCP считается, что функция $\sqrt{x}$ вогнута, см. рис.~\ref{ch1_dr2_cvx_error}).

Следует отметить, что в пакете CVX
по умолчанию используется некоммерческий
оптимизационный солвер {\bf SDPT3}~\cite{SDPT3},
но можно подключить и ряд других солверов:
{\bf SeDuMi}~\cite{Sedumi} (\url{http://sedumi.ie.lehigh.edu/},
наряду с SDPT3 входит в стандартный дистрибутив CVX), {\bf GLPK} (\url{http://www.gnu.org/software/glpk/},
экспериментальная поддержка, только для задач линейного и целочисленного программирования), а
также коммерческие {\bf MOSEK} и {\bf Gurobi}
(см. раздел~\ref{subsec_other_solvers}). Для
использования последних требуется
также и профессиональная версия (лицензия) CVX.
Такое разнообразие объясняется тем,
что каждый солвер имеет свой ограниченный набор возможностей (в частности, поддержку разных типов ограничений) и производительности.
Эксперименты авторов пакета CVX показали\footnote{\url{http://web.cvxr.com/cvx/doc/solver.html}},
что для большинства проблем SeDuMi находит решение быстрее,
чем SDPT3, но надёжность у последнего выше, именно поэтому он выбран солвером по умолчанию.

Оболочки, предоставляющие
возможность использовать CVX-подобный
синтаксис представления задач выпуклого
программирования, в
настоящее время (на конец 2020~года) выпущены
для следующих языков программирования (помимо MATLAB):
\begin{enumerate}
    \item Python, пакет {\bf CVXPY} (\url{https://www.cvxpy.org/})\footnote{Для CVXPY в 2019 году появилась также
библиотека {\bf cvxpylayers} (\url{https://github.com/cvxgrp/cvxpylayers})
с возможностью 
включать в вычислительные графы (в том числе искусственные
нейронные сети), создаваемые в PyTorch, JAX и
TensorFlow,
параметризованные задачи выпуклого программирования (т.е. объекты CVXPY), с целью далее использовать
для этих объектов автоматическое
дифференцирование.};
\item R, пакет  {\bf CVXR}~\cite{CVXR};
\item  Julia, пакет {\bf Convex.jl}~\cite{Convex.jl}, \url{https://jump.dev/Convex.jl/stable/}. 
Является частью экосистемы
{\bf JuMP}, набора пакетов для решения
задач математического программирования 
на языке Julia: \url{https://jump.dev/}.
\end{enumerate}

\subsubsection{Оптимизация для нейронных сетей}

Как известно, целью обучения практически любой модели машинного обучения, в том числе и модели с нейронной сетью~\cite{good, NikKad_DL}, является уменьшение заданной функции ошибки (функции потерь). Следовательно, и здесь возникает задача оптимизации,
для решения которой
необходимы эффективные
методы. В случае нейронных сетей~(HC) возникающая оптимизационная задача
обычно является задачей без ограничений (оптимальные веса нейронов чаще всего могут
быть произвольными вещественными числами). Но это всё равно не делает
задачу обучения НС проще, поскольку,
во-первых, возникающая целевая функция потерь является невыпуклой,
а, во-вторых,
даже не самая сложная архитектура НС может содержать миллионы параметров (весов), что превращает задачу нахождения
оптимальных параметров в критично требовательную к вычислительным ресурсам. Тем не менее в настоящее время для обучения НС уже создан ряд удобных библиотек,
которые позволяют за несколько строк кода создать 
и натренировать (обучить) нейронную сеть с любой современной архитектурой. К тому же
за одну строку кода можно добавить возможность тренировки НС на видеокарте~(GPU),
что значительно ускоряет этот процесс (учитывая сверхбольшую размерность задачи!).

Два наиболее популярных фреймворка для глубокого обучения (так
называют сейчас все эти библиотеки из-за основного применения)
были созданы конкурирующими корпорациями-гигантами. {\bf TensorFlow} был разработан во внутреннем подразделении ``Google Brain''
компании Google и активно
используется внутри Google
в их собственных продуктах,
работающих с 
применением возможностей
искусственного интеллекта (Gmail, Photos, поиск,\,\dots). 
{\bf PyTorch} был создан в AI Research lab компании Facebook.
Рассмотрим их немного подробнее.

TensorFlow~--- 
платформа
с открытым исходным кодом 
и экосистема инструментов,
библиотек и ресурсов,
предназначенная для решения задач машинного обучения. 
Cвоё название получила из-за основного примитива
платформы~--- многомерного массива, или {\it тензора}. 

\begin{defin}
Формально, тензором называется линейное
по каждому из аргументов отображение 
$$
f \, : \, \underbrace{V^* \times \ldots \times V^*}_{n} \times \underbrace{V \times \ldots V}_{m} \rightarrow \mathbb{R},
$$
где $V$~--- векторное пространство, $V^*$~--- сопряжённое пространство.
\end{defin}
Скаляр является тензором ($f\, : \, \mathbb{R} \rightarrow \mathbb{R}$), $n$-мерный вектор можно
представить как тензор $f\, : \, \mathbb{R}^n \rightarrow \mathbb{R}$, матрицу размера $n \times m$~--- как тензор $f\, : \, \mathbb{R}^n \times \mathbb{R}^m \rightarrow \mathbb{R}$,
и так далее. Именно поэтому тензоры
чаще всего представляются как многомерные массивы чисел.

Тензоры участвуют во всех вычислениях.
В TensorFlow определены функции над тензорами, производные
которых вычисляются автоматически
(здесь тоже работает {\it автоматическое дифференцирование}). Эти возможности
работать с функциями над тензорами и получать
их производные любых порядков~--- ключевое отличие обычных библиотек для работы с массивами
типа {\bf NumPy} от фреймворков глубокого обучения.

Все операции должны быть представлены
в виде {вычислительного графа}\footnote{См. стр.~\pageref{page_comp_graph}.},
в первых версиях TensorFlow статического
(без изменения состояния).
Тензор имеет начальное состояние, <<прогоняется>> через вычислительный граф операций,
и выходит с обновлёнными значениями. Таким образом, вершинами графа являются математические
операции, а в качестве рёбер можно
рассматривать тензоры с данными. Модель вычислений в виде графа оказывается очень
удобной по нескольким причинам: 
можно легко распределить вычислительную работу
на несколько процессоров; можно сохранить граф
вычислений для использования в другой
момент времени.

Одной из полезных особенностей TensorFlow 
является набор инструментов TensorBoard,
предназначенный для
визуализации экспериментов,
структуры нейронной сети,
процесса вычисления в TensorFlow,
и в том числе графического представления
вычислительного графа. Также Google
предоставляет удобную среду {\bf Colab notebook}
(\url{https://colab.research.google.com}) для программирования прямо в браузере,
без необходимости устанавливать что-либо на свой компьютер и с поддержкой всех возможностей TensorFlow.

Первый релиз TensorFlow
состоялся в ноябре 2015~года.
Обновлённая версия, которая сейчас называется TensorFlow~2.0, была выпущена в сентябре
2019~года. Примерно в то же
время вышел обзор <<The State of Machine Learning Frameworks in 2019>>~\cite{He},
подтвердивший тенденцию: всего за год
бывший экзотикой в 2018-м году PyTorch
стали использовать для экспериментов
в абсолютном большинстве академических статей,
в том числе и в статьях, прошедших на ICLR и ICML.
Самым значительным изменением в TensorFlow~2.0
стало появление ``Eager execution'' (<<Немедленного выполнения>>, Define-by-Run), когда при выполнении
вычислений есть возможность получения
немедленного результата
в динамическом вычислительном графе
вместо построения
статического вычислительного графа. Такой
режим выполнения кода обеспечивает
интуитивно-понятный интерфейс с использованием
исходных типов данных Python и его
управляющую логику,
а также расширенные возможности отладки. Непосредственно для нейросетей динамический вычислительный граф даёт возможность менять параметры нейросети во время тренировки.
В то же время режим ''Eager execution'' не подходит для внедрения и при больших размерах данных будет
работать медленнее.

Часто вместе с TensorFlow используется и библиотека
{\bf Keras} (\url{keras.io})~--- высокоуровневое API на языке Python для
вызова
реализаций 
композиции различных часто используемых
строительных блоков нейронных сетей (слоёв
всех существующих типов, функций активации,
оптимизационных алгоритмов).
До версии~2.3 Keras поддерживал и другие
фреймворки глубокого обучения
для построения графа вычислений
и проведения автоматического дифференцирования
(в том числе и долгое время занимавшую
лидирующие позиции среди таких фреймворков
библиотеку {\bf Theano}~\cite{Theano},
созданную в университете Монреаля, в группе
одного из ведущих мировых экспертов
в области искусственного интеллекта и первооткрывателей глубокого обучения Йошуа Бенджио (Yoshua Bengio)).
Начиная с версии~2.4, Keras поддерживает
только TensorFlow~2.0.

PyTorch~--- библиотека для машинного
обучения с открытым исходным
кодом, созданная в Facebook
на основе написанной на языке Lua библиотеки
Torch. Первая версия вышла осенью 2016~года. Несмотря на название,
в дополнение к основному интерфейсу на Python,
PyTorch имеет также ограниченную поддержку языков С++ и Java. 

Так же, как и TensorFlow, PyTorch оперирует с {\it тензорами} и рассматривает любую модель как ориентированный ациклический граф. Для получения
производных с целью последующей оптимизации
используется автоматическое дифференцирование.
Можно использовать все популярные библиотеки Python
вместе с PyTorch, в том числе {\bf Matplotlib} и {\bf Seaborn} для визуализации. 
Стоит упомянуть также {\bf PyTorch Lightning},
появившуюся в середине 2019~года библиотеку на Python с открытым исходным
кодом, которая создана для учёных, работающих в области искусственного интеллекта. PyTorch Lightning преобразует код на PyTorch таким
образом, чтобы отделить науку от инженерии,
то есть сделать код экспериментов с нейронными
сетями легкочитаемым и масштабируемым, а 
сами эксперименты~--- повторно воспроизводимыми,
в том числе и распределённо. Lightning имеет абстракции в коде
для превращения цикла обучения
в объектно-ориентированную структуру, с заменой
сложных, повторяющихся в любом
коде обучения НС, конструкций на стандартные
макросы. В конце
2019~года конференция NeurIPS 
рекомендовала использование PyTorch Lightning
в качестве стандарта для представления
кода к статьям, поданным на эту конференцию.

\begin{exercise}
Проведите экспериментальное сравнение 
представленных в Pytorch оптимизационных алгоритмов
на какой-либо задаче выпуклого программирования
из представленных в пособии. 
Пример кода можно найти, например, здесь:
\url{https://github.com/Karina1997/math-optimization-ellipsoids-method}. Обратите внимание,
что в представленных экспериментах
библиотечные методы сравниваются на практике с методом, разработанным автором репозитория. Подобный
подход возможен из-за того, что PyTorch позволяет
подключать сторонние методы в качестве оптимизационных, что удобно при  научных исследованиях.
\end{exercise}

К недостаткам PyTorch
можно отнести отсутствие поддержки статических вычислительных графов\footnote{Что в какой-то мере компенсируется возможной оптимизацией модели с помощью TorchScript.}. PyTorch поддерживает
только динамические вычислительные графы, в отличие от TensorFlow~2.0, в котором есть поддержка обоих типов графов.

В качестве оптимизационных алгоритмов в PyTorch
(на конец 2020 года) реализованы
\begin{itemize}
    \item Adadelta (метод с адаптивным подбором скорости обучения)~\cite{ada};
    \item Adagrad (адаптивный субградиентный метод для онлайн-обучения и стохастической оптимизации)~\cite{adagrad};
    \item несколько вариантов метода Adam (метод для стохастической оптимизации)~\cite{adam};
    \item  стохастический градиентный спуск, обычный, усреднённый и с ускорением;
    \item  L-BFGS (квазиньютоновский метод с ограниченной памятью)~\cite{lbfgs};
    \item RMSprop Д.~Хинтона (\url{https://www.cs.toronto.edu/~tijmen/csc321/slides/lecture_slides_lec6.pdf});
    \item Rprop (эвристический алгоритм устойчивого обратного распространения)~\cite{rprop}.
\end{itemize}

Очень важно понимать, что, несмотря на общее устоявшееся название~--- фреймворки для
глубокого обучения,  TensorFlow
и  PyTorch не являются исключительно библиотеками обучения нейронных сетей. Их
можно также использовать для выполнения
численных расчётов на графах,
в том числе и для 
выполнения автоматического дифференцирования
(см. раздел~\ref{subsect_autodiff}).

Перспективной для выполнения АД библиотекой машинного обучения является созданная
в Google Research Python-библиотека {\bf JAX}\footnote{Использование JAX
в DeepMind для исследований в ML: \url{https://deepmind.com/blog/article/using-jax-to-accelerate-our-research}} (\url{https://github.com/google/jax}).
По сути, это обновлённая версия библиотеки Autograd (\url{https://github.com/hips/autograd}),
соединённая с технологией ускоренных векторных вычислений XLA (Accelerated Linear Algebra)
для моделей TensorFlow.
JAX может автоматически дифференцировать функции Python, совместима
с NumPy и SciPy,
и поддерживает прямой и обратный режим дифференцирования.

Стоит упомянуть также {\bf Scikit-learn} ---
самую популярную на сегодняшний день библиотеку на языке Python, предназначенную для решения задач машинного обучения
и анализа данных.
Библиотека также включает несколько реализованных
моделей нейронных сетей. В частности, многослойный перцептрон (\url{https://scikit-learn.org/stable/modules/neural_networks_supervised.html})\footnote{
Многослойный перцептрон (МСП)~--- класс искусственных
нейронных сетей прямого распространения. 
МСП состоит по крайней мере из 3 слоёв узлов: входного, скрытого и выходного слоя. Узлы всех слоёв, за исключением входного, представляют собой модели искусственных нейронов с нелинейной функцией активации.} 
и ограниченную машину Больцмана (\url{https://scikit-learn.org/stable/modules/neural_networks_unsupervised.html#restricted-boltzmann-machines})\footnote{Ограниченная машина Больцмана (ОГМ)~---  генеративная стохастическая нейронная сеть, которая может обучиться определять вероятностное распределение по множеству входных данных. Впервые такая архитектура была получена в 1986 году, но расцвет этих сетей
пришёлся на середину двухтысячных, когда Джеффри Хинтон с коллегами разработали для них быстрые алгоритмы
обучения~\cite{Hinton_RBM}. В настоящее время ОГМ в практических приложениях почти полностью вытеснены генеративно-состязательными нейронными сетями~(GAN) и вариационными автоэнкодерами.}.

\subsubsection{И другие}\label{subsec_other_solvers}
Кроме вышеупомянутого программного обеспечения, следует
обратить внимание также на следующие пакеты:
\begin{enumerate}
\item 
{\bf SciPy.optimize}  ---
библиотека на Python (\url{https://docs.scipy.org/doc/scipy/reference/tutorial/optimize.html}), 
предназначенная для решения
задач минимизации, в том числе и с ограничениями. 
В библиотеку входят солверы для задач 
нелинейного выпуклого программирования,
глобальной оптимизации, линейного программирования,
регрессии, а также реализации метода наименьших квадратов,
и нахождения корней уравнений.
В частности, в библиотеке есть
реализации метода сопряжённых градиентов,
квазиньютоновских методов,
методов доверительных областей и многих других.
Кроме того, в оболочке \texttt{minimize} можно подключить и свой собственный метод
минимизации, реализованный на Python, что может быть
полезно для тестирования и сравнения производительности этого
метода с другими.

    \item IBM ILOG {\bf CPLEX} Optimization Studio (часто название сокращают до CPLEX)~--- 
    корпоративная аналитическая система
    для поддержки принятия решений, включающая
    средства для описания математических моделей,
    приводящих к оптимизационным задачам на
    специальном языке Optimization Programming Language (OPL) и высокопроизводительный пакет
    для решения задач линейного, смешанного целочисленного и квадратичного программирования большой размерности CPLEX (\url{https://www-01.ibm.com/software/commerce/optimization/cplex-optimizer/}). 
\item {\bf MOSEK} --- коммерческий солвер, предназначенный для
решения задач выпуклой оптимизации большой размерности (\url{https://www.mosek.com/}). Имеет интерфейсы для
языков программирования C, C++, Java, MATLAB, .NET, Python и R. Бесплатен для академического использования. Имеет демо-версию (\url{https://solve.mosek.com/web/index.html})
с лимитом в 1000 переменных и
1000 ограничений и временным лимитом на расчёт 20 секунд.
Широко используется для расчётных задач в финансовом секторе, энергетике и лесной промышленности. 

\item {\bf Gurobi}~--- один из наиболее мощных на текущий момент коммерческих солверов задач математического (линейного и квадратичного) программирования. Бесплатен для
академического использования. Есть интерфейсы
к языкам программирования C, C++, Java, Matlab, Python, R. На сайте выложен набор Jupyter-ноутбуков с примерами использования (\url{https://gurobi.github.io/modeling-examples/intro_to_modeling/introduction_to_modeling.html}, \url{https://www.gurobi.com/resource/discover-how-you-can-boost-your-mathematical-optimization-modeling-skills-with-python/}). 

\item {\bf Local Solver}~--- коммерческий солвер (\url{https://www.localsolver.com/home.html}), позиционирующий себя как универсальное решение для
любых бизнес-задач, в которых возможна какая-либо оптимизация, в том числе и невыпуклая. Бесплатен для
академического использования (обновляемая лицензия на один месяц). Есть интерфейсы
к языкам программирования C++, Java, Python. 
\end{enumerate}

\vspace{5mm}
\textbf{Системы алгебраического моделирования}
Для работы со многими из вышеупомянутых солверов необходимо подготовить описание оптимизационной задачи на одном из специальных языков
так называемого {\it алгебраического моделирования}.

Языки алгебраического моделирования занимают
среднюю позицию между математическим описанием оптимизационной задачи и созданием
алгоритма решения такой задачи на
языке программирования. Они позволяют
подготовить описание задачи на высокоуровневом
символьном языке и затем передать его солверу,
тем самым упрощая анализ и проверку модели.
В таких языках используется 
машиночитаемая адаптация традиционных математических выражений, таких как $x_i + y_i$,
$\sum_{i=1}^n a_{ij}$, $x_i \geq 0$, $x \, \in \, Q$ для описания оптимизационных задач как
задач минимизации или максимизации,
с учётом ограничений --- в виде равенств или неравенств. Эти языки сыграли и продолжают играть существенную
роль в развитии применения оптимизации
для таких дисциплин, как исследование операций
и многих других.

В настоящее время наиболее 
распространёнными системами (языками) алгебраического моделирования являются
следующие:
\begin{enumerate}
    \item {\bf AMPL}~\cite{AMPL} (A Modeling Language
for Mathematical Programming)~--- фактически,
стандарт де-факто в мире языков
алгебраического моделирования. Абсолютное большинство
существующих на рынке оптимизационных солверов
поддерживают AMPL.
    
    \item {\bf GAMS} (General Algebraic Modeling System), один из первых таких языков, в настоящее время широко используется в ряде отраслей промышленности. Есть онлайн-сервер
    с возможностью использования GAMS
    в Jupyter-ноутбуках (\url{https://www.gams.com/blog/2020/08/how-to-use-gams-in-jupyter-notebooks/}).
    
    \item {\bf AIMMS} (Advanced Integrated Multidimensional Modeling Software, \url{https://www.aimms.com/})~--- отличается от некоторых других языков тем,
    что имеет графический интерфейс и для разработчика, и для пользователя. Бесплатен
    для академического использования (вместе с набором солверов). 

\item {\bf LINGO}, язык от компании LINDO Systems, 
используется для решения задач линейного, целочисленного и нелинейного программирования
в программном обеспечении от данной компании. 

\item {\bf MPL} (Mathematical Programming Language), язык от компании Maximal Software, есть графический интерфейс и возможность
соединения с базами данных.
Имеется бесплатная версия для 
некоммерческого использования.
\end{enumerate}

Существуют и другие системы оптимизационного моделирования, основанные на представлениях 
данных, отличающихся от языков алгебраического моделирования. Такие альтернативные форматы
включают:
\begin{itemize}
    \item диаграммы в виде блок-схем для линейных ограничений;
\item диаграммы деятельности, в которых
переменные~--- это виды деятельности,
а ограничения~--- это эффекты,
влияющие на эти виды деятельности;
\item графовые или сетевые модели для моделирования потоков.
\end{itemize}
Эти форматы могут иметь некоторые преимущества для специальных классов задач,
но языки алгебраического моделирования являются
более универсальными.

Знание одного из языков алгебраического моделирования даёт возможность воспользоваться
бесплатным интернет-сервисом {\bf NEOS}
(\url{https://neos-server.org/neos/})
от университета Висконсина
для численного решения оптимизационных задач большой размерности и разных типов. На данный момент
NEOS Server предоставляет абсолютно
бесплатный доступ к более чем 60 современным
солверам, в том числе CPLEX, Gurobi,
MOSEK, Knitro, LINDOGlobal и др.

%
%
\subsection{Модель вычислений Блюм--Шуба--Смейла}\label{ch2_subsectBSSm}

Главной нерешённой проблемой в теории линейного программирования
по мнению выдающегося американского математика Стивена Смейла,
лауреата Филдсовской премии 1966~года за 
работы в области дифференциальной топологии\footnote{Премия была присуждена Смейлу за доказательство
обобщённой гипотезы Пуанкаре в пространствах 
размерности больше~$4$.}, является вопрос
о существовании полиномиального по времени алгоритма
{\bf над вещественными числами}, определяющего,
совместима ли линейная система неравенств $Ax \, \ge \, b$. 
\begin{remark}
Данная проблема является девятой по списку
из 18~великих нерешённых математических проблем XXI~века, описанных
Стивеном Смейлом~\cite{Math2005} по предложению 
занимавшего на рубеже XX и XXI~веков пост вице-президента Международного математического союза
В.И.~Арнольда. Смейл отобрал эти проблемы на основе следующих критериев: 1)~математически точная простая формулировка; 2)~личное знакомство с проблемой; 3)~вера
в то, что вклад в решение проблемы будет иметь большое значение для развития математики. В этом списке присутствует и гипотеза Римана, и уже доказанная гипотеза Пуанкаре, и проблема равенства P и NP.
\end{remark}

В 80-х годах XX~века Смейл, Ленор Блюм и Майк Шуб решили по-новому объединить теоретическую информатику
с численным анализом. Обеспечение основ теоретической информатики даёт классическая теория вычислений
Алана Тьюринга.
Но модель дискретных вычислений Тьюринга, в которой вычислимы
далеко не все вещественные числа~\cite{Kohlas_Turing} (ведь оперирует она только с целочисленным их
представлением) не является подходящей для анализа большинства современных компьютерных алгоритмов, 
которые предназначены для работы с вещественными числами.
В 1989~году Смейл, Блюм и Шуб представили свою модель вычислений
над произвольным кольцом или полем~$R$~\cite{BSS89, BSS_book}. 
В этой модели появилась возможность брать в расчёт
зависимость сложности вычислений от произвольной битовой
длины чисел (в отличие от модели алгебраической сложности,
для которой размерность входных данных зафиксирована\footnote{Т.е., если,
например, ставится задача перемножения матриц, то нужно отдельно
исследовать алгоритмы для {\bf каждой} комбинации размерностей
входных матриц, и в этом проблема модели алгебраической сложности.}). 

В случае, когда $R$ является кольцом из двух чисел $0$ и $1$,
модель Блюм--Шуба--Смейла (БШС) превращается в классическую
модель вычислений Тьюринга. Если $R$~--- поле вещественных чисел, то  
численные методы решения задач, оперирующие с такими числами, можно проанализировать в модели БШС. 

Опишем формально упомянутую в начале параграфа главную проблему линейного программирования,
поставленную Стивеном Смейлом.
Система $Ax \, \ge \, b$ имеет на входе вещественную матрицу $A$ размера $m \, \times \, n$ и вектор $b \, \in \, \mathbb{R}^m$. Существует ли полиномиальный
по числу арифметических операций алгоритм, который выполняет машина
над вещественными числами по модели БШС, определяющий,
существует ли такой $x \, \in \, \mathbb{R}^n$,
что $\sum\limits_{j=1}^n a_{ij} x_j \, \ge \, b_i$
для всех $i = 1, \, \ldots, \, m$? (подробнее см.~\cite[Глава~15]{BSS_book}).
Данная задача сводится к следующей задаче~ЛП:
даны матрица~$A$ размера $m \, \times \, n$, вектор~$b \, \in \, \mathbb{R}^m$ и вектор~$c \, \in \, \mathbb{R}^n$,
необходимо найти решение задачи 
$$
\max_{Ax \, \ge \, b} \, \langle c, \,  x \rangle,
$$
при $x \, \in \, \mathbb{R}^n$, если это решение существует. В терминах тьюринговой модели вычислений с использованием рациональных чисел $Q$ эта
проблема была решена Л.~Хачияном, но проблема над $\mathbb{R}$ остаётся открытой.
 
%
\subsubsection{Гипотеза Шуба--Смейла}\label{ch2_subsubsectSS_hyp}
К изложенному выше варианту теории сложности по Л.~Блюм, С.~Смейлу, М.~Шубу относится
и гипотеза Шуба--Смейла.
Оказывается, что эта гипотеза тесно связана со знаменитой проблемой математики и computer science о
равенстве или неравенстве классов $P$ и $NP$.
В случае, если гипотеза верна, Шубом
и Смейлом доказано неравенство $P$ и $NP$ над полем комплексных чисел.

Поясним, в чём состоит гипотеза Шуба--Смейла (далее в этом параграфе изложение частично
следует лекции Стивена Смейла, прочитанной
в Независимом Московском университете в 1999~году
и записанной в~\cite{Smail_mathpros}).

Рассматривается полином~$f(t)$ с целыми коэффициентами и последовательность $(1, \, t, \, u_1, \, \ldots, \, u_m = f)$, в которой каждый последующий член получается из некоторых двух предшествующих с помощью одной из трёх арифметических операций ($+$, $-$, $*$). Пусть число
$\tau(f)$ равно наименьшему возможному $m$.

В {\bf гипотезе Шуба--Смейла} утверждается, что {\it количество различных {\bf целых} корней многочлена $f$ не превосходит $\tau(f)^c$, где $c$~--- некоторая абсолютная константа}. 
Для комплексных и вещественных корней аналогичная
гипотеза неверна.

Рассмотрим многочлены $f_1(z_1, \, \ldots, \, z_n)$, 
$\ldots$, $f_k(z_1, \, \ldots, \, z_n)$ над полем
комплексных чисел~$\mathbb{C}$. Ответ на вопрос,
имеют ли эти многочлены общий нуль, даёт теорема
Гильберта о нулях (см. п.~\ref{ch1_subsect_Hilb}).
Общего нуля не существует тогда и только тогда, 
когда существуют многочлены $g_1, \, \ldots, \, g_k$ такие, что выполняется $\sum_{i = 1}^k g_i f_i = 1$.

Но теорема Гильберта о нулях даёт критерий существования общего нуля,
но не предлагает никакого алгоритма
нахождения $g_i$. Этот пробел
был закрыт в 1987~году Браунавеллом~\cite{BrBoundsNull}, который получил следующую неулучшаемую оценку
на степени многочленов $g_i$:
$$
\mbox{deg } g_i \leq \max \left \{3, \,
\max \deg f_i
\right \}^n.
$$
Таким образом, алгоритм сводится к решению системы линейных уравнений для коэффициентов многочленов~$g_i$.
Будем считать {\it временем} работы
алгоритма количество арифметических операций, а {\it размером} задачи~--- 
количество коэффициентов многочленов~$f_i$. В алгоритме Браунавелла время работы 
зависит от размера задачи экспоненциально.
Более того, по гипотезе Шуба--Смейла
полиномиального алгоритма определения,
существует ли общий нуль системы полиномиальных уравнений над $\mathbb{C}$, не существует в~принципе.

Здесь имеется в виду алгоритм не в смысле машины Тьюринга, а алгоритм над~$\mathbb{C}$ по модели БШС. Более конкретно, алгоритм~--- это ориентированный граф со входом и выходом
(в нашем случае на выходе один элемент,
который может принимать только два значения~--- <<да>> и <<нет>>). На вход
подаётся бесконечная в обе стороны последовательность комплексных чисел
$(\dots, \, 0, \, z_1 \neq 0, \, \dots, \, z_n \neq 0, \, 0, \, \dots)$.
В вычислительном узле графа производится
арифметическая операция с какими-то членами последовательности, и один член
последовательности заменяется на результат операции. Кроме того, можно умножить все члены последовательности на одно и то же число или произвести сдвиг последовательности. Ещё один тип узла графа~--- узел ветвления (проверка $z_i$ на
равенство 0). Граф может иметь циклы. Если называть {\it размером} входных данных число $n$, а {\it временем} работы
при данном входе ($T$)~--- длину пути от входа к выходу, то для алгоритма
с полиномиальным временем работы при всех входах должно
выполняться следующее неравенство:
$$
T \leq n^C,
$$
где $C$~--- некоторая константа.
После такого определения вопрос о том,
существует ли полиномиальный алгоритм 
для нашей задачи (об общем нуле системы
полиномиальных уравнений над $\mathbb{C}$), приобретает строгий математический смысл.

\subsection{О машинной арифметике}

Данный пункт подготовлен по материалам научно-популярной статьи Ю.\,В.~Матиясевича \cite{Mateyasevich}.

Таблицу умножения изучают с начальной школы. Даже простейший калькулятор легко перемножает шестизначные числа. Этого вполне достаточно не только для бытовых нужд, но и для большинства инженерных расчётов. В научных исследованиях приходится иметь дело и с более <<длинными>> числами~--- с десятками и сотнями десятичных знаков. Мы об этом не задумываемся, но в повседневной жизни есть области, в которых совершаются арифметические действия с такими числами. Для примера рассмотрим криптографию. Совершая покупку в Интернете, мы, естественно, не хотим, чтобы передаваемые нами данные банковской карты стали известны ещё кому-либо. Поэтому весь обмен информацией с банком шифруется. Многие другие сетевые сервисы в Интернете также используют более защищённый протокол HTTPS вместо обычного HTTP. А шифры часто основаны на больших простых числах и степень защиты информации напрямую зависит от величины используемых простых чисел. Ещё лет 10--15 назад простые числа, имеющие 128 цифр в двоичной записи, обеспечивали достаточную надёжность. Благодаря прогрессу как вычислительной техники, так и математических методов, на сегодняшний день такие коды успешно <<взламываются>>. Поэтому для надёжного шифрования приходится использовать существенно большие простые числа.

В настоящее время при работе с большими числами используют компьютеры. С~одной стороны, не так трудно <<обучить>> компьютер арифметическим действиям. С~другой стороны, важно, чтобы он выполнял эти операции не только правильно, но и быстро. В контексте оптимизационной тематики настоящего пособия заметим, что к примеру в задачах квадратичной оптимизации градиент целевой функции находится с помощью операции умножения матрицы на вектор. А эта операция естественным образом связана с необходимостью максимально быстро выполнять операцию умножения чисел.

Из четырёх арифметических действий сложение и вычитание представляются более простыми, чем умножение и деление. Но как сравнить их сложность (трудоёмкость) математически? Для этого можно посмотреть, как увеличивается количество выполняемых элементарных шагов с ростом длины чисел, над которыми производится действие. Шаг можно считать элементарным, если его трудность не зависит от длины чисел. Мерой трудоёмкости операции будем считать количество выполняемых элементарных шагов. Например:

1. Для операции сложения элементарным шагом можно считать \textbf{\emph{сложение двух цифр в одном разряде вместе с учётом переноса из предыдущего разряда и переноса в следующий разряд}} (в этом контексте зрительно можно представить сложение <<в столбик>>). Количество таких элементарных шагов пропорционально длине складываемых чисел. Скажем, при удвоении такой длины трудоёмкость всей операции сложения также удваивается.

2. Для операции умножения мерой трудоёмкости можно считать \textbf{\emph{количество обращений к таблице умножения для цифр}}. При <<школьном>> методе перемножения двух чисел <<в столбик>> количество обращений к таблице умножения будет равно произведению длин чисел (каждую цифру первого числа необходимо умножить на каждую цифру второго). При удвоении длин сомножителей необходимое количество обращений к таблице умножения увеличивается в 4 раза. Именно это отличие~--- возрастание трудоёмкости в 4, а не в 2 раза~--- делает этот алгоритм умножения намного более сложным при работе с <<длинными>> числами, чем алгоритм сложения. 

Однако стоит сказать, что умножение в столбик интересно, прежде всего, своей методической простотой. \ag{Известно} много других способов умножения натуральных чисел с лучшими оценками трудоёмкости (количества обращений к таблице умножения). Опишем один из таких методов, предложенный Анатолем Алексеевичем Карацубой в 1961 году. Поясним основную идею этого метода на примере умножения двух восьмизначных чисел $a$ и $b$. Представим их в виде
$$a=10^4 a_1+a_2,\;b=10^4 b_1+b_2,$$
где $a_1,\;a_2,\;b_1,\;b_2$~--- четырёхзначные числа. Тогда
$$a\times b = (10^4 a_1+a_2)(10^4 b_1+b_2)=$$
$$=10^8 a_1 b_1+10^4 (a_1 b_2+a_2 b_1)+a_2 b_2.$$
Заметим, что произведения $a_1 b_2$ и $a_2 b_1$ нам нужны не сами по себе, а только в сумме. Эту сумму можно вычислить, если мы уже знаем произведения $a_1 b_1$ и $a_2 b_2$, ценой только одного дополнительного перемножения двух четырёхзначных чисел и нескольких <<лёгких>> операций сложения/вычитания, поскольку
$$a_1 b_2+a_2 b_1=a_1 b_1+a_2 b_2-(a_1-a_2)(b_1-b_2).$$

\begin{example}
Пусть $a=73251846$ и $b=93725418$. Тогда
$$a_1=7325,\;a_2=1846,\;b_1=9372,\;b_2=5418.$$
Последовательно вычисляем:
$$a_1-a_2=7325-1846=5479,$$
$$b_1-b_2=9372-5418=3954,$$
$$a_1 b_1=7325\times9372=68649900,$$
$$a_2 b_2=1846\times5418=10001628,$$
$$(a_1-a_2)(b_1-b_2)=5479\times3954=21663966,$$
$$a_1 b_1+a_2 b_2=68649900+10001628=78651528.$$
По приведённой выше формуле произведение $a\times b$ равно
$$6864990000000000+786515280000+10001628=6865776525281628.$$
\end{example}

Сравним трудоёмкость метода А.А. Карацубы с трудоёмкостью традиционного метода <<в столбик>>, который изучают в школе.

1. При использовании школьного метода к таблице умножения требуется обратиться $8\times8=64$ раза.

2. В вычислениях методом А.А. Карацубы придётся найти три произведения четырёхзначных чисел. Если это делать школьным методом, то к таблице умножения придётся обратиться $3\times4\times4=48$ раз (умножение на степени числа 10~--- это просто дописывание нулей).

Но ведь <<по-новому>> можно перемножать и четырёхзначные числа, разбив каждое из них на два двузначных. По аналогии с разобранным выше примером, для перемножения четырёхзначных чисел достаточно найти три произведения двузначных чисел. Если эти произведения вычислять <<по-школьному>>, то потребуется $3\times(3\times2\times2)=36$ обращений к таблице умножения. Ещё одно применение той же идеи уменьшит количество обращений к таблице умножения до $3^3=27$. В итоге получаем, что для перемножения $8$-значных натуральных чисел достаточно 27 элементарных шагов вместо~64!

Можно показать, что в общем случае при умножении методом А.А. Карацубы двух $2^m$-значных чисел требуется провести $3^m$ обращений к таблице умножения. Следовательно, при удвоении длин сомножителей трудоёмкость увеличивается лишь в 3 раза, а не в 4, как при <<школьном>> способе. Общий выигрыш будет тем больше, чем длиннее перемножаемые числа.

Выше длина числа измерялась по его представлению в десятичной системе счисления. При переходе к двоичному представлению длина записи чисел возрастает примерно в 3,3 раза. По этой причине в двоичной системе счисления, в которой работают современные компьютеры, преимущества метода А.А. Карацубы проявляются, начиная с меньших по значению чисел, чем в десятичной системе. 

В 1961 году метод Карацубы был революционным прорывом. После того, как стало ясно, что школьный метод не оптимален, многие математики задумались~--- а нельзя ли перемножать большие числа ещё быстрее, чем методом Карацубы? Оказалось, что можно. Один из путей ускорения~--- разбиение каждого сомножителя не на две части, а на большее (фиксированное или растущее с длиной сомножителей) количество частей. К настоящему времени трудоёмкость перемножения многозначных чисел почти сравнялась с трудоёмкостью сложения. А именно, два $n$-значных числа можно перемножить, выполнив $O(n\log n)$ элементарных операций. Правдоподобной считается гипотеза, что два $n$-значных числа нельзя перемножить быстрее, чем за $Cn\log(n)$ шагов, где $C$~--- постоянная. Но пока никому это доказать не удалось. Такие результаты, называемые {\it нижними оценками} ({\it сложности вычисления}), обычно гораздо труднее верхних оценок. Для получения верхней оценки достаточно предъявить конкретный алгоритм и оценить его трудоёмкость, в то время как для получения нижней оценки нужно суметь обозреть все мыслимые способы вычисления искомой величины.


Отметим, что операция деления натуральных чисел не сильно усложняется по сравнению с умножением. Тогда оказывается, что его можно преобразовать в метод для деления чисел такой же длины, который будет требовать не более $5M(n)$ операций обращения к таблице умножения. Как говорилось ранее, сложение и вычитание~--- самые простые среди основных арифметических действий, быстро умножать мы также уже научились. Приведём для полноты \nk{изложения} схематическое описание аналога метода А.А. Карацубы для деления натуральных чисел. Следует различать две постановки задачи:

--- деление одного целого числа на другое с получением неполного частного и остатка;

--- получение приближённого значения отношения двух чисел, целых или десятичных.

Рассмотрим второй вариант. Прежде всего, заметим, что общую операцию нахождения частного $c=a/b$ можно свести к специальному случаю~--- нахождению обратного числа $d=1/b$. После этого искомый ответ получается как результат (быстрого) перемножения: $c=a\times d$. 

В свою очередь, нахождение обратного к $b$ числа означает необходимость решить уравнение $by=1$. Предположим, что выбрано начальное приближение~--- число $y_0$, примерно равное решению этого уравнения. Мы увидим, сколь точно это приближение, когда вычислим произведение $b\times y_0$. Оно должно быть близко к $1$, что можно представить формулой
$$b\times y_0=1-\varepsilon,$$
где $\varepsilon$~--- некоторое число с малой абсолютной величиной. Умножив обе части этого равенства на $1+\varepsilon$, получим
$$by_0\times (1+\varepsilon)=1-\varepsilon^2.$$
Это равенство можно интерпретировать так: если $|\varepsilon|<1$, то число
$$y_1=y_0(1+\varepsilon)$$
является лучшим приближением к $1/b$ по сравнению $y_0$ в силу
$\varepsilon^2<|\varepsilon|$. Заметим, $y_1$ можно найти также формуле $y_1=y_0(2-b y_0)$ без умножения на $1+\varepsilon$.
Процесс улучшения приближения можно продолжить, вычисляя последовательно
$$y_2=y_1(2-b y_1), \;\;\; y_3=y_2 (2-by_2), \text{ и т. д.}
$$
Легко проверить, что для всякого натурального числа $k$ верно
$$b\times y_k=1-\varepsilon^{2^k}.$$
Поэтому по мере увеличения $k$ получаются всё более и более точные приближения $y_k$ к $1/b$. Покажем пример.

\begin{example} Пусть $b=2,8062865$ и выбрано начальное приближение $y_0=0,4$. Последовательно находим:
$$y_1=0,40\times (2-2,8\times0,40)=0,3520,$$
$$y_2=0,3520\times (2-2,806\times 0,3520)=0,356325376,$$
$$y_3=0,356325376\times (2-2,8062865\times 0,356325376)=0,35634280\ldots$$
После трёх итераций мы получили 8 верных цифр числа $1/b$. На первом шаге использовались лишь 2 значащие цифры чисел $b$ и $y_0$, на втором шаге~--- лишь 4 значащие цифры чисел $b$ и $y_1$, и только на третьем шаге вычисления проводились с восьмизначными числами. Большая точность на начальных шагах не нужна. Длина используемых чисел должна соответствовать точности получаемых приближений. В итоге трудоёмкость последнего шага в таком итерационном процессе оказывается большей, чем трудоёмкость всех предыдущих шагов, вместе взятых. 
\end{example}

\begin{exercise}
С использованием рассуждений в конце предыдущего примера покажите, что работа по нахождению частного $a/b$ для $n$-значных (для наглядности можно положить $n = 2^q$ при некотором натуральном $q$) натуральных чисел оказывается не более сложной, чем $5$ операций умножения $n$-значных натуральных чисел.
\end{exercise}

Аналогичным образом осуществляется быстрое деление целых чисел с остатком.

\begin{exercise} По аналогии с описанными выше рассуждениями необходимо свести решение уравнения $by^2=1$ к последовательности выполнения операций умножения. Умножение полученного решения на $b$ даст, как легко понять, $\sqrt{b}$. \end{exercise}

\subsection[Аппроксимация трансцендентных условий\\ полуопределённо представимыми]{Аппроксимация трансцендентных условий\\ полуопределённо представимыми}

Некоторые задачи принципиально невозможно переписать в виде конической программы над симметричными конусами. Ясно, например, что трансцендентные ограничения невозможно представить линейными, конично-квадратичными или полуопределёнными условиями, поскольку все симметричные конуса задаются алгебраическими уравнениями. Однако точно так же, как конус Лоренца, можно аппроксимировать проекцией полиэдрального конуса (см.~п.~\ref{ch1_relaxations}), некоторые трансцендентные ограничения аппроксимируются условиями над конусами Лоренца или матричными конусами.

Ниже мы представим конструкцию, с помощью которой можно аппроксимировать ограничение $y \geq e^x$ конично-квадратичными условиями. Эта конструкция была предложена А.С.~Немировским \cite{Nemir_conrut}.
Пусть $p(x) = 1 + x + \frac{x^2}{2}$ обозначает квадратичный полином Тейлора функции $f(x) = e^x$ вокруг точки $x_0 = 0$. Тогда имеем
\[ e^x = \lim_{r \to \infty} \left( p\left(\frac{x}{2^r}\right) \right)^{2^r}.
\]
Сходимость не является равномерной по $x$, однако она равномерная и достаточно быстрая на компактных подмножествах (интервалах) ${I_T = \{ x \mid |x| < T \}}$. Чтобы аппроксимировать экспоненту с относительной точностью $\varepsilon$ на $I_T$, достаточно взять $r = O(T + \log\varepsilon^{-1})$. Но для каждого данного $r$, неравенство
\[ y \geq \left( p\left(\frac{x}{2^r}\right) \right)^{2^r}
\]
эквивалентно системе неравенств
\[ y \geq y_1^2,\quad y_1 \geq y_2^2,\quad\dots,\quad y_{r-1} \geq y_r^2,\quad y_r \geq p\left(\frac{x}{2^r}\right).
\]
Каждое из этих неравенств выпуклое и квадратичное, и поэтому представляется коническим условием над конусом Лоренца $L^3$ (см.~п.~\ref{ch1_sect_cone_repr}). В итоге мы получаем совокупность конических условий:
\[ \begin{pmatrix} 2y_1 \\ y-1 \\ y+1 \end{pmatrix},\ \begin{pmatrix} 2y_2 \\ y_1-1 \\ y_1+1 \end{pmatrix},\ \dots,\ \begin{pmatrix} 2y_r \\ y_{r-1}-1 \\ y_{r-1}+1 \end{pmatrix},\ \begin{pmatrix} 1+2^{-r}x \\ y_r-1 \\ y_r \end{pmatrix} \in L^3
\]
с $r$ дополнительными переменными. Для достижения относительной точности $10^{-12}$ для всех $|x| \leq 10$ достаточно взять $r = 24$, для достижения этой же точности для всех $|x| \leq 100$ достаточно взять $r = 29$.

\medskip

Похожие полуопределённые аппроксимации существуют и для других трансцендентных функций, например, (матричного) логарифма \cite{FawziSaundersonParrilo19} и эллиптического интеграла \cite{GlineurAGM}.


\chapter*{Заключение\markboth{Заключение}{Заключение}}\addcontentsline{toc}{chapter}{\hspace*{5.4mm}Заключение\dotfill}

Как уже отмечалось в аннотации к пособию, в образовании студентов-математиков МФТИ вот уже 50 лет особую роль играет курс <<\nk{Оптимизация}>>. В последние десятилетия роль оптимизации ещё больше возросла. Однако заметно сместились акценты. Если раньше большой акцент делался на решении задач управления, теоретико-игоровых задач, то современная оптимизация, конечно, в большой степени вдохновляется анализом данных. 

Важно подчеркнуть, что современная оптимизация и оптимизация \ag{70-х годов XX века}  довольно сильно отличаются не только основными приложениями, но и вычислительным аппаратом. Однако учебные пособия, которые доступны на русском языке, за редким исключением, во многом состоят из описания состояния развития данной области того времени. В~\nk{данной книге} предпринята попытка устранить отмеченный зазор. Как следствие, пособие местами может быть достаточно сложным для изучения. Заметная часть отобранных  материалов, по-видимому, излагается в учебной литературе впервые.

Однако не следует рассматривать данное пособие как исчерпывающий источник знаний по тому, какие задачи оптимизации сейчас решают и какими методами. Например, в пособии практически  ничего не говорится об очень популярных сейчас методах стохастического градиентного спуска и их приложениях к обучению глубоких нейронных сетей. Ничего не говорится про бурно развивающееся направление <<Распределённая оптимизация>>, которое играет чрезвычайно важную роль в решении современных задач оптимизации больших размеров. Мало место уделено невыпуклой оптимизации, а большинство практических задач оказываются невыпуклыми… С~другой стороны, все отмеченные и многие другие направления современной оптимизации, в том числе те, которые сильно востребованы на практике, едва ли возможно качественно освоить, не изучив {\it Выпуклую оптимизацию}. Мы~очень надеемся, что данное пособие станет достаточно хорошей стартовой площадкой, начиная с которой заинтересованные студенты в дальнейшем смогут  изучить более специальные разделы современной оптимизации, например, по книге \cite{Gas_Pos18} и цитированной там литературе.	
\clearpage

\addcontentsline{toc}{chapter}{\hspace*{5.4mm}Литература\dotfill}

\clearpage

\thispagestyle{empty}

{ \center \large  Учебное издание
	
	\vspace{26mm}

	\vspace{5mm}
	
	\textbf{Воронцова} Евгения Алексеевна\\ \textbf{Хильдебранд} Роланд Фалькович\\  \textbf{Гасников} Александр Владимирович\\  \textbf{Стонякин} Федор Сергеевич

	\vspace{15mm}

	{\bf\large ВЫПУКЛАЯ~ОПТИМИЗАЦИЯ%
	
 }

}

\vfill

{\parindent=0mm \small
	
	\small {Редакторы: \emph{В.\,А.~Дружинина}, \emph{И.\,А.~Волкова}, \emph{О.\,П.~Котова}. Корректор \emph{Н.\,Е.~Кобзева}}
	
	\small {Компьютерная верстка  \emph{Н.\,Е.~Кобзева}}
	
	\small {Дизайн обложки  \emph{Е.\,А.~Казённова}}

	\vspace{5mm}

	Подписано в печать 11.06.2021. Формат 60$\times$84 $^{1}\!/\!_{16}$. 
	
	Усл. печ. л.  22{,}75. Уч.-изд. л. 20{,}3. Тираж \ag{25}0~экз. Заказ №\,58.
	
	\vspace{5mm}

	Федеральное государственное автономное образовательное учреждение
	
	высшего  образования <<Московский физико-технический институт  
	
	(национальный исследовательский университет)>>
	
	141700, Московская обл., г. Долгопрудный, Институтский пер., 9
	
	Тел. (495) 408-58-22, e-mail: rio@mipt.ru
	
	\rule[2pt]{\textwidth}{0.2pt}
	
	Отпечатано в полном соответствии с предоставленным оригиналом-макетом\\
	ООО <<Печатный салон ШАНС>> \\
	127412, г.\,Москва, ул.\,Ижорская, д.\,13, стр.\,2. Тел. (495) 484-26-55


}


\end{document}